\documentclass[11.5pt]{article}
\usepackage{amsmath,amssymb}
\usepackage{amsthm}
\usepackage{multirow}
\usepackage{graphicx}
\usepackage{color}
\renewcommand{\theequation}{\thesection.\arabic{equation}}
\newtheorem{thm}{Theorem}[section] 
\newtheorem{lemma}{Lemma}[section] 
\newtheorem{cond}{Condition}[section]

\newtheorem{definition}{Definition}[section]
\newcommand{\bed}{\begin{definition}}
\newcommand{\eed}{\end{definition}}

\newcommand{\rom}[1]{\uppercase\expandafter{\romannumeral #1\relax}}

\newcommand{\eps}{\epsilon}

\newcommand{\bitem}{\begin{itemize}}
\newcommand{\eitem}{\end{itemize}}

\newcommand{\goto}{\rightarrow}

\newcommand{\beqn}{\begin{equation}}
\newcommand{\eeqn}{\end{equation}}
\newcommand{\balign}{\begin{align}}
\newcommand{\ealign}{\end{align}}

\newcommand{\s}{\sigma}

\newcommand{\tr}{\mathrm{tr}}
\usepackage{amssymb}
\newcommand{\beq}{\begin{equation}}
\newcommand{\eeq}{\end{equation}}

\newcommand{\diag}{\mathrm{diag}}

\newcommand{\teta}{\tilde{\eta}} 

\newcommand{\heta}{\hat\eta}
\allowdisplaybreaks

\newcommand\numberthis{\addtocounter{equation}{1}\tag{\theequation}}

\title{Optimal Adaptivity of Signed-Polygon Statistics for Network Testing}
\begin{document}
\author{Jiashun Jin$^*$, Zheng Tracy Ke$^\dagger$ and Shengming Luo$^*$ \\
Department of Statistics \\
Carnegie Mellon University$^*$ and Harvard University$^\dagger$} 
\maketitle

\begin{abstract} 
Given a symmetric social network, we are interested in testing whether it has only one community or multiple  communities.  The desired tests should (a) accommodate severe degree heterogeneity, (b) 
accommodate mixed-memberships, (c) have a tractable null distribution, and (d)  adapt automatically to different levels of sparsity, and achieve the optimal phase diagram.  
How to find such a test is a challenging problem. 

We propose the Signed Polygon as a class of new tests.  
Fixing $m \geq 3$, for each $m$-gon in the network,  define a score using the centered 
adjacency matrix. The sum of such scores is then the $m$-th order Signed Polygon statistic. The Signed Triangle (SgnT) and the Signed Quadrilateral (SgnQ) 
are special examples of the Signed Polygon.  

We show that both the SgnT and SgnQ tests satisfy (a)-(d), and especially, they work well for both very sparse and less sparse networks. Our proposed tests compare favorably with the existing tests. For example, the EZ and GC tests  behave unsatisfactorily in the less sparse case and do not achieve the optimal phase diagram. Also, many existing tests do not allow for severe heterogeneity or mixed-memberships, and they behave unsatisfactorily in our settings. 

The analysis of the SgnT and SgnQ tests is delicate and extremely tedious, and the main reason is that we need  a unified proof that covers  a wide range of sparsity levels and a wide range of degree heterogeneity. For lower bound theory, we use a phase transition framework, which includes the standard minimax argument, but is more informative.  The proof uses classical theorems on matrix scaling.

\end{abstract}

\section{Introduction} \label{sec:intro}  
Given a symmetrical social network, we are interested in the {\it global testing problem} where we use the adjacency matrix of the network to test whether it has only one community or multiple communities.  A good understanding of the problem is useful for discovering non-obvious social groups and patterns \cite{pattern2, pattern1},  measuring diversity of individual nodes \cite{diversity},   determining stopping time in a recursive community detection scheme \cite{li2018hierarchical, zhao2011community}.  It may also help  understand other related problems such as membership estimation \cite{JiZhuMM}, and estimating the number of the communities \cite{saldana2015many,wang2017likelihood}.    
  
Natural networks have several characteristics that are ubiquitously found: 
\begin{itemize} 
\item {\it Severe degree heterogeneity}. The distribution of the node degrees   usually has a power-law tail, implying severe degree heterogeneity.  
\item {\it Mixed-memberships}.  Communities are tightly woven clusters of nodes where we have more edges within than between \cite{GirvanNewman, Radicchi2658}.  Communities are rarely non-overlapping, and some nodes may belong to more than one community (and thus have mixed-memberships). 
\item {\it Sparsity}. Many networks are sparse. The sparsity levels may range significantly from one network to another, and may also range significantly from one node to another (due to severe degree heterogeneity).  
\end{itemize} 
Phase transition is a well-known optimality framework \cite{DJ04, Ingster2010, MW, Paul}. It is related to the minimax framework but can be more informative in many cases. 
Conceptually, for the global testing problem, in the two-dimensional phase space with the two axes calibrating the ``sparsity" and ``signal strength", respectively,  there is a ``Region of Possibility" and a ``Region of Impossibility". In  ``Region of Possibility", any alternative is separable from the null. In 
  ``Region of Impossibility", any alternative is inseparable from the null. 

If  a  test is able to automatically adapt to different levels of sparsity and is able to separate any given alternative in  the ``Region of Possibility" from the null, then we call it ``optimally adaptive".   

We are interested in finding tests that satisfy the following requirements. 
\begin{itemize}  
\item (R1). Applicable to networks with severe degree heterogeneity. 
\item (R2). Applicable to networks with mixed-memberships. 
\item (R3). The asymptotic null distribution is easy to track, so the rejection regions are  
easy to set.  
\item (R4). Optimally adaptive:  We desire a single test that is able to adapt to different levels of sparsity and is optimally adaptive. 
\end{itemize} 

\subsection{The DCMM model} 
\label{subsec:DCMM} 
We adopt the {\it Degree Corrected Mixed Membership (DCMM)} model  \cite{JiZhuMM, mixed-SCORE}. Denote the adjacency matrix by $A$, where 
\begin{equation} \label{model1a} 
A_{ij} = 
\left\{
\begin{array}{ll} 
1, &\qquad  \mbox{if node $i$ and node $j$ have an edge},  \\
0, &\qquad  \mbox{otherwise}.   \\
\end{array}
\right.  
\end{equation} 
Conventionally,  self-edges are not allowed  so all the diagonal entries of $A$ are $0$.   
In DCMM, we assume there are $K$ perceivable communities ${\cal C}_1, {\cal C}_2, \ldots, {\cal C}_K$,  and  each node is associated with a mixed-membership weight vector $\pi_i = (\pi_i(1),\pi_i(2),\ldots,\pi_i(K))'$ where for $1 \leq k \leq K$ and $1 \leq i \leq n$, 
\begin{equation} \label{model1b} 
\pi_i(k) = \mbox{the weight node $i$ puts in community $k$}.    
\end{equation}  
Moreover, for a $K \times K$ symmetric nonnegative matrix $P$ which models the community structure,  and positive parameters $\theta_1, \theta_2, \ldots, \theta_n$ which model  the degree heterogeneity,  
we assume the upper triangular entries of $A$ are independent Bernoulli variables satisfying    
\begin{equation} \label{model1c} 
\mathbb{P}(A_{ij}=1) = \theta_i\theta_j\cdot \pi_i'P\pi_j \equiv \Omega_{ij},   \qquad 1 \leq i < j \leq n,   
\end{equation} 
where $\Omega$ denotes the matrix $\Theta \Pi P \Pi' \Theta$, with $\Theta$ being the $n \times n$ diagonal matrix $\diag(\theta_1, \ldots, \theta_n)$ and $\Pi$ being the $n \times K$ matrix $[\pi_1, \pi_2, \ldots, \pi_n]'$. For identifiability (see \cite{mixed-SCORE} for more discussion), we assume 
\begin{equation} \label{condition1d} 
\mbox{all diagonal entries of $P$ are $1$}.   
\end{equation} 
When $K = 1$,  (\ref{condition1d}) implies $P = 1$,  and so  
$\Omega_{ij} = \theta_i \theta_j$, $1 \leq i, j \leq n$.   

Write for short $\diag(\Omega) = \diag(\Omega_{11}, \Omega_{22}, \ldots, \Omega_{nn})$, and let $W$ be the matrix where 
for $1 \leq i, j \leq n$,  $W_{ij} = A_{ij} -  \Omega_{ij}$ if $i \neq j$ and $W_{ij} = 0$ otherwise. In matrix form, we have    
\begin{equation} \label{DCMM-matrixform} 
A = \Omega - \diag(\Omega) + W,   \qquad \mbox{where} \;  \Omega = \Theta \Pi P \Pi' \Theta. 
\end{equation}

DCMM includes three models as special cases, each of which is well-known and has been 
studied extensively recently. 
\begin{itemize} 
\item {\it Degree Corrected Block Model (DCBM)} \cite{DCBM}.  If we do not allow mixed-memberships (i.e.,  each weight vector $\pi_i$ is degenerate with one entry being nonzero), then DCMM reduces to the DCBM.  
\item {\it Mixed Membership Stochastic Block Model (MMSBM)} \cite{airoldi2009mixed}.  If $\theta_1 = \theta_2 = \ldots = \theta_n$ and we denote the common value by $\alpha_n$, then $\Omega$ reduces to 
$\Omega = \alpha_n \Pi  P \Pi'$.  
For identifiability in this special case, (\ref{condition1d}) is too strong, so we relax it to that the average of the diagonals of $P$ is $1$.  
\item {\it Stochastic Block Model (SBM)} \cite{holland1983stochastic}.  MMSBM further reduces to the classical SBM  if additionally we do not allow mixed-memberships.   
\end{itemize} 
It is instructive to consider a special DCMM model, which is a  generalization from the symmetric SBM \cite{mossel2015reconstruction} to the case with degree heterogeneity. 

{\bf Example 1} {\it (A special DCMM)}. Let $e_1,\ldots,e_K$ be  the standard basis of $\mathbb{R}^K$. 
Fixing a positive vector $\theta\in\mathbb{R}^n$ and a scalar $b_n\in (0,1)$, we assume
\beq \label{example1} 
P = (1 - b_n) I_K + b_n 1_K 1_K',  
\qquad \mbox{$\pi_i$ are iid sampled from $e_1, \ldots, e_K$}.  
\eeq 
In this model,  $\|\theta\| $ measures the sparsity level, and $(1 - b_n)$ measures the ``similarity" between different communities.

Under DCMM,  the global testing problem is the problem of testing 
\beq \label{Problem}
H_0^{(n)}:  K=1   \qquad  \mbox{vs}. \qquad H_1^{(n)}: K\geq 2.  
\eeq
The seeming simplicity of the two hypotheses is deceiving, as both of them are highly 
composite, consisting of many different parameter configurations.  
\subsection{Phase transition: a preview of our main results} 
\label{subsec:preview} 
Let $\lambda_1, \lambda_2, \ldots, \lambda_K$ be the first $K$  eigenvalues of $\Omega$, arranged in the descending order in magnitude. We can view (a) $\sqrt{\lambda_1}$ both as the sparsity level and the noise level  \cite{SCORE} (i.e., spectral norm of the noise matrix $W$),  (b) $|\lambda_2|$ as the signal strength, so $|\lambda_2| / \sqrt{\lambda_1}$ is the Signal-to-Noise Ratio (SNR), and (c) $|\lambda_2|/\lambda_1$ as a measure for the similarity between different communities.

\begin{figure}[tb!]
\centering 
\includegraphics[width=.442\textwidth]{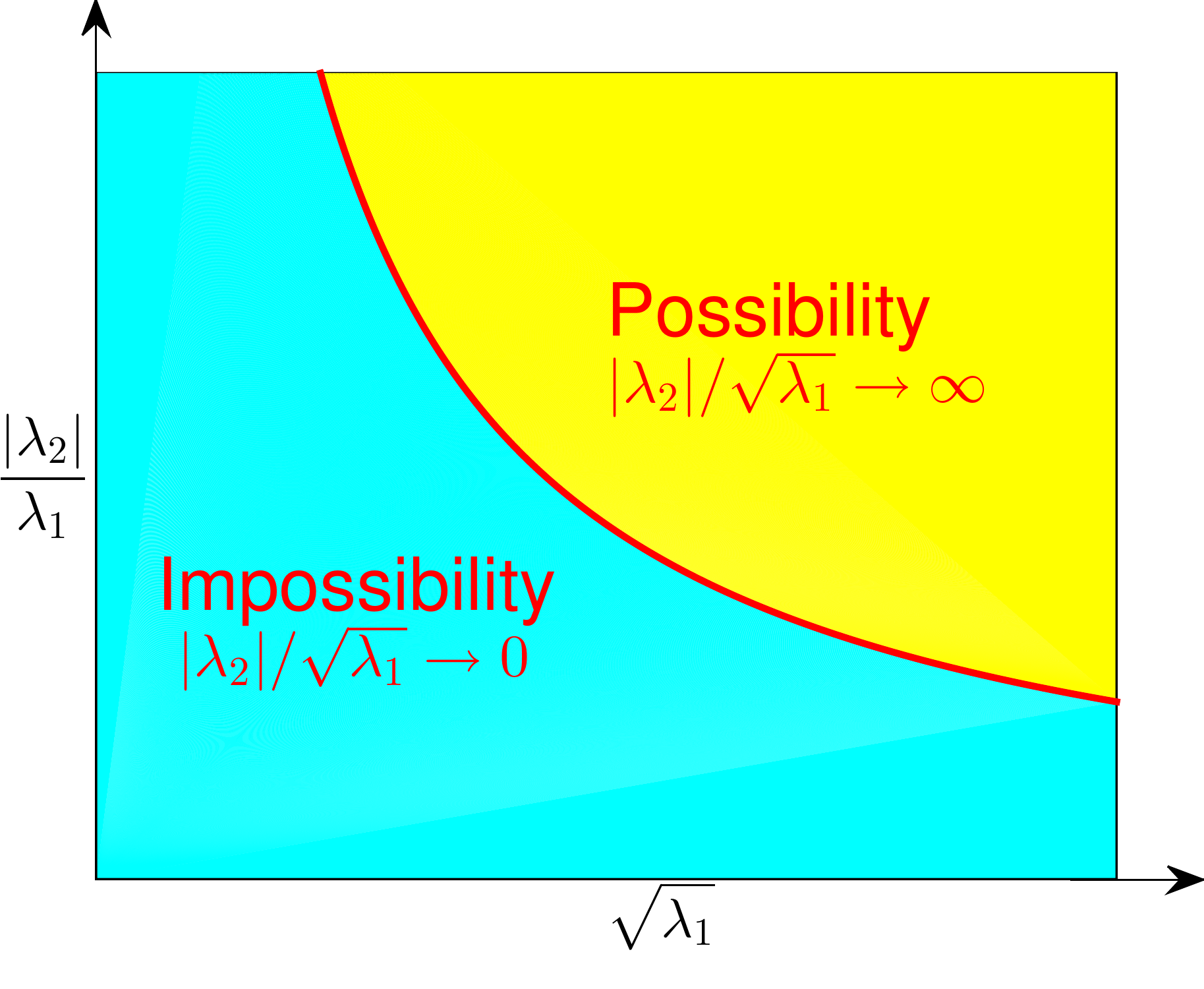} $\qquad$
\includegraphics[width=.442\textwidth]{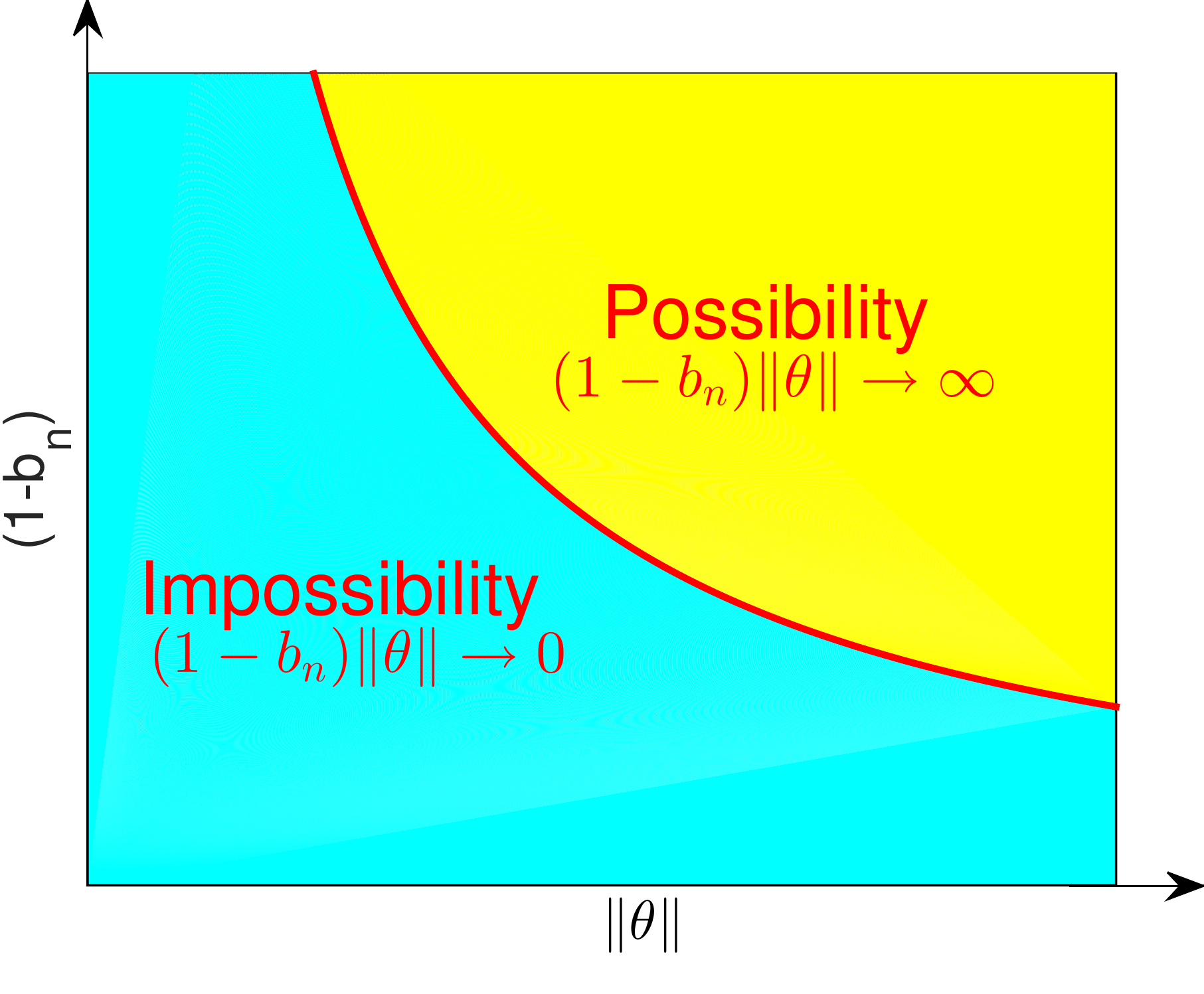}
\caption{Left: Phase transition. In Region of Impossibility, any alternative hypothesis is indistinguishable 
from a null hypothesis, provided that some mild conditions hold.  In Region of Possibility, the Signed Polygon test is able to 
separate any alternative hypothesis from a null hypothesis asymptotically. Right: Phase transition for 
the special DCMM model in Example 1, where $\sqrt{\lambda_1} \asymp \|\theta\|$, $|\lambda_2|/\lambda_1 \asymp (1 - b_n)$, and $|\lambda_2| / \sqrt{\lambda_1} \asymp (1 - b_n)\|\theta\| $.} 
\label{figure:phase1}
\end{figure}

Now, in the two-dimensional {\it phase space} where the $x$-axis is $\sqrt{\lambda_1}$ which measures the sparsity level, and the $y$-axis is $|\lambda_2|/\lambda_1$ which measures the community similarity, we have two regions. 
\begin{itemize} 
\item {\it Region of Possibility ($1 \ll \sqrt{\lambda_1} \ll \sqrt{n},\; |\lambda_2| / \sqrt{\lambda_1} \goto \infty$)}. For any alternative hypothesis in this region, it is possible to distinguish it from any null hypothesis, by the Signed Polygon tests to be introduced. 
\item {\it Region of Impossibility ($1 \ll \sqrt{\lambda_1} \ll \sqrt{n}, \;  |\lambda_2| / \sqrt{\lambda_1} \goto 0$)}. In this region, any alternative hypothesis is inseparable from the null hypothesis, provided that some mild conditions hold.  
\end{itemize} 
See Figure \ref{figure:phase1} (left panel).  The Signed Polygon test satisfies all requirements (R1)-(R4) aforementioned.  
Since the test is able to separate all alternatives (ranging from very sparse to less sparse)  in the  Region of Possibility from the null,   it is {\it optimally adaptive}.

To further elucidate, consider the special DCMM in Example 1, where    
\[
\lambda_1 \sim (1 + (K-1) b_n) \|\theta\|^2, \qquad \lambda_k \sim   (1 - b_n)\|\theta\|^2, \;\; k = 2, 3, \ldots, K. 
\] 
The sparsity level is $\sqrt{\lambda}_1 \asymp \|\theta\|$, and the SNR is $|\lambda_2| / 
\sqrt{\lambda_1}  \asymp \|\theta\|(1 - b_n)$, where $(1 - b_n)$ measures the community similarity.   
In this example,  the Region of Possibility and Region of Impossibility are defined by 
\[
\{1 \ll \|\theta\| \ll \sqrt{n}, \|\theta\| (1 - b_n) \goto \infty\}, \;  \mbox{and}\; \{1 \ll \|\theta\| \ll \sqrt{n}, \|\theta\|(1 - b_n) \goto 0\}, 
\] 
respectively. See Figure \ref{figure:phase1} (right panel).  

{\bf Remark 1}. As the phase transition is hinged on $\lambda_2/\sqrt{\lambda_1}$,  one may think 
that the statistic $\hat{\lambda}_2/ \sqrt{\hat{\lambda}_1}$ is optimally adaptive, where $\hat{\lambda}_k$ 
is the $k$-th eigenvalue of $A$, $1 \leq k \leq K$,  arranged in the descending order in magnitude. 
This is however not true, for the consistency of $\hat{\lambda}_2$  to $\lambda_2$
can not be guaranteed in our range of interest, unless with strong conditions on $\theta_{max}$ \cite{SCORE}.

\subsection{Literature review, the Signed Polygon and our contribution, and four surprises} 
\label{subsec:review}
Recently, the global testing problem has attracted much attention and many interesting approaches have been proposed. 
To name a few, Mossel et al. \cite{mossel2015reconstruction} and Banerjee and Ma \cite{banerjee2017optimal} (see also \cite{banks2016information}) considered a special case of the testing problem, where they assume a simple null of Erdos-Renyi random graph model and a special alternative which is an SBM with two equal-sized communities. They provided the asymptotic distribution of the log-likelihood ratio within the contiguous regime. Since the likelihood ratio test statistic is NP-hard to compute, \cite{banerjee2017optimal}  
introduced an approximation by linear spectral statistics. Lei \cite{lei2016goodness} also considered the SBM model and studied the problem of testing whether $K = K_0$ or $K > K_0$, where $K_0$ is the number of communities.  His approach is based on the Tracy-Widom law of extreme eigenvalues and requires delicate random matrix theory.  
Unfortunately, these work have been focused on the SBM (which allows neither severe degree heterogeneity nor  mixed membership). 
Therefore, despite the elegant theory in these works, it remains unclear  how to extend their ideas to our settings.

The approach by Gao and Lafferty \cite{gao2017testing} is probably the first that tackles the global testing problem in settings that allow severe degree heterogeneity. They showed that the EZ test has a 
null that is asymptotically $N(0,1)$, and has competitive powers in 
many interesting settings. 
However, they only considered a relatively idealized setting where the off-diagonal entries of $P$ are all equal and where $(\theta_i, \pi_i)$'s are $iid$ generated (see details therein), and whether their ideas continue to work in our setting remains unclear.

Jin et al. \cite{OGC} considered the problem in   much broader settings, with a very different theoretical framework. 
They suggested  a general recipe for constructing test statistics 
that have $N(0,1)$ as the asymptotic null distribution, and proposed a class of test statistics called the {\it graphlet counting (GC)}, which includes the EZ test as a special case.  
They explained why both GC and EZ tests are reasonable ideas (for settings much broader than that  of \cite{gao2017testing}) and showed that both tests have competitive power  in many cases.

At that time, our hope was that the GC test is the desired test:   
in our manuscript \cite{Ke2019}, we tried very hard to analyze the GC test,  hoping that it satisfies  (R1)-(R4). 
Unfortunately, after substantial time and efforts,  we found that in the less sparse case, the variance of the GC test  becomes unsatisfactorily large, and so the test loses power in many easy-to-test scenarios and is not optimally adaptive.

Fortunately, right at the moment of despair, we came to realize that 
\begin{itemize} 
\item  Especially for the less sparse case,  the key to constructing a powerful test is not how to  capture the signal, but to reduce the variance. 
\item The GC test is based on counts of {\it non-centered} cycles/paths. The variance can be much smaller if we count the {\it centered}  cycles instead.    
\end{itemize} 
Centered and non-centered cycles are defined on the centered and non-centered adjacency matrix, respectively. See details below. 

These insights motivate a class of new tests which we call {\it Signed Polygon}, including the Signed Triangle (SgnT) and the Signed Quadrilateral (SgnQ). The Signed Polygon statistics are related to 
the Signed Cycle statistics, first introduced by Bubeck et al. \cite{bubeck2016testing} and later generalized by Banerjee \cite{banerjee2018contiguity}.   

The Signed Polygon and the Signed Cycle
are cycle-counting approaches, both of which recognize the benefit of variance reduction by counting centered cycles instead of non-centered cycles, but there are some major differences.   
The study of the Signed Cycles has been focused on the SBM and similar models, where  under the null,   $\mathbb{P}(A_{ij}= 1)  = \alpha$, $1 \leq i\neq j \leq n$,  and $\alpha$ is the only unknown parameter. In this case, a natural approach to centering the adjacency matrix $A$ is to first  estimate $\alpha$ using the whole matrix $A$ (say, $\hat{\alpha}$),  and then subtract all off-diagonal entries of $A$ by $\hat{\alpha}$.  However, under the null of our setting, $\mathbb{P}(A_{ij} = 1)  = \theta_i \theta_j$, $1 \leq i\neq j \leq n$, and there are $n$ different unknown parameters $\theta_1, \theta_2, \ldots, \theta_n$. In this case, how to center the matrix $A$ is not only unclear but also {\it  worrisome}, especially when the network is very sparse, because we have to use limited data to estimate a large number of  unknown parameters. Also, for any approaches we may have, the analysis is seen to be  much harder than that of  the previous case.

Note that the ways how two statistics are defined over the centered adjacency matrix are also different; see Section \ref{subsec:SP} and \cite{bubeck2016testing}.   

In the Signed Polygon, we use a new approach to estimate $\theta_1, \theta_2, \ldots, \theta_n$ under the null, and use the estimates to center the matrix $A$. 
To our surprise, data limitation (though a challenge) does not ruin the idea, and even for  very sparse networks, the estimation errors of $\theta_1, \theta_2, \ldots, \theta_n$ only have a negligible effect. 
The main contributions of the paper are as follows. 
\begin{itemize} 
\item Identify the Region of Impossibility and the Region of Possibility in the phase space.
\item Propose the Signed Polygon as a class of new tests that are appropriate for networks with severe degree heterogeneity and mixed-memberships. 
\item Prove that the Signed Triangle and  Signed Quadrilateral tests satisfy all the requirements (R1)-(R4), and especially that they are optimally adaptive and perform well for all networks in the Region of Possibility, ranging from very sparse ones to the least sparse ones. 
\end{itemize} 
To show the success of the Signed Polygon test for the whole Region of Possibility is very subtle and extremely tedious. The main  reason is that we hope to cover the {\it whole spectrum} of degree heterogeneity and sparsity levels. Crude bounds may work in one case but not another, and many seemingly negligible terms turn out to be non-negligible (see Sections~\ref{subsec:SP} and \ref{sec:SgnedPolygon}).  
The lower bound argument  is also very subtle. Compared to work on SBM where there is only 
one unknown parameter under the null, our null has $n$ unknown parameters. 
The difference provides a lot of freedom in constructing inseparable hypothesis pairs, and so the Region of Impossibility in our setting is much wider than that for SBM. Our construction of inseparable hypothesis pairs uses  theorems on non-negative matrix scaling, a mathematical area pioneered by Sinkhorn \cite{DADSinkhorn} and Olkin \cite{DADOlkin} among others (e.g., \cite{DADBrualdi, DADJohnson}).

The search for optimally adaptive tests  has been quite a journey, during which we have encountered  an array of {\it surprises}.  
\begin{itemize} 
\item {\it Which case is harder to construct a good test, the more sparse or less sparse?}  Surprisingly, the latter is harder: A good test needs not only to capture the signal, 
but also to have a small variance; in the less sparse case, the latter is more important than the former. The variances of the EZ and the GC test statistics are unexpectedly large in the less sparse case,  which is the main reason why they 
perform well in the more sparse case but 
unsatisfactorily in the less sparse case.  
\item {\it Are the estimation errors for $\theta$ worrisome?} The Signed Polygon relies on  an estimate of $\theta$. One may worry that the estimation errors are hard to control, especially for the very sparse  case, as we have insufficient data to estimate an $n$-dimensional vector $\theta$.  Surprisingly, while the effect of the  estimation error is not always negligible, it is quite tractable, and the resultant tests work,  even in the very sparse case.   
\item {\it Is it really difficult to have test statistics with tractable null distributions?} Surprisingly, no.   In \cite{Ke2019}, we have identified a wide family of graphlet counting statistics, including GC, 
Signed Polygon, and the $m$-cycles and length-$m$ paths introduced there. These test statistics are 
all asymptotically normal, where the asymptotic means and variances can be approximated by a simple function  of $\|\theta\|_q$, $q = 1, 2, 3, 4$.  
As $\|\theta\|_q$ can be  accurately  estimated under the null, we can normalize the statistics 
so that $N(0,1)$ is the limiting null distribution.   
\item {\it Comparison with linear spectral statistics}. The power of our tests hindges on 
$\tr(\widetilde{\Omega}^m)$ ($\widetilde{\Omega}$ is closely related to $\Omega$;  to be defined). 
One may think a good estimate for $\tr(\widetilde{\Omega}^m)$ is the linear spectral statistic $\tr((A-\hat{\eta}\hat{\eta}')^m)$ ($(A-\hat{\eta}\hat{\eta}')$: a stochastic counterpart of $\widetilde{\Omega}$; to be defined), but Random Matrix Theory suggests that this is not always the case. 
Surprisingly, the Signed Polygon provides a better estimate for $\tr(\widetilde{\Omega}^m)$, and the key is that, many linear spectral statistics can also be viewed 
as graphlet counting statistics, but they include degenerate $m$-gons in the counts and have unsatisfactory variances;  the Signed Polygon statistics do not  count degenerate $m$-gons 
and have more satisfactory variances. 
\end{itemize} 
We shall explain these points with more details in Section~\ref{sec:main}. 

\subsection{The Signed Polygon statistic} 
\label{subsec:SP} 
Recall that $A$ is the adjacency matrix of the network. Introduce 
a vector $\hat{\eta}$ by (${\bf 1}_n$ denotes the vector of $1$'s)
\begin{equation} \label{Definehateta} 
\hat{\eta} = (1/\sqrt{V})\, A {\bf 1}_n , \qquad \mbox{where $V = {\bf 1}_n' A {\bf 1}_n$}.  
\end{equation} 
Fixing $m \geq 3$, the order-$m$ {\it Signed Polygon} statistic is defined by (notation: $(dist)$ is short for ``distinct", which means any two of $i_1,\ldots,i_m$ are unequal)
\begin{equation} \label{U-stat} 
U_n^{(m)} = \sum_{i_1, i_2, \ldots, i_m (dist)} (A_{i_1 i_2} - \hat{\eta}_{i_1} \hat{\eta}_{i_2}) (A_{i_2 i_3} - \hat{\eta}_{i_2} \hat{\eta}_{i_3}) \ldots (A_{i_m i_1} - \hat{\eta}_{i_m} \hat{\eta}_{i_1}). 
\end{equation} 
When $m = 3$, we call it the Signed-Triangle (SgnT) statistic:  
\begin{equation} \label{DefineST} 
T_n = \sum_{i_1,i_2,i_3 (dist)} (A_{i_1i_2} - \hat{\eta}_{i_1} \hat{\eta}_{i_2}) (A_{i_2 i_3} - \hat{\eta}_{i_2} \hat{\eta}_{i_3})  (A_{i_3 i_1} - \hat{\eta}_{i_3} \hat{\eta}_{i_1}).  
\end{equation}  
When $m = 4$, we call it the Signed-Quadrilateral (SgnQ) statistic: 
\begin{equation}  \label{DefienSQ} 
Q_n = \sum_{i_1,i_2,i_3,i_4 (dist)} (A_{i_1 i_2} - \heta_{i_1}  \heta_{i_2}) (A_{i_2i_3} - \heta_{i_2}  \heta_{i_3})  (A_{i_3i_4} - \heta_{i_3}  \heta_{i_4}) (A_{i_4 i_1} - \heta_{i_4} \heta_{i_1}).
\end{equation}  
For analysis, we focus  on $T_n$ and $Q_n$, but the theoretical framework is extendable to general $m$.  

The key to understanding and analyzing the Signed Polygon is the {\it Ideal Signed Polygon}.  
Introduce a {\it non-stochastic counterpart} of $\hat{\eta}$ by 
\begin{equation} \label{Defineeta*} 
\eta^* = \Omega {\bf 1}_n / \sqrt{v_0}, \qquad \mbox{where $v_0 = {\bf 1}_n' \Omega {\bf 1}_n$}.  
\end{equation} 
Define the order-$m$ {\it Ideal Signed Polygon} statistic by  
\begin{equation} \label{Ideal-SP} 
\widetilde{U}_n^{(m)} = \sum_{i_1, i_2, \ldots, i_m (dist)} (A_{i_1 i_2} - \eta_{i_1}^* \eta_{i_2}^*) (A_{i_2 i_3} - \eta_{i_2}^* \eta_{i_3}^*) \ldots (A_{i_m i_1} - \eta_{i_m}^*\eta_{i_1}^*). 
\end{equation} 
We expect to see that 
\[
\hat{\eta} \approx \mathbb{E}[\heta] \approx \eta^*.
\]   
We can view $\widetilde{U}_n^{(m)}$ as the oracle version of $U_n^{(m)}$, with $\eta^*$ given. We can also view $U_n^{(m)}$ as the {\it plug-in} version of $\widetilde{U}_n^{(m)}$, where we replace $\eta^*$ by $\hat{\eta}$.

For implementation, it is desirable to rewrite $T_n$ and $Q_n$ in  matrix forms, which allows us to avoid using a for loop and compute much faster (say, in MATLAB or R). 
For any two matrices $M, N \in \mathbb{R}^{n,n}$, let $\tr(M)$ be the trace of $M$,   $\diag(M) = \diag(M_{11}, M_{22}, \ldots, M_{nn})$, and $M \circ N$ be  
the Hadamard product of $M$ and $N$ (i.e., $M \circ N \in \mathbb{R}^{n, n}$, $(M \circ N)_{ij} = M_{ij} N_{ij}$). 
Denote $\widetilde{A} = A - \hat{\eta}\hat{\eta}'$.  
The following theorem is proved in the  appendix.  
\begin{thm} 
\label{thm:compu} 
We have $T_n = \tr(\widetilde{A}^3) - 3\tr(\widetilde{A}\circ \widetilde{A}^2) + 2\tr(\widetilde{A}\circ\widetilde{A}\circ\widetilde{A})$ and 
\begin{align*} 
Q_n  & =   \tr(\widetilde{A}^4) - 4\tr(\widetilde{A}\circ\widetilde{A}^3)  + 8\tr(\widetilde{A}\circ\widetilde{A}\circ\widetilde{A}^2)  - 6\tr(\widetilde{A}\circ\widetilde{A}\circ\widetilde{A}\circ\widetilde{A})   \\
 & - 2\tr(\widetilde{A}^2\circ\widetilde{A}^2) + 2\cdot 1_n'[\diag(\widetilde{A})(\widetilde{A} \circ \widetilde{A})\diag(\widetilde{A})]1_n 
 + 1_n'[ \widetilde{A}\circ\widetilde{A}\circ\widetilde{A}\circ\widetilde{A}]1_n.  
\end{align*}
The complexity of computing both $T_n$ and $Q_n$ is $O(n^2 \bar{d})$, where $\bar{d}$ is the 
average degree of the network.   
\end{thm}
Compared to the EZ and GC tests \cite{gao2017testing,OGC},  the computational complexity of SgnT and SgnQ is of the same order. 
  

{\bf Remark 2}. The computational complexity of  $U_n^{(m)}$ remains as $O(n^2 \bar{d})$ for larger  $m$. Similarly as that in Theorem~\ref{thm:compu},   the main complexity of $U_n^{(m)}$ comes from computing $\widetilde{A}^m$. Since we can compute $\widetilde{A}^m$ 
with $\widetilde{A}^m = \widetilde{A}^{m-1} \widetilde{A}$ 
and recursive matrix multiplications, each time with a complexity of $O(n^2 \bar{d})$,  the overall complexity is $O(n^2 \bar{d})$.

{\bf Remark 3} ({\it Connection to the Signed Cycle}).  In the more idealized MMSBM or SBM model,  we don't have degree heterogeneity, and 
$\Omega = \alpha_n {\bf 1}_n {\bf 1}_n'$ under the null, where $\alpha_n$ is the only unknown parameter. 
In this simple setting, it makes sense to estimate $\alpha_n$ by 
$\hat{\alpha}_n = \bar{d} /(n-1)$, where $\bar{d}$ is the average degree. This gives rise to the {\it Signed Cycle} statistics \cite{banerjee2018contiguity, bubeck2016testing}: 
\[
C_n^{(m)} = \sum_{i_1, i_2, \ldots, i_m (dist)} (A_{i_1 i_2} - \hat{\alpha}_n) (A_{i_2 i_3} - \hat{\alpha}_n) \ldots (A_{i_m i_1} - \hat{\alpha}_n). 
\] 
Bubeck et al. \cite{bubeck2016testing} first proposed $C_n^{(3)}$ for a global testing problem in a model similar to MMSBM. Although their test statistic is also called the Signed Triangle, it is different from our SgnT statistic \eqref{DefineST}, for their tests are only applicable to models without degree heterogeneity. The analysis of the  Signed Polygon is also much more delicate than that of the Singed Cycle, as
the error $(\hat{\alpha}_n-\alpha_n)$ is much smaller than the errors in $(\hat{\eta} - \eta^*)$.

It remains to understand   
\begin{itemize}
\item how the Signed Polygon manages to reduce variance, 
\item  what are the analytical challenges. 
\end{itemize}  
Consider the first question. We illustrate it with the Ideal Signed Polygon \eqref{Ideal-SP} and the null case. In this case,  $\Omega = \theta \theta'$. It is seen 
$\eta^* = \theta$, $A_{ij} - \eta_{i}^* \eta_j^*  =  A_{ij} - \Omega_{ij} =  W_{ij}$, for $i \neq j$ (see (\ref{DCMM-matrixform}) for definition of $W$), and so 
\[
\widetilde{U}^{(m)}_n =  \sum_{i_1, i_2, \ldots, i_m (dist)} W_{i_1 i_2} W_{i_2 i_3} \ldots W_{i_m i_1}.   
\] 
In the sum, each term is an $m$-product of independent centered Bernoulli variables, and two terms $W_{i_1 i_2} W_{i_2 i_3} \ldots W_{i_m i_1}$ and  $W_{i_1' i_2'} W_{i_2' i_3'} \ldots W_{i_m' i_1'}$ are correlated only when $\{i_1, i_2, \ldots, i_m\}$ and $\{i_1', i_2', \ldots, i_m'\}$ are the vertices of the same polygon. Such a construction is known to be efficient in variance reduction (e.g., \cite{bubeck2016testing}).

In comparison, the main term of an order-$m$ GC test statistic   \cite{OGC} is  
\[
N_n^{(m)} =  \sum_{i_1, i_2, \ldots, i_m (dist)} A_{i_1 i_2} A_{i_2 i_3} \ldots A_{i_m i_1}. 
\] 
Since here the Bernoulli variables are not centered, we can split $N_n^{(m)}$ into two uncorrelated terms: $N_n^{(m)} = \widetilde{U}_n^{(m)}+ (N_n^{(m)} - \widetilde{U}_n^{(m)})$. Compared to the Signed Polygon, the additional variance comes from the second term, which is undesirably large in the less sparse case \cite{Ke2019}. 

{\bf Remark 4}. The above argument also explains  why the order-$2$ Signed Polygon does not work well, so the Signed Triangle statistic is the lowest order 
Signed Polygon statistic.  To see the point, note that when $m = 2$, 
$\widetilde{U}^{(m)}_n = \sum_{i_1 \neq i_2} W_{i_1 i_2}^2$ under the null,   
which has an unsatisfactory variance due to the square of the $W$-terms. 

Consider the second question.  We discuss with the SgnQ statistic. 
Recall that $\eta^*$ is a non-stochastic proxy of $\hat{\eta}$.  For any $1 \leq i, j \leq n$ and $i \neq j$, we decompose  
$\eta_i^* \eta_j^*  - \hat{\eta}_i \hat{\eta}_j  = \delta_{ij} + r_{ij}$,   
where $\delta_{ij}$ is the main term, which is a linear function of $\hat{\eta}_i$ and $\heta_j$, and $r_{ij}$ is the remainder term. Introduce 
\begin{equation} \label{DefinetildeOmega} 
\widetilde{\Omega} = \Omega - \eta^* (\eta^*)'. 
\end{equation} 
We have $A_{ij} - \hat{\eta}_i \heta_j = \widetilde{\Omega}_{ij} + W_{ij} + \delta_{ij} + r_{ij}$.   
After inserting this into $Q_n$, each $4$-product is now the product of $4$ bracketed terms, where each bracketed term is the sum of 
$4$ terms. Expanding the brackets and re-organizing,  $Q_n$ splits into $4 \times 4 \times 4 \times 4 = 256$ {\it post-expansion} sums, each of the form
\[ 
\sum_{i_1, i_2, i_3, i_4 (dist)}  a_{i_1 i_2} b_{i_2 i_3} c_{i_3 i_4} d_{i_4 i_1},  
\] 
where $a$ is a generic term which can equal to either of the four terms $\widetilde{\Omega}$, $W$, 
$\delta$, and $r$; same for $b, c$ and $d$.  While some of these terms 
may equal to each other,   the symmetry we can exploit is limited, due to  
  (a) degree heterogeneity,  (b) mixed-memberships, and (c) the underlying polygon structure.  As a result, we still have   more than $50$ post-expansion sums to analyze.  

The analysis of a post-expansion sum with the presence of one or more $r$-term is the most tedious of all, where we need to further decompose each $r$-term into three different terms. 
This requires analysis of more than $100$ additional post-expansion sums.

At first glance, we may think most of the post-expansion sums are easy to control via a crude bound (e.g., the Cauchy-Schwartz inequality).  Unfortunately, this is not the case, and many seemingly negligible terms turn out to be non-negligible. Here are some of the reasons. 
\begin{itemize} 
\item Due to the scarcity of data, the estimation error $(\hat{\eta}_i - \eta_i)$ is not sufficiently small.  
Also, severe degree heterogeneity dictates that a crude bound may be enough for some $\heta_i$  but not for other $\heta_i$. 
\item We aim to cover all interesting sparsity levels: a crude bound may be enough for a specific range of sparsity levels, but not for others. 
\item We desire to have a {\it single} test that works for all levels of sparsity. Alternatively, we can find one test that works well for the more sparse case and another test that works well for the less sparse case,  but this is less appealing from a practical  viewpoint. 
\end{itemize} 
As a result, we have to analyze a large number of post-expansion sums, where the analysis is subtle, extremely tedious, and error-prone, involving delicate combinatorics, due to the underlying polygon structure. See Section~\ref{sec:SgnedPolygon}.

\subsection{Organization of the paper} 
Section \ref{sec:main} 
focuses on the Region of Possibility and contains the upper bound argument. 
Section \ref{sec:LB} focuses on the Region of Impossibility and contains the 
lower bound argument. Section \ref{sec:SgnedPolygon} presents the key proof ideas, with the proof of secondary lemmas deferred to the  appendix. Section \ref{sec:Simul} presents the numerical study, and Section \ref{sec:Discu} discusses extensions and connections.

For any $q > 0$ and $\theta \in \mathbb{R}^n$,  $\|\theta\|_q$ denotes the $\ell^q$-norm of $\theta$ (when $q = 2$, we drop the subscript for simplicity). Also,  $\theta_{min}$ and $\theta_{max}$ denote $\min\{\theta_1,  \ldots, \theta_n\}$ and $\max\{\theta_1,  \ldots, \theta_n\}$, respectively.  For any $n > 1$, ${\bf 1}_n  \in \mathbb{R}^n$ denotes the vector of $1$'s. 
For two positive sequences $\{a_n\}_{n = 1}^{\infty}$ and $\{b_n\}_{n = 1}^{\infty}$, we write $a_n \sim b_n$ if $\lim_{n \goto \infty} a_n/b_n = 1$, 
and we write $a_n \asymp b_n$ if for sufficiently large $n$, there are two constants 
$c_2 > c_1 > 0$ such that $c_1 \leq a_n/b_n \leq c_2$. We use $\sum_{i_1,i_2,\ldots,i_m (dist)}$ to denote the sum over all $(i_1,\ldots,i_m)$ such that $1\leq i_k\leq n$ and $i_k\neq i_\ell$ for $1\leq k\neq \ell\leq m$ (so the number of summands is $n(n-1)\cdots(n-m+1)$). 
\section{The Signed Polygon test and the upper bound} \label{sec:main} 
For reasons aforementioned, we focus our discussion on the SgnT statistic $T_n$ and the SgnQ statistic $Q_n$, but the ideas are extendable to general Signed Polygon statistics. In Section~\ref{subsec:normality}, we establish the asymptotic normality of two statistics. In 
Section~\ref{subsec:twotests}, we use two statistics to construct two tests, the SgnT test and the SgnQ test. 
In Section~\ref{subsec:power}, we discuss the power of the two tests.

In a DCMM model, $\Omega=\Theta\Pi P\Pi'\Theta$, where  $\Theta = \diag(\theta_1, \ldots, \theta_n)$,   and $\Pi$ is the $n \times K$ membership matrix $[\pi_1, \pi_2, \ldots, \pi_n]'$.  
We assume as $n \goto \infty$,     
\begin{equation} \label{cond-theta}
\|\theta\| \goto \infty, \;\;\;  \theta_{max} \goto 0,  \;\;\;  \mbox{and} \;\;\; (\|\theta\|^2/\|\theta\|_1) \sqrt{\log(\|\theta\|_1)} \goto 0. 
\end{equation} 
The first condition is necessary. In fact, if $\|\theta\| \goto 0$, then the alternative is indistinguishable from the null, as suggested by lower bounds in Section \ref{sec:LB}.  The second one is mild as the eligible range for $\theta_{max}$ is roughly $(n^{-1/2}, 1)$.  The last one is weaker than that of $\theta_{max} \sqrt{\log(n)} \goto 0$, and is very mild.  

Moreover, introduce $G=\|\theta\|^{-2}\Pi'\Theta^2\Pi\in\mathbb{R}^{K\times K}$. 
This matrix is properly scaled and it can be shown that $\|G\| \leq 1$ (Appendix C,  appendix). 
When the null is true, $K = P = G = 1$, and we don't need any additional condition. When the alternative is true,  we assume 
\beq \label{cond-balance}
\frac{\max_{1\leq k\leq K}\{\sum_{i=1}^n \theta_i\pi_i(k)\}}{\min_{1\leq k\leq K}\{\sum_{i=1}^n \theta_i\pi_i(k)\}}\leq C, \qquad \|G^{-1}\|\leq C, \qquad \|P\| \leq C.   
\eeq 
The  conditions are mild. Take the first two for example.  When there is no mixed membership, they only require  
the $K$ classes to be relatively balanced.

\subsection{Asymptotic normality of the null}  \label{subsec:normality} 
The following two theorems are proved in Section~\ref{subsec:shortmainpf}. 
\begin{thm}[Limiting null of the SgnT statistic] \label{thm:null-SgnT}  
Consider the testing problem \eqref{Problem} under the DCMM model \eqref{model1a}-\eqref{condition1d}, where the condition (\ref{cond-theta}) is satisfied.  Suppose the null hypothesis is true. As $n \goto \infty$,  
\[
\mathbb{E}[T_n]  = o(\|\theta\|^3),  \qquad \mbox{and} \qquad \mathrm{Var}(T_n)  \sim  6 \|\theta\|^6.   
\]
and 
\[
\frac{T_n - \mathbb{E}[T_n]}{\sqrt{\mathrm{Var}(T_n)}}  \;\;   \longrightarrow \;\;  N(0,1), \qquad \mbox{in law}.   
\] 
\end{thm} 
\begin{thm}[Limiting null of the SgnQ statistic] \label{thm:null-SgnQ}   
Consider the testing problem \eqref{Problem} under the DCMM model \eqref{model1a}-\eqref{condition1d}, where the condition \eqref{cond-theta} is satisfied. Suppose the null hypothesis is true. As $n \goto \infty$, 
\[
\mathbb{E}[Q_n]  = (2 + o(1)) \|\theta\|^4,  \qquad \mbox{and} \qquad  \mathrm{Var}(Q_n) \sim  8 \|\theta\|^8, 
\] 
and 
\[
\frac{Q_n - \mathbb{E}[Q_n]}{\sqrt{\mathrm{Var}(Q_n)}}  \;\;   \longrightarrow \;\;  N(0,1), \qquad \mbox{in law}. 
\] 
\end{thm}

Note that under the null, the limiting distributions of 
$T_n/\sqrt{\mathrm{Var}(T_n)}$ and $Q_n/\sqrt{\mathrm{Var}(Q_n)}$ are $N(0,1)$ and $N(1/\sqrt{2}, 1)$, respectively. To appreciate the difference, recall that the Signed Polygon can be viewed as 
a plug-in statistic, where we replace $\eta^*$ in the Ideal Signed Polygon by $\hat{\eta}$. 
Under the null, the effect of the plug-in is negligible for SgnT but not for SgnQ, so the two limiting distributions are different. See Section \ref{sec:SgnedPolygon} for details.


\subsection{The level-$\alpha$ SgnT and SgnQ tests} \label{subsec:twotests} 
By Theorems~\ref{thm:null-SgnT} and \ref{thm:null-SgnQ}, the null variances of the two statistics depend on $\|\theta\|^2$.  To use the two statistics as  tests, we need to estimate $\|\theta\|^2$. 
For $\hat{\eta}$ and $\eta^*$ defined in (\ref{Definehateta}) and (\ref{Defineeta*}), respectively, we have 
$\hat{\eta} \approx \eta^*$ and $\eta^* = \theta$ under the null. A reasonable estimator for 
$\|\theta\|^2$  under the null is therefore $\|\hat{\eta}\|^2$. We propose to estimate 
$\|\theta\|^2$ with $(\|\hat{\eta}\|^2-1)$,  which corrects the bias and is slightly more accurate than 
$\|\hat{\eta}\|^2$. 
The following lemma is proved in the  appendix.  
\begin{lemma}[Estimation of $\|\theta\|^2$] \label{lemma:heta} 
Consider the testing problem \eqref{Problem} under the DCMM model \eqref{model1a}-\eqref{condition1d}, where the condition \eqref{cond-theta} holds when either hypothesis is true and condition (\ref{cond-balance}) holds when 
the alternative is true.  Then, under both hypotheses, as $n\to\infty$
\[
(\|\hat{\eta}\|^2 -1) / \|\eta^*\|^2  \goto 1, \qquad \mbox{in probability},  
\] 
where 
\[
\|\eta^*\|^2 = ({\bf 1}_n' \Omega^2 {\bf1}_n)/({\bf1}_n' \Omega {\bf1}_n) 
\begin{cases}
= \|\theta\|^2,               & \text{under }H_0^{(n)},  \\
\asymp \|\theta\|^2,  &\text{under }H_1^{(n)}.  
\end{cases}
\] 
\end{lemma} 

Combining Lemma~\ref{lemma:heta} with Theorem~\ref{thm:null-SgnT} gives 
\begin{equation} \label{DefineSgnTtestA}
\frac{T_n}{\sqrt{6(\|\hat{\eta}\|^2 - 1)^3}} \;\;   \longrightarrow \;\;  N(0,1), \qquad \mbox{in law}. 
\end{equation} 
Fix $\alpha \in (0,1)$. 
We propose the following SgnT test,  which is a two-sided test where we reject the null hypothesis if and only if 
\begin{equation} \label{SgnTtest} 
|T_n| \geq z_{\alpha/2} \sqrt{6} (\|\hat{\eta}\|^2-1)^{3/2}, \;\;\;  \mbox{$z_{\alpha/2}$: upper $\frac{\alpha}{2}$ quantile of $N(0,1)$}.  
\end{equation} 
Similarly, combining Theorem~\ref{thm:null-SgnQ} and Lemma~\ref{lemma:heta}, we have
\begin{equation} \label{DefineSgnQtestA} 
\frac{Q_n - 2(\|\hat{\eta}\|^2-1)^2}{\sqrt{8(\|\hat{\eta}\|^2 -1)^4}} \;\;   \longrightarrow \;\;  N(0,1), \qquad \mbox{in law}. 
\end{equation} 
With the same $\alpha$, we propose the following SgnQ test, which is a one-sided test where we reject the null hypothesis if and only if 
\begin{equation} \label{SgnQtest} 
Q_n \geq \bigl(2 + z_{\alpha} \sqrt{8}\bigr) (\|\hat{\eta}\|^2-1)^2,  \;\;\;  \mbox{$z_{\alpha}$: upper $\alpha$ quantile of $N(0,1)$}. 
\end{equation} 
As a result, for both tests we just defined, the levels satisfy 
\[
\mathbb{P}_{H_0^{(n)}}(\mbox{Reject the null}) \goto \alpha, \qquad \mbox{as $n \goto \infty$}.  
\]

\begin{figure}[tb] 
\centering
\includegraphics[height = 2 in, width = 2.45 in]{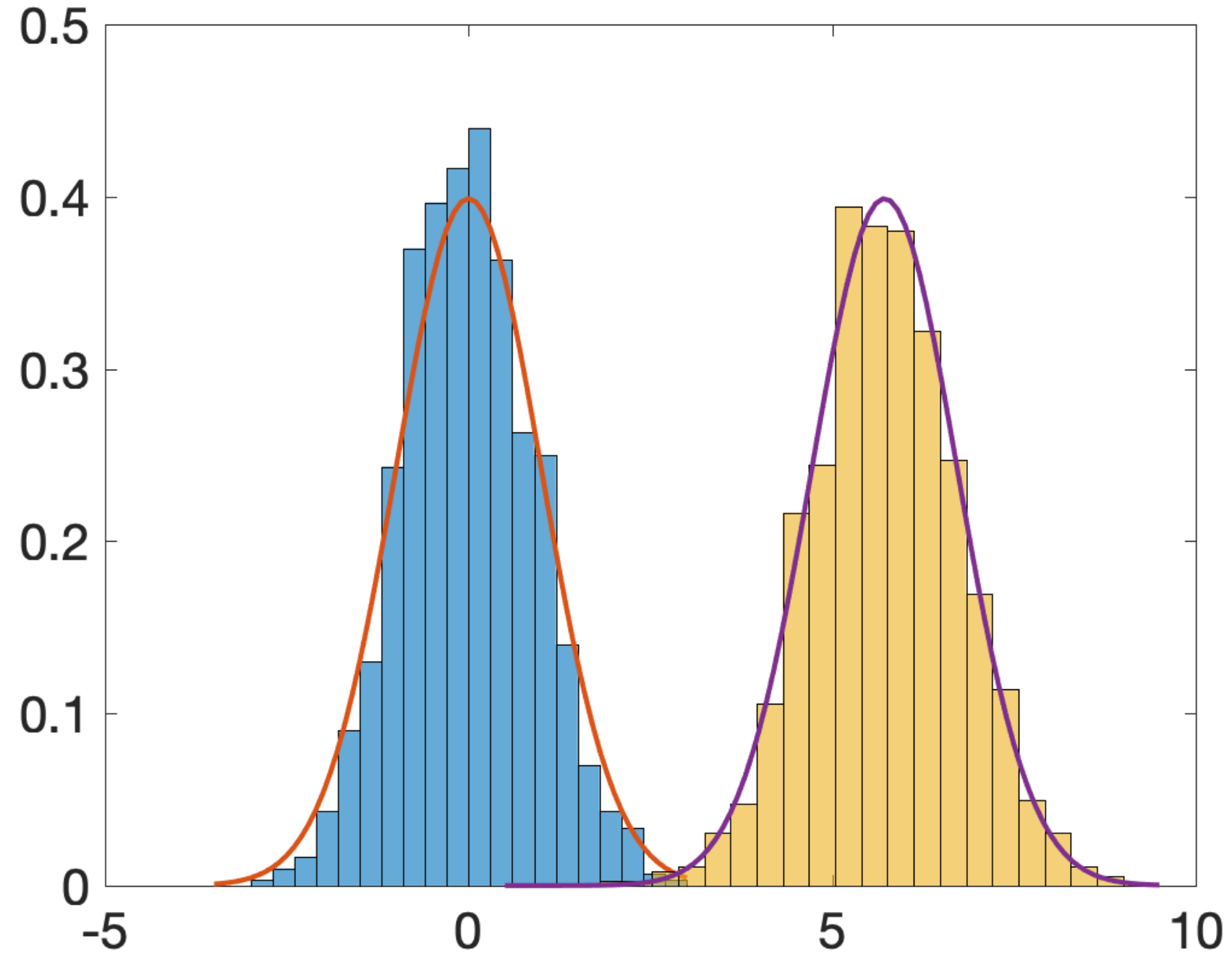} 
\includegraphics[height = 2 in, width = 2.45 in]{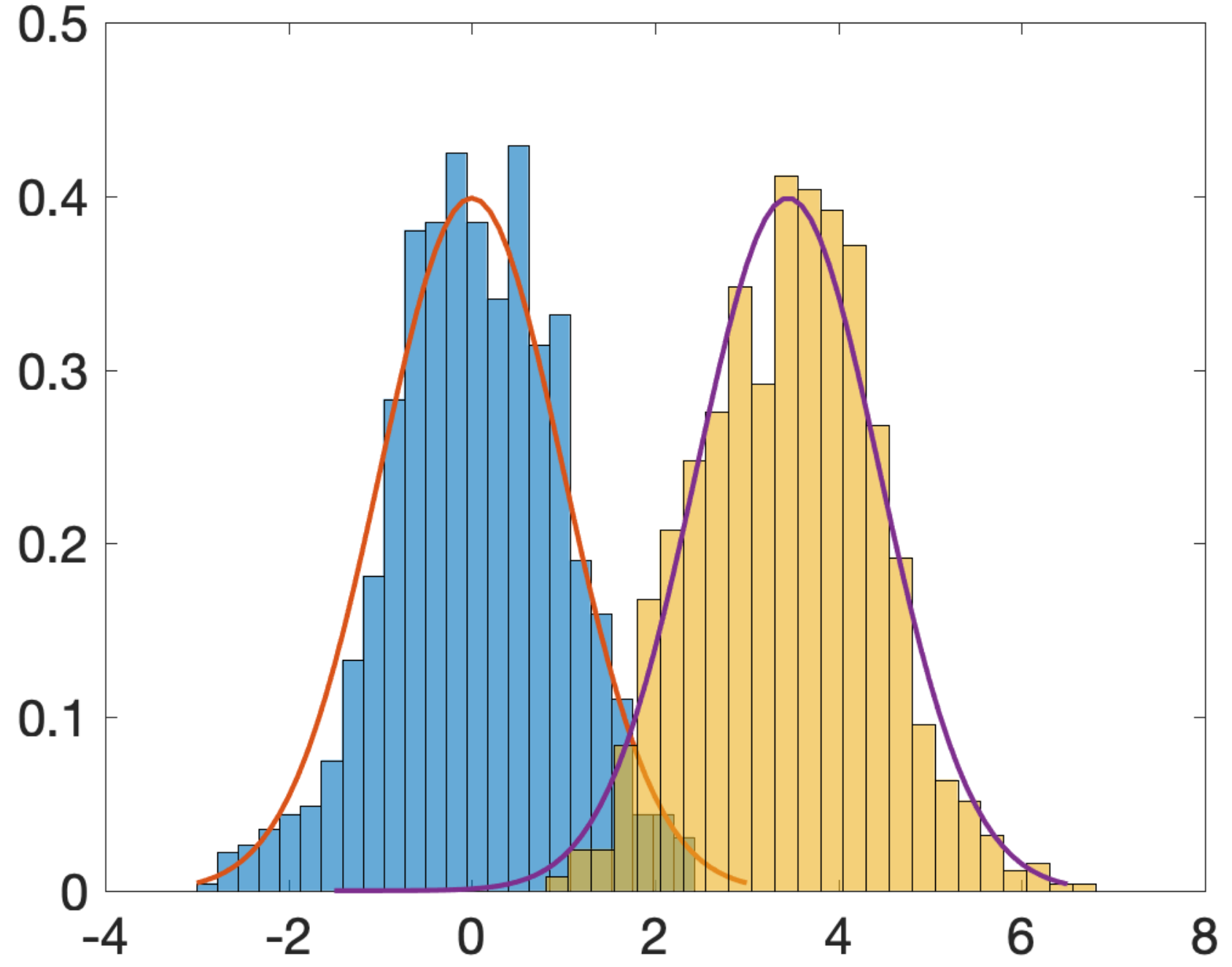} 
\caption{Left: histograms of the SgnT test statistics in (\ref{DefineSgnTtestA}) for the null (blue) and the alternative (yellow). Empirical mean and SD under the null: $0.04$ and $0.94$. Right:
same but for SgnQ test statistic in (\ref{DefineSgnQtestA}). Empirical mean and SD under the null: 
$-0.02$ and $0.92$. Repetition: $1000$ times.  See setting details in the main text.}  \label{fig:limitnull}
\end{figure}

Figure \ref{fig:limitnull} presents the histograms of $T_n/\sqrt{6 (\|\hat{\eta}\|^2 - 1)^3}$ (left) and 
$(Q_n - 2(\|\hat{\eta}\|^2-1)^2)/(\sqrt{8(\|\hat{\eta}\|^2 -1)^4})$ (right) 
under a null and an alternative setting simulated from DCMM. 
Recall that in DCMM, $\Omega = \theta \theta'$ under the null and $\Omega = \Theta \Pi P \Pi \Theta$, where $\Theta = \diag(\theta_1, \theta_2, \ldots, \theta_n)$. For the null,  
we take $n = 2000$ and draw $\theta_i$ from $\mathrm{Pareto}(12,3/8)$ and scale $\theta$ to have an $\ell^2$-norm of $8$. For the alternative, we let $(n, K) = (2000, 2)$,  $P$ be the matrix with $1$ on the diagonal and $0.6$ on the off-diagonal, rows of $\Pi$ equal to $\{1,0\}$ and $\{0, 1\}$ half by half, and with the same $\theta$ as in the null  but (to make it harder to separate from the null) rescaled to have an $\ell^2$-norm of $9$. The results confirm the limiting null of $N(0,1)$ for both tests.  
 
\subsection{Power analysis of the SgnT and SgnQ tests}  \label{subsec:power} 
The matrices $\Omega$ and $\widetilde{\Omega}$ play a key role in power analysis. 
Recall that $\Omega$ is defined in \eqref{model1c} where $\mathrm{rank}(\Omega) = K$, and  $\widetilde{\Omega} = \Omega - \eta^* (\eta^*)'$ is defined in (\ref{DefinetildeOmega}) with $\eta^* = \Omega {\bf1}_n / \sqrt{{\bf1}_n' \Omega {\bf1}_n}$  as in \eqref{Defineeta*}. Recall that  $\lambda_1, \lambda_2, \ldots, \lambda_K$ are the $K$ nonzero eigenvalues of $\Omega$. Let  $\xi_1, \xi_2, \ldots, \xi_K$ be the corresponding eigenvectors.  
The following theorems are proved in Section~\ref{subsec:shortmainpf}. 
\begin{thm}[Limiting behavior the SgnT statistic (altenative)] \label{thm:alt-SgnT}  
Consider the testing problem \eqref{Problem} under the DCMM model \eqref{model1a}-\eqref{condition1d}. 
Suppose the alternative hypothesis is true, and the conditions  \eqref{cond-theta}-\eqref{cond-balance} hold. As $n \goto \infty$,  
\[
\mathbb{E}[T_n]     = \tr(\widetilde{\Omega}^3) + o((|\lambda_2|/\lambda_1)^3 \|\theta\|^6)  + o(\|\theta\|^3),  
\] 
and 
\[ 
\mathrm{Var}(T_n)  
\leq C\bigl(\|\theta\|^6 + (\lambda_2/\lambda_1)^4 \|\theta\|^4 \|\theta\|_3^6\bigr).
\]  
\end{thm} 

\begin{thm}[Limiting behavior of the SgnQ statistic (alternative)] \label{thm:alt-SgnQ}   
Consider the testing problem \eqref{Problem} under the DCMM model \eqref{model1a}-\eqref{condition1d}. 
Suppose the alternative hypothesis is true and the conditions  \eqref{cond-theta}-\eqref{cond-balance} hold. As $n \goto \infty$, 
\[
\mathbb{E}[Q_n]  =  \tr(\widetilde{\Omega}^4) +    o((\lambda_2/\lambda_1)^4 \|\theta\|^8)  + o(\|\theta\|^4), 
\] 
and 
\[ 
\mathrm{Var}(Q_n) 
\leq C \bigl(\|\theta\|^8 + C (\lambda_2/\lambda_1)^6 \|\theta\|^8 \|\theta\|_3^6\bigr). 
\] 
\end{thm}  
We conjecture that both $T_n$ and $Q_n$ are asymptotically normal under the alternative. 
In fact, asymptotic normality is easy to establish for the Ideal SgnT and Ideal 
SgnQ. 
To establish results for the real SgnT and real SgnQ, we need very precise characterization of the 
plug-in effect. For reasons of space, we leave them to the future.

Consider the SgnT test \eqref{SgnTtest} first.  By Theorem~\ref{thm:alt-SgnT} and Lemma~\ref{lemma:heta},  under the alternative,  
\begin{equation} \label{SNRAnalysis1} 
\mbox{the mean and variance of $\frac{T_n}{\sqrt{6(\|\hat{\eta}\|^2-1)^3}}$ are $\frac{\tr(\widetilde{\Omega}^3)}{\sqrt{6\|\eta^*\|^6}}$ and $\sigma_n^2$}, 
\end{equation} 
respectively, 
where $\sigma_n^2$ denotes the asymptotic variance, which satisfies that
\begin{equation} \label{SNRAnalysis2} 
\sigma_n^2 \leq  \begin{cases}
C, &\mbox{if } |\lambda_2/\lambda_1| \ll \sqrt{\|\theta\|/\|\theta\|_3^3},\\
C (\lambda_2/\lambda_1)^4 \cdot (\|\theta\|_3^6/\|\theta\|^2), &\mbox{if } |\lambda_2/\lambda_1| \gg \sqrt{\|\theta\|/\|\theta\|_3^3}. 
\end{cases}
\end{equation} 
If we fix the degree heterogeneity vector $\theta$ and let 
$(\lambda_2/\lambda_1)$ range, there is a {\it phase change} in the variance. We shall call: 
\begin{itemize}
\item the case of $|\lambda_2/\lambda_1| \leq C \sqrt{\|\theta\|/\|\theta\|_3^3}$ as the {\it weak signal} case for SgnT.
\item the case of $|\lambda_2/\lambda_1| \gg \sqrt{\|\theta\|/\|\theta\|_3^3}$ as the {\it strong signal} case for SgnT.  
\end{itemize} 

It remains to derive a more explicit formula for $\tr(\widetilde{\Omega}^3)$. Recall that $\lambda_k$ and $\xi_k$ are the $k$-th eigenvalue and eigenvector of $\Omega$, $1\leq k\leq K$, respectively. 
When $\lambda_2, \ldots, \lambda_K$ have the same signs,  
$|\tr(\widetilde{\Omega}^3)| \geq C \sum_{k=2}^K|\lambda_k|^3$ (Lemma \ref{lem:tracetOmega1} below), so what we hope to show is that 
such an equality continues to hold in more general settings.  
Unfortunately, this is not true, and it is possible that the SgnT test suffers from the ``signal cancellation" in terms of 
\[
|\tr(\widetilde{\Omega}^3)| \ll \sum_{k = 2}^K |\lambda_k|^3. 
\] 
Normally, the ``signal cancellation" is found for odd-order moment-based statistics (e.g., $3rd$, $5th$, $\ldots$,  moment), but not for even-order moment methods (in fact, the SgnQ test won't experience such ``signal cancellation").

Fortunately,  ``signal cancellation" for the SgnT test  can be avoided if some conditions are imposed. 
Below, we list three conditions, and Lemma \ref{lem:tracetOmega1} shows that either of three conditions can ensure no  ``signal cancellation". 
Define $\Lambda\in\mathbb{R}^{(K-1)\times (K-1)}$ and $h \in\mathbb{R}^{K-1}$ by
\[
\Lambda = \diag(\lambda_2, \lambda_3, \ldots, \lambda_K), \qquad \mbox{and}
\qquad h_k = \frac{{\bf 1}_n' \xi_{k+1}}{{\bf 1}_n' \xi_1}, \;\; 1\leq k\leq K-1. 
\]
It can be shown that ${\bf 1}_n'\xi_1\neq 0$, so the vector $h$ is well-defined. Additionally, it can be shown that $\|h\|_\infty\leq C$.  See the proof of Lemma \ref{lem:tracetOmega1}. 
\begin{cond} \label{ExtraCondST} (a)  $\lambda_2, \lambda_3, \ldots, \lambda_K$ have the same signs, (b) $K = 2$, and (c) $|\lambda_2|/\lambda_1 \goto 0$, and $
|\tr(\Lambda^3) + 3 h' \Lambda^3 h + 3  (h' \Lambda h) (h' \Lambda^2 h) + (h' \Lambda h)^3| \geq C \sum_{k = 2}^K |\lambda_k|^3$. 
\end{cond}
The following is proved in Appendix C of the  appendix.  
\begin{lemma}[Analysis of $\tr(\widetilde{\Omega}^3)$] \label{lem:tracetOmega1}  
Suppose conditions of Theorem \ref{thm:alt-SgnT} hold. Under the alternative hypothesis,  
\begin{itemize}
\item If $|\lambda_2| /\lambda_1 \goto 0$,  then 
$\tr(\widetilde{\Omega}^3) = \tr(\Lambda^3) + 3 h' \Lambda^3 h + 3  (h' \Lambda h) (h' \Lambda^2 h) + (h' \Lambda h)^3 + o(|\lambda_2|^3)$.     
\item If $\lambda_2, \lambda_3, \ldots, \lambda_K$ have the same signs, then 
\[
|\tr(\widetilde{\Omega}^3)|  \geq   
\begin{cases} 
\sum_{k = 2}^K |\lambda_k|^3 + o(|\lambda_2|^3),  &\quad\mbox{if $|\lambda_2/\lambda_1| \goto 0$},  \\
C |\lambda_2|^3,  &\quad\mbox{if $|\lambda_2/\lambda_1| \geq C$}.  
\end{cases} 
\] 
\item In the special case where $K = 2$, the vector $h$ is a scalar, and  
\[
|\tr(\widetilde{\Omega}^3)| 
\begin{cases}
=  [(h^2+1)^3 + o(1)] \cdot  |\lambda_2|^3, &\quad\mbox{if $|\lambda_2|/\lambda_1 \goto 0$}, \\
\geq C |\lambda_2|^3,  &\quad \mbox{if $|\lambda_2/\lambda_1| \geq C$}.  
\end{cases} 
\] 
\end{itemize} 
As a result, when either one of  (a)-(c) holds, 
$|\tr(\widetilde{\Omega}^3)| \geq C \sum_{k = 2}^K |\lambda_k|^3$.  
\end{lemma}

It can be shown $\|\eta^*\| \asymp \sqrt{\lambda_1}  \asymp \|\theta\|$. Combining Lemma \ref{lem:tracetOmega1} with (\ref{SNRAnalysis1})-(\ref{SNRAnalysis2}), we have that  in the {\it weak signal} case,
\[
\frac{\mathbb{E}[T_n]}{\sqrt{\mathrm{Var}(T_n)}}\geq \frac{C\bigl(\sum_{k=2}^K|\lambda_k|^3\bigr)}{\|\theta\|^3}\geq C \Bigl(\lambda_1^{-\frac{3}{2}}\sum_{k=2}^K|\lambda_k|^3\Bigr),
\]
and that in the {\it strong signal} case, since $(\lambda_2/\lambda_1)^2\leq \lambda_1^{-2}(\sum_{k=2}^K|\lambda_k|^3)^{\frac{2}{3}}$, 
\[
\frac{\mathbb{E}[T_n]}{\sqrt{\mathrm{Var}(T_n)}}\geq \frac{C(\sum_{k=2}^K|\lambda_k|^3)}{\lambda_1^{-2}(\sum_{k=2}^K|\lambda_k|^3)^{\frac{2}{3}}\|\theta\|_3^3\|\theta\|^2}\geq \frac{C\|\theta\|^3}{\|\theta\|_3^3}\Bigl(\lambda_1^{-\frac{3}{2}}  \sum_{k=2}^K|\lambda_k|^3 \Bigr)^{\frac{1}{3}},   
\]
where it should be noted that in our setting, $\|\theta\|^3/\|\theta\|_3^3\to\infty$. As a result, in both cases, the power of the SgnT test $\to 1$ as long as $\lambda_1^{-3/2}\sum_{k=2}^K|\lambda_k|^3\to\infty$. This is validated in the following theorem, which is proved in Section \ref{subsec:shortmainpf}. 
\begin{thm}[Power of the SgnT test] \label{thm:SgnTtest} 
Under the conditions of Theorem~\ref{thm:alt-SgnT}, for any fixed $\alpha\in (0,1)$, consider the SgnT test in \eqref{SgnTtest}. Suppose one of the cases in Condition~\ref{ExtraCondST} holds. 
As $n \goto \infty$, if 
\[
\lambda_1^{-1/2} \cdot \Bigl(\sum_{k = 2}^K |\lambda_k|^3\Bigr)^{1/3}  \goto \infty, 
\] 
then the Type I error $\goto \alpha$, and the Type II error $\goto 0$. 
\end{thm}  

Next, we consider the SgnQ test \eqref{SgnQtest}. By Theorem~\ref{thm:alt-SgnQ} and Lemma~\ref{lemma:heta}, under the alternative hypothesis,
\[
\mbox{the mean and variance of $\frac{Q_n-2(\|\hat{\eta}\|^2-1)^2}{\sqrt{8(\|\hat{\eta}\|^2-1)^4}}$ are    $\frac{\tr(\widetilde{\Omega}^4)}{\sqrt{8\|\eta^*\|^8}}$ and $\sigma_n^2$}, 
\]
respectively, 
where $\sigma_n^2$ denotes the asymptotic variance, which satisfies that
\[
\sigma_n^2 \leq \begin{cases}
C, &\mbox{if } |\lambda_2/\lambda_1| \ll \|\theta\|^{-1}_3,\\
 C(\lambda_2/\lambda_1)^6 \cdot \|\theta\|_3^6, &\mbox{if } |\lambda_2/\lambda_1| \gg \|\theta\|^{-1}_3. 
\end{cases}
\]
Similar to the SgnT test, if we fix the degree heterogeneity vector $\theta$ and let 
$(\lambda_2/\lambda_1)$ range, there is a {\it phase change} in the variance. We shall call: 
\begin{itemize}
\item the case of $|\lambda_2/\lambda_1| \leq C\|\theta\|^{-1}_3$ as the {\it weak signal} case for SgnQ. 
\item the case of $|\lambda_2/\lambda_1| \gg \|\theta\|^{-1}_3$ as the {\it strong signal} case for SgnQ. 
\end{itemize}
We now analyze $\tr(\widetilde{\Omega}^4)$. The following lemma is proved in the  appendix. 
\begin{lemma}[Analysis of $\tr(\widetilde{\Omega}^4)$] \label{lem:tracetOmega2}  
Suppose the conditions of Theorem \ref{thm:alt-SgnQ} hold. Under the alternative hypothesis,   
\begin{itemize} 
\item If $|\lambda_2|/\lambda_1 \goto 0$, then $\tr(\widetilde{\Omega}^4) =  \tr(\Lambda^4) + (q' \Lambda q)^4 + 2 (h' \Lambda^2 h)^2 + 4 (h' \Lambda h)^2  (h' \Lambda^2 h) + 4 h' \Lambda^4 h 
  + 4 (h' \Lambda h) (h' \Lambda^3 h)   + o(\lambda_2^4) \gtrsim \sum_{k = 2}^4 \lambda_k^4$.  
\item If $|\lambda_2| / \lambda_1 \geq C$, then 
$\tr(\widetilde{\Omega}^4) \geq C  \sum_{k = 2}^K \lambda_k^4$.   
\item In the special case where $K = 2$, the vector $h$ is a scalar, and  
$\tr(\widetilde{\Omega}^4)  =  [(h^2+1)^4 + o(1)] \cdot \lambda_2^4$.   
\end{itemize} 
\end{lemma} 
As a result,  the SgnQ test has no issue of ``signal cancellation", and it always holds that  
$\tr(\widetilde{\Omega}^4) \geq C \sum_{k = 2}^K \lambda_k^4$. Then,  
in the {\it weak signal} case,
\[
\frac{\mathbb{E}[Q_n]}{\sqrt{\mathrm{Var}(Q_n)}}\geq \frac{C\bigl(\sum_{k=2}^K\lambda_k^4\bigr)}{\|\theta\|^4}\geq C \Bigl(\lambda_1^{-2}\sum_{k=2}^K \lambda_k^4\Bigr),
\]
In the {\it strong signal} case, since $(\lambda_2/\lambda_1)^3\leq \lambda_1^{-3}(\sum_{k=2}^K \lambda_k^4)^{\frac{3}{4}}$, 
\[
\frac{\mathbb{E}[Q_n]}{\sqrt{\mathrm{Var}(Q_n)}}\geq \frac{C\bigl(\sum_{k=2}^K\lambda_k^4\bigr)}{\lambda_1^{-3}(\sum_{k=2}^K\lambda_k^4)^{\frac{3}{4}} \|\theta\|_3^3\|\theta\|^4}\geq \frac{C\|\theta\|^3}{\|\theta\|_3^3}\Bigl(\lambda_1^{-2}  \sum_{k=2}^K \lambda_k^4 \Bigr)^{\frac{1}{4}},    
\]
where $\|\theta\|^3/\|\theta\|_3^3\to\infty$. So, in both cases, the power of the SgnQ test goes to $1$ if $\lambda_1^{-2}\sum_{k=2}^K\lambda_k^4\to\infty$. This is validated in Theorem \ref{thm:SgnQtest}, 
which is proved in Section \ref{subsec:shortmainpf}. 
\begin{thm}[Power of the SgnQ test] 
\label{thm:SgnQtest} 
Under the conditions of Theorem~\ref{thm:alt-SgnQ}, for any fixed $\alpha\in (0,1)$, consider the SgnQ test in \eqref{SgnQtest}. 
As $n \goto \infty$, if 
\[
\lambda_1^{-1/2} \Bigl(\sum_{k = 2}^K \lambda_k^4\Bigr)^{1/4} \goto \infty, 
\] 
then the Type I error $\goto \alpha$, and the Type II error $\goto 0$. 
\end{thm}  

In summary, Theorem~\ref{thm:SgnTtest} and Theorem~\ref{thm:SgnQtest} imply that as long as 
\begin{equation} \label{Maincondition} 
|\lambda_2| / \sqrt{\lambda_1} \goto \infty, 
\end{equation} 
the levels of SgnT and SgnQ tests tend to $\alpha$ as expected, and their powers tend to $1$. The SgnT test requires mild conditions to avoid ``signal cancellation", but the SgnQ test has no such issue. Our simulations further support that SgnQ may have better performance than SgnT.   See Section \ref{sec:Simul}.  

{\bf Remark 5}.  Practically, we prefer to fix $\alpha$, say, $\alpha = 5\%$.  If we allow the level $\alpha$ to change with $n$, then when (\ref{Maincondition}) holds, there is a sequence of $\alpha_n$ that tends to $0$ slowly enough such that $|\lambda_2| /(z_{\alpha_n/2} \cdot \sqrt{\lambda}_1)    \goto \infty$.  As a result, for either of the two tests, the Type I error $\goto 0$ and the power $\goto 1$, so the sum of Type I and Type II errors $\goto 0$.

{\bf Example 1 (contd)}. For this example,  $\lambda_1 \sim (1 + (K-1) b_n) \|\theta\|^2$, and $\lambda_k \sim   (1 - b_n)\|\theta\|^2, \;\; k = 2, 3, \ldots, K$. 
The condition \eqref{Maincondition} of $|\lambda_2| / \sqrt{\lambda_1} \goto \infty$ translates to  
$(1-b_n)\|\theta\|\to\infty$. 
See Section \ref{subsec:preview} and also Section \ref{sec:LB} for more discussion.

\section{Optimal adaptivity, lower bound,  and Region of Impossibility}    \label{sec:LB} 
We now focus on the Region of Impossibility, where  $|\lambda_2| / \sqrt{\lambda_1} \goto 0$.   
We first present a standard minimax lower bound, from which we can conclude that  
there is a sequence of hypothesis pairs (one alternative and one null) that are asymptotically indistinguishable. 
However, this does not answer the question whether {\it all alternatives} in the Region of Impossibility are indistinguishable from the null. To answer this question, we need much more sophisticated study; see Section \ref{subsec:RI}.

\subsection{Minimax lower bound}  \label{LB-minimax} 
Given an integer $K\geq 1$, a constant $c_0>0$, and two positive sequences $\{\alpha_n\}_{n = 1}^{\infty}$ and $\{\beta_n\}_{n = 1}^{\infty}$, we define a class of parameters for DCMM (recall that $\Omega = \Theta \Pi P \Pi' \Theta$, $G = \|\theta\|^{-2}\Pi'\Theta^2\Pi$ and is properly scaled, and $\lambda_k$ is the $k$-th largest eigenvalue of $\Omega$ in magnitude): 
\begin{align*}
& {\cal M}_n(K, c_0, \alpha_n, \beta_n)\cr
 = & \left\{
\begin{array}{rl}
(\theta, \Pi, P):& \theta_{\max}\leq \beta_n, \|\theta\|^{-1}\leq\beta_n, \|\theta\|^2\|\theta\|_1^{-1}\sqrt{\log(\|\theta\|_1)}\leq \beta_n,\\
& \frac{\max_k\{\sum_{i=1}^n \theta_i\pi_i(k)\}}{\min_k\{\sum_{i=1}^n \theta_i\pi_i(k)}\leq c_0, 
  \|G^{-1}\|\leq c_0,   |\lambda_2|/\sqrt{\lambda_1}\leq \alpha_n 
\end{array}
\right\}. 
\end{align*}
For the null case,  $K = P = \pi_i = 1$, and the above defines a class of $\theta$, which we write for short by 
\[
{\cal M}_n(1, c_0, \alpha_n, \beta_n)={\cal M}_n^*(\beta_n). 
\] 
The following theorem is proved in the  appendix: 
\begin{thm}[Minimax lower bound] \label{thm:minimax}
Fix $K\geq 2$, a constant $c_0>0$, and any sequences $\{\alpha_n\}_{n=1}^\infty$ and $\{\beta_n\}_{n=1}^\infty$ such that $\alpha_n\to 0$ and $\beta_n\to 0$ as $n\to\infty$. Then, as $n\to\infty$, 
\[
\inf_{\psi}\Bigl\{ \sup_{\theta\in {\cal M}^*_n(\beta_n)}\mathbb{P}(\psi=1) + \sup_{(\theta,\Pi,P)\in {\cal M}_n(K, c_0, \alpha_n, \beta_n)} \mathbb{P}(\psi=0)\Bigr\} \to 1,
\]
where the infimum is taken over all possible tests $\psi$. 
\end{thm}
The minimax theorem says that in the Region of Impossibility, {\it there exists a sequence of alternatives} that are inseparable from the null. This does not show what we desire, that is {\it any sequence in the Region of Impossibility} is inseparable from the null.  
This is discussed in the next section.  

\subsection{Region of Impossibility}  \label{subsec:RI} 
Recall that under DCMM, $\Omega = \Theta \Pi P \Pi' \Theta$ and $\Pi = [\pi_1, \pi_2, \ldots, \pi_n]'$.  Since our model is a mixed-membership latent variable model, in order to characterize the {\it least favorable configuration}, it is conventional to use a {\it random mixed-membership (RMM) model} for the matrix $\Pi$, while $(\Theta, P)$ are still non-stochastic. In detail, 
\begin{itemize} 
\item Let $V = \{x \in \mathbb{R}^K,  x_k \geq 0,  \sum_{k = 1}^K x_k = 1\}$.  
\item Let $V_0 = \{e_1, e_2, \ldots, e_K\}$, where $e_k$ is the $k$-th Euclidean basis vector. 
\end{itemize} 
In DCMM-RMM, we fix a distribution $F$ defined over $V$ and assume 
\[
\pi_i \stackrel{iid}{\sim}  F, \qquad \mbox{where $h \equiv \mathbb{E}[\pi_i]$}. 
\] 
If we further restrict that $F$ is defined over $V_0$, then 
the network has no mixed-membership, and DCMM-RMM reduces to DCBM-RMM. 

The desired result is to show that, for any given $P$ and $F$, there is a sequence of hypothesis pairs (a null and an alternative) 
\begin{equation} \label{LF-model} 
H_0^{(n)}: \quad \Omega = \theta \theta', \qquad\mbox{and}\qquad H_1^{(n)}: \quad \Omega = \widetilde{\Theta} \Pi P \Pi'  \widetilde{\Theta},   
\end{equation} 
where $\widetilde{\Theta} = \diag(\widetilde{\theta}_1, \widetilde{\theta}_1, \ldots, \widetilde{\theta}_n)$ and $\widetilde{\theta}_i$ can be different from $\theta_i$, 
such that the two hypotheses within each pair are asymptotically indistinguishable from each other, provided that under the alternative $|\lambda_2| / \sqrt{\lambda_1} \goto 0$.   

Here, since $\Omega$ depends on $\pi_i$, $\lambda_k$ is random, and it is more convenient to translate the condition of $|\lambda_2|/ \sqrt{\lambda_1} \goto 0$ to the condition of 
\begin{equation} \label{newcondition} 
\| \theta\| \cdot  |\mu_2(P)| \goto 0,  
\end{equation} 
where $\mu_k(P)$ is the $k$-th largest eigenvalue of $P$ in magnitude. The equivalence of 
two conditions are justified in Appendix D.1 of the  appendix. The regularity condition (\ref{cond-balance}) 
can also be ensured with high probability, by assuming that all entries of $\mathbb{E}[\pi_i]$  
are at the order of $O(1)$.

Under the DCBM, the desired result  can be proved satisfactorily. The key is the following lemma, which is in the line of Sinkhorn's beautiful work on scalable matrices \cite{DADSinkhorn}  (see also \cite{DADBrualdi, DADJohnson,  DADOlkin}) and is proved in the appendix. 
\begin{lemma} \label{lem:existence}  
Fix a matrix $A \in \mathbb{R}^{K,K}$ with strictly positive diagonal entries 
and non-negative off-diagonal entries, 
and a  strictly  positive vector $h \in \mathbb{R}^K$, 
there exists a diagonal matrix $D = \diag(d_1, d_2, \ldots, d_K)$ such that 
$D A D h    = 1_K$ and $d_k > 0$, $1 \leq k \leq K$. 
\end{lemma}

In detail, consider a DCBM-RMM setting where $\pi_i \stackrel{iid}{\sim} F$ and $F$ is supported over $V_0$  (with possibly unequal probabilities on the $K$ points).  Recall $h = \mathbb{E}[\pi_i]$. By Lemma \ref{lem:existence}, there is a unique diagonal matrix $D$ such that $D P D h  = 1_K$.   
Let 
\begin{equation} \label{config1} 
\widetilde{\theta}_i = d_k\cdot \theta_i, \qquad \mbox{if $\pi_i = e_k$},   \qquad 1 \leq i \leq n, \; 1 \leq k \leq K. 
\end{equation} 
The following theorem is proved in the  appendix. 

\begin{thm}[Region of Impossibility (DCBM)] \label{thm:RI-DCBM}    
Fix $K>1$ and a distribution $F$ defined over $V_0$. 
Consider a sequence of DCBM model pairs indexed by $n$: 
\[
H_0^{(n)}: \Omega = \theta \theta'  \quad\;\; \mbox{and}  
\quad\;\;
H_1^{(n)}: \Omega = \widetilde{\Theta} \Pi P \Pi' \widetilde{\Theta}, 
\]
where $\pi_i \stackrel{iid}{\sim} F$ and $\widetilde{\Theta} = \diag(\widetilde{\theta}_1, \widetilde{\theta}_2, \ldots, \widetilde{\theta}_n)$ with $\widetilde{\theta}_i$ defined as in (\ref{config1}).  If $\theta_{max} \leq c_0$ for a constant $c_0 < 1$, 
\[
\min_{1 \leq k \leq K} \{h_k\} \geq C, \qquad \mbox{and} \qquad \|\theta\| \cdot |\mu_2(P)| \goto 0,  
\]   
then for each pair of two hypotheses,  the $\chi^2$-distance between the two joint distributions tends to $0$, as $n \goto \infty$.   
\end{thm} 
 
We now generalize the result to DCMM. Fix a distribution $F$ defined over $V$. Given a set of $(\Theta, P, \Pi)$ with $\Theta = \diag(\theta_1, \theta_2, \ldots, \theta_n)$ and $\pi_i \stackrel{iid}{\sim} F$,  let $\widetilde{h}_D = \mathbb{E}[D^{-1} \pi_i / \| D^{-1} \pi_i\|_1]$ for any diagonal matrix $D\in\mathbb{R}^{K\times K}$ with positive diagonals. We assume that there exists  
$D$ such that 
\begin{equation} \label{config-condition} 
D P D \widetilde{h}_D = 1_K, \qquad \min_{1 \leq k \leq K} \{\widetilde{h}_{D,k}\} \geq C.  
\end{equation} 
When such a $D$ exists, we let 
\begin{equation} \label{config2} 
 \widetilde{\theta}_i = \theta_i/\|D^{-1} \pi_i\|_1, \qquad 1\leq i\leq n. 
\end{equation} 
When the support of $F$ is restricted to $V_0$, this reduces to the DCBM setting discussed above, in which (\ref{config-condition}) always holds, and $\widetilde{\theta}_i$ is the same as that in \eqref{config1}.  
When $K = 2$ (but the support of $F$ is not restricted to $V_0$), 
condition (\ref{config-condition}) also holds for all matrices $A$ in our setting. The proof is elementary so is omitted. 
 The following theorem is proved in the  appendix. 
\begin{thm}[Region of Impossibility (DCMM)] \label{thm:RI-DCMM}   
Fix $K>1$ and a distribution $F$ defined over $V$. 
Consider a sequence of DCMM model pairs indexed by $n$: 
\[
H_0^{(n)}: \Omega = \theta \theta'  \quad\;\; \mbox{and}  
\quad\;\;
H_1^{(n)}: \Omega = \widetilde{\Theta} \Pi P \Pi' \widetilde{\Theta}, 
\]
where $\pi_i \stackrel{iid}{\sim} F$ and $\widetilde{\Theta} = \diag(\widetilde{\theta}_1, \widetilde{\theta}_2, \ldots, \widetilde{\theta}_n)$ with $\widetilde{\theta}_i$ defined as in (\ref{config2}). If  (\ref{config-condition}) holds,  $\theta_{max} \leq c_0$ for a constant $c_0 < 1$,  and 
\[
\|\theta\| \cdot |\mu_2(P)| \goto 0,  
\]   
then for each pair of two hypotheses,  the $\chi^2$-distance between the two joint distributions tends to $0$, as $n \goto \infty$.   
\end{thm}

In Theorems \ref{thm:RI-DCBM} and \ref{thm:RI-DCMM}, we try to be as general as possible, where $F$ and $P$ are arbitrarily given,  and we seek for a $\Theta$-matrix in the alternative to make it most delicate to separate two hypotheses.   
We now consider a special case where $P$ is arbitrarily given, but $F$ is allowed to alter slightly. 
For any $P$ and $F$,  by Lemma \ref{lem:existence}, there is a unique positive diagonal matrix $D$ such that 
\begin{equation} \label{config3a} 
D P   D h = 1_K, \qquad \mbox{where} \qquad h = \mathbb{E}[\pi_i]. 
\end{equation} 
Let $\widetilde{\Pi} = [\widetilde{\pi}_1, \widetilde{\pi}_2, \ldots, \widetilde{\pi}_n]'$ and $\widetilde{\Theta} = \diag(\widetilde{\theta}_1, \widetilde{\theta}_2, \ldots, \widetilde{\theta}_n)$, with 
\begin{equation} \label{config3b} 
\widetilde{\pi}_i = D \pi_i/ \|D \pi_i\|_1,  \qquad \widetilde{\theta}_i = \|D \pi_i\|_1 \cdot \theta_i. 
\end{equation} 
The following theorem is proved in the  appendix. 
\begin{thm}[Region of Impossibility (DCMM with flexible $\Pi$)] \label{thm:RI-differentPi}    
Fix $K>1$ and a distribution $F$ defined over $V$. 
Consider a sequence of DCMM model pairs indexed by $n$: 
\[
H_0^{(n)}: \Omega = \theta \theta'  \quad\;\; \mbox{and}  
\quad\;\;
H_1^{(n)}: \Omega = \widetilde{\Theta} \widetilde{\Pi} P \widetilde{\Pi}' \widetilde{\Theta}, 
\]
where $\widetilde{\Pi}$ and $\widetilde{\Theta}$ are defined as in (\ref{config3a})-(\ref{config3b}). If   
$\theta_{max} \leq c_0$ for a constant $c_0 < 1$,   
\[
\min_{1 \leq k \leq K} \{h_k\} \geq C, \qquad \mbox{and} \qquad \|\theta\| \cdot |\mu_2(P)| \goto 0,  
\]   
then for each pair of two hypotheses,  the $\chi^2$-distance between the two joint distributions tends to $0$, as $n \goto \infty$.   
\end{thm} 

For completeness, one may wonder what happens if we require the null and the alternative have perfectly matching $\Theta$ matrix (up to an overall scaling). Such a scenario is natural when we focus on SBM or MMSBM,   
where degree heterogeneity is not allowed and so there is little freedom in choosing the $\Theta$ matrix. In this case, in order that the two hypotheses are indistinguishable, the expected node degrees under the alternative have to match those under the null. For each node $1 \leq i \leq n$, conditional on $\pi_i$ and neglecting the effect of no self edges, the expected degree equals to 
\[
\|\theta\|_1 \cdot \theta_i  \qquad \mbox{and} \qquad  \|\theta\|_1 \cdot (\pi_i'Ph) \cdot \theta_i,  
\] 
under the null and under the alternative, respectively, where $\{\pi_j\}_{j\neq i} \stackrel{iid}{\sim} F$ and $h = \mathbb{E}[\pi_j]$. 
For the expected degrees to match under any realized $\pi_i$, it is necessary that 
\begin{equation} \label{config4} 
P h = q_n 1_K, \qquad \mbox{for some scaling parameter $q_n  > 0$}. 
\end{equation} 
The following theorem is proved in the  appendix. 
\begin{thm}[Region of Impossibility (DCMM with matching $\Theta$)] \label{thm:RI-sameTheta}    
Fix $K>1$ and a distribution $F$ defined over $V$. 
Consider a sequence of DCMM model pairs indexed by $n$: 
\[
H_0^{(n)}: \Omega = q_n \cdot \theta \theta'  \quad\;\; \mbox{and}  
\quad\;\;
H_1^{(n)}: \Omega =  \Theta \Pi  P  \Pi' \Theta, 
\]
where $\Theta = \diag(\theta_1, \theta_2, \ldots, \theta_n)$, $\pi_i \stackrel{iid}{\sim}  F$, and $(P, h, q_n)$ satisfy  (\ref{config4}). If   $\theta_{max} \leq c_0$ for a constant $c_0 < 1$,    
\[
\min_{1 \leq k \leq K} \{h_k\} \geq C, \qquad \mbox{and} \qquad \|\theta\| \cdot |\mu_2(P)| \goto 0,  
\]   
then for each pair of two hypotheses,  the $\chi^2$-distance between the two joint distributions tends to $0$, as $n \goto \infty$.   
\end{thm} 

{\bf Example 1 (contd)}. In Example 1,  $\pi_i$ is drawn from $e_1, e_2, \ldots, e_K$ with equal probabilities, and 
$P=(1-b_n)I_K + b_n 1_K 1_K'$. 
Therefore, $h = \mathbb{E}[\pi_i] = (1/K) 1_K$.  In this case, all conditions of Theorem \ref{thm:RI-sameTheta}  hold, and especially,  $q_n = (1/K) + (K-1) b_n
 / K$, and $\mu_2(P) = (1 - b_n)$. 

{\bf Remark 6} (Least favorable configuration of LDA-DCMM).  The Dirichlet model is often used for mixed-memberships \cite{airoldi2009mixed}. 
Consider the model pairs 
\[
H_0^{(n)}: \Omega = q_n \theta\theta'  \qquad \mbox{and}  
\qquad 
H_1^{(n)}: \Omega = \Theta  \Pi P \Pi' \Theta, \;\;\; \pi_i\overset{iid}{\sim}\mathrm{Dir}(\alpha),
\]
where $\mathrm{Dir}(\alpha)$ is a Dirichlet distribution with parameters $\alpha=(\alpha_1,\ldots,\alpha_K)'$. 
By Theorem~\ref{thm:RI-sameTheta}, as long as $P\alpha\propto 1_K$, the null and alternative hypotheses are asymptotically indistinguishable if $(1-q_n)\|\theta\|\to 0$. One can easily construct $P$ such that $P\alpha\propto 1_K$. For example, $P=(1-q_n)MM'+q_n1_K1_K'$, where $M\in\mathbb{R}^{K\times(K-1)}$ is a matrix whose columns are from $Span^{\perp}(\alpha)$ and satisfy $\mathrm{diag}(MM')=I_K$.

 
\subsection{Optimal adaptivity}  \label{subsec:optimality} 
Recall that $\sqrt{\lambda_1}$, $|\lambda_2|/\lambda_1$, and $|\lambda_2|/\sqrt{\lambda_1}$ can be viewed as a measure for the sparsity, community similarity, and SNR, respectively.   
Combining Theorems~\ref{thm:null-SgnT}-\ref{thm:alt-SgnQ}, Theorems~\ref{thm:RI-DCBM}-\ref{thm:RI-sameTheta},  and Remark 5 in Section \ref{subsec:power},  
in the two-dimensional phase space where the $x$-axis is $\sqrt{\lambda_1}$ and the $y$-axis 
is the $|\lambda_2|/\lambda_1$, we have a partition to two regions, the Region of Possibility and the Region of Impossibility.  
\begin{itemize}
\item {\bf Region of Impossibility} ($1 \ll \sqrt{\lambda_1} \ll \sqrt{n}$, $|\lambda_2|/ \sqrt{\lambda_1}  = o(1)$).  
In this region, any DCBM alternative is asymptotically inseparable from the null, and up to a mild condition, any DCMM alternative is also asymptotically inseparable from the null. 
\item {\bf Region of Possibility} ($1 \ll \sqrt{\lambda_1} \ll \sqrt{n}$, $|\lambda_2|/ \sqrt{\lambda_1}  \goto \infty$). 
In this region, asymptotically, any alternative is completely separable form any null.   
\end{itemize} 
The SgnQ test is optimally adaptive:  for any alternative in the Region 
of Possibility, the test is able to separate it from the null with a sum of Type I and Type II errors tending to $0$.  
The SgnT test is also optimally adaptive, provided that some mild conditions hold to avoid signal cancellation.

To the best of our knowledge, the Signed Polygon is the only known test that 
is both applicable to general DCMM (where we allow severe degree heterogeneity and arbitrary mixed-memberships) and optimally adaptive. The EZ and GC tests are the only other tests we know 
that are applicable to general DCMM, but their variances are unsatisfactorily large for the less sparse case, so they are not optimally adaptive. See \cite{Ke2019} for details.

{\bf Remark 7}. Most lower bound results in the literature \cite{mossel2015reconstruction,banerjee2018contiguity,gao2017testing}
are in the standard minimax framework, where they focus on a particular sequence of alternative (e.g., the off-diagonals of $P$ are equal). In our case, the standard minimax theorem only implies  that in the Region of Impossibility, there is a sequence of alternative that are inseparable from the null. 
Our results (Theorems \ref{thm:RI-DCBM}-\ref{thm:RI-sameTheta}) are much stronger, implying that {\it any} alternative in the Region of Possibility is inseparable from the null.

{\bf Remark 8}.  Existing minimax lower bounds  \cite{mossel2015reconstruction,banks2016information,banerjee2018contiguity} 
have been largely focused on the SBM.  Though a least favorable scenario for SBM is also (one of the) least favorable scenario for DCMM, the former does not provide much insight on how the least favorable configurations and the separating boundary 
of the two regions (Possibility and Impossibility) depend on the degree heterogeneity and mixed-memberships. 
Moreover, our results suggest that $\|\theta\|$, not $\|\theta\|_1$, determines the separating boundary.  
In the SBM case,  $\theta_1 =  \ldots = \theta_n$ and $\|\theta\|_1 = \sqrt{n} \|\theta\|$, so it is hard to tell 
which of the two norms decides the boundary. In DCMM,  there is no simple 
relationship between $\|\theta\|_1$ and $\|\theta\|$, and we can tell this clearly.

\section{The behavior of the SgnT and SgnQ test statistics}  \label{sec:SgnedPolygon} 
Recall that the SgnT   statistic $T_n$ and the SgnQ statistic $Q_n$  are defined as 
\begin{align*} 
T_n & = \sum_{i_1, i_2, i_3 (dist)} (A_{i_1 i_2} - \hat{\eta}_{i_1} \hat{\eta}_{i_2}) (A_{i_2 i_3} - \eta_{i_2} \hat{\eta}_{i_3})  (A_{i_3 i_1} - \hat{\eta}_{i_3} \hat{\eta}_{i_1}),\cr
Q_n & = \sum_{i_1, i_2,i_3, i_4 (dist)} (A_{i_1 i_2} - \heta_{i_1} \heta_{i_2}) (A_{i_2 i_3} - \heta_{i_2} \heta_{i_3})  (A_{i_3 i_4} - \heta_{i_3} \heta_{i_4}) (A_{i_4 i_1} - \heta_{i_4} \heta_{i_1}),  
\end{align*}  
where 
\[
\hat{\eta} = A {\bf 1}_n /\sqrt{V}, \qquad \mbox{where} \;   V = {\bf 1}_n' A {\bf 1}_n. 
\]
In Section \ref{subsec:SP}, we have introduced the following non-stochastic proxy of $\heta$:  
\[
\eta^* = \Omega {\bf 1}_n / \sqrt{v_0}, \qquad \mbox{where} \; v_0 = {\bf 1}_n \Omega {\bf 1}_n.    
\]
We now introduce another non-stochastic proxy $\tilde{\eta}$ by  
\beq
\tilde{\eta} = A {\bf 1}_n / \sqrt{v}, \qquad \mbox{where $v = \mathbb{E}[{\bf 1}_n'  A {\bf 1}_n] = {\bf 1}_n (\Omega - \diag(\Omega)) {\bf 1}_n$}.   
\eeq 
Denoting the mean of $\tilde{\eta}$ by $\eta$, it is seen that 
\beq
\eta  = ([\Omega - \diag(\Omega)]  {\bf 1}_n) / \sqrt{{\bf 1}_n'  (\Omega - \diag(\Omega)) {\bf 1}_n}.  
\eeq
Here, $\eta$ and $\eta^*$ are close to each other but $\eta^*$ has a more explicit form. For example,  
under the null hypothesis, $\Omega = \theta \theta'$, and it is seen that $\eta^* = \theta$. Recall that 
\[
A = \Omega - \diag(\Omega) + W, \qquad \mbox{and} \qquad \widetilde{\Omega} = \Omega - \eta^* (\eta^*)'. 
\]
Fix $1 \leq i, j \leq n$ and $i \neq j$. First, we write 
\[
A_{ij} - \hat{\eta}_i \hat{\eta}_j = (A_{ij} - \eta^*_i \eta^*_j) + (\eta^*_i \eta_j^* - \hat{\eta}_i \hat{\eta}_j) = \widetilde{\Omega}_{ij} + W_{ij} + (\eta^*_i \eta_j^* - \hat{\eta}_i \hat{\eta}_j).  
\] 
Second, we write 
\[ 
\eta_i^*  \eta_j^* - \hat{\eta}_i \hat{\eta}_j   =    \delta_{ij} + r_{ij}, 
\] 
where  
\begin{equation} \label{Define-delta} 
\delta_{ij}  = \eta_i (\eta_j - \tilde{\eta}_j) + \eta_j (\eta_i - \tilde{\eta}_i) 
\end{equation} 
is  the linear approximation term of $(\eta_i^*  \eta_j^* - \hat{\eta}_i \hat{\eta}_j)$  and $
r_{ij} \equiv (\eta_i^* \eta_j^* - \hat{\eta}_i \hat{\eta}_j') - \delta_{ij}$  
is the remainder term. By definition and elementary algebra,  
\beq  \label{Define-r} 
r_{ij} = (\eta_i^* \eta_j^* - \eta_i \eta_j)  - (\eta_i - \tilde{\eta}_i) (\eta_j - \tilde{\eta}_j) + (1 - \frac{v}{V}) \tilde{\eta}_i \tilde{\eta}_j. 
\eeq 
It is seen that $r_{ij}$ is of a smaller order than that of $\delta_{ij}$.  The remainder term can be shown to  have a negligible effect over $T_n$ and $Q_n$, in terms of the  variances of $T_n$ and $Q_n$, respectively;   see Theorems \ref{thm:RealSgnT}-\ref{thm:RealSgnQ}.

Let $X$ be the symmetric matrix where all diagonal entries are $0$ and for $1 \leq i, j \leq n$ but $i \neq j$,  $X_{ij} = A_{ij} - \heta_i \heta_j$, or equivalently,  
\begin{equation} \label{DefineX}
X_{ij} = \widetilde{\Omega}_{ij} + W_{ij} + \delta_{ij}  + r_{ij}.  
\end{equation} 
If we omit the remainder term, then we have a proxy of $X$, denoted by $X^*$, where 
all diagonal entries of $X^*$ are $0$, and 
for $1 \leq i, j \leq n$ but $i \neq j$, 
\begin{equation} \label{DefineX*} 
X^*_{ij} = \widetilde{\Omega}_{ij} + W_{ij} + \delta_{ij}. 
\end{equation}  
If we further omit the $\delta$ term, then we have another proxy of $X$, denoted by $\widetilde{X}$, where all diagonal entries of $\widetilde{X}$ are $0$, and for $1 \leq i, j \leq n$ but $i \neq j$, 
\begin{equation} \label{DefinetX} 
\widetilde{X}_{ij} = \widetilde{\Omega}_{ij} + W_{ij}. 
\end{equation} 
With the above notations, we can rewrite $T_n$ and $Q_n$ as follows:
\[ 
T_n = \sum_{i_1, i_2, i_3 (dist)} X_{i_1i_2}X_{i_2i_3}X_{i_3i_1}, \qquad Q_n =  \sum_{i_1, i_2,i_3, i_4 (dist)} X_{i_1i_2}X_{i_2i_3}X_{i_3i_4}X_{i_4i_1}. 
\]
For the Ideal Signed Polygon in (\ref{Ideal-SP}), we have the 
{\it Ideal SgnT test statistic} 
\[ 
\widetilde{T}_n = \sum_{i_1, i_2, i_3 (dist)} (\widetilde{\Omega}_{i_1 i_2}  + W_{i_1 i_2}) 
(\widetilde{\Omega}_{i_2 i_3}  + W_{i_2 i_3}) (\widetilde{\Omega}_{i_3 i_1}  + W_{i_3 i_1}),  
\]
where with the notations above, 
\beq \label{DefineIdealSgnT} 
\widetilde{T}_n = \sum_{i_1, i_2, i_3 (dist)}  \widetilde{X}_{i_1 i_2}   \widetilde{X}_{i_2 i_3}   \widetilde{X}_{i_3 i_1}. 
\eeq
Similarly, we have the {\it Ideal SgnQ test statistic} 
\begin{equation} \label{DefineIdealSQ} 
\widetilde{Q}_n = \sum_{i_1, i_2, i_3, i_4 (dist)}  \widetilde{X}_{i_1 i_2} \widetilde{X}_{i_2 i_3} \widetilde{X}_{i_3 i_4}  \widetilde{X}_{i_4 i_1}.   
\end{equation}  
The Ideal SgnT and the Ideal SgnQ test statistics can be viewed as  proxies of the SgnT and the SgnQ test statistics, respectively, but such proxies are frequently not accurate enough. Therefore, 
we introduce another pair of proxies for SgnT and SgnQ, which we call the Proxy SgnT and Proxy SgnQ test statistics, respectively. Recall that $X_{ij}^* = \widetilde{\Omega}_{ij} + W_{ij} + \delta_{ij}$. 
\begin{definition}  \label{def:ST*SQ*}
The {\it Proxy SgnT test statistic} is
\[ 
T_n^*= \sum_{i_1, i_2, i_3 (dist)}  X_{i_1 i_2}^* X_{i_2 i_3}^* X_{i_3 i_1}^*. 
\]
The {\it Proxy SgnQ} test statistic is
\[
Q_n^*= \sum_{i_1, i_2, i_3, i_4 (dist)} X_{i_1 i_2}^* X_{i_2 i_3}^* X_{i_3 i_4}^* X_{i_4 i_1}^*.   
\]
\end{definition}
\noindent 
By these notations, we can partition SgnT as 
\[
T_n = \widetilde{T}_n + (T_n^* - \widetilde{T}_n) + (T_n - T_n^*), 
\] 
where the three terms on the right hand side are {\it the Ideal SgnT}, {\it the difference between the Ideal SgnT and the Proxy SgnT}, and {\it  the differnece between the Proxy SgnT and the real SgnT}, 
respectively.  Similarly, we can partition SgnQ by 
\[ 
\qquad Q_n =  \widetilde{Q}_n + (Q_n^* - \widetilde{Q}_n) + (Q_n - Q_n^*). 
\]
Below, first in  Section  \ref{subsec:IdealSgn}, we  analyze the Ideal SgnT and the Ideal SgnQ test statistics. Then in Section \ref{subsec:ProxySgn}, we analyze the difference between the Ideal SgnT and the Proxy SgnT, and that for the Ideal SgnQ and the Proxy SgnQ.   
Last, in Section \ref{subsec:RealSgn}, we analyze the difference between the Proxy SgnT and the real SgnT, and that for the Proxy SgnQ and the real SgnQ.

\subsection{The behavior of the Ideal SgnT and Ideal SgnQ test statistics}  \label{subsec:IdealSgn} 
The Ideal SgnT  test statistic is defined as  
\beq \label{IdealSgnT2}
\widetilde{T}_n = \sum_{i_1, i_2, i_3 (dist)} \widetilde{X}_{i_1 i_2}  \widetilde{X}_{i_2 i_3} \widetilde{X}_{i_3 i_1}, 
\eeq  
where $\widetilde{X}_{ij} = \widetilde{\Omega}_{ij} + W_{ij}$, $i \neq j$. 
Consider the null hypothesis first. When the null hypothesis holds, 
$\Omega = \theta \theta'$,    $\eta^* = \theta$, and 
\[
\widetilde{\Omega} = {\bf 0}_{n\times n}.  
\] 
In this case, the Ideal SgnT test statistic reduces to 
\[
\widetilde{T}_n = 
\sum_{i_1, i_2, i_3 (dist)}  W_{i_1i_2} W_{i_2 i_3} W_{i_3i_1},  
\]
which is the sum of a large number of uncorrelated terms, with each term being a $3$-product of independent centered-Bernoulli variables.  
It can be shown that 
\[
\mathbb{E}[\widetilde{T}_n] = 0, \qquad \mathrm{Var}(\widetilde{T}_n)  \sim 6 \|\theta\|^6,   
\] 
and 
\[
\widetilde{T}_n  / \sqrt{\mathrm{Var}(\widetilde{T}_n)} \; \longrightarrow \;  N(0,1), \qquad \mbox{in law}. 
\]

Consider the alternative hypothesis next.  In this case, the right hand side of (\ref{IdealSgnT2}) is the sum of many $3$-product, where each factor in the product is the sum of a $\widetilde{\Omega}$ term and a $W$ term. 
Expanding the bracket and re-arranging, we have 
\[
2 \times 2 \times 2  = 8 
\]  
post-expansion sums, 
each having the form of 
\begin{equation} \label{abc} 
\sum_{i_1, i_2, i_3  (dist)} a_{i_1 i_2} b_{i_2 i_3} c_{i_3 i_1},     
\end{equation} 
where $a$ is a generic notation which may either equal to $\widetilde{\Omega}$ or  $W$; same for $b$ and $c$.  For example, one of the $8$ post-expansion sums is 
\[ 
\sum_{i_1, i_2, i_3  (dist)} W_{i_1 i_2} \widetilde{\Omega}_{i_2 i_3} W_{i_3 i_1}, 
\]
which can be written in the form of (\ref{abc}) with  $b = \widetilde{\Omega}$, and $a = c  =W$. 
Note that each of  $8$ post-expansion sums is the sum of many $3$-product, 
where the number  of the $\widetilde{\Omega}$ factors in each product is the same; denote this number (which can be $0$, $1$, $2$, or $3$) by $N_{\widetilde{\Omega}}$. 
Similarly, the number  of the  $W$ factors in each product are also the same. Denote it by $N_W$, we have $N_{\widetilde{\Omega}}  + N_W = 3$. For the example above,  
$(N_{\widetilde{\Omega}}, N_W) = (1, 2)$. 

According to $(N_{\widetilde{\Omega}}, N_W)$, we can group the $8$ post-expansion sums into $4$ different types. Table \ref{tab:IdealSgnTsum} presents the mean and variance of each type 
(Recall that $\lambda_1, \lambda_2, \ldots, \lambda_K$ are the $K$ eigenvalues of $\Omega$, arranged in   descending order in magnitude. In Table~\ref{tab:IdealSgnTsum},  
$\alpha = |\lambda_2| / \lambda_1$.  When the alternative holds, we assume $|\lambda_2| /\sqrt{\lambda_1} \goto \infty$, and it translates to $\alpha \|\theta\| \goto \infty$ since $\sqrt{\lambda_1} \asymp \|\theta\|$).  
\begin{table}[htb!]
\centering
\caption{The $4$ types of the $8$ post-expansion sums for $\widetilde{T}_n$ ($\|\theta\|_q$ is the $\ell^q$-norm of $\theta$ (the subscript is dropped when $q = 2$). In our setting,  $\alpha \|\theta\| \goto \infty$, and $\|\theta\|_4^4   \ll     \|\theta\|_3^3   \ll     \|\theta\|^2   \ll   \|\theta\|_1$.}      \label{tab:IdealSgnTsum} 
\scalebox{0.97}{
\begin{tabular}{lcclcr}
Type  & $\#$ &  ($N_{\widetilde{\Omega}}, N_W)$ &     Examples  & Mean & Variance \\
\hline 
I   & $1$  &(0, 3)   & $\sum_{i, j, k (dist)} W_{ij} W_{jk} W_{ki} $ & 0 & $\asymp \|\theta\|^6$    \\  
II  & $3$  &(1, 2)   & $\sum_{i, j, k (dist)} \widetilde{\Omega}_{ij} W_{jk} W_{ki}$ & 0 & $\leq C\alpha^2 \|\theta\|^2 \|\theta\|_3^6 = o(\|\theta\|^6)$    \\ 
III & $3$ &(2, 1)   & $\sum_{i, j, k (dist)} \widetilde{\Omega}_{ij} \widetilde{\Omega}_{jk} W_{ki}$ & 0 & $\leq C\alpha^4\|\theta\|^4\|\theta\|_3^6$  \\   
IV & 1 & (3, 0)   & $\sum_{i, j, k (dist)} \widetilde{\Omega}_{ij} \widetilde{\Omega}_{jk} \widetilde{\Omega}_{ki}$    & $\sim \tr(\widetilde{\Omega}^3)$ &  0   \\  
\hline 
\end{tabular} 
}  
\end{table}

From the table, it can be concluded that among all $8$ post-expansion sums, 
the total mean is 
\[
\sim \tr(\widetilde{\Omega}^3), 
\] 
with Type IV sum being the only contributor, and 
the total variance is  
\[
\leq C \|\theta\|^6 + C(|\lambda_2|/\lambda_1)^4\|\theta\|^4\|\theta\|_3^6,  
\]  
with Type I sum and Type III sum being the major contributors. 
The following theorem is proved in the  appendix. 
\begin{thm}[Ideal SgnT test statistic] \label{thm:IdealSgnT} 
Consider the testing problem \eqref{Problem} under the DCMM model \eqref{model1a}-\eqref{condition1d}, where the condition \eqref{cond-balance} is satisfied under the alternative hypothesis. Suppose $\theta_{\max}\to 0$ and $\|\theta\| \goto \infty$ as $n \goto \infty$, and suppose $|\lambda_2|/\sqrt{\lambda_1}\to\infty$ under the alternative hypothesis. Then, under the null hypothesis, as $n \goto \infty$, 
\[
\mathbb{E}[\widetilde{T}_n]=0, \qquad \mathrm{Var}(\widetilde{T}_n)=6\|\theta\|^6\cdot [1+o(1)], 
\]
and
\[
\frac{\widetilde{T}_n - \mathbb{E}[\widetilde{T}_n]}{\sqrt{\mathrm{Var}(\widetilde{T}_n)}}  \;\;   \longrightarrow \;\;  N(0,1), \qquad \mbox{in law}.
\] 
Furthermore, under the alternative hypothesis, as $n\to\infty$, 
\[
\mathbb{E}[\widetilde{T}_n] = \tr(\widetilde{\Omega}^3) + o(\|\theta\|^3),  \qquad  
\mathrm{Var}(\widetilde{T}_n)  \leq C \|\theta\|^6 + C(|\lambda_2|/\lambda_1)^4\|\theta\|^4\|\theta\|_3^6. 
\] 
\end{thm}

Next, we discuss the Ideal SgnQ test statistic, which is defined as 
\begin{equation} \label{IdealSgnQ2} 
\widetilde{Q}_n = \sum_{i_1, i_2, i_3, i_4 (dist)}  \widetilde{X}_{i_1 i_2}  
\widetilde{X}_{i_2 i_3}  \widetilde{X}_{i_3 i_4}   \widetilde{X}_{i_4i_1},
\end{equation}  
where for any $i \neq j$,  $\widetilde{X}_{ij} = \widetilde{\Omega}_{ij}  + W_{ij}$.
Under the null, since $\widetilde{\Omega}$ is a zero matrix, the statistic reduces to 
\[
\widetilde{Q}_n = \sum_{i_1, i_2, i_3, i_4 (dist)}  W_{i_1 i_2} W_{i_2 i_3}   W_{i_3i_4} W_{i_4i_1}.   
\] 
Similarly, it can be shown that the statistic is asymptotically normal, with  
\[
\mathbb{E}[\widetilde{Q}_n] = 0, \qquad \mbox{and} \qquad \mathrm{Var}(\widetilde{Q}_n) \sim 8 \|\theta\|^8.   
\] 
Under the alternative, similarly, we obtain 
\[
2 \times 2 \times 2 \times 2 = 16 
\] 
post-expansion sums, and divide them into $6$ different types, according to $(N_{\widetilde{\Omega}}, N_W$). See Table \ref{tab:IdealSgnQsums}, where we recall  
$\alpha = |\lambda_2|/\lambda_1$.   
\begin{table}[htb!]   
\caption{The $6$ different types of the $16$ post-expansion sums of $\widetilde{Q}_n$. Notations:  same as Table \ref{tab:IdealSgnTsum}.}  
\centering
\scalebox{0.85}{
\begin{tabular}{lcclcr}
Type  & $\#$ &  ($N_{\widetilde{\Omega}}, N_W)$ &     Examples  &  Mean & Variance \\
\hline 
I   & $1$  &(0, 4)   & $\sum_{i, j, k, \ell (dist)} W_{ij} W_{jk} W_{k \ell} W_{\ell i}$ & 0 & $\asymp  \|\theta\|^8$    \\  
II  & $4$  &(1, 3)   & $\sum_{i, j, k, \ell (dist)} \widetilde{\Omega}_{ij} W_{jk} W_{k \ell} W_{\ell i}$ & 0 & $\leq C\alpha^2   \|\theta\|^4 \|\theta\|_3^6 = o(\|\theta\|^8)$    \\ 
IIIa & $4$ &(2, 2)   & $\sum_{i, j, k, \ell (dist)} \widetilde{\Omega}_{ij}  \widetilde{\Omega}_{jk} W_{k \ell} W_{\ell i}$ & 0 & $\leq C\alpha^4  \|\theta\|^6 \|\theta\|_3^6=o( \alpha^6\|\theta\|^8  \|\theta\|_3^6 )$  \\   
IIIb             & 2  & (2, 2)  &  $\sum_{i, j, k, \ell (dist)} \widetilde{\Omega}_{ij} W_{jk} \widetilde{\Omega}_{k \ell} W_{\ell i}$ & 0 &  
$\leq C  \alpha^4 \|\theta\|_3^{12} = o(\|\theta\|^8)$ \\  
IV & 4 & (3, 1)   & $\sum_{i, j, k, \ell (dist)} \widetilde{\Omega}_{ij} \widetilde{\Omega}_{jk} \widetilde{\Omega}_{k \ell} W_{\ell i}$    & 0 & $\leq \alpha^6\|\theta\|^8  \|\theta\|_3^6$   \\  
V & 1 & (4, 0)  & $\sum_{i, j, k, \ell (dist)} \widetilde{\Omega}_{ij} \widetilde{\Omega}_{jk} \widetilde{\Omega}_{k \ell} \widetilde{\Omega}_{\ell i}$ & $\sim \tr(\widetilde{\Omega}^4)$  & 0 \\ 
\hline 
\end{tabular} 
} 
 \label{tab:IdealSgnQsums}
\end{table} 

From the table, among all $16$ post-expansion sums, 
the total mean is 
\[
\sim \tr(\widetilde{\Omega}^4), 
\] 
with Type V sum being the only contributor, and 
the total variance   
\[
\leq C \|\theta\|^8 + C(|\lambda_2|/\lambda_1)^6\|\theta\|^8\|\theta\|_3^6,  
\]   
with Type I sum and Type IV sum being the major contributors.   
The following theorem is proved in the  appendix. 
\begin{thm}[Ideal SgnQ test statistic] \label{thm:IdealSgnQ} 
Consider the testing problem \eqref{Problem} under the DCMM model \eqref{model1a}-\eqref{condition1d}, where the condition \eqref{cond-balance} is satisfied under the alternative hypothesis. Suppose $\theta_{\max}\to 0$ and $\|\theta\| \goto \infty$ as $n \goto \infty$, and suppose $|\lambda_2|/\sqrt{\lambda_1}\to\infty$ under the alternative hypothesis. Then, under the null hypothesis, as $n \goto \infty$, 
\[
\mathbb{E}[\widetilde{Q}_n]=0, \qquad \mathrm{Var}(\widetilde{Q}_n)=8\|\theta\|^8\cdot [1+o(1)], 
\]
and
\[
\frac{\widetilde{Q}_n - \mathbb{E}[\widetilde{Q}_n]}{\sqrt{\mathrm{Var}(\widetilde{Q}_n)}}  \;\;   \longrightarrow \;\;  N(0,1), \qquad \mbox{in law}.
\] 
Furthermore, under the alternative hypothesis, as $n\to\infty$, 
\[
\mathbb{E}[\widetilde{Q}_n]  = \tr(\widetilde{\Omega}^4)  + o(\|\theta\|^4),  \qquad  
\mathrm{Var}(\widetilde{T}_n)  \leq C[\|\theta\|^8 + (|\lambda_2|/\lambda_1)^6\|\theta\|^8\|\theta\|_3^6]. 
\] 
\end{thm}

\subsection{The behavior of $(T_n^* - \widetilde{T}_n)$ and $(Q_n^* - \widetilde{Q}_n)$}  
\label{subsec:ProxySgn} 
Consider $(T_n^* - \widetilde{T}_n)$ first. The Proxy SgnT test statistic $T_n^*$ is defined as 
\begin{equation} \label{ProxySgnT2} 
T^*_n = \sum_{i_1, i_2, i_3 (dist)} X_{i_1 i_2}^* X_{i_2 i_3}^* X_{i_3 i_1}^*,  
\end{equation} 
where for all $i \neq j$, $X_{ij}^* = \widetilde{\Omega}_{ij} + W_{ij} + \delta_{ij}$, and 
$\delta_{ij} = \eta_i (\eta_j - \tilde{\eta}_j) + \eta_j (\eta_i - \tilde{\eta}_i)$.  
Similar to Section \ref{subsec:IdealSgn},  we first expand every bracket in (\ref{ProxySgnT2}) and obtain 
\[
3 \times 3 \times 3 = 27 
\]  
different post-expansion sums corresponding to $T^*_n$.  Out of these $27$ post-expansion sums, $2 \times 2 \times 2 = 8$ of them do not involve any $\delta$ term. The sum of these $8$ 
post-expansion sums equals to 
$\widetilde{T}_n$,  and the sum the remaining $19$ post-expansion sums equals to 
$(T_n^* - \widetilde{T}_n)$. 

For each of these $19$ post-expansion sums, we are summing over many $3$-products, 
where each of them has the same number  of $\widetilde{\Omega}$ factors,  $W$ factors, and $\delta$ factors,  which we denote  by 
$N_{\widetilde{\Omega}}, N_W$, and $N_{\delta}$, respectively.  According to $(N_{\widetilde{\Omega}}, N_W,  N_{\delta})$,  
 we divide the $19$ post-expansion sums into $6$ different types; see Table 
\ref{tab:ProxySgnTsum}.  

\begin{table}[htb!]  
\caption{The $6$ types of the $19$ post-expansion sums for $(T_n^* - \widetilde{T}_n)$. 
Notations: same as Table \ref{tab:IdealSgnTsum}.}
\centering
\scalebox{.8}{
\begin{tabular}{lcclcr}
Type  & $\#$ &  ($N_{\delta}, N_{\widetilde{\Omega}}, N_W)$ &     Examples  & Abs. Mean &  Variance \\
\hline 
Ia   & $3$  &(1, 0, 2)   & $\sum_{\substack{i, j, k\\ (dist)}} \delta_{ij} W_{jk} W_{k i}$ & 0 & $\leq C \frac{\|\theta\|^4\|\theta\|_3^3}{\|\theta\|_1}=o(\|\theta\|^6)$   \\  
Ib  & $6$  &(1, 1, 1)   & $\sum_{\substack{i, j, k\\ (dist)}} \delta_{ij}\widetilde{\Omega}_{jk}  W_{k i}$ &   
$\leq C\alpha \|\theta\|^4$=$o(\alpha^3\|\theta\|^6)$ & $\leq \frac{C\alpha^2\|\theta\|^6\|\theta\|_3^3}{\|\theta\|_1}=o(\|\theta\|^6)$  \\ 
Ic & $3$ &(1, 2, 0)   & $\sum_{\substack{i, j, k\\ (dist)}} \delta_{ij}  \widetilde{\Omega}_{jk} \widetilde{\Omega}_{ki}$ & 0 & $\leq \frac{C\alpha^4\|\theta\|^8\|\theta\|_3^3}{\|\theta\|_1} = O(\alpha^4\|\theta\|^4\|\theta\|_3^6)$\\   
IIa & 3 & (2, 0, 1)   & $\sum_{\substack{i, j, k\\ (dist)}} \delta_{ij} \delta_{jk} W_{k i}$    & $\leq C\|\theta\|^2$=$o(\|\theta\|^3)$ & $\leq C\|\theta\|_3^6=o(\|\theta\|^6)$ \\  
IIb & 3 & (2, 1, 0) & $\sum_{\substack{i, j, k\\ (dist)}} \delta_{ij} \delta_{jk} \widetilde{\Omega}_{k i}$    &  $\leq \frac{C\alpha\|\theta\|^6}{\|\theta\|_1^2}$=$o(\|\theta\|^3)$ & $\leq \frac{C\alpha^2\|\theta\|^{10}}{\|\theta\|_1^2}=o(\|\theta\|^6)$\\  
III & 1 &  (3,  0,  0)  & $\sum_{\substack{i, j, k\\ (dist)}} \delta_{ij} \delta_{jk} \delta_{ki}$ & $\leq \frac{C\|\theta\|^4}{\|\theta\|_1^2}$=$o(\|\theta\|^3)$ & $\leq \frac{C\|\theta\|^4\|\theta\|_3^3}{\|\theta\|_1}=o(\|\theta\|^6)$ \\ 
\hline 
\end{tabular} 
} 
  \label{tab:ProxySgnTsum} 
\end{table}

Consider the null hypothesis first. When the null hypothesis holds, $\widetilde{\Omega}$ is a zero matrix, so only Type Ia, Type IIa and Type III post-expansion sums are nonzero. From the table, it is seen that  
\[
| \mathbb{E}[T_n^*-\widetilde{T}_n]| =  o(\|\theta\|^3), 
\] 
and that 
\[
\mathrm{Var}(T_n^* - \widetilde{T}_n)  = o(\|\theta\|^6). 
\] 
Note that $\|\theta\|^6$ is the order of $\mathrm{Var}(\widetilde{T}_n)$ under the null. These show that under the null,  the difference between $T_n^*$ and $\widetilde{T}_n$ is negligible.

Consider the alternative hypothesis next. First, from the table,  
\[
|\mathbb{E}[T_n^*-\widetilde{T}_n]| \leq  C(|\lambda_2|/\lambda_1)\|\theta\|^4 +  o(\|\theta\|^3),    
\] 
where the main contribution comes from Type Ib post-expansion sum. 
On the right hand side, the second part is smaller than $\|\theta\|^3$ (the standard deviation of $T_n^*$ under the null) and is thus negligible. The first part can be much larger than $\|\theta\|^3$, 
but under our range of interest, $|\lambda_2| / \sqrt{\lambda_1} \goto \infty$ and $\lambda_1 \asymp 
\|\theta\|^2$, it is seen that this part 
is $o(|\lambda_2|^3)$,  where $|\lambda_2|^3$ is the order of $\mathbb{E}[\widetilde{T}_n]$ 
and $\tr(\widetilde{\Omega}^3)$; see Lemma~\ref{lem:tracetOmega1} and Theorem~\ref{thm:IdealSgnT}.  
Therefore, $\mathbb{E}[T_n^*]$ only has a negligible difference when compared to 
$\mathbb{E}[\widetilde{T}_n]$. 
Second, from the table, 
\[
\mathrm{Var}(T_n^* - \widetilde{T}_n)  \leq  \frac{C(|\lambda_2|/\lambda_1)^4\|\theta\|^8\|\theta\|_3^3}{\|\theta\|_1} + o(\|\theta\|^6), 
\] 
where the second term is smaller than $\mathrm{Var}(\widetilde{T}_n)$, and the first term is no greater than 
\[ 
C(|\lambda_2|/\lambda_1)^4\|\theta\|^4\|\theta\|_3^6,   
\] 
where we have used $\|\theta\|^4 \leq \|\theta\|_1 \|\theta\|_3^3$ (from the Cauchy-Schwartz inequality).  
Note that $\|\theta\|^6+(|\lambda_2|/\lambda_1)^4\|\theta\|^4\|\theta\|_3^6$ is the order of $\mathrm{Var}(\widetilde{T}_n)$; see Theorem \ref{thm:IdealSgnT}. 
Therefore, $\mathrm{Var}(T_n^*) \leq C \mathrm{Var}(\widetilde{T}_n)$. 
Combining these, the SNR of $\widetilde{T}_n$ and 
$T_n^*$ are at the same order. 

The above results are summarized in the following theorem. It is proved in the  
 appendix. 
\begin{thm}[Proxy SgnT test statistic] \label{thm:ProxySgnT} 
Consider the testing problem \eqref{Problem} under the DCMM model \eqref{model1a}-\eqref{condition1d}, where the condition \eqref{cond-balance} is satisfied under the alternative hypothesis. Suppose $\theta_{\max}\to 0$ and $\|\theta\| \goto \infty$ as $n \goto \infty$, and suppose $|\lambda_2|/\sqrt{\lambda_1}\to\infty$ under the alternative hypothesis. Then, under the null hypothesis, as $n \goto \infty$, 
\[
\mathbb{E}[T_n^*-\widetilde{T}_n] = o(\|\theta\|^3), \qquad \mathrm{Var}(T_n^* - \widetilde{T}_n)= o(\|\theta\|^6).  
\]
Furthermore, under the alternative hypothesis, 
\begin{align*}
& \mathbb{E}[T_n^*-\widetilde{T}_n] = o((|\lambda_2|/\lambda_1)^3\|\theta\|^6), \\
& \mathrm{Var}(T_n^* - \widetilde{T}_n) \leq C(|\lambda_2|/\lambda_1)^4\|\theta\|^4\|\theta\|_3^6 + o(\|\theta\|^6). 
\end{align*}
\end{thm}

Next, we consider $(Q_n^* - \widetilde{Q}_n)$, which is defined as
\begin{equation} \label{ProxySgnQ2} 
\widetilde{Q}_n = \sum_{i_1, i_2, i_3, i_4 (dist)}  X_{i_1 i_2}^* X_{i_2 i_3}^* X_{i_3 i_4}^* X_{i_4 i_1}^*.  
\end{equation} 
Similarly, if we expand the bracket of all individual terms and re-organize, we have 
\[
3 \times 3 \times  3 \times 3 = 81 
\] 
post-expansion sums. Out of the $81$ post-expansion sums, $2\times 2\times 2\times 2 = 16$ of them do not depend on $\delta$, the  
sum of which equals to $\widetilde{Q}_n$. These leave us with $65$ post-expansion sums, the total sum of which is $(Q_n^* - \widetilde{Q}_n)$.  Similarly, according to 
$(N_{\widetilde{\Omega}}, N_W, N_{\delta})$, we divide these $65$ sums into $10$ types. 
See Table \ref{tab:ProxySgnQsum}, where we recall that $\alpha=|\lambda_2|/\lambda_1$. 

\begin{table}[htb!]  
\caption{The $10$ types of the post-expansion sums for $(Q_n^* - \widetilde{Q}_n)$. Notations: same as in Table \ref{tab:IdealSgnTsum}.}     \label{tab:ProxySgnQsum} 
\centering
\scalebox{.83}{
\begin{tabular}{lcclcr}
Type  & $\#$ &  ($N_{\delta}$,$N_{\widetilde{\Omega}}$,$N_W$) &     Examples  &  Abs. Mean &  Variance \\
\hline 
Ia   & 4  &(1, 0, 3)   & $\sum_{\substack{i,j,k,\ell\\ (dist)}} \delta_{ij} W_{jk} W_{k \ell} W_{\ell i}$ & 0  & $\leq C\|\theta\|^2\|\theta\|_3^6=o(\|\theta\|^8)$    \\  
Ib  &  8  &(1, 1, 2)   & $\sum_{\substack{i, j, k, \ell\\ (dist)}} \delta_{ij} \widetilde{\Omega}_{jk}  W_{k \ell}  W_{\ell i}$& 0 &  $\leq C\alpha^2\|\theta\|^4\|\theta\|_3^6=o(\|\theta\|^8)$  \\ 
   & 4  &    & $\sum_{\substack{i, j, k, \ell\\ (dist)}} \delta_{ij} W_{jk}   \widetilde{\Omega}_{k \ell}  W_{\ell i}$& 0 & $\leq C\alpha^2\|\theta\|^4\|\theta\|_3^6=o(\|\theta\|^8)$   \\ 
Ic & 8 &(1, 2, 1)   & $\sum_{\substack{i, j, k, \ell\\ (dist)}} \delta_{ij}  \widetilde{\Omega}_{jk} \widetilde{\Omega}_{k \ell}W_{\ell i}$ & $\leq C\alpha^2\|\theta\|^6$=$o(\alpha^4\|\theta\|^8)$ & $\leq \frac{C\alpha^4\|\theta\|^{10}\|\theta\|_3^3}{\|\theta\|_1}=o(\alpha^6\|\theta\|^8\|\theta\|_3^6)$  \\   
  & 4 &  & $\sum_{\substack{i, j, k, \ell\\ (dist)}} \delta_{ij}  \widetilde{\Omega}_{jk} W_{k \ell}  \widetilde{\Omega}_{\ell i}$ & 0 & $\leq \frac{C\alpha^4\|\theta\|^4\|\theta\|_3^9}{\|\theta\|_1}=o(\|\theta\|^8)$  \\   
Id & 4  & (1, 3, 0) & $\sum_{\substack{i, j, k, \ell\\ (dist)}} \delta_{ij}  \widetilde{\Omega}_{jk} \widetilde{\Omega}_{k \ell}  \widetilde{\Omega}_{\ell i}$  & 0  & $\leq \frac{C\alpha^6\|\theta\|^{12}\|\theta\|_3^3}{\|\theta\|_1}=O(\alpha^6\|\theta\|^8\|\theta\|_3^6)$  \\
IIa & 4 & (2, 0, 2)   & $\sum_{\substack{i, j, k, \ell\\ (dist)}} \delta_{ij} \delta_{jk} W_{k \ell}  W_{\ell i}$  &  $\leq C\|\theta\|^4$=$o(\alpha^4\|\theta\|^8)$ &  $\leq C \|\theta\|^2\|\theta\|_3^6=o(\|\theta\|^8)$   \\  
     & 2 &     & $\sum_{\substack{i, j, k, \ell\\ (dist)}} \delta_{ij} W_{jk} \delta_{k \ell}  W_{\ell i}$    &  $\leq C\|\theta\|^4$=$o(\alpha^4\|\theta\|^8)$  &  $\leq \frac{C\|\theta\|^6\|\theta\|_3^3}{\|\theta\|_1}=o(\|\theta\|^8)$    \\ 
IIb & 8 & (2, 1, 1) & $\sum_{\substack{i, j, k, \ell\\ (dist)}} \delta_{ij} \delta_{jk} \widetilde{\Omega}_{k \ell} W_{\ell i}$    & $0$ &  $\leq C\alpha^2\|\theta\|^4\|\theta\|_3^6=o(\|\theta\|^8)$   \\ 
& 4 &  & $\sum_{\substack{i, j, k, \ell\\ (dist)}} \delta_{ij}  \widetilde{\Omega}_{j k} \delta_{k \ell}  W_{\ell i}$    & $\leq C\alpha\|\theta\|^4$=$o(\alpha^4\|\theta\|^8)$  &  $\leq \frac{C\alpha^2\|\theta\|^8\|\theta\|_3^3}{\|\theta\|_1} = o(\|\theta\|^8)$   \\
IIc &  4 &  (2,  2,  0)  & $\sum_{\substack{i, j, k, \ell\\ (dist)}} \delta_{ij} \delta_{jk} \widetilde{\Omega}_{k \ell} \widetilde{\Omega}_{\ell i}$ &  $\leq C\alpha^2\|\theta\|^6$=$o(\alpha^4\|\theta\|^8)$ & $\leq\frac{C\alpha^4\|\theta\|^{14}}{\|\theta\|_1^2}=o(\alpha^6\|\theta\|^8\|\theta\|_3^6)$ \\ 
    &  2 &    & $\leq \sum_{\substack{i, j, k, \ell\\ (dist)}} \delta_{ij} \widetilde{\Omega}_{jk} \delta_{k \ell} \widetilde{\Omega}_{\ell i}$ &  $\frac{C\alpha^2\|\theta\|^8}{\|\theta\|_1^2}$=$o(\alpha^4\|\theta\|^8)$  & $\leq \frac{C\alpha^4\|\theta\|^8\|\theta\|_3^6}{\|\theta\|_1^2}=o(\|\theta\|^8)$  \\ 
IIIa & 4 &  (3,  0, 1)  & $\sum_{\substack{i, j, k, \ell\\ (dist)}} \delta_{ij} \delta_{jk} \delta_{k\ell} W_{\ell i}$ &  $\leq C\|\theta\|^4$=$o(\alpha^4\|\theta\|^8)$ & $\leq \frac{C\|\theta\|^6 \|\theta\|_3^3}{\|\theta\|_1}=o(\|\theta\|^8)$  \\ 
IIIb & 4 &  (3,  1,  0)  & $\leq \sum_{\substack{i, j, k, \ell\\ (dist)}} \delta_{ij} \delta_{jk} \delta_{k\ell} \widetilde{\Omega}_{\ell i}$ &
$\leq \frac{C\alpha\|\theta\|^6}{\|\theta\|_1^3}$=$o(\alpha^4\|\theta\|^8)$   & $\leq \frac{C\alpha^2 \|\theta\|^8 \|\theta\|_3^3}{\|\theta\|_1}=o(\|\theta\|^8)$  \\ 
IV & 1 &  (4,  0,  0)  & $\sum_{\substack{i, j, k, \ell\\ (dist)}} \delta_{ij} \delta_{jk} \delta_{k\ell} \delta_{\ell i}$ &  $\leq C\|\theta\|^4$=$o(\alpha^4\|\theta\|^8)$  & $\leq \frac{C\|\theta\|^{10}}{\|\theta\|_1^2}=o(\|\theta\|^8)$ \\ 
\hline 
\end{tabular} 
} 
\end{table}

Consider the null hypothesis first. Under the null, $\widetilde{\Omega}$ is a zero matrix, so the nonzero post-expansion sums only include Type Ia, Type IIa, Type IIIa, and Type IV. It is seen that  
\[
|\mathbb{E}[Q_n^* - \widetilde{Q}_n] | \leq C\|\theta\|^4, 
\] 
and that
\[
\mathrm{Var}(Q_n^*- \widetilde{Q}_n) = o(\|\theta\|^8). 
\]
Note that $\|\theta\|^8$ is the order of $\mathrm{Var}(\widetilde{Q}_n)$ under the null. The difference between the variance of $Q^*_n$ and the variance of $\widetilde{Q}_n$ is negligible, but the difference between the mean of $Q_n^*$ and the mean of $\widetilde{Q}_n$ is non-negligible. With lengthy calculations (see the  appendix), we can show that
\[
\mathbb{E}[Q_n^* - \widetilde{Q}_n] \sim  2\|\theta\|^4.
\]
Therefore, $(Q_n^*-2\|\theta\|^4)$ and $\widetilde{Q}_n$ have a negligible difference under the null. 

Consider the alternative hypothesis next. From Table~\ref{tab:ProxySgnQsum}, 
\[
|\mathbb{E}[Q_n^* - \widetilde{Q}_n] | \leq C(|\lambda_2|/\lambda_1)^2\|\theta\|^6,
\]
where the major contribution is from Type Ic and Type IIc post-expansion sums. Under our assumptions for the alternative, $|\lambda_2|/\sqrt{\lambda_1}\to\infty$ and $\lambda_1 \asymp\|\theta\|^4$. It is easy to see that
$|\mathbb{E}[Q_n^* - \widetilde{Q}_n] | = o(\lambda_2^4)$, 
where $\lambda_2^4$ is the order of $\tr(\widetilde{\Omega}^4)$ and $\mathbb{E}[\widetilde{Q}_n]$; see Lemma~\ref{lem:tracetOmega2} and Theorem~\ref{thm:IdealSgnQ}. Additionally, $\|\theta\|^4=O(\lambda_1^2)=o(\lambda_2^4)$, which is also of a smaller order of $\mathbb{E}[\widetilde{Q}_n]$. We conclude that
\[
\bigl|\mathbb{E}[Q_n^* -\widetilde{Q}_n-2\|\theta\|^4]\bigr| = o(\mathbb{E}[\widetilde{Q}_n]). 
\]
From the table, 
\[
\mathrm{Var}(Q_n^* - \widetilde{Q}_n)  \leq \frac{C(|\lambda_2|/\lambda_1)^6\|\theta\|^{12}\|\theta\|_3^3}{\|\theta\|_1} + o(\|\theta\|^8),
\]  
with the major contribution from Type Id. Here, the second term is smaller than $\mathrm{Var}(\widetilde{Q}_n)$, and the first term is upper bounded by (using the universal inequality of $\|\theta\|^4\leq \|\theta\|_1\|\theta\|_3^3$)
\[
C(|\lambda_2|/\lambda_1)^6\|\theta\|^8\|\theta\|_3^6,
\]  
which has a comparable order as $\mathrm{Var}(\widetilde{Q}_n)$. It follows that
\[
\mathrm{Var}(Q_n^* -\widetilde{Q}_n-2\|\theta\|^4) = \mathrm{Var}(Q_n^* - \widetilde{Q}_n)\leq C\mathrm{Var}(\widetilde{Q}_n). 
\]
Combining the above, we obtain that the SNR of $(Q_n^*-2\|\theta\|^4)$ and $\widetilde{Q}_n$ are at the same order. 

These results are summarized in the following theorem, which is proved in the  
 appendix.

\begin{thm}[Proxy SgnQ test statistic] \label{thm:ProxySgnQ} 
Consider the testing problem \eqref{Problem} under the DCMM model \eqref{model1a}-\eqref{condition1d}, where the condition \eqref{cond-balance} is satisfied under the alternative hypothesis. Suppose $\theta_{\max}\to 0$ and $\|\theta\| \goto \infty$ as $n \goto \infty$, and suppose $|\lambda_2|/\sqrt{\lambda_1}\to\infty$ under the alternative hypothesis. Then, under the null hypothesis, as $n \goto \infty$, 
\[
\mathbb{E}[(Q_n^*-2\|\theta\|^4)-\widetilde{Q}_n] =o(\|\theta\|^4), \qquad \mathrm{Var}(Q_n^* - \widetilde{Q}_n)= o(\|\theta\|^8).  
\]
Furthermore, under the alternative hypothesis, 
\begin{align*}
& \mathbb{E}[(Q_n^*-2\|\theta\|^4)-\widetilde{Q}_n] = o((|\lambda_2|/\lambda_1)^4\|\theta\|^8), \\
& \mathrm{Var}(Q_n^* - \widetilde{Q}_n) \leq C(|\lambda_2|/\lambda_1)^6\|\theta\|^8\|\theta\|_3^6 + o(\|\theta\|^8). 
\end{align*}
\end{thm} 

\subsection{The behavior of $(T_n - T_n^*)$ and $(Q_n - Q_n^*)$} 
\label{subsec:RealSgn} 
The SgnT statistic and SgnQ statistic we introduce in Section~\ref{subsec:SP} are defined as
\[
T_n  =\sum_{i_1,i_2,i_3 (dist)}X_{i_1i_2}X_{i_2i_3}X_{i_3i_1}, \qquad Q_n  =\sum_{i_1,i_2,i_3,i_4 (dist)}X_{i_1i_2}X_{i_2i_3}X_{i_3i_4}X_{i_4i_1},
\] 
where $X_{ij}=\widetilde{\Omega}_{ij}+W_{ij}+\delta_{ij}+r_{ij}$ for any $i\neq j$. Similar to Sections~\ref{subsec:IdealSgn}-\ref{subsec:ProxySgn},  we first expand every bracket in the definitions and obtain 
\[
4 \times 4 \times 4 = 64, \qquad\mbox{and}\qquad 4\times 4\times 4\times 4=256
\]  
different post-expansion sums in $T_n$ and $Q_n$, respectively. Out of the $64$ post-expansion sums in $T_n$, $3 \times 3 \times 3 = 27$ of them do not involve any $r$ term  and are contained in $T_n^*$; this leaves a total of 
\[
64 - 27 = 37
\] 
different post-expansion sums in $(T_n-T_n^*)$. 
Out of the $256$ post-expansion sums in $Q_n$, $3 \times 3 \times 3\times 3 = 81$ of them do not involve any $r$ term  and are contained in $Q_n^*$; this leaves a total of 
\[
256 - 81 = 175
\]
different post-expansion sums in $(Q_n-Q_n^*)$. In the  appendix, we investigate the order of mean and variance of each of $37$ post-expansion sums in $(T_n-T_n^*)$ and each of $175$ post-expansion sums in $(Q_n-Q_n^*)$. The calculations are very tedious: although we expect these post-expansion sums to be of a smaller order than the post-expansion sums in Sections~\ref{subsec:IdealSgn}-\ref{subsec:ProxySgn}, it is impossible to prove this argument rigorously using only some crude bounds (such as Cauchy-Schwarz inequality). Instead, we still need to do calculations for each post-expansion sum; details are in the  appendix. 
\begin{thm}[Real SgnT test statistic] \label{thm:RealSgnT}
Consider the testing problem \eqref{Problem} under the DCMM model \eqref{model1a}-\eqref{condition1d}, where the condition \eqref{cond-balance} is satisfied under the alternative hypothesis. Suppose $\theta_{\max}\to 0$ and $\|\theta\| \goto \infty$ as $n \goto \infty$, and suppose $|\lambda_2|/\sqrt{\lambda_1}\to\infty$ under the alternative hypothesis. Then, under the null hypothesis, as $n \goto \infty$,  
\[
|\mathbb{E}[T_n - T_n^*]| = o(\|\theta\|^3),  
\qquad  \mbox{and} \qquad  
\mathrm{Var}(T_n - T_n^*) =o(|\theta\|^6). 
\] 
Under the alternative hypothesis, as $n\to\infty$, 
\begin{align*}
& |\mathbb{E}[T_n - T_n^*]| = o((|\lambda_2|/\lambda_1)^3\|\theta\|^6), \\
&  \mathrm{Var}(T_n - T_n^*) =o((|\lambda_2|/\lambda_1)^4\|\theta\|^4|\theta\|_3^6) +o(\|\theta\|^6 ). 
\end{align*}
\end{thm} 

Note that $\|\theta\|^3$ is the order of the standard deviation under the null. The theorem implies that the difference between $T_n$ and $T_n^*$ is negligible. In Section~\ref{subsec:ProxySgn}, we have seen the difference between $T_n^*$ and $\widetilde{T}_n$ is negligible. Therefore, $T_n$ has a negligible difference from $\widetilde{T}_n$, and so $T_n$ enjoys the same asymptotic normality under the null as $\widetilde{T}_n$ does. Under the alternative, $(|\lambda_2|/\lambda_1)^3\|\theta\|^6$ is the order of $\mathbb{E}[\widetilde{T}_n]$, and $(|\lambda_2|/\lambda_1)^4\|\theta\|^4|\theta\|_3^6+\|\theta\|^6$ is the order of $\mathrm{Var}(\widetilde{T}_n)$. Therefore, the above theorem implies that $(T_n-T_n^*)$ has a negligible effect to the SNR. 
\begin{thm}[Real SgnQ test statistic] \label{thm:RealSgnQ}
Consider the testing problem \eqref{Problem} under the DCMM model \eqref{model1a}-\eqref{condition1d}, where the condition \eqref{cond-balance} is satisfied under the alternative hypothesis. Suppose $\theta_{\max}\to 0$ and $\|\theta\| \goto \infty$ as $n \goto \infty$, and suppose $|\lambda_2|/\sqrt{\lambda_1}\to\infty$ under the alternative hypothesis. Then, under the null hypothesis, as $n \goto \infty$,  
\[
|\mathbb{E}[Q_n - Q_n^*]| = o(\|\theta\|^4),  
\qquad  \mbox{and} \qquad  
\mathrm{Var}(Q_n - Q_n^*) =o(|\theta\|^8). 
\] 
Under the alternative hypothesis, as $n\to\infty$, 
\begin{align*}
& |\mathbb{E}[Q_n - Q_n^*]| = o((|\lambda_2|/\lambda_1)^4\|\theta\|^8), \\
&  \mathrm{Var}(Q_n - Q_n^*) =o((|\lambda_2|/\lambda_1)^6\|\theta\|^8|\theta\|_3^6) +o(\|\theta\|^8). 
\end{align*}
\end{thm}

Similarly, we can conclude that $(Q_n-Q_n^*)$ has a negligible effect to both the asymptotic normality under the null and the SNR under the alternative.   

\subsection{Proof of the main theorems}  \label{subsec:shortmainpf} 
We now prove Theorems \ref{thm:null-SgnT}-\ref{thm:SgnQtest}.  Consider Theorems \ref{thm:null-SgnT}-\ref{thm:alt-SgnQ}  first. Since the proofs are similar, 
we only show Theorems \ref{thm:null-SgnQ} and \ref{thm:alt-SgnQ}. 

Consider Theorem \ref{thm:null-SgnQ}. In this theorem, we assume the null is true. 
First, by Theorems \ref{thm:ProxySgnQ} and \ref{thm:RealSgnQ} and elementary statistics,   
\begin{equation} \label{mainpfA1a} 
 \mathbb{E}[Q_n^* - \widetilde{Q}_n]  \sim 2 \|\theta\|^4,  \qquad 
|\mathbb{E}[Q_n - Q_n^*]|  = o(\|\theta\|^4), 
\end{equation} 
and 
\begin{equation} \label{mainpfA1b} 
\mathrm{Var}(Q_n^* - \widetilde{Q}_n) = o(\|\theta\|^8), \qquad 
\mathrm{Var}(Q_n    - Q_n^*) = o(\|\theta\|^8).  
\end{equation} 
It follows that 
\begin{equation} \label{mainpfA1} 
\mathbb{E}[Q_n]  - \mathbb{E}[\widetilde{Q}_n] = (2 + o(1)) \|\theta\|^4, 
 \qquad  \mathrm{Var}(Q_n - \widetilde{Q}_n) = o(\|\theta\|^8).  
\end{equation} 
By Theorem \ref{thm:IdealSgnQ}.  
\begin{equation} \label{mainpfA2} 
\mathbb{E}[\widetilde{Q}_n] = o(\|\theta\|^4), \;\;\;   \mathrm{Var}(\widetilde{Q}_n) \sim 8 \|\theta\|^8, \;\;\;  
\frac{\widetilde{Q}_n - \mathbb{E}[\widetilde{Q}_n]}{\sqrt{\mathrm{Var}(\widetilde{Q}_n)}} \goto N(0,1). 
\end{equation} 
Since for any random variables $X$ and $Y$, $\mathrm{Var}(X + Y) \leq (1 + a_n) \mathrm{Var}(X) + (1 + \frac{1}{a_n}) \mathrm{Var}(Y)$ for any number $a_n > 0$, combining the above and letting $a_n$ tend  to $0$ appropriately slow,  
\begin{equation} \label{mainpfA3}
\mathbb{E}[Q_n] \sim 2 \|\theta\|^4, \qquad 
\mathrm{Var}(Q_n) \sim 8 \|\theta\|^8.   
\end{equation}  
Moreover, write 
\[ 
\frac{Q_n - \mathbb{E}[Q_n]}{\sqrt{\mathrm{Var}(Q_n)}} = \sqrt{\frac{\mathrm{Var}(\widetilde{Q}_n)}{\mathrm{Var}(Q_n)}} \cdot   \biggl[ \frac{(Q_n - \widetilde{Q}_n)}{\sqrt{\mathrm{Var}(\widetilde{Q}_n)}} +  \frac{\widetilde{Q}_n - \mathbb{E}[\widetilde{Q}_n]}{\sqrt{\mathrm{Var}(\widetilde{Q}_n)}}  +  \frac{\mathbb{E}[\widetilde{Q}_n] - \mathbb{E}[Q_n]}{\sqrt{\mathrm{Var}(\widetilde{Q}_n)}} \biggr]. 
\] 
On the right hand side,  by (\ref{mainpfA1})-(\ref{mainpfA3}), as $n \goto \infty$, the term outside the bracket $\goto 1$, and for the three terms in the bracket, the first one has a mean and variance that tend to $0$ so it tends to $0$ in probability, the second one weakly converges to $N(0,1)$, and the last one $\goto 0$. Combining these,  
\begin{equation} \label{mainpfA4} 
\frac{Q_n - \mathbb{E}[Q_n]}{\sqrt{\mathrm{Var}(Q_n)}} \goto N(0,1), \qquad \mbox{in law}. 
\end{equation} 
Combining (\ref{mainpfA3}) and (\ref{mainpfA4}) proves Theorem \ref{thm:null-SgnQ}. 

Next, we consider Theorem \ref{thm:alt-SgnQ}, where we assume the alternative is true. First, similarly, by 
Theorems \ref{thm:ProxySgnQ} and \ref{thm:RealSgnQ},  
\[
\mathbb{E}[Q_n^* - \widetilde{Q}_n]   = (2 + o(1)) \|\theta\|^4 + o((|\lambda_2|/\lambda_1)^4 \|\theta\|^8), 
\] 
and
\[ 
\mathrm{Var}(Q_n - \widetilde{Q}_n) \leq 
C(\lambda_2 / \lambda_1)^6 \|\theta\|^8 \|\theta\|_3^6 + o(\|\theta\|^8).   
\] 
Second, by Theorems \ref{thm:IdealSgnT} and \ref{thm:IdealSgnQ},  
\[
\mathbb{E}[\widetilde{Q}_n] = \tr(\widetilde{\Omega}^4)  +  o(\|\theta\|^4),    \qquad 
\mathrm{Var}(\widetilde{Q}_n) \leq C[\|\theta\|^8 + (\lambda_2/\lambda_1)^6 \|\theta\|^8 \|\theta\|_3^6]. 
\] 
Combining these proves Theorem \ref{thm:alt-SgnQ}. 

Last, we consider Theorems \ref{thm:SgnTtest}-\ref{thm:SgnQtest}. 
Since the proofs are similar, we only show Theorem \ref{thm:SgnQtest}. 
First,  by Theorem \ref{thm:null-SgnQ} and Lemma \ref{lemma:heta}, under the null,  
\[
\frac{Q_n  - 2(\|\hat{\eta}\|^2-1)^2}{\sqrt{8 (\|\hat{\eta}\|^2 -1)^4}} \goto N(0,1), 
\]
so the Type I error is 
\[
\mathbb{P}_{H_0^{(n)}}\biggl(Q_n \geq (2 + z_{\alpha} \sqrt{8})   (|\hat{\eta}\|^2 - 1)^2 \biggr) = P\biggl(\frac{Q_n - 2 (\|\hat{\eta}\|^2 - 1)^2}{\sqrt{8 (\|\hat{\eta}\|^2 - 1)^4}} \geq  z_{\alpha} \biggr)   = \alpha + o(1). 
\] 
Second, fixing $0 < \eps < 1$, let $A_{\eps}$ be the event  $\{(\|\hat{\eta}\|^2 - 1) \leq (1 + \eps)  \|\eta^*\|^2\}$.  
By Lemma \ref{lemma:heta} and definitions,  on one hand,  over the event $A_{\eps}$, $(\|\hat{\eta}\|^2 - 1) \leq (1 + \eps) \|\eta^*\|^2 \leq C \|\theta\|^2$, and on the other hand, 
 $\mathbb{P}(A_{\eps}^c) = o(1)$.   Therefore,  the Type II error 
\begin{align*} 
& \mathbb{P}_{H_1^{(n)}}\biggl(Q_n \leq (2 + z_{\alpha} \sqrt{8}) (\|\hat{\eta}\|^2 - 1)^2 \biggr)  \\
\leq & \mathbb{P}_{H_1^{(n)}}\biggl(Q_n \leq (2 + z_{\alpha} \sqrt{8}) (\|\hat{\eta}\|^2 - 1)^2, A_{\eps}\biggr) + \mathbb{P}(A_{\eps}^c)  \\
\leq & \mathbb{P}_{H_1^{(n)}} \biggl(Q_n \leq C  (2 + z_{\alpha} \sqrt{8}) \|\theta\|^4 \biggr) + o(1),  
\end{align*} 
where by Chebyshev's inequality, the first term in the last line 
\begin{equation} \label{mypfC1}
\leq  [\mathbb{E}(Q_n) - C  (2 + z_{\alpha} \sqrt{8}) \|\theta\|^4]^{-2} \cdot \mathrm{Var}(Q_n). 
\end{equation} 
By Lemma D.2 of the  appendix and our assumptions,   $\lambda_1 \asymp \|\theta\|^2$,  $|\lambda_2|/\sqrt{\lambda_1} \goto \infty$, and $\|\theta\| \goto \infty$. 
Using Lemma \ref{lem:tracetOmega2}  
$\mathbb{E}[Q_n] \geq C \lambda_2^4  \gg \lambda_1^2$,    
and it follows that 
$\mathbb{E}(Q_n)  \gg C  (2 + z_{\alpha} \sqrt{8}) \|\theta\|^4$,   
so for sufficiently large $n$, 
\[
\mathbb{E}(Q_n) - C  (2 + z_{\alpha} \sqrt{8}) \|\theta\|^4  \geq \frac{1}{2}  \mathbb{E}[Q_n] \geq C \lambda_2^4. 
\] 
At the same time, by Theorem \ref{thm:alt-SgnQ}, 
\[
\mathrm{Var}(Q_n) \leq C (\|\theta\|^8 + (\lambda_2 / \lambda_1)^6 \|\theta\|^8 \|\theta\|_3^6).  
\] 
Combining these, the right hand side of (\ref{mypfC1}) does not exceed 
\begin{equation} \label{mypfC2} 
C \frac{\|\theta\|^8 + (\lambda_2/\lambda_1)^6 \|\theta\|^8 \|\theta\|_3^6}{\lambda_2^8}  =   (I) + (II), 
\end{equation} 
where  
$(I) = C \lambda_2^{-8} \|\theta\|^8$ and $(II) = C \lambda_2^{-8} (\lambda_2/\lambda_1)^6 \|\theta\|^8 \|\theta\|_3^6$. 
Now, first, since $\lambda_1 \asymp \|\theta\|^2$ and $|\lambda_2|/\sqrt{\lambda_1} \goto 0$, 
$(I)  \leq C (\lambda_2 / \sqrt{\lambda_1})^{-8} \goto 0$.  
Second, since $\lambda_1 \asymp \|\theta\|^2$ and $\|\theta\|_3^6 \leq \|\theta\|^4$,  
$(II) = C \lambda_2^{-2}  \lambda_1^{-6} \|\theta\|^8 \|\theta\|_3^6 \leq C \lambda_2^{-2}$.    
As  $|\lambda_2|/\sqrt{\lambda_1}  \goto \infty$, $\sqrt{\lambda_1} \asymp \|\theta\|$ with $\|\theta\| \goto \infty$, 
$|\lambda_2| \goto \infty$ and  $(II) \goto 0$. 
Inserting these into (\ref{mypfC2}), the Type II error $\goto 0$ and the claim follows.  \qed

\section{Simulations} \label{sec:Simul}
We investigate the numerical performance of two Signed Polygon tests, the SgnT test \eqref{SgnTtest} and the SgnQ test \eqref{SgnQtest}. We also include the EZ test  \cite{gao2017testing} and the GC test \cite{OGC} for comparison. For reasons mentioned in \cite{OGC}, we use a two-sided rejection region for EZ and a one-sided rejection region for GC. 

Given $(n,K)$, a scalar $\beta_n>0$ that controls $\|\theta\|$, a symmetric nonnegative matrix $P\in\mathbb{R}^{K\times K}$, a distribution $f(\theta)$ on $\mathbb{R}_+$, and a distribution $g(\pi)$ on the standard simplex of $\mathbb{R}^K$, we generate two network adjacency matrices $A^{null}$ and $A^{alt}$, under the null and the alternative, respectively, as follows: 
\begin{itemize}
\item Generate $\tilde{\theta}_1,\tilde{\theta}_2,\ldots,\tilde{\theta}_n$ $iid$ from $f(\theta)$. Let $\theta_i=\beta_n\cdot \tilde{\theta}_i/\|\tilde{\theta}\|$, $1\leq i\leq n$.
\item Generate $\pi_1,\pi_2,\ldots,\pi_n$ $iid$ from $g(\pi)$.
\item Let $\Omega^{alt}=\Theta\Pi P\Pi'\Theta'$, where $\Theta=\mathrm{diag}(\theta_1,\cdots,\theta_n)$ and $\Pi=[\pi_1,\pi_2,\ldots,\pi_n]'$. Generate $A^{alt}$ from $\Omega^{alt}$ according to Model \eqref{model1a}.
\item Let $\Omega^{null}= (a' P a) \cdot \theta \theta'$, where $a = \mathbb{E}_g\pi\in\mathbb{R}^K$ is the mean vector of $g(\pi)$. Generate $A^{null}$ from $\Omega^{null}$ according to Model \eqref{model1a}.
\end{itemize}
The pair $(\Omega^{null},\Omega^{alt})$ is constructed in a way such that the corresponding networks have approximately the same expected average degree. This is the most subtle case for distinguishing two hypotheses (see Section~\ref{sec:LB}). 

It is of interest to explore different sparsity levels and also to focus on the parameter settings 
where the SNR is neither too large or too small. Therefore, for most of the experiments, we let 
$\beta_n = \|\theta\|$ range but fix the SNR at a more or less the same level. See details below.

For each parameter setting, we generate $200$ networks under the null hypothesis and $200$ networks under the alternative hypothesis, run all the four tests with a targeting level $\alpha=5\%$, and then record the sum of percent of type I errors and percent of type II errors.  

We consider three experiments (and a total of $8$ sub-experiments), exploring different sets of $n$, $K$, $\theta$, $\Pi$, and $P$, etc. 

\begin{figure}[tb!]
\includegraphics[width=.32\textwidth, height = .32\textwidth]{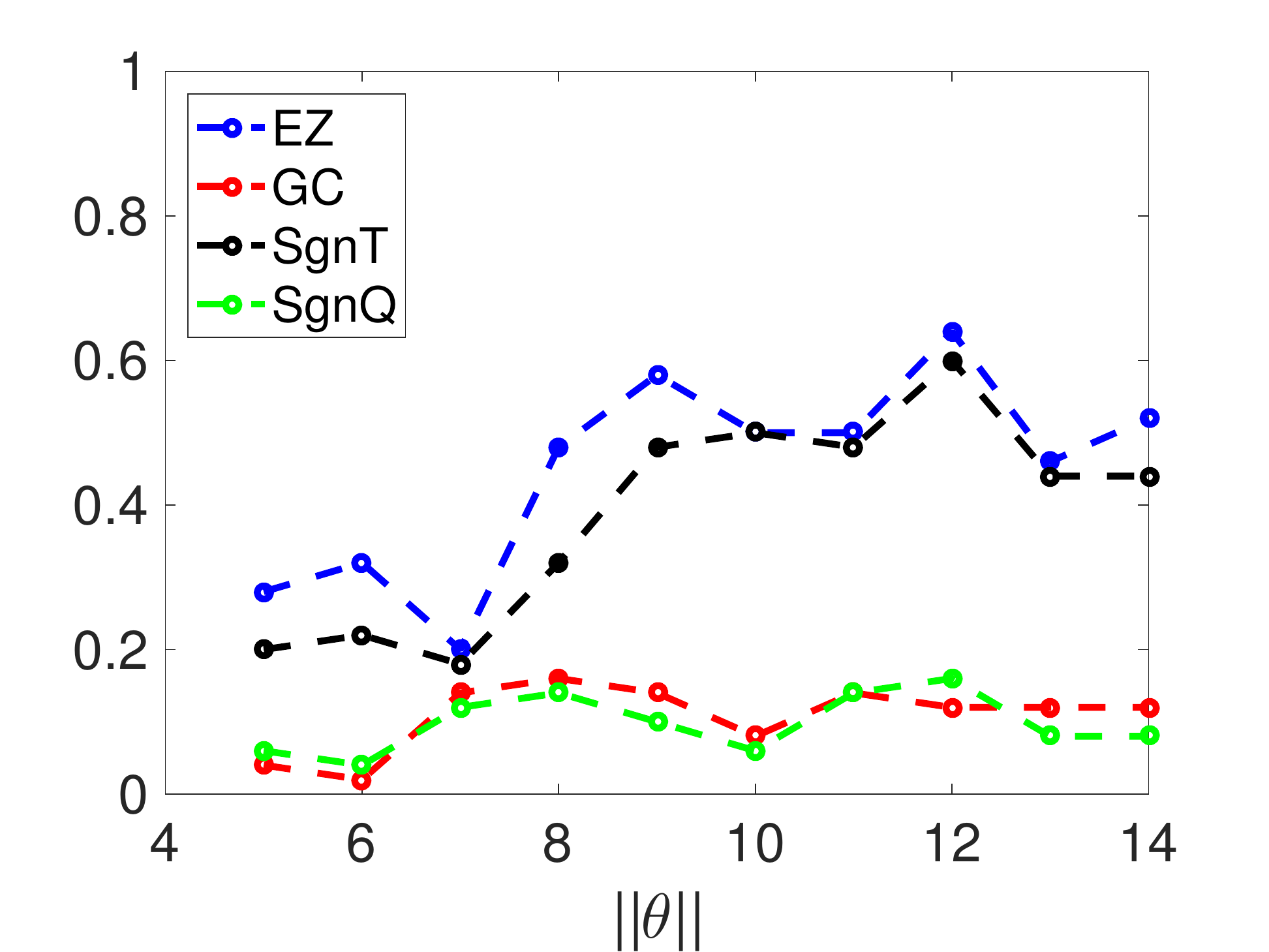}
\includegraphics[width=.32\textwidth, height = .32\textwidth]{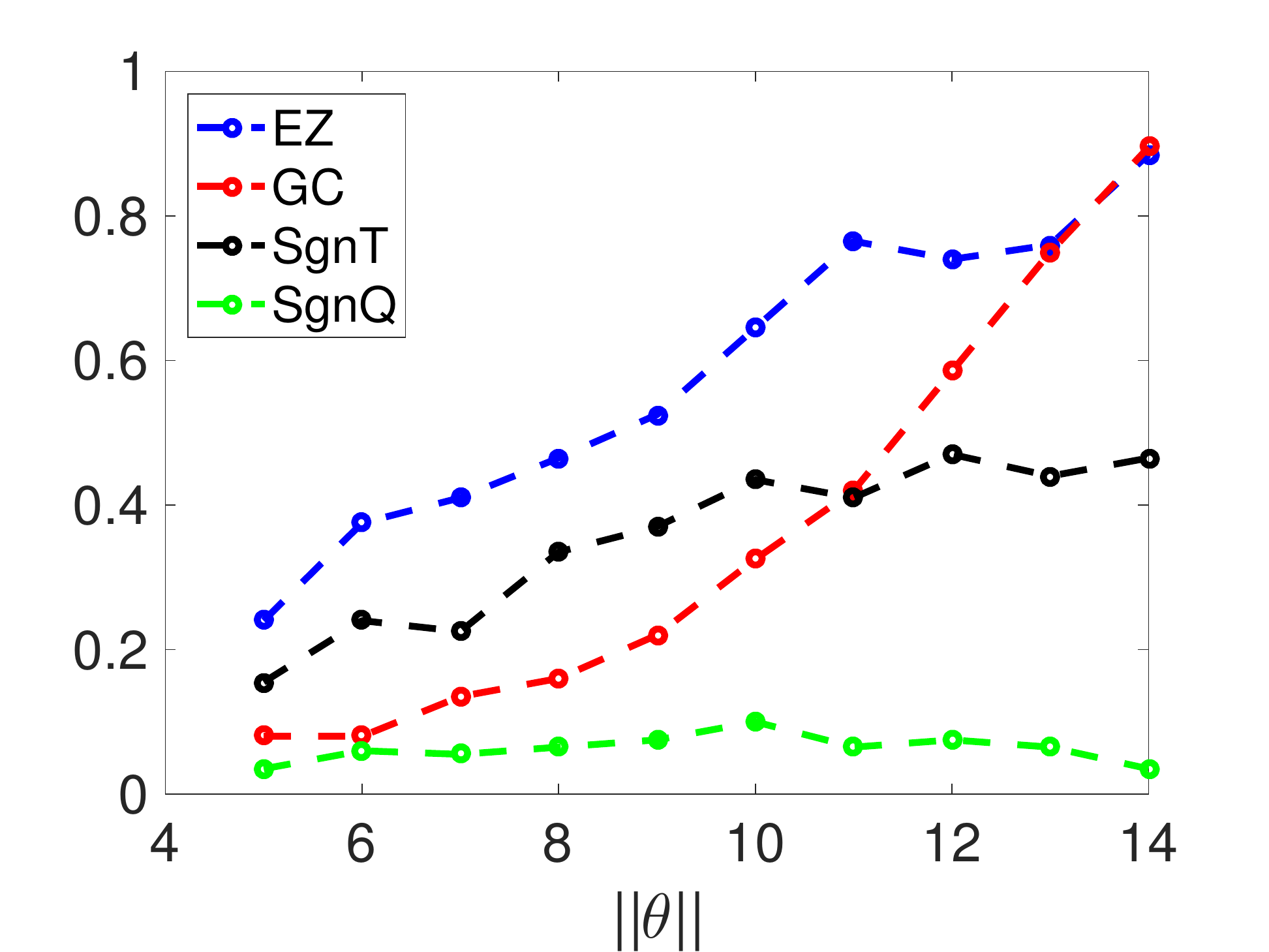}
\includegraphics[width=.32\textwidth, height = .32\textwidth]{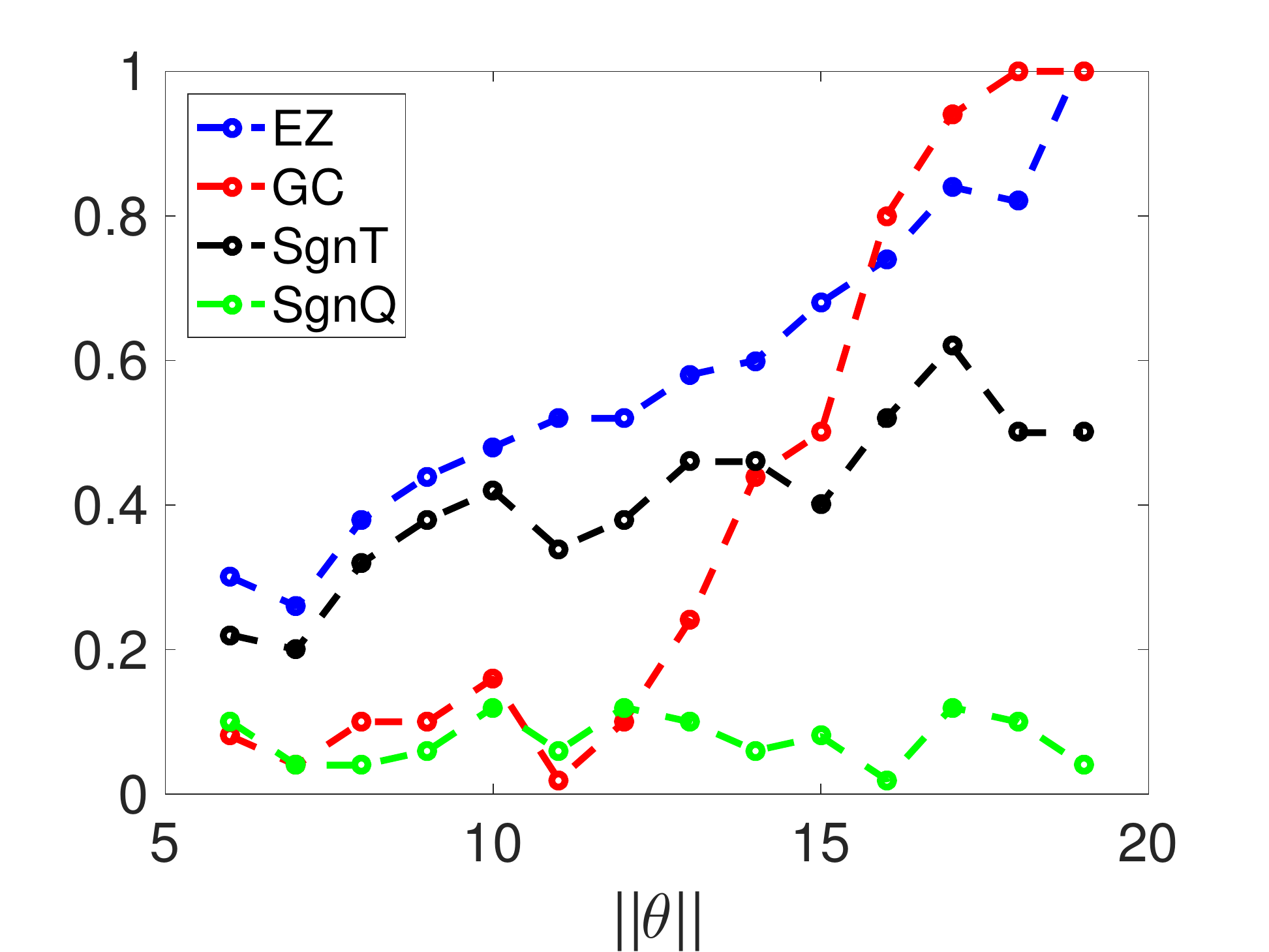}
\caption{From left to right: Experiment 1a, 1b, and 1c. The $y$-axis are the sum of Type I and Type II errors (testing level is fixed at $5\%$). The $x$-axis are $\|\theta\|$ or sparsity levels.   Results are based on $200$ repetitions.} \label{fig:Exp1}
\end{figure}

\paragraph{Experiment 1} We study the role of  degree heterogeneity.  
Fix $(n,K)=(2000,2)$. Let $P$ be a $2\times 2$ matrix with unit diagonal entries and all off-diagonal entries  equal to $b_n$. Let $g(\pi)$ be the uniform distribution on $\{(0,1), (1,0)\}$. We consider three sub-experiments, Exp 1a-1c, where respectively we take $f(\theta)$ to be the following:    
\begin{itemize}
\item Uniform distribution $U(2,3)$. 
\item Two-point distribution $0.95\delta_1 + 0.05\delta_3$, where $\delta_a$ is a point mass at $a$.   
\item Pareto distribution $\mathrm{Pareto}(10, 0.375)$, where $10$ is the shape parameter and $0.375$ is the scale parameter. 
\end{itemize}
The degree heterogeneity is moderate in the Exp 1a-1b, but  more severe Exp 1c.  In such a setting, SNR is at the order of $\|\theta\| (1 - b_n)$.  
Therefore, for each sub-experiment, we  let $\beta_n=\|\theta\|$ vary while fixing the SNR to be 
\[
\|\theta\| (1 - b_n)  = 3.2.
\]
The sum of Type I and Type II errors are displayed in Figure~\ref{fig:Exp1}. 

First, both the SgnQ test and the GC test are based on the counts of $4$-cycles, but the $GC$ test counts {\it non-centered} cycles and the SgnQ test counts {\it centered} cycles. As we pointed out in  Section~\ref{sec:intro},  counting {\it centered} cycles may have much smaller variances 
than counting {\it non-centered} cycles, 
especially in the less sparse case, and thus improves the testing power.  This is confirmed by numerical results here, where the SgnQ test is consistently better than the GC test, significantly so in the less sparse case. Similarly,  both the SgnT test and the EZ test are based on the counts of $3$-cycles, but the $EZ$ test counts {\it non-centered} cycles and the SgnQ test counts {\it centered} cycles, and 
we expect the GC to significantly improve the EZ, especially in the less sparse case. This is also confirmed in the experiment. 

Second,  SgnQ and GC are order-$4$ test graphlet counting statistics, and 
SgnT and EZ are order-$3$ graphlet counting statistics. In comparison, SgnQ significantly outperforms  
SgnT, and GC significantly outperforms  EZ (in the more sparse case; see discussion for the 
less sparse case below). A possible explanation is that higher-order graphlet counting statistics 
have larger SNR. Investigation on this is interesting, and we leave this to the future study.  Note that SgnQ is the best among all four tests. 

Last, GC outperforms EZ in the more sparse case, but  underperforms in the less sparse case. The reason for the latter is as follows. The biases of both tests are negligible in the more sparse case, but are non-negligible in the less sparse case, with that of GC is much larger. 
In \cite{Ke2019}, we propose a bias correction 
method, where the performance of GC is significantly improved. 
However, GC continues to underperform SgnQ, because even with the bias corrected, 
it still has a variance that is unsatisfactorily large.

\begin{figure}[tb] 
\centering
\vspace{-0mm}
\includegraphics[width=.42\textwidth, height = .35\textwidth]{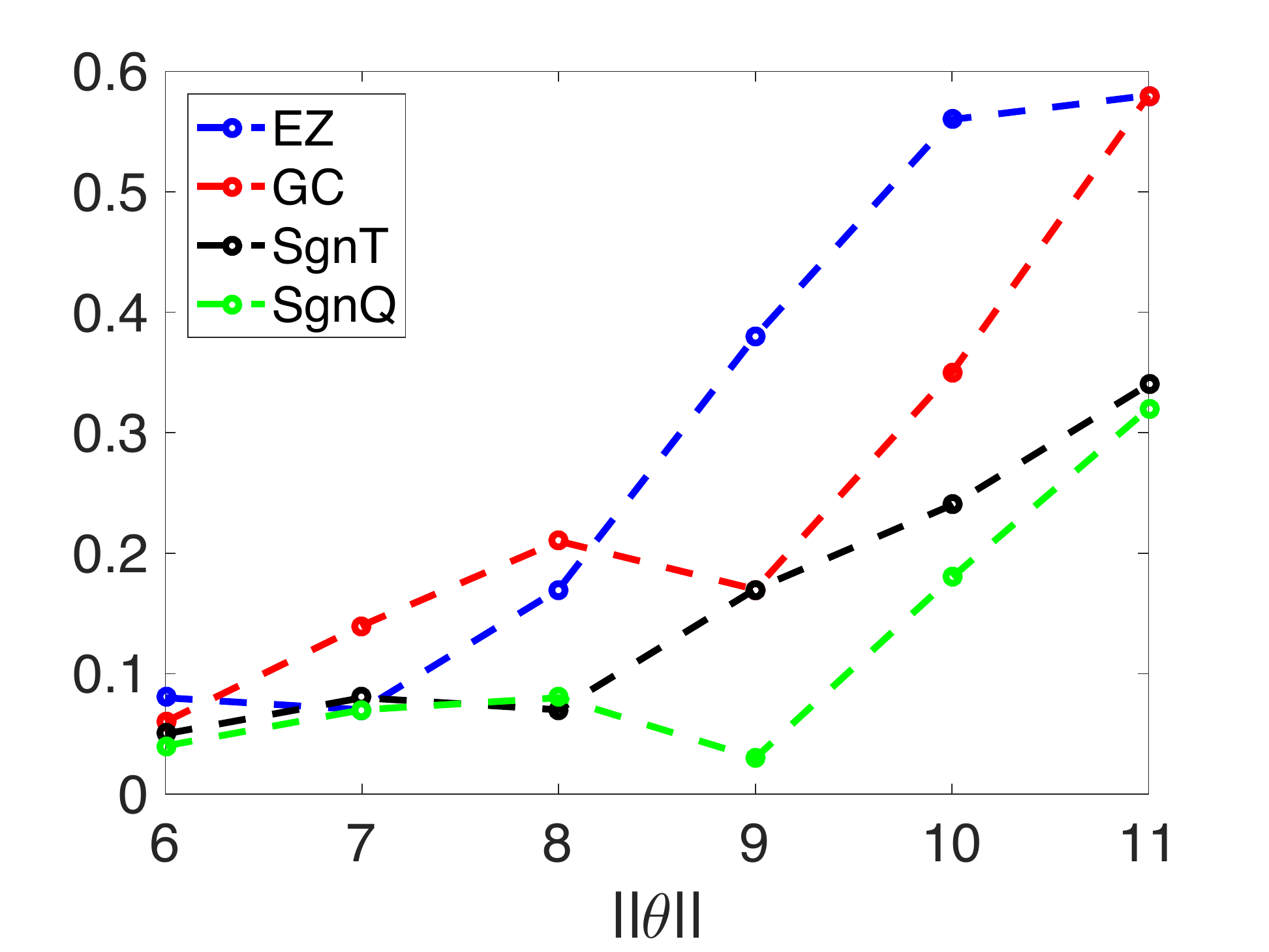}
\includegraphics[width=.42\textwidth, height = .35\textwidth]{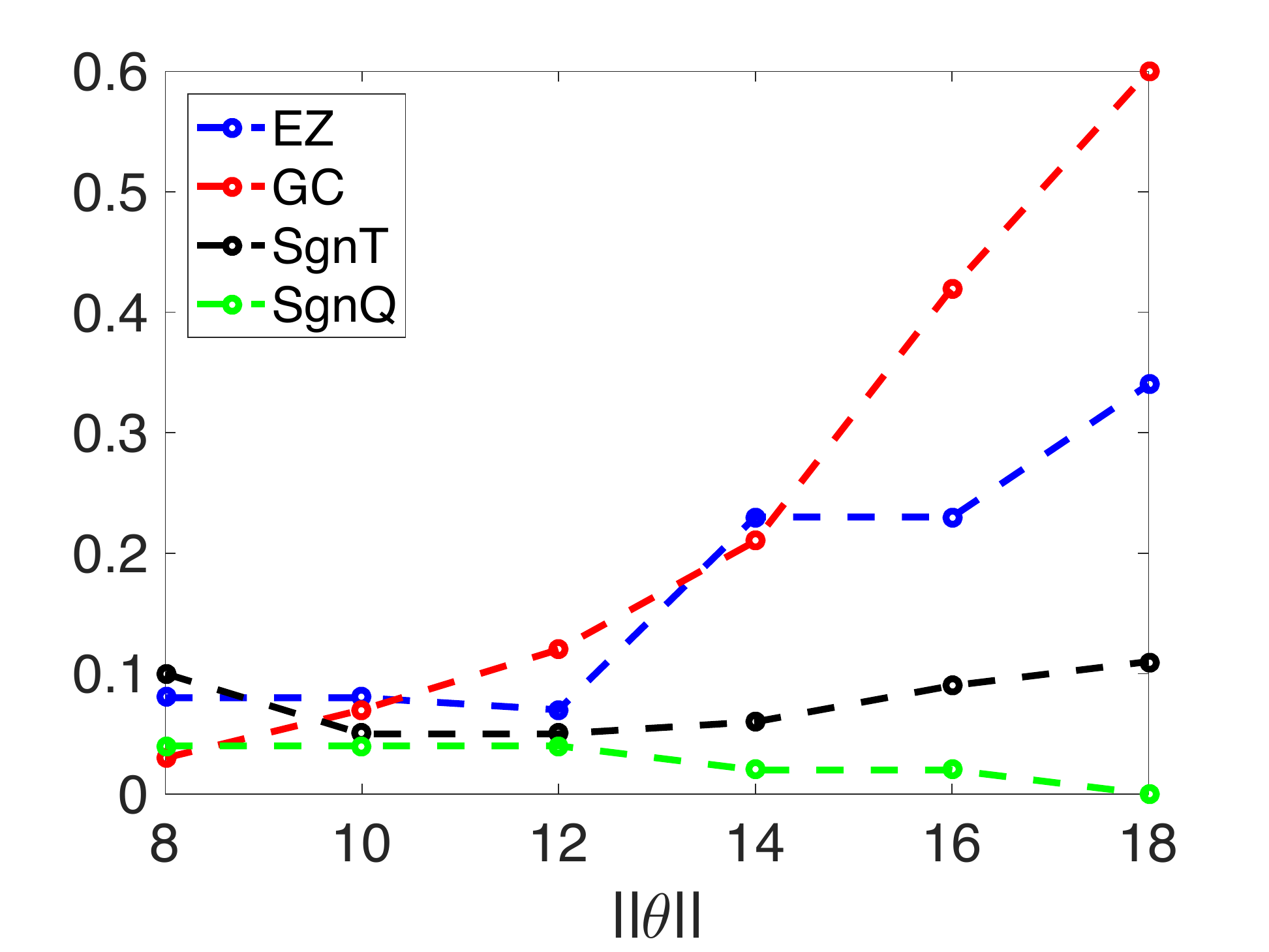}
\caption{From left to right: Experiment 2a and 2b. The $y$-axis are the sum of Type I and Type II errors (testing level is fixed at $5\%$). The $x$-axis are $\|\theta\|$ or sparsity levels.   Results are based on $200$ repetitions.}  \label{fig:Experiment2}
\end{figure}

\paragraph{Experiment 2}  We study the cases with larger $K$ and more complicate matrix $P$. For a $b_n \in (0,1)$, let $\eps_n = \frac{1}{6}\min(1-b_n, b_n)$,  and  let $P$ be the matrix with $1$ on the diagonal but 
the off-diagonal entries are iid drawn from $\mathrm{Unif}(b_n-\eps_n, b_n + \eps_n)$; once the $P$ matrix is drawn, it is fixed throughout different repetitions. 
We consider two sub-experiments, Exp 2a and 2b.   In Exp 2a, we take $(n, K) = (1000, 5)$,  
$f(\theta)$  to be Pareto$(10,0.375)$, and $g(\pi)$ to be the uniform distribution on $\{e_1, \cdots, e_K\}$ (the standard basis vectors of $\mathbb{R}^K$). We let $\beta_n$ range but $\|\theta\| (1 - b_n )$ is fixed at $4.5$, so the SNR will not change drastically. In Exp 2b,  we take $(n,K) = (3000, 10)$, $f(\theta)$ to be $0.95\delta_{1}+0.05\delta_{3}$, and $g(\pi) = 0.1 \sum_{k=1}^2\delta_{e_k} + 0.15 \sum_{k=3}^6\delta_{e_k} + 0.05 \sum_{k=7}^{10}\delta_{e_k}$ (so to have unbalanced community sizes). Similarly, we let $\beta_n$ range but  fix  $\|\theta\| (1 - b_n )  = 5.2$. 
The sum of Type I and Type II errors are displayed in Figure~\ref{fig:Experiment2}.

In these examples,  EZ and GC underperform SgnT and SgnQ, especially in the less sparse case, and the performances of the SgnT and SgnQ are more similar to each other, compared to  those in Experiment 1. In these examples, we have larger $K$, more complicate $P$, and unbalanced community sizes, and the performance of the SgnT and SgnQ test statistics suggest that they are relatively robust.

\paragraph{Experiment 3} We investigate the role of mixed-membership. We have three sub-experiments, Exp 3a-3c.  where the memberships are not-mixed, lightly mixed, and significantly mixed,  respectively. 
For all sub-experiments, we take $(n,K) = (2000, 3)$ and $f(\theta)$ to be $\mathrm{Unif}(2,3)$. For Exp 3a,  we let $g_1(\pi) = 0.4 \delta_{e_1} + 0.3 \delta_{e_2} + 0.3 \delta_{e_3}$. In 
Exp 3b,  we let  $g_2(\pi) = 0.3 \sum_{k=1}^3 \delta_{e_k} + 0.1 \cdot \text{Dirichlet}$,  and in Exp 3c, we let $g_3(\pi) = 0.25 \sum_{k=1}^3 \delta_{e_k} + 0.25 \cdot \text{Dirichlet}$, where $\text{Dirichlet}$ represents the symmetric $K$-dimensional Dirichlet distribution. In Exp 3a-3b,  we 
let $\beta_n$ range while $(1 -b_n)\|\theta\|$ is fixed at  $4.2$ so the SNR's are roughly the same. In Exp 3c,  we also let $\beta_n$ range but $(1 -b_n)\|\theta\| = 4.5$ (the SNR's need to be slightly larger to counter the effect of mixed-membership, which makes the testing problem harder).

\begin{figure}[tb] 
\centering
\vspace{-0mm}
\includegraphics[width=.32\textwidth, height = .32\textwidth]{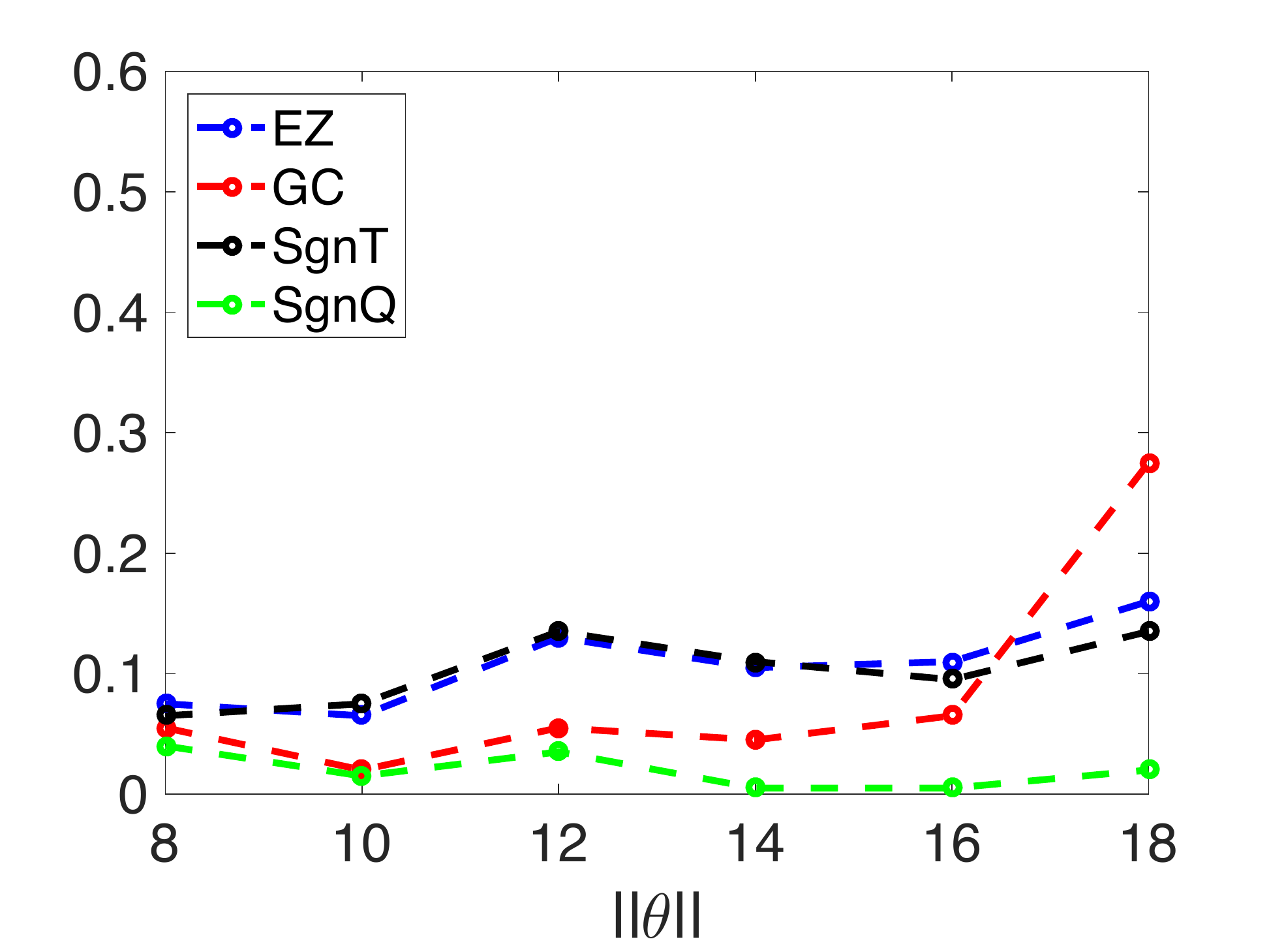}
\includegraphics[width=.32\textwidth, height = .32\textwidth]{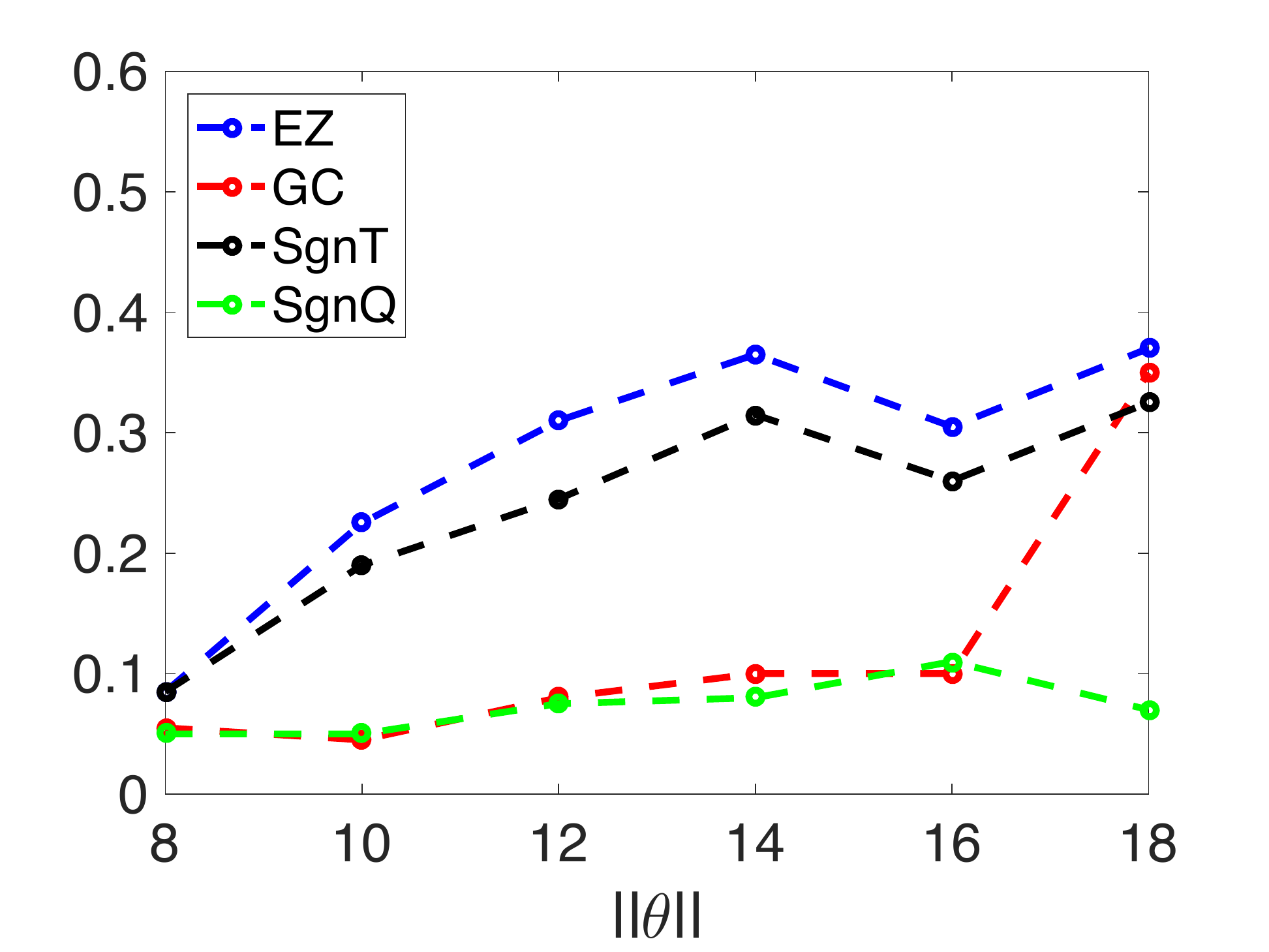}
\includegraphics[width=.32\textwidth, height = .32\textwidth]{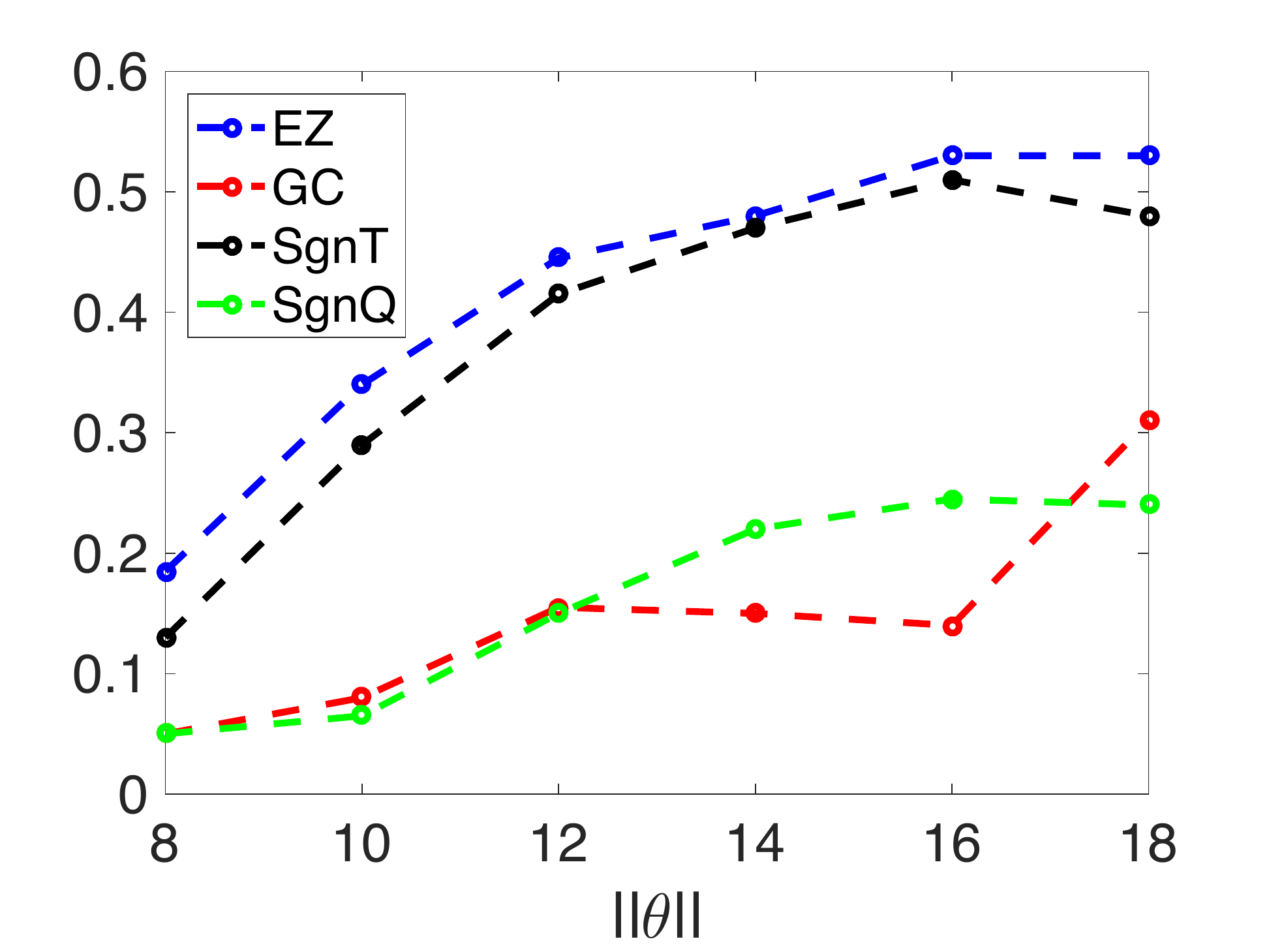}
\caption{From left to right: Experiment 3a,  3b, and 3c. The $y$-axis are the sum of Type I and Type II errors (testing level is fixed at $5\%$). The $x$-axis are $\|\theta\|$ or sparsity levels.   Results are based on $200$ repetitions.}  \label{fig:Experiment3}
\end{figure}

The sum of Type I and Type II errors are presented in Figure~\ref{fig:Experiment3}. 
First, the results confirm that mixed-memberships make the testing problem harder. 
For example, the value of $\|\theta\| (1 - b_n)$ in Exp 3c is higher than that of Exp 3a-3b, 
but the testing errors are higher, due to that the memberships in Exp 3c are more mixed. 
Second,  SgnQ consistently outperforms EZ and SgnT. Third, GC is comparable with SgnQ in the 
more sparse case, but performs unsatisfactorily in the less sparse case, for reasons explained before. Last, in these settings, SgnT is uniformly better than EZ, and more so when the memberships  become more mixed.

\section{Discussions}  
\label{sec:Discu} 
A closely related idea is to use $\|A - \hat{\eta} \hat{\eta}'\|$ as the 
test statistics. To see why this is a reasonable choice, consider the proxy test statistic 
$\|A - \eta^* (\eta^*)'\|$, where we recall $\eta^* = \theta$ under the null; see (\ref{Defineeta*}).  
Therefore,  
\[
A - \eta^* (\eta^*)' = \left\{
\begin{array}{ll} 
W,  &\qquad \mbox{under the null}, \\ 
(\Omega -(\eta^* (\eta^*)')  + W,  &\qquad \mbox{under the alternative},     \\ 
\end{array} 
\right.
\]  
and the test has reasonable power, as $\|A - \eta^* (\eta^*)'\|$ is expected to be bigger in the alternative than in the null. Another related idea is to extend the Signed Polygon to address the problem of 
testing whether $K = k_0$ vs. $K > k_0$, where $k_0 > 1$ is a prescribed integer. Let $\hat{\Omega} = \sum_{k = 1}^{k_0} \hat{\lambda}_k \hat{\xi}_k \hat{\xi}_k'$, where $\hat{\lambda}_k$ are the $k$-th eigenvalue of $A$, arranged in the descending order in magnitude, and $\hat{\xi}_k$ is the corresponding eigenvector. 
The Signed Polygon test statistic can then be extended to 
\[
U_{n,k_0}^{(m)} = \sum_{i_1, i_2, \ldots, i_m (dist)} (A_{i_1 i_2} - \hat{\Omega}_{i_1 i_2}) (A_{i_2 i_3} - \hat{\Omega}_{i_2 i_3}) \ldots (A_{i_m i_1} - \hat{\Omega}_{i_m i_1}). 
\] 
See \cite{JLWtesting2019} for more discussion. It remains unclear whether these test statistics are optimally adaptive, and we leave the study to the future. 

The testing problem is also closely related to the problem of estimating $K$. 
In fact, we can cast the estimation problem as a sequential testing problem where we test 
$K = k_0$ vs. $K > k_0$ for $k_0 = 1, 2, 3, \ldots$, and estimate $K$ to be the smallest $k_0$ where 
we accept the null.  

Note also the lower bound argument for the global testing problem sheds useful insight  for many other problems (e.g., estimating $K$, community detection, mixed-membership). Take the problem of estimating $K$ for example. Given an alternative setting, if we can not distinguish it from some null setting, then the underlying parameter $K$ is not estimable.

In a high level, these ideas, together with the Signed Polygon, are related to the ideas in \cite{Hu2017TestEntry} on testing $K = k_0$ vs. $K > k_0$,  in \cite{lei2016goodness} on goodness of fit, and in \cite{le2015estimating} on estimating $K$. However, the focus of these works are on the more idealized model where  we don't have degree heterogeneity, and how to extend their ideas to the current setting remains unclear.

   \newpage

\begin{center}
{\LARGE Appendix for ``Optimally Adaptivity of Signed Polygon Statistics for Network Testing"}
\end{center}

\bigskip\bigskip

\noindent
This appendix contains technical proofs of the main article \cite{utest}. Appendix~\ref{supp:MatrixForm} derives matrix forms of signed-polygon statistics and proves Theorem~\ref{thm:compu}. Appendix~\ref{supp:eta-norm} studies the estimation error of $\|\theta\|^2$ and proves Lemma~\ref{lemma:heta}. Appendix~\ref{supp:Spectral} contains spectral analysis for $\Omega$ and $\widetilde{\Omega}$ and proves Lemmas~\ref{lem:tracetOmega1}-\ref{lem:tracetOmega2}. Appendix~\ref{supp:LB} analyzes the region of impossibility and proves Lemma~\ref{lem:existence} and Theorems~\ref{thm:minimax}-\ref{thm:RI-sameTheta}. Appendix~\ref{supp:Var} calculates the mean and variance of signed-polygon statistics and proves the results in Tables~\ref{tab:IdealSgnTsum}-\ref{tab:ProxySgnQsum} and Theorems~\ref{thm:IdealSgnT}-\ref{thm:RealSgnQ}.

 \appendix

\section{Matrix forms of Signed-Polygon statistics} \label{supp:MatrixForm}
We prove Theorem~\ref{thm:compu}. Recall that $\widetilde{A} = A - \hat{\eta}\hat{\eta}$. By definition, 
\begin{align*}
T_n & = \tr(\widetilde{A}^3) - \sum_{\substack{\text{at least two of}\\ i,j,k \text{ are equal}}}\widetilde{A}_{ij}\widetilde{A}_{jk}\widetilde{A}_{ki},\cr
Q_n & = \tr(\widetilde{A}^4) - \sum_{\substack{\text{at least two of}\\ i,j,k,\ell \text{ are equal}}}\widetilde{A}_{ij}\widetilde{A}_{jk}\widetilde{A}_{k\ell}\widetilde{A}_{\ell i}. 
\end{align*}

First, we derive the matrix form of $T_n$. If at least two of $\{i,j,k\}$ are equal, there are four cases: (a) $i=j$, $k\neq i$, (b) $j=k$, $i\neq j$, (c) $k=i$, $j\neq k$, (d) $i=j=k$. The first three cases are similar. It follows that
\begin{align*}
T_n &= \tr(\widetilde{A}^3) - 3\sum_{i, k (dist)} \widetilde{A}_{ii}\widetilde{A}_{ik}^2 - \sum_{i} \widetilde{A}_{ii}^3\cr
&= \tr(\widetilde{A}^3) - 3\Bigl(\sum_{i, k} \widetilde{A}_{ii}\widetilde{A}_{ik}^2 - \sum_{i}\widetilde{A}_{ii}^3\Bigr) - \sum_{i} \widetilde{A}_{ii}^3\cr
&=\tr(\widetilde{A}^3) - 3\tr(\widetilde{A}\circ \widetilde{A}^2) + 2\tr(\widetilde{A}\circ\widetilde{A}\circ\widetilde{A}). 
\end{align*}
This gives the desired expression of $T_n$. 

\begin{table}[!bt]
\centering
\caption{Decomposition of $\tr(\widetilde{A}^4)$. We note that the last column sums to $n^4$.}\label{tb:Qterms}
\scalebox{1}{\begin{tabular}{ccccccc}
\hline
Pattern & Variations & Summand & Sum & \#Summands \\
\hline
$\{i,j,k,\ell\}$ & $(i,j,k,\ell)$ & $\widetilde{A}_{ij}\widetilde{A}_{jk}\widetilde{A}_{k\ell}\widetilde{A}_{\ell i}$ & $Q_n$ & $n(n-1)(n-2)(n-3)$\\
\hline
\multirow{6}{*}{$\{i,i,j,k\}$} & $(i,j,k,i)$ & $\widetilde{A}_{ij}\widetilde{A}_{jk}\widetilde{A}_{ki}\widetilde{A}_{ii}$ & $S_1$ & \multirow{6}{*}{$6n(n-1)(n-2)$}\\
& $(i,j,k,j)$ & $\widetilde{A}_{ij}\widetilde{A}_{jk}\widetilde{A}_{kj}\widetilde{A}_{ji}$ & $S_2$\\
& $(i,j,k,k)$ & $\widetilde{A}_{ij}\widetilde{A}_{jk}\widetilde{A}_{kk}\widetilde{A}_{ki}$ & $S_1$\\
& $(i,j,i,k)$ & $\widetilde{A}_{ij}\widetilde{A}_{ji}\widetilde{A}_{ik}\widetilde{A}_{ki}$ & $S_2$ \\
& $(i,j,j,k)$ & $\widetilde{A}_{ij}\widetilde{A}_{jj}\widetilde{A}_{jk}\widetilde{A}_{ki}$ & $S_1$\\
& $(i,i,j,k)$ & $\widetilde{A}_{ii}\widetilde{A}_{ij}\widetilde{A}_{jk}\widetilde{A}_{ki}$ & $S_1$\\
\hline
\multirow{4}{*}{$\{i,i,i,j\}$} & $(i,j,i,i)$& $\widetilde{A}_{ij}\widetilde{A}_{ji}\widetilde{A}_{ii}\widetilde{A}_{ii}$ & $S_3$ & \multirow{4}{*}{$4n(n-1)$}\\
& $(i,j,j,j)$ & $\widetilde{A}_{ij}\widetilde{A}_{jj}\widetilde{A}_{jj}\widetilde{A}_{ji}$ & $S_3$\\
& $(i,i,j,i)$ & $\widetilde{A}_{ii}\widetilde{A}_{ij}\widetilde{A}_{ji}\widetilde{A}_{ii}$ & $S_3$\\
& $(i,i,i,j)$ & $\widetilde{A}_{ii}\widetilde{A}_{ii}\widetilde{A}_{ij}\widetilde{A}_{ji}$ & $S_3$ \\
\hline
\multirow{3}{*}{$\{i,i,j,j\}$}  & $(i,j,i,j)$ & $\widetilde{A}_{ij}\widetilde{A}_{ji}\widetilde{A}_{ij}\widetilde{A}_{ji}$ & $S_4$ &\multirow{3}{*}{$3n(n-1)$}\\
& $(i,j,j,i)$ & $\widetilde{A}_{ij}\widetilde{A}_{jj}\widetilde{A}_{ji}\widetilde{A}_{ii}$  &  $S_5$ \\
& $(i,i,j,j)$ & $\widetilde{A}_{ii}\widetilde{A}_{ij}\widetilde{A}_{jj}\widetilde{A}_{ji}$ & $S_5$ \\
\hline
$\{i,i,i,i\}$ & $(i,i,i,i)$ & $\widetilde{A}_{ii}\widetilde{A}_{ii}\widetilde{A}_{ii}\widetilde{A}_{ii}$ & $S_6$&$n$\\
\hline
\end{tabular}}
\end{table} 

Next, we derive the matrix form of $Q_n$. 
When at least two of $\{i,j,k,\ell\}$ are equal, depending on how many indices are equal, we have four patterns: $\{i,i,i,i\}$, $\{i,i,i,j\}$, $\{i,i,j,j\}$, $\{i,i,j,k\}$, where $(i,j,k)$ are distinct. For each pattern, depending on the appearing locations of the next distinct indices,  there are a few variations. Take the pattern $\{i,i,j,k\}$ for example: (a) when a new distinct index appears at location 2 and at location 3, the variations are $(i,j,k,i)$, $(i,j,k,j)$, $(i,j,k,k)$; (b) when a new distinct index appears at location 2 and at location 4, the variations are $(i,j,i,k)$, $(i,j,j,k)$; (c) when a new distinct index appears at location 3 and location 4, the variation is $(i,i,j,k)$. Using similar arguments, we can find all variations of each pattern. They are summarized in Table~\ref{tb:Qterms}.
Define
\begin{align*}
& S_1= \sum_{i,j,k (dist)}\widetilde{A}_{ii}\widetilde{A}_{ij}\widetilde{A}_{jk}\widetilde{A}_{ki}, &&
S_2 = \sum_{i, j, k (dist)}\widetilde{A}^2_{ij}\widetilde{A}^2_{ik},\cr
&S_3 = \sum_{i, j (dist)} \widetilde{A}^2_{ii}\widetilde{A}^2_{ij}, &&S_4 = \sum_{i, j (dist)} \widetilde{A}_{ij}^4,  \cr
& S_5 = \sum_{i, j (dist)}\widetilde{A}_{ii}\widetilde{A}^2_{ij}\widetilde{A}_{jj},&& S_6 = \sum_i \widetilde{A}_{ii}^4. 
\end{align*}
It follows from Table~\ref{tb:Qterms} that 
\beq \label{Q-compu}
Q_n = \tr(\widetilde{A}^4) - 4S_1-2S_2-4S_3 - S_4 - 2S_5-S_6. 
\eeq
What remains is to derive the matrix form of $S_1$-$S_6$. By direct calculations,
\begin{align*}
S_1 &= \sum_i \widetilde{A}_{ii}\biggl[ \sum_{j\neq i,k\neq i}\widetilde{A}_{ij}\widetilde{A}_{jk}\widetilde{A}_{ki} -\sum_{j\neq i} \widetilde{A}_{ij}\widetilde{A}_{jj}\widetilde{A}_{ji}\biggr]\cr
&= \sum_i \widetilde{A}_{ii}\biggl[ \Bigl( \sum_{j,k}\widetilde{A}_{ij}\widetilde{A}_{jk}\widetilde{A}_{ki}- 2\sum_{j}\widetilde{A}^2_{ij}\widetilde{A}_{ii} + \widetilde{A}_{ii}^3 \Bigr)-\Bigl(\sum_{j} \widetilde{A}^2_{ij}\widetilde{A}_{jj} - \widetilde{A}^3_{ii}\Bigr)\biggr]\cr
&=\sum_{i,j,k}\widetilde{A}_{ii}\widetilde{A}_{ij}\widetilde{A}_{jk}\widetilde{A}_{ki} -2\sum_{i,j}\widetilde{A}_{ii}^2\widetilde{A}_{ij}^2 -\sum_{i,j}\widetilde{A}_{ii}\widetilde{A}_{ij}^2\widetilde{A}_{jj} + 2\sum_{i}\widetilde{A}^{4}_{ii}\cr
&=\tr(\widetilde{A}\circ\widetilde{A}^3) - 2\tr(\widetilde{A}\circ\widetilde{A}\circ\widetilde{A}^2) - 1_n'[\diag(\widetilde{A})(\widetilde{A} \circ \widetilde{A})\diag(\widetilde{A})]1_n+2S_6. 
\end{align*}
Moreover, we can derive that
\begin{align*}
S_2 &= \sum_i \biggl[ \sum_{j\neq i, k\neq i}\widetilde{A}^2_{ij}\widetilde{A}^2_{ik} -  \sum_{j\neq i}\widetilde{A}^4_{ij}\biggr]\cr
&= \sum_i \biggl[ \Bigl( \sum_{j,k}\widetilde{A}^2_{ij}\widetilde{A}^2_{ik} -2\sum_{j}\widetilde{A}_{ij}^2\widetilde{A}_{ii}^2 + \widetilde{A}^4_{ii}\Bigr)-\Bigl(\sum_j\widetilde{A}^4_{ij}-\widetilde{A}^4_{ii}\Bigr) \biggr]\cr
&=\sum_{i,j,k}\widetilde{A}^2_{ij}\widetilde{A}^2_{ik} -2\sum_{i,j}\widetilde{A}_{ij}^2\widetilde{A}^2_{ii}-\sum_{i,j}\widetilde{A}^4_{ij} + 2\sum_i\widetilde{A}^4_{ii}\cr
&= \tr(\widetilde{A}^2\circ\widetilde{A}^2) - 2\tr(\widetilde{A}\circ\widetilde{A}\circ\widetilde{A}^2)-1_n'[\widetilde{A}\circ\widetilde{A}\circ\widetilde{A}\circ\widetilde{A}]1_n + 2S_6. 
\end{align*}
It is also easy to see that
\begin{align*}
S_3 &= \sum_{i,j} \widetilde{A}^2_{ii}\widetilde{A}^2_{ij} - \sum_{i}\widetilde{A}^{4}_{ii}  =\tr(\widetilde{A}\circ\widetilde{A} \circ \widetilde{A}^2)- S_6,\cr
S_4 &= \sum_{i,j}\widetilde{A}_{ij}^4 - \sum_{i}\widetilde{A}^{4}_{ii} = 1_n'[ \widetilde{A}\circ \widetilde{A}\circ\widetilde{A}\circ\widetilde{A} ]1_n  - S_6,\cr
S_5 &= \sum_{i,j}\widetilde{A}_{ii}\widetilde{A}^2_{ij}\widetilde{A}_{jj} -S_6  = 1_n'[\diag(\widetilde{A})(\widetilde{A} \circ \widetilde{A})\diag(\widetilde{A})]1_n - S_6,\cr
S_6 &= \tr( \widetilde{A}\circ \widetilde{A}\circ\widetilde{A}\circ\widetilde{A}). 
\end{align*}
Plugging the matrix forms of $S_1$-$S_6$ into \eqref{Q-compu}, we obtain
\begin{align*}
Q_n =& \tr(\widetilde{A}^4) - 4\tr(\widetilde{A}\circ\widetilde{A}^3) - 2\tr(\widetilde{A}^2\circ\widetilde{A}^2) + 8\tr(\widetilde{A}\circ\widetilde{A}\circ\widetilde{A}^2)- 6\tr(\widetilde{A}\circ\widetilde{A}\circ\widetilde{A}\circ\widetilde{A})\cr
& + 2\cdot 1_n'[\diag(\widetilde{A})(\widetilde{A} \circ \widetilde{A})\diag(\widetilde{A})]1_n + 1_n'[ \widetilde{A}\circ\widetilde{A}\circ\widetilde{A}\circ\widetilde{A}]1_n.  
\end{align*}
This gives the desired expression of $Q_n$.

Last, we discuss the complexity of computing $T_n$ and $Q_n$. It involves the following operations:
\begin{itemize}
\item Compute the matrix $\widetilde{A}=A-\hat{\eta}\hat{\eta}'$.
\item Compute the Hadamard product of finitely many matrices.
\item Compute the trace of a matrix.
\item Compute the matrix $DMD$ for a matrix $M$ and a diagonal matrix $D$.
\item Compute $1_n'M1_n$ for a matrix $M$.
\item Compute the matrices $\widetilde{A}^k$, for $k=2,3,4$. 
\end{itemize}
Excluding the last operation, the complexity is $O(n^2)$. For the last operation, since we can compute $\widetilde{A}^k$ recursively from $\widetilde{A}^k=\widetilde{A}^{k-1}\widetilde{A}$, it suffices to consider the complexity of computing $B\widetilde{A}$, for an arbitrary $n\times n$ matrix $B$. Write
\[
B\widetilde{A}=BA - B\heta (\heta)'. 
\]
Consider computing $BA$. 
The $(i,j)$-th entry of $BA$ is $\sum_{\ell: A_{\ell j}\neq 0}B_{i\ell}A_{\ell j}$, where the total number of nonzero $A_{\ell j}$ equals to $d_j$, the degree of node $j$. Hence, the complexity of computing the $(i,j)$-th entry of $BA$ is $O(d_j)$. It follows that the complexity of computing $BA$ is $O(\sum_{i,j=1}^n d_j)=O(n^2\bar{d})$. Consider computing $B\heta (\heta)'$. We first compute the vector $v=B\heta$ and then compute $v(\heta)'$, where the complexity of both steps is $O(n^2)$. Combining the above, the complexity of computing $B\widetilde{A}$ is $O(n^2\bar{d})$. We have seen that this is the dominating step in computing $T_n$ and $Q_n$, so the complexity of the latter is also $O(n^2\bar{d})$.

\newpage
\section{Estimation of $\|\theta\|$}  \label{supp:eta-norm}
We prove Lemma~\ref{lemma:heta}. 
First, we show that 
\[
\|\eta^*\|^2   \left\{ 
\begin{array}{ll} 
= \|\theta\|^2, &\qquad \mbox{under the null}, \\
\asymp \|\theta\|^2, &\qquad \mbox{under the alternative}.  \\
\end{array} 
\right. 
\] 
Recall that $\eta^*=(1/\sqrt{1_n'\Omega 1_n})\Omega 1_n$. Hence, 
\beq \label{lem-heta-1}
\|\eta^*\|^2 = (1_n'\Omega^2 1_n)/(1_n'\Omega 1_n). 
\eeq
Under the null, $\Omega=\theta\theta'$, and the claim follows by direct calculations. Under the alternative, $\Omega=\sum_{k=1}^K\lambda_k\xi_k\xi_k'$, so
\[
1_n'\Omega 1_n = \sum_{k=1}^K\lambda_k (1_n'\xi_k)^2, \qquad 1_n'\Omega^2 1_n = \sum_{k=1}^K\lambda_k^2 (1_n'\xi_k)^2. 
\] 
By Lemma~\ref{lem:Omega}, $\lambda_1\asymp \|\theta\|^2$. 
By Lemma~\ref{lem:vectorh}, $1_n'\xi_1\asymp \|\theta\|^{-1}\|\theta\|_1$ and $|1_n'\xi_k| = O(\|\theta\|^{-1}\|\theta\|_1)$. It follows that $1_n'\Omega^2 1_n\geq \lambda_1^2(1_n'\xi_1)^2\geq C\|\theta\|_1^2\|\theta\|^2$ and $1_n'\Omega^21_n\leq \lambda_1^2\sum_{k=1}^K(1_n'\xi_k)^2\leq C\|\theta\|_1^2\|\theta\|^2$. We conclude that
\beq \label{lem-heta-2}
1_n'\Omega^2 1_n \asymp \|\theta\|_1^2\|\theta\|^2. 
\eeq
Moreover, $1_n'\Omega 1_n\leq |\lambda_1|\sum_{k=1}^K (1_n'\xi_k)^2\leq C\|\theta\|_1^2$, and by Lemma~\ref{lem:1Omega1}, $1_n'\Omega 1_n\geq C\|\theta\|_1^2$. It follows that
\beq \label{lem-heta-3}
1_n'\Omega 1_n \asymp \|\theta\|_1^2. 
\eeq
Plugging \eqref{lem-heta-2}-\eqref{lem-heta-3} into \eqref{lem-heta-1} gives the claim.

Next, we show $(\|\hat{\eta}\|^2 - 1)/\|\eta^*\|^2 \goto 1$ in probability. Since $\|\eta^*\|\asymp\|\theta\| \goto \infty$ as $n \goto \infty$, it suffices to show $\|\hat{\eta}\|^2 / \|\eta^*\|^2 \goto 1$ in probability.  
By definition, 
\[
\|\hat{\eta}\|^2 = \frac{1_n' A^2 1_n}{1_n' A 1_n}. 
\] 
Compare this with (\ref{lem-heta-1}), all we need to show is that in probability, 
\begin{equation} \label{mypf0} 
\frac{1_n' A 1_n}{1_n' \Omega 1_n} \goto 1, \qquad \mbox{and} \qquad \frac{1_n' A^2 1_n}{1_n' \Omega^2 1_n} \goto 1. 
\end{equation} 
Since the proofs are similar, we only show the second one. 
By elementary probability, it is sufficient to show that as $n \goto \infty$,  
\begin{equation} \label{mypf1}
\frac{\mathbb{E}[1_n' A^2 1_n]}{1_n' \Omega^2 1_n} \goto 1,  \qquad 
\frac{\mathrm{Var}(1_n' A^2 1_n)}{(1_n' \Omega^2 1_n)^2} \goto 0. 
\end{equation} 

We now prove (\ref{mypf1}). Consider the first claim. Write 
\beq  \label{lem-heta-7}
1_n'A^21_n = \sum_{i,j,k}A_{ij}A_{jk} = \sum_{i\neq j}A^2_{ij} + \sum_{i,j,k (dist)} A_{ij}A_{jk}.  
\eeq
It follows that
\[
\mathbb{E}[1_n'A^21_n] =\sum_{i\neq j}\Omega_{ij}+\sum_{i,j,k (dist)}\Omega_{ij}\Omega_{jk}. 
\]
Since $\Omega_{ij}\leq \theta_i\theta_j$ under both hypotheses, we have
\begin{align*}
\bigl| \mathbb{E}[1_n'A^21_n]-1_n'\Omega 1_n-1_n'\Omega^2 1_n \bigr|&\leq 
\Bigl|\sum_{i}\Omega_{ii}+ \sum_{\substack{(i,j,k)\text{ are}\\\text{not distinct}}} \Omega_{ij}\Omega_{jk}\Bigr|\cr
& \leq \sum_i\theta_i^2+C\sum_{i,j}\theta_i^2\theta_j^2 + C\sum_{i,k}\theta_i^3\theta_k\cr
&\leq C\|\theta\|^2 +C\|\theta\|^4 +C\|\theta\|_3^3\|\theta\|_1\cr
&\leq C\|\theta\|_3^3\|\theta\|_1,
\end{align*}
where we have used the universal inequality $\|\theta\|^4\leq \|\theta\|_3^3\|\theta\|_1$. Since $\|\theta\|_3^3\leq \theta_{\max}^2\|\theta\|_1=o(\|\theta\|_1)$, the right hand side is $o(\|\theta\|_1^2)=o(1_n'\Omega 1_n)$. So, 
\beq   \label{lem-heta-8}
\mathbb{E}[1_n'A^21_n] =1_n'\Omega^21_n +1_n'\Omega 1_n + o(1_n'\Omega 1_n). 
\eeq
Combining this with (\ref{lem-heta-2})-(\ref{lem-heta-3}) gives 
\[
\Bigl|\frac{\mathbb{E}[1_n' A^2 1_n]}{1_n' \Omega^2 1_n} - 1\Bigr| \lesssim  \frac{1_n' \Omega 1_n}{1_n' \Omega^2 1_n} \asymp \frac{1}{\|\theta\|^2},  
\] 
and the claim follows by $\|\theta\| \goto \infty$.

Consider the second claim. By \eqref{lem-heta-7}, 
\beq \label{lem-heta-add1}
\mathrm{Var}(1_n'A^2 1_n)\leq 2\mathrm{Var}\Bigl(\sum_{i\neq j}A^2_{ij}\Bigr)+ 2\mathrm{Var}\Bigl( \sum_{i,j,k (dist)}A_{ij}A_{jk}\Bigr). 
\eeq
We re-write $\sum_{i\neq j}A^2_{ij}=\sum_{i\neq j}A_{ij}=2\sum_{i<j}A_{ij}$. The variables $\{A_{ij}\}_{1\leq i<j\leq n}$ are mutually independent. It follows that
\beq \label{lem-heta-add2}
\mathrm{Var}\Bigl( \sum_{i\neq j}A^2_{ij} \Bigr)=4\sum_{i<j}\mathrm{Var}(A_{ij})\leq C\sum_{i,j}\Omega_{ij}\leq C\|\theta\|_1^2. 
\eeq
Moreover, since $A_{ij}A_{jk}=(\Omega_{ij}+W_{ij})(\Omega_{jk}+W_{jk})$, we have
\begin{align*}
\sum_{i,j,k (dist)}A_{ij}A_{jk} &= \sum_{i,j,k (dist)} \Omega_{ij}\Omega_{jk} + 2\sum_{i,j,k (dist)}\Omega_{ij}W_{jk} + \sum_{i,j,k (dist)}W_{ij}W_{jk}\cr
&\equiv   \sum_{i,j,k (dist)} \Omega_{ij}\Omega_{jk}  + X_1 + X_2. 
\end{align*}
By elementary probability,
\[
\mathrm{Var}\Bigl( \sum_{i,j,k (dist)}A_{ij}A_{jk} \Bigr)\leq 2\mathrm{Var}(X_1)+2\mathrm{Var}(X_2). 
\]
To compute the variance of $X_1$, we note that
\[
X_1 = 4\sum_{j<k}\beta_{jk}W_{jk}, \qquad \beta_{jk}=\sum_{i\notin\{j,k\}}\Omega_{ij}. 
\]
The variables $\{W_{jk}\}_{1\leq j<k\neq n}$ are mutually independent, and $|\beta_{jk}|\leq C\sum_{i}\theta_i\theta_j\leq C\|\theta\|_1\theta_j$. It follows that 
\[
\mathrm{Var}(X_1)\leq C\sum_{j,k}(\|\theta\|_1\theta_j)^2(\theta_j\theta_k)\leq C\|\theta\|_1^3\|\theta\|_3^3. 
\]
To compute the variance of $X_2$, we note that 
\[
\mathrm{Var}(X_2) = \sum_{i,j,k (dist)}\sum_{i',j',k' (dist)}\mathbb{E}[W_{ij}W_{jk}W_{i'j'}W_{j'k'}]. 
\] 
The summand is nonzero only when the two variables $\{W_{i'j'}, W_{j'k'}\}$ are the same as the two variables $\{W_{ij}, W_{jk}\}$. This can only happen if $(i,j,k)=(i',j',k')$ or $(i,j,k)=(k',j',i')$, where in either case the summand equals to $\mathbb{E}[W^2_{ij}W^2_{jk}]$. It follows that
\[
\mathrm{Var}(X_2)=\sum_{i,j,k (dist)} 2\mathbb{E}[W^2_{ij}W^2_{jk}] \leq C\sum_{i,j,k}\theta_i\theta_j^2\theta_k \leq C\|\theta\|^2\|\theta\|_1^2. 
\]
Combining the above gives
\beq \label{lem-heta-add3}
\mathrm{Var}\Bigl( \sum_{i,j,k (dist)}A_{ij}A_{jk} \Bigr)\leq C\|\theta\|_1^3\|\theta\|_3^3 + C\|\theta\|^2\|\theta\|_1^2\leq C\|\theta\|_1^3\|\theta\|_3^3,
\eeq
where we have used the fact that $\|\theta\|_1\|\theta\|_3^3\geq \|\theta\|^4$ (Cauchy-Schwarz inequality) and $\|\theta\|\to\infty$. Plugging \eqref{lem-heta-add2}-\eqref{lem-heta-add3} into \eqref{lem-heta-add1} gives
\beq   \label{lem-heta-9}
\mathrm{Var}(1_n'A^21_n)\leq C\|\theta\|_1^3\|\theta\|_3^3. 
\eeq 
Comparing this with (\ref{lem-heta-2}) and using $\|\theta\|_3^3\leq\theta^2_{\max}\|\theta\|_1$, we obtain
\[
\frac{\mathrm{Var}(1_n'A^21_n)}{(1_n'\Omega^2 1_n)^2}\leq \frac{C\|\theta\|_1^3\|\theta\|_3^3}{\|\theta\|_1^4\|\theta\|^4}\leq \frac{C\theta_{\max}^2}{\|\theta\|^4},
\]
and the claim follows by $\|\theta\|\to\infty$.

\newpage
\section{Spectral analysis for $\Omega$ and $\widetilde{\Omega}$}   \label{supp:Spectral}
We state and prove some useful results about eigenvalues and eigenvectors of $\Omega$ and $\widetilde{\Omega}$. In Section \ref{subsec:twolemmas}, we prove Lemma  \ref{lem:tracetOmega1} and 
\ref{lem:tracetOmega2} of the main file.

For $1\leq k\leq K$, let $\lambda_k$ be the $k$-th largest (in absolute value) eigenvalue of $\Omega$ and let $\xi_k\in\mathbb{R}^n$ be the corresponding unit-norm eigenvector. We write 
\[
\Xi = [\xi_1, \xi_2, \ldots, \xi_K] = [u_1, u_2, \ldots, u_n]',  
\] 
so that $u_i$ is the $i$-th row of $\Xi$. Recall that $G$ is the $K\times K$ matrix $\|\theta\|^{-2}(\Pi'\Theta^2\Pi)$.

\subsection{Spectral analysis of $\Omega$}
The following lemma relates $\lambda_k$ and $\xi_k$ to the eigenvalues and eigenvectors of the  $K\times K$ matrix $G^{\frac{1}{2}}PG^{\frac{1}{2}}$.  
\begin{lemma} \label{lem:connection}
Consider the DCMM model. Let $d_k$ be the $k$-th largest (in absolute value) eigenvalue of $G^{\frac{1}{2}}PG^{\frac{1}{2}}$ and let $\beta_k\in\mathbb{R}^K$ be the associated eigenvector, $1\leq k\leq K$. Then
under the null, 
\[
\lambda_1 = \|\theta\|^2, \qquad \xi_1 =  \pm \theta/\|\theta\|. 
\]
Under the alternative, for $1 \leq k \leq K$,   
\[
\lambda_k = d_k\|\theta\|^2, \qquad \xi_k = \|\theta\|^{-1}[ \theta\circ (\Pi G^{-\frac{1}{2}}\beta_k)]. 
\]
\end{lemma}

Under the alternative hypothesis, we further have the following lemma:
\begin{lemma} \label{lem:Omega}
Under the DCMM model, as $n\to\infty$, suppose \eqref{cond-balance} holds. As  $n\to\infty$, under the alternative hypothesis,  
\[
\lambda_1\asymp \|\theta\|^2, \qquad   \|u_i\| \leq C\|\theta\|^{-1}\theta_i,\quad \mbox{for all }1\leq i\leq n. 
\]
\end{lemma}

The quantities $(1_n'\xi_k)$ play key roles in the analysis of the Signed Polygon tests. 
By Lemma \ref{lem:connection},  
\[
\xi_1 = (\|\theta\|)^{-1} \Theta \Pi G^{-1/2} \beta_1, 
\] 
where $\beta_1$ is the first eigenvector of $G^{1/2} P G^{1/2}$, corresponding to the 
largest eigenvalue of $G^{1/2} P G^{1/2}$. It is seen $G^{-1/2} \beta_1$ is the eigenvector of the 
matrix $P G$ associated with the largest eigenvalue of $GP$, which is the same as the largest eigenvalue of $G^{1/2} P G^{1/2}$.  Since $PG$ is a non-negative matrix, by Perron's theorem, we can assume all entries of $G^{-1/2} \beta_1$ are non-negative. As a result, all entries of $\xi_1$ are non-negative, and 
\[
1_n' \xi_1 > 0. 
\] 
The following lemma is proved in Section~\ref{supp:SpectralProof}.  
\begin{lemma} \label{lem:vectorh} 
Under the DCMM model, as $n\to\infty$, suppose \eqref{cond-balance} holds. As $n\to\infty$,  
\[
\max_{1\leq k\leq K}|1_n'\xi_k|\leq C\|\theta\|^{-1}\|\theta\|_1, \qquad 1_n'\xi_1 \geq C\|\theta\|^{-1}\|\theta\|_1. 
\]
and so for any $2 \leq k \leq K$,  
\[
|1_n' \xi_k| \leq C |1_n' \xi_1|
\]  
\end{lemma}

We also have a lower bound for $1_n'\Omega 1_n$. The following lemma is proved in Section~\ref{supp:SpectralProof}. 
\begin{lemma} \label{lem:1Omega1}
Under the DCMM model, as $n\to\infty$, suppose \eqref{cond-balance} holds. As $n\to\infty$, both under the null hypothesis and the alternative hypothesis,  
\[
1_n'\Omega 1_n\geq C \|\theta\|_1^2. 
\]
\end{lemma}

\subsection{Spectral analysis of $\widetilde{\Omega}$}
Recall that 
\[
\widetilde{\Omega}=\Omega - (\eta^*)(\eta^*)', \qquad \mbox{where }\eta^*=(1/\sqrt{1_n'\Omega 1_n})\Omega 1_n, 
\]
and $\lambda_1,\ldots, \lambda_K$ are the $K$ nonzero eigenvalues of $\Omega$, arranged in the descending order in magnitude, and $\xi_1,\ldots,\xi_K$ are the corresponding unit-norm eigenvectors of $\Omega$
The following lemma is proved in Section~\ref{supp:SpectralProof}. 
\begin{lemma} \label{lem:tildeOmega}
Under the DCMM model, as $n\to\infty$, suppose \eqref{cond-balance} holds. Then,  
\[
|\lambda_2|  \leq \|\widetilde{\Omega}\|\leq C|\lambda_2|. 
\]
Moreover, for any fixed integer $m\geq 1$, 
\[
|(\widetilde{\Omega}^m)_{ij}|\leq C|\lambda_2|^m\cdot \|\theta\|^{-2}\theta_i\theta_j, \;\;\;\mbox{for all $1\leq i,j\leq n$}. 
\] 
\end{lemma}

Recall that $d_1,\ldots,d_K$ are the nonzero eigenvalues of $G^{\frac{1}{2}}PG^{\frac{1}{2}}$. Introduce 
\[
D = \diag(d_1, d_2, \ldots, d_K), \qquad \widetilde{D} = \diag(d_2, d_3, \ldots, d_K), 
\] 
and
\[
h = \Bigl(\frac{1_n' \xi_2}{1_n' \xi_1}, \frac{1_n' \xi_3}{1_n' \xi_1}, \ldots, \frac{1_n' \xi_K}{1_n' \xi_1}\Bigr)', 
\qquad u_0 = \sum_{k=2}^K \frac{d_k(1_n'\xi_k)^2}{d_1 (1_n'\xi_1)^2}. 
\] 
By Lemma~\ref{lem:vectorh}, $1_n'\xi_1>0$, so $h$ and $u_0$ are both well-defined. 
Write $\Xi = [\xi_1, \xi_2, \ldots, \xi_K]$. 
The following lemma gives an alternative expression of $\widetilde{\Omega}$. 
\begin{lemma} \label{lem:connection2}
Under the DCMM model,  
\[
\widetilde{\Omega} = \|\theta\|^2 \cdot \Xi M\Xi',  
\]
where $M$ is a $K \times K$ matrix satisfying 
\[
M =   
\left[ 
\begin{array}{ll} 
(1 + u_0)^{-1}  h' \widetilde{D}  h & - (1 + u_0)^{-1} h' \widetilde{D} \\ 
- (1 + u_0)^{-1}  \widetilde{D} h  &   \widetilde{D}  - (d_1 (1 + u_0))^{-1} \widetilde{D}  h h' \widetilde{D}  \\ 
\end{array} 
\right]. 
\] 
If additionally $|\lambda_2|/\lambda_1 \goto 0$, then for the matrix $\widetilde{M} \in \mathbb{R}^{K, K}$,  
\[
\widetilde{M} = 
\|\theta\|^2 \cdot 
\left[ 
\begin{array}{ll} 
h' \widetilde{D}  h & -   h' \widetilde{D} \\ 
-  \widetilde{D} h  &   \widetilde{D}    \\ 
\end{array} 
\right], 
\] 
we have 
\[
|M_{ij} - \widetilde{M}_{ij}|  \leq C \lambda_2^2/\lambda_1,  \qquad \mbox{for all $1 \leq i, j \leq K$}.   
\]  
\end{lemma}

We now study $\tr(\widetilde{\Omega}^3)$ and $\tr(\widetilde{\Omega}^4)$. They are related to the power of the SgnT test and SgnQ test, respectively. 
We discuss the two cases $|\lambda_2|/\lambda_1 \goto 0$ and 
$|\lambda_2|/\lambda_1 \geq c_0$ separately. 
Consider the case of $|\lambda_2| / \lambda_1 = o(1)$. 
Since $\widetilde{\Omega} = \Xi M \Xi'$, where $\Xi' \Xi = I_K$, we have 
\[
\tr(\widetilde{\Omega}^3) = \tr(M^3), \qquad \mbox{and} \qquad 
\tr(\widetilde{\Omega}^4) = \tr(M^4).  
\] 
The following lemma is proved in Section~\ref{supp:SpectralProof}.
\begin{lemma} \label{lem:tracetOmega}  
Consider the DCMM model, where \eqref{cond-balance} holds.
As $n\to\infty$, if $|\lambda_2|/\lambda_1 \goto 0$, then 
\begin{equation} \label{tracetOmega1topf1} 
|\tr(\widetilde{\Omega}^3) - \tr(\widetilde{M}^3)| \leq o(|\lambda_2|^3),   \qquad |\tr(\widetilde{\Omega}^4) - \tr(\widetilde{M}^4)| \leq o(|\lambda_2|^3),   
\end{equation} 
Moreover,  
\[ 
\tr(\widetilde{M}^3) = \tr(\widetilde{D}^3) + 3 h' \widetilde{D}^3  h + 3  (h' \widetilde{D} h) (h' \widetilde{D}^2 h) + (h' \widetilde{D} h)^3,  
\] 
and 
\begin{align*} 
\tr(\widetilde{M}^4) & = \tr(\widetilde{D}^4) + (h' \widetilde{D} h)^4 + 4 (h' \widetilde{D}^2 h)^2 + 4 (h' \widetilde{D} h)^2  (h' \widetilde{D}^2 h) + 4 h' \widetilde{D}^4 h + 4 (h' \widetilde{D} h) (h' \widetilde{D}^3 h) \cr   
&\geq  \tr(\widetilde{D}^4) + (h' \widetilde{D} h)^4 + 2[(h' \widetilde{D}^2 h)^2 +   (h' \widetilde{D} h)^2
(h' \widetilde{D}^2 h) +   h' \widetilde{D}^4 h] \label{trMpf0}    \cr
&\geq \tr(\widetilde{D}^4). 
\end{align*}  
\begin{itemize} 
\item In the special case where $\lambda_2, \lambda_3, \ldots, \lambda_K$ have the same signs, 
\[
|\tr(\widetilde{M}^3)| \geq |\sum_{k = 2}^K \lambda_k^3| = \sum_{k = 2}^K |\lambda_k|^3, 
\] 
and so 
\[
|\tr(\widetilde{\Omega}^3)|  \geq   \sum_{k = 2}^K |\lambda_k|^3 + o(|\lambda_2|^3). 
\] 
\item In the special case where $K = 2$, the vector $h$ is a scalar, and  
\[
\tr(\widetilde{M}^3) = (1 + h^2)^3 \lambda_2^3, \qquad \tr(\widetilde{M}^4) = (1 + h^2)^4 \lambda_2^4,  
\] 
and so 
\[
\tr(\widetilde{\Omega}^3)  =  [(1 + h^2)^3 + o(1)] \lambda_2^3, \qquad 
\tr(\widetilde{\Omega}^4)  =  [(1 + h^2)^4 + o(1)]  \lambda_2^4.  
\] 
\end{itemize} 
\end{lemma} 

We now consider the case $|\lambda_2/\lambda_1| \geq c_0$. In this case, $\widetilde{M}$ is not a good proxy for $M$ any more, so we can not derive a simple formula for $\tr(\widetilde{\Omega}^3)$ or 
$\tr(\widetilde{\Omega}^4)$ as above. However, for $\tr(\widetilde{\Omega}^4)$, since
\[
\tr(\widetilde{\Omega}^4)\geq \|\widetilde{\Omega}\|^4,
\]
by Lemma~\ref{lem:tildeOmega}, we immediately have
\begin{equation} \label{tildeOmegaPFAdd} 
\tr(\widetilde{\Omega}^4) \geq C\lambda_2^4\geq C(\sum_{k=2}^K\lambda_k^4)/(K-1)\geq C\sum_{k=2}^K\lambda_k^4. 
\end{equation} 


\subsection{Proof of Lemmas~\ref{lem:connection}-\ref{lem:tracetOmega}} \label{supp:SpectralProof}
 
\subsubsection{Proof of Lemma~\ref{lem:connection}}
The proof for the null case is straightforward, so we only prove the lemma for the alternative case.  
Consider the spectral decomposition 
\[
G^{1/2} P G^{1/2} = B D B'. 
\] 
where 
\[ 
D =\diag(d_1,\ldots,d_K) \qquad \mbox{and}  
\qquad  B =[\beta_1,\ldots,\beta_K].
\] 
Combining this with $\Omega=\Theta\Pi P\Pi'\Theta$ gives  
\begin{align*} 
\Omega & = \Theta \Pi G^{-\frac{1}{2}}(G^{\frac{1}{2}}PG^{\frac{1}{2}}) G^{-\frac{1}{2}}\Pi'\Theta\cr
&=  \Theta \Pi G^{-\frac{1}{2}}( B D B') G^{-\frac{1}{2}}\Pi'\Theta\cr
&= (\|\theta\|^{-1}\Theta\Pi G^{-\frac{1}{2}}B)( \|\theta\|^2 D)(\|\theta\|^{-1}\Theta\Pi G^{-\frac{1}{2}}B)'  \\
& = H (\|\theta\|^2 D) H', 
\end{align*}
where 
\[
H =  \|\theta\|^{-1}\Theta\Pi G^{-\frac{1}{2}} B.  
\] 
Recalling that $G = (\|\theta\|^2)^{-1} \cdot \Pi' \Theta^2 \Pi$, it is seen 
\begin{equation} 
H'H = \|\theta\|^{-2} B'G^{-\frac{1}{2}}(\Pi'\Theta^2\Pi) G^{-\frac{1}{2}} B  = B' B =I_K,  
\end{equation}
Therefore,  
\[
\Omega = H (\|\theta\|^2 D) H' 
\] 
is the spectral decomposition of $\Omega$.  Since $(\widetilde{D}_k, \xi_k)$ are the $k$-th eigenvalue of 
$\Omega$ and unit-norm eigenvector respectively, we have 
\[
\xi_k = \pm 1 \cdot  \mbox{the $k$-th column of H} =  \pm (\|\theta\|)^{-1} \Theta \Pi G^{-1/2} \beta_k. 
\] 
This proves the claim. \qed 
 
\subsubsection{Proof of Lemma~\ref{lem:Omega}}
Consider the first claim.  By Lemma~\ref{lem:connection}, $\lambda_1 = d_1 \|\theta\|^2$, where $d_1$ is the maximum eigenvalue of $G^{\frac{1}{2}}PG^{\frac{1}{2}}$. It suffices to show that $d_1 \asymp 1$. Since all entries of $P$ are upper bounded by constants, we have
\[
\|P\| \leq C. 
\]
Additionally, since $G$ is a nonnegative symmetric matrix, 
\beq \label{lem-Omega-1}
\|G\|\leq \|G\|_{\max} =\max_{1\leq k\leq K}\sum_{\ell=1}^KG(k,\ell)= \|\theta\|^{-2}\max_{1\leq k\leq K}\sum_{\ell=1}^K \sum_{i=1}^n\pi_i(k)\pi_i(\ell)\theta_i^2\leq 1. 
\eeq
It follows that 
\begin{equation} \label{d1UB} 
d_1\leq \|G\| \|P\| \leq C.
\end{equation}  
At the same time,  for any unit-norm non-negative vector $x \in \mathbb{R}^K$, since all entries of $P$ are non-negative and all diagonal entries of $P$ are $1$, 
\[
x' P x  \geq x' x = 1. 
\] 
It follows that
\[
d_1= \|G^{\frac{1}{2}}PG^{\frac{1}{2}}\|  \geq \frac{(G^{-\frac{1}{2}}x)'(G^{\frac{1}{2}}PG^{\frac{1}{2}})(G^{-\frac{1}{2}}x)}{\|(G^{-\frac{1}{2}}x)\|^2}=\frac{x'Px}{x'G^{-1}x}\geq \frac{1}{\|G^{-1}\|}. 
\]
Combining it with the assumption \eqref{cond-balance} gives 
\begin{equation} \label{d1LB} 
d_1\geq C. 
\end{equation} 
where we note $C$ denotes a generic constant which may vary from occurrence to occurrence. 
Combining (\ref{d1UB})-(\ref{d1LB}) gives the claim. 

Consider the second claim.  Let $B = [\beta_1, \beta_2, \ldots, \beta_K]$ and $D = \diag(d_1, d_2, \ldots, d_K)$ as in the proof of Lemma \ref{lem:connection}, where we note $B$ is orthonormal. 
By Lemma~\ref{lem:connection} and definitions, 
\[
u_i'  = \|\theta\|^{-1}\theta_i \pi_i'G^{-\frac{1}{2}}B.
\]
It follows that
\[
\|u_i\|\leq \|\theta\|^{-1}\theta_i \cdot \|\pi_i\|\|G^{-\frac{1}{2}}\|\|B\| \leq (\|\theta\|)^{-1} \theta_i \|G^{-1/2}\|, 
\]
where we have used $\|B\|=1$ and $\|\pi_i\|  =  [\sum_{k = 1}^K \pi_i(k)^2]^{1/2}  \leq 1$.  
Finally,  by the assumption \eqref{cond-balance},  $\|G^{-1}\|  \leq C$ and so $\|G^{-1/2}\| \leq C$. 
Combining these gives the claim.  \qed

\subsubsection{Proof of Lemma \ref{lem:vectorh}}  
It is sufficient to show the first two claims. Consider the first claim. 
By Lemma~\ref{lem:Omega}, for all $1 \leq k \leq K$ and $1 \leq i \leq n$,  
\[
|\xi_k(i)|\leq C\|\theta\|^{-1}\theta_i.
\]
It follows that 
\beq \label{lem-tOmega-4}
|1_n'\xi_k| \leq C\sum_{i=1}^n \|\theta\|^{-1}\theta_i \leq C\|\theta\|^{-1}\|\theta\|_1, \qquad \mbox{for all }1\leq k\leq K, 
\eeq
and the claim follows. 

Consider the second claim.  
By Lemma~\ref{lem:connection}, 
\begin{equation} \label{vectorhpf1} 
\xi_1 = \|\theta\|^{-1}\Theta \Pi (G^{-\frac{1}{2}}\beta_1), 
\end{equation} 
where $\beta_1$ is the (unit-norm) eigenvector of $G^{\frac{1}{2}}PG^{\frac{1}{2}}$ associated with $\lambda_1$, which is the largest eigenvalue of $G^{1/2} P G^{1/2}$. By basic algebra, $\lambda_1$ is 
also the largest eigenvalue of the matrix $PG$, with $G^{-1/2} \beta_1$ being the corresponding eigenvector. Since $PG$ is a nonnegative matrix, $G^{-\frac{1}{2}}\beta_1$ is a nonnegative vector (e.g., \cite[Theorem 8.3.1]{HornJohnson}). Denote for short by 
\[
h = G^{-1/2} \beta_1. 
\] 
It follows from (\ref{vectorhpf1}) that 
\beq \label{lem-tOmega-6}
1_n'\xi_1 = (\|\theta\|)^{-1} \cdot 1_n' \Theta \Pi h = \|\theta\|^{-1}  \cdot \sum_{k=1}^K \Bigl(\sum_{i=1}^n\pi_i(k)\theta_i\Bigr)h_k . 
\eeq
We note that $\sum_{k=1}^K \bigl(\sum_{i=1}^n \pi_i(k)\theta_i\bigr)=\|\theta\|_1$. Combining it with the assumption \eqref{cond-balance} yields 
\[
\min_{1\leq k\leq K}\Bigl\{ \sum_{i=1}^n \pi_i(k)\theta_i \Bigr\}\geq C \|\theta\|_1. 
\]
Inserting this into (\ref{lem-tOmega-6}) gives 
\begin{equation} \label{lem-tOmega-6b} 
1_n' \xi_1 \geq C (\|\theta\|)^{-1} \|\theta\|_1 \cdot \|h\|_1. 
\end{equation} 
We claim that $\|h\| \geq 1$. Otherwise, if $\|h\|  <  1$, then 
every entry of $h$ is no greater than $1$ in magnitude, and so 
\[
\|h\|_1 \geq \|h\|^2 = \|G^{-1}  \beta_1 \|^2. 
\]  
Since $\|G^{-1}\|=\|G\|^{-1}\geq 1$ (see \eqref{lem-Omega-1}) and $\|\beta_1\|=1$, 
\[
\|G^{-\frac{1}{2}}\beta_1\|\geq 1.  
\]
and so it follows $\|h\| \geq 1$. The contradiction show that $\|h\| \geq 1$. 
The claim follows by combining this with (\ref{lem-tOmega-6b}).   \qed

\subsubsection{Proof of Lemma \ref{lem:1Omega1}} 
For $1 \leq k \leq K$, let 
\[
c = (\|\theta\|_1)^{-1}  \Pi' \Theta 1_n = (\|\theta\|_1)^{-1} (1_n' \Theta \Pi)'.  
\]
Since $\Omega = \Theta \Pi  P \Pi' \Theta$ and all entries of $P$ are non-negative,    
\begin{equation} \label{1Omega1pf1} 
1_n' \Omega 1_n =  \|\theta\|_1^2 (c' P c )\geq \|\theta\|^2 \Bigl(\sum_{k = 1}^K c_k^2\Bigr).   
\end{equation} 
Note that,   first, $c_k \geq 0$, and second,  
$\|\theta\|_1\sum_{k=1}^K c_k = 1_n'\Pi\Theta 1_n = 1_n' \Theta 1_n$, where the last term is $\|\theta\|_1$, and so 
\[ 
\sum_{k = 1}^K  c_k  = 1.  
\] 
Together with the Cauchy-Schwartz inequality, we have  
\[
\sum_{k = 1}^K c_k^2 \geq (\sum_{k = 1}^K c_k)^2/K = 1/K. 
\]  
Combining this with (\ref{1Omega1pf1}) gives the claim.  \qed

\subsubsection{Proof of Lemma~\ref{lem:tildeOmega}}
Consider the first claim. We first derive a lower bound for $\|\widetilde{\Omega}\|$. By Lemma~\ref{lem:connection2}, 
\beq \label{lem-tOmega-1}
\|\widetilde{\Omega}\| =  \|\theta\|^2 \cdot \|M\|, 
\eeq
where with the same notations as in the proof of Lemma \ref{lem:connection2},  $M= D - (1 + u_0)^{-1}  vv'$. Let $M_0$ be the top left $2\times 2$ block of $M$. Let $D_0 =\diag(d_1, d_2)$, and let $v_0$ be the sub-vector of $v$ in \eqref{lem-connection2-2} restricted to the first two coordinates. By \eqref{lem-connection2-2}, 
\[
M_0 = D_0 - (1 + u_0)^{-1}  v_0 v_0' = D_0^{\frac{1}{2}} \Bigl(I_2 -  (1+u_0)^{-1}  D_0^{-1/2}  v_0 v_0' D_0^{-\frac{1}{2}}\Bigr) D_0^{\frac{1}{2}}, 
\]
and so by $\|D_0^{-1/2}\| = |d_2|^{-1/2}$ we have  
\begin{equation} \label{lemeigtOmegapf1} 
\|\Bigl(I_2 -  (1+u_0)^{-1}  D_0^{-1/2}  v_0 v_0' D_0^{-\frac{1}{2}}\Bigr) \| 
\leq \| D_0^{-1/2} M_0 D_0^{-1/2} \|  \leq |d_2|^{-1} \cdot \|M_0\|.  
\end{equation} 
At the same time,  since $(1 + u_0)^{-1} D_0^{-1/2} v_0 v_0' D_0^{-1/2}$ is a rank-$1$ matrix, there is an orthonormal matrix and a number $b$  such that 
\[
Q (1 + u_0)^{-1} D_0^{-1/2} v_0 v_0' D_0^{-1/2}   Q' = \diag(b, 0). 
\] 
It follows 
\[
\|\Bigl(I_2 -  (1+u_0)^{-1}  D_0^{-1/2}  v_0 v_0' D_0^{-\frac{1}{2}}\Bigr) \| 
  = \|I_2 - \diag(b, 0)\| = \max\{|1-b|, 1\} \geq 1. 
\] 
Inserting this into (\ref{lemeigtOmegapf1}) gives 
\[
\|M_0\| \geq |d_2|, 
\] 
Note that $\|M\| \geq \|M_0\|$. Combining this with (\ref{lem-tOmega-1}) gives 
\[
\|\widetilde{\Omega}\| \geq  |d_2| \|\theta\|^2.  
\]

Next, we derive an upper bound for $\|\widetilde{\Omega}\|$.  By Lemma \ref{lem:vectorh}, 
\beq \label{lem-tOmega-3}
\max_{1\leq k\leq K}|1_n'\xi_k|\leq C\|\theta\|^{-1}\|\theta\|_1, \qquad 1_n'\xi_1 \geq C\|\theta\|^{-1}\|\theta\|_1. 
\eeq
By (\ref{lem-tOmega-3}),  all the entries of $M$ are upper bounded by $C|\lambda_2|$, which implies $\|M\|\leq C|d_2|$. Plugging it into \eqref{lem-tOmega-1} gives
\begin{equation} \label{lem-tOmega-5} 
\|\widetilde{\Omega}\|\leq \frac{C}{|1+u_0|}|d_2|\|\theta\|^2, 
\end{equation} 
and all remains to show is 
\[
1 + u_0 \geq C > 0. 
\] 
Now, recalling that $\Omega = \sum_{k = 1}^K \lambda_k \xi_k \xi_k'$ and $\lambda_k = d_k \|\theta\|^2$,   by definitions, 
\[
d_1(1_n'\xi_1)^2 (1 + u_0) =\sum_{k=1}^K d_k(1_n'\xi_k)^2=\|\theta\|^{-2}1_n'\Omega 1_n.
\]
By Lemma~\ref{lem:1Omega1} which gives $1_n'\Omega 1_n\geq C \|\theta\|_1^2$. 
It follows that 
\[
1 + u_0 \geq \frac{\|\theta\|^{-2} 1_n'\Omega 1_n}{d_1(1_n'\xi_1)^2}\geq  C \frac{\|\theta\|^{-2}\cdot   \|\theta\|_1^2}{ \|\theta\|^{-2}\cdot \|\theta\|^2_1}\geq C, 
\]
where in the second inequality we have used \eqref{lem-tOmega-3} and $d_1  = (\|\theta\|)^{-2} \cdot  \lambda_1  \leq 1$ (see Lemma \ref{lem:Omega}).  Inserting this into (\ref{lem-tOmega-5}) gives the claim.

Consider the second claim.   
By Lemma~\ref{lem:connection2},    
\[
\widetilde{\Omega} = \Xi M \Xi',
\]
where $\Xi$ and $M$ are the same there.  Write 
 \[
\Xi = [\xi_1, \xi_2, \ldots, \xi_K] = [u_1, u_2, \ldots, u_n]'.   
\] 
Note that $\widetilde{\Omega}$ and $M$ have the same spectral norm. 
It follows that 
\[
\widetilde{\Omega}^m = \Xi  M^m \Xi',  
\] 
and 
\[
|(\widetilde{\Omega}^m)_{ij}| = |u_i' M^m u_j|\leq \|u_i\|\|M\|^m \|u_j\|. 
\] 
By Lemma~\ref{lem:Omega}, $\|u_i\| \|u_j\|\leq C\|\theta\|^{-2}\theta_i\theta_j$, and by the first part of the current lemma, 
\[
\|M \|  = \|\tilde{\Omega} \|  \leq  C|d_2| \|\theta\|^2.
\] 
It follows  that
\[
|(\widetilde{\Omega}^m)_{ij}| \leq C |d_2|^m\|\theta\|^{2m-2} \theta_i\theta_j. 
\]
This proves the claim. \qed

\subsubsection{Proof of Lemma~\ref{lem:connection2}}
Consider the first claim.  By definitions, 
\beq \label{lem-connection2-0}
\widetilde{\Omega} = \Omega -(\eta^*)(\eta^*)', \qquad \mbox{where}\;\; \eta^* =\frac{1}{\sqrt{1_n'\Omega 1_n}}\Omega 1_n. 
\eeq
Recalling $\widetilde{D}_k = d_k \|\theta\|^2$ and  $\Xi = [\xi_1, \xi_2, \ldots, \xi_K]$, we have 
\[
\Omega =\sum_{k=1}^K \widetilde{D}_k  \xi_k\xi_k' = \|\theta\|^2 \cdot  \Xi D \Xi'. 
\] 
It follows  that 
\[
1_n'\Omega 1_n=\|\theta\|^2 \sum_{k=1}^K d_k(1_n'\xi_k)^2, 
\]
and
\begin{align*}
\eta^*   =\frac{\|\theta\|}{\sqrt{\sum_{s=1}^K d_s(1_n'\xi_s)^2}} \sum_{k=1}^K d_k(1_n'\xi_k)\xi_k  = \frac{\|\theta\|}{\sqrt{(1+u_0)}}\biggl[ \sqrt{d_1}\, \xi_1 +\sum_{k=2}^K\frac{d_k(1_n'\xi_k)}{\sqrt{d_1}(1_n'\xi_1)} \xi_k\biggr],  
\end{align*}
where the vector in the big bracket on the right is $\Xi v$, if we let 
\[ 
v =   (\sqrt{d_1},  
\frac{d_2(1_n'\xi_2)}{\sqrt{d_1}(1_n'\xi_1)}, \ldots,  \frac{d_K(1_n'\xi_K)}{\sqrt{d_1}(1_n'\xi_1)})'.  
\] 
Combining these gives 
\[
\widetilde{\Omega} = \|\theta\|^2  \Xi  D \Xi' - \frac{\|\theta\|^2}{1 + u_0} \Xi v v' \Xi.  
\] 
Plugging it into \eqref{lem-connection2-0} gives 
\beq \label{lem-connection2-2}
\widetilde{\Omega} = \|\theta\|^2  \Xi  D \Xi' - \frac{\|\theta\|^2}{1 + u_0} \Xi v v' \Xi  = \|\theta\|^2 \Xi (D - (1 + u_0)^{-1} v v') \Xi'. 
\end{equation} 
By definitions, 
\[
D = \diag(d_1, d_2, \ldots, d_K), \qquad \mbox{and} \qquad v = d_1^{-1/2} \cdot (d_1, h' \widetilde{D})'. 
\] 
It follows 
\[
D - (1 + u_0)^{-1} vv' = \left[ 
\begin{array}{cc} 
(1 + u_0)^{-1}  d_1 u_0  &  - (1 + u_0)^{-1} h' \widetilde{D} \\ 
-(1 + u_0)^{-1} \widetilde{D} h  & \widetilde{D}  - (d_1(1 + u_0))^{-1}  \widetilde{D} h h' \widetilde{D} \\ 
\end{array} 
\right],  
\]
where we note that 
\[
d_1 u_0 = \sum_{s = 2}^K d_s \frac{(1_n' \xi_s)^2}{(1_n' \xi_1)^2} = h' \widetilde{D} h, 
\] 
Combining these gives the claim.

Consider the second claim. By definitions,  
\[
M - \widetilde{M} = \|\theta\|^2 \cdot  \left[
\begin{array}{cc}  
[(1 + u_0)^{-1} - 1] d_1 u_0 & (1 - (1 + u_0)^{-1}) h' \widetilde{D} \\ 
(1 - (1 + u_0)^{-1}) \widetilde{D} h  & - (d_1(1+u_0))^{-1} \widetilde{D} h h' \widetilde{D} \\ 
\end{array} 
\right].  
\] 
Note that 
\[
|1 - (1 + u_0)^{-1}| \leq C |u_0| \leq C |\widetilde{D}_2|/\widetilde{D}_1,   
\] 
and that by Lemma \ref{lem:vectorh}, 
\[
|(1_n' \xi_k)| \leq C 1_n' \xi_1,
\] 
and so each entry of $\widetilde{D} h$ does not exceed $C|d_2|$.  
It follows that for all $2 \leq i, j \leq K$,   
\[
|M_{1i} - \widetilde{M}_{1i}| \leq C\|\theta\|^2 (|\widetilde{D}_2|/\widetilde{D}_1)  d_2^2  \leq C \widetilde{D}_2^2/\widetilde{D}_1, 
\] 
and 
\[
|M_{ij} - \widetilde{M}_{ij}| \leq C\|\theta\|^2 d_1^{-1}   d_2^2  \leq C \widetilde{D}_2^2/\widetilde{D}_1.  
\] 
Finally,  
\[
d_1 u_0^2   = d_1^{-1}  (\sum_{s = 2} d_2 \frac{(1_n' \xi_s)^2}{(1_n' \xi_1)^2})^2 \leq C d_2^2/d_1, 
\] 
so 
\[
|M_{11} - \widetilde{M}_{11}| \leq C \|\theta\|^2  d_2^2/d_1  \leq C  \widetilde{D}_2^2/\widetilde{D}_1.   
\] 
Combining these gives the claim.  \qed

\subsubsection{Proof of Lemma \ref{lem:tracetOmega}} 
It is sufficient to show (\ref{tracetOmega1topf1}). In fact, 
once (\ref{tracetOmega1topf1}) is proved, other claims follow  
by direct calculations, except for the first inequality regarding $\tr(\widetilde{\Omega}^4)$, we have used 
\[
|(h' \widetilde{D} h) (h' \widetilde{D}^3 h)|  \leq  |h' \widetilde{D} h|    \sqrt{(h' \widetilde{D}^2 h) (h' \widetilde{D}^4 h)} \leq \frac{1}{2} \biggl[ (h' \widetilde{D} h)^2 (h' \widetilde{D}^2 h) + h' \widetilde{D}^4 h  \biggr].  
\] 

We now show (\ref{tracetOmega1topf1}). 
Since $\tr(\widetilde{\Omega}^m) = \tr(\widetilde{M}^m)$, for $m = 3, 4$, it is sufficient to show 
\begin{equation} \label{tracetOmegatopf2}
|\tr(M^3) - \tr(\widetilde{M}^3)| \leq C \lambda_2^4/\lambda_1),  \qquad 
|\tr(M^4) - \tr(\widetilde{M}^4)| \leq  C |\lambda_2|^5/\lambda_1. 
\end{equation}
Since the proofs are similar, we only show the first one.  
By basic algebra,  
\[
\tr(M^3 - \widetilde{M}^3) = \tr((M - \widetilde{M})^3) + 3 \tr(\widetilde{M} (M - \widetilde{M})^2) + 
3 \tr(\widetilde{M}^2 (M - \widetilde{M})).  
\] 
By Lemma \ref{lem:connection2},   for all $1 \leq i, j \leq K$,  
\[
|M_{ij} - \widetilde{M}_{ij}| \leq C\lambda_2^2/\lambda_1.   
\] 
Also, by Lemma \ref{lem:vectorh}, all entries of $h$ are bounded, so for all $1 \leq i, j \leq K$,  
\[
|\widetilde{M}_{ij}|  \leq |\lambda_2|. 
\]  
It follows 
\[
|\tr((M - \widetilde{M})^3| \leq C (\lambda_2^2/\lambda_1)^3,  
\]  
\[
|\tr(\widetilde{M} (M - \widetilde{M})^2)| \leq C |\lambda_2|  (\lambda^2/\lambda_1)^2 \leq C |\lambda_2|^5/\lambda_1^2,   
\] 
and 
\[
|\tr(\widetilde{M}^2 (M - \widetilde{M})| \leq C \lambda_2^2  (\lambda^2/\lambda_1) \leq C \lambda_2^4/\lambda_1.    
\] 
where we note that $\lambda_2/ \lambda_1 = o(1)$. 
Combining these gives the claim.

\subsection{Proof of Lemmas \ref{lem:tracetOmega1} and \ref{lem:tracetOmega2}} 
\label{subsec:twolemmas} 
Lemma \ref{lem:tracetOmega1}  follows directly from Lemma \ref{lem:tracetOmega} of this appendix. 
Consider Lemma \ref{lem:tracetOmega2}. The second bullet point is a direct result of (\ref{tildeOmegaPFAdd}), and the other two bullet points follow directly from Lemma \ref{lem:tracetOmega} of this appendix.

\newpage
\section{Lower bounds, Region of Impossibility}  \label{supp:LB} 
We study the Region of Impossibility by considering a DCMM with random mixed memberships. First, in Section~\ref{subsec:equivalence}, we establish the equivalence between regularity conditions for a DCMM with non-random mixed memberships and those for a DCMM with random mixed memberships. Next, we prove Lemma~\ref{lem:existence}, which is key to the construction of inseparable hypothesis pairs. Last, we prove Theorem~\ref{thm:minimax}-\ref{thm:RI-sameTheta} of the main article. 

\subsection{Equivalence of regularity conditions} \label{subsec:equivalence}
Let $\mu_1, \mu_2, \ldots, \mu_K$ be the eigenvalues of $P$, 
arranged in the descending order in magnitude. Recall that $\lambda_1,\lambda_2,\ldots,\lambda_K$ are the eigenvalues of $\Omega$. The following lemma is proved in Section~\ref{subsec:Proof-of-Lemma-M}. 
\begin{lemma}[Equivalent definition of Region of Impossibility] \label{lem:eigP} 
Consider the DCMM model \eqref{model1a}-\eqref{condition1d}, where the alternative is true and the condition \eqref{cond-balance} holds. Suppose $\theta_{\max}\to 0$ and $\|\theta\| \goto \infty$ as $n \goto \infty$. Then, as $n \goto \infty$,  
\[
\mu_1 \asymp 1, \qquad   \frac{|\mu_2|}{\mu_1} \asymp \frac{|\lambda_2|}{\lambda_1}, \qquad \max_{1 \leq i, j \leq K} |P_{ij} - 1| \leq C(|\lambda_2|/\lambda_1). 
\] 
As a result, $|\lambda_2|/\sqrt{\lambda_1}\to 0$ if and only if $\|\theta\|\cdot|\mu_2(P)|\to 0$. 
\end{lemma}

We now consider DCMM with random mixed memberships: Given $(\Theta, P)$ and a distribution $F$ over $V$ (the standard simplex in $\mathbb{R}^K$), let
\beq \label{DCMM-random}
\Omega = \Theta \Pi P\Pi'\Theta, \qquad \Pi=[\pi_1,\pi_2,\ldots,\pi_n]', \qquad \pi_i\overset{iid}{\sim}F. 
\eeq
We notice that the conclusion of Lemma~\ref{lem:eigP} holds provided that the regularity condition \eqref{cond-balance} is satisfied. The next lemma shows that \eqref{cond-balance} holds with high probability. It is proved in Section~\ref{subsec:Proof-of-Lemma-M}.
\begin{lemma}[Equivalence of regularity conditions] \label{lem:Regu}  
Consider the model \eqref{DCMM-random}. Let $h=\mathbb{E}[\pi_i]$ and $\Sigma=\mathrm{Cov}(\pi_i)$. Suppose $\|P\|\leq C$, $\min_{1\leq k\leq K}\{h_k\}\geq C$ and $\|\Sigma^{-1}\|\leq C$. Suppose $\theta_{\max}\to 0$, $\|\theta\| \goto \infty$, and $(\|\theta\|^2/\|\theta\|_1)\sqrt{\log(\|\theta\|_1)}\to 0$, as $n \goto \infty$. Then, as $n \goto \infty$,  with probability $1-o(1)$, the condition \eqref{cond-balance} is satisfied, i.e.,  
\[
\frac{\max_{1\leq k\leq K}\{\sum_{i=1}^n \theta_i\pi_i(k)\}}{\min_{1\leq k\leq K}\{\sum_{i=1}^n \theta_i\pi_i(k)\}}\leq C_0, \qquad \|G^{-1}\|\leq C_0,
\] 
for a constant $C_0>0$ and $G=\|\theta\|^{-2}(\Pi'\Theta^2\Pi)$. 
\end{lemma}

\subsection{Proof of Lemma~\ref{lem:existence}} 
Let  $M = \diag(\mu_1, \mu_2, \ldots, \mu_K)$. It is seen 
$\mu = M 1_K$ and so the desired result is to find a $D$ such that 
\[
D A D M 1_K = 1_K    \Longleftrightarrow   
M DAD M 1_K = M 1_K = \mu    \Longleftrightarrow    D (M A M) D 1_K = \mu. 
\] 
Since $M A M$ has strictly positive entries, it is sufficient to show that for any matrix $A$ ($MAM$ in our case; a slight misuse notation here)  
with strictly positive entries, there is a unique diagonal matrix $D$ with strictly positive diagonal entries such that 
\begin{equation} \label{DADexistence} 
D A D 1_k = \mu_K. 
\end{equation} 
We now show the existence and uniqueness separately. 

For existence, we follow the proof in \cite{DADOlkin}. Consider $d' A d$ for a vector $d \in \mathbb{R}^K$ with strictly positive entries.  It is shown there that $d' A d$ can be minimized using Lagrange multiplier: 
\[
\frac{1}{2} d' A d -  \lambda \sum_{k = 1}^K \mu_k \log(d_k). 
\] 
Differentiating with respect to $d$ and set the derivative to $0$ gives 
\begin{equation} \label{DADexistence1} 
A d = \lambda  \sum_{k = 1}^K \mu_k / d_k, 
\end{equation} 
where $\lambda = d' A d/(\sum_{k = 1}^K \mu_k) > 0$. Letting $D = \lambda^{-1/2} \diag(d_1, d_2, \ldots, d_K)$. 
It is seen that (\ref{DADexistence1}) can be rewritten as 
\[ 
D A D 1_K =  \mu,  
\] 
and the claim follows. 

For uniqueness, we adapt the proof in \cite{DADJohnson} to our case. Suppose there are two different eligible diagonal matrices  $D_1$ and $D_2$ satisfying (\ref{DADexistence}). Let $d_1 = D_1 1_K$ and $d_2 = D_2 1_K$, and let 
$M = \diag(\mu_1, \mu_2, \ldots, \mu_K)$. It follows that 
\[
D_2 D_1 A d_1 = D_2  D_1 A D_1  1_K = D_2 \mu  = M  d_2,  
\] 
and so 
\[
M^{-1} D_2 D_1 A d_1 = d_2. 
\] 
Now, for a diagonal matrix $S$ with strictly positive diagonal entries to be determined, we have 
\[
S^{-1} M^{-1} D_2 D_1 A S S^{-1}  d_1 = S^{-1}  d_2. 
\] 
We pick $S$ such that 
\[
S^{-1} M^{-1} D_2 D_1 = S, 
\]  
and denote such an $S$ by $S_0$. It follows  
\[
S_0 A S_0  (S_0^{-1} d_1) = S_0^{-1} d_2. 
\] 
or equivalently, if we let $\tilde{d}_1 = S_0^{-1} d_1$ and $\tilde{d}_2 = S_0^{-1} d_2$, 
\begin{equation} \label{DADunique1} 
S_0 A S_0 \tilde{d}_1= \tilde{d}_2. 
\end{equation} 

Similarly, by switching the places of $D_1$ and $D_2$, we have 
\begin{equation} \label{DADunique2} 
S_0 A S_0  \tilde{d}_2 = \tilde{d}_1. 
\end{equation} 
Combining (\ref{DADunique1}) and (\ref{DADunique2}) gives 
\[
S_0  A S_0 (\tilde{d}_1 +  \tilde{d}_2) = (\tilde{d}_1 + \tilde{d}_2),   \qquad   \mbox{and} \qquad S_0 A S_0 (\tilde{d}_1 -  \tilde{d}_2) = - (\tilde{d}_1 - \tilde{d}_2). 
\] 
This implies that $1$ and $-1$ are the two eigenvalues of $S_0 A S_0$, with $\tilde{d}_1 + \tilde{d}_2$ 
and $\tilde{d}_1 - \tilde{d}_2$ being the corresponding eigenvectors, respectively, 
where we note that especially, $\tilde{d}_1 + \tilde{d}_2$ has all strictly positive entries. 
By Perron's theorem \cite{HornJohnson}, since $S_0 A S_0$ have all  strictly positive entries, 
the eigenvector corresponding to the largest eigenvalue (i.e., the Perron root)   have  all strictly positive entries. 
As for any symmetric matrix, we can only have one eigenvector that has all strictly positive entries, so $1$ must be the Perron root of $S_0 A S_0$.  Using Perron's Theorem again,  all eigenvalues of $S_0 A S_0$ except the Perron root itself  should be strictly smaller than $1$ in magnitude. 
This contradicts with the fact that $-1$ is an eigenvalue of $S_0 A S_0$. The contradiction proves the uniqueness. \qed

\subsection{Proof of Theorem~\ref{thm:minimax}}
This theorem follows easily from Theorems~\ref{thm:RI-DCBM}-\ref{thm:RI-sameTheta}. Fix $(\Theta, P, F)$ such that $\theta\in{\cal M}^*_n(\beta_n/\log(n))$ and $\|\theta\|\cdot|\mu_2(P)|\leq \alpha_n/\log(n)$. Consider a sequence of hypotheses indexed by $n$, where $\Omega=\theta\theta'$ under $H_0^{(n)}$, and $\Omega$ follows the construction in any of  Theorems~\ref{thm:RI-DCBM}-\ref{thm:RI-sameTheta} under $H_1^{(n)}$. Let $P_0^{(n)}$ and $P_1^{(n)}$ be the probability measures associated with two hypotheses, respectively. By those theorems, the $\chi^2$-distance satisfy
\[
{\cal D}(P_0^{(n)}, P_1^{(n)})=o(1), \qquad \mbox{as }n\to\infty. 
\]
By connection between $L^1$-distance and $\chi^2$-distance, it follows that  
\[
\|P_0^{(n)}-P_1^{(n)}\|_1=o(1), \qquad \mbox{as }n\to\infty.  
\]
We now slightly modify the alternative hypothesis. Let $\Pi_0$ be a non-random membership matrix such that $(\theta,\Pi_0, P)\in {\cal M}_n(K, c_0, \alpha_n, \beta_n)$. In the modified alternative hypothesis $\widetilde{H}_1^{(n)}$, 
\[
\Pi = \begin{cases}
\widetilde{\Pi}, & \mbox{if }(\theta,\widetilde{\Pi}, P)\in {\cal M}_n(K, c_0, \alpha_n, \beta_n),\\
\Pi_0, & \mbox{otherwise}, 
\end{cases}
\qquad\mbox{where}\quad \widetilde{\pi}_i\overset{iid}{\sim}F. 
\]
Let $\widetilde{P}_1^{(n)}$ be the probability measure associated with $\widetilde{H}_1^{(n)}$.  
By Lemmas~\ref{lem:eigP}-\ref{lem:Regu}, $\Pi=\widetilde{\Pi}$, except for vanishing probability. It follows that
\[
\|P_1^{(n)}-\widetilde{P}_1^{(n)}\|_1=o(1), \qquad \mbox{as }n\to\infty. 
\]
Under $\widetilde{H}_1^{(n)}$, all realizations $(\theta,\Pi,P)$ are in the class ${\cal M}_n(K, c_0, \alpha_n, \beta_n)$. By Neymann-Pearson lemma and elementary inequalities,
\begin{align*}
 &\inf_{\psi}\Bigl\{ \sup_{\theta\in {\cal M}^*_n(\beta_n)}\mathbb{P}(\psi=1) + \sup_{(\theta,\Pi,P)\in {\cal M}_n(K, c_0, \alpha_n, \beta_n)} \mathbb{P}(\psi=0)\Bigr\}\cr
\geq\; &1- \inf_{f_0\in{\cal M}^*_n(\beta_n),f_1\in {\cal M}_n(K, c_0, \alpha_n, \beta_n)}\{ \|f_0-f_1\|_1\}\cr
\geq\;& 1-\|P_0^{(n)}-\widetilde{P}_1^{(n)}\|_1\cr
\geq\; & 1-\|P_0^{(n)}-P_1^{(n)}\|_1 - \|P_1^{(n)}-\widetilde{P}_1^{(n)}\|_1\cr
\geq\;& 1-o(1),
\end{align*}
where in the second line we have mis-used the notation $f\in {\cal M}_n(K, c_0, \alpha_n, \beta_n)$ to denote the probability density for a DCMM with non-random mixed memberships whose parameters are in the class ${\cal M}_n(K, c_0, \alpha_n, \beta_n)$.

%
%

\subsection{Proof of Theorems~\ref{thm:RI-DCBM}-\ref{thm:RI-sameTheta}}
We note that Theorem~\ref{thm:RI-DCBM}, Theorem~\ref{thm:RI-differentPi} and Theorem~\ref{thm:RI-sameTheta} can be deduced from Theorem~\ref{thm:RI-DCMM}. To see this, recall that Theorem~\ref{thm:RI-DCMM} assumes there exists a positive diagonal matrix $D$ such that 
\beq \label{LF-new-0}
D P D \widetilde{h}_D = 1_K, \qquad \min_{1 \leq k \leq K} \{\widetilde{h}_{D,k}\} \geq C, 
\eeq
where $\widetilde{h}_D = \mathbb{E}[D^{-1} \pi_i / \| D^{-1} \pi_i\|_1]$. We show that the condition \eqref{LF-new-0} is implied by conditions of other theorems. Theorem~\ref{thm:RI-DCBM} assumes $\pi_i\in\{e_1,e_2,\ldots,e_K\}$. It follows that $D^{-1} \pi_i / \| D^{-1} \pi_i\|_1=\pi_i$, and so $\widetilde{h}_D=h$. By Lemma~\ref{lem:existence}, there exists $D$ such that $DPDh=1_K$, hence, \eqref{LF-new-0} is satisfied. Theorem~\ref{thm:RI-differentPi} constructs the alternative hypothesis using $\widetilde{\pi}_i=D\pi_i/\|D\pi_i\|_1$. Equivalently, $D^{-1}\widetilde{\pi}_i/\|D^{-1}\pi_i\|_1=\pi_i$, and so $\widetilde{h}_D$ becomes $h$. Since $DPDh=1_K$, condition \eqref{LF-new-0} holds. Theorem~\ref{thm:RI-sameTheta} assumes $Ph=q_n 1_K$. Let $D=q_n^{-1/2}I_K$. Then, $\widetilde{h}_D=h$ and $DPDh =q_n^{-1}Ph= 1_K$. Again, \eqref{LF-new-0} is satisfied. 

We only need to prove Theorem~\ref{thm:RI-DCMM}. Let $P_0^{(n)}$ and $P_1^{(n)}$ be the probability measure associated with $H_0^{(n)}$ and $H_1^{(n)}$, respectively. Let ${\cal D}(P_0^{(n)}, P_1^{(n)})$ be the chi-square distance between two probability measures. By elementary probability, 
\[
{\cal D}(P_0^{(n)}, P_1^{(n)})=\int \biggl[\frac{dP_1^{(n)}}{dP_0^{(n)}}\biggr]^2dP_0^{(n)}-1. 
\]
It suffices to show that, when $\|\theta\|\cdot\mu_2(P)\to 0$, 
\beq \label{LB0}
\int \biggl[\frac{dP_1^{(n)}}{dP_0^{(n)}}\biggr]^2dP_0^{(n)}  = 1+o(1). 
\eeq
  
Let $p_{ij}$ and $q_{ij}(\Pi)$ be the corresponding $\Omega_{ij}$ under the null and the alternative, respectively.  It is seen that
\[
dP_0^{(n)} = \prod_{i<j} p_{ij}^{A_{ij}}(1-p_{ij})^{1-A_{ij}}, \qquad  
dP_1^{(n)} = \mathbb{E}_{\Pi}\Bigl[ \prod_{i<j} [q_{ij}(\Pi)]^{A_{ij}} [1-q_{ij}(\Pi)]^{1-A_{ij}} \Bigr]. 
\]
Let $\tilde{\Pi}$ be an independent copy of $\Pi$. Then,
\begin{align*}
\biggl[\frac{dP_1^{(n)}}{dP_0^{(n)}}\biggr]^2 & = \mathbb{E}_\Pi\biggl[ \prod_{i<j} \Bigl( \frac{q_{ij}(\Pi)}{p_{ij}}\Bigr)^{A_{ij}} \Bigl( \frac{1-q_{ij}(\Pi)}{1-p_{ij}}\Bigr)^{1-A_{ij}}  \biggr]\cdot \mathbb{E}_{\tilde{\Pi}}\biggl[ \prod_{i<j} \Bigl( \frac{q_{ij}(\tilde{\Pi})}{p_{ij}}\Bigr)^{A_{ij}} \Bigl( \frac{1-q_{ij}(\tilde{\Pi})}{1-p_{ij}}\Bigr)^{1-A_{ij}}  \biggr]\cr
& = \mathbb{E}_{\Pi,\tilde{\Pi}}\biggl [\underbrace{\prod_{i<j} \Bigl( \frac{q_{ij}(\Pi)q_{ij}(\tilde{\Pi})}{p^2_{ij}}\Bigr)^{A_{ij}} \Bigl( \frac{[1-q_{ij}(\Pi)][1-q_{ij}(\tilde{\Pi})]}{[1-p_{ij}]^2}\Bigr)^{1-A_{ij}}}_{S(A, \Pi,\tilde{\Pi})}  \biggr].  
\end{align*}
It follows that
\begin{align*}
\int \biggl[\frac{dP_1^{(n)}}{dP_0^{(n)}}\biggr]^2dP_0^{(n)} &= \mathbb{E}_A \biggl[\frac{dP_1^{(n)}}{dP_0^{(n)}}\biggr]^2\cr
 & =\mathbb{E}_{A, \Pi,\tilde{\Pi}}[ S(A,\Pi,\tilde{\Pi})]\cr
 & =\mathbb{E}_{\Pi,\tilde{\Pi}}\bigl\{\mathbb{E}_A\bigl[S(A,\Pi,\tilde{\Pi})|\Pi,\tilde{\Pi}\bigr]\bigr\}, 
\end{align*}
where the distribution of $A|(\Pi,\tilde{\Pi})$ is under the null hypothesis. Under the null hypothesis, $A$ is independent of $(\Pi,\tilde{\Pi})$, the upper triangular entries of $A$ are independent of each other, and $A_{ij}\sim \mathrm{Bernoulli}(p_{ij})$. It follows that
\begin{align*} 
\mathbb{E}_A\bigl[S(A,\Pi,\tilde{\Pi})|\Pi,\tilde{\Pi}\bigr] & = \prod_{i<j}\mathbb{E}_{A}\biggl[ \Bigl( \frac{q_{ij}(\Pi)q_{ij}(\tilde{\Pi})}{p^2_{ij}}\Bigr)^{A_{ij}} \Bigl( \frac{[1-q_{ij}(\Pi)][1-q_{ij}(\tilde{\Pi})]}{[1-p_{ij}]^2}\Bigr)^{1-A_{ij}}\bigg|\Pi,\tilde{\Pi}\biggr]\cr
& = \prod_{i<j}\biggl\{ p_{ij}\frac{q_{ij}(\Pi)q_{ij}(\tilde{\Pi})}{p^2_{ij}} + (1-p_{ij})\frac{[1-q_{ij}(\Pi)][1-q_{ij}(\tilde{\Pi})]}{[1-p_{ij}]^2}\biggr\}\cr
& = \prod_{i<j}\biggl\{ \frac{q_{ij}(\Pi)q_{ij}(\tilde{\Pi})}{p_{ij}} + \frac{[1-q_{ij}(\Pi)][1-q_{ij}(\tilde{\Pi})]}{1-p_{ij}}\biggr\}.
\end{align*}
Let $\Delta_{ij}=q_{ij}(\Pi)-p_{ij}$ and $\tilde{\Delta}_{ij}=q_{ij}(\tilde{\Pi})-p_{ij}$. By direct calculations,
\[
\frac{q_{ij}(\Pi)q_{ij}(\tilde{\Pi})}{p_{ij}} + \frac{[1-q_{ij}(\Pi)][1-q_{ij}(\tilde{\Pi})]}{1-p_{ij}} =1+\frac{\Delta_{ij}\tilde{\Delta}_{ij}}{p_{ij}(1-p_{ij})}.
\]
Combining the above gives
\beq \label{LF-new-1}
\int \biggl[\frac{dP_1^{(n)}}{dP_0^{(n)}}\biggr]^2dP_0^{(n)} = \mathbb{E}_{\Pi,\tilde{\Pi}}\biggl[ \prod_{i<j}\Bigl(1 + \frac{\Delta_{ij}\tilde{\Delta}_{ij}}{p_{ij}(1-p_{ij})} \Bigr)\biggr]. 
\eeq

We then plug in the expressions of $\Delta_{ij}$ and $\tilde{\Delta}_{ij}$ from the model. 
Let $D$ be the matrix in \eqref{LF-new-0}. Introduce $M=DPD-{\bf 1}_K{\bf 1}_K'$. We re-write
\[
DPD = {\bf 1}_K{\bf 1}_K' + M. 
\]
It is seen that $M\widetilde{h}_D ={\bf 0}_K$. 
The following lemma is proved in Section~\ref{subsec:Proof-of-Lemma-M}.
\begin{lemma} \label{lem:Mnorm}
Under the conditions of Theorem~\ref{thm:RI-DCMM}, $\|M\|\leq C|\mu_2(P)|$. 
\end{lemma}
Write for short $\pi_i^{D} = \frac{1}{\|D^{-1}\pi_i\|_1}D^{-1}\pi_i$ and $y_i= \pi_i^D - \mathbb{E}[\pi_i^D]=\pi_i^D-\widetilde{h}_D$. 
Under the alternative hypothesis, 
\begin{align*}
q_{ij}(\Pi)& =\theta_i\theta_j\|D^{-1}\pi_i\|_1\|D^{-1}\pi_j\|_1\cdot\pi_i'P\pi_j\cr
&= \theta_i\theta_j\cdot (\pi_i^D)'(DPD)(\pi_j^D)\cr
&=\theta_i\theta_j\cdot (\pi_i^D)'({\bf 1}_K{\bf 1}_K'+M)(\pi_j^D)\cr
&= \theta_i\theta_j\cdot \bigl[ 1 + (\pi_i^D)'M(\pi_j^D)\bigr]\cr
&=\theta_i\theta_j\cdot \bigl[1+(\widetilde{h}_D+y_i)'M(\widetilde{h}_D+y_j)]\cr
&= \theta_i\theta_j\cdot (1 + y_i'My_j).
\end{align*}
Here, the fourth line is due to ${\bf 1}_K'\pi_i=1$ and the last line is due to $M\widetilde{h}_D={\bf 0}_K$. 
Under the null hypothesis, $p_{ij}=\theta_i\theta_j$. As a result,
\[ 
\Delta_{ij} =\theta_i\theta_j \cdot y_i'My_j, \qquad y_i\equiv \pi^D_i-\mathbb{E}[\pi^D_i]. 
\]
Similarly, $\tilde{\Delta}_{ij}=\theta_i\theta_j\cdot \tilde{y}_i'M\tilde{y}_j$, with $\tilde{y}_i=\tilde{\pi}_i^D-\mathbb{E}[\tilde{\pi}_i^D]$. 
We plug them into \eqref{LF-new-1} and use $p_{ij}=\theta_i\theta_j$. It gives
\beq \label{LF-new-2}
\int \biggl[\frac{dP_1^{(n)}}{dP_0^{(n)}}\biggr]^2dP_0^{(n)} = \mathbb{E}\biggl[ \prod_{i<j}\Bigl(1 + \frac{\theta_i\theta_j}{1-\theta_i\theta_j}(y_i'My_j)(\tilde{y}_i'M\tilde{y}_j) \Bigr)\biggr],
\eeq
where $\{y_i, \tilde{y}_i\}_{i=1}^n$ are $iid$ random vectors with $\mathbb{E}[y_i]={\bf 0}_K$. 

We bound the right hand side of \eqref{LF-new-2}. Since $1+x\leq e^x$ for all $x\in\mathbb{R}$, 
\[
{\cal D}(P_0^{(n)}, P_1^{(n)})\leq \mathbb{E}[\exp(S)], \qquad \mbox{where}\quad S \equiv \sum_{i<j}\frac{\theta_i\theta_j}{1-\theta_i\theta_j}(y_i'My_j)(\tilde{y}_i'M\tilde{y}_j). 
\]
Let $M=\sum_{k=1}^K\delta_k b_kb_k'$ be the eigen-decomposition of $M$. Then,
\[
(y_i'My_j)(\tilde{y}_i'M\tilde{y}_j) =\sum_{1\leq k,\ell\leq K}\delta_k\delta_\ell (b_k'y_i)(b_k'y_j)(b_\ell'\tilde{y}_i)(b_{\ell}'\tilde{y}_j).
\]
This allows us to decompose 
\[
S = \frac{1}{K^2}\sum_{1\leq k,\ell\leq K}S_{k\ell}, \qquad\mbox{where}\;\; S_{k\ell}=K^2\delta_k\delta_\ell \sum_{i<j}\frac{\theta_i\theta_j}{1-\theta_i\theta_j}(b_k'y_i)(b_k'y_j)(b_\ell'\tilde{y}_i)(b_{\ell}'\tilde{y}_j). 
\]
By Jensen's inequality, $\exp(\frac{1}{K^2}\sum_{k,\ell}S_{k\ell})\leq \frac{1}{K^2}\sum_{k,\ell}\exp(S_{k\ell})$. It follows that
\beq \label{LF-new-3}
\int \biggl[\frac{dP_1^{(n)}}{dP_0^{(n)}}\biggr]^2dP_0^{(n)} \leq \mathbb{E}[\exp(S)]\leq \max_{1\leq k,\ell\leq K}\mathbb{E}[\exp(S_{k\ell})]. 
\eeq

We now fix $(k,\ell)$ and derive a bound for $\mathbb{E}[\exp(S_{k\ell})]$. For $n$ large enough, $\theta_{\max}\leq 1/2$ and $K^4\|M\|^2\|\theta\|^2\leq 1/9$. By Taylor expansion of $(1-\theta_i\theta_j)^{-1}$,
\begin{align*}
S_{k\ell} &= K^2\delta_k\delta_\ell \sum_{i<j}\sum_{m=1}^\infty \theta^m_i\theta^m_j(b_k'y_i)(b_k'y_j)(b_\ell'\tilde{y}_i)(b_{\ell}'\tilde{y}_j)\cr
&\equiv \sum_{m=1}^\infty X_m, \qquad \mbox{where}\quad X_m \equiv K^2\delta_k\delta_\ell \sum_{i<j}\theta_i^m\theta_j^m(b_k'y_i)(b_k'y_j)(b_\ell'\tilde{y}_i)(b_{\ell}'\tilde{y}_j). 
\end{align*}
Since $|X_m|\leq C\|M\|^2 \|\theta\|_{m}^{2m}\leq C\|M\|\|\theta\|_1^2\theta_{\max}^{2(m-1)}$, where $\sum_{m=1}^\infty\theta_{\max}^{2(m-1)}<\infty$, the random variable $\sum_{m=1}^\infty X_m$ is always well-defined. For $m\geq 1$, let $a_m=\theta_{\max}^{2(m-1)}(1-\theta_{\max}^2)$. Then, $\sum_{m=1}^\infty a_m=1$. By Jenson's inequality,
\[
\exp\Bigl(\sum_{m=1}^\infty X_m\Bigr)= \exp\Bigl( \sum_{m=1}^\infty a_m\cdot a_m^{-1} |X_m| \Bigr) 
\leq  \sum_{m=1}^\infty a_m\cdot \exp(a_m^{-1} X_m). 
\]
Using Fatou's lemma, we have
\beq \label{LF-new-4}
\mathbb{E}[\exp(S_{k\ell})] \leq \sum_{m=1}^\infty a_m\cdot \mathbb{E}\bigl[ \exp(a_m^{-1}X_m) \bigr]. 
\eeq
By definition of $X_m$, 
\begin{align*}
X_m &= K^2\delta_k\delta_\ell \biggl\{ \Big[\sum_i\theta_i^m(b_k'y_i)(b_\ell'\tilde{y}_i)\Bigr]^2 - \sum_i \theta_i^{2m}(b_k'y_i)^2(b_\ell'\tilde{y}_i)^2\biggr\}. 
\end{align*}
Note that $\max_{i}\{\|y_i\|,\|\tilde{y}_i\|\}\leq \sqrt{K}$ and $\max_k|\delta_k|=\|M\|$. Therefore,
\[
|X_m|\leq K^2\|M\|^2\Big[\sum_i\theta_i^m(b_k'y_i)(b_\ell'\tilde{y}_i)\Bigr]^2 + K^4\|M\|^2\|\theta\|_{2m}^{2m}. 
\] 
Write $Y=\sum_i \theta_i^m(b_k'y_i)(b_\ell'\tilde{y}_i)$. We see that $Y$ is sum of independent, mean-zero random variables. Since $|(b_k'y_i)(b_\ell'\tilde{y}_i)|\leq K$, by  Hoeffding's inequality, 
\[
\mathbb{P}(|Y|>t)\leq 2\exp\Bigl(-\frac{t^2}{4K^2\|\theta\|_{2m}^{2m}}\Bigr), \qquad \mbox{for any $t>0$}. 
\]
Since $\|\theta\|_{2m}^{2m}\leq \|\theta\|^2\theta_{\max}^{2(m-1)}\leq 2a_m\|\theta\|^2$, we have $a_m^{-1}K^4\|M\|^2\|\theta\|_{2m}^{2m}\leq 2 K^4\|M\|^2\|\theta\|^2$. 
Note that $K^4\|M\|^2\|\theta\|^2\leq 1/9$. 
By direct calculations,
\begin{align*} 
\mathbb{E}\bigl[\exp(a_m^{-1}|X_m|)\bigr]&\leq e^{a_m^{-1}K^4\|M\|^2\|\theta\|_{2m}^{2m}}\cdot\mathbb{E}\bigl[e^{a_m^{-1}K^2\|M\|^2Y^2}\bigr]\cr
&\leq e^{2K^4\|M\|^2\|\theta\|^2} \cdot \mathbb{E}\bigl[e^{a_m^{-1}K^2\|M\|^2Y^2}\bigr]\cr
&= e^{2K^4\|M\|^2\|\theta\|^2} \Bigl[ 1+ \int_0^{\infty}e^{t}\cdot\mathbb{P}\bigl(a_m^{-1}K^2\|M\|^2Y^2>t\bigr)dt\Bigr]\cr
&\leq e^{2K^4\|M\|^2\|\theta\|^2} \Bigl[ 1+  \int_0^{\infty} e^t\cdot e^{-\frac{t}{8 K^4\|M\|^2\|\theta\|^2}}dt\Bigr]\cr
&\leq e^{K^4\|M\|^2\|\theta\|^2}\cdot(1+72 K^4\|M\|^2\|\theta\|^2). 
\end{align*}
We plug it into \eqref{LF-new-4} and notice that $\sum_{m=1}^\infty a_m=1$. It gives
\beq \label{LF-new-7}
\mathbb{E}[\exp(S_{k\ell})]\leq e^{K^4\|M\|^2\|\theta\|^2}\cdot(1+72 K^4\|M\|^2\|\theta\|^2). 
\eeq

Combining \eqref{LF-new-3} and \eqref{LF-new-7} gives
\[
\int \biggl[\frac{dP_1^{(n)}}{dP_0^{(n)}}\biggr]^2dP_0^{(n)} \leq e^{K^4\|M\|^2\|\theta\|^2}\cdot(1+72 K^4\|M\|^2\|\theta\|^2). 
\]
We recall that $\|\theta\|\cdot \|M\|\leq C\|\theta\|\cdot |\mu_2(P)|\to 0$. Hence, the right hand side is $1+o(1)$.  This proves \eqref{LB0}.

\subsection{Proof of Lemmas~\ref{lem:eigP}-\ref{lem:Mnorm}} \label{subsec:Proof-of-Lemma-M}
\subsubsection{Proof of Lemma~\ref{lem:eigP}}
The first claim follows by our assumptions on $P$, so we omit the proof. 
Consider the second claim. Recall that $G = \|\theta\|^{-2} \Pi' \Theta^2 \Pi$ and 
 $d_1, d_2, \ldots, d_K$ are the eigenvalues of 
$G^{1/2} P G^{1/2}$, arranged in the descending order in magnitude. 
By Lemmas D.1 and D.2, $\lambda_k = \|\theta\|^2 d_k$,  $1 \leq k \leq K$, and $d_1 \asymp 1$. 
Combining these,  it suffices to show 
\[
|\mu_2| \asymp |d_2|. 
\] 
 
We now prove for the cases where $P$ is non-singular and singular, separately. Consider the first case.  
Since $1/d_k$ and $1/\mu_K$ are the largest eigenvalue of $G^{-1/2} P^{-1/2} G^{-1/2}$ and $P^{-1}$  in magnitude, respectively,  and $\|G\| \leq C$ and $\|G^{-1}\| \leq C$, it is seen that 
$|\mu_K| \asymp |d_K|$. To show the claim, it sufficient to show that for any $m \geq 2$, if $|\mu_k| \asymp |d_k|$ for $k = m+1, \ldots, K$, then $|\mu_m| \asymp |d_m|$. 

We now fix $m \geq 2$, and assume  $|\mu_k| \asymp |d_k|$ for $k = m+1, \ldots, K$. The goal is to show 
$|\mu_m| \asymp |d_m|$. 
By symmetry, it is sufficient to show that 
\begin{equation} \label{eigPadd} 
|d_m| \leq C |\mu_m|. 
\end{equation} 
Let $P = V \diag(d_1, d_2,\ldots, d_K)  V'$ be the SVD of $P$, where $V \in \mathbb{R}^{K, K}$ is orthonormal, and let $V_m$ be the sub-matrix of $V$ consisting the first $m$ columns of $V$. Introduce 
\[
\widetilde{P}_m = V_m D_m V_m', \qquad \mbox{where} \; D_m = \diag(d_1, d_2, \ldots, d_m).  
\] 
Let $\mu_1^*, \mu_2^*, \ldots, \mu_m^*$ and $d_1^*, d_2^*, \ldots, d_m^*$ 
be the first $m$ eigenvalues of $\widetilde{P}_m$ and 
  $G^{1/2} P_m G^{1/2}$, respectively, arranged in the descending order in magnitude. 
Since $\|G\| \leq C$, we have 
\[
\|P - P_m\| \leq C  |\mu_{m+1}|, \qquad  \|G^{1/2} (P - P_m) G^{1/2}\| \leq C |\mu_{m+1}|. 
\] 
By Theorem \cite[Theorem A.46]{RMTbookBai},   
\begin{equation} \label{eigPadd1A} 
|\mu_m - \mu_m^*| \leq C \|P - P_m\| \leq  |\lambda_{m+1}|, 
\end{equation}
and 
\begin{equation} \label{eigPadd1B}
|d_m - d_m^*| \leq  \|G^{1/2} (P-P_m) G^{1/2}\| \leq  C |\mu_{m+1}|. 
\end{equation} 
and 
At the same time,  note that the nonzero eigenvalues of $G^{1/2} P_m G^{1/2}$ are the same as the nonzero eigenvalues of $ D_m V_m'GV_m$, and also the same as those of $(V_m'GV_m)^{1/2} D_m (V_m'GV_m)^{1/2}$. Since  
$\|G\| \leq C$ and $\|G^{-1}\| \leq C$,  it is seen $\|V_m' G V_m\| \leq C$ and $\|V_m' G V_m)^{-1}\| \leq C$. 
Therefore, by similar arguments, 
\begin{equation} \label{eigPadd2} 
|\mu_m^*| \asymp |d_m^*|. 
\end{equation} 
Combining (\ref{eigPadd1A}), (\ref{eigPadd1B}), and (\ref{eigPadd2}) gives 
\begin{align*} 
&  |\mu_m| \leq |\mu_m^*| + |\mu_m - \mu_m^*|   \leq C(|d_m^*|+ |d_{m+1}|)  \\
\leq & C[(|d_m|  + |d_m - d_m^*|)+ |d_{m+1}|]  \leq C |d_m|. 
\end{align*} 
This proves (\ref{eigPadd}) and the claim follows.

We now consider the case where $P$ is singular, say, $rank(P) = r  < K$, and the nonzero eigenvalues are 
$\mu_1, \mu_2, \ldots, \mu_r$. Let $P = U D U'$ be the SVD, where $U \in \mathbb{R}^{n,r}$ and 
$D = \diag(\mu_1, \mu_2, \ldots, \mu_r)$.  By similar argument, the nonzero eigenvalues of  
 $G^{1/2} P G^{1/2}$ are the same as $(U' G U)^{1/2} D (U' G U)^{1/2}$, where 
 $\|U' G U\| \leq C$ and $\|(U' G U)^{-1}\| \leq C$. The remaining part of the proof is similar so is omitted.

Consider the last claim. Let $\widetilde{P} = \eta \eta'$, where $\eta$ is the first eigenvector of $P$, scaled to have a $\ell^2$-norm of $\sqrt{\mu_1}$.   Write 
\begin{equation} \label{eigP3A}
|P_{ij} - 1| = |P_{ij} - \eta_i \eta_j| + |\eta_i \eta_j - 1|. 
\end{equation} 
Now, first, by definitions and elementary algebra, for $1 \leq i, j \leq K$, 
\begin{equation} \label{eigP3B}
|P_{ij}  - \eta_i \eta_j| \leq |P_{ij} - \widetilde{P}_{ij}| \leq \|P - \widetilde{P}\|  \leq \mu_2, 
\end{equation} 
where by the second claim, $\mu_2 = o(1)$. Note that for $1 \leq i, j \leq K$, $P_{ii} = 1$ and $P_{ij} \geq 0$. 
It is seen that $|\eta_i| = 1 + o(1)$ and all $\eta_i$ must have the positive sign. 
It follows $|\eta_i - 1| =  (1 + \eta_i)^{-1} (1 - \eta_i^2) \leq \mu_2$, and so 
\begin{equation} \label{eigP3C}
|1 - \eta_i \eta_j| \leq |(1 - \eta_i) (1 - \eta_j)| + |1 - \eta_i| + |1 - \eta_j| \leq C \mu_2.   
\end{equation}
Combining (\ref{eigP3A})-(\ref{eigP3C}) gives the claim.  \qed 

\subsubsection{Proof of Lemma~\ref{lem:Regu}}
Consider the first claim about $\sum_i\theta_i\pi_i(k)$. 
Write $X=\sum_{i=1}^n\theta_i(\pi_i(k)-h_k)$. It is seen that $X$ is sum of independent mean-zero random variables, where $\theta_i |\pi_i(k)-h_k|\leq C\theta_{\max}$ and $\sum_{i=1}^n\mathrm{Var}(\theta_i(\pi_i(k)-h_k))\leq C\|\theta\|^2$. By Bernstein's inequality, for any $t>0$,
\[
\mathbb{P}(|X|>t)\leq \exp\Bigl(- \frac{t^2}{C\|\theta\|^2+C\theta_{\max}t}\Bigr). 
\]
It follows that, with probability $1-\|\theta\|_1^{-1}$, 
\[
\Bigl| \sum_i\theta_i\pi_i(k) - h_k\|\theta\|_1\Bigr|=|X|\leq C\|\theta\|\sqrt{\log(\|\theta\|_1)}+C\theta_{\max}\log(\|\theta\|_1). 
\]
Since $\|\theta\|\to\infty$, $\theta_{\max}\to 0$, and $(\|\theta\|^2/\|\theta\|_1)\sqrt{\log(\|\theta\|_1)}\to 0$, the right hand side is $o(\|\theta\|_1)$. Combining it with the assumption of $\min_k\{h_k\}\geq C$, we have
\[
\sum_i\theta_i\pi_i(k)\geq C\|\theta\|_1, \qquad\mbox{with probability }1-\|\theta\|^{-1}=1-o(1). 
\]
Additionally, since $\pi_i(k)\leq 1$, $\sum_i\theta_i\pi_i(k)\leq \|\theta\|_1$. Therefore, with probability $1-o(1)$, each $\sum_{i}\theta_i\pi_i(k)$ is at the order of $\|\theta\|_1$. This proves the first claim. 

Consider the second claim about $G$. Let $y_i=\pi_i-h$. We can write 
\[
 \|\theta\|^2 G =\sum_{i=1}^n \theta_i^2\pi_i\pi_i'=\|\theta\|^2(hh') + \underbrace{\sum_{i=1}^n \theta_i^2y_iy_i'}_{\equiv Y}+\underbrace{\sum_{i=1}^n\theta_i^2hy_i'}_{Z_1} + \underbrace{\sum_{i=1}^n\theta_i^2y_ih'}_{\equiv Z_2}. 
\] 
Note that $\mathbb{E}[y_iy_i']=\Sigma$. Then, $Y-\|\theta\|^2\Sigma=\sum_{i}\theta_i^2(y_iy_i'-\Sigma)$ is sum of independent, mean-zero random matrices, where $\theta_i^2\|y_iy_i'-\Sigma\|\leq C\theta_i^2$. Using the matrix Hoeffding inequality \cite{tropp2012user}, $
\mathbb{P}\bigl(\|Y-\|\theta\|^2\Sigma\|>t\bigr)\leq \exp\bigl(-\frac{t^2}{C\|\theta\|_4^4}\bigr)$, for any $t>0$.  
With $t=\|\theta\|^{-1}$, we have
\[
\bigl\|Y-\|\theta\|^2\Sigma\bigr\|\leq C\|\theta\|_4^2\sqrt{\log(\|\theta\|)}, \qquad \mbox{with probability }1-\|\theta\|^{-1}. 
\]
Similarly, we can apply matrix Hoeffding inequality to $Z_1$ and $Z_2$. It gives
\[
\|Z_1+Z_2\|\leq C\|\theta\|_4^2\sqrt{\log(\|\theta\|)}, \qquad \mbox{with probability }1-\|\theta\|^{-1}. 
\]
Since $\|\theta\|_4^2\leq \theta_{\max}\|\theta\|\ll\|\theta\|^2$, it follows that, with probability $1-o(1)$,
\[
\bigl\|Y+Z_1+Z_2-\|\theta\|^2\Sigma\bigr\|=o(\|\theta\|^2).
\]
At the same time, $\lambda_{\min}(\|\theta\|^2\Sigma)=\|\theta\|^2\cdot \|\Sigma^{-1}\|^{-1}\geq C\|\theta\|^2$. We thus have, with probability $1-o(1)$, 
\[
\lambda_{\min}(\|\theta\|^2G)\geq \lambda_{\min}(Y+Z_1+Z_2)\geq \lambda_{\min}(\|\theta\|^2\Sigma)-\bigl\|Y+Z_1+Z_2-\|\theta\|^2\Sigma\bigr\|\geq C\|\theta\|^2. 
\]
This guarantees $\|G^{-1}\|\leq C$. \qed

\subsubsection{Proof of Lemma~\ref{lem:Mnorm}}
Let $Q=P-1_K1_K'$, and introduce $d\in\mathbb{R}^K$ such that $D=\mathrm{diag}(d)$. By Lemma~\ref{lem:eigP}, $\|Q\|\leq C|\mu_2|$. With these notations, 
\beq \label{lem-M-1}
DPD - 1_K1_K' = dd' + DQD - 1_K1_K'. 
\eeq
Using the same notations, the assumption $DPD\widetilde{h}_D=1_K$ can be written as $D(1_K1_K'+Q)D\widetilde{h}_D =1_K$. It implies
\beq \label{lem-M-2}
1_K = (d'\widetilde{h}_D)d + DQD\widetilde{h}_D. 
\eeq 
We multiply $\widetilde{h}_D'$ on both sides and notice that $1_K'\widetilde{h}_D=1$. It gives
\beq \label{lem-M-3}
(d'\widetilde{h}_D)^2 = 1 - \widetilde{h}_D'DQD\widetilde{h}_D. 
\eeq
Combining \eqref{lem-M-2}-\eqref{lem-M-3} gives
\begin{align*}
dd' - 1_K 1_K' &= [1-(d'\widetilde{h}_D)^2] dd' - (d'\widetilde{h}_D) ( DQD\widetilde{h}_D d+d\widetilde{h}_D DQD) - DQD\widetilde{h}_D\widetilde{h}_D'DQD\cr
&= (\widetilde{h}_D'DQD\widetilde{h}_D)\cdot dd' - (d'\widetilde{h}_D) ( DQD\widetilde{h}_D d+d\widetilde{h}_D DQD) - DQD\widetilde{h}_D\widetilde{h}_D'DQD. 
\end{align*}
Since $\|\widetilde{h}_D\|\leq C$ and $\|d\|\leq C$, we immediately have
\[
\|dd' - 1_K 1_K'\|\leq C\|Q\| \leq C|\mu_2|. 
\]
Plugging it into \eqref{lem-M-1} gives
\[
\|DPD - 1_K1_K'\| \leq C\|Q\|\leq C|\mu_2|. 
\]
\qed

 \newpage
\section{Properties of Signed Polygon statistics}  \label{supp:Var}
We prove Tables~\ref{tab:IdealSgnTsum}-\ref{tab:ProxySgnQsum} and Theorem~\ref{thm:IdealSgnT}-\ref{thm:RealSgnQ}. The analysis of $T_n$ and $Q_n$ is very similar. To save space, we only present the proof for results of $Q_n$. The proof for results of $T_n$ (Tables~\ref{tab:IdealSgnTsum}, \ref{tab:ProxySgnTsum}, and Theorems~\ref{thm:IdealSgnT}, \ref{thm:ProxySgnT}, \ref{thm:RealSgnT}) is omitted. 


We recall the following notations:
\begin{align*}
& \widetilde{\Omega}=\Omega-(\eta^*)(\eta^*)', \qquad \mbox{where}\quad \eta^* = \frac{1}{\sqrt{v_0}}\Omega {\bf 1}_n, \;\; v_0= {\bf 1}_n'\Omega {\bf 1}_n;\cr
& \delta_{ij} = \eta_i(\eta_j-\teta_j)+\eta_j(\eta_i-\teta_i), \qquad \mbox{where}\quad \eta=\frac{1}{\sqrt{v}}(\mathbb{E}A){\bf 1}_n,\;\; \teta = \frac{1}{\sqrt{v}}A{\bf 1}_n,\;\; v= {\bf 1}_n'(\mathbb{E}A){\bf 1}_n;\cr
&r_{ij} = (\eta_i^*\eta_j^*-\eta_i\eta_j) - (\eta_i-\teta_i)(\eta_j-\teta_j)+(1-\frac{v}{V})\teta_i\teta_j, \qquad \mbox{where}\;\; V = {\bf 1}_n'A{\bf 1}_n. 
\end{align*}
Then, the Ideal SgnQ statistic equals to 
\[
\widetilde{Q}_n = \sum_{i,j,k,\ell (dist)} (\widetilde{\Omega}_{ij}  + W_{ij}) 
(\widetilde{\Omega}_{jk}  + W_{jk}) (\widetilde{\Omega}_{k\ell}  + W_{k\ell}) (\widetilde{\Omega}_{\ell i} + W_{\ell i}),
\]
the Proxy SgnQ statistic equals to 
\[
Q_n^*= \sum_{i,j,k,\ell (dist)} (\widetilde{\Omega}_{ij}  + W_{ij} + \delta_{ij})(\widetilde{\Omega}_{jk}  + W_{jk} + \delta_{jk}) (\widetilde{\Omega}_{k\ell}  + W_{k\ell} + \delta_{k\ell}) (\widetilde{\Omega}_{\ell i} + W_{\ell i} + \delta_{\ell i}),
\]
and the SgnQ statistic equals to 
\[
Q_n = \sum_{i,j,k,\ell (dist)} (\widetilde{\Omega}_{ij}  + W_{ij}+\delta_{ij}+r_{ij}) 
(\widetilde{\Omega}_{jk}  + W_{jk}+\delta_{jk}+r_{jk}) (\widetilde{\Omega}_{k\ell}  + W_{k\ell }+\delta_{k\ell}+r_{k\ell}) (\widetilde{\Omega}_{\ell i} + W_{\ell i}+\delta_{\ell i}+r_{\ell i}). 
\]

As explained in Section~\ref{sec:SgnedPolygon}, each of $\widetilde{Q}_n, Q_n^*, Q_n$ is the sum of a finite number of post-expansion sums, each having the form
\begin{equation} \label{post-expansion-sum} 
\sum_{i,j,k,\ell(dist)} a_{ij} b_{jk} c_{k\ell} d_{\ell i},     
\end{equation} 
where $a_{ij}$ equals to one of $\{\widetilde{\Omega}_{ij},W_{ij}, \delta_{ij}, r_{ij}\}$; same for $b_{ij}$, $c_{ij}$ and $d_{ij}$. Let $N_{\widetilde{\Omega}}$ be the (common) number of $\widetilde{\Omega}$ terms in each product; similarly, we define $N_W, N_\delta, N_r$. These numbers satisfy $N_{\widetilde{\Omega}}+N_W + N_\delta + N_r=4$. For example, for the post-expansion sum $\sum_{i,j,k,\ell(dist)}\widetilde{\Omega}_{ij}W_{jk}W_{k\ell}W_{\ell i}$, $(N_{\widetilde{\Omega}}, N_W, N_\delta, N_r)=(1,3,0,0)$. In Section~\ref{suppC1}, we study $\widetilde{Q}_n$, and it involves these post-expansion sums such that 
\[
N_\delta=N_r=0,
\]
In Section~\ref{suppC2}, we study $(Q_n^*- \widetilde{Q}_n)$, which involves post-expansion sums such that 
\[
N_\delta>0, \;\; \mbox{and}\;\; N_r=0,
\]
In Section~\ref{suppC3}, we study $(Q_n-Q_n^*)$, which is related to the sums such that 
\[
N_r>0.
\]

\subsection{Analysis of Table~\ref{tab:IdealSgnQsums}, proof of Theorem~\ref{thm:IdealSgnQ}} \label{suppC1}
Define
\begin{align*}
& X_1 = \sum_{i, j, k, \ell (dist)} W_{ij} W_{jk} W_{k \ell} W_{\ell i}, \qquad X_2 =  \sum_{i, j, k, \ell (dist)} \widetilde{\Omega}_{ij} W_{jk} W_{k \ell} W_{\ell i},\cr
& X_3 = \sum_{i, j, k, \ell (dist)} \widetilde{\Omega}_{ij}  \widetilde{\Omega}_{jk} W_{k \ell} W_{\ell i}, \qquad\; X_4 = \sum_{i, j, k, \ell (dist)} \widetilde{\Omega}_{ij} W_{jk} \widetilde{\Omega}_{k \ell} W_{\ell i},\cr
& X_5 = \sum_{i, j, k, \ell (dist)} \widetilde{\Omega}_{ij} \widetilde{\Omega}_{jk} \widetilde{\Omega}_{k \ell} W_{\ell i}, \qquad\;\, X_6 = \sum_{i, j, k, \ell (dist)} \widetilde{\Omega}_{ij} \widetilde{\Omega}_{jk} \widetilde{\Omega}_{k \ell} \widetilde{\Omega}_{\ell i}.  
\end{align*}
We first consider the null hypothesis. Since $\widetilde{\Omega}$ is a zero matrix, it is not hard to see that
\[
\widetilde{Q}_n = X_1. 
\]
The following lemmas are proved in Section~\ref{subsec:Calculations}. 
\begin{lemma} \label{lem:IdealSgnQ-null}
Suppose the conditions of Theorem~\ref{thm:IdealSgnQ} hold. Under the null hypothesis,  as $n \goto \infty$, $\mathbb{E}[\widetilde{Q}_n] =0$ and $\mathrm{Var}(\widetilde{Q}_n) = 8\|\theta\|^8\cdot [1+o(1)]$.
\end{lemma}

\begin{lemma} \label{lem:Normality}
Suppose the conditions of Theorem~\ref{thm:IdealSgnQ} hold. Under the null hypothesis,  as $n \goto \infty$, 
\[
\frac{\widetilde{Q}_n - E[\widetilde{Q}_n]}{\sqrt{\mathrm{Var}(\widetilde{Q}_n)}}  \;\;   \longrightarrow \;\;  N(0,1), \qquad \mbox{in law}. 
\] 
\end{lemma}

We then consider the alternative hypothesis. By elementary algebra,
\[
\widetilde{Q}_n = X_1 + 4X_2 + 4X_3 + 2X_4 + 4X_5 + X_6. 
\]
The following lemma characterizes the asymptotic mean and variance of $X_1$-$X_6$ under the alternative hypothesis. It gives rise to Columns 5-6 of Table~\ref{tab:IdealSgnQsums}. 
\begin{lemma}[Table~\ref{tab:IdealSgnQsums}] \label{lem:IdealSgnQ-alt}
Suppose conditions of Theorem~\ref{thm:IdealSgnQ} hold. Write $\alpha=|\lambda_2|/\lambda_1$. Under the alternative hypothesis,  as $n \goto \infty$,
\begin{itemize} \itemsep +5pt
\item $\mathbb{E}[X_k]=0$ for $1\leq k\leq 5$, and $\mathbb{E}[X_6]=\tr(\widetilde{\Omega}^4)\cdot[1+o(1)]$. 
\item $C^{-1}\|\theta\|^8\leq \mathrm{Var}(X_1)\leq C\|\theta\|^8$. 
\item $\mathrm{Var}(X_2)\leq C\alpha^2 \|\theta\|^4\|\theta\|_3^6=o(\|\theta\|^8)$. 
\item $\mathrm{Var}(X_3)\leq C\alpha^4 \|\theta\|^6\|\theta\|_3^6=o( \alpha^6\|\theta\|^8  \|\theta\|_3^6 )$. 
\item $\mathrm{Var}(X_4)\leq C\alpha^4\|\theta\|_3^{12}=o(\|\theta\|^8)$. 
\item $\mathrm{Var}(X_5)\leq C \alpha^6\|\theta\|^8  \|\theta\|_3^6$. 
\end{itemize}  
As a result, $\mathbb{E}[\widetilde{Q}_n]\sim \tr(\widetilde{\Omega}^4)$ and $\mathrm{Var}(\widetilde{Q}_n)\leq C(\|\theta\|^8 + \alpha^6\|\theta\|^8  \|\theta\|_3^6)$. 
\end{lemma}

Theorem~\ref{thm:IdealSgnQ} follows directly from Lemmas~\ref{lem:IdealSgnQ-null}-\ref{lem:IdealSgnQ-alt}.

\subsection{Analysis of Table~\ref{tab:ProxySgnQsum}, proof of Theorem~\ref{thm:ProxySgnQ}} \label{suppC2} 
We introduce $U_a$, $U_b$ and $U_c$ such that 
\[
Q_n^* - \widetilde{Q}_n = U_a + U_b + U_c,
\]
where $U_a$, $U_b$ and $U_c$ contain post-expansion sums \eqref{post-expansion-sum} with $N_\delta=1$, $N_\delta=2$, and $N_\delta\geq 3$, respectively. 

First, we consider the post-expansion sums with $N_\delta=1$. Define
\beq \label{def-Ua}
U_a = 4Y_1+8Y_2 + 4Y_3 + 8Y_4 + 4Y_5 + 4Y_6, 
\eeq
where 
\begin{align*}
& Y_1 = \sum_{i, j, k, \ell (dist)} \delta_{ij} W_{jk} W_{k \ell} W_{\ell i}, \qquad\;\;Y_2 =  \sum_{i, j, k, \ell (dist)} \delta_{ij} \widetilde{\Omega}_{jk}  W_{k \ell}  W_{\ell i},\cr
& Y_3 = \sum_{i, j, k, \ell (dist)} \delta_{ij} W_{jk}   \widetilde{\Omega}_{k \ell}  W_{\ell i}, \qquad\;\;\; Y_4 = \sum_{i, j, k, \ell (dist)} \delta_{ij}  \widetilde{\Omega}_{jk} \widetilde{\Omega}_{k \ell}W_{\ell i},\cr
& Y_5 = \sum_{i, j, k, \ell (dist)} \delta_{ij}  \widetilde{\Omega}_{jk} W_{k \ell}  \widetilde{\Omega}_{\ell i}, \qquad\quad Y_6 = \sum_{i, j, k, \ell (dist)} \delta_{ij}  \widetilde{\Omega}_{jk} \widetilde{\Omega}_{k \ell}  \widetilde{\Omega}_{\ell i}.  
\end{align*}
Under the null hypothesis, only $Y_1$ is nonzero, and 
\[
U_a=4Y_1.
\] 
\begin{lemma} \label{lem:ProxySgnQ(a)-null}
Suppose the conditions of Theorem~\ref{thm:IdealSgnQ} hold. Under the null hypothesis,  as $n \goto \infty$, $\mathbb{E}[U_a] =0$ and $\mathrm{Var}(U_a) \leq C\|\theta\|^2\|\theta\|_3^6=o(\|\theta\|^8)$.
\end{lemma}

Under the alternative hypothesis, the following lemma characterizes the asymptotic means and variances of $Y_1$-$Y_6$. It gives rise to Rows 1-6 of Table~\ref{tab:ProxySgnQsum} and is proved in Section~\ref{subsec:Calculations}. 
\begin{lemma}[Table~\ref{tab:ProxySgnQsum}, Rows~1-6] \label{lem:ProxySgnQ(a)-alt}
Suppose the conditions of Theorem~\ref{thm:IdealSgnQ} hold. Let $\alpha=|\lambda_2|/\lambda_1$. Under the alternative hypothesis,  as $n \goto \infty$,
\begin{itemize} \itemsep +5pt
\item $\mathbb{E}[Y_k]=0$ for $k\in\{1,2,3,5,6\}$, and $|\mathbb{E}[Y_4]|\leq C\alpha^2\|\theta\|^6=o(\alpha^4\|\theta\|^8)$. 
\item $\mathrm{Var}(Y_1)\leq C\|\theta\|^2\|\theta\|_3^6=o(\|\theta\|^8)$. 
\item $\mathrm{Var}(Y_2)\leq C\alpha^2\|\theta\|^4\|\theta\|_3^6=o(\|\theta\|^8)$. 
\item $\mathrm{Var}(Y_3)\leq C\alpha^2\|\theta\|^4\|\theta\|_3^6=o(\|\theta\|^8)$. 
\item $\mathrm{Var}(Y_4)\leq \frac{C\alpha^4\|\theta\|^{10}\|\theta\|_3^3}{\|\theta\|_1}=o(\alpha^6\|\theta\|^8\|\theta\|_3^6)$. 
\item $\mathrm{Var}(Y_5)\leq \frac{C\alpha^4\|\theta\|^4\|\theta\|_3^9}{\|\theta\|_1}=o(\|\theta\|^8)$.
\item $\mathrm{Var}(Y_6)\leq \frac{C\alpha^6\|\theta\|^{12}\|\theta\|_3^3}{\|\theta\|_1}=O(\alpha^6\|\theta\|^8\|\theta\|_3^6)$. 
\end{itemize}  
As a result, $\mathbb{E}[U_a]=o(\alpha^4 \|\theta\|^8)$ and $\mathrm{Var}(U_a)\leq C\alpha^6\|\theta\|^8\|\theta\|_3^6 +o(\|\theta\|^8)$.
\end{lemma}

Next, we consider the post-expansion sums with $N_\delta=2$. Define
\beq  \label{def-Ub}
U_b = 4Z_1+2Z_2 + 8Z_3 + 4Z_4 + 4Z_5 + 2Z_6, 
\eeq
where 
\begin{align*}
& Z_1 = \sum_{i, j, k, \ell (dist)} \delta_{ij} \delta_{jk} W_{k \ell} W_{\ell i}, \qquad\;\;Z_2 =  \sum_{i, j, k, \ell (dist)} \delta_{ij} W_{jk}  \delta_{k \ell}  W_{\ell i},\cr
& Z_3 = \sum_{i, j, k, \ell (dist)} \delta_{ij} \delta_{jk}   \widetilde{\Omega}_{k \ell}  W_{\ell i}, \qquad\;\;\; Z_4 = \sum_{i, j, k, \ell (dist)} \delta_{ij}  \widetilde{\Omega}_{jk} \delta_{k \ell}W_{\ell i},\cr
& Z_5 = \sum_{i, j, k, \ell (dist)} \delta_{ij} \delta_{jk} \widetilde{\Omega}_{k \ell}  \widetilde{\Omega}_{\ell i}, \qquad\quad Z_6 = \sum_{i, j, k, \ell (dist)} \delta_{ij}  \widetilde{\Omega}_{jk} \delta_{k \ell}  \widetilde{\Omega}_{\ell i}.  
\end{align*}
Under the null hypothesis, only $Z_1$ and $Z_2$ are nonzero, and 
\[
U_b = 4Z_1 + 2Z_2. 
\]
\begin{lemma} \label{lem:ProxySgnQ(b)-null}
Suppose the conditions of Theorem~\ref{thm:IdealSgnQ} hold. Under the null hypothesis,  as $n \goto \infty$, 
\begin{itemize} \itemsep +5pt
\item $\mathbb{E}[Z_1] =\|\theta\|^4\cdot[1+o(1)]$, and $\mathrm{Var}(Z_1)\leq C \|\theta\|^2\|\theta\|_3^6=o(\|\theta\|^8)$.
\item $\mathbb{E}[Z_2] =2\|\theta\|^4\cdot[1+o(1)]$, and $\mathrm{Var}(Z_2)\leq \frac{C\|\theta\|^6\|\theta\|_3^3}{\|\theta\|_1}=o(\|\theta\|^8)$.
\end{itemize}
As a result, $\mathbb{E}[U_b]\sim 8\|\theta\|^4$ and $\mathrm{Var}(U_b)=o(\|\theta\|^8)$. 
\end{lemma}

Under the alternative hypothesis, the following lemma provides the asymptotic means and variances of $Z_1$-$Z_6$. It gives rise to Rows 7-12 of Table~\ref{tab:ProxySgnQsum}:
\begin{lemma}[Table~\ref{tab:ProxySgnQsum}, Rows~7-12] \label{lem:ProxySgnQ(b)-alt}
Suppose conditions of Theorem~\ref{thm:IdealSgnQ} hold. Write $\alpha=|\lambda_2|/\lambda_1$. Under the alternative hypothesis,  as $n \goto \infty$,
\begin{itemize}  \itemsep +5pt
\item $|\mathbb{E}[Z_1]|\leq C\|\theta\|^4=o(\alpha^4\|\theta\|^8)$, and $\mathrm{Var}(Z_1)\leq C \|\theta\|^2\|\theta\|_3^6=o(\|\theta\|^8)$. 
\item $|\mathbb{E}[Z_2]|\leq C\|\theta\|^4=o(\alpha^4\|\theta\|^8)$, and $\mathrm{Var}(Z_2)\leq \frac{C\|\theta\|^6\|\theta\|_3^3}{\|\theta\|_1}=o(\|\theta\|^8)$. 
\item $\mathbb{E}Z_3 = 0$, and $\mathrm{Var}(Z_3)\leq C\alpha^2\|\theta\|^4\|\theta\|_3^6=o(\|\theta\|^8)$.  
\item $|\mathbb{E}[Z_4]|\leq C\alpha\|\theta\|^4=o(\alpha^4\|\theta\|^8)$, and $\mathrm{Var}(Z_4)\leq \frac{C\alpha^2\|\theta\|^8\|\theta\|_3^3}{\|\theta\|_1} = o(\|\theta\|^8)$. 
\item $|\mathbb{E}[Z_5]|\leq C\alpha^2\|\theta\|^6=o(\alpha^4\|\theta\|^8)$, and $\mathrm{Var}(Z_5)\leq \frac{C\alpha^4\|\theta\|^{14}}{\|\theta\|_1^2} = o(\alpha^6\|\theta\|^8\|\theta\|_3^6)$. 
\item $|\mathbb{E}[Z_6]|\leq \frac{C\alpha^2\|\theta\|^8}{\|\theta\|_1^2}=o(\alpha^4\|\theta\|^8)$, and $\mathrm{Var}(Z_6)\leq \frac{C\alpha^4\|\theta\|^8\|\theta\|_3^6}{\|\theta\|_1^2}=o(\|\theta\|^8)$. 
\end{itemize}  
As a result, $\mathbb{E}[U_b] =o(\alpha^4\|\theta\|^8)$ and $\mathrm{Var}(U_b)=o(\|\theta\|^8 + \alpha^6\|\theta\|^8\|\theta\|_3^6) $. 
\end{lemma}

Last, we consider the post-expansion sums with $N_\delta\geq 3$. Define
\beq   \label{def-Uc}
U_c = 4T_1 + 4T_2 + F, 
\eeq
where
\begin{align*}
&T_1= \sum_{i,j,k,\ell (dist)} \delta_{ij}\delta_{jk}\delta_{k\ell}W_{\ell i}, \qquad\;\; T_2= \sum_{i,j,k,\ell (dist)} \delta_{ij}\delta_{jk}\delta_{k\ell}\widetilde{\Omega}_{\ell i},\cr
& F = \sum_{i,j,k,\ell (dist)} \delta_{ij}\delta_{jk}\delta_{k\ell}\delta_{\ell i}. 
\end{align*}
Under the null hypothesis, only $T_1$ and $F$ are nonzero, and 
\[
U_b = 4T_1 + F. 
\]
\begin{lemma} \label{lem:ProxySgnQ(c)-null}
Suppose the conditions of Theorem~\ref{thm:IdealSgnQ} hold. Under the null hypothesis,  as $n \goto \infty$, 
\begin{itemize} \itemsep +5pt
\item $\mathbb{E}[T_1] =-2\|\theta\|^4\cdot [1+o(1)]$, and $\mathrm{Var}(T_1)\leq \frac{C\|\theta\|^6 \|\theta\|_3^3}{\|\theta\|_1}=o(\|\theta\|^8)$.
\item $|\mathbb{E}[F]| =2\|\theta\|^4\cdot [1+o(1)]$, and $\mathrm{Var}(F)\leq \frac{C\|\theta\|^{10}}{\|\theta\|_1^2}=o(\|\theta\|^8)$.
\end{itemize}
As a result, $\mathbb{E}[U_c] \sim -6\|\theta\|^4$ and $\mathrm{Var}(U_c)=o(\|\theta\|^8)$. 
\end{lemma}

Under the alternative hypothesis, the next lemma studies the asymptotic means and variances of $T_1$, $T_2$ and $F$. It gives rise to Rows 13-15 of Table~\ref{tab:ProxySgnQsum}:
\begin{lemma}[Table~\ref{tab:ProxySgnQsum}, Rows~13-15] \label{lem:ProxySgnQ(c)-alt}
Suppose conditions of Theorem~\ref{thm:IdealSgnQ} hold. Write $\alpha=|\lambda_2|/\lambda_1$. Under the alternative hypothesis,  as $n \goto \infty$,
\begin{itemize} \itemsep +5pt
\item $|\mathbb{E}[T_1]|\leq C\|\theta\|^4=o(\alpha^4\|\theta\|^8)$, and $\mathrm{Var}(T_1)\leq \frac{C\|\theta\|^6 \|\theta\|_3^3}{\|\theta\|_1}=o(\|\theta\|^8)$. 
\item $|\mathbb{E}[T_2]|\leq \frac{C\alpha\|\theta\|^6}{\|\theta\|_1^3}=o(\alpha^4\|\theta\|^8)$, and $\mathrm{Var}(T_2)\leq \frac{C\alpha^2 \|\theta\|^8 \|\theta\|_3^3}{\|\theta\|_1}=o(\|\theta\|^8)$. 
\item $|\mathbb{E}[F]|\leq C\|\theta\|^4=o(\alpha^4\|\theta\|^8)$, and $\mathrm{Var}(F)\leq \frac{C\|\theta\|^{10}}{\|\theta\|_1^2}=o(\|\theta\|^8)$.  
\end{itemize}  
As a result, $\mathbb{E}|U_c| =o(\alpha^4\|\theta\|^8)$ and $\mathrm{Var}(U_c)=o(\|\theta\|^8)$. 
\end{lemma}

We now prove Theorem~\ref{thm:ProxySgnQ}. Since $Q_n^*-\widetilde{Q}_n=U_a+U_b+U_c$, we have
\begin{align*}
& \mathbb{E}[Q_n^*-\widetilde{Q}_n]=\mathbb{E}[U_a]+\mathbb{E}[U_b]+\mathbb{E}[U_c], \cr
& \mathrm{Var}(Q_n^*-\widetilde{Q}_n)\leq 3\mathrm{Var}(U_a)+3\mathrm{Var}(U_b)+3\mathrm{Var}(U_c). 
\end{align*}
Consider the null hypothesis. By Lemmas~\ref{lem:ProxySgnQ(a)-null}, \ref{lem:ProxySgnQ(b)-null}, \ref{lem:ProxySgnQ(c)-null},  
\[
\mathbb{E}[Q_n^*- \widetilde{Q}_n] = 0 + 8\|\theta\|^4 -6\|\theta\|^4 +o(\|\theta\|^4)\sim 2\|\theta\|^4, 
\]
and 
\[
\mathrm{Var}(Q_n^* - \widetilde{Q}_n)\leq C\|\theta\|^2\|\theta\|_3^6 + \frac{C\|\theta\|^6\|\theta\|_3^3}{\|\theta\|_1} +\frac{C\|\theta\|^{10}}{\|\theta\|_1^2}.  
\]
Using the universal inequality $\|\theta\|^4\leq \|\theta\|_1\|\theta\|_3^3$, we further have
\[
\mathrm{Var}(Q_n^* - \widetilde{Q}_n)\leq C\|\theta\|^2\|\theta\|_3^6 = o(\|\theta\|^8),
\]
where $\|\theta\|_3^3=o(\|\theta\|^2)$ and $\|\theta\|\to\infty$ in our range of interest. This proves claims for the null hypothesis. Consider the alternative hypothesis. By Lemmas~\ref{lem:ProxySgnQ(a)-alt}, \ref{lem:ProxySgnQ(b)-alt}, \ref{lem:ProxySgnQ(c)-alt}, 
\[
\bigl|\mathbb{E}[Q_n^*- \widetilde{Q}_n]\bigr| \leq C\alpha^2\|\theta\|^6, 
\]
where the main contributors are $Y_4$ and $Z_5$. Since $\alpha\|\theta\|\to\infty$ in our range of interest, the above is $o(\alpha^4\|\theta\|^8)$. By Lemmas~\ref{lem:ProxySgnQ(a)-alt}, \ref{lem:ProxySgnQ(b)-alt}, \ref{lem:ProxySgnQ(c)-alt}, 
\[
\mathrm{Var}(Q_n^* - \widetilde{Q}_n)\leq \frac{C\alpha^6\|\theta\|^{12}\|\theta\|_3^3}{\|\theta\|_1},
\]
where the main contributor is $Y_6$. Using the universal inequality of $\|\theta\|^4\leq \|\theta\|_1\|\theta\|_3^3$, the above is $O(\alpha^6\|\theta\|^8\|\theta\|_3^6)$.
This proves claims for the alternative hypothesis.

\subsection{Analysis of $(Q_n-Q_n^*)$, proof of Theorem~\ref{thm:RealSgnQ}} \label{suppC3}
By definition, $(Q_n-Q^*_n)$ expands to the sum of $175$ post-expansion sums, where each has the form \eqref{post-expansion-sum} and satisfies $N_r>0$. Recall that 
\[
r_{ij} = (\eta_i^*\eta_j^*-\eta_i\eta_j) - (\eta_i-\teta_i)(\eta_j-\teta_j)+(1-\frac{v}{V})\teta_i\teta_j. 
\]
Since $\delta_{ij}=\eta_i(\eta_j-\teta_j)+\eta_j(\eta_i-\teta_i)$, we have $\teta_i\teta_j =\eta_i\eta_j-\delta_{ij}+(\teta_i-\eta_i)(\teta_j-\eta_j)$. Inserting it into the definition of $r_{ij}$ gives
\beq
r_{ij} = (\eta_i^*\eta_j^*-\eta_i\eta_j)+(1-\frac{v}{V})\eta_i\eta_j-(1-\frac{v}{V})\delta_{ij} - \frac{v}{V}(\teta_i-\eta_i)(\teta_j-\eta_j).
\eeq
Define
\[
\tilde{r}_{ij} = - \frac{v}{V}(\teta_i-\eta_i)(\teta_j-\eta_j),\qquad \epsilon_{ij} = (\eta_i^*\eta_j^*-\eta_i\eta_j)+(1-\frac{v}{V})\eta_i\eta_j-(1-\frac{v}{V})\delta_{ij}.
\]
Then, we can write 
\beq  \label{r-decompose}
r_{ij}=\tilde{r}_{ij} + \epsilon_{ij}. 
\eeq
Using this notation, we re-write
\[
Q_n=\sum_{i,j,k,\ell (dist)}M_{ij}M_{jk}M_{k\ell}M_{\ell i},\qquad \mbox{where}\;\; M_{ij}=\widetilde{\Omega}_{ij}+W_{ij}+\delta_{ij}+\tilde{r}_{ij}+\epsilon_{ij},
\]
and 
\[
Q_n^* = \sum_{i,j,k,\ell (dist)} M^*_{ij}M^*_{jk}M^*_{k\ell}M^*_{\ell i}, \qquad \mbox{where}\;\; M^*_{ij}\equiv \widetilde{\Omega}_{ij}+W_{ij}+\delta_{ij}.
\]
We then introduce an intermediate variable: 
\beq \label{Define-intermediateSgnQ}
\widetilde{Q}^*_n = \sum_{i,j,k,\ell (dist)}\widetilde{M}^*_{ij}\widetilde{M}^*_{jk}\widetilde{M}^*_{k\ell}\widetilde{M}^*_{\ell i}, \quad \mbox{where}\; \widetilde{M}^*_{ij}= \widetilde{\Omega}_{ij}+W_{ij}+\delta_{ij}+\tilde{r}_{ij}. 
\eeq
As a result, $(Q_n-Q_n^*)$ decomposes into
\beq \label{remainder-decompose}
Q_n-Q_n^* = (\widetilde{Q}^*_n-Q_n^*) + (Q_n - \widetilde{Q}_n^*). 
\eeq

We note that $Q_n$ can be expanded to the sum of $5^4=625$ post-expansion sums, each with the form
\[
\sum_{i,j,k,\ell (dist)} a_{ij}b_{jk}c_{k\ell}d_{\ell i},
\]
where each of $a_{ij}, b_{ij}, c_{ij}, d_{ij}$ takes values in $\{ \widetilde{\Omega}_{ij}, W_{ij}, \delta_{ij}, \tilde{r}_{ij}, \epsilon_{ij} \}$. Let $N_{\widetilde{\Omega}}$ be the (common) number of $\widetilde{\Omega}$ terms in each product and define $N_W, N_\delta, N_{\tilde{r}}, N_{\epsilon}$ similarly. Among the $625$ post-expansion sums, 
\begin{itemize}
\item $3^4=81$ of them are contained in $Q_n^*$, 
\item $4^4-3^4=175$ of them are contained in $(\widetilde{Q}^*_n-Q_n^*)$, 
\item and $5^4-4^4=369$ of them are contained in $(Q_n - \widetilde{Q}_n^*)$.
\end{itemize}
We shall study $(\widetilde{Q}_n^*-Q_n^*)$ and $(Q_n-\widetilde{Q}_n^*)$, separately. 

In our analysis, one challenge is to deal with the random variable $V$ that appears in the denominator in the expression of $r_{ij}$.  
The following lemma is useful and proved in Section~\ref{subsec:Calculations}. 

\begin{lemma} \label{lem:event}
Suppose conditions of Theorem~\ref{thm:RealSgnQ} hold. As $n\to\infty$, for any sequence $x_n$ such that $\sqrt{\log(\|\theta\|_1)}\ll x_n\ll \|\theta\|_1$,  
\[
\mathbb{E}\bigl[(\widetilde{Q}_n-Q_n)^2\cdot I\{|V-v|>\|\theta\|_1x_n\}\bigr] \to 0. 
\]
\end{lemma}

The next two lemmas are proved in Section~\ref{subsec:Calculations}. 
\begin{lemma} \label{lem:remainder1}
Suppose conditions of Theorem~\ref{thm:RealSgnQ} hold. Write $\alpha=|\lambda_2|/\lambda_1$. As $n\to\infty$, 
\begin{itemize} \itemsep +5pt
\item Under the null hypothesis, $|\mathbb{E}[\widetilde{Q}_n^*-Q_n^*]|=o(\|\theta\|^4)$ and $\mathrm{Var}(\widetilde{Q}^*_n-Q_n^*)=o(\|\theta\|^8)$. 
\item Under the alternative hypothesis, $|\mathbb{E}[\widetilde{Q}^*_n-Q_n^*]|=o(\alpha^4\|\theta\|^8)$ and $\mathrm{Var}(\widetilde{Q}_n^*-Q_n^*)=o(\|\theta\|^8+\alpha^6\|\theta\|^8\|\theta\|_3^6)$. 
\end{itemize}
\end{lemma}
\begin{lemma} \label{lem:remainder2}
Suppose conditions of Theorem~\ref{thm:RealSgnQ} hold. Write $\alpha=|\lambda_2|/\lambda_1$. As $n\to\infty$, 
\begin{itemize} \itemsep +5pt
\item Under the null hypothesis, $|\mathbb{E}[Q_n - \widetilde{Q}^*_n]|=o(\|\theta\|^4)$ and $\mathrm{Var}(Q_n - \widetilde{Q}^*_n)=o(\|\theta\|^8)$. 
\item Under the alternative hypothesis, $|\mathbb{E}[Q_n - \widetilde{Q}^*_n]|=o(\alpha^4\|\theta\|^8)$ and $\mathrm{Var}(\widetilde{Q}_n^*-Q_n^*)=O(\|\theta\|^8+\alpha^6\|\theta\|^8\|\theta\|_3^6)$. 
\end{itemize}
\end{lemma}

Theorem~\ref{thm:RealSgnQ} follows directly from \eqref{remainder-decompose} and Lemmas~\ref{lem:remainder1}-\ref{lem:remainder2}.

\subsection{Proof of Lemmas~\ref{lem:IdealSgnQ-null}-\ref{lem:remainder2}} \label{subsec:Calculations}

\subsubsection{Proof of Lemma~\ref{lem:IdealSgnQ-null}}
Under the null hypothesis, 
\[
\widetilde{Q}_n = X_1 = \sum_{i, j, k, \ell (dist)} W_{ij} W_{jk} W_{k \ell} W_{\ell i}. 
\]
For mutually distinct indices $(i,j,k,\ell)$, $(W_{ij}, W_{jk}, W_{k \ell}, W_{\ell i})$ are independent of each other, each with mean zero. So $\mathbb{E}[W_{ij} W_{jk} W_{k \ell} W_{\ell i}]=0$. It follows that
\[
\mathbb{E}[\widetilde{Q}_n]=0.
\] 
We now calculate the variance of $\widetilde{Q}_n$. Under the null hypothesis, $\Omega_{ij}=\theta_i\theta_j$; hence,  $\mathrm{Var}(W_{ij})=\Omega_{ij}(1-\Omega_{ij})=\theta_i\theta_j-\theta_i^2\theta_j^2=\theta_i\theta_j[1+O(\theta^2_{\max})]$. It follows that
\begin{align} \label{LemmaC1-eq1}
\mathrm{Var}(W_{ij} W_{jk} W_{k \ell} W_{\ell i})&=\theta^2_i\theta^2_j\theta^2_k\theta^2_\ell\cdot[1+O(\theta_{\max}^2)]^4\cr
& = \theta^2_i\theta^2_j\theta^2_k\theta^2_\ell\cdot[1+O(\theta_{\max}^2)]. 
\end{align}
Note that each $(i,j,k,\ell)$ corresponds to a $4$-cycle in a complete graph of $n$ nodes. For $(i,j,k,\ell)$ and $(i',j',k',\ell')$, we can write $W_{ij} W_{jk} W_{k \ell} W_{\ell i}\cdot W_{i'j'} W_{j'k'} W_{k' \ell'} W_{\ell' i'}$ in the form of $\prod_{t}(W_{i_tj_t})^{m_t}$, where $\{W_{i_tj_t}\}$ are mutually distinct with each other and $m_t$ is the number of times that $W_{i_tj_t}$ appears in this product. If the two $4$-cycles corresponding to $(i,j,k,\ell)$ and $(i',j',k',\ell')$ are not exactly overlapping, then at least two of $m_t$ equals to $1$. As a result, the mean of $\prod_{t}(W_{i_tj_t})^{m_t}$ is zero. In other words, we have argued that
\beq  \label{LemmaC1-eq2}
\begin{array}{l}
\mathrm{Cov}(W_{ij} W_{jk} W_{k \ell} W_{\ell i},\; W_{i'j'} W_{j'k'} W_{k' \ell'} W_{\ell' i'})=0 \mbox{ if the }\\
\mbox{two cycles corresponding to $(i,j,k,\ell)$ and $(i',j',k',\ell')$}\\
\mbox{are not exactly overlapping}. 
\end{array}
\eeq
In the sum over all distinct $(i,j,k,\ell)$, each $4$-cycle is repeatedly counted by 8 times
\[
\begin{array}{l}
(i, j,k,\ell), (j,k,\ell,i), (k,\ell,i,j),(\ell, i,j,k),\\
(\ell, k, j, i), (k, j,i, \ell), (j, i, \ell,k), (i, \ell, k, j).
\end{array}
\]
It follows that
\begin{align}  \label{LemmaC1-eq3}
\mathrm{Var}(\widetilde{Q}_n)&=\mathrm{Var}\biggl( 8 \sum_{\substack{\text{unique}\\\text{4-cycles}}} W_{ij}W_{jk}W_{k\ell}W_{\ell i}\biggr)\cr
&= 64\cdot \mathrm{Var}\biggl(\sum_{\substack{\text{unique}\\\text{4-cycles}}}W_{ij}W_{jk}W_{k\ell}W_{\ell i}\biggr)\cr
&= 64\sum_{\substack{\text{unique}\\\text{4-cycles}}} \mathrm{Var}\bigl( W_{ij}W_{jk}W_{k\ell}W_{\ell i} \bigr)\cr
&= 8 \sum_{i,j,k,\ell (dist)} \mathrm{Var}\bigl( W_{ij}W_{jk}W_{k\ell}W_{\ell i} \bigr)\cr
&= [1+O(\theta_{\max}^2)]\cdot 8 \sum_{i,j,k,\ell (dist)}\theta^2_i\theta^2_j\theta^2_k\theta^2_\ell, 
\end{align}
where the third line is from \eqref{LemmaC1-eq2} and the last line is from \eqref{LemmaC1-eq1}. 
We then compute the right hand side of \eqref{LemmaC1-eq3}. Note that
\[
\sum_{i,j,k,\ell (dist)}\theta^2_i\theta^2_j\theta^2_k\theta^2_\ell =\|\theta\|^8 - \sum_{i,j,k,\ell (not\; dist)}\theta^2_i\theta^2_j\theta^2_k\theta^2_\ell,
\]
where
\[
\sum_{i,j,k,\ell (not\; dist)}\theta^2_i\theta^2_j\theta^2_k\theta^2_\ell \leq {4\choose 2}\sum_{i,j,k}\theta^2_i\theta^2_j\theta^4_k\leq C\|\theta\|^4\|\theta\|_4^4=\|\theta\|^8\cdot O\Bigl(\frac{\|\theta\|_4^4}{\|\theta\|^4}\Bigr). 
\]
Combining the above gives
\beq   \label{LemmaC1-eq4}
\sum_{i,j,k,\ell (dist)}\theta^2_i\theta^2_j\theta^2_k\theta^2_\ell  = \|\theta\|^8\cdot \Bigl[ 1+ O\Bigl(\frac{\|\theta\|_4^4}{\|\theta\|^4}\Bigr)\Bigr]. 
\eeq
We combine \eqref{LemmaC1-eq3}-\eqref{LemmaC1-eq4} and note that $\theta_{\max}=o(1)$ and $\|\theta\|_4^4/\|\theta\|^4\leq (\|\theta\|^2\theta_{\max}^2)/\|\theta\|^4=o(1)$. So,
\[
\mathrm{Var}(\widetilde{Q}_n) = 8\|\theta\|^8\cdot [1+o(1)]. 
\]
This completes the proof.

\subsubsection{Proof of Lemma~\ref{lem:Normality}}
Under the null hypothesis, 
\[
\widetilde{Q}_n = X_1 = \sum_{i,j,k,\ell (dist)} W_{ij}W_{jk}W_{k\ell}W_{\ell i}. 
\]
In the proof of Theorem 3.2 of \cite{OGC}, it has been shown that $X_1/\sqrt{\mathrm{Var}(X_1)}\to N(0,1)$ in law (in the proof there, $X_1/\sqrt{\mathrm{Var}(X_1)}$ is denoted as $S_{n,n}$). Since $\mathbb{E}[X_1]=0$, we can directly quote their results to get the desired claim.

\subsubsection{Proof of Lemma~\ref{lem:IdealSgnQ-alt}}
We shall study the mean and variance of each of $X_1$-$X_6$ and then combine those results. 

Consider $X_1$. We have analyzed this term under the null hypothesis. Under the alternative hypothesis, the difference is that we no longer have $\Omega_{ij}=\theta_i\theta_j$. Instead, we have an upper bound $\Omega_{ij}= \theta_i\theta_j(\pi_i'P\pi_j)\leq C\theta_i\theta_j$. Using similar proof as that for the null hypothesis, we can derive that
\beq \label{LemmaC2-X1}
\mathbb{E}[X_1]=0, \qquad \mathrm{Var}(X_1)\leq C\|\theta\|^8. 
\eeq
To get a lower bound for $\mathrm{Var}(X_1)$, we notice that $\mathrm{Var}(W_{ij})=\Omega_{ij}(1-\Omega_{ij})\geq \Omega_{ij}[1-O(\theta_{\max}^2)]\geq \Omega_{ij}/2$; this inequality is true even when $\Omega_{ij}=0$. It follows that
\[
\mathrm{Var}(W_{ij} W_{jk} W_{k \ell} W_{\ell i}) \geq \frac{1}{16}\Omega_{ij}\Omega_{jk}\Omega_{k\ell}\Omega_{\ell i}. 
\]
Note that the second last line of \eqref{LemmaC1-eq3} is still true. As a result, 
\begin{align*}
\mathrm{Var}(X_1)& =8\sum_{i,j,k,\ell (dist)} \mathrm{Var}\bigl( W_{ij}W_{jk}W_{k\ell}W_{\ell i} \bigr)\cr
&\geq \frac{1}{2}\sum_{i,j,k,\ell (dist)} \Omega_{ij}\Omega_{jk}\Omega_{k\ell}\Omega_{\ell i}\cr
&= \frac{1}{2}\tr(\Omega^4) - \frac{1}{2}\sum_{i,j,k,\ell (not\; dist)} \Omega_{ij}\Omega_{jk}\Omega_{k\ell}\Omega_{\ell i}\cr
&\geq \frac{1}{2}\tr(\Omega^4) - C\sum_{i,j,k,\ell (not\; dist)}\theta^2_i\theta^2_j\theta^2_k\theta^2_\ell\cr
&\geq \frac{1}{2}\tr(\Omega^4) - o(\|\theta\|^8), 
\end{align*}
where the last inequality is due to \eqref{LemmaC1-eq4}. Recall that $\lambda_1,\ldots,\lambda_K$ denote the $K$ nonzero eigenvalues of $\Omega$. By Lemma~\ref{lem:Omega}, $\lambda_1\geq C^{-1}\|\theta\|^2$. It follows that
\[
\tr(\Omega^4)=\sum_{k=1}^K\lambda_k^4\geq \lambda_1^4\geq C^{-1}\|\theta\|^8. 
\]
Combining the above gives
\beq \label{LemmaC2-X1add}
\mathrm{Var}(X_1)\geq C^{-1}\|\theta\|^8. 
\eeq
So far, we have proved all claims about $X_1$. 

Consider $X_2$. Recall that
\[
X_2 = \sum_{i, j, k, \ell (dist)} \widetilde{\Omega}_{ij} W_{jk} W_{k \ell} W_{\ell i}. 
\]
It is easy to see that $\mathbb{E}[X_2]=0$. Below, we bound its variance. 
Each index choice $(i,j,k,\ell)$ defines a undirected path $j$-$k$-$\ell$-$i$ in the complete graph of $n$ nodes. If the two paths $j$-$k$-$\ell$-$i$ and $j'$-$k'$-$\ell'$-$i'$ are not exactly overlapping, then $W_{jk}W_{k\ell}W_{\ell i}\cdot W_{j'k'} W_{k' \ell'} W_{\ell' i'}$ have mean zero. In the sum above, each unique path $j$-$k$-$\ell$-$i$ is counted twice as $(i,j,k,\ell)$ and $(j,i,\ell,k)$. Mimicking the argument in \eqref{LemmaC1-eq3}, we immediately have
\begin{align*}
\mathrm{Var}(X_2) & = 2\sum_{i,j,k,\ell (dist)} \mathrm{Var}\bigl(  \widetilde{\Omega}_{ij} W_{jk} W_{k \ell} W_{\ell i} \bigr)\cr
&= 2\sum_{i,j,k,\ell (dist)} \widetilde{\Omega}_{ij}^2\cdot \mathrm{Var}\bigl( W_{jk} W_{k \ell} W_{\ell i} \bigr).
\end{align*}
By Lemma~\ref{lem:tildeOmega}, $|\widetilde{\Omega}_{ij}|\leq |\lambda_2|\|\theta\|^{-2}\theta_i\theta_j$. In our notations, $\alpha=|\lambda_2|/\lambda_1$; additionally, by Lemma~\ref{lem:Omega}, $\lambda_1\leq C\|\theta\|^2$. Combining them gives
\beq \label{tOmega(ij)}
|\widetilde{\Omega}_{ij}|\leq C\alpha\theta_i\theta_j. 
\eeq
Moreover, $\mathrm{Var}( W_{jk} W_{k \ell} W_{\ell i})\leq \Omega_{jk}\Omega_{k\ell}\Omega_{\ell i}\leq C\theta_j\theta_k^2\theta_\ell^2\theta_i$. It follows that
\begin{align*}
\mathrm{Var}(X_2)& \leq C\sum_{i,j,k,\ell (dist)}(\alpha\theta_i\theta_j)^2\cdot \theta_j\theta_k^2\theta_\ell^2\theta_i\cr
&\leq C\alpha^2 \sum_{i,j,k,\ell }\theta^3_i\theta_j^3\theta_k^2\theta_\ell^2\cr
&\leq C\alpha^2 \|\theta\|^4\|\theta\|_3^6. 
\end{align*}
Since $\|\theta\|_3^3\leq \theta_{\max}\sum_i\theta_i^2=\theta_{\max}\|\theta\|^2$, the right hand side is $\leq C\alpha^2\|\theta\|^8\theta_{\max}^2$. Note that $\alpha\leq 1$ and $\theta_{\max}\to 0$. So, this term is $o(\|\theta\|^8)$. 
We have proved all claims about $X_2$. 

Consider $X_3$. Recall that
\[
X_3 = \sum_{i, j, k, \ell (dist)} \widetilde{\Omega}_{ij}  \widetilde{\Omega}_{jk} W_{k \ell} W_{\ell i}=\sum_{i, k,\ell (dist)}\Bigl(\sum_{j\notin\{i,k,\ell\}} \widetilde{\Omega}_{ij}  \widetilde{\Omega}_{jk} \Bigr)W_{k \ell} W_{\ell i}. 
\]
It is easy to see that $\mathbb{E}[X_3]=0$. We then study its variance. We note that for $W_{k\ell}W_{\ell i}$ and $W_{k'\ell'}W_{\ell' i'}$ to be correlated, we must have that $(k',\ell',i')=(k,\ell, i)$ or $(k',\ell',i')=(i,\ell, k)$; in other words, the two underlying paths $k$-$\ell$-$i$ and $k'$-$\ell'$-$i'$ have to be equal. Mimicking the argument in \eqref{LemmaC1-eq3}, we have
\begin{align*}
\mathrm{Var}(X_3) &\leq C\sum_{i, k,\ell (dist)}\mathrm{Var}\Bigl[\Bigl(\sum_{j\notin\{i,k,\ell\}} \widetilde{\Omega}_{ij}  \widetilde{\Omega}_{jk} \Bigr)W_{k \ell} W_{\ell i}\Bigr]\cr
&\leq  C\sum_{i, k,\ell (dist)}\Bigl(\sum_{j\notin\{i,k,\ell\}} \widetilde{\Omega}_{ij}  \widetilde{\Omega}_{jk} \Bigr)^2\cdot\mathrm{Var}(W_{k \ell} W_{\ell i}). 
\end{align*}
By \eqref{tOmega(ij)},
\[
\Bigl|\sum_{j\notin\{i,k,\ell\}} \widetilde{\Omega}_{ij}  \widetilde{\Omega}_{jk} \Bigr|\leq C\sum_{j}\alpha^2\theta_i\theta_j^2\theta_k\leq C\alpha^2\|\theta\|^2\cdot \theta_i\theta_k. 
\]
Combining the above gives
\begin{align*}
\mathrm{Var}(X_3) & \leq C\sum_{i,k,\ell}(\alpha^2\|\theta\|^2\theta_i\theta_k)^2\cdot \theta_k\theta_\ell^2\theta_i \cr
&\leq C\alpha^4\|\theta\|^4\sum_{i,k,\ell}\theta_i^3\theta_k^3\theta_\ell^2\cr
&\leq C\alpha^4\|\theta\|^6\|\theta\|_3^6. 
\end{align*}
Since $\|\theta\|\to \infty$, the right hand side is $o(\alpha^4\|\theta\|^8\|\theta\|_3^6)$. We have proved all claims about $X_3$. 

Consider $X_4$. Recall that 
\begin{align*}
X_4 & = \sum_{i, j, k, \ell (dist)} \widetilde{\Omega}_{ij} W_{jk} \widetilde{\Omega}_{k \ell} W_{\ell i}=\sum_{i, j, k, \ell (dist)} \widetilde{\Omega}_{ij} \widetilde{\Omega}_{k \ell} W_{jk} W_{\ell i}. 
\end{align*}
It is easy to see that $\mathbb{E}[X_4]=0$. To calculate its variance, note that $W_{jk}W_{\ell i}$ and $W_{j'k'}W_{\ell' i'}$ are uncorrelated unless (i) $\{j',k'\}=\{j,k\}$ and $\{\ell',i'\}=\{\ell,i\}$ or (ii) $\{j',k'\}=\{\ell ,i\}$ and $\{\ell',i'\}=\{j,k\}$. Mimicking the argument in \eqref{LemmaC1-eq3}, we immediately have
\begin{align*}
\mathrm{Var}(X_4) &\leq C \sum_{i,j,k,\ell (dist)}\mathrm{Var}\bigl( \widetilde{\Omega}_{ij} \widetilde{\Omega}_{k \ell} W_{jk} W_{\ell i} \bigr)\cr
&\leq C \sum_{i,j,k,\ell (dist)} \widetilde{\Omega}^2_{ij} \widetilde{\Omega}^2_{k \ell}\cdot \mathrm{Var}( W_{jk} W_{\ell i})\cr
&\leq C\sum_{i,j,k,\ell} (\alpha\theta_i\theta_j)^2(\alpha\theta_k\theta_\ell)^2\cdot \theta_j\theta_k\theta_\ell\theta_i\cr
&\leq C\alpha^4\sum_{i,j,k,\ell}\theta_i^3\theta_j^3\theta_k^3\theta_\ell^3\cr
&\leq C\alpha^4\|\theta\|_3^{12}. 
\end{align*}
Since $\|\theta\|_3^3\leq \theta_{\max}\|\theta\|^2=o(\|\theta\|^2)$, the right hand side is $o(\|\theta\|^8)$. This proves the claims of $X_4$. 

Consider $X_5$. Recall that
\[
X_5 = \sum_{i, j, k, \ell (dist)} \widetilde{\Omega}_{ij} \widetilde{\Omega}_{jk} \widetilde{\Omega}_{k \ell} W_{\ell i} =2 \sum_{i< \ell}\Bigl( \sum_{\substack{j,k\notin\{i,\ell\}\\j\neq k}}\widetilde{\Omega}_{ij} \widetilde{\Omega}_{jk} \widetilde{\Omega}_{k \ell} \Bigr)W_{\ell i}. 
\]
It is easily seen that $\mathbb{E}[X_5]=0$. Furthermore, we have
\beq  \label{LemmaC2-X5}
\mathrm{Var}(X_5) = 2\sum_{i<\ell }\Bigl( \sum_{\substack{j,k\notin\{i,\ell\}\\j\neq k}}\widetilde{\Omega}_{ij} \widetilde{\Omega}_{jk} \widetilde{\Omega}_{k \ell} \Bigr)^2\cdot \mathrm{Var}(W_{\ell i}). 
\eeq
By \eqref{tOmega(ij)}, 
\[
\Bigl|\sum_{\substack{j,k\notin\{i,\ell\}\\j\neq k}}\widetilde{\Omega}_{ij} \widetilde{\Omega}_{jk} \widetilde{\Omega}_{k \ell}\Bigr|\leq C\sum_{j,k}\alpha^3\theta_i\theta_j^2\theta_k^2\theta_\ell \leq C\alpha^3\|\theta\|^4\cdot \theta_i\theta_\ell
\]
We plug it into \eqref{LemmaC2-X5} and use $\mathrm{Var}(W_{\ell i})\leq \Omega_{\ell i}\leq C\theta_\ell\theta_i$. It yields that 
\begin{eqnarray} \label{LemmaC2-X5-UB}
\mathrm{Var}(X_5) & \leq & C\sum_{\ell,i (dist)}(\alpha^3\|\theta\|^4\theta_i\theta_\ell)^2\cdot \theta_\ell\theta_i \cr
&\leq & C\alpha^6\|\theta\|^8\sum_{\ell,i}\theta_i^3\theta_\ell^3\cr
&\leq &  C\alpha^6\|\theta\|^8\|\theta\|_3^6. 
\end{eqnarray}
This proves the claims of $X_5$.  

Consider $X_6$. Recall that
\[
X_6 = \sum_{i, j, k, \ell (dist)} \widetilde{\Omega}_{ij} \widetilde{\Omega}_{jk} \widetilde{\Omega}_{k \ell} \widetilde{\Omega}_{\ell i}=\tr(\widetilde{\Omega}^4) - \sum_{i, j, k, \ell (not\; dist)} \widetilde{\Omega}_{ij} \widetilde{\Omega}_{jk} \widetilde{\Omega}_{k \ell} \widetilde{\Omega}_{\ell i}.
\]
This is a non-stochastic number, so its variance is zero and its mean is $X_6$ itself. By Lemma~\ref{lem:tildeOmega}, $|\lambda_2|\leq \|\widetilde{\Omega}\|\leq C |\lambda_2|$. Since $\|\widetilde{\Omega}\|^4 \leq \tr(\widetilde{\Omega}^4)\leq K\|\widetilde{\Omega}\|^4$, we immediately have $\tr(\widetilde{\Omega}^4)\asymp\|\widetilde{\Omega}\|^4\asymp |\lambda_2|^4$. Additionally, $|\lambda_2|=\alpha\lambda_1$ in our notation, and $\lambda_1\asymp \|\theta\|^2$ by Lemma~\ref{lem:Omega}. It follows that
\[
\tr(\widetilde{\Omega}^4)\asymp |\lambda_2|^4 \asymp \alpha^4\|\theta\|^8. 
\]
At the same time, by \eqref{tOmega(ij)}, $|\widetilde{\Omega}_{ij} \widetilde{\Omega}_{jk} \widetilde{\Omega}_{k \ell} \widetilde{\Omega}_{\ell i}|\leq C\alpha^4\theta_i^2\theta_j^2\theta_k^2\theta_\ell^2$. 
We thus have
\begin{align*}
|X_6 - \tr(\widetilde{\Omega}^4)| &\leq C\alpha^4 \sum_{i, j, k, \ell (not\; dist)} \theta_i^2\theta_j^2\theta_k^2\theta_\ell^2\cr
& \leq C\alpha^4\sum_{i,j,k}\theta_i^2\theta_j^2\theta_k^4\cr
&\leq C\alpha^4\|\theta\|^4\|\theta\|_4^4=o(\alpha^4\|\theta\|^8), 
\end{align*}
where the last equality is due to $\|\theta\|_4^4\leq\theta_{\max}^2\|\theta\|^2=o(\|\theta\|^4)$. 
Combining the above gives
\[
X_6 = \tr(\widetilde{\Omega}^4)\cdot [1+o(1)]. 
\]
This proves the claims of $X_6$. 

Last, we combine the results for $X_1$-$X_6$ to study $\widetilde{Q}_n$. Note that
\[
\widetilde{Q}_n = X_1 + 4X_2 + 4X_3 + 2X_4 + 4X_5 + X_6. 
\]
Only $X_6$ has a nonzero mean. So,  
\[
\mathbb{E}[\widetilde{Q}_n]=\mathbb{E}[X_6]=\tr(\widetilde{\Omega}^4)\cdot [1+o(1)].
\]
At the same time, given random variables $Z_1,Z_2,\ldots,Z_m$, $\mathrm{Var}(\sum_{k=1}^m Z_k)=\sum_k \mathrm{Var}(Z_k)+\sum_{k\neq \ell}\mathrm{Cov}(Z_k, Z_\ell)\leq \sum_k \mathrm{Var}(Z_k)+\sum_{k\neq \ell}\sqrt{\mathrm{Var}(Z_k)\mathrm{Var}(Z_\ell)}\leq m^2 \max_{k}\{\mathrm{Var}(Z_k)\}$. We thus have
\[
\mathrm{Var}(\widetilde{Q}_n)\leq C\max_{1\leq k\leq 6}\mathrm{Var}(X_k)\leq C\bigl(\|\theta\|^8+\alpha^6\|\theta\|^8\|\theta\|_3^6\bigr). 
\]
The proof of this lemma is now complete.

\subsubsection{Proof of Lemma~\ref{lem:ProxySgnQ(a)-null}}
Recall that $U_a = 4Y_1 = 4\sum_{i, j, k, \ell (dist)} \delta_{ij} W_{jk} W_{k \ell} W_{\ell i}$. 
By definition, $\delta_{ij}=\eta_i(\eta_j-\tilde{\eta}_j)+\eta_j(\eta_i-\tilde{\eta}_i)$.
It follows that 
\[
U_a =4\sum_{i, j, k, \ell (dist)} \eta_i(\eta_j - \tilde{\eta}_j)W_{jk} W_{k \ell} W_{\ell i}+ 4\sum_{i, j, k, \ell (dist)} \eta_j(\eta_i - \tilde{\eta}_i)W_{jk} W_{k \ell} W_{\ell i}. 
\]
In the second sum, if we relabel $(i,j,k,\ell)=(j',i',\ell',k')$, it becomes
\[
4\sum_{i', j', k', \ell' (dist)} \eta_{i'}(\eta_{j'} - \tilde{\eta}_{j'})W_{i'\ell'} W_{\ell' k'} W_{k' j'}=4\sum_{i, j, k, \ell (dist)} \eta_{i}(\eta_{j} - \tilde{\eta}_{j})W_{i\ell} W_{\ell k} W_{k j},
\]
which is the same as the first term. It follows that
\[
U_a =8\sum_{i, j, k, \ell (dist)} \eta_i(\eta_j - \tilde{\eta}_j)W_{jk} W_{k \ell} W_{\ell i}. 
\]
By definition, $\eta_j=\frac{1}{\sqrt{v}}\sum_{s\neq j}\mathbb{E}A_{js}$ and $\teta_j=\frac{1}{\sqrt{v}}\sum_{s\neq j}A_{js}$. Hence,  
\beq \label{tilde-eta}
\teta_j - \eta_j = \frac{1}{\sqrt{v}}\sum_{s\neq j}W_{js}.
\eeq 
We then re-write
\begin{eqnarray*}
U_a & =&-8\sum_{i, j, k, \ell (dist)} \eta_i\Bigl(\frac{1}{\sqrt{v}}\sum_{s\neq j}W_{js}\Bigr)W_{jk} W_{k \ell} W_{\ell i}\cr
&=& -\frac{8}{\sqrt{v}}\sum_{\substack{i, j, k, \ell (dist)\\s\neq j}} \eta_iW_{js}W_{jk} W_{k \ell} W_{\ell i}. 
\end{eqnarray*}
In the summand, $(i,j,k,\ell)$ are distinct, but $s$ is only required to be distinct from $j$. We consider two different cases: (a) the case of $s=k$, where the summand becomes $W_{jk}^2W_{k\ell}W_{\ell i}$, and (b) the case of $s\neq k$. Correspondingly, we write
\begin{eqnarray} \label{proof-Ua-nullmean}
U_a  &=& - \frac{8}{\sqrt{v}}\sum_{i, j, k, \ell (dist)} \eta_iW^2_{jk} W_{k \ell} W_{\ell i}-  \frac{8}{\sqrt{v}}\sum_{\substack{i, j, k, \ell (dist)\\s\notin\{ j,k\}}} \eta_iW_{js}W_{jk} W_{k \ell} W_{\ell i}\cr
&\equiv& U_{a1}+U_{a2}. 
\end{eqnarray}
It is easy to see that the summands in both sums have mean zero. Therefore,
\[ 
\mathbb{E}[U_a]=0. 
\]

Next, we bound the variance of $U_a$. Since $\mathrm{Var}(U_a)\leq 2\mathrm{Var}(U_{a1})+2\mathrm{Var}(U_{a2})$, it suffices to bound the variances of $U_{a1}$ and $U_{a2}$. Consider $U_{a1}$. Note that
\beq \label{proof-Ua-null1}
\mathrm{Var}(U_{a1}) =\frac{64}{v}\sum_{\substack{i, j, k, \ell (dist)\\i', j',k',\ell' (dist)}} \eta_i\eta_{i'}\cdot \mathbb{E}[W^2_{jk} W_{k \ell} W_{\ell i}W^2_{j'k'} W_{k'\ell'} W_{\ell' i'}]. 
\eeq
By definition, $v=1_n'(\mathbb{E}A)1_n=1_n'\Omega 1_n-\sum_i\Omega_{ii}$. Since $\Omega_{ii}\leq \theta_i^2$, it implies  $v =1_n'\Omega 1_n-O(\|\theta\|^2)= 1_n'\Omega 1_n+o(\|\theta\|_1^2)$. Moreover, we note that $1_n'\Omega 1_n\leq C\sum_{i,j}\theta_i\theta_j\leq C\|\theta\|_1^2$, and by Lemma~\ref{lem:1Omega1}, $1_n'\Omega 1_n\geq C^{-1}\|\theta\|^2_1$. Combining these results gives
\beq \label{v-order}
C^{-1}\|\theta\|_1^2\leq v\leq C\|\theta\|_1^2. 
\eeq
Moreover, $\eta_i = \frac{1}{\sqrt{v}}\sum_{s\neq i}\Omega_{is}\leq \frac{C}{\|\theta\|_1}\sum_{s}\theta_i\theta_s$. This gives
\beq \label{eta-bound}
0\leq \eta_i\leq C\theta_i, \qquad \mbox{for all }1\leq i\leq n. 
\eeq
We plug \eqref{v-order}-\eqref{eta-bound} into \eqref{proof-Ua-null1} and find out that
\[
\mathrm{Var}(U_{a1})\leq \frac{C}{\|\theta\|_1^2}\sum_{\substack{i, j, k, \ell (dist)\\i', j',k',\ell' (dist)}} \theta_i\theta_{i'}\cdot \mathbb{E}[W^2_{jk} W_{k \ell} W_{\ell i}W^2_{j'k'} W_{k'\ell'} W_{\ell' i'}]. 
\]
In order for the summand to be nonzero, all $W$ terms have to be perfectly paired. By elementary calculations,
\[
\theta_i\theta_{i'} \mathbb{E}[W^2_{jk} W_{k \ell} W_{\ell i}W^2_{j'k'} W_{k'\ell'} W_{\ell' i'}]=
\begin{cases}
\theta_i^2\,\mathbb{E}[W^2_{jk}W^2_{k\ell}W^2_{\ell i}W^2_{j'k}], & \mbox{if }(\ell',k',i')\text{=}(\ell,k,i);\\
\theta_i\theta_k\, \mathbb{E}[W^2_{jk}W^2_{k\ell}W^2_{\ell i}W^2_{j'i}], & \mbox{if }(\ell', k',i')\text{=}(\ell, i,k);\\
\theta_i\theta_j\, \mathbb{E}[W^3_{jk}W^2_{ k\ell}W^3_{\ell i}], & \mbox{if } (j',k')\text{=}(i,\ell),\, (i', \ell')\text{=}(j,k);\\
0, & \mbox{otherwise}. 
\end{cases}
\]
Here, $(i,j,k,\ell)$ are distinct. In the second case above, $(W^2_{jk}, W^2_{k\ell}, W^2_{\ell i}, W^2_{j'i})$ are independent of each other, no matter $j=j'$ or $j\neq j'$ (we remark that $j'\neq \ell$, because $j'\notin\{i',k',\ell'\}=\{i,k,\ell\}$). It follows that $\mathbb{E}[W^2_{jk}W^2_{k\ell}W^2_{\ell i}W^2_{j'i}]\leq \Omega_{jk}\Omega_{k\ell}\Omega_{\ell i}\Omega_{j'i}\leq C\theta^2_i\theta_j\theta_k^2\theta^2_\ell\theta_{j'}$. In the first case, when $j\neq j'$, $\mathbb{E}[W^2_{jk}W^2_{k\ell}W^2_{\ell i}W^2_{j'k}]\leq \Omega_{jk}\Omega_{k\ell}\Omega_{\ell i}\Omega_{j'k}\leq C\theta_i\theta_j\theta_k^3\theta^2_\ell\theta_{j'}$; when $j=j'$, it holds that $\mathbb{E}[W^2_{jk}W^2_{k\ell}W^2_{\ell i}W^2_{j'k}]= \mathbb{E}[W^4_{jk}W^2_{k\ell}W^2_{\ell i}]\leq C\theta_i\theta_j\theta_k^2\theta^2_\ell$. In the third case, $(W^3_{jk}, W^2_{k\ell}, W^3_{\ell i})$ are mutually independent, so $\mathbb{E}[W^2_{jk}W^2_{k\ell}W^2_{\ell i}]\leq \Omega_{jk}\Omega_{k\ell}\Omega_{\ell i}\leq C\theta_i\theta_j\theta_k^2\theta^2_\ell$. We then have
\[
\theta_i\theta_{i'} \mathbb{E}[W^2_{jk} W_{k \ell} W_{\ell i}W^2_{j'k'} W_{k'\ell'} W_{\ell' i'}]
\leq  
\begin{cases}
C\theta^3_i\theta_j\theta_k^2\theta^2_\ell, & \mbox{if }(\ell',k',i')=(\ell,k,i),\; j'=j;\\
C\theta_i^3\theta_j\theta_k^3\theta^2_\ell\theta_{j'}, & \mbox{if }(\ell',k',i')=(\ell,k,i),\; j'\neq j;\\
C\theta_i^3\theta_j\theta_k^3\theta^2_\ell\theta_{j'}, & \mbox{if }(\ell', k',i')=(\ell, i,k);\\
C\theta^2_i\theta^2_j\theta_k^2\theta^2_\ell, & \mbox{if } (j',k')\text{=}(i,\ell),\, (i', \ell')\text{=}(j,k);\\
0, & \mbox{otherwise}. 
\end{cases}
\]
It follows that
\begin{eqnarray} \label{proof-Ua-null2}
\mathrm{Var}(U_{a1}) & \leq& \frac{C}{\|\theta\|_1^2}\Bigl( \sum_{i,j,k,\ell}\theta^3_i\theta_j\theta_k^2\theta^2_\ell +\sum_{i,j,k,\ell,j'}\theta_i^3\theta_j\theta_k^3\theta^2_\ell\theta_{j'} + \sum_{i,j,k,\ell}\theta_i^2\theta_j^2\theta_k^2\theta_\ell^2\Bigr)\cr
&\leq&  \frac{C}{\|\theta\|^2_1}\bigl( \|\theta\|^4\|\theta\|_3^3\|\theta\|_1 + \|\theta\|^2\|\theta\|_3^6\|\theta\|_1^2 +\|\theta\|^8\bigr) \cr
&\leq & C\|\theta\|^2\|\theta\|_3^6, 
\end{eqnarray}
where we obtain the last inequality as follows: By Cauchy-Schwarz inequality, $\|\theta\|^4 = (\sum_i\theta_i^{1/2}\cdot\theta^{3/2})^2\leq (\sum_i\theta_i)(\sum_i\theta_i^3)\leq \|\theta\|_1\|\theta\|_3^3$; therefore, $\|\theta\|^8\leq \|\theta\|^4\|\theta\|_3^3\|\theta\|_1\leq \|\theta\|_3^6\|\theta\|_1^2$. 
We then consider $U_{a2}$. Define
\[
{\cal P}_5^*=\biggl\{\begin{array}{ll}
\mbox{path $i$-$\ell$-$k$-$j$-$s$ in a complete}: &\mbox{nodes $i,j,k,\ell$ are distinct,}\\
\mbox{graph with $n$ nodes} & \mbox{and node $s$ is different from $j$, $k$}
\end{array}
\biggr\}. 
\]
Fix a path $i$-$\ell$-$k$-$j$-$s$ in ${\cal P}_5^*$. If $s\notin\{i,\ell\}$, then this path is counted twice in the definition of $U_{a2}$, as $i$-$\ell$-$k$-$j$-$s$ and $s$-$j$-$k$-$\ell$-$i$, respectively. If $s\in\{i,\ell\}$, then it is counted only once in the definition of $U_{a2}$. 
Hence, we can re-write
\[
U_{a2} = - \frac{8}{\sqrt{v}}\sum_{\substack{\text{path in ${\cal P}_5^*$}\\s\notin\{i,\ell\}}} (\eta_i+\eta_s) W_{sj}W_{jk} W_{k \ell} W_{\ell i}- \frac{8}{\sqrt{v}}\sum_{\substack{\text{path in ${\cal P}_5^*$}\\s\in\{i,\ell\}}} \eta_i W_{sj}W_{jk} W_{k \ell} W_{\ell i}.
\]
For two distinct paths in ${\cal P}_5^*$, the corresponding summands are uncorrelated with each other. It follows that
\begin{align}  \label{proof-Ua-null3}
\mathrm{Var}(U_{a2}) &= \frac{64}{v}\sum_{\substack{\text{path in ${\cal P}_5^*$}\\s\notin\{i,\ell\}}} (\eta_i+\eta_s)^2\,\mathrm{Var}(W_{sj}W_{jk} W_{k \ell} W_{\ell i})\cr
&\qquad + \frac{64}{v}\sum_{\substack{\text{path in ${\cal P}_5^*$}\\s\in\{i,\ell\}}} \eta^2_i \,\mathrm{Var}(W_{sj}W_{jk} W_{k \ell} W_{\ell i})\cr
&\leq \frac{C}{v}\sum_{i,j,k,\ell,s}(\eta_i^2+\eta_s^2)\cdot\theta_i\theta^2_j\theta^2_k\theta^2_\ell\theta_s\cr
&\leq \frac{C}{\|\theta\|_1^2}\sum_{i,j,k,\ell,s}(\theta^3_i\theta^2_j\theta^2_k\theta^2_\ell\theta_s+\theta_i\theta^2_j\theta^2_k\theta^2_\ell\theta^3_s)\cr
&\leq \frac{C\|\theta\|^6\|\theta\|_3^3}{\|\theta\|_1}. 
\end{align}
By Cauchy-Schwarz inequality, $\|\theta\|^4\leq \|\theta\|_1\|\theta\|_3^3$, so
the right hand side of \eqref{proof-Ua-null3} is $\leq C\|\theta\|^2\|\theta\|_3^6$. Combining it with \eqref{proof-Ua-null2} gives
\[
\mathrm{Var}(U_a)\leq C \|\theta\|^2\|\theta\|_3^6=o(\|\theta\|^8). 
\]
This proves the claim.

\subsubsection{Proof of Lemma~\ref{lem:ProxySgnQ(a)-alt}}
It suffices to prove the claims for each of $Y_1$-$Y_6$. Consider $Y_1$. We have analyzed this term under the null hypothesis. Using similar proof, we can easily derive that
\[
\mathbb{E}[Y_1]=0, \qquad \mathrm{Var}(Y_1)\leq C\|\theta\|^2\|\theta\|_3^6=o(\|\theta\|^8). 
\]

Consider $Y_2$. Using the definition of $Y_2$ and the expression of $\teta_i$ in \eqref{tilde-eta}, we have
\begin{align*} 
Y_2 & =  \sum_{i, j, k, \ell (dist)} \delta_{ij} \widetilde{\Omega}_{jk}  W_{k \ell}  W_{\ell i}\cr
&= \sum_{i, j, k, \ell (dist)} \eta_i(\eta_j-\teta_j) \widetilde{\Omega}_{jk}  W_{k \ell}  W_{\ell i} + \sum_{i, j, k, \ell (dist)} \eta_j(\eta_i-\teta_i) \widetilde{\Omega}_{jk}  W_{k \ell}  W_{\ell i}\cr
&= \frac{1}{\sqrt{v}} \sum_{i, j, k, \ell (dist)}\eta_i\Bigl(- \sum_{s\neq j}W_{js}\Bigr) \widetilde{\Omega}_{jk}  W_{k \ell}  W_{\ell i} + 
 \frac{1}{\sqrt{v}} \sum_{i, j, k, \ell (dist)} \eta_j \Bigl(-\sum_{s\neq i}W_{is}\Bigr)\widetilde{\Omega}_{jk}  W_{k \ell}  W_{\ell i}\cr
&=  - \frac{1}{\sqrt{v}} \sum_{\substack{i, j, k, \ell (dist)\\s\neq j}}\eta_i \widetilde{\Omega}_{jk}  W_{js} W_{k \ell}  W_{\ell i} - 
 \frac{1}{\sqrt{v}} \sum_{\substack{i, k, \ell (dist)\\s\neq i}} \Bigl(\sum_{j\notin\{i,k,\ell\}}\eta_j \widetilde{\Omega}_{jk}\Bigr)  W_{is}W_{k \ell}  W_{\ell i}.
 \end{align*}
In the second sum above, we further separate two cases, $s=\ell$ and $s\neq \ell$. It then gives rise to three terms: 
 \begin{eqnarray}  \label{proof-Y2-decompose}
Y_2 &=&  -\frac{1}{\sqrt{v}} \sum_{\substack{i, j, k, \ell (dist)\\s\neq j}}\eta_i \widetilde{\Omega}_{jk}  W_{js} W_{k \ell}  W_{\ell i}\cr
&& - \frac{1}{\sqrt{v}} \sum_{i, k, \ell (dist)} \Bigl(\sum_{j\notin\{i,k,\ell\}}\eta_j \widetilde{\Omega}_{jk}\Bigr) W^2_{i\ell} W_{k \ell}\cr
 && -   \frac{1}{\sqrt{v}} \sum_{\substack{i, k, \ell (dist)\\s\notin\{ i,\ell\}}} \Bigl(\sum_{j\notin\{i,k,\ell\}}\eta_j \widetilde{\Omega}_{jk}\Bigr)  W_{is}W_{k \ell}  W_{\ell i}\cr
&\equiv& Y_{2a} + Y_{2b} + Y_{2c}. 
\end{eqnarray}
Since $(i,j,k,\ell)$ are distinct, it is easy to see that all three terms have mean zero. We thus have 
\[
\mathbb{E}[Y_2] = 0. 
\]
Below, we calculate the variances. First, we bound the variance of $Y_{2a}$. Each $(i, j, k, \ell, s)$ is associated with a length-3 path $i$-$k$-$\ell$ and an edge $j$-$s$ in the complete graph. For $(i, j, k, \ell, s)$ and $(i', j', k', \ell', s')$, if the associated path and edge are the same, then we group them together. Given a length-3 path $i$-$k$-$\ell$ and an edge $j$-$s$ (such that the edge is not in the path), they are counted four times in the definition of $Y_{2a}$, as (i) $i$-$k$-$\ell$ and $j$-$s$, (ii) $i$-$k$-$\ell$ and $s$-$j$, (iii) $\ell$-$k$-$i$ and $j$-$s$, (iv) $\ell$-$k$-$i$ and $s$-$j$, so we group these four summands together. After grouping the summands, we re-write
\[
Y_{2a} = -\frac{1}{\sqrt{v}}\sum_{\substack{\text{length-3}\\\text{path}}}\, \sum_{\substack{\text{edge not}\\\text{in the path}}}\bigl( \eta_i \widetilde{\Omega}_{jk} + \eta_i \widetilde{\Omega}_{sk}+\eta_k \widetilde{\Omega}_{ji}+\eta_k \widetilde{\Omega}_{si}\bigr)W_{js}W_{k\ell}W_{\ell i}. 
\]
In this new expression of $Y_{2a}$, two summands are correlated only when the underlying path\&edge pairs are exactly the same. Additionally, 
by \eqref{tOmega(ij)} and \eqref{eta-bound}, 
\[
\bigl| \eta_i \widetilde{\Omega}_{jk} + \eta_i \widetilde{\Omega}_{sk}+\eta_k \widetilde{\Omega}_{ji}+\eta_k \widetilde{\Omega}_{si}\bigr|\leq C\alpha(\theta_j+\theta_s)\theta_i\theta_k. 
\]
It follows that
\begin{align} \label{proof-Y2-result1}
\mathrm{Var}(Y_{2a})& \leq  \frac{C}{v}\sum_{i,j,k,\ell,s} \alpha^2(\theta_j+\theta_s)^2\theta_i^2\theta_k^2\cdot \mathrm{Var}(W_{js} W_{k \ell}  W_{\ell i})\cr
&\leq \frac{C}{\|\theta\|_1^2}\sum_{i,j,k,\ell,s} \alpha^2(\theta_j+\theta_s)^2\theta_i^2\theta_k^2\cdot \theta_i\theta_j\theta_k\theta^2_\ell\theta_s\cr
&\leq \frac{C\alpha^2}{\|\theta\|_1^2}\sum_{i,j,k,\ell,s} (\theta^3_i\theta^3_j\theta^3_k\theta^2_\ell\theta_s+ \theta^3_i\theta_j\theta^3_k\theta^2_\ell\theta^3_s)\cr
&\leq \frac{C\alpha^2\|\theta\|^2\|\theta\|_3^9}{\|\theta\|_1}. 
\end{align}
Second, we bound the variance of $Y_{2b}$. Write $\beta_{ik\ell}
= \sum_{j\notin\{i,k,\ell\}}\eta_j \widetilde{\Omega}_{jk}$. By \eqref{tOmega(ij)} and \eqref{eta-bound}, $|\beta_{ik\ell}|\leq C\sum_j\theta_j\cdot \alpha\theta_j\theta_k\leq C\alpha\|\theta\|^2\theta_k$. Using this notation, 
\[
Y_{2b} = \frac{1}{v}\sum_{i,j,k,\ell (dist)}\beta_{ik\ell}W^2_{i\ell} W_{k \ell}, \qquad \mbox{where}\quad |\beta_{ik\ell}|\leq C\alpha\|\theta\|^2\theta_k. 
\]
It follows that 
\begin{align*}
\mathrm{Var}(Y_{2b}) =\mathbb{E}[Y_{2b}^2] &\leq \frac{C}{v}\sum_{\substack{i,k,\ell (dist)\\i',k', \ell' (dist)}} \beta_{ik\ell}\beta_{i'k'\ell'}\cdot \mathbb{E}[ W^2_{i\ell} W_{k \ell} W^2_{i'\ell'} W_{k' \ell'}]\cr
&\leq\frac{C\alpha^2\|\theta\|^4}{\|\theta\|_1^2}\sum_{\substack{i,k,\ell (dist)\\i',k', \ell' (dist)}} \theta_k\theta_{k'}\cdot \mathbb{E}[ W^2_{i\ell} W_{k \ell} W^2_{i'\ell'} W_{k' \ell'}].
\end{align*}
The summand is nonzero only when the two variables $W_{k\ell}$ and $W_{k'\ell'}$ equal to each other or when each of them equals to some other squared variables. By elementary calculations,
\begin{align*}
& \theta_k\theta_{k'}\cdot \mathbb{E}[ W^2_{i\ell} W_{k \ell} W^2_{i'\ell'} W_{k' \ell'}]\cr
=\; & 
\begin{cases}
\theta^2_k\, \mathbb{E}[ W^4_{i\ell} W^2_{k \ell}]\leq C\theta_i\theta^3_k\theta^2_\ell, &\mbox{if }(k',\ell')=(k,\ell),\; i'=i;\\
\theta^2_k\, \mathbb{E}[ W^2_{i\ell} W^2_{k \ell} W^2_{i'\ell}]\leq C\theta_i\theta^3_k\theta^3_\ell\theta_{i'}, &\mbox{if }(k',\ell')=(k,\ell),\; i'\neq i;\\
\theta_k\theta_{\ell}\, \mathbb{E}[ W^2_{i\ell} W^2_{k \ell} W^2_{i'k}]\leq C\theta_i\theta^3_k\theta^3_\ell\theta_{i'}, &\mbox{if }(k',\ell') = (\ell, k);\\
\theta^2_k\, \mathbb{E}[ W^3_{i\ell} W^3_{k \ell}]\leq C\theta_i\theta^3_k\theta^2_\ell, & \mbox{if }\ell'=\ell,\; (i', k')=(i, k);\\
\theta_k\theta_i \, \mathbb{E}[ W^3_{i\ell} W^3_{k \ell}]\leq C\theta^2_i\theta^2_k\theta^2_\ell, & \mbox{if }\ell'=\ell,\; (i', k')=(k, i);\\
0, & \mbox{otherwise}. 
\end{cases}
\end{align*}
As a result,
\begin{align} \label{proof-Y2-result2}
\mathrm{Var}(Y_{2b})& \leq \frac{C\alpha^2\|\theta\|^4}{\|\theta\|_1^2}\Bigl(\sum_{i,k,\ell}\theta_i\theta^3_k\theta^2_\ell+\sum_{i,k,\ell,i'}\theta_i\theta^3_k\theta^3_\ell\theta_{i'}  +\sum_{i, k, \ell}\theta_i^2\theta_k^2\theta_\ell^2 \Bigr)\cr
& \leq \frac{C\alpha^2\|\theta\|^4}{\|\theta\|_1^2}\bigl( \|\theta\|_3^3\|\theta\|^2\|\theta\|_1 + \|\theta\|_3^6\|\theta\|_1^2 + \|\theta\|^6\bigr)\cr 
&\leq C\alpha^2\|\theta\|^4\|\theta\|_3^6, 
\end{align}
where to get the last inequality we have used $\|\theta\|^6\ll \|\theta\|^8\leq (\|\theta\|_1\|\theta\|_3^3)^2$ and $\|\theta\|_3^3\|\theta\|^2\|\theta\|_1\ll \|\theta\|_3^3\|\theta\|^4\|\theta\|_1\leq (\|\theta\|_1\|\theta\|_3^3)^2$. Last, we bound the variance of $Y_{2c}$. Let $\beta_{ik\ell}
= \sum_{j\notin\{i,k,\ell\}}\eta_j \widetilde{\Omega}_{jk}$ be the same as above. We write
\[
Y_{2c} = \frac{1}{\sqrt{v}} \sum_{\substack{i, k, \ell (dist)\\s\notin\{ i,\ell\}}} \beta_{ik\ell}  W_{is}W_{k \ell}  W_{\ell i}, \qquad \mbox{where}\quad |\beta_{ik\ell}|\leq C\alpha\|\theta\|^2\theta_k. 
\]
For $\mathbb{E}[W_{is}W_{k\ell}W_{\ell i}\cdot W_{i's'}W_{k'\ell'}W_{\ell' i'}]$ to be nonzero, it has to be the case that $(W_{is}, W_{k\ell}, W_{\ell i})$ and $(W_{i's'}, W_{k'\ell'}, W_{\ell' i'})$ are the same set of variables, up to an order permutation. For each fixed $(i,k,\ell, s)$, there are only a constant number of $(i', k', \ell', s')$ such that the above is satisfied. As we have argued many times before (e.g., see \eqref{LemmaC1-eq3}), it is true that
\begin{align} \label{proof-Y2-result3}
\mathrm{Var}(Y_{2c})& \leq \frac{C}{v} \sum_{\substack{i, k, \ell (dist)\\s\notin\{ i,\ell\}}} \beta^2_{ik\ell}\cdot\mathrm{Var}(W_{is}W_{k \ell}  W_{\ell i})\cr
&\leq \frac{C}{\|\theta\|_1^2}\sum_{i,k,\ell,s}(\alpha\|\theta\|^2\theta_k)^2\cdot \theta_i^2\theta_k\theta_\ell^2\theta_s\cr
&\leq \frac{C\alpha^2\|\theta\|^8\|\theta\|_3^3}{\|\theta\|_1}. 
\end{align}
We now combine the variances of $Y_{2a}$-$Y_{2c}$. Since $\|\theta\|_3^3\leq \theta_{\max}^2\|\theta\|_1\ll \|\theta\|_1$, the right hand side is \eqref{proof-Y2-result1} is $o(\alpha^2\|\theta\|^2\|\theta\|_3^6)=o(\alpha^2\|\theta\|^4\|\theta\|_3^6)$. Since $\|\theta\|^4\leq \|\theta\|_1\|\theta\|_3^3$, the right hand side is \eqref{proof-Y2-result3} is $\leq C\alpha^2\|\theta\|^4\|\theta\|_3^6$. It follows that
\[
\mathrm{Var}(Y_2)\leq C\alpha^2\|\theta\|^4\|\theta\|_3^6=o(\|\theta\|^8). 
\]
This proves the claims of $Y_2$. 

Consider $Y_3$. By definition,
\[
Y_3 = \sum_{i, j, k, \ell (dist)}\eta_i(\eta_j-\tilde{\eta}_j) W_{jk}\widetilde{\Omega}_{k\ell}  W_{\ell i} + 
\sum_{i, j, k, \ell (dist)} \eta_j (\eta_i-\tilde{\eta}_i)W_{jk}\widetilde{\Omega}_{k \ell}  W_{\ell i}. 
\]
In the second sum, if we relabel $(i, j, k, \ell)=(j', i', \ell', k')$, it can be written as $\sum_{i', j', k', \ell' (dist)} \eta_{i'} (\eta_{j'}-\tilde{\eta}_{j'})W_{i'\ell'}\widetilde{\Omega}_{\ell' k'}  W_{k' j'}$. This shows that the second sum is indeed equal to the first sum. As a result,
\begin{align}\label{proof-Y3-decompose}
Y_3 &= 2\sum_{i, j, k, \ell (dist)}\eta_i(\eta_j-\tilde{\eta}_j) W_{jk}\widetilde{\Omega}_{k\ell}  W_{\ell i}\cr
& = 2\sum_{i, j, k, \ell (dist)}\eta_i\Bigl(-\frac{1}{\sqrt{v}}\sum_{s\neq j}W_{js}\Bigr) W_{jk}\widetilde{\Omega}_{k\ell}  W_{\ell i}\cr
& = - \frac{2}{\sqrt{v}} \sum_{\substack{i, j, k, \ell (dist)\\s\neq j}}\eta_i \widetilde{\Omega}_{k\ell}  W_{js} W_{jk}  W_{\ell i}\cr
 &= - \frac{2}{\sqrt{v}} \sum_{i, j, k, \ell (dist)}\eta_i \widetilde{\Omega}_{k\ell}  W^2_{jk}  W_{\ell i} - \frac{2}{\sqrt{v}} \sum_{\substack{i, j, k, \ell (dist)\\s\notin\{j,k\}}}\eta_i \widetilde{\Omega}_{k\ell}  W_{js} W_{jk}  W_{\ell i}\cr
 & \equiv Y_{3a} + Y_{3b},
\end{align}
where the second line is from \eqref{tilde-eta} and the second last line is from dividing all summands into two cases of $s=k$ and $s\neq k$. 
Both terms have mean zero, so 
\[
\mathbb{E}[Y_3] = 0. 
\]
Below, first, we calculate the variance of $Y_{3a}$. 
\[
\mathrm{Var}(Y_{3a}) = \frac{4}{v} \sum_{\substack{i,j,k,\ell (dist)\\i',j',k',\ell' (dist)}} (\eta_i \widetilde{\Omega}_{k\ell} \eta_{i'} \widetilde{\Omega}_{k'\ell'})\cdot\mathbb{E}[W^2_{jk}  W_{\ell i}W^2_{j'k'}W_{\ell' i'}]. 
\]
The summand is nonzero only if either the two variables $W_{\ell i}$ and $W_{\ell' i'}$ are the same, or each of the two variables $W_{\ell i}$ and $W_{\ell' i'}$ equals to another squared $W$ term. By \eqref{tOmega(ij)}, \eqref{eta-bound},  and elementary calculations,
\begin{align*} 
& (\eta_i \widetilde{\Omega}_{k\ell} \eta_{i'} \widetilde{\Omega}_{k'\ell'})\cdot\mathbb{E}[W^2_{jk}  W_{\ell i}W^2_{j'k'}W_{\ell' i'}]\cr
&\leq C\alpha^2 \theta_i\theta_k\theta_\ell \theta_{i'}\theta_{k'}\theta_{\ell'}\cdot \mathbb{E}[W^2_{jk}  W_{\ell i}W^2_{j'k'}W_{\ell' i'}]\cr
=\; & \begin{cases}
C\alpha^2 \theta^2_i\theta^2_\ell\theta^2_k \, \mathbb{E}[W^4_{jk}  W^2_{\ell i}]\leq C\alpha^2\theta_i^3\theta_j\theta_k^3\theta^3_\ell, &\mbox{if }\{\ell',i'\}=\{\ell,i\},\, (j', k')= (j, k);\\
C\alpha^2 \theta^2_i\theta^2_\ell\theta_k\theta_j \, \mathbb{E}[W^4_{jk}  W^2_{\ell i}]\leq C\alpha^2\theta_i^3\theta^2_j\theta_k^2\theta^3_\ell, &\mbox{if }\{\ell',i'\}=\{\ell,i\}, \, (j', k')= (k, j);\\
C\alpha^2 \theta^2_i\theta^2_\ell\theta_k\theta_{k'}\, \mathbb{E}[W^2_{jk}  W^2_{\ell i}W^2_{j'k'}]\leq C\alpha^2\theta_i^3\theta_j\theta_k^2\theta^3_\ell\theta_{j'}\theta^2_{k'}, &\mbox{if }\{\ell',i'\}=\{\ell,i\},\, \{j', k'\}\neq \{j, k\};\\
C\alpha^2 \theta^2_i\theta_\ell\theta_j\theta_k^2\, \mathbb{E}[W^3_{jk}  W^3_{\ell i}]\leq C\alpha^2\theta^3_i\theta^2_j\theta^3_k\theta^2_\ell , &\mbox{if }\{\ell', i'\}=\{j,k\},\, (j', k')=(\ell, i);\\
C\alpha^2 \theta_i\theta^2_\ell\theta_j\theta_k^2\, \mathbb{E}[W^3_{jk}  W^3_{\ell i}]\leq C\alpha^2\theta^2_i\theta^2_j\theta^3_k\theta^3_\ell , &\mbox{if }\{\ell', i'\}=\{j,k\},\, (j', k')=(i, \ell);\\
0, & \mbox{otherwise}. 
\end{cases}
\end{align*}
There are only three different cases in the bounds. 
It follows that
\begin{eqnarray} \label{proof-Y3-result1}
\mathrm{Var}(Y_{3a})& \leq& \frac{C\alpha^2}{\|\theta\|_1^2}\Bigl(\sum_{i,j,k,\ell}\theta_i^3\theta_j\theta_k^3\theta^3_\ell  + \sum_{i,j,k,\ell}\theta_i^3\theta^2_j\theta_k^2\theta^3_\ell + \sum_{i,j,k,\ell,j',k'} \theta_i^3\theta_j\theta_k^2\theta^3_\ell\theta_{j'}\theta^2_{k'} \Bigr)\cr
&\leq& \frac{C\alpha^2}{\|\theta\|_1^2}\bigl( \|\theta\|_1\|\theta\|_3^9 + \|\theta\|^4\|\theta\|_3^6+\|\theta\|^4\|\theta\|_1^2\|\theta\|_3^6 \bigr)\cr
&\leq&  C\alpha^2\|\theta\|^4\|\theta\|_3^6, 
\end{eqnarray}
where in the last line we have used $\|\theta\|_3^9\leq \|\theta\|_3^6 (\theta_{\max}\|\theta\|^2)=o(\|\theta\|^2\|\theta\|_3^6)$ and $\|\theta\|_1\geq \theta_{\max}^{-1}\|\theta\|^2\to\infty$. 
Next, we calculate the variance of $Y_{3b}$. We mimic the argument in \eqref{proof-Y2-result1} and group summands according to the underlying path $s$-$j$-$k$ and edge $\ell$-$i$ in a complete graph. It yields 
\[
Y_{3b} = -\frac{2}{\sqrt{v}}\sum_{\substack{\text{length-3}\\\text{path}}}\, \sum_{\substack{\text{edge not}\\\text{in the path}}}\bigl( \eta_i \widetilde{\Omega}_{k\ell} + \eta_\ell \widetilde{\Omega}_{ki}+\eta_i \widetilde{\Omega}_{s\ell}+\eta_\ell \widetilde{\Omega}_{si}\bigr)W_{sj}W_{jk}W_{\ell i}, 
\]
where
\[
\bigl| \eta_i \widetilde{\Omega}_{k\ell} + \eta_\ell \widetilde{\Omega}_{ki}+\eta_i \widetilde{\Omega}_{s\ell}+\eta_\ell \widetilde{\Omega}_{si}\bigr|\leq C\alpha(\theta_k+\theta_s)\theta_i\theta_\ell. 
\]
It follows that
\begin{align} \label{proof-Y3-result2}
\mathrm{Var}(Y_{3b})& \leq  \frac{C}{v}\sum_{i,j,k,\ell,s} \alpha^2(\theta_k+\theta_s)^2\theta_i^2\theta_\ell^2\cdot \mathrm{Var}(W_{sj} W_{jk}  W_{\ell i})\cr
&\leq \frac{C\alpha^2}{\|\theta\|_1^2}\sum_{i,j,k,\ell,s} (\theta^3_i\theta^2_j\theta^3_k\theta^3_\ell\theta_s+ \theta^3_i\theta^2_j\theta_k\theta^3_\ell\theta^3_s)\cr
&\leq \frac{C\alpha^2\|\theta\|^2\|\theta\|_3^9}{\|\theta\|_1}. 
\end{align}
Since $\|\theta\|_3^9\leq \|\theta\|_3^6(\theta_{\max}\|\theta\|_1)=o(\|\theta\|_1\|\theta\|_3^6)$, so the right hand side of \eqref{proof-Y3-result2} is much smaller than the right hand side of \eqref{proof-Y3-result1}. Together, we have
\[
\mathrm{Var}(Y_3)\leq C\alpha^2\|\theta\|^4\|\theta\|_3^6 = o(\|\theta\|^8). 
\]
This proves the claims of $Y_3$. 

Consider $Y_4$. We plug in $\delta_{ij}=\eta_i(\eta_j-\teta_j)+\eta_j(\eta_i-\teta_i)$ and the expression \eqref{tilde-eta}. It gives
\begin{align*}
Y_4 & = \sum_{i,j,k,\ell (dist)} \eta_i(\eta_j-\teta_j)\widetilde{\Omega}_{jk}\widetilde{\Omega}_{k\ell}W_{\ell i} + \sum_{i,j,k,\ell (dist)} \eta_j (\eta_i-\teta_i)\widetilde{\Omega}_{jk}\widetilde{\Omega}_{k\ell}W_{\ell i}\cr
&=\sum_{i,j,k,\ell (dist)} \eta_i\Bigl(-\frac{1}{\sqrt{v}}\sum_{s\neq j}W_{js}\Bigr)\widetilde{\Omega}_{jk}\widetilde{\Omega}_{k\ell}W_{\ell i} + \sum_{i,j,k,\ell (dist)}\eta_j \Bigl(-\frac{1}{\sqrt{v}}\sum_{s\neq i}W_{is}\Bigr) \widetilde{\Omega}_{jk}\widetilde{\Omega}_{k\ell}W_{\ell i}\cr
&= -\frac{1}{\sqrt{v}}\sum_{\substack{i,j,\ell (dist)\\s\neq j}}\Bigl(\sum_{k\notin\{i,j,\ell\}} \eta_i\widetilde{\Omega}_{jk}\widetilde{\Omega}_{k\ell}\Bigr) W_{js} W_{\ell i} - \frac{1}{\sqrt{v}} \sum_{\substack{i,\ell (dist)\\s\neq i}} \Bigl(\sum_{j,k\notin\{ i,\ell\}}\eta_j \widetilde{\Omega}_{jk}\widetilde{\Omega}_{k\ell}\Bigr) W_{is}W_{\ell i}\cr
&\equiv Y_{4a} + Y_{4b}. 
\end{align*}
First, we analyze $Y_{4a}$. When $(i,j,\ell)$ are distinct, $W_{js} W_{\ell i}$ has a mean zero. Therefore,
\[
\mathbb{E}[Y_{4a}] = 0. 
\]
To calculate the variance, we rewrite
\[
Y_{4a} =  -\frac{1}{\sqrt{v}}\sum_{\substack{i,j,\ell (dist)\\s\neq j}}\beta_{ij\ell} W_{js} W_{\ell i}, \qquad \mbox{where}\quad \beta_{ij\ell} = \sum_{k\notin\{i,j,\ell\}} \eta_i\widetilde{\Omega}_{jk}\widetilde{\Omega}_{k\ell}
\]
By \eqref{tOmega(ij)} and \eqref{eta-bound}, $|\beta_{ij\ell}|\leq C\sum_k\alpha^2\theta_i\theta_j\theta_k^2\theta_\ell\leq C\alpha^2 \|\theta\|^2\theta_i\theta_j\theta_\ell$. Also, for $W_{js}W_{\ell i}$ and $W_{j's'}W_{\ell' i'}$ to be correlated, there are only two cases: $(W_{js}, W_{\ell i})=(W_{j's'}, W_{\ell'i'})$ or $(W_{js}, W_{\ell i})=(W_{\ell'i'}, W_{j's'})$. Mimicking the argument in \eqref{proof-Y2-result1} or \eqref{proof-Y3-result2}, we can easily obtain that
\begin{align}  \label{proof-Y4-result1}
\mathrm{Var}(Y_{4a}) & \leq  \frac{C}{v}\sum_{\substack{i,j,\ell (dist)\\s\neq j}}\beta^2_{ij\ell}\cdot \mathrm{Var}( W_{js} W_{\ell i}) \cr
&\leq \frac{C}{\|\theta\|_1^2}\sum_{i,j,\ell,s}(\alpha^2\|\theta\|^2\theta_i\theta_j\theta_\ell)^2\cdot \theta_i\theta_j\theta_\ell\theta_s\cr
&\leq  \frac{C\alpha^4\|\theta\|^4\|\theta\|_3^9}{\|\theta\|_1}. 
\end{align}
Next, we analyze $Y_{4b}$. We re-write
\[
Y_{4b} = -\frac{1}{\sqrt{v}} \sum_{\substack{i,\ell (dist)\\s\neq i}} \beta_{i\ell} W_{is}W_{\ell i}, \qquad \mbox{where}\quad \beta_{i\ell} = \sum_{j,k\notin\{ i,\ell\}}\eta_j \widetilde{\Omega}_{jk}\widetilde{\Omega}_{k\ell}. 
\]
By separating the case of $s=\ell$ from the case of $s\neq \ell$, we have
\[
Y_{4b} = -\frac{1}{\sqrt{v}} \sum_{i,\ell (dist)} \beta_{i\ell} W^2_{\ell i}- \frac{1}{\sqrt{v}} \sum_{\substack{i,\ell (dist)\\s\notin\{ i,\ell\}}} \beta_{i\ell} W_{is}W_{\ell i}\equiv \widetilde{Y}_{4b} + Y_{4b}^*. 
\]
Only $\widetilde{Y}_{4b}$ has a nonzero mean. 
By \eqref{tOmega(ij)} and \eqref{eta-bound}, 
\[
|\beta_{i\ell}|\leq C\sum_{j,k}\alpha^2\theta^2_j\theta_k^2\theta_\ell\leq C\alpha^2\|\theta\|^4\theta_\ell.
\]
It follows that
\beq \label{proof-Y4-result2}
|\mathbb{E}[Y_{4b}]|= |\mathbb{E}[\widetilde{Y}_{4b}]|\leq \frac{C}{\|\theta\|_1}\sum_{i,\ell}(\alpha^2\|\theta\|^4\theta_\ell)\theta_i\theta_\ell \leq C\alpha^2\|\theta\|^6. 
\eeq
We now bound the variances of $\widetilde{Y}_{4b}$ and $Y^*_{4b}$. By direct calculations, 
\begin{align*}
\mathrm{Var}(\widetilde{Y}_{4b}) &= \frac{2}{v} \sum_{i,\ell (dist)} \beta^2_{i\ell}\cdot\mathrm{Var}( W^2_{i\ell}) \leq \frac{C}{\|\theta\|^2_1}\sum_{i,\ell} (\alpha^2\|\theta\|^4\theta_\ell)^2\cdot \theta_i\theta_\ell \leq \frac{C\alpha^4\|\theta\|^8\|\theta\|_3^3}{\|\theta\|_1},\cr
\mathrm{Var}(Y^*_{4b}) &\leq \frac{C}{v} \sum_{\substack{i,\ell (dist)\\s\notin\{ i,\ell\}}} \beta^2_{i\ell} \cdot \mathrm{Var}( W_{is}W_{\ell i}) \leq \frac{C}{\|\theta\|_1^2}\sum_{i,\ell,s} (\alpha^2\|\theta\|^4\theta_\ell)^2\cdot\theta_i^2\theta_\ell\theta_s\leq \frac{C\alpha^4\|\theta\|^{10}\|\theta\|_3^3}{\|\theta\|_1}. 
\end{align*}
Together, we have
\beq \label{proof-Y4-result3}
\mathrm{Var}(Y_{4b})\leq 2\mathrm{Var}(\widetilde{Y}_{4b})+ 2\mathrm{Var}(Y^*_{4b})  \leq \frac{C\alpha^4\|\theta\|^{10}\|\theta\|_3^3}{\|\theta\|_1}. 
\eeq
We combine the results of $Y_{4a}$ and $Y_{4b}$. Since $\|\theta\|_3^6\leq (\theta_{\max}\|\theta\|^2)^2=o(\|\theta\|^4)$, the right hand side of \eqref{proof-Y4-result2} dominates the right hand side of \eqref{proof-Y4-result1}. It follows that
\[
|\mathbb{E}[Y_4]|\leq C\alpha^2\|\theta\|^6=o(\alpha^4\|\theta\|^8), \quad \mathrm{Var}(Y_4)\leq \frac{C\alpha^4\|\theta\|^{10}\|\theta\|_3^3}{\|\theta\|_1}=o(\alpha^6\|\theta\|^8\|\theta\|_3^6). 
\]
Here, we explain the equalities. The first one is due to $\alpha^2\|\theta\|^2\to\infty$. To get the second equality, we compare $\mathrm{Var}(Y_4)$ with the order of $\alpha^6\|\theta\|^8\|\theta\|_3^6$. Note that $\frac{\|\theta\|^{10}\|\theta\|_3^3}{\|\theta\|_1}=\frac{\|\theta\|^6\|\theta\|_3^3}{\|\theta\|_1}\|\theta\|^4\leq \frac{\|\theta\|^6\|\theta\|_3^3}{\|\theta\|_1}\|\theta\|_1\|\theta\|_3^3\leq \|\theta\|^6\|\theta\|_3^6$. It follows that $\mathrm{Var}(Y_4)\leq C\alpha^4\|\theta\|^6\|\theta\|_3^6\ll C\alpha^6\|\theta\|^8\|\theta\|_3^6$, where the last inequality is due to $\alpha^2\|\theta\|^2\to\infty$. So far, we have proved all claims about $Y_4$. 

Consider $Y_5$. Recall that
\[
Y_5 = \sum_{i, j, k, \ell (dist)}\eta_i(\eta_j - \tilde{\eta}_j) \widetilde{\Omega}_{jk} W_{k \ell}  \widetilde{\Omega}_{\ell i} + \sum_{i, j, k, \ell (dist)} (\eta_i - \tilde{\eta}_i)\eta_j  \widetilde{\Omega}_{jk} W_{k \ell}  \widetilde{\Omega}_{\ell i}. 
\]
With relabeling of $(i,j,k,\ell)=(j', i', \ell', k')$, the second sum can be written as $\sum_{i', j', k', \ell' (dist)} (\eta_{j'}- \tilde{\eta}_{j'})\eta_{i'}  \widetilde{\Omega}_{i'\ell'} W_{\ell'k'}  \widetilde{\Omega}_{k' j'}$. This suggests that it is actually equal to the first sum above. Hence, 
\begin{align*}
Y_5 & = 2\sum_{i, j, k, \ell (dist)}\eta_i(\eta_j - \tilde{\eta}_j) \widetilde{\Omega}_{jk} W_{k \ell}  \widetilde{\Omega}_{\ell i}\cr
& = \sum_{i, j, k, \ell (dist)}\eta_i\Bigl(-\frac{2}{\sqrt{v}}\sum_{s\neq j}W_{js}\Bigr) \widetilde{\Omega}_{jk} W_{k \ell}  \widetilde{\Omega}_{\ell i}\cr
&= -\frac{2}{\sqrt{v}}\sum_{\substack{j,k,\ell (dist)\\s\neq j}}\Bigl(\sum_{i\notin\{j,k,\ell\}} \eta_i\widetilde{\Omega}_{jk}\widetilde{\Omega}_{\ell i}\Bigr)W_{js}W_{k\ell}\cr
&\equiv -\frac{2}{\sqrt{v}}\sum_{\substack{j,k,\ell (dist)\\s\neq j}}\beta_{jk\ell} W_{js}W_{k\ell}, \qquad \mbox{where}\quad \beta_{jk\ell} \equiv \sum_{i\notin\{j,k,\ell\}} \eta_i\widetilde{\Omega}_{jk}\widetilde{\Omega}_{\ell i}. 
\end{align*}
It is easy to see that $\mathbb{E}[W_{js}W_{k\ell}]=0$ when $(j,k,\ell)$ are distinct. Hence,
\[
\mathbb{E}[Y_5]=0. 
\]
By \eqref{tOmega(ij)} and \eqref{eta-bound}, $|\beta_{jk\ell}|\leq  C\sum_i \theta_i\cdot \alpha^2\theta_j\theta_k\theta_\ell\theta_i\leq C\alpha^2\|\theta\|^2\theta_j\theta_k\theta_\ell$. Similar to the argument in \eqref{proof-Y2-result1} or \eqref{proof-Y3-result2}, we can show that 
\begin{align*}
\mathrm{Var}(Y_5)&\leq \frac{C}{v}\sum_{\substack{j,k,\ell (dist)\\s\neq j}}\beta^2_{jk\ell}\cdot\mathrm{Var}( W_{js}W_{k\ell})\cr
&\leq \frac{C}{\|\theta\|_1^2}\sum_{j,k,\ell,s}(\alpha^2\|\theta\|^2 \theta_j\theta_k\theta_\ell)^2\theta_j\theta_s\theta_k\theta_\ell\cr
&\leq \frac{C\alpha^4\|\theta\|^4\|\theta\|_3^9}{\|\theta\|_1}. 
\end{align*}
Since $\|\theta\|_3^9=(\|\theta\|_3^3)^2\|\theta\|_3^3\leq (\theta_{\max}\|\theta\|^2)^2(\theta_{\max}^2\|\theta\|_1)=o(\|\theta\|^4\|\theta\|_1)$, the right hand side is $o(\|\theta\|^8)$. This proves the claims of $Y_5$. 

Consider $Y_6$. By definition and elementary calculations, 
\begin{align*}
Y_6 &= \sum_{i, j, k, \ell (dist)}\eta_i(\eta_j-\tilde{\eta}_j) \widetilde{\Omega}_{jk} \widetilde{\Omega}_{k \ell}  \widetilde{\Omega}_{\ell i} + \sum_{i, j, k, \ell (dist)} \eta_j (\eta_i - \tilde{\eta}_i)\widetilde{\Omega}_{jk} \widetilde{\Omega}_{k \ell}  \widetilde{\Omega}_{\ell i}\cr
&= 2\sum_{i, j, k, \ell (dist)}\eta_i(\eta_j - \tilde{\eta}_j) \widetilde{\Omega}_{jk} \widetilde{\Omega}_{k \ell}  \widetilde{\Omega}_{\ell i}\cr
&= 2\sum_{i, j, k, \ell (dist)}\eta_i\Bigl(-\frac{1}{\sqrt{v}}\sum_{s\neq j}W_{js}\Bigr) \widetilde{\Omega}_{jk} \widetilde{\Omega}_{k \ell}  \widetilde{\Omega}_{\ell i}\cr
&= -\frac{2}{\sqrt{v}}\sum_{j,s (dist)}\Bigl(\sum_{i,k,\ell (dist) \notin\{j\}} \eta_i\widetilde{\Omega}_{jk}\widetilde{\Omega}_{k\ell}\widetilde{\Omega}_{\ell i}\Bigr)W_{js}. 
\end{align*}
Here, to get the second line above, we relabeled $(i,j,k,\ell)=(j', i', \ell', k')$ in the second sum and found out the two sums are equal; the third line is from \eqref{tilde-eta}. 
We immediately see that
\[
\mathrm{E}[Y_6]=0. 
\]
By \eqref{tOmega(ij)} and \eqref{eta-bound}, 
\begin{align*}
\Bigl| \sum_{i,k,\ell (dist) \notin\{j\}} \eta_i\widetilde{\Omega}_{jk}\widetilde{\Omega}_{k\ell}\widetilde{\Omega}_{\ell i} \Bigr|
& \leq \sum_{i,k,\ell}C\theta_i\cdot \alpha^3 \theta_j\theta_k^2\theta_\ell^2\theta_i\leq C\alpha^3\|\theta\|^6\theta_j. 
\end{align*}
It follows that
\begin{align*}
\mathrm{Var}(Y_6)& = \frac{8}{v}\sum_{j,s (dist)}\Bigl(\sum_{i,k,\ell (dist) \notin\{j\}} \eta_i\widetilde{\Omega}_{jk}\widetilde{\Omega}_{k\ell}\widetilde{\Omega}_{\ell i}\Bigr)^2\cdot \mathrm{Var}(W_{js})\cr
&\leq \frac{C}{\|\theta\|_1^2}\sum_{j,s}(\alpha^3\|\theta\|^6\theta_j)^2\theta_j\theta_s\cr
&\leq \frac{C\alpha^6\|\theta\|^{12}\|\theta\|_3^3}{\|\theta\|_1}. 
\end{align*}
Since $\|\theta\|^4\leq \|\theta\|_1\|\theta\|_3^3$, the variance is bounded by $C\alpha^6\|\theta\|^8\|\theta\|_3^6$. This proves the claims of $Y_6$.

\subsubsection{Proof of Lemma~\ref{lem:ProxySgnQ(b)-null}}
It suffices to prove the claims for each of $Z_1$ and $Z_2$; then, the claims of $U_b$ follow immediately.

We first analyze $Z_1$. 
Plugging $\delta_{ij}=\eta_i(\eta_j - \tilde{\eta}_j)+\eta_j(\eta_i - \tilde{\eta}_i)$ into the definition of $Z_1$ gives  
\begin{align*}
Z_1 &= \sum_{i,j,k,\ell (dist)}\eta_i (\eta_j - \tilde{\eta}_j)\eta_j (\eta_k - \tilde{\eta}_k) W_{k\ell}W_{\ell i}+ \sum_{i,j,k,\ell (dist)}\eta_i (\eta_j - \tilde{\eta}_j)^2\eta_k W_{k\ell}W_{\ell i}\cr
&\quad +\sum_{i,j,k,\ell (dist)}(\eta_i - \tilde{\eta}_i)\eta^2_j (\eta_k - \tilde{\eta}_k) W_{k\ell}W_{\ell i} +  \sum_{i,j,k,\ell (dist)}(\eta_i - \tilde{\eta}_i)\eta_j (\eta_j-\tilde{\eta}_j)\eta_k W_{k\ell}W_{\ell i}.
\end{align*}
In the last term above, if we relabel $(i,j,k,\ell)=(k',j', i',\ell')$, it becomes $\sum_{i',j',k',\ell' (dist)}(\eta_{k'} - \tilde{\eta}_{k'})\eta_{j'} (\eta_{j'}-\tilde{\eta}_{j'})\eta_{i'} W_{i'\ell'}W_{\ell' k'}$. This shows that the last sum equals to the first sum. Therefore, 
\begin{eqnarray}  \label{proof-Z1-decompose}
Z_1 &=& \sum_{i,j,k,\ell (dist)}\eta_i (\eta_j - \tilde{\eta}_j)^2\eta_k W_{k\ell}W_{\ell i} \cr
&&+ 2\sum_{i,j,k,\ell (dist)}\eta_i (\eta_j - \tilde{\eta}_j)\eta_j (\eta_k - \tilde{\eta}_k) W_{k\ell}W_{\ell i}\cr
&&+ \sum_{i,j,k,\ell (dist)}(\tilde{\eta}_i-\eta_i)\eta^2_j (\tilde{\eta}_k-\eta_k) W_{k\ell}W_{\ell i}\cr
&\equiv& Z_{1a} + Z_{1b} + Z_{1c}. 
\end{eqnarray}
Below, we compute the means and variances of $Z_{1a}$-$Z_{1c}$. 

First, we study $Z_{1a}$. When $(i,j,k,\ell)$ are distinct, $W_{k\ell}W_{\ell i}$ has a mean zero and is independent of $(\tilde{\eta}_j-\eta_j)^2$, so $\mathbb{E}[ (\eta_j - \tilde{\eta}_j)^2W_{k\ell}W_{\ell i} ]=0$. It follows that 
\[
\mathbb{E}[Z_{1a}]=0. 
\]
To bound the variance of $Z_{1a}$, we use \eqref{tilde-eta} to re-write
\begin{eqnarray*}
Z_{1a} &=& \sum_{i,j,k,\ell (dist)}\eta_i \Bigl(-\frac{1}{\sqrt{v}}\sum_{s\neq j}W_{js}\Bigr)\Bigl(-\frac{1}{\sqrt{v}} \sum_{t\neq j}W_{jt}\Bigr) \eta_k W_{k\ell}W_{\ell i}\cr
&=& \frac{1}{v} \sum_{\substack{i,j,k,\ell (dist)\\s,t\notin\{j\}}}\eta_i \eta_k W_{js}W_{jt}W_{k\ell}W_{\ell i}\cr
&=& \frac{1}{v} \sum_{\substack{i,j,k,\ell (dist)\\s\notin\{j\}}}\eta_i \eta_k W^2_{js}W_{k\ell}W_{\ell i} + \frac{1}{v} \sum_{\substack{i,j,k,\ell (dist)\\s,t (dist)\notin\{j\}}}\eta_i \eta_k W_{js}W_{jt}W_{k\ell}W_{\ell i}\cr
&\equiv& \widetilde{Z}_{1a} + Z^*_{1a}. 
\end{eqnarray*}
We first bound the variance of $\widetilde{Z}_{1a}$. It is seen that
\[
\mathrm{Var}(\widetilde{Z}_{1a}) = \frac{1}{v^2} \sum_{\substack{i,j,k,\ell (dist), s\notin\{j\}\\i',j',k',\ell' (dist), s'\notin\{j'\}}}\eta_i \eta_k \eta_{i'}\eta_{k'}\cdot \mathbb{E}[W^2_{js}W_{k\ell}W_{\ell i}\cdot W^2_{j's'}W_{k'\ell'}W_{\ell' i'}]. 
\]
The summand is nonzero only if $\ell'=\ell$ and $\{k',i'\}=\{k,i\}$. We also note that, if we switch $i'$ and $k'$, the summand remains unchanged. So, it suffices to consider the case of $\ell'=\ell$ and $(k',i')=(k,i)$. 
By \eqref{eta-bound} and elementary calculations,
\begin{align*}
& \eta_i \eta_k \eta_{i'}\eta_{k'}\cdot \mathbb{E}[W^2_{js}W_{k\ell}W_{\ell i}\cdot W^2_{j's'}W_{k'\ell'}W_{\ell' i'}]\cr
=\; & \begin{cases}
\eta^2_i \eta^2_k\, \mathbb{E}[W^4_{js}W^2_{k\ell}W^2_{\ell i}]\leq C\theta^3_i\theta_j\theta^3_k\theta^2_\ell\theta_s, &\mbox{if }(\ell',k',i')=(\ell,k,i),\, \{j',s'\}= \{j,s\};\\
 \eta^2_i \eta^2_k\, \mathbb{E}[W^2_{js}W^2_{k\ell}W^2_{\ell i}W^2_{j's'}]\leq C\theta^3_i\theta_j\theta^3_k\theta^2_\ell\theta_s \theta_{j'}\theta_{s'}, &\mbox{if }(\ell',k',i')=(\ell,k,i),\, \{j',s'\}\neq \{j,s\};\\
 0, &\mbox{otherwise}. 
\end{cases}
\end{align*}
It follows that
\begin{align*}
\mathrm{Var}(\widetilde{Z}_{1a})&\leq \frac{C}{\|\theta\|_1^4}\Bigl(\sum_{i,j,k,\ell,s}\theta^3_i\theta_j\theta^3_k\theta^2_\ell\theta_s + \sum_{i,j,k,\ell,s,j',s'}\theta^3_i\theta_j\theta^3_k\theta^2_\ell\theta_s \theta_{j'}\theta_{s'}\Bigr)\cr
&\leq \frac{C}{\|\theta\|_1^4}\bigl( \|\theta\|^2\|\theta\|_3^6\|\theta\|_1^2 + \|\theta\|^2\|\theta\|_3^6\|\theta\|_1^4  \bigr)\cr
&\leq C\|\theta\|^2\|\theta\|_3^6. 
\end{align*}
We then bound the variance of $Z^*_{1a}$. Note that
\begin{align*}
& \eta_i \eta_k \eta_{i'}\eta_{k'}\cdot \mathbb{E}[W_{js}W_{jt}W_{k\ell}W_{\ell i}\cdot W_{j's'}W_{j't'}W_{k'\ell'}W_{\ell' i'}]\cr
=\; & \begin{cases}
\eta^2_i \eta^2_k\, \mathbb{E}[W^2_{js}W^2_{jt}W^2_{k\ell}W^2_{\ell i}]\leq C\theta^3_i\theta^2_j\theta^3_k\theta^2_\ell\theta_s\theta_t, &\mbox{if }(j', \ell')=(j,\ell), \{s', t'\}=\{ s, t\}, \{k',i'\}=\{k,i\};\\
 \eta_i \eta_k\eta_s\eta_t\, \mathbb{E}[W^2_{js}W^2_{jt}W^2_{k\ell}W^2_{\ell i}] \leq C\theta^2_i\theta^2_j\theta^2_k\theta^2_\ell\theta^2_s\theta^2_t, &\mbox{if }(j', \ell')=(\ell,j), \{s', t'\}=\{k,i\}, \{k', i'\}=\{s,t\};\\
 0, &\mbox{otherwise}. 
\end{cases}
\end{align*}
It follows that
\begin{align*}
\mathrm{Var}(Z^*_{1a})&\leq \frac{C}{\|\theta\|_1^4}\Bigl(\sum_{i,j,k,\ell,s,t}\theta^3_i\theta^2_j\theta^3_k\theta^2_\ell\theta_s\theta_t + \sum_{i,j,k,\ell,s,t}\theta^2_i\theta^2_j\theta^2_k\theta^2_\ell\theta^2_s\theta^2_t\Bigr)\cr
&\leq \frac{C}{\|\theta\|_1^4}\bigl( \|\theta\|^4\|\theta\|_3^6\|\theta\|_1^2 + \|\theta\|^{12}  \bigr)\cr
&\leq \frac{C\|\theta\|^4\|\theta\|_3^6}{\|\theta\|_1^2},
\end{align*}
where the last inequality is because of $\|\theta\|^{12}=\|\theta\|^4(\|\theta\|^4)^2\leq \|\theta\|^4(\|\theta\|_1\|\theta\|_3^3)^2=\|\theta\|^4\|\theta\|_3^6\|\theta\|_1^2$. Combining the above gives
\beq \label{proof-Z1-result1}
\mathrm{Var}(Z_{1a})\leq 2\mathrm{Var}(\widetilde{Z}_{1a}) + 2\mathrm{Var}(Z^*_{1a})\leq C\|\theta\|^2\|\theta\|_3^6. 
\eeq

Second, we study $Z_{1b}$. Since $(\eta_j-\teta_j)$, $(\eta_k-\teta_k)W_{k\ell}$ and $W_{\ell i}$ are independent of each other, each summand in $Z_{1b}$ has a zero mean. It follows that
\[
\mathbb{E}[Z_{1b}]=0. 
\]
We now compute its variance. By direct calculations, 
\begin{align*}
Z_{1b} &= 2\sum_{i,j,k,\ell (dist)}\eta_i \Bigl(-\frac{1}{\sqrt{v}}\sum_{s\neq j}W_{js}\Bigr)\eta_j \Bigl(-\frac{1}{\sqrt{v}} \sum_{t\neq k}W_{kt}\Bigr) W_{k\ell}W_{\ell i}\cr
&= \frac{2}{v}\sum_{\substack{i,j,k,\ell (dist)\\s\neq j,t\neq k}}\eta_i \eta_j W_{js}W_{kt}W_{k\ell}W_{\ell i}\cr
&= \frac{2}{v}\sum_{\substack{i,j,k,\ell (dist)\\s\neq j}}\eta_i \eta_j W_{js}W^2_{k\ell}W_{\ell i} + \frac{2}{v}\sum_{\substack{i,j,k,\ell (dist)\\s\neq j,t\notin\{ k,\ell\}}}\eta_i \eta_jW_{js}W_{kt}W_{k\ell}W_{\ell i}\cr
&\equiv \widetilde{Z}_{1b}+ Z^*_{1b}. 
\end{align*}
We first bound the variance of $\widetilde{Z}_{1b}$. Note that
\[
\mathrm{Var}(\widetilde{Z}_{1b}) = \frac{4}{v}\sum_{\substack{i,j,k,\ell (dist), s\neq j\\i',j',k',\ell' (dist), s'\neq j'}}\eta_i \eta_j\eta_{i'}\eta_{j'}\cdot\mathbb{E}[W_{js}W^2_{k\ell}W_{\ell i}\cdot W_{j's'}W^2_{k'\ell'}W_{\ell' i'}]. 
\]
For this summand to be nonzero, there are only two cases. In the first case, $(W_{js}, W_{\ell i})$ are paired with $(W_{j's'}, W_{\ell' i'})$. It follows that 
\[
\eta_i \eta_j\eta_{i'}\eta_{j'}\cdot\mathbb{E}[W_{js}W^2_{k\ell}W_{\ell i}W_{j's'}W^2_{k'\ell'}W_{\ell' i'}] = \eta_i \eta_j\eta_{i'}\eta_{j'}\cdot 
\mathbb{E}[W^2_{js}W^2_{k\ell}W^2_{\ell i}W^2_{k'\ell'}]. 
\]
This happens only if (i) $\{j',s'\}=\{j,s\}$ and $\{\ell',i'\}=\{\ell,i\}$, or (ii) $\{j',s'\}=\{\ell, i\}$ and $\{\ell',i'\}=\{j,s\}$. By \eqref{eta-bound} and elementary calculations,
\begin{align*}
& \eta_i \eta_j\eta_{i'}\eta_{j'}\cdot\mathbb{E}[W_{js}W^2_{k\ell}W_{\ell i}\cdot W_{j's'}W^2_{k'\ell'}W_{\ell' i'}]\cr
=\; & \begin{cases}
\eta^2_i \eta^2_j\cdot\mathbb{E}[W^2_{js}W_{\ell i}^2 W^2_{k\ell}W^2_{k'\ell}]\leq C\theta^3_i\theta^3_j\theta_k\theta^3_\ell\theta_s\theta_{k'}, &\mbox{if }(j',s')=(j,s), (\ell',i')=(\ell,i);\\
\eta_i\eta^2_j\eta_\ell \cdot\mathbb{E}[W^2_{js}W_{\ell i}^2 W^2_{k\ell}W^2_{k'i}]\leq C\theta^3_i\theta^3_j\theta_k\theta^3_\ell\theta_s\theta_{k'}, &\mbox{if }(j',s')=(j,s), (\ell',i')=(i,\ell);\\
\eta^2_i\eta_j\eta_s\cdot\mathbb{E}[W^2_{js}W_{\ell i}^2 W^2_{k\ell}W^2_{k'\ell}]\leq C\theta^3_i\theta^2_j\theta_k\theta^3_\ell\theta^2_s\theta_{k'}, &\mbox{if }(j',s')=(s,j), (\ell',i')=(\ell,i);\\
\eta_i\eta_j\eta_\ell\eta_s\cdot\mathbb{E}[W^2_{js}W_{\ell i}^2 W^2_{k\ell}W^2_{k'i}]\leq C\theta^3_i\theta^2_j\theta_k\theta^3_\ell\theta^2_s\theta_{k'}, &\mbox{if }(j',s')=(s,j), (\ell',i')=(i,\ell);\\
\eta_i \eta_j\eta_\ell\eta_s\cdot\mathbb{E}[W^2_{js}W_{\ell i}^2 W^2_{k\ell}W^2_{k'j}]\leq C\theta^2_i\theta^3_j\theta_k\theta^3_\ell\theta^2_s\theta_{k'}, &\mbox{if }(j',s')=(\ell,i), (\ell',i')=(j,s);\\
\eta_i \eta^2_j\eta_\ell  \cdot\mathbb{E}[W^2_{js}W_{\ell i}^2 W^2_{k\ell}W^2_{k's}]\leq C\theta^2_i\theta^3_j\theta_k\theta^3_\ell\theta^2_s\theta_{k'}, &\mbox{if }(j',s')=(\ell,i), (\ell',i')=(s,j);\\
\eta^2_i \eta_j\eta_s\cdot\mathbb{E}[W^2_{js}W_{\ell i}^2 W^2_{k\ell}W^2_{k'j}]\leq C\theta^3_i\theta^3_j\theta_k\theta^2_\ell\theta^2_s\theta_{k'}, &\mbox{if }(j',s')=(i, \ell), (\ell',i')=(j,s);\\
\eta^2_i \eta^2_j\cdot\mathbb{E}[W^2_{js}W_{\ell i}^2 W^2_{k\ell}W^2_{k's}]\leq C\theta^3_i\theta^3_j\theta_k\theta^2_\ell\theta^2_s\theta_{k'}, &\mbox{if }(j',s')=(i, \ell), (\ell',i')=(s,j);\\
0, &\mbox{otherwise}. 
\end{cases}
\end{align*}
The upper bound on the right hand side only has two types $C\theta^3_i\theta^3_j\theta_k\theta^3_\ell\theta_s\theta_{k'}$ and $C\theta^3_i\theta^2_j\theta_k\theta^3_\ell\theta^2_s\theta_{k'}$. The contribution of this case to $\mathrm{Var}(\widetilde{Z}_{1b})$ is
\begin{align*}
&\leq \frac{C}{v^2}\Bigl( \sum_{i,j,k,\ell,s,k'}\theta^3_i\theta^3_j\theta_k\theta^3_\ell\theta_s\theta_{k'} + \sum_{i,j,k,\ell,s,k'}\theta^3_i\theta^2_j\theta_k\theta^3_\ell\theta^2_s\theta_{k'}\Bigr)\cr
&\leq \frac{C}{\|\theta\|_1^4}\bigl(\|\theta\|_3^9\|\theta\|^3_1 + \|\theta\|^4\|\theta\|_3^6\|\theta\|^2_1 \bigr)\cr
&\leq \frac{C\|\theta\|_3^9}{\|\theta\|_1}.
\end{align*}
In the second case, $\{W_{js}, W_{k\ell}, W_{\ell i}\}$ and  $\{W_{j's'}, W_{k'\ell'}, W_{\ell' i'}\}$ are two sets of same variables. Then,
\[
\eta_i \eta_j\eta_{i'}\eta_{j'}\cdot\mathbb{E}[W_{js}W^2_{k\ell}W_{\ell i}W_{j's'}W^2_{k'\ell'}W_{\ell' i'}] = \eta_i \eta_j\eta_{i'}\eta_{j'}\cdot 
\mathbb{E}[W^3_{js}W^3_{k\ell}W^3_{\ell i}]. 
\]
This can only happen if $\ell'=\ell$, $\{i',k'\}=\{i,k\}$, and $\{j',s'\}=\{j,s\}$. By \eqref{eta-bound} and elementary calculations,
\begin{align*}
& \eta_i \eta_j\eta_{i'}\eta_{j'}\cdot\mathbb{E}[W_{js}W^2_{k\ell}W_{\ell i}\cdot W_{j's'}W^2_{k'\ell'}W_{\ell' i'}]\cr
=\; & \begin{cases}
\eta^2_i \eta^2_j\cdot\mathbb{E}[W^3_{js}W^3_{\ell i} W^3_{k\ell}]\leq C\theta^3_i\theta^3_j\theta_k\theta^2_\ell\theta_s, &\mbox{if }\ell'=\ell,\, (i',k')=(i,k),\, (j',s')=(j,s);\\
\eta^2_i \eta_j\eta_s\cdot\mathbb{E}[W^3_{js}W^3_{\ell i} W^3_{k\ell}]\leq C\theta^3_i\theta^2_j\theta_k\theta^2_\ell\theta^2_s, &\mbox{if }\ell'=\ell,\, (i',k')=(i,k),\, (j',s')=(s,j);\\
\eta_i\eta_k \eta^2_j\cdot\mathbb{E}[W^3_{js}W^3_{\ell i} W^3_{k\ell}]\leq C\theta^2_i\theta^3_j\theta^2_k\theta^2_\ell\theta_s, &\mbox{if }\ell'=\ell,\, (i',k')=(k,i),\, (j',s')=(j,s);\\
\eta_i\eta_k \eta_j\eta_s\cdot\mathbb{E}[W^3_{js}W^3_{\ell i} W^3_{k\ell}]\leq C\theta^2_i\theta^2_j\theta^2_k\theta^2_\ell\theta^2_s, &\mbox{if }\ell'=\ell,\, (i',k')=(i,k),\, (j',s')=(s,j);\\
0, &\mbox{otherwise}. 
\end{cases}
\end{align*}
The upper bound on the right hand side has three types, and the contribution of this case to $\mathrm{Var}(\widetilde{Z}_{1b})$ is
\begin{align*}
&\leq \frac{C}{v^2}\Bigl( \sum_{i,j,k,\ell,s}\theta^3_i\theta^3_j\theta_k\theta^2_\ell\theta_s + \sum_{i,j,k,\ell,s}\theta^3_i\theta^2_j\theta_k\theta^2_\ell\theta^2_s + \sum_{i, j,k,\ell,s} \theta^2_i\theta^2_j\theta^2_k\theta^2_\ell\theta^2_s \Bigr)\cr
&\leq \frac{C}{\|\theta\|_1^4}\bigl(\|\theta\|^2\|\theta\|_3^6\|\theta\|^2_1 + \|\theta\|^6\|\theta\|_3^3\|\theta\|_1+\|\theta\|^{10} \bigr)\cr
&\leq \frac{C\|\theta\|^2\|\theta\|_3^6}{\|\theta\|^2_1},
\end{align*}
where we use $\|\theta\|^4\leq \|\theta\|_1\|\theta\|_3^3$ (Cauchy-Schwarz) in the last line. 
It is seen that the contribution of the first case is dominating, and so
\[
\mathrm{Var}(\widetilde{Z}_{1b})\leq \frac{C\|\theta\|_3^9}{\|\theta\|_1}. 
\]
We then bound the variance of $Z^*_{1b}$. Note that
\[
\mathrm{Var}(Z^*_{1b}) = \frac{4}{v^2}\sum_{\substack{i,j,k,\ell (dist), s\neq j,t\notin\{ k,\ell\}\\ i',j',k',\ell' (dist), s'\neq j',t'\notin\{ k',\ell'\}}}\eta_i \eta_j\eta_{i'}\eta_{j'}\cdot\mathbb{E}[W_{js}W_{kt}W_{k\ell}W_{\ell i}\cdot W_{j's'}W_{k't'}W_{k'\ell'}W_{\ell' i'}]. 
\]
For the summand to be nonzero, all $W$ terms have to be perfectly matched, so that the expectation in the summand becomes 
\[
\mathbb{E}[W_{js}W_{kt}W_{k\ell}W_{\ell i}\cdot W_{j's'}W_{k't'}W_{k'\ell'}W_{\ell' i'}]= \mathbb{E}[W^2_{js}W^2_{kt}W^2_{k\ell}W^2_{\ell i} ]\leq C\theta_i\theta_j\theta_k^2\theta_\ell^2\theta_s\theta_t. 
\] 
For this perfect match to happen, we need $(t', k', \ell', i')=(t, k, \ell, i)$ or $(t', k', \ell', i')=(i, \ell, k, t)$, as well as $\{j', s'\}=\{j,s\}$. This implies that, $i'$ can only take values in $\{i, t\}$ and $j'$ can only take values in $\{j, s\}$. It follows that $\eta_i \eta_j\eta_{i'}\eta_{j'}$ belongs to one of the following cases: 
\begin{align*}
&\eta_i\eta_j(\eta_i\eta_j)\leq C\theta_i^2\theta_j^2, &&\eta_i\eta_j(\eta_i\eta_s)=C\theta^2_i\theta_j\theta_s,\\
& \eta_i\eta_j(\eta_t\eta_j)\leq C\theta_i\theta_j^2\theta_t, && \eta_i\eta_j(\eta_t\eta_s)\leq C\theta_i\theta_j\theta_t\theta_s.
\end{align*}
Combining the above gives
\begin{align*}
\mathrm{Var}(Z^*_{1b}) &\leq  \frac{C}{v^2}\sum_{i,j,k,\ell,s,t}(\theta_i^2\theta_j^2+ \theta^2_i\theta_j\theta_s + \theta_i\theta_j^2\theta_t+\theta_i\theta_j\theta_t\theta_s)\cdot \theta_i\theta_j\theta_k^2\theta_\ell^2\theta_s\theta_t\cr
&\leq \frac{C}{\|\theta\|_1^4}\bigl( \|\theta\|^4 \|\theta\|^6_3\|\theta\|_1^2+2\|\theta\|^8 \|\theta\|^3_3\|\theta\|_1+\|\theta\|^{12} \bigr)\cr
&\leq \frac{C\|\theta\|^4 \|\theta\|^6_3}{\|\theta\|_1^2}.  
\end{align*}
We combine the variances of $\widetilde{Z}_{1b}$ and $Z^*_{1b}$. Since $\|\theta\|^4\leq \|\theta\|_1\|\theta\|_3^3$, the variance of $\widetilde{Z}_{1b}$ dominates. It follows that
\beq \label{proof-Z1-result2}
\mathrm{Var}(Z_{1b})\leq 2\mathrm{Var}(\widetilde{Z}_{1b}) + 2\mathrm{Var}(Z^*_{1b})\leq \frac{C\|\theta\|_3^9}{\|\theta\|_1}.  
\eeq

Third, we study $Z_{1c}$. It is seen that
\begin{align*}
Z_{1c} &= \sum_{i,j,k,\ell (dist)}\Bigl(-\frac{1}{\sqrt{v}}\sum_{s\neq i}W_{is}\Bigr)\eta^2_j \Bigl(-\frac{1}{\sqrt{v}} \sum_{t\neq k}W_{kt}\Bigr) W_{k\ell}W_{\ell i}\cr
&= \frac{1}{v}\sum_{\substack{i,k,\ell (dist)\\ s\neq i, t\neq k}}\Bigl(\sum_{j\notin\{i,k,\ell\}}\eta^2_j\Bigr) W_{is}W_{kt}W_{k\ell}W_{\ell i}\cr
&\equiv \frac{1}{v}\sum_{\substack{i,k,\ell (dist)\\ s\neq i, t\neq k}}\beta_{ik\ell}W_{is}W_{kt}W_{k\ell}W_{\ell i},
\end{align*}
where 
\beq  \label{proof-Z1-beta}
\beta_{ik\ell}\equiv\sum_{j\notin\{i,k,\ell\}}\eta^2_j\leq C\sum_j\theta_j^2\leq C\|\theta\|^2. 
\eeq
We divide all summands into four groups: (i) $s=t=\ell$; (ii) $s=\ell$, $t\neq \ell$; (iii) $s\neq \ell$, $t=\ell$; (iv) $s\neq \ell$, $t\neq \ell$. It yields that
\begin{align*}
Z_{1c}&= \frac{1}{v}\sum_{i,k,\ell (dist)}\beta_{ik\ell}  W^2_{k\ell}W^2_{\ell i} + \frac{1}{v}\sum_{\substack{i,k,\ell (dist)\\t\neq \{k,\ell\}}}\beta_{ik\ell}  W_{kt}W_{k\ell}W^2_{\ell i} \cr
&\qquad + \frac{1}{v}\sum_{\substack{i,k,\ell (dist)\\s\notin \{i,\ell\}}}\beta_{ik\ell}  W_{is}W^2_{k\ell}W_{\ell i}+ \frac{1}{v}\sum_{\substack{i,k,\ell (dist)\\s\notin\{ i,\ell\},t\notin\{ k,\ell\}}} \beta_{ik\ell}    W_{is}W_{kt}W_{k\ell}W_{\ell i}.
\end{align*}
In the third sum, if we relabel $(i,k,\ell,s)=(k',i',\ell',t')$, it then has the form of $\sum_{i',k',\ell' (dist), t'\notin \{k',\ell'\}}\beta_{k'i'\ell'}  W_{k't'}W^2_{i'\ell'}W_{\ell'k'}$. This shows that this sum equals to the second sum. We thus have
\begin{align*}
Z_{1c}&= \frac{1}{v}\sum_{i,k,\ell (dist)}\beta_{ik\ell}  W^2_{k\ell}W^2_{\ell i} + \frac{2}{v}\sum_{\substack{i,k,\ell (dist)\\t\neq \{k,\ell\}}}\beta_{ik\ell}  W_{kt}W_{k\ell}W^2_{\ell i} \cr
&\qquad + \frac{1}{v}\sum_{\substack{i,k,\ell (dist)\\s\notin\{ i,\ell\}, t\notin\{ k,\ell\}}} \beta_{ik\ell}    W_{is}W_{kt}W_{k\ell}W_{\ell i}\cr
&\equiv \widetilde{Z}_{1c} + Z^*_{1c} + Z^{\dag}_{1c}. 
\end{align*}
Among all three terms, only $\widetilde{Z}_{1c}$ has a nonzero mean. It follows that
\begin{align*}
\mathbb{E}[Z_{1c}] = \mathbb{E}[\widetilde{Z}_{1c}] &= \frac{1}{v}\sum_{i,k,\ell (dist)}\beta_{ik\ell}\Omega_{k\ell}(1-\Omega_{k\ell})\Omega_{\ell i}(1-\Omega_{\ell i})\cr
&=  \frac{1}{v}\sum_{i,k,\ell (dist)}\beta_{ik\ell}\Omega_{k\ell}\Omega_{\ell i}[1+O(\theta^2_{\max})]. 
\end{align*}
Under the null hypothesis, $\Omega_{ij}=\theta_i\theta_j$. It follows that $\eta_j=\frac{\theta_j}{\sqrt{v}}\sum_{i:i\neq j}\theta_i =[1+o(1)] \frac{\theta_j \|\theta\|_1}{\sqrt{v}}$ and that $\beta_{ik\ell}=\sum_{j\notin\{i,k,\ell\}}\eta_j^2= [1+o(1)] \frac{\|\theta\|^2_1}{v}\sum_{j\notin\{i,k,\ell\}}\theta_j^2=[1+o(1)]\frac{\|\theta\|^2_1\|\theta\|^2}{v}$. Additionally, $v=\sum_{i\neq j}\theta_i\theta_j=\|\theta\|_1^2\cdot [1+o(1)]$. As a result,
\begin{align} \label{proof-Z1-result3}
\mathbb{E}[Z_{1c}] & =\frac{1}{v}\sum_{i,k,\ell (dist)}[1+o(1)]\frac{\|\theta\|_1^2\|\theta\|^2}{v}\cdot \theta_k\theta^2_\ell\theta_i \cr
&=[1+o(1)]\cdot \frac{\|\theta\|_1^2\|\theta\|^2}{v^2}\sum_{i,k,\ell (dist)}\theta_k\theta^2_\ell\theta_i\cr
&=[1+o(1)]\cdot \frac{\|\theta\|_1^2\|\theta\|^2}{\|\theta\|_1^4} \bigl[\|\theta\|_1^2\|\theta\|^2 - O(\|\theta\|^4+\|\theta\|_1\|\theta\|_3^3)\bigr]\cr
&= [1+o(1)]\cdot\|\theta\|^4,
\end{align}
where in the last line we have used $\|\theta\|^2=o(\|\theta\|_1)$, $\|\theta\|_3^3=o(\|\theta\|_1)$ and $\|\theta\|_1\to\infty$. 
We then bound the variance of $Z_{1c}$ by studying the variance of each of the three variables, $\widetilde{Z}_{1c}$, $Z^*_{1c}$ and $Z^{\dag}_{1c}$. Consider $\widetilde{Z}_{1c}$ first. For $W^2_{k\ell}W^2_{\ell i}$ and $W^2_{k'\ell'}W^2_{\ell' i'}$ to be correlated, it has to be the case of either $\{k',\ell'\}=\{k,\ell\}$ or $\{i',\ell'\}=\{i,\ell\}$. By symmetry between $k$ and $i$ in the expression, it suffices to consider $\{k',\ell'\}=\{k,\ell\}$. Direct calculations show that 
\[
\mathrm{Cov}(W^2_{k\ell}W^2_{\ell i},\, W^2_{k'\ell'}W^2_{\ell' i'})\leq
\begin{cases}
\mathbb{E}[W^4_{k\ell}W^4_{\ell i}]\leq C\theta_k\theta_\ell^2\theta_i, &\mbox{if }(k',\ell')=(k,\ell),\, i'=i;\\
\mathbb{E}[W^4_{k\ell}W^2_{\ell i}W^2_{\ell i'}]\leq C\theta_k\theta^3_\ell\theta_i\theta_{i'}, &\mbox{if }(k',\ell')=(k,\ell),\, i'\neq i;\\
\mathbb{E}[W^4_{k\ell}W^2_{\ell i}W^2_{ki}] \leq C\theta^2_k\theta_\ell^2\theta^2_i, &\mbox{if }(k',\ell')=(\ell,k),\, i'=i;\\
\mathbb{E}[W^4_{k\ell}W^2_{\ell i}W^2_{k i'}]\leq C\theta^2_k\theta^2_\ell\theta_i\theta_{i'}, &\mbox{if }(k',\ell')=(\ell,k),\, i'\neq i;\\
0, &\mbox{otherwise}. 
\end{cases}
\]
Combining it with \eqref{proof-Z1-beta} and the fact of $v\geq C^{-1}\|\theta\|_1^2$, we have
\begin{align*}
\mathrm{Var}(\widetilde{Z}_{1c})&\leq \frac{C\|\theta\|^4}{\|\theta\|_1^4}\Bigl(\sum_{i,k,\ell} \theta_k\theta_\ell^2\theta_i + \sum_{i,k,\ell,i'}\theta_k\theta^3_\ell\theta_i\theta_{i'} + \sum_{i,k,\ell}\theta^2_k\theta_\ell^2\theta^2_i + \sum_{i,k,\ell,i'} \theta^2_k\theta^2_\ell\theta_i\theta_{i'}\Bigr)\cr
&\leq \frac{C\|\theta\|^4}{\|\theta\|_1^4}\bigl( \|\theta\|^2\|\theta\|_1^2 + \|\theta\|_3^3\|\theta\|_1^3 + \|\theta\|^6 + \|\theta\|^4\|\theta\|_1^2 \bigr)\cr
&\leq\frac{C\|\theta\|^4\|\theta\|_3^3}{\|\theta\|_1}.
\end{align*}
Consider $Z^*_{1c}$. By direct calculations,
\begin{align*}
& \mathbb{E}[W_{kt}W_{k\ell}W^2_{\ell i}W_{k't'}W_{k'\ell'}W^2_{\ell' i'}]\cr
=\; & \begin{cases}
\mathbb{E}[W^2_{kt}W^2_{k\ell}W^4_{\ell i}]\leq C\theta_i\theta_k^2\theta_\ell^2\theta_t, & \mbox{if }(k', t',\ell')=(k, t,\ell),\, i= i';\\
\mathbb{E}[W^2_{kt}W^2_{k\ell}W^2_{\ell i}W^2_{\ell i'}]\leq C\theta_i\theta_k^2\theta_\ell^3\theta_t\theta_{i'}, & \mbox{if } (k', t',\ell')=(k, t,\ell),\, i\neq i';\\
\mathbb{E}[W^2_{kt}W^2_{k\ell}W^2_{\ell i}W^2_{t i'}]\leq C\theta_i\theta_k^2\theta_\ell^2\theta^2_t\theta_{i'}, & \mbox{if } (k', t',\ell')=(k, \ell,t);\\
\mathbb{E}[W^3_{kt}W^2_{k\ell}W^3_{\ell i}]\leq C\theta_i\theta^2_k\theta^2_\ell\theta_t, & \mbox{if }(k', t', \ell', i')=(\ell, i, k, t);\\
0, & \mbox{otherwise}. 
\end{cases}
\end{align*}
We combine it with \eqref{proof-Z1-beta} and find that
\begin{align*}
\mathrm{Var}(Z^*_{1c}) &= \frac{4}{v^2}\sum_{ \substack{i,k,\ell (dist), t\neq \{k,\ell\}\\i',k',\ell' (dist), t'\neq \{k',\ell'\}}} \beta_{ik\ell}\beta_{i'k'\ell'} \cdot\mathbb{E}[ W_{kt}W_{k\ell}W^2_{\ell i}W_{k't'}W_{k'\ell'}W^2_{\ell' i'}]\cr
&\leq \frac{C\|\theta\|^4}{\|\theta\|_1^4}\Bigl( \sum_{i,k,\ell,t}\theta_i\theta_k^2\theta_\ell^2\theta_t + \sum_{i,k,\ell,t,i'}\theta_i\theta_k^2\theta_\ell^3\theta_t\theta_{i'} + \sum_{i,k,\ell,t,i'}\theta_i\theta_k^2\theta_\ell^2\theta^2_t\theta_{i'}  \Bigr)\cr
&\leq \frac{C\|\theta\|^4}{\|\theta\|_1^4}\bigl( \|\theta\|^4\|\theta\|_1^2 + \|\theta\|^2\|\theta\|_3^3\|\theta\|_1^3 + \|\theta\|^6\|\theta\|_1^2   \bigr)\cr
&\leq \frac{C\|\theta\|^6\|\theta\|_3^3}{\|\theta\|_1}. 
\end{align*}
Consider $Z_{1c}^{\dag}$. Re-write
\[
Z_{1c}^{\dag} = \frac{1}{v}\sum_{i,k,\ell (dist)} \beta_{ik\ell}    W_{ik}^2 W_{k\ell}W_{\ell i} + \frac{1}{v}\sum_{\substack{i,k,\ell (dist)\\s\notin\{ i,\ell\}, t\notin\{ k,\ell\}\\(s,t)\neq (k,i)}} \beta_{ik\ell}    W_{is}W_{kt}W_{k\ell}W_{\ell i}. 
\]
Regarding the first term, by direct calculations,
\begin{align*}
& \mathbb{E}[W_{ik}^2 W_{k\ell}W_{\ell i}\cdot W_{i'k'}^2 W_{k'\ell'}W_{\ell' i'}]\cr
= & \begin{cases}
\mathbb{E}[W_{ik}^4 W^2_{k\ell}W^2_{\ell i}]\leq C\theta_i^2\theta_k^2\theta^2_{\ell}, &\mbox{if }\ell'=\ell, \, \{i',k'\}=\{i,k\};\\
\mathbb{E}[W_{ik}^3 W^2_{k\ell}W^3_{\ell i}] \leq C\theta_i^2\theta_k^2\theta^2_{\ell}, &\mbox{if }(\ell', k')=(k, \ell), \, i'=i;\\
0, &\mbox{otherwise}. 
\end{cases}
\end{align*}
Combining it with \eqref{proof-Z1-beta} gives
\begin{align*}
\mathrm{Var}\Bigl( \frac{1}{v}\sum_{i,k,\ell (dist)} \beta_{ik\ell}    W_{ik}^2 W_{k\ell}W_{\ell i} \Bigr)\leq \frac{C\|\theta\|^4}{\|\theta\|_1^4}\sum_{i,j,k,\ell}\theta_i^2\theta_k^2\theta_\ell^2\leq \frac{C\|\theta\|^{10}}{\|\theta\|_1^4}. 
\end{align*}
Regarding the second term, for $W_{is}W_{kt}W_{k\ell}W_{\ell i}$ and $W_{i's'}W_{k't'}W_{k'\ell'}W_{\ell' i'}$ to be correlated, all $W$ terms have to be perfectly matched. For each fixed $(i,k,\ell,s,t)$, there are only a constant number of $(i',k',\ell',s',t')$ so that the above is satisfied. Mimicking the argument in \eqref{LemmaC1-eq3}, we have
\begin{align*}
\mathrm{Var}\Bigl(\frac{1}{v}\sum_{\substack{i,k,\ell (dist)\\s\notin\{ i,\ell\}, t\notin\{ k,\ell\}\\(s,t)\neq (k,i)}} &\beta_{ik\ell}    W_{is}W_{kt}W_{k\ell}W_{\ell i}\Bigr) \leq \frac{C}{v^2}\sum_{\substack{i,k,\ell (dist)\\s\notin\{ i,\ell\}, t\notin\{ k,\ell\}\\(s,t)\neq (k,i)}}\beta^2_{ik\ell}  \cdot \mathrm{Var}( W_{is}W_{kt}W_{k\ell}W_{\ell i})\cr
&\leq \frac{C}{\|\theta\|_1^4} \sum_{i,k,\ell,s,t} \|\theta\|^4\cdot \theta_i^2\theta_k^2\theta_\ell^2\theta_s\theta_t \leq \frac{C\|\theta\|^{10}}{\|\theta\|_1^2}.  
\end{align*}
It follows that
\[
\mathrm{Var}(Z^{\dag}_{1c})\leq \frac{C\|\theta\|^{10}}{\|\theta\|_1^2}. 
\]
Combining the above results and noticing that $\|\theta\|^4\leq \|\theta\|_1\|\theta\|_3^3$, we immediately have
\beq   \label{proof-Z1-result4}
\mathrm{Var}(Z_{1c})\leq 3\mathrm{Var}(\widetilde{Z}_{1c}) + 3\mathrm{Var}(Z^*_{1c})  + 3\mathrm{Var}(Z^{\dag}_{1c})\leq \frac{C\|\theta\|^6\|\theta\|_3^3}{\|\theta\|_1}. 
\eeq

We now combine \eqref{proof-Z1-result1}, \eqref{proof-Z1-result2}, \eqref{proof-Z1-result3}, and \eqref{proof-Z1-result4}. Since $Z_1=Z_{1a}+Z_{1b}+Z_{1c}$, it follows that
\[
\mathbb{E}[Z_1] = \|\theta\|^4\cdot [1+o(1)], \qquad \mathrm{Var}(Z_1)\leq C\|\theta\|^2\|\theta\|_3^6=o(\|\theta\|^8). 
\]
This proves the claims of $Z_1$.

Next, we analyze $Z_2$. Since $\delta_{ij}=\eta_i(\eta_j - \teta_j)+\eta_j(\eta_i - \teta_i)$, by direct calculations,  
\begin{align*}
Z_2 &= \sum_{i,j,k,\ell (dist)}\eta_i (\eta_j - \tilde{\eta}_j)W_{jk}\eta_k (\eta_\ell - \tilde{\eta}_\ell)W_{\ell i} + \sum_{i,j,k,\ell (dist)}\eta_i (\eta_j - \tilde{\eta}_j)W_{jk}(\eta_k - \tilde{\eta}_k)\eta_\ell W_{\ell i}  \cr
&\qquad  +\sum_{i,j,k,\ell (dist)}(\eta_i - \tilde{\eta}_i)\eta_j W_{jk}\eta_k (\eta_\ell - \tilde{\eta}_\ell) W_{\ell i}  +  \sum_{i,j,k,\ell (dist)}(\eta_i - \tilde{\eta}_i)\eta_j W_{jk}(\eta_k - \tilde{\eta}_k)\eta_\ell W_{\ell i}. 
\end{align*}
By relabeling the indices, we find out that the first and last sums are equal and that the second and third sums are equal. It follows that
\begin{eqnarray} \label{proof-Z2-decompose}
Z_2 &=& 2\sum_{i,j,k,\ell (dist)}\eta_i (\eta_j - \tilde{\eta}_j)W_{jk}\eta_k (\eta_\ell - \tilde{\eta}_\ell)W_{\ell i}\cr
&& + 2\sum_{i,j,k,\ell (dist)}\eta_i (\eta_j - \tilde{\eta}_j)W_{jk}(\eta_k - \tilde{\eta}_k)\eta_\ell W_{\ell i}\cr
&\equiv& Z_{2a} + Z_{2b}. 
\end{eqnarray}

First, we study $Z_{2a}$. It is seen that
\begin{align*}
Z_{2a} & = 2\sum_{i,j,k,\ell (dist)}\eta_i \Bigl(-\frac{1}{\sqrt{v}}\sum_{s\neq j}W_{js}\Bigr)W_{jk}\eta_k \Bigl(-\frac{1}{\sqrt{v}}\sum_{t\neq \ell }W_{\ell t}\Bigr)W_{\ell i}\cr
&= \frac{2}{v}\sum_{\substack{i,j,k,\ell (dist)\\s\neq j,t\neq \ell}} \eta_i\eta_k W_{js}W_{jk}W_{\ell t}W_{\ell i}.
\end{align*}
We divide summands into four groups: (i) $s=k$ and $t=i$, (ii) $s=k$ and $t\neq i$, (iii) $s\neq k$ and $t=i$, (iv) $s\neq k$ and $t\neq i$. By symmetry between $(j, k, s)$ and $(\ell, i, t)$, the sum of group (ii) and group (iii) are equal. We end up with 
\begin{align*}
Z_{2a}&= \frac{2}{v}\sum_{i,j,k,\ell (dist)} \eta_i\eta_k W^2_{jk}W^2_{\ell i} + \frac{4}{v}\sum_{\substack{i,j,k,\ell (dist)\\s\notin\{j,k\}}} \eta_i\eta_k W_{js}W_{jk}W^2_{\ell i}\cr
&\qquad + \frac{2}{v}\sum_{\substack{i,j,k,\ell (dist)\\s\notin\{j,k\}, t\notin\{ \ell,i\}}} \eta_i\eta_k W_{js}W_{jk}W_{\ell t}W_{\ell i}\cr
&\equiv \widetilde{Z}_{2a}+Z^*_{2a}+Z^{\dag}_{2a},
\end{align*}
Only $\widetilde{Z}_{2a}$ has a nonzero mean. It follows that
\[
\mathbb{E}[Z_{2a}] =\mathbb{E}[\widetilde{Z}_{2a}] =  \frac{2}{v}\sum_{i,j,k,\ell (dist)} \eta_i\eta_k \Omega_{jk}(1-\Omega_{jk})\Omega_{\ell i}(1-\Omega_{\ell i}). 
\] 
Under the null hypothesis, $\Omega_{ij}=\theta_i\theta_j$. Hence, $\Omega_{jk}(1-\Omega_{jk})\Omega_{\ell i}(1-\Omega_{\ell i})=\theta_j\theta_k\theta_\ell\theta_i\cdot [1+O(\theta^2_{\max})]$. Additionally, in the proof of \eqref{proof-Z1-result3}, we have seen that $v=[1+o(1)]\cdot \|\theta\|_1^2$ and $\eta_j=[1+o(1)]\cdot \theta_j$. Combining these results gives
\begin{eqnarray} \label{proof-Z2-result1}
\mathbb{E}[Z_{2a}] &=&\frac{2[1+o(1)]}{\|\theta\|_1^2}\sum_{i,j,k,\ell (dist)}(\theta_i\theta_k)(\theta_j\theta_k\theta_\ell\theta_i)\cr
&=&\frac{2[1+o(1)]}{\|\theta\|_1^2}\Bigl[ \sum_{i,j,k,\ell}\theta^2_i\theta_j\theta^2_k\theta_\ell - \sum_{\substack{i,j,k,\ell\\\ \text{(not dist)}}}\theta^2_i\theta_j\theta^2_k\theta_\ell\Bigr]\cr
&=&\frac{2[1+o(1)]}{\|\theta\|_1^2}\Big[\|\theta\|^4\|\theta\|_1^2  - O\bigl(\|\theta\|_4^4\|\theta\|_1^2+\|\theta\|_3^3\|\theta\|^2\|\theta\|_1+\|\theta\|^6\bigr)\Bigr]\cr
&=& \frac{2[1+o(1)]}{\|\theta\|_1^2}\cdot \|\theta\|^4\|\theta\|_1^2 [1+o(1)] \cr
&=&[1+o(1)]\cdot 2\|\theta\|^4. 
\end{eqnarray}
We then bound the variance of $Z_a$. Consider $\widetilde{Z}_{2a}$ first. Note that $W_{jk}^2W_{\ell i}^2$ and $W_{j'k'}^2W_{\ell' i'}^2$ are correlated only if either $\{j',k'\}=\{j,k\}$ or $\{j',k'\}=\{\ell,i\}$. By symmetry, it suffices to consider $\{j', k'\}=\{j,k\}$. Direct calculations show that
\begin{align*}
& \mathrm{Cov}(\eta_i\eta_k W^2_{jk}W^2_{\ell i},\, \eta_{i'}\eta_{k'} W^2_{j'k'}W^2_{\ell' i'})\cr
\leq\; &
\begin{cases}
\eta_k^2\eta^2_i \mathbb{E}[W^4_{jk}W^4_{\ell i}]\leq C\theta^3_i\theta_j\theta^3_k\theta_\ell, &\mbox{if }(j', k')=(j,k),\, i=i', \, \ell=\ell';\\
\eta_k^2\eta^2_i \mathbb{E}[W^4_{jk}W^2_{\ell i}W^2_{\ell' i}]\leq C\theta^4_i\theta_j\theta^3_k\theta_\ell\theta_{\ell'}, &\mbox{if }(j', k')=(j,k),\, i=i',\, \ell \neq \ell';\\
\eta_k^2\eta_i\eta_{i'} \mathbb{E}[W^4_{jk}W^2_{\ell i}W^2_{\ell'i'}]\leq C\theta^2_i\theta_j\theta^3_k\theta_\ell\theta^2_{i'}\theta_{\ell'}, &\mbox{if }(j', k')=(j,k),\, i\neq i';\\
\eta_j\eta_k\eta^2_i \mathbb{E}[W^4_{jk}W^4_{\ell i}]\leq C\theta^3_i\theta^2_j\theta^2_k\theta_\ell, &\mbox{if }(j', k')=(k, j),\, i=i', \, \ell=\ell';\\
\eta_j\eta_k\eta^2_i \mathbb{E}[W^4_{jk}W^2_{\ell i}W^2_{\ell' i}]\leq C\theta^4_i\theta^2_j\theta^2_k\theta_\ell\theta_{\ell'}, &\mbox{if }(j', k')=(k, j),\, i=i',\, \ell \neq \ell';\\
\eta_j\eta_k \eta_i\eta_{i'} \mathbb{E}[W^4_{jk}W^2_{\ell i}W^2_{\ell'i'}]\leq C\theta^2_i\theta^2_j\theta^2_k\theta_\ell\theta^2_{i'}\theta_{\ell'}, &\mbox{if }(j', k')=(k, j),\, i\neq i';\\
0, &\mbox{otherwise}. 
\end{cases}
\end{align*}
As a result,
\begin{align*}
\mathrm{Var}(\widetilde{Z}_{2a}) &= \frac{4}{v^2}\sum_{\substack{i,j,k,\ell (dist)\\i',j',k',\ell' (dist)}} \mathrm{Cov}(\eta_i\eta_k W^2_{jk}W^2_{\ell i},\, \eta_{i'}\eta_{k'} W^2_{j'k'}W^2_{\ell' i'})\cr
&\leq \frac{C}{\|\theta\|_1^4} \bigl( \|\theta\|_3^6 \|\theta\|_1^2 + \|\theta\|_4^4\|\theta\|_3^3 \|\theta\|_1^3 +  \|\theta\|_3^3\|\theta\|^4 \|\theta\|_1^3 \cr
&\qquad + \|\theta\|_3^3\|\theta\|^4 \|\theta\|_1 +  \|\theta\|_4^4\|\theta\|^4 \|\theta\|_1^2 + \|\theta\|^8\|\theta\|_1^2  \bigr)\cr
&\leq \frac{C\|\theta\|^4\|\theta\|_3^3}{\|\theta\|_1},
\end{align*}
where the last line is obtained as follows: There are six terms in the brackets; since $\|\theta\|^4\leq \|\theta\|_1\|\theta\|_3^3$, the last three terms are dominated by the first three terms; for the first three terms, since $\|\theta\|_3^3\leq \theta^2_{\max}\|\theta\|_1=o(\|\theta\|_1)$ and $\|\theta\|_4^4\leq \theta_{\max}^2\|\theta\|^2=o(\|\theta\|^2)$, the third term dominates. Consider $Z^*_{2a}$ next. We note that for 
\[
\mathbb{E}[W_{js}W_{jk}W^2_{\ell i}\cdot W_{j's'}W_{j'k'}W^2_{\ell' i'}]
\]
to be nonzero, it has to be the case of either $(W_{j's'}, W_{j'k'})=(W_{js}, W_{jk})$ or $(W_{j's'}, W_{j'k'})=(W_{jk}, W_{js})$. This can only happen if $(j', s',k')=(j,s,k)$ or $(j',s',k')=(j,k,s)$. By elementary calculations,
\begin{align*}
& \eta_i\eta_k\eta_{i'}\eta_{k'}\cdot\mathbb{E}[W_{js}W_{jk}W^2_{\ell i}\cdot W_{j's'}W_{j'k'}W^2_{\ell' i'}]\cr
=\; &\begin{cases}
 \eta^2_i\eta^2_k\, \mathbb{E}[W^2_{js}W^2_{jk}W^4_{\ell i}]\leq C\theta^3_i\theta^2_j\theta^3_k\theta_\ell\theta_s, &\mbox{if }(j', s', k')=(j, s, k),\; i'= i,\, \ell'=\ell;\cr
 \eta^2_i\eta^2_k\, \mathbb{E}[W^2_{js}W^2_{jk}W^2_{\ell i}W^2_{\ell' i}]\leq C\theta^4_i\theta^2_j\theta^3_k\theta_\ell\theta_s\theta_{\ell'}, &\mbox{if }(j', s', k')=(j, s, k),\; i'= i,\, \ell'\neq \ell;\cr
 \eta_i\eta_{i'}\eta^2_k\, \mathbb{E}[W^2_{js}W^2_{jk}W^2_{\ell i}W^2_{\ell' i'}]\leq C\theta^2_i\theta^2_j\theta^3_k\theta_\ell\theta_s\theta^2_{i'}\theta_{\ell'}, &\mbox{if }(j', s', k')=(j, s, k),\; i\neq i';\cr
 \eta^2_i\eta_k\eta_s\, \mathbb{E}[W^2_{js}W^2_{jk}W^4_{\ell i}]\leq C\theta^3_i\theta^2_j\theta^2_k\theta_\ell\theta^2_s, &\mbox{if }(j', s', k')=(j, k, s),\; i'= i,\, \ell'=\ell;\cr
 \eta^2_i\eta_k\eta_s\, \mathbb{E}[W^2_{js}W^2_{jk}W^2_{\ell i}W^2_{\ell' i}]\leq C\theta^4_i\theta^2_j\theta^2_k\theta_\ell\theta^2_s\theta_{\ell'}, &\mbox{if }(j', s', k')=(j, k, s),\; i'= i,\, \ell'\neq \ell;\cr
 \eta_i\eta_{i'}\eta_k\eta_s\, \mathbb{E}[W^2_{js}W^2_{jk}W^2_{\ell i}W^2_{\ell' i'}]\leq C\theta^2_i\theta^2_j\theta^2_k\theta_\ell\theta^2_s\theta^2_{i'}\theta_{\ell'}, &\mbox{if }(j', s', k')=(j, k, s),\; i\neq i';\cr
 0, &\mbox{otherwise}. 
\end{cases}
\end{align*}
It follows that
\begin{align*}
\mathrm{Var}(Z^*_{2a}) &= \frac{16}{v^2}\sum_{\substack{i,j,k,\ell (dist)\\i',j',k',\ell' (dist)}} \eta_i\eta_k\eta_{i'}\eta_{k'}\cdot\mathbb{E}[W_{js}W_{jk}W^2_{\ell i}\cdot W_{j's'}W_{j'k'}W^2_{\ell' i'}] \cr
&\leq \frac{C}{\|\theta\|_1^4} \bigl( \|\theta\|_3^6 \|\theta\|^2 \|\theta\|_1^2 + \|\theta\|_4^4\|\theta\|_3^3\|\theta\|^2 \|\theta\|_1^3 +  \|\theta\|_3^3\|\theta\|^6 \|\theta\|_1^3 \cr
&\qquad + \|\theta\|_3^3\|\theta\|^6 \|\theta\|_1 +  \|\theta\|_4^4\|\theta\|^6 \|\theta\|_1^2 + \|\theta\|^{10}\|\theta\|_1^2  \bigr)\cr
&\leq \frac{C\|\theta\|^6\|\theta\|_3^3}{\|\theta\|_1},
\end{align*} 
where the last inequality is obtained similarly as in the calculation of $\mathrm{Var}(\widetilde{Z}_{2a})$. Last, consider $Z^{\dag}_{2a}$. Write
\beq \label{Z2a-dag}
Z_{2a}^{\dag} = \frac{2}{v}\sum_{i,j,k,\ell (dist)} \eta_i\eta_k W_{j\ell}^2 W_{jk}W_{\ell i} + \frac{2}{v}\sum_{\substack{i,j,k,\ell (dist)\\s\notin\{j,k\}, t\notin\{ \ell,i\}\\(s,t)\neq (\ell,j)}} \eta_i\eta_k W_{js}W_{jk}W_{\ell t}W_{\ell i}
\eeq
Regarding the first term, we note that
\begin{align*}
&\eta_i\eta_k\eta_{i'}\eta_{k'}\cdot\mathbb{E}[W_{j\ell}^2 W_{jk}W_{\ell i}\cdot W_{j'\ell'}^2 W_{j'k'}W_{\ell' i'}]\cr
= &\begin{cases}
\eta_i^2\eta_k^2\, \mathbb{E}[W^2_{jk}W^2_{\ell i} W^4_{j\ell}]\leq C\theta_i^3\theta^2_j\theta_k^3\theta_{\ell}^2, &\mbox{if }(j',k')=(j,k),\, (i',\ell')=(i,\ell);\\
\eta_i\eta_k^2\eta_\ell\, \mathbb{E}[W^2_{jk}W^2_{\ell i} W^2_{j\ell}W^2_{ji}]\leq C\theta_i^3\theta^3_j\theta_k^3\theta_{\ell}^3,  &\mbox{if }(j',k')=(j,k),\, (i',\ell')=(\ell,i);\\
\eta_i^2\eta_k\eta_\ell\, \mathbb{E}[W^2_{jk}W^2_{\ell i} W^2_{j\ell}W^2_{k\ell}]\leq C\theta_i^3\theta^2_j\theta_k^3\theta_{\ell}^4, &\mbox{if }(j',k')=(k,j),\, (i',\ell')=(i,\ell);\\
\eta_i\eta_k\eta_\ell\eta_j\, \mathbb{E}[W^2_{jk}W^2_{\ell i} W^2_{j\ell}W^2_{k i}]\leq C\theta_i^3\theta^3_j\theta_k^3\theta_{\ell}^3, &\mbox{if }(j',k')=(k,j),\, (i',\ell')=(\ell,i);\\
0, &\mbox{otherwise}. 
\end{cases} 
\end{align*}
It follows that
\begin{align*}
& \mathrm{Var}\Bigl( \frac{2}{v}\sum_{i,j,k,\ell (dist)} \eta_i\eta_k W_{j\ell}^2 W_{jk}W_{\ell i} \Bigr)\cr
\leq &\frac{C}{\|\theta\|_1^4}\sum_{i,j,k,\ell}(\theta_i^3\theta^2_j\theta_k^3\theta_{\ell}^2 + \theta_i^3\theta^3_j\theta_k^3\theta_{\ell}^3 + \theta_i^3\theta^2_j\theta_k^3\theta_{\ell}^4)\cr
\leq &\frac{C}{\|\theta\|_1^4}\bigl( \|\theta\|_3^6\|\theta\|^4 + \|\theta\|_3^{12} + \|\theta\|_4^4\|\theta\|_3^6\|\theta\|^2 \bigr)\cr
\leq &\frac{C\|\theta\|_3^6\|\theta\|^4}{\|\theta\|_1^4}.
\end{align*}
Regarding the second term in \eqref{Z2a-dag}. We note that, for $\eta_i\eta_kW_{js}W_{jk}W_{\ell t}W_{\ell i}$ and $\eta_{i'}\eta_{k'}W_{j's'}W_{j'k'}W_{\ell' t'}W_{\ell' i'}$ to be correlated, all the $W$ terms have to be perfectly paired. It turns out that
\[
\mathbb{E}[W_{js}W_{jk}W_{\ell t}W_{\ell i}\cdot W_{j's'}W_{j'k'}W_{\ell' t'}W_{\ell' i'}] = \mathbb{E}[W^2_{js}W^2_{jk}W^2_{\ell t}W^2_{\ell i}]. 
\]
To perfectly pair the $W$ terms, there are two possible cases: (i) $(j', \ell')=(j,\ell)$, $\{s',k'\}=\{s,k\}$, $\{\ell', i'\}=\{\ell,i\}$. (ii) $(j',\ell')=(\ell,j)$, $\{s',k'\}=\{\ell,i\}$, $\{\ell',i'\}=\{s,k\}$. As a result, $\eta_i\eta_k\eta_{i'}\eta_{k'}$ only has the following possibilities:
\[
\begin{array}{llll}
\eta_i\eta_k(\eta_i\eta_k)= \eta^2_i\eta^2_k,& \eta_i\eta_k(\eta_i\eta_s)=\eta^2_i\eta_k\eta_s,&  \eta_i\eta_k(\eta_\ell \eta_k)=\eta_i\eta^2_k\eta_\ell,&\eta_i\eta_k(\eta_\ell \eta_s)=\eta_i\eta_k\eta_\ell\eta_s,\\
\eta_i\eta_k(\eta_k\eta_i)=\eta_i^2\eta_k^2, & \eta_i\eta_k(\eta_k\eta_\ell)=\eta_i\eta_k^2\eta_\ell, & \eta_i\eta_k(\eta_s \eta_i)=\eta_i^2\eta_k\eta_s, & \eta_i\eta_k(\eta_s \eta_\ell)=\eta_i\eta_k\eta_\ell\eta_s. 
\end{array}
\]
By symmetry, there are only three different types: $\eta^2_i\eta^2_k$, $\eta^2_i\eta_k\eta_s$, and $\eta_i\eta_k\eta_\ell\eta_s$. It follows that
\begin{align*}
& \mathrm{Var}\Bigl( \frac{2}{v}\sum_{\substack{i,j,k,\ell (dist)\\s\notin\{j,k\}, t\notin\{ \ell,i\}, (s,t)\neq (\ell,j)}} \eta_i\eta_k W_{js}W_{jk}W_{\ell t}W_{\ell i} \Bigr)\cr
\leq & \frac{C}{\|\theta\|_1^4}\sum_{i,j,k,\ell,s,t}(\theta^2_i\theta^2_k+\theta^2_i\theta_k\theta_s+\theta_i\theta_k\theta_\ell\theta_s)\cdot \theta^2_j\theta_s\theta_k\theta^2_\ell\theta_t\theta_i \cr
\leq & \frac{C}{\|\theta\|_1^4}\sum_{i,j,k,\ell,s,t}(\theta^3_i\theta^2_j\theta^3_k\theta^2_\ell \theta_s\theta_t +\theta^3_i\theta^2_j\theta^2_k\theta^2_\ell \theta^2_s\theta_t  +\theta^2_i\theta^2_j\theta^2_k\theta^3_\ell \theta^2_s\theta_t)\cr
\leq & \frac{C}{\|\theta\|_1^4}\bigl( \|\theta\|_3^6\|\theta\|^4\|\theta\|_1^2+\|\theta\|_3^3\|\theta\|^8\|\theta\|_1\bigr)
\; \leq  \frac{C\|\theta\|^4\|\theta\|_3^6}{\|\theta\|_1^2}. 
\end{align*}
It follows that
\[
\mathrm{Var}(Z^{\dag}_{2a})\leq \frac{C\|\theta\|^4\|\theta\|_3^6}{\|\theta\|_1^2}.
\]
Comparing the variances of $\widetilde{Z}_{2a}$, $Z^*_{2a}$ and $Z^{\dag}_{2a}$, we find out that the variance of $Z^*_{2a}$ dominates. As a result,
\beq \label{proof-Z2-result2}
\mathrm{Var}(Z_{2a})\leq 3\mathrm{Var}(\widetilde{Z}_{2a}) +3\mathrm{Var}(Z^*_{2a})+3\mathrm{Var}(Z^{\dag}_{2a})\leq \frac{C\|\theta\|^6\|\theta\|_3^3}{\|\theta\|_1}. 
\eeq

Second, we study $Z_{2b}$. It is seen that
\begin{align*}
Z_{2b} &= 2\sum_{i,j,k,\ell (dist)}\eta_i \Bigl(-\frac{1}{\sqrt{v}}\sum_{s\neq j}W_{js}\Bigr) W_{jk}\Bigl(-\frac{1}{\sqrt{v}}\sum_{t\neq k}W_{kt}\Bigr)\eta_\ell W_{\ell i}\cr
&=\frac{2}{v}\sum_{\substack{i,j,k,\ell (dist)\\s\neq j, t\neq k}}\eta_i\eta_\ell W_{js}W_{jk}W_{kt}W_{\ell i}. 
\end{align*}
We divide summands into four groups: (i) $s=k$ and $t=j$, (ii) $s=k$ and $t\neq j$, (iii) $s\neq k$ and $t=j$, (iv) $s\neq k$ and $t\neq j$. By index symmetry, the sums of group (ii) and group (iii) are equal. We end up with 
\begin{align*}
Z_{2b} &= \frac{2}{v}\sum_{i,j,k,\ell (dist)}\eta_i\eta_\ell W^3_{jk}W_{\ell i} + \frac{4}{v}\sum_{i,j,k,\ell (dist), t\notin\{ k,j\}}\eta_i\eta_\ell W^2_{jk}W_{kt}W_{\ell i}\cr
&\qquad + \frac{2}{v}\sum_{i,j,k,\ell (dist), s\neq \{j,k\}, t\neq \{j,k\}}\eta_i\eta_\ell W_{js}W_{jk}W_{kt}W_{\ell i}\cr
&\equiv \widetilde{Z}_{2b} + Z^*_{2b} + Z^{\dag}_{2b}. 
\end{align*}
It is easy to see that all three terms have mean zero. Therefore,
\beq \label{proof-Z2-result3}
\mathbb{E}[Z_{2b}]= 0.
\eeq
We then bound the variances. Consider $\widetilde{Z}_{2b}$ first. By direct calculations,
\begin{align*}
&\eta_i\eta_\ell\eta_{i'}\eta_{\ell'}\cdot \mathbb{E}[ W^3_{jk}W_{\ell i}\cdot W^3_{j'k'}W_{\ell' i'}]\cr
=\; &\begin{cases}
\eta^2_i\eta^2_\ell\cdot \mathbb{E}[ W^6_{jk}W^2_{\ell i}]\leq C\theta^3_i\theta_j\theta_k\theta^3_\ell, &\mbox{if }\{j', k'\}=\{j,k\},\, \{\ell',i'\}=\{\ell,i\};\\
\eta_i\eta_\ell\eta_j\eta_k \cdot \mathbb{E}[ W^4_{jk}W^4_{\ell i}]\leq C\theta^2_i\theta^2_j\theta^2_k\theta^2_\ell, &\mbox{if }\{j', k'\}=\{\ell,i\},\, \{\ell',i'\}=\{j,k\};\\
0, &\mbox{otherwise}. 
\end{cases}
\end{align*}
It follows that
\begin{align*}
\mathrm{Var}(\widetilde{Z}_{2b})&\leq \frac{C}{\|\theta\|_1^4}\Bigl(\sum_{i,j,k,\ell}\theta^3_i\theta_j\theta_k\theta^3_\ell + \sum_{i,j,k,\ell}\theta^2_i\theta^2_j\theta^2_k\theta^2_\ell  \Bigr)\cr
&\leq \frac{C}{\|\theta\|_1^4}\bigl( \|\theta\|_3^6\|\theta\|_1^2 + \|\theta\|^8 \bigr)\cr
&\leq  \frac{C\|\theta\|_3^6}{\|\theta\|_1^2}. 
\end{align*}
Consider $Z^*_{2b}$ next. By direct calculations,
\begin{align*}
&\eta_i\eta_\ell\eta_{i'}\eta_{\ell'}\cdot \mathbb{E}[W^2_{jk}W_{kt}W_{\ell i}\cdot W^2_{j'k'}W_{k't'}W_{\ell' i'}]\cr
=\; &\begin{cases}
\eta^2_i\eta^2_\ell\, \mathbb{E}[ W^4_{jk}W^2_{kt}W^2_{\ell i}]\leq C\theta^3_i\theta_j\theta^2_k\theta^3_\ell\theta_t, &\mbox{if }(k', t')=(k,t),\, \{\ell',i'\}=\{\ell,i\},\, j'= j;\\
\eta^2_i\eta^2_\ell\, \mathbb{E}[ W^2_{jk}W^2_{kt}W^2_{\ell i}W^2_{j'k}]\leq C\theta^3_i\theta_j\theta^3_k\theta^3_\ell\theta_t\theta_{j'}, &\mbox{if }(k', t')=(k,t),\, \{\ell',i'\}=\{\ell,i\},\, j'\neq j;\\
\eta^2_i\eta^2_\ell\, \mathbb{E}[ W^2_{jk}W^2_{kt}W^2_{\ell i}W^2_{j't}]\leq C\theta^3_i\theta_j\theta^2_k\theta^3_\ell\theta^2_t\theta_{j'}, &\mbox{if }(k', t')=(t,k),\, \{\ell',i'\}=\{\ell,i\};\\
\eta_i\eta_\ell\eta_k\eta_t \, \mathbb{E}[ W^2_{jk}W^2_{kt}W^4_{\ell i} ]\leq C\theta^2_i\theta_j\theta^3_k\theta^2_\ell\theta^2_t, &\mbox{if }(k', t')=(\ell,i),\, \{\ell',i'\}=\{k, t\},\, j'= i;\\
\eta_i\eta_\ell\eta_k\eta_t \, \mathbb{E}[ W^2_{jk}W^2_{kt}W^2_{\ell i}W^2_{j'\ell} ]\leq C\theta^2_i\theta_j\theta^3_k\theta^3_\ell\theta^2_t\theta_{j'}, &\mbox{if }(k', t')=(\ell,i),\, \{\ell',i'\}=\{k, t\},\, j'\neq i;\\
\eta_i\eta_\ell\eta_k\eta_t \, \mathbb{E}[ W^2_{jk}W^2_{kt}W^4_{\ell i}]\leq C\theta^2_i\theta_j\theta^3_k\theta^2_\ell\theta^2_t, &\mbox{if }(k', t')=(i, \ell),\, \{\ell',i'\}=\{k, t\},\, j'= \ell;\\
\eta_i\eta_\ell\eta_k\eta_t \, \mathbb{E}[ W^2_{jk}W^2_{kt}W^2_{\ell i}W^2_{j'i} ]\leq C\theta^3_i\theta_j\theta^3_k\theta^2_\ell\theta^2_t\theta_{j'}, &\mbox{if }(k', t')=(i, \ell),\, \{\ell',i'\}=\{k, t\},\, j'\neq \ell;\\
\eta^2_i\eta^2_\ell \, \mathbb{E}[ W^3_{jk}W^3_{kt}W^2_{\ell i}]\leq C\theta^3_i\theta_j\theta^2_k\theta^3_\ell\theta_t, & \mbox{if }(k', t', j') = (k, j, t), \, \{i', \ell'\}=\{i,\ell\};\\
0, &\mbox{otherwise}. 
\end{cases}
\end{align*}
There are only two four types on the right hand side. It follows that
\begin{align*}
\mathrm{Var}(Z^*_{2b})&\leq \frac{C}{\|\theta\|_1^4}\Bigl(\sum_{i,j,k,\ell,t,j'}\theta^3_i\theta_j\theta^3_k\theta^3_\ell\theta_t\theta_{j'} + \sum_{i,j,k,\ell,t,j'}\theta^3_i\theta_j\theta^2_k\theta^3_\ell\theta^2_t\theta_{j'} \cr
&\qquad + \sum_{i,j,k,\ell,t}\theta^3_i\theta_j\theta^2_k\theta^3_\ell\theta_t + \sum_{i,j,k,\ell,t}\theta^2_i\theta_j\theta^3_k\theta^2_\ell\theta^2_t  \Bigr)\cr
&\leq \frac{C}{\|\theta\|_1^4}\bigl( \|\theta\|_3^9\|\theta\|_1^3 +  \|\theta\|_3^6 \|\theta\|^4 \|\theta\|_1^2+ \|\theta\|_3^6 \|\theta\|^2 \|\theta\|_1^2 + \|\theta\|_3^3 \|\theta\|^6\|\theta\|_1\bigr)\cr
&\leq  \frac{C\|\theta\|_3^9}{\|\theta\|_1}. 
\end{align*}
Last, consider $Z^{\dag}_{2b}$. By direct calculations,
\begin{align*}
& \eta_i\eta_\ell \eta_{i'}\eta_{\ell'}\cdot \mathbb{E}[W_{js}W_{jk}W_{kt}W_{\ell i}\cdot W_{j's'}W_{j'k'}W_{k't'}W_{\ell' i'}]\cr
=\; & \begin{cases}
\eta_i^2\eta_\ell^2\, \mathbb{E}[W^2_{js}W^2_{jk}W^2_{kt}W^2_{\ell i}]\leq C\theta^3_i\theta^2_j\theta^2_k\theta^3_\ell\theta_s\theta_t, &\mbox{if }(j',s')=(j,s),\, (k',t')=(k,t),\, \{\ell', i'\}=\{\ell, i\};\\
\eta_i^2\eta_\ell^2\, \mathbb{E}[W^2_{js}W^2_{jk}W^2_{kt}W^2_{\ell i}]\leq C\theta^3_i\theta^2_j\theta^2_k\theta^3_\ell\theta_s\theta_t, &\mbox{if }(j',s')=(k, t),\, (k',t')=(j,s),\, \{\ell', i'\}=\{\ell, i\};\\
0, &\mbox{otherwise}. 
\end{cases} 
\end{align*}
It follows that
\[
\mathrm{Var}(Z^{\dag}_{2b})\leq \frac{C}{\|\theta\|_1^4}\sum_{i,j,k,\ell,s,t} \theta^3_i\theta^2_j\theta^2_k\theta^3_\ell\theta_s\theta_t\leq \frac{C\|\theta\|^4\|\theta\|_3^6}{\|\theta\|_1^2}. 
\]
Since $\|\theta\|_1\|\theta\|_3^3\geq\|\theta\|^4\to\infty$, the variance of $Z^*_{2b}$ dominates the variances of $\widetilde{Z}_{2b}$ and $Z^{\dag}_{2b}$. We thus have
\beq \label{proof-Z2-result4}
\mathrm{Var}(Z_{2b})\leq 3\mathrm{Var}(\widetilde{Z}_{2b}) + 3\mathrm{Var}(Z^*_{2b}) + 3\mathrm{Var}(Z^{\dag}_{2b})\leq \frac{C\|\theta\|_3^9}{\|\theta\|_1}. 
\eeq

We now combine \eqref{proof-Z2-result1}, \eqref{proof-Z2-result2}, \eqref{proof-Z2-result3}, and \eqref{proof-Z2-result4}. Since $\|\theta\|_3^6\leq \theta_{\max}^2\|\theta\|^4\ll \|\theta\|^6$, the right hand side of \eqref{proof-Z2-result4} is much smaller than the right hand side of \eqref{proof-Z2-result2}. It yields that
\[
\mathbb{E}[Z_2] = 2\|\theta\|^4\cdot [1+o(1)], \qquad \mathrm{Var}(Z_2)\leq \frac{C\|\theta\|^6\|\theta\|_3^3}{\|\theta\|_1}=o(\|\theta\|^8). 
\]
This proves the claims of $Z_2$.

\subsubsection{Proof of Lemma~\ref{lem:ProxySgnQ(b)-alt}}
It suffices to prove the claims for each of $Z_1$-$Z_6$. We have analyzed $Z_1$-$Z_2$ under the null hypothesis. The proof for the alternative hypothesis is similar and omitted. We obtain that 
\begin{align*}
& \bigl|\mathbb{E}[Z_1]\bigr|\leq C\|\theta\|^4, \qquad \mathrm{Var}(Z_1)\leq C\|\theta\|^2\|\theta\|_3^6=o(\|\theta\|^8),\cr
&\bigl|\mathbb{E}[Z_2]\bigr| \leq C\|\theta\|^4, \qquad \mathrm{Var}(Z_2) \leq \frac{C\|\theta\|^6\|\theta\|_3^3}{\|\theta\|_1}=o(\|\theta\|^8).  
\end{align*}

First, we analyze $Z_3$. Since $\delta_{ij}=\eta_i(\eta_j-\teta_j)+\eta_j(\eta_i-\teta_i)$, we have
\begin{align} \label{proof-Z3-decompose}
Z_3 =& \sum_{\substack{i, j, k, \ell\\ (dist)}}\eta_i(\eta_j-\teta_j)\eta_j(\eta_k-\teta_k)   \widetilde{\Omega}_{k \ell}  W_{\ell i} +
\sum_{\substack{i, j, k, \ell\\ (dist)}}\eta_i(\eta_j-\teta_j)^2\eta_k   \widetilde{\Omega}_{k \ell}  W_{\ell i} \cr
&\hspace{-30pt}+\sum_{\substack{i, j, k, \ell\\ (dist)}}(\eta_i-\teta_i)\eta^2_j(\eta_k-\teta_k)   \widetilde{\Omega}_{k \ell}  W_{\ell i} +
\sum_{\substack{i, j, k, \ell\\ (dist)}}(\eta_i-\teta_i)\eta_j (\eta_j-\teta_j)\eta_k   \widetilde{\Omega}_{k \ell}  W_{\ell i}\cr
&\equiv Z_{3a}+Z_{3b} + Z_{3c}+Z_{3d}. 
\end{align}

First, we study $Z_{3a}$. By direct calculations,
\begin{align*}
Z_{3a} &= \sum_{i, j, k, \ell (dist)}\eta_i\Bigl(-\frac{1}{\sqrt{v}}\sum_{s\neq j}W_{js}\Bigr)\eta_j\Bigl(-\frac{1}{\sqrt{v}}\sum_{t\neq k}W_{kt}\Bigr)  \widetilde{\Omega}_{k \ell}  W_{\ell i}\cr
&=\frac{1}{v}\sum_{\substack{ i, j, k, \ell (dist)\\s\neq j,t\neq k}} \beta_{ijk\ell}  W_{js}W_{kt}W_{\ell i}, \qquad \mbox{where}\;\; \beta_{ijk\ell}=\eta_i\eta_j\widetilde{\Omega}_{k\ell}. 
\end{align*}
Since $(i,j,k,\ell)$ are distinct, all summands have mean zero. Hence,
\beq \label{proof-Z3-result1}
\mathbb{E}[Z_{3a}]=0.
\eeq
To bound its variance, re-write
\begin{align*}
Z_{3a} &=\frac{1}{v}\sum_{i, j, k, \ell (dist)} \beta_{ijk\ell}  W^2_{jk}W_{\ell i} +   \frac{1}{v}\sum_{\substack{ i, j, k, \ell (dist)\\s\neq j,t\neq k, (s,t)\neq (k,j)}} \beta_{ijk\ell} W_{js}W_{kt}W_{\ell i}\cr
&\equiv \widetilde{Z}_{3a} + Z^*_{3a}. 
\end{align*}
We note that $|\beta_{ijk\ell}|\leq C\alpha\theta_i\theta_j\theta_k\theta_\ell$ by \eqref{tOmega(ij)} and \eqref{eta-bound}. Consider the variance of $\widetilde{Z}_{3a}$. By direct calculations, 
\begin{align*}
&\beta_{ijk\ell}\beta_{i'j'k'\ell'}\cdot\mathrm{Cov}(W^2_{jk}W_{\ell i},\, W^2_{j'k'}W_{\ell' i'})\cr
=& \begin{cases}
C\alpha^2\theta_i^2\theta_j^2\theta_k^2\theta_\ell^2\,\mathbb{E}[W^4_{jk}W^2_{\ell i}]\leq C\alpha^2\theta^3_i\theta^3_j\theta^3_k\theta^3_\ell, &\mbox{if }\{\ell',i'\}=\{\ell,i\},\,\{j',k'\}=\{j,k\};\\
C\alpha^2\theta_i^2\theta_j\theta_k\theta_\ell^2\theta_{j'}\theta_{k'}\,\mathbb{E}[W^2_{jk}W^2_{j'k'}W^2_{\ell i}]\leq C\alpha^2\theta^3_i\theta^2_j\theta^2_k\theta^3_\ell\theta^2_{j'}\theta^2_{k'}, &\mbox{if }\{\ell',i'\}=\{\ell,i\},\, \{j',k'\}\neq \{j,k\};\\
C\alpha^2\theta_i^2\theta_j^2\theta_k^2\theta_\ell^2\,\mathbb{E}[W^3_{jk}W^3_{\ell i}]\leq C\alpha^2\theta^3_i\theta^3_j\theta^3_k\theta^3_\ell, &\mbox{if }\{j',k'\}=\{\ell,i\}, \,\{\ell',i'\}=\{j,k\};\\
0, &\mbox{otherwise}. 
\end{cases}
\end{align*}
It follows that 
\begin{align*}
\mathrm{Var}(\widetilde{Z}_{3a})& \leq \frac{C\alpha^2}{\|\theta\|_1^4}\Bigl( \sum_{i,j,k,\ell}\theta^3_i\theta^3_j\theta^3_k\theta^3_\ell + \sum_{i,j,k,\ell,j',k'}\theta^3_i\theta^2_j\theta^2_k\theta^3_\ell\theta^2_{j'}\theta^2_{k'} \Bigr)\cr
&\leq \frac{C\alpha^2}{\|\theta\|_1^4}\bigl( \|\theta\|_3^{12} + \|\theta\|^8\|\theta\|_3^6 \bigr)\cr
&\leq \frac{C\alpha^2\|\theta\|_3^{12}}{\|\theta\|_1^2}. 
\end{align*}
Consider the variance of $Z^*_{3a}$. For $W_{js}W_{kt}W_{\ell i}$ and $W_{j's'}W_{k't'}W_{\ell'i'}$ to be correlated, all $W$ terms have to be perfectly paired. 
By symmetry across indices, it reduces to three cases: (i) $(\ell',i')=(\ell,i)$, $(j',s')=(j,s)$, $(k',t')=(k,t)$; (ii) $(\ell', i')=(j,s)$, $(j',s')=(\ell,i)$, $(k',t')=(k,t)$; (iii) $(\ell', i')=(j,s)$, $(j', s')=(k,t)$, $(k',t')=(\ell, i)$. It follows that
\begin{align*}
& \beta_{ijk\ell}\beta_{i'j'k'\ell'}\cdot \mathbb{E}[ W_{js}W_{kt}W_{\ell i}\cdot W_{j's'}W_{k't'}W_{\ell'i'} ]\cr
\leq \; &C\alpha^2(\theta_i\theta_j\theta_k\theta_\ell)(\theta_{i'}\theta_{j'}\theta_{k'}\theta_{\ell'})\cdot \mathbb{E}[ W^2_{js}W^2_{kt}W^2_{\ell i}]\cr
\leq \; &\begin{cases}
C\alpha^2\theta_i^2\theta_j^2\theta_k^2\theta_\ell^2 \mathbb{E}[ W^2_{js}W^2_{kt}W^2_{\ell i}]\leq C\alpha^2\theta^3_i\theta^3_j\theta^3_k\theta^3_\ell\theta_s\theta_t, &\mbox{case (i)}\\
C\alpha^2(\theta_i\theta_j\theta_k\theta_\ell)(\theta_s \theta_\ell \theta_k\theta_j)\mathbb{E}[ W^2_{js}W^2_{kt}W^2_{\ell i}]\leq C\alpha^2\theta^2_i\theta^3_j\theta^3_k\theta^3_\ell \theta^2_s\theta_t, &\mbox{case (ii)}\\
C\alpha^2(\theta_i\theta_j\theta_k\theta_\ell)(\theta_s \theta_k \theta_\ell\theta_j)\mathbb{E}[ W^2_{js}W^2_{kt}W^2_{\ell i}]\leq C\alpha^2\theta^2_i\theta^3_j\theta^3_k\theta^3_\ell \theta^2_s\theta_t, &\mbox{case (iii)}\\
0, &\mbox{otherwise}. 
\end{cases}
\end{align*}
As a result,
\begin{align*} 
\mathrm{Var}(Z^*_{3a}) &\leq \frac{C}{\|\theta\|_1^4}\Bigl(\sum_{i,j,k,\ell,s,t}\alpha^2\theta^3_i\theta^3_j\theta^3_k\theta^3_\ell\theta_s\theta_t + \sum_{i,j,k,\ell,s,t}\alpha^2\theta^2_i\theta^3_j\theta^3_k\theta^3_\ell \theta^2_s\theta_t \Bigr)\cr
&\leq \frac{C\alpha^2}{\|\theta\|_1^4}\bigl(\|\theta\|_3^{12}\|\theta\|_1^2+\|\theta\|^4\|\theta\|_3^9\|\theta\|_1 \bigr)\cr
&\leq \frac{C\alpha^2\|\theta\|_3^{12}}{\|\theta\|_1^2}. 
\end{align*}
Combining the variance of $\widetilde{Z}_{3a}$ and $Z^*_{3a}$ gives
\beq \label{proof-Z3-result2}
\mathrm{Var}(Z_{3a})\leq \frac{C\alpha^2\|\theta\|_3^{12}}{\|\theta\|_1^2}.
\eeq

Second, we study $Z_{3b}$. It is seen that 
\begin{align*}
Z_{3b} &= \sum_{i, j, k, \ell (dist)}\eta_i\Bigl(-\frac{1}{\sqrt{v}}\sum_{s\neq j}W_{js}\Bigr)\Bigl(-\frac{1}{\sqrt{v}}\sum_{t\neq j}W_{jt}\Bigr)\eta_k   \widetilde{\Omega}_{k \ell}  W_{\ell i}\cr
&= \frac{1}{v}\sum_{\substack{i,j,\ell (dist)\\s\neq j,t\neq j}} \Bigl( \sum_{k\notin\{i,j,\ell\}} \eta_i\eta_k\widetilde{\Omega}_{k\ell} \Bigr)    W_{js}W_{jt}W_{\ell i}\cr
&\equiv \frac{1}{v}\sum_{\substack{i,j,\ell (dist)\\s\neq j,t\neq j}} \beta_{ij\ell}   W_{js}W_{jt}W_{\ell i},
\end{align*}
where by \eqref{tOmega(ij)} and \eqref{eta-bound},
\beq\label{proof-Z3-beta}
|\beta_{ij\ell}|\leq \sum_{k\notin\{i,j,\ell\}}|\eta_i\eta_k\widetilde{\Omega}_{k\ell}|\leq \sum_{k}C\alpha\theta_i\theta_k^2\theta_\ell\leq C\alpha\|\theta\|^2\cdot\theta_i\theta_\ell. 
\eeq
We further decompose $Z_{3b}$ into 
\[
Z_{3b} = \frac{1}{v}\sum_{\substack{i,j,\ell (dist)\\s\neq j}} \beta_{ij\ell}   W^2_{js}W_{\ell i} + \frac{1}{v}\sum_{\substack{i,j,\ell (dist)\\s, t (dist) \notin \{j\}}} \beta_{ij\ell}   W_{js}W_{jt}W_{\ell i} \equiv \widetilde{Z}_{3b} + Z^*_{3b}. 
\]
It is easy to see that both terms have mean zero. It follows that
\beq  \label{proof-Z3-result3}
\mathbb{E}[Z_{3b}] = 0. 
\eeq
To calculate the variance of $\widetilde{Z}_{3b}$, we note that
\begin{align*}
& \beta_{ij\ell} \beta_{i'j'\ell'} \cdot\mathbb{E}[ W^2_{js}W_{\ell i}\cdot  W^2_{j's'}W_{\ell' i'}]\cr
\leq\; &C\alpha^2\|\theta\|^4\theta_i\theta_{i'}\theta_\ell\theta_{\ell'}\cdot \mathbb{E}[ W^2_{js}W_{\ell i}\cdot  W^2_{j's'}W_{\ell' i'}]\cr
\leq \; &\begin{cases}
C\alpha^2\|\theta\|^4\theta_i^2\theta_\ell^2\cdot \mathbb{E}[W^4_{js}W^2_{\ell i}]\leq C\alpha^2\|\theta\|^4\theta_i^3\theta_j\theta_\ell^3\theta_s &\mbox{if }\{\ell', i'\}=\{\ell, i\},\; \{j',s'\}= \{j,s\}\\
C\alpha^2\|\theta\|^4\theta_i^2\theta_\ell^2\cdot \mathbb{E}[W^2_{js}W^2_{\ell i}W^2_{j's'}]\leq C\alpha^2\|\theta\|^4\theta_i^3\theta_j\theta_\ell^3\theta_s\theta_{j'}\theta_{s'}, &\mbox{if }\{\ell', i'\}=\{\ell, i\},\; \{j',s'\}\neq \{j,s\};\\
C\alpha^2\|\theta\|^4\theta_i\theta_\ell\theta_j\theta_s \cdot \mathbb{E}[ W^3_{js}W^3_{\ell i}]\leq C\alpha^2\|\theta\|^4\theta_i^2\theta^2_j\theta^2_\ell\theta^2_s, & \mbox{if }\{\ell', i'\}=\{j, s\},\; \{j',s'\}= \{\ell,i\};\\
0, &\mbox{otherwise}. 
\end{cases}
\end{align*}
It follows that
\begin{align*}
\mathrm{Var}(\widetilde{Z}_{3b})&\leq \frac{C\alpha^2\|\theta\|^4}{\|\theta\|_1^4}\Bigl( \sum_{i,j,\ell,s}\theta_i^3\theta_j\theta_\ell^3\theta_s + \sum_{i,j,\ell,s,j',s'}\theta_i^3\theta_j\theta_\ell^3\theta_s\theta_{j'}\theta_{s'} +\sum_{i,j,\ell,s,j',s'} \theta_i^2\theta^2_j\theta^2_\ell\theta^2_s  \Bigr)\cr
&\leq \frac{C\alpha^2\|\theta\|^4}{\|\theta\|_1^4}\bigl(  \|\theta\|_3^6\|\theta\|_1^2 + \|\theta\|_3^6\|\theta\|_1^4 + \|\theta\|^8 \bigr)\cr
&\leq  C\alpha^2\|\theta\|^4\|\theta\|_3^6. 
\end{align*}
To calculate the variance of $Z^*_{3b}$, we note that $\mathrm{E}[W_{js}W_{jt}W_{\ell i}\cdot W_{j's'}W_{j't'}W_{\ell' i'}]$ is nonzero only if $j'=j$, $\{s',t'\}=\{s,t\}$ and $\{\ell', i'\}=\{\ell,i\}$. Combining it with \eqref{proof-Z3-beta} gives
\begin{align*}
\mathrm{Var}(Z^*_{3b})&\leq \frac{C}{v^2}\sum_{\substack{i,j,\ell (dist)\\s, t (dist) \notin \{j\}}} \beta^2_{ij\ell} \cdot \mathbb{E}[ W^2_{js}W^2_{jt}W^2_{\ell i}]\cr
&\leq \frac{C}{\|\theta\|_1^4}\sum_{i,j,\ell,s,t} (\alpha\|\theta\|^2\theta_i\theta_\ell)^2\cdot\theta_j^2\theta_s\theta_t\theta_\ell\theta_i\cr
&\leq \frac{C\alpha^2\|\theta\|^4}{\|\theta\|_1^4}\sum_{i,j,\ell,s,t}\theta_i^3\theta_j^2\theta^3_{\ell}\theta_s\theta_t\cr
&\leq  \frac{C\alpha^2\|\theta\|^6\|\theta\|_3^6}{\|\theta\|_1^2}. 
\end{align*}
Since $\|\theta\|^6\leq \|\theta\|^4\|\theta\|^2\ll\|\theta\|^4\|\theta\|_1$, the variance of $\widetilde{Z}_{3b}$ dominates the variance of $Z^*_{3b}$. 
Combining the above gives
\beq  \label{proof-Z3-result4}
\mathrm{Var}(Z_{3b})\leq 2\mathrm{Var}(\widetilde{Z}_{3b}) + 2\mathrm{Var}(Z^*_{3b})\leq C\alpha^2\|\theta\|^4\|\theta\|_3^6. 
\eeq

Third, we study $Z_{3c}$. It is seen that
\begin{align*}
Z_{3c}&= \sum_{i, j, k, \ell (dist)}\Bigl(-\frac{1}{\sqrt{v}}\sum_{s\neq i}W_{is}\Bigr)\eta^2_j\Bigl(-\frac{1}{\sqrt{v}}\sum_{t\neq k}W_{kt}\Bigr)   \widetilde{\Omega}_{k \ell}  W_{\ell i}\cr
&= \frac{1}{v}\sum_{\substack{i,k,\ell (dist)\\s\neq i, t\neq k}}\Bigl( \sum_{j\notin\{i,k,\ell\}} \eta_j^2\widetilde{\Omega}_{k\ell} \Bigr)   W_{is}W_{kt}W_{\ell i}\cr
&\equiv \frac{1}{v}\sum_{\substack{i,k,\ell (dist)\\s\neq i, t\neq k}}\beta_{ik\ell}W_{is}W_{kt}W_{\ell i}, 
\end{align*}
where by \eqref{tOmega(ij)} and \eqref{eta-bound},
\beq \label{proof-Z3-beta}
|\beta_{ik\ell}|\leq \sum_{j\notin\{i,k,\ell\}}|\eta^2_j \widetilde{\Omega}_{k\ell}|\leq \sum_j C\alpha\theta_j^2\theta_k\theta_\ell\leq C\alpha\|\theta\|^2\theta_k\theta_\ell. 
\eeq
There are two cases for the indices: $i=\ell$ and $i\neq \ell$.  
We further decompose $Z_{3c}$ into
\[
Z_{3c} = \frac{1}{v}\sum_{\substack{i,k,\ell (dist)\\t\neq k}}\beta_{ik\ell}W^2_{i\ell}W_{kt} + \frac{1}{v}\sum_{\substack{i,k,\ell (dist)\\s\notin\{ i,\ell\}, t\neq k}}\beta_{ik\ell}W_{is}W_{kt}W_{\ell i}\equiv \widetilde{Z}_{3c} + Z^*_{3c}. 
\]
It is easy to see that both terms have zero mean. Hence,
\beq \label{proof-Z3-result5}
\mathbb{E}[Z_{3c}]=0. 
\eeq
To calculate the variance of $\widetilde{Z}_{3c}$, we note that $W^2_{i\ell}W_{kt}$ and $W^2_{i'\ell'}W_{k't'}$ are correlated only when (i) $\{k',t'\}=\{k,t\}$ or (ii) $\{k',t'\}=\{i,\ell\}$ and $\{i',\ell'\}=\{k,t\}$. By direct calculations,
\begin{align*}
& \beta_{ik\ell}\beta_{i'k'\ell'}\cdot \mathbb{E}[W^2_{i\ell}W_{kt}\cdot W^2_{i'\ell'}W_{k't'}]\cr
\leq \; & C\alpha^2\|\theta\|^4\theta_k\theta_{k'}\theta_\ell\theta_{\ell'}\cdot \mathbb{E}[W^2_{i\ell}W_{kt}\cdot W^2_{i'\ell'}W_{k't'}]\cr
\leq \; &\begin{cases}
C\alpha^2\|\theta\|^4\theta_k^2\theta_\ell^2\, \mathbb{E}[W^4_{i\ell}W^2_{kt}]\leq C\alpha^2\|\theta\|^4 \theta_i\theta^3_k\theta^3_\ell\theta_t, &\mbox{if }(k',t')=(k,t),\; (i',\ell')=(i,\ell);\\
C\alpha^2\|\theta\|^4\theta^2_k\theta_\ell\theta_i \, \mathbb{E}[W^4_{i\ell}W^2_{kt}]\leq C\alpha^2\|\theta\|^4 \theta^2_i\theta^3_k\theta^2_\ell\theta_t, &\mbox{if }(k',t')=(k,t),\; (i',\ell')=(\ell,i);\\
C\alpha^2\|\theta\|^4\theta_k\theta_\ell^2\theta_t\, \mathbb{E}[W^4_{i\ell}W^2_{kt}]\leq C\alpha^2\|\theta\|^4 \theta_i\theta^2_k\theta^3_\ell\theta^2_t, &\mbox{if }(k',t')=(t,k),\; (i',\ell')=(i,\ell);\\
C\alpha^2\|\theta\|^4\theta_k\theta_t\theta_\ell\theta_i\, \mathbb{E}[W^4_{i\ell}W^2_{kt}]\leq C\alpha^2\|\theta\|^4 \theta^2_i\theta^2_k\theta^2_\ell\theta^2_t, &\mbox{if }(k',t')=(t,k),\; (i',\ell')=(\ell,i);\\
C\alpha^2\|\theta\|^4\theta^2_k\theta_\ell\theta_{\ell'}\, \mathbb{E}[W^2_{i\ell}W^2_{kt}W^2_{i'\ell'}]\leq C\alpha^2\|\theta\|^4 \theta_i\theta^3_k\theta^2_\ell\theta_t\theta_{i'}\theta^2_{\ell'}, &\mbox{if }(k',t')=(k,t),\; \{i',\ell'\}\neq \{i,\ell\};\\
C\alpha^2\|\theta\|^4\theta_k\theta_t\theta_\ell\theta_{\ell'}\, \mathbb{E}[W^2_{i\ell}W^2_{kt}W^2_{i'\ell'}]\leq C\alpha^2\|\theta\|^4 \theta_i\theta^2_k\theta^2_\ell\theta^2_t\theta_{i'}\theta^2_{\ell'}, &\mbox{if }(k',t')=(t,k),\; \{i',\ell'\}\neq \{i,\ell\};\\
C\alpha^2\|\theta\|^4\theta_k\theta_i\theta_\ell\theta_t\, \mathbb{E}[W^3_{i\ell}W^3_{kt}]\leq C\alpha^2\|\theta\|^4 \theta^2_i\theta^2_k\theta^2_\ell\theta^2_t, &\mbox{if }(k',t')=(i,\ell),\; (i',\ell')=(k,t);\\
C\alpha^2\|\theta\|^4\theta^2_k\theta_i\theta_\ell\, \mathbb{E}[W^3_{i\ell}W^3_{kt}]\leq C\alpha^2\|\theta\|^4 \theta^2_i\theta^3_k\theta^2_\ell\theta_t, &\mbox{if }(k',t')=(i,\ell),\; (i',\ell')=(t,k);\\
C\alpha^2\|\theta\|^4\theta_k\theta^2_\ell\theta_t\, \mathbb{E}[W^3_{i\ell}W^3_{kt}]\leq C\alpha^2\|\theta\|^4 \theta_i\theta^2_k\theta^3_\ell\theta^2_t, &\mbox{if }(k',t')=(\ell,i),\; (i',\ell')=(k,t);\\
C\alpha^2\|\theta\|^4\theta^2_k\theta^2_\ell\, \mathbb{E}[W^3_{i\ell}W^3_{kt}]\leq C\alpha^2\|\theta\|^4 \theta_i\theta^3_k\theta^3_\ell\theta_t, &\mbox{if }(k',t')=(\ell,i),\; (i',\ell')=(t,k);\\
0, &\mbox{otherwise}. 
\end{cases}
\end{align*}
There are only five types on the right hand side. It follows that
\begin{align*}
\mathrm{Var}(\widetilde{Z}_{3c})&\leq \frac{C\alpha^2\|\theta\|^4}{\|\theta\|_1^4}\Bigl(\sum_{i,k,\ell,t}\theta_i\theta^3_k\theta^3_\ell\theta_t+\sum_{i,k,\ell,t} \theta^2_i\theta^3_k\theta^2_\ell\theta_t + \sum_{i,k,\ell,t}\theta^2_i\theta^2_k\theta^2_\ell\theta^2_t\cr
&\qquad +\sum_{i,k,\ell,t,i',\ell'} \theta_i\theta^3_k\theta^2_\ell\theta_t\theta_{i'}\theta^2_{\ell'}+\sum_{i,k,\ell,t,i',\ell'}\theta_i\theta^2_k\theta^2_\ell\theta^2_t\theta_{i'}\theta^2_{\ell'}\Bigr)\cr
&\leq \frac{C\alpha^2\|\theta\|^4}{\|\theta\|_1^4}
\bigl( \|\theta\|_3^6\|\theta\|_1^2 + \|\theta\|^4\|\theta\|_3^3\|\theta\|_1 + \|\theta\|^8 + \|\theta\|^4\|\theta\|_3^3\|\theta\|^3_1 +\|\theta\|^8\|\theta\|^2_1\bigr)\cr
&\leq \frac{C\alpha^2\|\theta\|^8\|\theta\|_3^3}{\|\theta\|_1},
\end{align*}
where the last inequality is obtained as follows: Among the five terms in the brackets, the first and third terms are dominated by the last term, and the second term is dominated by the fourth term; it remains to compare the fourth term and the last term, where the fourth term dominated because $ \|\theta\|^4\leq \|\theta\|_1\|\theta\|_3^3$. To calculate the variance of $Z_{3c}^*$, we write
\[
Z_{3c}^*=  \frac{1}{v}\sum_{i,k,\ell (dist)}\beta_{ik\ell}W_{ik}^2W_{\ell i}+ \frac{1}{v}\sum_{\substack{i,k,\ell (dist)\\s\notin\{ i,\ell\}, t\neq k, (s,t)\neq (k,i)}}\beta_{ik\ell}W_{is}W_{kt}W_{\ell i}. 
\]
Regarding the first term, we note that
\begin{align*}
& \beta_{ik\ell}\beta_{i'k'\ell'}\cdot\mathbb{E}[W_{ik}^2W_{\ell i}\cdot W_{i'k'}^2W_{\ell' i'}]\cr
\leq &C\alpha^2\|\theta\|^4\theta_k\theta_\ell\theta_{k'}\theta_{\ell'}\cdot \mathbb{E}[W_{ik}^2W_{\ell i}\cdot W_{i'k'}^2W_{\ell' i'}]\cr
\leq &
\begin{cases}
C\alpha^2\|\theta\|^4\theta^2_k\theta^2_\ell \, \mathbb{E}[W_{ik}^4W^2_{\ell i}]\leq C\alpha^2\|\theta\|^4\theta^2_i\theta^3_k\theta^3_\ell, &\mbox{if }(\ell',i')=(\ell,i),\, k'= k;\\
C\alpha^2\|\theta\|^4\theta_k\theta^2_\ell\theta_{k'}\, \mathbb{E}[W_{ik}^2W^2_{\ell i}W_{ik'}^2]\leq C\alpha^2\|\theta\|^4\theta^3_i\theta^2_k\theta^3_\ell\theta^2_{k'}, &\mbox{if }(\ell',i')=(\ell,i),\, k'\neq k;\\
C\alpha^2\|\theta\|^4\theta_i \theta_k\theta_\ell\theta_{k'}\, \mathbb{E}[W_{ik}^2W^2_{\ell i}W_{\ell k'}^2]\leq C\alpha^2\|\theta\|^4\theta^3_i\theta^2_k\theta^3_\ell\theta^2_{k'}, &\mbox{if }(\ell',i')=(i, \ell);\\
C\alpha^2\|\theta\|^4\theta^2_k\theta^2_\ell\, \mathbb{E}[W_{ik}^3W^3_{\ell i}]\leq C\alpha^2\|\theta\|^4\theta_i^2\theta_k^3\theta^3_\ell, & \mbox{if }(\ell', i')=(k, i),\, k'=\ell;\\
0, &\mbox{otherwise}. 
\end{cases}
\end{align*}
It follows that 
\begin{align*}
\mathrm{Var}\Bigl(  & \frac{1}{v}\sum_{i,k,\ell (dist)}\beta_{ik\ell}W_{ik}^2W_{\ell i} \Bigr)\leq \frac{C\alpha^2\|\theta\|^4}{\|\theta\|_1^4}\Bigl(\sum_{i,k,\ell}\theta^2_i\theta^3_k\theta^3_\ell+\sum_{i,k,\ell,k'}\theta^3_i\theta^2_k\theta^3_\ell\theta^2_{k'}\Bigr)\cr
&\leq \frac{C\alpha^2\|\theta\|^4}{\|\theta\|_1^4}\bigl( \|\theta\|^2\|\theta\|_3^6 + \|\theta\|^4\|\theta\|_3^6 \bigr)\;
\leq \frac{C\alpha^2\|\theta\|^8\|\theta\|_3^6}{\|\theta\|_1^4}. 
\end{align*}
Regarding the second term, we note that 
\begin{align*}
& \beta_{ik\ell}\beta_{i'k'\ell'}\cdot \mathbb{E}[W_{is}W_{kt}W_{\ell i}\cdot W_{i's'}W_{k't'}W_{\ell' i'}]\cr
\leq \; & C\alpha^2\|\theta\|^4\theta_k\theta_{k'}\theta_\ell\theta_{\ell'}\cdot \mathbb{E}[W_{is}W_{kt}W_{\ell i}\cdot W_{i's'}W_{k't'}W_{\ell' i'}]\cr
\leq \; &\begin{cases}
C\alpha^2\|\theta\|^4\theta^2_k\theta^2_\ell\, \mathbb{E}[W^2_{is}W^2_{kt}W^2_{\ell i}]\leq C\alpha^2\|\theta\|^4\theta^2_i\theta^3_k\theta^3_\ell\theta_s\theta_t, & \mbox{if }(i', s', \ell' )= (i,s,\ell), \; (k', t')=(k,t);\\
C\alpha^2\|\theta\|^4\theta_k\theta_t\theta^2_\ell\, \mathbb{E}[W^2_{is}W^2_{kt}W^2_{\ell i}]\leq C\alpha^2\|\theta\|^4\theta^2_i\theta^2_k\theta^3_\ell\theta_s\theta^2_t, & \mbox{if }(i', s', \ell' )= (i,s,\ell), \; (k', t')=(t,k);\\
C\alpha^2\|\theta\|^4\theta^2_k\theta_\ell\theta_s\, \mathbb{E}[W^2_{is}W^2_{kt}W^2_{\ell i}]\leq C\alpha^2\|\theta\|^4\theta^2_i\theta^3_k\theta^2_\ell\theta^2_s\theta_t, & \mbox{if }(i', s', \ell' )= (i,\ell, s), \; (k', t')=(k,t);\\
C\alpha^2\|\theta\|^4\theta_k\theta_t \theta_\ell\theta_s\, \mathbb{E}[W^2_{is}W^2_{kt}W^2_{\ell i}]\leq C\alpha^2\|\theta\|^4\theta^2_i\theta^2_k\theta^2_\ell\theta^2_s\theta^2_t, & \mbox{if }(i', s', \ell' )= (i,\ell, s), \; (k', t')=(t,k);\\
0, &\mbox{otherwise}. 
\end{cases}
\end{align*}
It follows that 
\begin{align*}
\mathrm{Var}\Bigl(\frac{1}{v}\sum_{\substack{i,k,\ell (dist)\\s\notin\{ i,\ell\}, t\neq k,\\ (s,t)\neq (k,i)}}&\beta_{ik\ell}W_{is}W_{kt}W_{\ell i}\Bigr) \leq \frac{C\alpha^2\|\theta\|^4}{\|\theta\|_1^4} \sum_{\substack{i,k,\ell,\\s, t}}(\theta^2_i\theta^3_k\theta^3_\ell\theta_s\theta_t+\theta^2_i\theta^2_k\theta^3_\ell\theta_s\theta^2_t + \theta^2_i\theta^2_k\theta^2_\ell\theta^2_s\theta^2_t)\cr
&\leq \frac{C\alpha^2\|\theta\|^4}{\|\theta\|_1^4}\bigl( \|\theta\|^2\|\theta\|_3^6\|\theta\|_1^2 + \|\theta\|^6\|\theta\|_3^3\|\theta\|_1 + \|\theta\|^{10} \bigr)\cr
&\leq  \frac{C\alpha^2\|\theta\|^6\|\theta\|_3^6}{\|\theta\|_1^2}. 
\end{align*}
We plug the above results into $Z_{3c}^*$. Since $\|\theta\|^2\leq \|\theta\|_1\theta_{\max}\ll\|\theta\|_1^2$, we have $\frac{C\alpha^2\|\theta\|^8\|\theta\|_3^6}{\|\theta\|_1^4}\ll \frac{C\alpha^2\|\theta\|^6\|\theta\|_3^6}{\|\theta\|_1^2}$. It follows that
\[
\mathrm{Var}(Z^*_{3c})\leq  \frac{C\alpha^2\|\theta\|^6\|\theta\|_3^6}{\|\theta\|_1^2}. 
\]
Since $\|\theta\|_3^6\ll\|\theta\|_3^3\|\theta\|_1$, the variance of $Z^*_{3c}$ is dominated by the variance of $\widetilde{Z}_{3c}$. It follows that 
\beq  \label{proof-Z3-result6}
\mathrm{Var}(Z_{3c})\leq 2\mathrm{Var}(\widetilde{Z}_{3c})+2\mathrm{Var}(Z^*_{3c})\leq \frac{C\alpha^2\|\theta\|^8\|\theta\|_3^3}{\|\theta\|_1}. 
\eeq

Last, we study $Z_{3d}$. In the definition of $Z_{3d}$, if we switch the two indices $(j,k)$, then it becomes
\[
Z_{3d} = \sum_{\substack{i,j,k,\ell\\ (dist)}}(\eta_i-\teta_i)\eta_k (\eta_k-\teta_k)\eta_j   \widetilde{\Omega}_{j \ell}  W_{\ell i}=
\sum_{\substack{i,j,k,\ell\\ (dist)}}(\eta_k\eta_j\widetilde{\Omega}_{j\ell})(\eta_i-\teta_i)(\eta_k-\teta_k). 
\]
At the same time, we recall that
\[
Z_{3c} = \sum_{\substack{i,j,k,\ell\\ (dist)}}(\eta_i-\teta_i)\eta_j^2(\eta_k-\teta_k)\widetilde{\Omega}_{k\ell}W_{\ell i}
=\sum_{\substack{i,j,k,\ell\\ (dist)}}(\eta_j^2\widetilde{\Omega}_{k\ell})(\eta_i-\teta_i)(\eta_k-\teta_k). 
\]
Here, $Z_{3d}$ has a similar structure as $Z_{3c}$. Moreover, in deriving the bound for $\mathrm{Var}(Z_{3c})$, we have used $|\eta_j^2\widetilde{\Omega}_{k\ell}|\leq C\alpha\theta_j^2\theta_k\theta_\ell$. In the expression of $Z_{3d}$ above, we also have $|\eta_k\eta_j\widetilde{\Omega}_{j\ell}|\leq C\alpha\theta_j^2\theta_k\theta_\ell$. Therefore, we can use \eqref{proof-Z3-result5} and \eqref{proof-Z3-result6} to directly get
\beq \label{proof-Z3-result7} 
\mathbb{E}[Z_{3d}]=0, \qquad \mathrm{Var}(Z_{3d})\leq \frac{C\alpha^2\|\theta\|^8\|\theta\|_3^3}{\|\theta\|_1}
\eeq

Now, we combine \eqref{proof-Z3-result1}, \eqref{proof-Z3-result3}, \eqref{proof-Z3-result5} and \eqref{proof-Z3-result6} to get
\[
\mathbb{E}[Z_3]=0.
\]
We also combine \eqref{proof-Z3-result2}, \eqref{proof-Z3-result4}, \eqref{proof-Z3-result6}-\eqref{proof-Z3-result7}. Since $\|\theta\|^4\leq \|\theta\|_1\|\theta\|_3^3$, the right hand side of \eqref{proof-Z3-result6}-\eqref{proof-Z3-result7} is dominated by the right hand side of \eqref{proof-Z3-result4}; since $\|\theta\|_3^6\ll\|\theta\|_1^2$, the right hand side of \eqref{proof-Z3-result2} is negligible to the right hand side of \eqref{proof-Z3-result4}. It follows that
\[
\mathrm{Var}(Z_3)\leq C\alpha^2\|\theta\|^4\|\theta\|_3^6. 
\]
This proves the claims of $Z_3$. 

Next, we analyze $Z_4$. Since $\delta_{ij}=\eta_i(\eta_j-\teta_j)+\eta_j(\eta_i-\teta_i)$, 
\begin{align*} 
Z_4 &= \sum_{i, j, k, \ell(dist)} \eta_i(\eta_j-\teta_j)  \widetilde{\Omega}_{jk} \eta_k (\eta_\ell - \teta_\ell )W_{\ell i}
+  \sum_{i, j, k, \ell(dist)} \eta_i(\eta_j-\teta_j)  \widetilde{\Omega}_{jk}(\eta_k - \teta_k )\eta_\ell W_{\ell i}\cr
&+ \sum_{i, j, k, \ell(dist)} (\eta_i-\teta_i)\eta_j  \widetilde{\Omega}_{jk} \eta_k (\eta_\ell - \teta_\ell )W_{\ell i}
+ \sum_{i, j, k, \ell(dist)} (\eta_i-\teta_i)\eta_j \widetilde{\Omega}_{jk}(\eta_k - \teta_k )\eta_\ell W_{\ell i}.
\end{align*}
If we relabel $(i,j,k,\ell)$ as $(\ell',k', j',i')$ in the last sum, it is equal to the first sum. Therefore, 
\begin{align} \label{proof-Z4-decompose}
Z_4 &= 2\sum_{i, j, k, \ell(dist)}  \eta_i(\eta_j-\teta_j)  \widetilde{\Omega}_{jk} \eta_k (\eta_\ell - \teta_\ell )W_{\ell i}\cr
&+ \sum_{i, j, k, \ell(dist)} \eta_i(\eta_j-\teta_j)  \widetilde{\Omega}_{jk}(\eta_k - \teta_k )\eta_\ell W_{\ell i}\cr
& +  \sum_{i, j, k, \ell(dist)} (\eta_i-\teta_i)\eta_j  \widetilde{\Omega}_{jk} \eta_k (\eta_\ell - \teta_\ell )W_{\ell i}\cr
&\equiv Z_{4a}+Z_{4b}+Z_{4c}. 
\end{align}

First, we study $Z_{4a}$ and $Z_{4b}$. We show that they have the same structure as $Z_{3c}$ and $Z_{3a}$, respectively. In $Z_{4a}$, by relabeling $(i,j,k,\ell)$ as $(\ell,k,j,i)$, we have
\[
Z_{4a} = 2\sum_{\substack{i, j, k, \ell\\(dist)}} \eta_\ell (\eta_k-\teta_k)  \widetilde{\Omega}_{kj} \eta_j (\eta_i - \teta_i )W_{i \ell} = 2\sum_{\substack{i,j,k,\ell\\(dist)}}(\eta_j\eta_{\ell}\widetilde{\Omega}_{kj})(\eta_i - \teta_i )(\eta_k-\teta_k)W_{\ell i}.
\]
At the same time, we recall the definition of $Z_{3c}$ in \eqref{proof-Z3-decompose}:
\[
Z_{3c} = \sum_{\substack{i, j, k, \ell\\(dist)}}(\eta_i-\teta_i)\eta_j^2(\eta_k-\teta_k)\widetilde{\Omega}_{k\ell}W_{\ell i} = \sum_{\substack{i,j,k,\ell\\(dist)}}(\eta^2_j\widetilde{\Omega}_{k\ell})(\eta_i - \teta_i )(\eta_k-\teta_k)W_{\ell i}. 
\]
It is seen that $Z_{4a}$ has a similar structure as $Z_{3c}$ does. Also, by \eqref{tOmega(ij)} and \eqref{eta-bound}, in the expression of  $Z_{4a}$, we have $|\eta_j\eta_{\ell}\widetilde{\Omega}_{kj}|\leq C\alpha\theta_j^2\theta_k\theta_\ell$, while in the expression of $Z_{3d}$, we have $|\eta^2_j\widetilde{\Omega}_{k\ell}|\leq C\alpha\theta_j^2\theta_k\theta_{\ell}$. As a result, if we use similar calculation as before, we will get the same conclusion for $Z_{4a}$ and $Z_{3d}$. Hence, we use \eqref{proof-Z3-result5}-\eqref{proof-Z3-result6} to conclude that
\beq \label{proof-Z4-result1}
\mathbb{E}[Z_{4a}]=0, \qquad \mathrm{Var}(Z_{4a})\leq \frac{C\alpha^2\|\theta\|^8 \|\theta\|_3^3}{\|\theta\|_1}.
\eeq 
For $Z_{4b}$, we note that 
\[
Z_{4b} = \sum_{\substack{i, j, k, \ell\\ (dist)}} \eta_i(\eta_j-\teta_j)  \widetilde{\Omega}_{jk}(\eta_k - \teta_k )\eta_\ell W_{\ell i}
= \sum_{\substack{i, j, k, \ell\\ (dist)}}(\eta_i\eta_\ell \widetilde{\Omega}_{jk})(\eta_j-\teta_j)(\eta_k-\teta_k)W_{\ell i},  
\]
where $|\eta_i\eta_\ell \widetilde{\Omega}_{jk}|\leq C\alpha\theta_i\theta_j\theta_k\theta_\ell$. 
At the same time, we recall the definition of $Z_{3a}$ in \eqref{proof-Z3-decompose}:
\[
Z_{3a} = \sum_{\substack{i, j, k, \ell\\ (dist)}} \eta_i(\eta_j-\teta_j) \eta_j (\eta_k - \teta_k )\widetilde{\Omega}_{k\ell} W_{\ell i}
= \sum_{\substack{i, j, k, \ell\\ (dist)}}(\eta_i\eta_j \widetilde{\Omega}_{k\ell})(\eta_j-\teta_j)(\eta_k-\teta_k)W_{\ell i},
\]
where $|\eta_i\eta_j \widetilde{\Omega}_{k\ell}|\leq C\alpha\theta_i\theta_j\theta_k\theta_\ell$. Therefore, we can quote the results for $Z_{3a}$ in \eqref{proof-Z3-result1}-\eqref{proof-Z3-result2} to get 
\beq \label{proof-Z4-result2}
\mathbb{E}[Z_{4b}]=0, \qquad \mathrm{Var}(Z_{4b})\leq \frac{C\alpha^2\|\theta\|_3^{12}}{\|\theta\|^2_1}.
\eeq 

Second, we study $Z_{4c}$. It is seen that
\begin{align*}
Z_{4c} &= \sum_{i, j, k, \ell (dist)} \Bigl(-\frac{1}{\sqrt{v}}\sum_{s\neq i}W_{is}\Bigr)\eta_j  \widetilde{\Omega}_{jk} \eta_k \Bigl( -\frac{1}{\sqrt{v}}\sum_{t\neq \ell }W_{\ell t} \Bigr)W_{\ell i}\cr
&=\frac{1}{v}\sum_{\substack{i, \ell (dist)\\s\neq i, t\neq \ell}} \Bigr( \sum_{j, k (dist)\notin\{i,\ell\}} \eta_j\eta_k\widetilde{\Omega}_{jk}\Bigr)   W_{is}W_{\ell t}W_{\ell i}\cr
&\equiv \frac{1}{v}\sum_{\substack{i, \ell (dist)\\s\neq i, t\neq \ell}} \beta_{i\ell}  W_{is}W_{\ell t}W_{\ell i},
\end{align*}
where
\beq \label{proof-Z4-beta2}
|\beta_{i\ell}|\leq \sum_{j, k (dist)\notin\{i,\ell\}} |\eta_j\eta_k\widetilde{\Omega}_{jk}|\leq \sum_{j,k} C\alpha\theta_j^2\theta_k^2\leq C\alpha\|\theta\|^4. 
\eeq
We divide the summands into four groups: (i) $s=\ell$, $t=i$; (ii) $s=\ell$, $t\neq i$; (iii) $s\neq \ell$, $t=i$; (iv) $s\neq \ell$, $t\neq i$. By symmetry, the sum of group (ii) and the sum of group (iii) are equal. It yields that 
\begin{align*}
Z_{4c} &
= \frac{1}{v}\sum_{i, \ell (dist)} \beta_{i\ell}  W^3_{\ell i}  + \frac{2}{v}\sum_{\substack{i, \ell (dist)\\s\notin\{i,\ell\}}} \beta_{i\ell}  W_{is}W^2_{\ell i} + \frac{1}{v}\sum_{\substack{i, \ell (dist)\\s\notin\{ i,\ell\}, t\notin\{\ell, i\}}} \beta_{i\ell}  W_{is}W_{\ell t}W_{\ell i}\cr
&\equiv \widetilde{Z}_{4c} + Z^*_{4c} + Z^{\dag}_{4c}. 
\end{align*}
Only $\widetilde{Z}_{4c}$ has a nonzero mean. By \eqref{v-order} and \eqref{proof-Z4-beta2},
\beq \label{proof-Z4-result5}
\bigl|\mathbb{E}[Z_{4c}]\bigr| = \bigl|\mathbb{E}[\widetilde{Z}_{4c}]\bigr| 
\leq \frac{C}{\|\theta\|_1^2} \sum_{i,\ell} \alpha\|\theta\|^4\theta_i\theta_\ell \leq C\alpha\|\theta\|^4. 
\eeq
We now compute the variances of these terms. It is seen that
\begin{align*}
\mathrm{Var}(\widetilde{Z}_{4c}) &\leq \frac{C}{v^2}\sum_{i,\ell (dist)}\beta^2_{i\ell}\mathrm{Var}(W^3_{i\ell})
\leq \frac{C\alpha^2\|\theta\|^8}{\|\theta\|_1^4}\sum_{i,\ell}\theta_i\theta_\ell \leq \frac{C\alpha^2\|\theta\|^8}{\|\theta\|_1^2}. 
\end{align*}
For $Z_{4c}^*$, by direct calculations, 
\begin{align*}
& \beta_{i\ell}\beta_{i'\ell'}\cdot\mathbb{E}[W_{is}W^2_{\ell i}\cdot W_{i's'}W^2_{\ell'i'}]\cr
\leq \; & C\alpha^2\|\theta\|^8\cdot \mathbb{E}[W_{is}W^2_{\ell i}\cdot W_{i's'}W^2_{\ell'i'}]\cr
\leq \; &\begin{cases}
C\alpha^2\|\theta\|^8\cdot\mathbb{E}[W^2_{is}W^4_{\ell i}]\leq C\alpha^2\|\theta\|^8\theta_i^2\theta_\ell\theta_s, &\mbox{if }i'=i, \,s'=s,\, \ell'= \ell;\\
C\alpha^2\|\theta\|^8\cdot\mathbb{E}[W^2_{is}W^2_{\ell i}W^2_{\ell'i}]\leq C\alpha^2\|\theta\|^8\theta_i^3\theta_\ell\theta_s\theta_{\ell'}, &\mbox{if }i'=i, \; s'=s,\; \ell'\neq \ell;\\
C\alpha^2\|\theta\|^8\cdot\mathbb{E}[W^3_{is}W^3_{\ell i}]\leq C\alpha^2\|\theta\|^8\theta_i^2\theta_\ell\theta_s, &\mbox{if }i'=i, \, s'=\ell,\; \ell'= s;\\
0, &\mbox{otherwise}. 
\end{cases}
\end{align*}
It follows that
\begin{align*}
\mathrm{Var}(Z_{4c}^*)&\leq \frac{C\alpha^2\|\theta\|^8}{\|\theta\|_1^4}\Bigl(\sum_{i, \ell, s} \theta_i^2\theta_\ell\theta_s + \sum_{i,\ell, s, \ell'} \theta_i^3\theta_\ell\theta_s\theta_{\ell'} \Bigr)\cr
&\leq \frac{C\alpha^2\|\theta\|^8}{\|\theta\|_1^4}\bigl( \|\theta\|^2\|\theta\|_1^2 + \|\theta\|_3^3\|\theta\|_1^3 \bigr)\cr
&\leq \frac{C\alpha^2\|\theta\|^8\|\theta\|_3^3}{\|\theta\|_1},
\end{align*}
where, to get the last line, we have used $\|\theta\|^2\ll\|\theta\|^4\leq \|\theta\|_1\|\theta\|_3^3$. Regarding the variance of $Z^{\dag}_{4c}$, we note that $W_{is}W_{\ell t}W_{\ell i}$ and $W_{i's'}W_{\ell' t'}W_{\ell' i'}$ are correlated only when the two undirected paths $s$-$i$-$\ell$-$t$ and $s'$-$i'$-$\ell'$-$t'$ are the same. Mimicking the argument in \eqref{proof-Y2-result1} or \eqref{proof-Y3-result2}, we can derive that
\begin{align*}
\mathrm{Var}(Z^{\dag}_{4c}) &\leq \frac{C}{v^2}\sum_{\substack{i, \ell (dist)\\s\notin\{ i,\ell\}, t\notin\{\ell, i\}}} \beta^2_{i\ell}\cdot\mathrm{Var}( W_{is}W_{\ell t}W_{\ell i})\cr
&\leq \frac{C\alpha^2\|\theta\|^8}{\|\theta\|_1^4}\sum_{i,\ell,s,t} \theta_i^2\theta_\ell^2\theta_s\theta_t\cr
&\leq  \frac{C\alpha^2\|\theta\|^{12}}{\|\theta\|_1^2}. 
\end{align*}
Since $\|\theta\|^4\leq \|\theta\|_1\|\theta\|_3^3$, the variance of $Z_{4c}^{\dag}$ is dominated by the variance of $Z^*_{4c}$. Since $\|\theta\|\to\infty$, we have $\|\theta\|^3_3\gg1/\|\theta\|_1$; it follows that the variance of $\widetilde{Z}_{4c}$ is dominated by the variance of $Z_{4c}^*$. Combining the above gives
\beq \label{proof-Z4-result6}
\mathrm{Var}(Z_{4c})\leq 3\mathrm{Var}(\widetilde{Z}_{4c}) + 3\mathrm{Var}(Z_{4c}^*) + 3\mathrm{Var}(Z^{\dag}_{4c}) \leq  \frac{C\alpha^2\|\theta\|^8\|\theta\|_3^3}{\|\theta\|_1}. 
\eeq

We combine \eqref{proof-Z4-result1}, \eqref{proof-Z4-result2} and \eqref{proof-Z4-result5} to get
\[
\bigl|\mathbb{E}[Z_4]\bigr|\leq C\alpha\|\theta\|^4 = o(\alpha^4\|\theta\|^8). 
\]
We then combine \eqref{proof-Z4-result1}, \eqref{proof-Z4-result2} and \eqref{proof-Z4-result6}. Since $\|\theta\|_3^6\leq (\theta_{\max}^2\|\theta\|_1)(\theta_{\max}\|\theta\|^2)=o(\|\theta\|_1\|\theta\|^2)$, the variance of $Z_{4b}$ is negligible compared to the variances of $Z_{4a}$ and $Z_{4c}$. It follows that
\[
\mathrm{Var}(Z_4)\leq \frac{C\alpha^2\|\theta\|^8\|\theta\|_3^3}{\|\theta\|_1} = o(\|\theta\|^8). 
\]
This proves the claims of $Z_4$.

Next, we analyze $Z_5$. By plugging in the definition of $\delta_{ij}$, we have
\begin{align} \label{proof-Z5-decompose}
Z_5 &= \sum_{i, j, k, \ell (dist)} \eta_i(\eta_j-\teta_j)\eta_j(\eta_k-\teta_k) \widetilde{\Omega}_{k \ell}  \widetilde{\Omega}_{\ell i} + \sum_{i, j, k, \ell (dist)} \eta_i(\eta_j-\teta_j)^2 \eta_k \widetilde{\Omega}_{k \ell}  \widetilde{\Omega}_{\ell i}\cr
&\qquad  +\sum_{i, j, k, \ell (dist)} (\eta_i-\teta_i)\eta^2_j(\eta_k-\teta_k) \widetilde{\Omega}_{k \ell}  \widetilde{\Omega}_{\ell i} + \sum_{i, j, k, \ell (dist)} (\eta_i-\teta_i)\eta_j (\eta_j-\teta_j) \eta_k \widetilde{\Omega}_{k \ell}  \widetilde{\Omega}_{\ell i}\cr
&= 2\sum_{i, j, k, \ell (dist)} \eta_i(\eta_j-\teta_j)\eta_j(\eta_k-\teta_k) \widetilde{\Omega}_{k \ell}  \widetilde{\Omega}_{\ell i} + \sum_{i, j, k, \ell (dist)} \eta_i(\eta_j-\teta_j)^2 \eta_k \widetilde{\Omega}_{k \ell}  \widetilde{\Omega}_{\ell i}\cr
&\qquad  +\sum_{i, j, k, \ell (dist)} (\eta_i-\teta_i)\eta^2_j(\eta_k-\teta_k) \widetilde{\Omega}_{k \ell}  \widetilde{\Omega}_{\ell i} \cr
&\equiv Z_{5a} + Z_{5b} + Z_{5c}. 
\end{align}

First, we study $Z_{5a}$. By definition, $(\teta_i-\eta_i)$ has the expression in \eqref{tilde-eta}. It follows that
\begin{align*}
Z_{5a} &= 2\sum_{i, j, k, \ell (dist)} \eta_i\Bigl(-\frac{1}{\sqrt{v}}\sum_{s\neq j}W_{js}\Bigr)\eta_j\Bigl(-\frac{1}{\sqrt{v}}\sum_{t\neq k}W_{kt}\Bigr) \widetilde{\Omega}_{k \ell}  \widetilde{\Omega}_{\ell i}\cr
&= \frac{2}{v} \sum_{\substack{j, k (dist)\\s\neq j, t\neq k}} \Bigl( \sum_{i,\ell (dist)\notin\{j,k\}}\eta_i\eta_j\widetilde{\Omega}_{k \ell}  \widetilde{\Omega}_{\ell i}\Bigr)  W_{js}W_{kt}\cr
&\equiv \frac{2}{v}\sum_{\substack{j, k (dist)\\s\neq j, t\neq k}}\beta_{jk}W_{js}W_{kt}, 
\end{align*}
where
\beq\label{proof-Z5-beta}
|\beta_{jk}| \leq   \sum_{i,\ell (dist)\notin\{j,k\}}|\eta_i\eta_j\widetilde{\Omega}_{k \ell}  \widetilde{\Omega}_{\ell i}|\leq \sum_{i,\ell}(C\theta_i\theta_j)(C\alpha^2\theta_k\theta_\ell^2\theta_i)\leq C\alpha^2\|\theta\|^4\theta_j\theta_k. 
\eeq
In $Z_{5a}$, the summand has a nonzero mean only if $(s, t)=(k,j)$. We further decompose $Z_{5a}$ into 
\[
Z_{5a} =  \frac{2}{v}\sum_{j, k (dist)}\beta_{jk}W_{jk}^2  + \frac{2}{v}\sum_{\substack{j, k (dist)\\s\neq j, t\neq k,\\(s,t)\neq (k,j)}}\beta_{jk}W_{js}W_{kt} \equiv \widetilde{Z}_{5a} + Z^*_{5a}. 
\]
Only the first term has a nonzero mean. By \eqref{v-order} and \eqref{proof-Z5-beta}, we have
\beq \label{proof-Z5-result1}
\bigl|\mathbb{E}[Z_{5a}]\bigr| = \bigl| \mathbb{E}[\widetilde{Z}_{5a}] \bigr|
\leq \frac{C}{\|\theta\|_1^2}\sum_{j,k}(\alpha^2\|\theta\|^4\theta_j\theta_k)(\theta_j\theta_k)\leq \frac{C\alpha^2\|\theta\|^8}{\|\theta\|_1^2}. 
\eeq
We then compute the variances. In each of $\widetilde{Z}_{5a}$ and $Z^*_{5a}$, two summands are uncorrelated unless they are exactly the same variables (e.g., when $(j', k')=(k, j)$ in $\widetilde{Z}_{5a}$). Mimicking the argument in \eqref{proof-Y2-result1} or \eqref{proof-Y3-result2}, we can derive that
\begin{align*}
\mathrm{Var}(\widetilde{Z}_{5a}) &\leq \frac{C}{v^2}\sum_{j,k (dist)}\beta_{jk}^2\,\mathrm{Var}(W^2_{jk})\leq \frac{C\alpha^4\|\theta\|^8}{\|\theta\|_1^4}\sum_{j,k}(\theta_j^2\theta_k^2)\theta_j\theta_k \leq \frac{C\alpha^4\|\theta\|^8\|\theta\|_3^6}{\|\theta\|_1^4},\cr
\mathrm{Var}(Z^*_{5a})&\leq \frac{C}{v^2}\sum_{\substack{j, k (dist)\\s\neq j, t\neq k,\\(s,t)\neq (k,j)}} \beta^2_{jk}\, \mathrm{Var}(W_{js}W_{kt})\leq \frac{C\alpha^4\|\theta\|^8}{\|\theta\|_1^4} \sum_{j,k}(\theta_j^2\theta_k^2)\theta_j\theta_s\theta_k\theta_t\leq  \frac{C\alpha^4\|\theta\|^8\|\theta\|_3^6}{\|\theta\|_1^2}.
\end{align*}
It immediately leads to
\beq  \label{proof-Z5-result2}
\mathrm{Var}(Z_{5a}) \leq 2\mathrm{Var}(\widetilde{Z}_{5a}) + 2\mathrm{Var}(Z^*_{5a})\leq  \frac{C\alpha^4\|\theta\|^8\|\theta\|_3^6}{\|\theta\|_1^2}. 
\eeq

Second, we study $Z_{5b}$. It is seen that
\begin{align*}
Z_{5b} &= \sum_{i, j, k, \ell (dist)} \eta_i\Bigl(-\frac{1}{\sqrt{v}}\sum_{s\neq j}W_{js}\Bigr)\Bigl(-\frac{1}{\sqrt{v}}\sum_{t\neq j}W_{jt}\Bigr) \eta_k \widetilde{\Omega}_{k \ell}  \widetilde{\Omega}_{\ell i}\cr
&=\frac{1}{v}\sum_{j, s\neq j,t\neq j} \Bigl(\sum_{i,k,\ell (dist)\notin\{j\}}\eta_i\eta_k \widetilde{\Omega}_{k \ell}  \widetilde{\Omega}_{\ell i} \Bigr) W_{js}W_{jt}\cr
& \equiv  \frac{1}{v}\sum_{j, s\neq j,t\neq j} \beta_j W_{js}W_{jt},
\end{align*}
where
\beq\label{proof-Z5-beta2}
|\beta_j|\leq  \sum_{i,k,\ell (dist)\notin\{j\}}|\eta_i\eta_k \widetilde{\Omega}_{k \ell}  \widetilde{\Omega}_{\ell i}|\leq \sum_{i, k, \ell}(C\theta_i\theta_k)(C\alpha^2\theta_i\theta_k\theta_{\ell}^2)\leq C\alpha^2\|\theta\|^6. 
\eeq
In $Z_{5b}$, the summand has a nonzero mean only if $s=t$. We further decompose $Z_{5b}$ into
\[
Z_{5b} = \frac{1}{v}\sum_{j, s (dist)} \beta_j W^2_{js} + \frac{1}{v}\sum_{\substack{j\\ s, t (dist) \notin\{j\}}} \beta_j W_{js}W_{jt}\equiv \widetilde{Z}_{5b} + Z^*_{5b}. 
\]
Only $\widetilde{Z}_{5b}$ has a nonzero mean. By \eqref{v-order} and \eqref{proof-Z5-beta2}, 
\beq \label{proof-Z5-result3}
\bigl|\mathbb{E}[Z_{5b}]\bigr| = \bigl|\mathbb{E}[\widetilde{Z}_{5b}]\bigr|\leq \frac{C}{\|\theta\|_1^2} \sum_{j,s} (\alpha^2\|\theta\|^6) \theta_j\theta_s \leq C\alpha^2\|\theta\|^6. 
\eeq
To compute the variance, we note that in each of $\widetilde{Z}_{5b}$ and $Z^*_{5b}$, two summands are uncorrelated unless they are exactly the same random variables (e.g., when $\{j', s'\}=\{s, j\}$ in $\widetilde{Z}_{5b}$, and when $j' = j$ and $\{s', t'\}=\{s,t\}$ in $Z^*_{5b}$). Mimicking the argument in \eqref{proof-Y2-result1} or \eqref{proof-Y3-result2}, we can derive that
\begin{align*}
\mathrm{Var}(\widetilde{Z}_{5b}) &\leq \frac{C}{v^2}\sum_{j, s (dist)} \beta^2_j\,\mathrm{Var}( W^2_{js})\leq \frac{C\alpha^4\|\theta\|^{12}}{\|\theta\|_1^4}\sum_{j,s}\theta_j\theta_s\leq \frac{C\alpha^4\|\theta\|^{12}}{\|\theta\|_1^2},\cr
\mathrm{Var}(Z^*_{5b}) &\leq \frac{C}{v^2}\sum_{\substack{j\\ s, t (dist) \notin\{j\}}}  \beta^2_j\,\mathrm{Var}( W_{js}W_{jt})\leq \frac{C\alpha^4\|\theta\|^{12}}{\|\theta\|_1^4}\sum_{j,s,t}\theta^2_j\theta_s\theta_t \leq \frac{C\alpha^4\|\theta\|^{14}}{\|\theta\|_1^2}. 
\end{align*}
Combining the above gives
\beq  \label{proof-Z5-result4}
\mathrm{Var}(Z_{5b}) \leq 2\mathrm{Var}(\widetilde{Z}_{5b}) + 2\mathrm{Var}(Z^*_{5b})\leq  \frac{C\alpha^4\|\theta\|^{14}}{\|\theta\|_1^2}. 
\eeq

Third, we study $Z_{5c}$. If we relabel $(i, j,k,\ell)=(j,i,k,\ell)$, then $Z_{5c}$ becomes
\[
Z_{5c} = \sum_{\substack{i, j, k, \ell\\ (dist)}} (\eta_j-\teta_j)\eta^2_i(\eta_k-\teta_k) \widetilde{\Omega}_{k \ell}  \widetilde{\Omega}_{\ell j} = \sum_{\substack{i, j, k, \ell\\ (dist)}}(\eta_i^2\widetilde{\Omega}_{k\ell}\widetilde{\Omega}_{\ell j})(\eta_j-\teta_j)(\eta_k-\teta_k),
\]
where $|\eta_i^2\widetilde{\Omega}_{k\ell}\widetilde{\Omega}_{\ell j}|\leq C\alpha^2\theta^2_i\theta_j\theta_k\theta_\ell^2$. At the same time, we recall that
\[
Z_{5a} = 2\sum_{\substack{i, j, k, \ell\\ (dist)}} \eta_i (\eta_j-\teta_j)\eta_j(\eta_k-\teta_k) \widetilde{\Omega}_{k \ell}  \widetilde{\Omega}_{\ell i} = \sum_{\substack{i, j, k, \ell\\ (dist)}}(\eta_i\eta_j\widetilde{\Omega}_{k\ell}\widetilde{\Omega}_{\ell i})(\eta_j-\teta_j)(\eta_k-\teta_k),
\]
where $|\eta_i\eta_j\widetilde{\Omega}_{k\ell}\widetilde{\Omega}_{\ell i}|\leq C\alpha^2\theta^2_i\theta_j\theta_k\theta^2_\ell$. It is easy to see that $Z_{5c}$ has a similar structure as $Z_{5c}$. As a result, from \eqref{proof-Z5-result1}-\eqref{proof-Z5-result2}, we immediately have
\beq  \label{proof-Z5-result5}
\bigl|\mathbb{E}[Z_{5c}]\bigr| \leq \frac{C\alpha^2\|\theta\|^8}{\|\theta\|_1^2}, \qquad  \mathrm{Var}(Z_{5c}) \leq \frac{C\alpha^4\|\theta\|^8\|\theta\|_3^6}{\|\theta\|_1^2}. 
\eeq

We now combine the results for $Z_{5a}$-$Z_{5c}$. Since $\|\theta\|^2\leq \theta_{\max}\|\theta\|_1\ll\|\theta\|_1^2$, $\mathbb{E}[Z_{5a}]$ and $\mathbb{E}[Z_{5c}]$ are of a smaller order than the the right hand side of \eqref{proof-Z5-result3}. Since $\|\theta\|_3^6\leq \theta_{\max}^2\|\theta\|^4\ll\|\theta\|^6$, $\mathrm{Var}(Z_{5a})$ and $\mathrm{Var}(Z_{5c})$ are of a smaller order than the right hand side of \eqref{proof-Z5-result4}. It follows that
\[
\bigl|\mathbb{E}[Z_5]\bigr|\leq C\alpha^2\|\theta\|^6=o(\alpha^4\|\theta\|^8), \qquad \mathrm{Var}(Z_5)\leq \frac{C\alpha^4\|\theta\|^{14}}{\|\theta\|_1^2} = o(\alpha^6\|\theta\|^8\|\theta\|_3^6). 
\]
We briefly explain why $\mathrm{Var}(Z_5) = o(\alpha^6\|\theta\|^8\|\theta\|_3^6)$: since $\|\theta\|^4\leq \|\theta\|_1\|\theta\|_3^3$, we immediately have $\|\theta\|^{14}\leq \|\theta\|^6(\|\theta\|_1\|\theta\|_3^3)^2$; it follows that the bound for $\mathrm{Var}(Z_5)$ is $\leq C\alpha^4\|\theta\|^6\|\theta\|_3^6$; note that $\alpha\|\theta\|\to\infty$, we immediately have $\alpha^4\|\theta\|^6\|\theta\|_3^6=o(\alpha^6\|\theta\|^8\|\theta\|_3^6)$. This proves the claims of $Z_5$.

Last, we analyze $Z_6$. Plugging in the definition of $\delta_{ij}$, we have
\begin{align*}
Z_6 &=  \sum_{i, j, k, \ell (dist)} \eta_i(\eta_j-\teta_j)  \widetilde{\Omega}_{jk}\eta_k(\eta_\ell-\teta_\ell)  \widetilde{\Omega}_{\ell i} +  \sum_{i, j, k, \ell (dist)} \eta_i(\eta_j-\teta_j)  \widetilde{\Omega}_{jk}(\eta_k-\teta_k)\eta_\ell  \widetilde{\Omega}_{\ell i} \cr
&\qquad + \sum_{i, j, k, \ell (dist)} (\eta_i-\teta_i)\eta_j  \widetilde{\Omega}_{jk}\eta_k(\eta_\ell-\teta_\ell)  \widetilde{\Omega}_{\ell i} +  \sum_{i, j, k, \ell (dist)} (\eta_i-\teta_i)\eta_j  \widetilde{\Omega}_{jk}(\eta_k-\teta_k)\eta_\ell  \widetilde{\Omega}_{\ell i} \cr
&=  2\sum_{i, j, k, \ell (dist)} \eta_i(\eta_j-\teta_j)  \widetilde{\Omega}_{jk}\eta_k(\eta_\ell-\teta_\ell)  \widetilde{\Omega}_{\ell i} +  2\sum_{i, j, k, \ell (dist)} \eta_i(\eta_j-\teta_j)  \widetilde{\Omega}_{jk}(\eta_k-\teta_k)\eta_\ell  \widetilde{\Omega}_{\ell i}\cr
&\equiv Z_{6a} + Z_{6b}. 
\end{align*}
By relabeling $(i,j,k,\ell)$ as $(i,j,\ell,k)$, we can write
\[
Z_{6a} = 2\sum_{\substack{i, j, k, \ell\\ (dist)}}  \eta_i(\eta_j-\teta_j)  \widetilde{\Omega}_{j\ell}\eta_\ell(\eta_k-\teta_k)  \widetilde{\Omega}_{k i} = \sum_{\substack{i, j, k, \ell\\ (dist)}}(\eta_i\eta_\ell \widetilde{\Omega}_{j\ell}\widetilde{\Omega}_{ki})(\eta_j-\teta_j)(\eta_k-\teta_k),
\]
where $|\eta_i\eta_\ell \widetilde{\Omega}_{j\ell}\widetilde{\Omega}_{ki}|\leq C\alpha^2\theta^2_i\theta_j\theta_k\theta^2_\ell$. Also, we write
\[
Z_{6b} = 2\sum_{\substack{i, j, k, \ell\\ (dist)}} \eta_i(\eta_j-\teta_j)  \widetilde{\Omega}_{jk}(\eta_k-\teta_k)\eta_\ell  \widetilde{\Omega}_{\ell i} = 2\sum_{\substack{i, j, k, \ell\\ (dist)}} (\eta_i\eta_\ell \widetilde{\Omega}_{jk}\widetilde{\Omega}_{\ell i})(\eta_j-\teta_j)(\eta_k-\teta_k). 
\]
where $|\eta_i\eta_\ell \widetilde{\Omega}_{jk}\widetilde{\Omega}_{\ell i}|\leq C\alpha^2\theta^2_i\theta_j\theta_k\theta^2_\ell$. At the same time, we recall that 
\[
Z_{5a} = 2\sum_{\substack{i, j, k, \ell\\ (dist)}} \eta_i (\eta_j-\teta_j)\eta_j(\eta_k-\teta_k) \widetilde{\Omega}_{k \ell}  \widetilde{\Omega}_{\ell i} = \sum_{\substack{i, j, k, \ell\\ (dist)}}(\eta_i\eta_j\widetilde{\Omega}_{k\ell}\widetilde{\Omega}_{\ell i})(\eta_j-\teta_j)(\eta_k-\teta_k),
\]
where $|\eta_i\eta_j\widetilde{\Omega}_{k\ell}\widetilde{\Omega}_{\ell i}|\leq C\alpha^2\theta^2_i\theta_j\theta_k\theta^2_\ell$. It is clear that both $Z_{6a}$ and $Z_{6b}$ have a similar structure as $Z_{5a}$. From \eqref{proof-Z5-result1}-\eqref{proof-Z5-result2}, we immediately have
\[
\bigl|\mathbb{E}[Z_6]\bigr| \leq \frac{C\alpha^2\|\theta\|^8}{\|\theta\|_1^2}=o(\alpha^4\|\theta\|^8), \qquad  \mathrm{Var}(Z_6) \leq \frac{C\alpha^4\|\theta\|^8\|\theta\|_3^6}{\|\theta\|_1^2}=o(\|\theta\|^8).
\]
This proves the claims of $Z_6$.

\subsubsection{Proofs of Lemmas \ref{lem:ProxySgnQ(c)-null} and \ref{lem:ProxySgnQ(c)-alt}}
Recall that $\lambda_1, \lambda_2, \ldots, \lambda_K$ are all the nonzero eigenvalues of $\Omega$, arranged in the descending order in magnitude. Write for short $\alpha=|\lambda_2|/|\lambda_1|$.  We shall repeatedly use the following results, which are proved in  \eqref{tOmega(ij)}, \eqref{v-order}, and \eqref{eta-bound}:
\[
v \asymp \|\theta\|_1^2, \qquad 0 < \eta_i < C \theta_i, \qquad |\widetilde{\Omega}_{ij}| \leq C\alpha \theta_i \theta_j. 
\]
Recall that $U_c = 4T_1+F$, under the null hypothesis; $U_c=4T_1+4T_2+F$ under the alternative hypothesis. By definition, 
\begin{align*}
T_1  &= \sum_{i_1, i_2, i_3, i_4 (dist)} \delta_{i_1i_2}\delta_{i_2i_3}\delta_{i_3i_4}W_{i_4i_1},\cr  
T_2 &=\sum_{i_1, i_2, i_3, i_4 (dist)} \delta_{i_1i_2}\delta_{i_2i_3}\delta_{i_3i_4}\widetilde{\Omega}_{i_4i_1}, \cr 
F &= \sum_{i_1, i_2, i_3, i_4 (dist)} \delta_{i_1i_2}\delta_{i_2i_3}\delta_{i_3i_4}\delta_{i_4i_1},
\end{align*}
where $\delta_{ij} = \eta_i(\eta_j - \widetilde{\eta}_j) + \eta_j (\eta_i - \widetilde{\eta}_i)$, for $1 \leq i, j \leq n$, $i \neq j$. 
By symmetry and elementary algebra, we further write 
\begin{equation} \label{T1-split} 
T_1  = 2T_{1a} + 2T_{1b} + 2T_{1c} + 2T_{1d},  
\end{equation} 
where
\begin{align*}
T_{1a} &= \sum_{i_1, i_2, i_3, i_4(dist)}\eta_{i_2}\eta_{i_3}\eta_{i_4}\big[(\eta_{i_1} - \tilde{\eta}_{i_1}) (\eta_{i_2} - \tilde{\eta}_{i_2})  (\eta_{i_3} - \tilde{\eta}_{i_3})    \big]\cdot  W_{i_4i_1},\cr
T_{1b} &=   \sum_{i_1, i_2, i_3, i_4(dist)}\eta_{i_2}\eta_{i_3}^2\big[(\eta_{i_1} - \tilde{\eta}_{i_1}) (\eta_{i_2} - \tilde{\eta}_{i_2}) (\eta_{i_4} - \tilde{\eta}_{i_4})     \big]\cdot  W_{i_4i_1},\cr
T_{1c} &=  \sum_{i_1, i_2, i_3, i_4(dist)}\eta_{i_1}\eta_{i_3}\eta_{i_4}\big[ (\eta_{i_2} - \tilde{\eta}_{i_2})^2  (\eta_{i_3} - \tilde{\eta}_{i_3})    \big]\cdot  W_{i_4i_1},\cr
T_{1d} &=   \sum_{i_1, i_2, i_3, i_4(dist)}\eta_{i_1}\eta_{i_3}^2\big[(\eta_{i_2} - \tilde{\eta}_{i_2})^2 (\eta_{i_4} - \tilde{\eta}_{i_4})     \big]\cdot  W_{i_4i_1}.
\end{align*}
Similarly, we write  
\begin{equation} \label{T2-split} 
T_2  = 2T_{2a} + 2T_{2b} + 2T_{2c} + 2T_{2d},  
\end{equation} 
where 
\begin{align*}
T_{2a} &=  \sum_{i_1, i_2, i_3, i_4(dist)}\eta_{i_2}\eta_{i_3}\eta_{i_4}\big[(\eta_{i_1} - \tilde{\eta}_{i_1}) (\eta_{i_2} - \tilde{\eta}_{i_2})  (\eta_{i_3} - \tilde{\eta}_{i_3})    \big]\cdot  \widetilde{\Omega}_{i_4i_1},\cr
T_{2b} &=  \sum_{i_1, i_2, i_3, i_4(dist)}\eta_{i_2}\eta_{i_3}^2\big[(\eta_{i_1} - \tilde{\eta}_{i_1}) (\eta_{i_2} - \tilde{\eta}_{i_2}) (\eta_{i_4} - \tilde{\eta}_{i_4})     \big]\cdot  \widetilde{\Omega}_{i_4i_1}, \cr
T_{2c} &= \sum_{i_1, i_2, i_3, i_4(dist)}\eta_{i_1}\eta_{i_3}\eta_{i_4}\big[ (\eta_{i_2} - \tilde{\eta}_{i_2})^2  (\eta_{i_3} - \tilde{\eta}_{i_3})    \big]\cdot  \widetilde{\Omega}_{i_4i_1},\cr
T_{2d} &= \sum_{i_1, i_2, i_3, i_4(dist)}\eta_{i_1}\eta_{i_3}^2\big[(\eta_{i_2} - \tilde{\eta}_{i_2})^2 (\eta_{i_4} - \tilde{\eta}_{i_4})     \big]\cdot  \widetilde{\Omega}_{i_4i_1}. 
\end{align*}
Also, similarly, we have  
\begin{equation} \label{F-split}  
F  = 2F_{a} + 12F_{b} + 2F_{c},  
\end{equation} 
where 
\begin{align*}
F_{a} &=  \sum_{i_1,i_2,i_3,i_4(dist)}\eta_{i_1}\eta_{i_2}\eta_{i_3}\eta_{i_4}\big[(\eta_{i_1} - \tilde{\eta}_{i_1}) (\eta_{i_2} - \tilde{\eta}_{i_2})  (\eta_{i_3} - \tilde{\eta}_{i_3})  (\eta_{i_4} - \tilde{\eta}_{i_4})    \big], \cr
F_{b} &=   \sum_{i_1,i_2,i_3,i_4 (dist)}\eta_{i_2}\eta_{i_3}^2\eta_{i_4}\big[(\eta_{i_1} - \tilde{\eta}_{i_1})^2 (\eta_{i_2} - \tilde{\eta}_{i_2}) (\eta_{i_4} - \tilde{\eta}_{i_4})     \big], \cr
F_{c} &=  \sum_{i_1,i_2,i_3,i_4  (dist)}\eta_{i_2}^2\eta^2_{i_4}\big[ (\eta_{i_1} - \tilde{\eta}_{i_1})^2  (\eta_{i_3} - \tilde{\eta}_{i_3})^2    \big]. 
\end{align*}
To show the lemmas, it is sufficient to show the following $11$ items (a)-(k), corresponding to 
$T_{1a}, T_{1b}, T_{1c}, T_{1d}, T_{2a}, T_{2b}, T_{2c}, T_{2d}, F_a, F_b, F_c$, respectively. 
Item (a) claims that both under the null and the alternative,   
\begin{equation} \label{T11varpf} 
|\mathbb{E}[T_{1a}]| \leq C\|\theta\|^6/\|\theta\|_1^2, \qquad \mathrm{Var}(T_{1a}) \leq C \|\theta\|^4\|\theta\|_3^6/\|\theta\|_1^2.   
\end{equation} 
Item (b) claims that both under the null and the alternative, 
\begin{equation} \label{T12varpf} 
|\mathbb{E}[T_{1b}]| \leq C\|\theta\|^6/\|\theta\|_1^2, , \qquad \mathrm{Var}(T_{1b}) \leq C \|\theta\|^6 \|\theta\|_3^3/\|\theta\|_1.   
\end{equation} 
Item (c) claims that both under the null and the alternative, 
\begin{equation} \label{T13varpf} 
\mathbb{E}[T_{1c}] = 0, \qquad \mathrm{Var}(T_{1c}) \leq C \|\theta\|_3^9/\|\theta\|_1,  
\end{equation} 
Item (d) claims that 
\begin{align} \label{T14meanpf} 
& \mathbb{E}[T_{1d}] \asymp  - \|\theta\|^4 \;  \mbox{under the null},\cr
& |\mathbb{E}[T_{1d}]| \leq C \|\theta\|^4  \;  \mbox{under the alternative}, 
\end{align} 
and that both under the null and the alternative, 
\begin{equation} \label{T14varpf} 
\mathrm{Var}(T_{1d}) \leq C \|\theta\|^6 \|\theta\|_3^3/\|\theta\|_1.   
\end{equation} 
Next, for item (e)-(h), we recall that under the null, $T_2 = 0$, and correspondingly $T_{2a} = T_{2b} = T_{2c} = T_{2d} = 0$, so we only need to consider the alternative.  Recall that 
$\alpha = |\lambda_2/\lambda_1|$.  
Item (e) claims that under the alternative,  
\begin{equation} \label{T21varpf} 
\mathbb{E}[T_{2a}] = 0, \qquad \mathrm{Var}(T_{2a}) \leq C \alpha^2\cdot \|\theta\|^{4} \|\theta\|_3^9/\|\theta\|_1^3.  
\end{equation} 
Item (f) claims that under the alternative, 
\begin{equation} \label{T22varpf} 
\mathbb{E}[T_{2b}] = 0,\qquad \mathrm{Var}(T_{2b}) \leq C \alpha^2\cdot \|\theta\|^{12}\|\theta\|_3^3/\|\theta\|_1^5,  
\end{equation} 
Item (g) claims that under the alternative,  
\begin{equation} \label{T23varpf} 
|\mathbb{E}[T_{2c}]| \leq C\alpha\|\theta\|^6/\|\theta\|_1^3, \qquad \mathrm{Var}(T_{2c}) \leq C \alpha^2\cdot \|\theta\|^{8}\|\theta\|_3^3/\|\theta\|_1.  
\end{equation}  
Item (h) claims that both under the null and the alternative, 
\begin{equation} \label{T24varpf} 
|\mathbb{E}[T_{2d}]| \leq C\alpha\|\theta\|^6/\|\theta\|_1^3, 
\qquad \mathrm{Var}(T_{2d}) \leq C \alpha^2\cdot \|\theta\|^{8} \|\theta\|_3^3/\|\theta\|_1.  
\end{equation} 
Finally, for items (i)-(k). Item (i) claims that both under the null and the alternative,  
\begin{equation} \label{F1varpf} 
|\mathbb{E}[F_a]|  \leq C \|\theta\|^8 / \|\theta\|_1^4,\qquad \mathrm{Var}(F_a)  \leq C \|\theta\|_3^{12}/\|\theta\|_1^4.   
\end{equation} 
Item (j) claims that both under the null and the alternative, 
\begin{equation} \label{F2varpf} 
|\mathbb{E}[F_b]|  \leq C\|\theta\|^6/\|\theta\|_1^2,   \qquad \mathrm{Var}(F_b)   \leq C\|\theta\|^4 \|\theta\|_3^6 /\|\theta\|_1^2.   
\end{equation} 
Item (k) claims that 
\begin{align} \label{F3meanpf} 
& \mathbb{E}[F_c] \asymp \|\theta\|^4 \; \mbox{under the null},\cr
&  |\mathbb{E}[F_c]| \leq C \|\theta\|^4 \; 
\mbox{under the alternative},  
\end{align} 
and that under both under the null and the alternative, 
\begin{equation} \label{F3varpf} 
\mathrm{Var}(F_3)  \leq C\|\theta\|^{10}/\|\theta\|_1^2.     
\end{equation} 

We now show Lemmas \ref{lem:ProxySgnQ(a)-null} and  \ref{lem:ProxySgnQ(a)-alt}   follow once (a)-(k) are proved. In detail, first, we note that $\|\theta\|^6/\|\theta\|_1^2 = o(\|\theta\|^4)$.   
Inserting (\ref{T14meanpf}) and the first equation in each of (\ref{T11varpf})-(\ref{T13varpf}) into (\ref{T1-split}) gives that  
\[ 
\mathbb{E}[T_1] \asymp - 2 \|\theta\|^4 \;  \mbox{under the null}, \qquad |\mathbb{E}[T_1]|  \leq C \|\theta\|^4  \; \mbox{under the alternative},  
\] 
and inserting (\ref{T14varpf}) and the second equation in each of (\ref{T11varpf})-(\ref{T13varpf}) into (\ref{T1-split}) gives that both under the null and the alternative,  
\[
\mathrm{Var}(T_1) \leq C[\|\theta\|^4\|\theta\|_3^6/\|\theta\|_1^2 +  
\|\theta\|^6 \|\theta\|_3^3/\|\theta\|_1 +  \|\theta\|_3^9/\|\theta\|_1 + \|\theta\|^6 \|\theta\|_3^3/\|\theta\|_1],  
\] 
where since $\|\theta\|_3^3/\|\theta\|^2  = o(1)$ and $\|\theta\|^2/\|\theta\|_1 = o(1)$, the right hand side 
\[
\leq  C[\|\theta\|^6 \|\theta\|_3^3/\|\theta\|_1^2  + \|\theta\|^6 \|\theta\|_3^3/\|\theta\|_1] \leq C \|\theta\|^6 \|\theta\|_3^3 /\|\theta\|_1. 
\] 
Second, inserting the first equation in each of (\ref{T21varpf})-(\ref{T24varpf})  into  (\ref{T2-split}) gives that under the alternative (recall that $T_2 = 0$ under the null), 
\[
|\mathbb{E}[T_2]|  \leq  C\alpha\|\theta\|^6/\|\theta\|_1^3,   
\] 
and inserting the second equation in each of (\ref{T21varpf})-(\ref{T24varpf}) into (\ref{T2-split}) gives 
\[
\mathrm{Var}(T_2) \leq C\alpha^2 [\|\theta\|^8\|\theta\|_3^3/\|\theta\|_1  + \|\theta\|^{12} \|\theta\|_3^3/\|\theta\|_1^5] \leq C\alpha^2 \|\theta\|^8 \|\theta\|_3^3 / \|\theta\|_1,    
\]
where we have used $\|\theta\|^2 = o(\|\theta\|_1^2)$. 
Third, note that $\|\theta\|^8 / \|\theta\|_1^4 = o(\|\theta\|^4)$ and $\|\theta\|^6 / \|\theta\|_1^2 = o(\|\theta\|^4)$.  Inserting (\ref{F3meanpf}) and the first equation in each of (\ref{F1varpf})-(\ref{F2varpf}) into (\ref{F-split}) gives 
\[
\mathbb{E}[F] \sim 2 \|\theta\|^4 \; \mbox{under the null},  \qquad 
|\mathbb{E}[F]| \leq C \|\theta\|^4 \; \mbox{under the alternative},  
\]
and inserting (\ref{F3varpf}) and the second equation in each of (\ref{F1varpf})-(\ref{F2varpf}) into (\ref{F-split}) gives that both under the null and the alternative,  
\[ 
\mathrm{Var}(F) \leq C[   \|\theta\|_3^{12}/\|\theta\|_1^4 + \|\theta\|^4 \|\theta\|_3^6 /\|\theta\|_1^2 + \|\theta\|^{10}/\|\theta\|_1^2] \leq C \|\theta\|^{10} /\|\theta\|_1^2,   
\] 
where we have used $\|\theta\|_3^3\ll\theta\|^2\ll\|\theta\|_1$ and $\|\theta\|_3^3 / \|\theta\|^2 = o(1)$. 

We now combine the above results for $T_1$, $T_2$ and $F$. First, since that $U_c = 4T_1+F$ under the null, 
it follows that under the null, 
\[
\mathbb{E}[U_c]  \sim -6 \|\theta\|^4, 
\]  
and 
\[
\mathrm{Var}(U_c) \leq C[ \|\theta\|^6 \|\theta\|_3^3 /\|\theta\|_1 +  \|\theta\|^{10} /\|\theta\|_1^2] \leq C \|\theta\|^6 \|\theta\|_3^3 / \|\theta\|_1,   
\] 
where we have used $\|\theta\|^4  \leq \|\theta\|_1 \|\theta\|_3^3$ (a direct use of Cauchy-Schwartz inequality).     
Second, since $U_c=4T_1+4T_2+F$ under the alternative, it follows that under the alternative, 
\[ 
|\mathbb{E}[U_c]| \leq C \|\theta\|^4, 
\]
and 
\[ 
\mathrm{Var}(U_c) \leq C [\|\theta\|^6 \|\theta\|_3^3 /\|\theta\|_1 +\alpha^2 \|\theta\|^8 \|\theta\|_3^3 / \|\theta\|_1 +\|\theta\|^{10} /\|\theta\|_1^2] \leq C \|\theta\|^6\|\theta\|_3^3  (\alpha^2\|\theta\|^2   + 1)/ \|\theta\|_1,  
\] 
where we have used $\|\theta\|^4  \leq \|\theta\|_1 \|\theta\|_3^3$  and basic algebra. 
Combining the above gives all the claims in Lemmas \ref{lem:ProxySgnQ(a)-null} and \ref{lem:ProxySgnQ(a)-alt}. 
 
It remains to show the $11$ items (a)-(k).  We consider them separately. 

Consider Item (a). The goal is to show (\ref{T11varpf}). Recall that   
\[
T_{1a} = \sum_{i_1, i_2, i_3, i_4(dist)}\eta_{i_2}\eta_{i_3}\eta_{i_4}\big[(\eta_{i_1} - \tilde{\eta}_{i_1}) (\eta_{i_2} - \tilde{\eta}_{i_2})  (\eta_{i_3} - \tilde{\eta}_{i_3})    \big]\cdot  W_{i_4i_1},
\]
and that 
\begin{equation} \label{etadiff} 
\widetilde{\eta}  - \eta = v^{-1/2} W 1_n. 
\end{equation} 
Plugging (\ref{etadiff}) into $T_{11}$ gives 
\begin{align*}
T_{1a} &=  - \frac{1}{v^{3/2}} \sum_{i_1, i_2,i_3, i_4 (dist)}\eta_{i_2}\eta_{i_3}\eta_{i_4} \Bigl(\sum_{j_1, j_1\neq i_1}W_{i_1j_1}\Bigr)\Bigl(\sum_{j_2, j_2\neq i_2}W_{i_2j_2}\Bigr)\Bigl(\sum_{j_3, j_3\neq i_3}W_{i_3j_3}\Bigr)W_{i_4 i_1}\cr
&= -\frac{1}{v^{3/2}}\sum_{\substack{i_1, i_2,i_3, i_4  (dist) \\ j_1\neq i_1, j_2\neq i_2, j_3\neq i_3}}\eta_{i_2}\eta_{i_3}\eta_{i_4} W_{i_1j_1}W_{i_2j_2}W_{i_3j_3}W_{i_1i_4}.
\end{align*}
By basic combinatorics and careful observations, we have 
\beq \label{2}
W_{i_1j_1}W_{i_2j_2}W_{i_3j_3}W_{i_1i_4} = 
\begin{cases}
W^2_{i_1i_4}W^2_{i_2i_3},&\mbox{if } j_1 = i_4, (j_2, j_3) = (i_3, i_2),   \\
W^2_{i_1i_4}W_{i_2j_2}W_{i_3j_3}, & \mbox{if } j_1 = i_4, (j_2, j_3) \neq (i_3, i_2),   \\
W^2_{i_2i_3}W_{i_1j_1} W_{i_1i_4}, &\mbox{if }j_1\neq i_4, (j_2, j_3) = (i_3, i_2),    \\
W_{i_1i_2}^2W_{i_3j_3}W_{i_1i_4}, & \mbox{if }(j_1, j_2) = (i_2, i_1),     \\
W^2_{i_1i_3}W_{i_2j_2}W_{i_1i_4}, & \mbox{if }(j_1, j_3) = (i_3, i_1),    \\
W_{i_1j_1}W_{i_2j_2}W_{i_3j_3}W_{i_1i_4}, & \mbox{otherwise}. 
\end{cases}
\eeq
This allows us to further split $T_{11}$ into $6$ different terms: 
\begin{equation} \label{T11-split} 
T_{1a} = X_a + X_{b1} + X_{b2} + X_{b3} + X_{b4}  + X_c, 
\end{equation} 
where 
\begin{align*}
X_{a} & = -\frac{1}{v^{3/2}}\sum_{i_1, i_2, i_3, i_4 (dist)}\eta_{i_2}\eta_{i_3}\eta_{i_4} W^2_{i_1i_4}W^2_{i_2i_3}, \cr
X_{b1} &= -\frac{1}{v^{3/2}}\sum_{i_1, i_2, i_3, i_4 (dist)}  
\sum_{\substack{j_2, j_3  \\ (j_2,j_3) \neq \{i_3,i_2\}}} 
\eta_{i_2}\eta_{i_3}\eta_{i_4} W^2_{i_1i_4}W_{i_2j_2}W_{i_3j_3},\cr
X_{b2}  & = -\frac{1}{v^{3/2}} \sum_{i_1, i_2, i_3, i_4 (dist)} \;   \sum_{j_1 (j_1 \neq i_4)} \eta_{i_2}\eta_{i_3}\eta_{i_4} W^2_{i_2i_3}W_{i_1j_1} W_{i_1i_4},\cr
X_{b3} &= -\frac{1}{v^{3/2}}\sum_{i_1, i_2, i_3, i_4 (dist)} \; 
\sum_{j_3 (j_3 \neq i_3)} \eta_{i_2}\eta_{i_3}\eta_{i_4}W_{i_1i_2}^2W_{i_3j_3}W_{i_1i_4},\cr
X_{b4} &=  -\frac{1}{v^{3/2}}\sum_{i_1, i_2, i_3, i_4(dist)} \; 
\sum_{j_2 (j_2 \neq i_2)}  
\eta_{i_2}\eta_{i_3}\eta_{i_4} W^2_{i_1i_3}W_{i_2j_2}W_{i_1i_4},\cr
X_c &=-\frac{1}{v^{3/2}}\sum_{i_1, i_2, i_3, i_4 (dist)}
\sum_{\substack{j_1, j_2, j_3 \\ j_1\notin\{ i_1, i_4\}, (j_2,j_3)\neq (i_3,i_2) \\  (j_1,j_2)\neq (i_2,i_1), (j_1,j_3)\neq (i_3,i_1)}} 
\eta_{i_2}\eta_{i_3}\eta_{i_4} W_{i_1j_1}W_{i_2j_2}W_{i_3j_3}W_{i_1 i_4}. 
\end{align*}

We now show (\ref{T11varpf}).  Consider the first claim of (\ref{T11varpf}). It is seen that out of the $6$ terms on the right hand side of 
(\ref{T11-split}), the mean of all terms are $0$, except for the first term. Note that for any $1 \leq i, j \leq n$, $i \neq j$, $\mathbb{E}[W_{ij}^2] = \Omega_{ij} (1 -\Omega_{ij})$, where 
$\Omega_{ij}$ are upper bounded by $o(1)$ uniformly for all such $i, j$. It follows 
\begin{align*}
\mathbb{E}[X_a] & = - v^{-3/2}\sum_{i_1, i_2, i_3, i_4 (dist)} \eta_{i_2} \eta_{i_3} \eta_{i_4} \mathbb{E}[W_{i_1i_4}^2] \mathbb{E}[W_{i_2i_3}^2] \\
& = -(1+o(1))\cdot v^{-3/2} \sum_{i_1, i_2, i_3, i_4(dist)} \eta_{i_2} \eta_{i_3} \eta_{i_4} \Omega_{i_1 i_4} \Omega_{i_2i_3}. 
\end{align*}
Since for any $1 \leq i, j\leq n$, $i \neq j$, $0 < \eta_i \leq C \theta_i$, $\Omega_{ij} \leq C \theta_i \theta_j$ and $v \asymp \|\theta\|_1^2$, 
\[
|\mathbb{E}[X_a]|  \leq C  (\|\theta\|_1)^{-3} \sum_{i_1, i_2, i_3, i_4 (dist)} \theta_{i_1} \theta_{i_2}^2 \theta_{i_3}^2 \theta_{i_4}^2 \leq C\|\theta\|^6 / \|\theta\|_1^2.  
\] 
Inserting these into (\ref{T11-split}) gives 
\begin{equation} \label{mean-T11} 
|\mathbb{E}[T_{1a}] |  \leq C \|\theta\|^6 / \|\theta\|_1^2, 
\end{equation} 
and the first claim of (\ref{T11varpf}) follows. 

Consider the second claim of  (\ref{T11varpf}) next. 
By (\ref{T11-split}) and Cauchy-Schwartz inequality,  
\begin{align}  
\mathrm{Var}(T_{1a}) & \leq C \mathrm{Var}(X_{a}) + \mathrm{Var}(X_{b1}) + \mathrm{Var}(X_{b2})  +  \mathrm{Var}(X_{b3})  + \mathrm{Var}(X_{b4}) + \mathrm{Var}(X_{c}))  \nonumber \\ 
& \leq  C(\mathrm{Var}(X_{a})+ \mathbb{E}[X_{b1}^2] + \mathbb{E}[X_{b2}^2] + \mathbb{E}[X_{b3}^2] + \mathbb{E}[X_{b4}^2] + 
\mathbb{E}[X_{c}^2]). \label{T11-var-split} 
\end{align} 
We now consider $ \mathrm{Var}(X_{a})$, $\mathbb{E}[X_{b1}^2] + \mathbb{E}[X_{b2}^2] + \mathbb{E}[X_{b3}^2] + \mathbb{E}[X_{b4}^2]$, and $\mathbb{E}[X_{c}^2]$, separately. 

Consider $\mathrm{Var}(X_{a})$.  Write $\mathrm{Var}(X_{a})$ as
\begin{align}  
& v^{-3}  \sum_{\substack{i_1, \cdots, i_4 (dist) \\ 
i_1', \cdots, i_4' (dist)} } \eta_{i_2} \eta_{i_3} \eta_{i_4} \eta_{i_2'} \eta_{i_3'} \eta_{i_4'}  \nonumber  \\
&\qquad  \qquad \qquad \mathbb{E}\big[(W_{i_1 i_4}^2 W_{i_2i_3}^2- \mathbb{E}[W_{i_1i_4}^2W_{i_2i_3}^2]) (W_{i_1' i_4'}^2 W_{i_2'i_3'}^2- \mathbb{E}[W_{i_1'i_4'}^2W_{i_2'i_3'}^2])\big].  \label{varXa-add}   
\end{align} 
In the sum, a term is nonzero only when one of the following cases happens.  
\begin{itemize}
\item (A). $\{W_{i_1 i_4}, W_{i_2 i_3}, W_{i_1'i_4'}, W_{i_2'i_3'}\}$ has $2$ distinct random variables. 
\item (B). $\{W_{i_1 i_4},  W_{i_2, i_3}, W_{i_1'i_4'}, W_{i_2'i_3'}\}$ has $3$ distinct random variables. This has $4$ sub-cases: (B1)  $W_{i_1 i_4} = W_{i_1' i_4'}$, (B2) $W_{i_1i_4} = W_{i_2'i_3'}$, (B3) $W_{i_2i_3} = W_{i_1' i_4'}$, and  (B4) $W_{i_2i_3} = W_{i_2' i_3'}$.  
\end{itemize} 
For Case (A),  the two sets $\{i_1, i_2, i_3, i_4\}$ and $\{i_1', i_2', i_3', 
i_4'\}$ are identical. By basic statistics and independence between $W_{i_1i_4}$ and $W_{i_2 i_3}$,  
\begin{align} 
& \mathbb{E}[(W_{i_1 i_4}^2 W_{i_2i_3}^2- \mathbb{E}[W_{i_1i_4}^2W_{i_2i_3}^2]) (W_{i_1' i_4'}^2 W_{i_2'i_3'}^2- \mathbb{E}[W_{i_1'i_4'}^2W_{i_2'i_3'}^2])]  \nonumber   \\ 
=&  \mathbb{E}[(W_{i_1i_4}^2 W_{i_2i_3}^2 - \mathbb{E}[W_{i_1i_4}^2  W_{i_2i_3}^2])^2] \nonumber  \\ \nonumber
= &  \mathbb{E}[W_{i_1i_4}^4] \mathbb{E}[W_{i_2i_3}^4]  - (\mathbb{E} [W_{i_1i_4}^2])^2 (\mathbb{E}[W_{i_2i_3}^2])^2 \\ 
\leq & \mathbb{E}[W_{i_1i_4}^4] \mathbb{E}[W_{i_2i_3}^4],        \label{varXa-addeq1} 
\end{align} 
where by basic statistics and that $\Omega_{ij} \leq C \theta_i \theta_j$ for all $1 \leq i, j \leq n$, $i \leq j$,  the  right hand side 
\[ 
\leq C  \Omega_{i_1i_4} \Omega_{i_2 i_3} \leq C \theta_{i_1} \theta_{i_2} \theta_{i_3} \theta_{i_4}. 
\] 
Combining these with  (\ref{varXa-add}) and noting that $v \sim \|\theta\|_1^2$ and that $0 < \eta_i \leq C \theta_i$ for all $1 \leq i \leq n$, the contribution of this case to $\mathrm{Var}(X_a)$ is no more than 
\begin{equation} \label{varXa-add1} 
C (\|\theta\|_1)^{-6}  \sum_{i_1, \cdots, i_4 (dist)} \sum_{a}  \theta_{i_1}^{a_1+1} \theta_{i_2}^{a_2+2} \theta_{i_3}^{a_3 + 2} \theta_{i_4}^{a_4 + 2}, 
\end{equation} 
where $a = (a_1, a_2, a_3, a_4)$ and each $a_i$ is either $0$ and $1$, satisfying $a_1 + a_2 + a_3 + a_4 = 3$. 
Note that the right hand side of (\ref{varXa-add1}) is no greater than 
\[
C (\|\theta\|_1)^{-6} \max\{\|\theta\|_1 \|\theta\|_3^9, \|\theta\|^4 \|\theta\|_3^6\} \leq C \|\theta\|_3^9 / \|\theta\|_1^5, 
\] 
where we have used $\|\theta\|^4 \leq \|\theta\|_1 \|\theta\|_3^3$. 

Next, consider (B1).  By independence between $W_{i_1i_4}$, $W_{i_2i_3}$, and $W_{i_2'i_3'}$ and basic algebra,  
\begin{align} 
& \mathbb{E}[(W_{i_1 i_4}^2 W_{i_2i_3}^2- \mathbb{E}[W_{i_1i_4}^2W_{i_2i_3}^2]) (W_{i_1' i_4'}^2 W_{i_2'i_3'}^2- \mathbb{E}[W_{i_1'i_4'}^2W_{i_2'i_3'}^2])] \nonumber \\
 = & \mathbb{E}[(W_{i_1i_4}^2 W_{i_2i_3}^2 - \mathbb{E}[W_{i_1i_4}^2 W_{i_2i_3}^2]) (W_{i_1i_4}^2 W_{i_2'i_3'}^2 - \mathbb{E}[W_{i_1i_4}^2 W_{i_2'i_3'}^2])]  \nonumber \\
 = & \mathbb{E}[W_{i_1i_4}^4] \mathbb{E}[W_{i_2i_3}^2] \mathbb{E}[W_{i_2'i_3'}^2] - (\mathbb{E}[W_{i_1i_4}^2])^2 \mathbb{E}[W_{i_2i_3}^2] \mathbb{E}[W_{i_2'i_3'}^2]  \nonumber  \\
 = & \mathrm{Var}(W_{i_1i_4}^2) \mathbb{E}[W_{i_2i_3}^2] \mathbb{E}[W_{i_2'i_3'}^2],    \label{varXa-addeq2} 
\end{align} 
where by similar arguments, the last term 
\[
\leq C \Omega_{i_1i_4} \Omega_{i_2i_3} \Omega_{i_2'i_3'} \leq C \theta_{i_1} \theta_{i_2}\theta_{i_3} \theta_{i_4} \theta_{i_2'} \theta_{i_3'}. 
\] 
Combining this with (\ref{varXa-add}) and using similar arguments,   the contribution of this case to $\mathrm{Var}(X_a)$ 
\begin{equation} \label{varXa-add2} 
\leq C  (\|\theta\|_1)^{-6}  \sum_{\substack{i_1, i_2, i_3, i_4 (dist) \\ i_2', i_3'  (dist)}} 
C \theta_{i_1}^{b_1 + 1} \theta_{i_2}^2  \theta_{i_3}^2 \theta_{i_4}^{b_2 + 2}  \theta_{i_2'}^2  \theta_{i_3'}^2,  
\end{equation} 
where similarly $b_1, b_2$ are either $0$ or $1$ and $b_1 + b_2 = 1$.  By similar argument, the right hand side 
\[
\leq C \|\theta\|_1 \|\theta\|^8 \|\theta\|_3^3 / \|\theta\|_1^6 = C \|\theta\|^8 \|\theta\|_3^3 /\|\theta\|_1^5. 
\] 

The discussion for (B2), (B3), and (B4) are similar so is omitted, and their contribution to $\mathrm{Var}(X_a)$ are respectively 
\begin{equation} \label{varXa-add3}
\leq C  \|\theta\|^8 \|\theta\|_3^3 /\|\theta\|_1^5,  
\end{equation} 
\begin{equation} \label{varXa-add4} 
\leq C \|\theta\|^8 \|\theta\|_3^3 /\|\theta\|_1^5,  
\end{equation} 
and 
\begin{equation} \label{varXa-add5} 
\leq C  \|\theta\|^4  \|\theta\|_3^6 /\|\theta\|_1^4. 
\end{equation} 
Finally, inserting (\ref{varXa-add1}), (\ref{varXa-add2}), (\ref{varXa-add3}), (\ref{varXa-add4}), and (\ref{varXa-add5}) into (\ref{varXa-add}) gives   
\begin{equation} \label{T11-Vara} 
\mathrm{Var}(X_a) \leq  C[\|\theta\|_3^9 / \|\theta\|_1^5 +  \|\theta\|^{8} \|\theta\|_3^3 / \|\theta\|_1^5 + \|\theta\|^4 \|\theta\|_3^6 / \|\theta\|_1^4] \leq  C \|\theta\|^4 \|\theta\|_3^6/\|\theta\|_1^4,  
\end{equation} 
where we have used $\|\theta\|_3^3 \ll \|\theta\|^2$ and $\|\theta\|^4 \leq \|\theta\|_1 \|\theta\|_3^3$.

Consider $\mathbb{E}[X_{b1}^2] + \mathbb{E}[X_{b21}^2] + \mathbb{E}[X_{b3}^2] + \mathbb{E}[X_{b4}^2]$.  
We claim that both under the null and the alternative,  
\begin{align} 
\mathbb{E}[X_{b1}^2] & \leq C \|\theta\|^4 \|\theta\|_3^6/\|\theta\|_1^2    \label{T11-Varb-1},  \\
\mathbb{E}[X_{b2}^2] & \leq C \|\theta\|^8 \|\theta\|_3^3/\|\theta\|_1^3  \label{T11-Varb-2},  \\ 
\mathbb{E}[X_{b3}^2] & \leq  C \|\theta\|^6 \|\theta\|_3^6/\|\theta\|_1^4 \label{T11-Varb-3},  \\
\mathbb{E}[X_{b4}^2] & \leq  C \|\theta\|^6 \|\theta\|_3^6/\|\theta\|_1^4 \label{T11-Varb-4},     
\end{align} 
where the last two terms are seen to be negligible compared to the other two. Therefore, 
\begin{equation} \label{T11-Varb} 
\mathbb{E}[X_{b1}^2] + \mathbb{E}[X_{b2}^2] + \mathbb{E}[X_{b3}^2] + \mathbb{E}[X_{b4}^2] \leq C[\|\theta\|^4 \|\theta\|_3^6/\|\theta\|_1^2 + \|\theta\|^8 \|\theta\|_3^3 / \|\theta\|_1^3], 
\end{equation}  
where since $\|\theta\|^4 \leq \|\theta\|_1 \|\theta\|_3^3$ (Cauchy-Schwartz inequality) the right hand side 
\[  
\leq C[\|\theta\|^4 \|\theta\|_3^6/\|\theta\|_1^2.   
\] 
 
We now prove (\ref{T11-Varb-1})-(\ref{T11-Varb-4}).  
Since the study for $\mathbb{E}[X_{b1}^2], \mathbb{E}[X_{b2}^2], \mathbb{E}[X_{b3}^2]$ and $\mathbb{E}[X_{b4}^2]$ are similar, we only present the proof for $\mathbb{E}[X_{b1}^2]$. Write  $\mathbb{E}[X_{b1}^2]$ as 
\begin{align*} 
  v^{-3} \sum_{\substack{i_1,i_2,i_3,i_4 (dist)  \\ i_1',i_2',i_3',i_4'(dist)}}   
\sum_{\substack{j_2, j_3  \\    (j_2,j_3) \neq (i_3,i_2)}} 
\sum_{\substack{j_2', j_3'  \\    (j_2',j_3') \neq (i_3',i_2')}} 
\eta_{i_2}\eta_{i_3}\eta_{i_4}\eta_{i_2'}\eta_{i_3'}\eta_{i_4'}   W^2_{i_1i_4} W_{i_2j_2}W_{i_3j_3} W^2_{i_1'i_4'} W_{i_2'j_2'}W_{i_3'j_3'}. 
\end{align*} 
Consider the term
\[
W^2_{i_1i_4} W_{i_2j_2}W_{i_3j_3} W^2_{i_1'i_4'} W_{i_2'j_2'}W_{i_3'j_3'}.
\] 
In order for the mean to be nonzero, we have two cases 
\begin{itemize} 
\item Case A. The two sets of random variables $\{ W_{i_1i_4}, W_{i_2j_2}, W_{i_3j_3}\}$ and 
$\{ W_{i_1'i_4'}, W_{i_2'j_2'}, W_{i_3'j_3'}\}$ are identical. 
\item Case B. The two sets $\{W_{i_2j_2}, W_{i_3j_3}\}$ and $\{W_{i_2'j_2'}, W_{i_3'j_3'}\}$ are identical. 
\end{itemize} 
Consider Case A. In this case, $\{i_2', i_3', i_4'\}$ are three distinct indices in $\{ i_1, i_2, i_3, i_4, j_2, j_3\}$, and for some integers satisfying  
$0 \leq a_1, a_2, \ldots, a_6 \leq 1$, $a_1 + a_2 + \ldots +  a_6 = 3$, 
\[
\eta_{i_2}\eta_{i_3}\eta_{i_4}\eta_{i_2'}\eta_{i_3'}\eta_{i_4'} = \eta_{i_1}^{a_1}  \eta_{i_2}^{1 + a_2}  \eta_{i_3}^{1 + a_3}  \eta_{i_4}^{1 + a_4}  \eta_{j_2}^{a_5}  \eta_{j_3}^{a_6}   
\] 
and for some integers satisfying  $0 \leq b_1, b_2, b_3\leq 1$,  and $b_1 + b_2 + b_3 = 1$,  
\[
W^2_{i_1i_4} W_{i_2j_2}W_{i_3j_3} W^2_{i_1'i_4'} W_{i_2'j_2'}W_{i_3'j_3'}
= W_{i_1i_4}^{b_1 + 3}  W_{i_2j_2}^{b_2 + 2} W_{i_3j_3}^{b_3+2}.
\] 
Similarly, by $v \sim \|\theta\|_1^2$, $0 < \eta_i \leq C \theta_i$, and uniformly for all $b_1, b_2, b_3$ above, 
\[
0 < \mathbb{E}[W_{i_1i_4}^{b_1 + 3}  W_{i_2j_2}^{b_2 + 2} W_{i_3j_3}^{b_3+2}] \leq C \Omega_{i_1i_4} \Omega_{i_2j_2} \Omega_{i_3j_3} \leq C \theta_{i_1}\theta_{i_2} \theta_{i_3}\theta_{i_4} \theta_{j_2}\theta_{j_3}. 
\] 
Therefore under both the null and the alternative, the contribution of Case $A$ to the variance is  
\begin{equation} \label{T11-Varb1-A} 
\leq C (\|\theta\|_1)^{-6}   \sum_{i_1, i_2, i_3, i_4 (dist)}   
\sum_{\substack{j_2, j_3  \\  j_2 \neq i_2, j_3 \neq i_3, (j_2,j_3) \neq (i_3,i_2)}}  [\sum_a  \theta_{i_1}^{a_1+1}  \theta_{i_2}^{a_2 + 2}  \theta_{i_3}^{a_3 + 2}  \theta_{i_4}^{a_4 + 2} \theta_{j_2}^{a_5 + 1} \theta_{j_3}^{a_6 + 1}],  
\end{equation} 
where $a = (a_1, a_2, \ldots, a_6)$ and $a_i$ satisfies the above constraints. Note that the right hand size 
\[
\leq C (\|\theta\|_1)^{-6} \cdot \max\{\|\theta\|_1^3  \|\theta\|_3^{9},  \|\theta\|_1^2 \|\theta\|^4 \|\theta\|_3^6, \|\theta\|_1 \|\theta\|^8 \|\theta\|_3^3, \|\theta\|^{12}\} \leq C \|\theta\|_3^9 / \|\theta\|_1^3.  
\] 
Here in the last inequality we have used $\|\theta\|^2 \leq \sqrt{\|\theta\|_1 \|\theta\|_3^3}$.  

Consider Case B. In this case, $\{i_2, i_3, j_2, j_3\} =\{i_2', i_3', j_2', j_3'\} $,  
and for some integers $0 \leq c_1, c_2, c_3, c_4 \leq 1$, $c_1 + c_2 + c_3 + c_4 = 2$, 
\[
\eta_{i_2} \eta_{i_3} \eta_{i_4} \eta_{i_2'} \eta_{i_3'} \eta_{i_4'}  = \eta_{i_2}^{c_1 + 1}  \eta_{i_3}^{c_2 + 1}   \eta_{i_4}    \eta_{j_2}^{c_3}  \eta_{j_3}^{c_4}   \eta_{i_4'}, 
\] 
and 
\[
W^2_{i_1i_4} W_{i_2j_2}W_{i_3j_3} W^2_{i_1'i_4'} W_{i_2'j_2'}W_{i_3'j_3'} =
W^2_{i_1i_4} W^2_{i_2j_2}W^2_{i_3j_3} W^2_{i_1'i_4'},  
\]
where the four $W$ terms on the right are independent of each other. 
Similarly, by $v \sim \|\theta\|_1^2$, $0 < \eta_i \leq C \theta_i$, 
\[
0 < \mathbb{E}[W^2_{i_1i_4} W^2_{i_2j_2}W^2_{i_3j_3} W^2_{i_1'i_4'}] \leq C \Omega_{i_1i_4} \Omega_{i_2j_2} \Omega_{i_3j_3} \Omega_{i_1'i_4'}  \leq C \theta_{i_1} \theta_{i_2} \theta_{i_3} \theta_{i_4} \theta_{j_2} \theta_{j_3} \theta_{i_1'} \theta_{i_4'},  
\] 
we have that under both the null and the alternative, the contribution of Case $B$ to the variance 
\[
\leq C (\|\theta\|_1)^{-6} \sum_{\substack{i_1, i_2, i_3, i_4 (dist) \\ i_1', i_4' (dist)}}  
\sum_{\substack{j_2, j_3  \\    (j_2,j_3) \neq (i_3,i_2)}}  \theta_{i_1} \theta_{i_2}^{c_1 + 2}  \theta_{i_3}^{c_2 + 2}  \theta_{i_4}^2 \theta_{j_2}^{c_3 + 1}  \theta_{j_3}^{c_4 + 1} \theta_{i_1'} \theta_{i_4'}^2, 
\] 
where the right hand size 
\begin{equation} \label{T11-Varb1-B} 
\leq C (\|\theta\|_1)^{-6} \cdot \|\theta\|_1^2 \|\theta\|^4 \cdot  \max\{\|\theta\|_1^2 \|\theta\|_3^6, \|\theta\|_1  \|\theta\|^4 \|\theta\|_3^3,   \|\theta\|^8\} \leq C \|\theta\|^4 \|\theta\|_3^6 /\|\theta\|_1^2.   
\end{equation} 
Here we have again used $\|\theta\|^2 \leq \sqrt{\|\theta\|_1 \|\theta\|_3^3}$. 

Finally, combining (\ref{T11-Varb1-A}) and (\ref{T11-Varb1-B})  gives 
\[
\mathbb{E}[X_{b1}^2] \leq C (\|\theta\|_3^9/\|\theta\|_1^3 + \|\theta\|^4 \|\theta\|_3^6/\|\theta\|_1^2) \leq C \|\theta\|^4 \|\theta\|_3^6/\|\theta\|_1^2, 
\] 
which proves (\ref{T11-Varb-1}).

Consider $\mathbb{E}[X_c^2]$. Consider the terms in the sum, 
\[
\eta_{i_2} \eta_{i_3}\eta_{i_4}   W_{i_1j_1} W_{i_2j_2} W_{i_3j_3} W_{i_1i_4}, \qquad \mbox{and} \qquad \eta_{i_2'} \eta_{i_3'}\eta_{i_4'}   W_{i_1'j_1'} W_{i_2'j_2'} W_{i_3'j_3'} W_{i_1'i_4'}. 
\] 
Each term has a mean $0$, and two terms are uncorrelated with each other if only if the two sets of random variables $\{ W_{i_1j_1} , W_{i_2j_2} , W_{i_3j_3}, W_{i_1i_4}\}$ and $\{ W_{i_1'j_1'} , W_{i_2'j_2'} , W_{i_3'j_3'}, W_{i_1'i_4'}\}$ are identical (however, it is possible that $W_{i_1j_1}$ does not equal to $W_{i_1j_1'}$ but equals to $W_{i_2'j_2'}$, say).  Additionally, the indices $i_2', i_3', i_4' \in \{i_1, i_2, i_3, i_4, j_1, j_2, j_3\}$, and it follows
\begin{align*}
\mathbb{E}[X_c^2] \leq C  v^{-3}    \sum_{i_1, i_2, i_3, i_4 (dist)  } 
& \sum_{\substack{j_1, j_2, j_3 \\ j_1\notin\{ i_1,i_4\}, (j_1,j_3)\neq (i_3,i_1) \\ (j_2,j_3)\neq (i_3,i_2), (j_2,j_1)\neq (i_2,i_1)}} \\
& [\sum_{a} \eta_{i_1}^{a_1} \eta_{i_2}^{a_2 +1}  \eta_{i_3}^{a_3+1} \eta_{i_4}^{a_{4}+1}\eta^{a_5}_{j_1}\eta^{a_6}_{j_2}\eta^{a_7}_{j_3}] \cdot   \mathbb{E}[W_{i_1j_1}^2 W_{i_2j_2}^2 W_{i_3j_3}^2 W_{i_1j_1}^2],  
\end{align*}
where $a = (a_1, a_2, \cdots, a_7)$ and the power $0 \leq a_1,  a_2, \cdots, a_7 \leq 1$, 
and $a_1 + a_2 + \cdots + a_7 = 3$. 
Note that $W_{i_1j_1},  W_{i_2j_2},  W_{i_3j_3}$ and $W_{i_1i_4}$ are independent and $\mathbb{E}(W_{ij}^2) \leq \Omega_{ij} \leq C \theta_i \theta_j$,  $1 \leq i, j \leq n$, $i \neq j$,  
\[
\mathbb{E}[W_{i_1j_1}^2 W_{i_2j_2}^2 W_{i_3j_3}^2 W_{i_1 i_4}^2] \leq 
\Omega_{i_1j_1} \Omega_{i_2j_2} \Omega_{i_3j_3} \Omega_{i_1i_4} \leq C \theta_{i_1}^2 \theta_{i_2} \theta_{i_3} \theta_{i_4} \theta_{j_1} \theta_{j_2} \theta_{j_3}. 
\] 
Also, recall that both under the null and the alternative,  $v \asymp \|\theta\|_1^2$ and  $0 < \eta_i \leq C \theta_i$, $1 \leq i \leq n$. 
Combining these gives  
\begin{align*}
\mathbb{E}[X_c^2] & \leq C (\|\theta\|_1)^{-6}  \sum_{i_1, i_2, i_3, i_4 (dist)} 
\sum_{\substack{j_1, j_2, j_3 \\ j_1\notin\{ i_1,i_4\}, (j_1,j_3)\neq (i_3,i_1) \\ (j_2,j_3)\neq (i_3,i_2), (j_2,j_1)\neq (i_2,i_1) }} \\ 
& \qquad\qquad  \qquad \qquad \qquad [\sum_a   \eta_{i_1}^{a_1+2} \eta_{i_2}^{a_2 +2}  \eta_{i_3}^{a_3+2} \eta_{i_4}^{a_{4}+2}\eta^{a_5+1}_{j_1}\eta^{a_6+1}_{j_2}\eta^{a_7+1}_{j_3}],  
\end{align*}
where the last term 
\[ 
\leq C \sum_{a} \|\theta\|_{a_1 + 2}^{a_1 + 2}  \cdot \|\theta\|_{a_2 + 2}^{a_2 + 2}  \cdot \|\theta\|_{a_3 + 2}^{a_3 + 2} \cdot  \|\theta\|_{a_4 + 2}^{a_4 + 2}\|\theta\|_{a_5 +1}^{a_5 + 1}\|\theta\|_{a_6 + 1}^{a_6 + 1}\|\theta\|_{a_7 + 1}^{a_7 + 1}/\|\theta\|_1^6.   
\] 
Since $a_1, a_2, \cdots, a_7$ have to take values from $\{0, 1\}$ and their sum is $3$,  the above term 
\[
\leq C \|\theta\|^2 \|\theta\|_3^9 / \|\theta\|_1^3 = o(\|\theta\|_3^3),  
\] 
where we have used  $\|\theta\|_3^3 \ll \|\theta\|_2^2 \ll \|\theta\|_1$. Combining these gives 
\begin{equation} \label{T11-Varc} 
\mathbb{E}[X_c^2] \leq C \|\theta\|^2 \|\theta\|_3^9 / \|\theta\|_1^3.  
\end{equation} 

Finally, inserting (\ref{T11-Vara}), (\ref{T11-Varb}), and (\ref{T11-Varc}) into (\ref{T11-split}) gives that both under the null and the alternative, 
\[
\mathrm{Var}(T_{11}) \leq C[\|\theta\|^{8} / \|\theta\|_1^4 +  \|\theta\|^4 \|\theta\|_3^6/\|\theta\|_1^2 + \|\theta\|^2 \|\theta\|_3^9 / \|\theta\|_1^3] \leq C \|\theta\|^4 \|\theta\|_3^6 /\|\theta\|_1^2, 
\]  
where we have used $\|\theta\|^4 \leq \|\theta\|_1 \|\theta\|_3^3$ and $\|\theta\|_3^3/\|\theta\|_1 = o(1)$. 
This gives (\ref{T11varpf}) and completes the proof for Item (a). 


Consider Item (b). The goal is to show (\ref{T12varpf}). Recall that \[
T_{1b} =   \sum_{i_1, i_2, i_3, i_4(dist)}\eta_{i_2}\eta_{i_3}^2\big[(\eta_{i_1} - \tilde{\eta}_{i_1}) (\eta_{i_2} - \tilde{\eta}_{i_2}) (\eta_{i_4} - \tilde{\eta}_{i_4})     \big]\cdot  W_{i_4i_1},
\]
and that 
\[
\widetilde{\eta}  - \eta = v^{-1/2}  W 1_n. 
\]
Plugging this into $T_{1b}$ gives 
\begin{align*}
T_{1b} &=  - v^{-3/2} \sum_{i_1,i_2,i_3,i_4 (dist)} \eta_{i_2}\eta^2_{i_3} \Bigl(\sum_{j_1\neq i_1}W_{i_1j_1}\Bigr)\Bigl(\sum_{j_2\neq i_2}W_{i_2j_2}\Bigr)\Bigl(\sum_{j_4\neq i_4}W_{i_4j_4}\Bigr)W_{i_1i_4}\cr
&= -\frac{1}{v^{3/2}}\sum_{\substack{i_1, i_2, i_3, i_4 (dist)\\ j_1\neq i_1, j_2\neq i_2, j_4\neq i_4}}\eta_{i_2}\eta^2_{i_3} W_{i_1j_1}W_{i_2j_2}W_{i_4j_4}W_{i_1i_4}.
\end{align*}
By basic combinatorics and careful observations, we have 
\beq \label{2}
W_{i_1j_1}W_{i_2j_2}W_{i_4j_4}W_{i_1i_4} = 
\begin{cases}
W^3_{i_1i_4}W_{i_2j_2},&\mbox{if } j_1=i_4, j_4 = i_1,    \\
W^2_{i_1i_2}W^2_{i_1i_4},&\mbox{if } j_1=i_2, j_2 = i_1, j_4 = i_1,   \\
W^2_{i_1i_4}W^2_{i_2i_4},&\mbox{if } j_1=i_4, j_2 = i_4, j_4 = i_2,  \\
W^2_{i_1i_2}W_{i_4j_4}W_{i_1i_4},&\mbox{if } j_1=i_2, j_2 = i_1,    \\
W_{i_1i_4}^2W_{i_1j_1}W_{i_2j_2},&\mbox{if } j_4=i_1,    \\
W^2_{i_1i_4}W_{i_2j_2}W_{i_4j_4},&\mbox{if } j_1=i_4, \{i_2, j_2\}\neq\{i_4, j_4\},     \\
W_{i_2i_4}^2W_{i_1j_1}W_{i_1i_4},&\mbox{if } j_2=i_4, j_4=i_2,     \\
 W_{i_1j_1}W_{i_2j_2}W_{i_4j_4}W_{i_1i_4}, & \mbox{otherwise}. 
\end{cases}
\eeq
This allows us to further split $T_{1b}$ into $8$ different terms: 
\begin{equation} \label{T12-split} 
T_{1b} = Y_{a1}  + Y_{a2} + Y_{a3} + Y_{b1} + Y_{b2} + Y_{b3}  + Y_{b4}  + Y_{c}, 
\end{equation} 
where 
\begin{align*}
Y_{a1} & = -\frac{1}{v^{3/2}}\sum_{i_1, i_2, i_3,i_4(dist)}\sum_{j_2 (j_2 \neq i_2)}\eta_{i_2}\eta_{i_3}^2W^3_{i_1i_4}W_{i_2j_2}, \cr
Y_{a2} & = -\frac{1}{v^{3/2}}\sum_{i_1, i_2, i_3,i_4 (dist)}\eta_{i_2}\eta_{i_3}^2W^2_{i_1i_2}W^2_{i_1i_4}, \cr
Y_{a3} & = -\frac{1}{v^{3/2}}\sum_{i_1, i_2, i_3,i_4 (dist)}\eta_{i_2}\eta_{i_3}^2W^2_{i_1i_4}W^2_{i_2i_4}, \cr
Y_{b1} &= -\frac{1}{v^{3/2}}\sum_{i_1, i_2, i_3,i_4 (dist)}  
\sum_{\substack{j_4( j_4\neq i_4)}} 
\eta_{i_2}\eta_{i_3}^2 W^2_{i_1i_2}W_{i_4j_4}W_{i_1i_4},\cr
Y_{b2} &=  -\frac{1}{v^{3/2}} \sum_{i_1, i_2, i_3,i_4 (dist)} \; 
\sum_{\substack{j_1( j_1\neq i_1), j_2( j_2\neq i_2)  \\ \{i_1, j_1\} \neq \{i_2, j_2\}}} 
\eta_{i_2}\eta_{i_3}^2  W_{i_1i_4}^2 W_{i_1j_1} W_{i_2j_2},\cr
Y_{b3} &= -\frac{1}{v^{3/2}}\sum_{i_1, i_2, i_3,i_4 (dist)}  
\sum_{\substack{j_2( j_2\neq i_2), j_4( j_4\neq i_4)\\ \{i_2, j_2\}\neq\{i_4, j_4\}}} 
\eta_{i_2}\eta_{i_3}^2 W^2_{i_1i_4}W_{i_2j_2}W_{i_4j_4},\cr
Y_{b4} &= -\frac{1}{v^{3/2}}\sum_{i_1, i_2, i_3,i_4 (dist)}  
\sum_{\substack{j_1( j_1\neq i_1)}} 
\eta_{i_2}\eta_{i_3}^2  W_{i_2i_4}^2 W_{i_1j_1}W_{i_1i_4},\cr
Y_{c} &=-\frac{1}{v^{3/2}} \sum_{i_1, i_2, i_3,i_4 (dist)}
\sum_{\substack{ j_1, j_2, j_4 \\ j_1\notin\{i_2, i_4\}, j_2\notin\{i_1, i_4\}, j_4\notin\{i_1, i_2\}   }} 
\eta_{i_2}\eta_{i_3}^2  W_{i_1j_1}W_{i_2j_2}W_{i_4j_4}W_{i_1 i_4}.  
\end{align*}
We now show the two claims in (\ref{T12varpf}) separately. 

Consider the first claim of (\ref{T12varpf}). It is seen that out of the $8$ terms on the right hand side of 
(\ref{T12-split}), the mean of all terms are $0$, except that of  the $Y_{a2}$ and $Y_{a3}$. Note that for any $1 \leq i, j \leq n$, $i \neq j$, $\mathbb{E}[W_{ij}^2] = \Omega_{ij} (1 -\Omega_{ij})$, where 
$\Omega_{ij}$ are upper bounded by $o(1)$ uniformly for all such $i, j$. It follows 
\begin{align*}
\mathbb{E}[Y_{a2}] & =  -\frac{1}{v^{3/2}}\sum_{i_1, i_2, i_3,i_4 (dist)}\eta_{i_2}\eta_{i_3}^2\mathbb{E}[W^2_{i_1i_2}]\mathbb{E}[W^2_{i_1i_4}] \\
& = -(1+o(1))\cdot v^{-3/2} \sum_{i_1, i_2, i_3, i_4 (dist)} \eta_{i_2} \eta_{i_3}^2 \Omega_{i_1i_2} \Omega_{i_1i_4}. 
\end{align*}
Since for any $1 \leq i, j\leq n$, $i \neq j$, $0 < \eta_i \leq C \theta_i$, $\Omega_{ij} \leq C \theta_i \theta_j$ and $v \asymp \|\theta\|_1^2$, 
\[
|\mathbb{E}[Y_{a2}]|  \leq C  (\|\theta\|_1)^{-3} \sum_{i_1, i_2, i_3,i_4 (dist)}\theta_{i_1}^2 \theta_{i_2}^2 \theta_{i_3}^2 \theta_{i_4}  \leq C\|\theta\|^6 / \|\theta\|_1^2.  
\] 
Therefore, 
\begin{equation} \label{mean-T12} 
|\mathbb{E}[Y_{a2}] |  \leq C \|\theta\|^6 / \|\theta\|_1^2. 
\end{equation} 
By symmetry, we similarly find 
\begin{equation} \label{mean-T12-2} 
|\mathbb{E}[Y_{a3}] |  \leq C \|\theta\|^6 / \|\theta\|_1^2. 
\end{equation} 
Combining \eqref{mean-T12}  and \eqref{mean-T12-2} gives \[
\mathbb{E}[|T_{1b}|] \leq C \|\theta\|^6 / \|\theta\|_1^2. 
\]
This completes the proof of the first claim of (\ref{T12varpf}). 

We now show the second claim of (\ref{T12varpf}) .  
By Cauchy-Schwartz inequality,  
\begin{align}  
\mathrm{Var}(T_{1b}) & \leq C (\mathrm{Var}(Y_{a1}) + \mathrm{Var}(Y_{a2}) +  \mathrm{Var}(Y_{a3}) + \sum_{s=1}^4\mathrm{Var}(Y_{bs})  + \mathrm{Var}(Y_{c}))  \nonumber \\ 
& \leq  C(\mathrm{Var}(Y_{a1}) + \mathrm{Var}(Y_{a2}) +  \mathrm{Var}(Y_{a3}) +  \sum_{s=1}^4\mathbb{E}[Y^2_{bs}] +  
\mathbb{E}[Y_{c}^2]). \label{T12-var-split} 
\end{align} 
We now show $\mathrm{Var}(Y_{a1})$, $\mathrm{Var}(Y_{a2})$, $\mathrm{Var}(Y_{a3})$,  $\sum_{s=1}^4\mathbb{E}[Y_{bs}^2]$, and $\mathbb{E}[Y_{c}^2]$, separately. 

Consider $\mathrm{Var}(Y_{a1})$. Recalling $\mathbb{E}[Y_{a1}] = 0$, we write $\mathrm{Var}(Y_{a1})$ as 
\beq\label{varYa-add}
v^{-3}  \sum_{\substack{i_1, i_2, i_3, i_4 (dist) \\ 
i_1', i_2', i_3', i_4'  (dist) }}  \sum_{j_2 (j_2 \neq i_2)} \sum_{j_2' (j_2' \neq i_2')} \eta_{i_2} \eta_{i_3}^2 \eta_{i_2'} \eta_{i_3'}^2 \mathbb{E}\big[W_{i_1i_4}^3 W_{i_2j_2} W_{i_1'i_4'}^3 W_{i_2'j_2'}\big].  
\eeq
In the sum, a term is nonzero only when one of the following cases happens.  
\begin{itemize}
\item (A). $\{W_{i_1 i_4}, W_{i_2 j_2}, W_{i_1'i_4'}, W_{i_2'j_2'}\}$ has $2$ distinct random variables. 
\item (B). $\{W_{i_1 i_4},  W_{i_2j_2}, W_{i_1'i_4'}, W_{i_2'j_2'}\}$ has $3$ distinct random variables. While it may seem we have $4$ possibilities in this case, but the only one that has a nonzero mean is when $W_{i_2j_2} = W_{i_2' j_2'}$. 
\end{itemize} 
For Case (A),  the two sets $\{i_1, i_2, i_4, j_2\}$ and $\{i_1', i_2', 
i_4', j_2' \}$ are identical, 
and so for two integers $0 \leq b_1, b_2 \leq 1$ and $b_1 + b_2 = 1$,  
\[
 W_{i_1i_4}^3 W_{i_2j_2} W_{i_1'i_4'}^3 W_{i_2'j_2'} = W_{i_1i_4}^{4+2b_1}   W_{i_2j_2}^{2+2b_2},         
\] 
and so 
\[
\mathbb{E}[W_{i_1i_4}^3 W_{i_2j_2} W_{i_1'i_4'}^3 W_{i_2'j_2'}]  = \mathbb{E}[W_{i_1i_4}^{4+2b_1}   W_{i_2j_2}^{2+2b_2}] = \mathbb{E}[W_{i_1i_4}^{4+2b_1}] \mathbb{E}[W_{i_2j_2}^{2+2b_2}], 
\] 
Note that for any integer $2 \leq b \leq 6$,   
\[
0 < \mathbb{E}[W_{ij}^b]  \leq C \Omega_{ij}, 
\]
where note that $\Omega_{ij} \leq C \theta_i \theta_j$ for all $1 \leq i, j \leq n$, $i \leq j$. Recall that $v \sim \|\theta\|_1^2$,  and that $0 < \eta_i \leq C \theta_i$ for all $1 \leq i \leq n$. Combining these that, the contribution of Case (A) to $\mathrm{Var}(Y_{a1})$ is no more than 
\begin{equation} \label{varYa-add1} 
C (\|\theta\|_1)^{-6}  \sum_{i_1, \cdots, i_4 (dist)} \sum_{i_3', j_2}\sum_{a}  \theta_{i_1}^{a_1+1} \theta_{i_2}^{a_2+2} \theta_{i_3}^{ 2} \theta_{i_4}^{a_3 + 1}\theta_{i_3'}^2\theta_{j_2}^{a_4+1}, 
\end{equation} 
where $a = (a_1, a_2, a_3, a_4)$ and each $a_i$ is either $0$ and $1$, satisfying $a_1 + a_2 + a_3 + a_4 = 1$. 
Note that the right hand side of (\ref{varYa-add1}) is no greater than 
\[
C (\|\theta\|_1)^{-6} \max\{ \|\theta\|_1^3 \|\theta\|^4 \|\theta\|_3^3, \|\theta\|_1^2\|\theta\|^{8} \} \leq C\|\theta\|^4 \|\theta\|_3^3/\|\theta\|_1^3, 
\] 
where we have used $\|\theta\|^4 \leq \|\theta\|_1 \|\theta\|_3^3$. 

Next, consider Case (B).  In this case, $\{i_2, j_2\} = \{i_2', j_2'\}$ and 
\[
W_{i_1 i_4}^3 W_{i_2j_2}  W_{i_1' i_4'}^3 W_{i_2'j_2'} = W_{i_1i_4}^3 W_{i_2j_2}^2 W_{i_1' i_4'}^3, 
\] 
and by similar argument,  
\begin{align} 
0 <  \mathbb{E}[W_{i_1 i_4}^3 W_{i_2j_2}^2 W_{i_1' i_4'}^3] \leq C \Omega_{i_1 i_4} \Omega_{i_2 j_2} \Omega_{i_1'i_4'}. 
\end{align} 
Recall that $\Omega_{ij} \leq C \theta_i \theta_j$ for all $1 \leq i, j \leq n$, $i \leq j$, that $v \sim \|\theta\|_1^2$,  and that $0 < \eta_i \leq C \theta_i$ for all $1 \leq i \leq n$. 
Combining this with (\ref{varYa-add}),    the contribution of this case to $\mathrm{Var}(Y_{a1})$ 
\begin{equation} \label{varYa-add2} 
\leq C  (\|\theta\|_1)^{-6}  \sum_{\substack{i_1, i_2, i_3, i_4 (dist) \\ i_1', i_3', i_4'  (dist)}} 
\sum_{j_2} C  \theta_{i_1} \theta_{i_2}^{2+b_1}\theta_{i_3}^2 \theta_{i_4} \theta_{i_1'}\theta_{i_3'}^2  \theta_{i_4'}\theta_{j_2}^{1+b_2},  
\end{equation} 
where similarly $b_1, b_2$ are either $0$ or $1$ and $b_1 + b_2 = 1$.  By similar argument, the right hand side 
\[
\leq C\|\theta\|_1^{-6}\cdot[ \|\theta\|_1^5 \|\theta\|^4\|\theta\|_3^3 + \|\theta\|_1^4 \|\theta\|^8]   \leq C\|\theta\|^4\|\theta\|_3^3 /\|\theta\|_1,
\] 
where we've used Cauchy-Schwartz inequality that $\|\theta\|^4\leq \|\theta\|_1\|\theta\|_3^3$.

Now, inserting (\ref{varYa-add1}) and (\ref{varYa-add2}) into (\ref{varYa-add}) gives   
\begin{equation} \label{T12-Vara} 
\mathrm{Var}(Y_{a1}) \leq  C[ \|\theta\|^4 \|\theta\|_3^3/\|\theta\|_1^3 +  \|\theta\|^4  \|\theta\|_3^3 /\|\theta\|_1] \leq  C \|\theta\|^4 \|\theta\|_3^3/\|\theta\|_1,  
\end{equation} 
where we have used $\|\theta\|_1 \goto \infty$ and $\|\theta\|^4 \leq \|\theta\|_1 \|\theta\|_3^3$.   
This shows 
\begin{equation} \label{varYa1-add} 
\mathrm{Var}(Y_{a1}) \leq C\|\theta\|^4\|\theta\|_3^3/\|\theta\|_1. 
\end{equation} 

Next, we consider $\mathrm{Var}(Y_{a2})$ and $\mathrm{Var}(Y_{a3})$. The proofs are similar 
to that of $\mathrm{Var}(X_a)$ of Item (a), so we skip the detail, but claim that 
\begin{equation} \label{varYa2-add} 
\mathrm{Var}(Y_{a2}) \leq C\|\theta\|^4\|\theta\|_3^6/\|\theta\|_1^4,  
\end{equation} 
and 
\begin{equation} \label{varYa3-add} 
\mathrm{Var}(Y_{a3}) \leq C\|\theta\|^4\|\theta\|_3^6/\|\theta\|_1^4.  
\end{equation} 
Combining (\ref{varYa1-add}), (\ref{varYa2-add}), and (\ref{varYa3-add}) gives 
\begin{equation} \label{T12-Vara} 
\mathrm{Var}(Y_{a1}) +  \mathrm{Var}(Y_{a2}) + \mathrm{Var}(Y_{a3}) 
\leq C[\|\theta\|^4\|\theta\|_3^3/\|\theta\|_1 + \|\theta\|^4\|\theta\|_3^6/\|\theta\|_1^4 ] \leq C\|\theta\|^4\|\theta\|_3^3/\|\theta\|_1, 
\end{equation} 
where we have used the universal inequality that $\|\theta\|_3^3\leq \|\theta\|_1^3$. 

Next, consider $\sum_{s=1}^4\mathbb{E}[Y_{bs}^2]$.  For each $1 \leq s \leq 4$, the study of $\mathbb{E}[Y_{bs}^2]$ is   similar to that of $\mathbb{E}[X_{b1}^2]$ in Item (a), so we skip the details. We have that both under the null and the alternative, 
\begin{align} 
\mathbb{E}[Y_{b1}^2] & \leq C \|\theta\|^{12}  /\|\theta\|^4_1   \label{T12-Varb-1},  \\
\mathbb{E}[Y_{b2}^2] & \leq  C \|\theta\|^6 \|\theta\|_3^3 /\|\theta\|_1 \label{T12-Varb-2},  \\ 
\mathbb{E}[Y_{b3}^2] & \leq C \|\theta\|^6 \|\theta\|_3^3 /\|\theta\|_1 \label{T12-Varb-3},  \\
\mathbb{E}[Y_{b4}^2] & \leq  C \|\theta\|^{12}  /\|\theta\|^4_1 \label{T12-Varb-4}.  
\end{align} 
Therefore,  
\begin{equation} \label{T12-Varb} 
\sum_{s = 1}^4 \mathbb{E}[Y_{bs}^2] \leq C[\|\theta\|^6 \|\theta\|_3^3/\|\theta\|_1 +   \|\theta\|^{12}  /\|\theta\|^4_1]  
\leq C \|\theta\|^6 \|\theta\|_3^3 / \|\theta\|_1. 
\end{equation} 

Third, we consider $\mathbb{E}[Y_c^2]$. The proof is very similar to that of $\mathbb{E}[X_c^2]$ and we have that both under the null and the alternative, 
\begin{equation} \label{T12-Varc} 
\mathbb{E}[Y_c^2] \leq C \|\theta\|^8 \|\theta\|_3^3 / \|\theta\|_1^3. 
\end{equation} 

Finally, combining (\ref{T12-Vara}), (\ref{T12-Varb}), and (\ref{T12-Varc}) with (\ref{T12-var-split})  gives 
\begin{equation} \label{T12-Varsum} 
\mathrm{Var}(T_{1b}) \leq C[ \|\theta\|^4 \|\theta\|_3^3 / \|\theta\|_1 +  \|\theta\|^6 \|\theta\|_3^3 / \|\theta\|_1 + \|\theta\|^8 \|\theta\|_3^3 / \|\theta\|_1^3] \leq C\|\theta\|^6 \|\theta\|_3^3 / \|\theta\|_1 , 
\end{equation} 
where we have used $\|\theta\|\goto\infty$ and $\|\theta\|^2\ll\|\theta\|_1$. 
This completes the proof of (\ref{T12varpf}).  

Consider Item (c). The goal is to show (\ref{T13varpf}). Recall that \[
T_{1c} =  \sum_{i_1, i_2, i_3, i_4(dist)}\eta_{i_1}\eta_{i_3}\eta_{i_4}\big[ (\eta_{i_2} - \tilde{\eta}_{i_2})^2  (\eta_{i_3} - \tilde{\eta}_{i_3})    \big]\cdot  W_{i_4i_1},
\]
and that 
\[
\widetilde{\eta}  - \eta = v^{-1/2} W 1_n. 
\]
Plugging this into $T_{1c}$ gives 
\begin{align*}
T_{1c} &=  -\frac{1}{v^{3/2}} \sum_{i_1,i_2,i_3,i_4 (dist)} \eta_{i_1}\eta_{i_3}\eta_{i_4}  \Bigl(\sum_{j_2\neq i_2}W_{i_2j_2}\Bigr)\Bigl(\sum_{\ell_2\neq i_2}W_{i_2\ell_2}\Bigr)  
\Bigl(\sum_{j_3\neq i_3}W_{i_3j_3}\Bigr)W_{i_1i_4}\cr
&= -\frac{1}{v^{3/2}}\sum_{\substack{i_1, i_2, i_3, i_4 (dist)\\ j_2\neq i_2, \ell_2\neq i_2, j_3\neq i_3}}\eta_{i_1}\eta_{i_3}\eta_{i_4} W_{i_2j_2}W_{i_2\ell_2}W_{i_3j_3}W_{i_1i_4}.
\end{align*}
By basic combinatorics and careful observations, we have 
\beq \label{2}
W_{i_2j_2}W_{i_2\ell_2}W_{i_3j_3}W_{i_1i_4} = 
\begin{cases}
W^3_{i_2i_3}W_{i_1i_4},&\mbox{if } j_2=\ell_2 = i_3, j_3 = i_2, \\
W_{i_2j_2}^2W_{i_3j_3}W_{i_1i_4},&\mbox{if } j_2 = \ell_2, (j_3, j_2) \neq (i_2, i_3) , \\
W_{i_2i_3}^2W_{i_2\ell_2}W_{i_1i_4},&\mbox{if } j_2=i_3, j_3 = i_2, \ell_2 \neq i_3,  \\
W_{i_2i_3}^2W_{i_2j_2}W_{i_1i_4},&\mbox{if } \ell_2=i_3, j_3 = i_2, j_2 \neq i_3, \\
W_{i_2j_2}W_{i_2\ell_2}W_{i_3j_3}W_{i_1i_4}, & \mbox{otherwise}. 
\end{cases}
\eeq
This allows us to further split $T_{1c}$ into $5$ different terms: 
\begin{equation} \label{T13-split} 
T_{1c} = Z_{a}  + Z_{b1} + Z_{b2} + Z_{b3} + Z_{c}, 
\end{equation} 
where 
\begin{align*}
Z_{a} & = -\frac{1}{v^{3/2}}\sum_{i_1, i_2, i_3,i_4 (dist)}\eta_{i_1}\eta_{i_3}\eta_{i_4}W^3_{i_2i_3}W_{i_1i_4}, \cr
Z_{b1} &= -\frac{1}{v^{3/2}}\sum_{i_1, i_2, i_3,i_4 (dist)}  
\sum_{\substack{ j_2, (j_3, j_2) \neq (i_2, i_3)}} 
\eta_{i_1}\eta_{i_3}\eta_{i_4}W^2_{i_2j_2}W_{i_3j_3}W_{i_1i_4} ,\cr
Z_{b2} &= -\frac{1}{v^{3/2}}\sum_{i_1, i_2, i_3,i_4 (dist)}  
\sum_{\substack{ j_2=i_3, j_3 = i_2\\ \ell_2 \neq i_3}} 
\eta_{i_1}\eta_{i_3}\eta_{i_4} W_{i_2i_3}^2W_{i_2\ell_2}W_{i_1i_4},\cr
Z_{b3} &= -\frac{1}{v^{3/2}}\sum_{i_1, i_2, i_3,i_4 (dist)}  
\sum_{\substack{ \ell_2=i_3, j_3 = i_2\\ j_2 \neq i_3}} 
\eta_{i_1}\eta_{i_3}\eta_{i_4} W_{i_2i_3}^2W_{i_2j_2}W_{i_1i_4},\cr
Z_{c} &=-\frac{1}{v^{3/2}} \sum_{i_1, i_2, i_3,i_4 (dist)}
\sum_{\substack{ j_2,\ell_2, j_3\\ j_2\neq \ell_2, j_2, \ell_2\neq i_3, j_3\neq i_2}} 
\eta_{i_1}\eta_{i_3}\eta_{i_4} W_{i_2j_2}W_{i_2\ell_2}W_{i_3j_3}W_{i_1i_4} .  
\end{align*}

We now show  the two claims in (\ref{T13varpf}) separately. The proof of the first claim is trivial, so we only show the second claim of (\ref{T13varpf}). 

Consider the second claim of (\ref{T13varpf}). 
By Cauchy-Schwartz inequality,  
\begin{align}  
\mathrm{Var}(T_{1c}) & \leq C (\mathrm{Var}(Z_{a}) + \mathrm{Var}(Z_{b1}) +  \mathrm{Var}(Z_{b2}) + \mathrm{Var}(Z_{b3})  + \mathrm{Var}(Z_{c}))  \nonumber \\ 
& \leq  C(\mathbb{E}[Z^2_{a}] +  \sum_{s=1}^3\mathbb{E}[Z^2_{bs}] +  
\mathbb{E}[Z_{c}^2]). \label{T13-var-split} 
\end{align}  
Note that 
\begin{itemize} 
\item The proof of $\mathrm{Var}(Z_a)$ is similar to that of $\mathrm{Var}(Y_a)$ in Item (b). 
\item The proof of $\sum_{s=1}^3\mathbb{E}[Z_{bs}^2]$ is similar to that of $\sum_{s = 1}^4 \mathbb{E}[X_{bs}^2]$ in Item (a). 
\item The proof of $\mathbb{E}[Z_c^2]$ is similar to that of $\mathbb{E}[X_c^2]$ in Item (a). 
\end{itemize} 
For these reasons, we omit the proof details and only state the claims. We have that under both the null and the 
alternative, 
\begin{equation} \label{T13-Vara} 
\mathrm{Var}(Z_a) \leq C \|\theta\|^4 \|\theta\|_3^6/\|\theta\|_1^4, 
\end{equation} 
\begin{equation} \label{T13-Varb} 
\sum_{s = 1}^3 \mathbb{E}[Z^2_{bs}] \leq C\|\theta\|_3^9 /\|\theta\|_1, 
\end{equation} 
and 
\begin{equation} \label{T13-Varc} 
\mathbb{E}[Z_c^2]  \leq C  \|\theta\|^2\|\theta\|_3^9 / \|\theta\|_1^3.  
\end{equation} 
Inserting (\ref{T13-Vara}), (\ref{T13-Varb}), and (\ref{T13-Varc}) into (\ref{T13-var-split}) gives 
\[
\mathrm{Var}(T_{1c}) \leq  C[\|\theta\|^4 \|\theta\|_3^6/\|\theta\|_1^4 + \|\theta\|_3^9 /\|\theta\|_1 + \|\theta\|^2\|\theta\|_3^9 / \|\theta\|_1^3] \leq C\|\theta\|_3^9 /\|\theta\|_1, 
\] 
where we have used $\|\theta\|_3^3\ll \|\theta\|^2 \ll\|\theta\|_1$, $\|\theta\|^4\leq \|\theta\|_1\|\theta\|_3^3$ and $\|\theta\|_1\goto \infty$. This proves (\ref{T13varpf}).


Consider Item (d). The goal is to show (\ref{T14meanpf}) and (\ref{T14varpf}). Recall that \[
T_{1d} = -\frac{1}{v^{3/2}} \sum_{i_1, i_2, i_3, i_4(dist)}\eta_{i_1}\eta_{i_3}^2\big[(\eta_{i_2} - \tilde{\eta}_{i_2})^2 (\eta_{i_4} - \tilde{\eta}_{i_4})     \big]\cdot  W_{i_4i_1}. 
\]
and that
\[
\widetilde{\eta}  - \eta = v^{-1/2} W 1_n. 
\]
Plugging this into $T_{1d}$ gives 
\begin{align*}
T_{1d} &= -\frac{1}{v^{3/2}} \sum_{i_1,i_2,i_3,i_4 (dist)} \eta_{i_1}\eta^2_{i_3} \Bigl(\sum_{j_2\neq i_2}W_{i_2j_2}\Bigr)\Bigl(\sum_{\ell_2\neq i_2}W_{i_2\ell_2}\Bigr)\Bigl(\sum_{j_4\neq i_4}W_{i_4j_4}\Bigr)W_{i_1i_4}\cr
&= -\frac{1}{v^{3/2}}\sum_{\substack{i_1, i_2, i_3, i_4 (dist) \\ j_2\neq i_2, \ell_2\neq i_2, j_4\neq i_4}}\eta_{i_1}\eta^2_{i_3}  W_{i_2j_2}W_{i_2\ell_2}W_{i_4j_4}W_{i_1i_4}.
\end{align*}
By basic combinatorics and careful observations, we have 
\beq \label{2}
W_{i_2j_2}W_{i_2\ell_2}W_{i_4j_4}W_{i_1i_4} = 
\begin{cases}
W^3_{i_2i_4}W_{i_1i_4},&\mbox{if } j_2=\ell_2=i_4, j_4 = i_2,  \\
W^2_{i_2j_2}W^2_{i_1i_4},&\mbox{if } j_2=\ell_2, j_4 = i_1,  \\
W^2_{i_2j_2}W_{i_4j_4}W_{i_1i_4},&\mbox{if } j_2=\ell_2, j_4\neq i_1, (j_2, j_4)\neq(i_4, i_2),  \\
W_{i_2j_2}W_{i_2i_4}^2W_{i_1i_4},&\mbox{if } \ell_2=i_4, j_4 = i_2, j_2\neq i_4,   \\
W_{i_2\ell_2}W_{i_2i_4}^2W_{i_1i_4},&\mbox{if } j_2=i_4, j_4 = i_2, \ell_2\neq i_4, \\
W_{i_2j_2}W_{i_2\ell_2}W_{i_1i_4}^2,&\mbox{if } j_4=i_1, j_2\neq \ell_2, \\
W_{i_2j_2}W_{i_2\ell_2}W_{i_4j_4}W_{i_1i_4} , & \mbox{otherwise}. 
\end{cases}
\eeq
This allows us to further split $T_{14}$ into $7$ different terms: 
\begin{equation} \label{T12-split} 
T_{1d} = U_{a1}  + U_{a2}  + U_{b1} + U_{b2} + U_{b3}  + U_{b4}  + U_{c}, 
\end{equation} 
where 
\begin{align*}
U_{a1} & = -\frac{1}{v^{3/2}}\sum_{i_1, i_2, i_3,i_4 (dist)}\eta_{i_1}\eta_{i_3}^2W^3_{i_2i_4}W_{i_1i_4}, \cr
U_{a2} & = -\frac{1}{v^{3/2}}\sum_{i_1, i_2, i_3,i_4 (dist)} \sum_{j_2}     \eta_{i_1}\eta_{i_3}^2W^2_{i_2j_2}W^2_{i_1i_4}, \cr
U_{b1} &= -\frac{1}{v^{3/2}}\sum_{i_1, i_2, i_3,i_4 (dist)}  
\sum_{\substack{j_2( j_2\neq i_2), j_4( j_4\neq i_4)\\ j_4\neq i_1, (j_2, j_4)\neq(i_4, i_2) }} 
\eta_{i_1}\eta_{i_3}^2 W^2_{i_2j_2}W_{i_4j_4}W_{i_1i_4},\cr
U_{b2} &= -\frac{1}{v^{3/2}}\sum_{i_1, i_2, i_3,i_4 (dist)}  
\sum_{\substack{j_2( j_2\neq i_4)}} 
\eta_{i_1}\eta_{i_3}^2 W_{i_2j_2}W_{i_2i_4}^2W_{i_1i_4},\cr
U_{b3} &= -\frac{1}{v^{3/2}}\sum_{i_1, i_2, i_3,i_4 (dist)}  
\sum_{\substack{\ell_2( \ell_2\neq i_4)}} 
\eta_{i_1}\eta_{i_3}^2 W_{i_2\ell_2}W_{i_2i_4}^2W_{i_1i_4},\cr
U_{b4} &=  -\frac{1}{v^{3/2}} \sum_{i_1, i_2, i_3,i_4 (dist)} \; 
\sum_{\substack{j_2\neq \ell_2}} 
\eta_{i_1}\eta_{i_3}^2 W_{i_2j_2}W_{i_2\ell_2}W_{i_1i_4}^2,\cr
U_{c} &=-\frac{1}{v^{3/2}} \sum_{i_1, i_2, i_3,i_4 (dist)}
\sum_{\substack{ j_2, \ell_2, j_4, W \text{dist}}} 
\eta_{i_1}\eta_{i_3}^2 W_{i_2j_2}W_{i_2\ell_2}W_{i_4j_4}W_{i_1i_4}.  
\end{align*}
We now show (\ref{T14meanpf}) and (\ref{T14varpf}) separately.

Consider (\ref{T14meanpf}). It is seen that out of the $7$ terms on the right hand side of 
(\ref{T13-split}), all terms are mean $0$, except for the second  term $U_{a2}$. Note that for any $1 \leq i, j \leq n$, $i \neq j$, $\mathbb{E}[W_{ij}^2] = \Omega_{ij} (1 -\Omega_{ij})$, where 
$\Omega_{ij}$ are upper bounded by $o(1)$ uniformly for all such $i, j$. It follows 
\begin{align*}
\mathbb{E}[U_{a2}] &=  -\frac{1}{v^{3/2}}\sum_{i_1, i_2, i_3,i_4 (dist)}\sum_{j_2}\eta_{i_1}\eta_{i_3}^2\mathbb{E}[W^2_{i_2j_2}]\mathbb{E}[W^2_{i_1i_4}] \\
& = -(1+o(1))\cdot v^{-3/2} \sum_{\substack{ i_1, i_2, i_3, i_4 (dist)\\ j_2}} \eta_{i_1} \eta^2_{i_3}\Omega_{i_2j_2} \Omega_{i_1i_4}. 
\end{align*}
Under null,  for any $1 \leq i, j\leq n$, $i \neq j$, $\eta_i = (1+o(1)) \theta_i$, $\Omega_{ij} = (1+o(1)) \theta_i \theta_j$ and $v \asymp \|\theta\|_1^2$, 
\[
\mathbb{E}[U_{a2}]  =  (\|\theta\|_1)^{-3} \sum_{i_1, i_2, i_3,i_4(dist)}\sum_{j_2}\theta^2_{i_1} \theta_{i_2} \theta_{i_3}^2 \theta_{i_4}\theta_{j_2} = - (1+o(1)) \|\theta\|^4,
\] 
and under alternative, a similar arguments yields
\begin{equation} \label{mean-T14} 
|\mathbb{E}[U_{a1}] |  \leq C \|\theta\|^4. 
\end{equation} 
This proves (\ref{T14meanpf}).

We now consider (\ref{T14varpf}). 
By Cauchy-Schwartz inequality,  
\begin{align}  
\mathrm{Var}(T_{1d}) & \leq C (\mathrm{Var}(U_{a1}) + \mathrm{Var}(U_{a2})  + \sum_{s=1}^4\mathrm{Var}(U_{bs})  + \mathrm{Var}(U_{c}))  \nonumber \\ 
& \leq  C(\mathrm{Var}(U_{a1}) + \mathrm{Var}(U_{a2}) +  \sum_{s=1}^4\mathbb{E}[U^2_{bs}] +  
\mathbb{E}[U_{c}^2]). \label{T14-var-split} 
\end{align} 
Note that 
\begin{itemize} 
\item The proof of $U_{a1}$ is similar to that of $Y_{a1}$ in Item (b). 
\item The proof of $U_{a2}$ is similar to that of $X_{a1}$ in Item (a). 
\item The proof of $U_{bs}$, $1 \leq s \leq 4$, is similar to that of $X_{b1}$ in Item (a).
\item The proof of $U_c$ is similar to that of $X_c$ in Item (a).  
\end{itemize} 
For these reasons, we omit the proof details, and claim that 
\begin{equation} \label{T14-Vara-1}
\mathrm{Var}(U_{a1}) \leq C\|\theta\|^4\|\theta\|_3^6/\|\theta\|_1^4,
\end{equation} 
\begin{equation} \label{T14-Vara-2}
\mathrm{Var}(U_{a2}) \leq  C\|\theta\|^4\|\theta\|_3^3/\|\theta\|_1,
\end{equation} 
\begin{equation} \label{T14-Varb} 
\sum_{s = 1}^4 \mathbb{E}[U_{bs}^2] \leq C\|\theta\|^6\|\theta\|_3^3/\|\theta\|_1,  
\end{equation} 
and 
\begin{equation} \label{T14-Varc} 
\mathrm{Var}(U_c)  \leq C\|\theta\|^8\|\theta\|_3^3/\|\theta\|_1^3,  
\end{equation} 
Inserting (\ref{T14-Vara-1}), (\ref{T14-Vara-2}), (\ref{T14-Varb}), and (\ref{T14-Varc}) into 
(\ref{T14-var-split}) gives 
\begin{align} 
\mathrm{Var}(T_{1d}) & \leq C[\|\theta\|^4\|\theta\|_3^6/\|\theta\|_1^4  + \|\theta\|^4\|\theta\|_3^3/\|\theta\|_1 + \|\theta\|^6\|\theta\|_3^3/\|\theta\|_1 + \|\theta\|^8\|\theta\|_3^3/\|\theta\|_1^3] \\
&  \leq C\|\theta\|^6\|\theta\|_3^3/\|\theta\|_1,   
\end{align}  
where we have used $\|\theta\|\goto\infty$ and $\|\theta\|_3^3\leq\|\theta\|_1^3$. This proves (\ref{T14varpf}). 


We now consider Item (e) and Item (f). Since the proof is similar, we only prove Item (e). The goal is to show  (\ref{T21varpf}). Recall that 
\begin{equation} \label{def:T21-add} 
T_{2a} =  \sum_{i_1, i_2, i_3, i_4(dist)}\eta_{i_2}\eta_{i_3}\eta_{i_4}\big[(\eta_{i_1} - \tilde{\eta}_{i_1}) (\eta_{i_2} - \tilde{\eta}_{i_2})  (\eta_{i_3} - \tilde{\eta}_{i_3})    \big]\cdot  \widetilde{\Omega}_{i_4i_1},
\end{equation} 
and 
\begin{equation} \label{etadiff-new} 
\widetilde{\eta}  - \eta = v^{-1/2} W 1_n. 
\end{equation} 
Plugging (\ref{etadiff-new}) into (\ref{def:T21-add})  gives 
\begin{align*}
T_{2a} &= -\frac{1}{v^{3/2}} \sum_{i_1, i_2,i_3, i_4 (dist)}\eta_{i_2}\eta_{i_3}\eta_{i_4} \Bigl( \sum_{j_1\neq i_1}W_{i_1j_1}\Bigr)\Bigl(\sum_{j_2\neq i_2}W_{i_2j_2}\Bigr)\Bigl(\sum_{j_3\neq i_3}W_{i_3j_3}\Bigr)\widetilde{\Omega}_{i_4 i_1}\cr
&= -\frac{1}{v^{3/2}}\sum_{\substack{i_1, i_2,i_3, i_4 (dist) \\ j_1\neq i_1, j_2\neq i_2, j_3\neq i_3}}\eta_{i_2}\eta_{i_3}\eta_{i_4} W_{i_1j_1}W_{i_2j_2}W_{i_3j_3}\widetilde{\Omega}_{i_1i_4}.
\end{align*}
By basic combinatorics and careful observations, we have 
\beq \label{2}
W_{i_1j_1}W_{i_2j_2}W_{i_3j_3} = 
\begin{cases}
W^2_{i_1i_2}W_{i_3j_3},&\mbox{if } j_1 = i_2,  j_2 = i_1, \\
W^2_{i_1i_3}W_{i_2j_2}, & \mbox{if }  j_1 = i_3,  j_3 = i_1, \\
W^2_{i_2i_3}W_{i_1j_1}, & \mbox{if }  j_2 = i_3,  j_3 = i_2, \\
W_{i_1j_1}W_{i_2j_2}W_{i_3j_3}, & \mbox{otherwise}. 
\end{cases}
\eeq
This allows us to further split $T_{2a}$ into $4$ different terms: 
\begin{equation} \label{T21-split} 
T_{2a} = X_{a1} + X_{a2} + X_{a3} + X_{b}, 
\end{equation} 
where 
\begin{align*}
X_{a1} & = -\frac{1}{v^{3/2}}\sum_{i_1, i_2, i_3, i_4 (dist)}\sum_{j_3\neq i_3}\eta_{i_2}\eta_{i_3}\eta_{i_4} W^2_{i_1i_2}W_{i_3j_3}\widetilde{\Omega}_{i_1i_4}, \cr
X_{a2} & = -\frac{1}{v^{3/2}}\sum_{i_1, i_2, i_3, i_4 (dist)}\sum_{j_2\neq i_2}\eta_{i_2}\eta_{i_3}\eta_{i_4} W^2_{i_1i_3}W_{i_2j_2}\widetilde{\Omega}_{i_1i_4}, \cr
X_{a3} & = -\frac{1}{v^{3/2}}\sum_{i_1, i_2, i_3, i_4 (dist)}\sum_{j_1\neq i_1}\eta_{i_2}\eta_{i_3}\eta_{i_4} W^2_{i_2i_3}W_{i_1j_1}\widetilde{\Omega}_{i_1i_4}, \cr
X_b &=-\frac{1}{v^{3/2}}\sum_{i_1, i_2, i_3, i_4 (dist)}
\sum_{\substack{j_1, j_2, j_3 \\ j_k\neq i_\ell, k, \ell = 1, 2, 3}} 
\eta_{i_2}\eta_{i_3}\eta_{i_4} W_{i_1j_1}W_{i_2j_2}W_{i_3j_3}\widetilde{\Omega}_{i_1i_4}. 
\end{align*}

We now consider the two claims of (\ref{T21varpf}) separately. Since the mean of $X_{a1}, X_{a2}, X_{a3}, X_b$ are all $0$, the first claim of (\ref{T21varpf}) follows trivially, so all remains to show is the second claim of (\ref{T21varpf}). 

We now consider  the second claim of (\ref{T21varpf}).  
By Cauchy-Schwartz inequality,  
\begin{align}  
\mathrm{Var}(T_{2a}) & \leq C \mathrm{Var}(X_{a1}) + \mathrm{Var}(X_{a2}) + \mathrm{Var}(X_{a3})  + \mathrm{Var}(X_{b}))  \nonumber \\ 
& \leq  C( \mathbb{E}[X_{a1}^2] + \mathbb{E}[X_{a2}^2] + \mathbb{E}[X_{a3}^2] + 
\mathbb{E}[X_{b}^2]). \label{T21-var-split} 
\end{align} 
We now consider $\mathbb{E}[X_{a1}^2] + \mathbb{E}[X_{a2}^2] + \mathbb{E}[X_{a3}^2]$, and $\mathbb{E}[X_{b}^2]$, separately.

Consider $\mathbb{E}[X_{a1}^2] + \mathbb{E}[X_{a2}^2] + \mathbb{E}[X_{a3}^2]$.  
We claim that both under the null and the alternative,  
\begin{align} 
\mathbb{E}[X_{a1}^2] & \leq C\alpha^2 \|\theta\|^{12} \|\theta\|_3^3/\|\theta\|_1^5 \label{T21-Varb-1},  \\
\mathbb{E}[X_{a2}^2] & \leq  C\alpha^2 \|\theta\|^{12} \|\theta\|_3^3/\|\theta\|_1^5  \label{T21-Varb-2},  \\ 
\mathbb{E}[X_{a3}^2] & \leq   C\alpha^2 \|\theta\|^{12} \|\theta\|_3^3/\|\theta\|_1^5 \label{T21-Varb-3}. 
\end{align} 
Combining these gives that both under the null and the alternative,    
\begin{equation} \label{T21-Vara} 
\mathbb{E}[X_{a1}^2] + \mathbb{E}[X_{a2}^2] + \mathbb{E}[X_{a3}^2]  \leq  C\alpha^2 \|\theta\|^{12} \|\theta\|_3^3/\|\theta\|_1^5. 
\end{equation}

It remains to show (\ref{T21-Varb-1})-(\ref{T21-Varb-3}).  
Since the proofs are similar, we only prove (\ref{T21-Varb-1}). Write  
\begin{align*} 
\mathbb{E}[X_{a1}^2] & = v^{-3} \sum_{\substack{i_1,i_2,i_3,i_4 (dist) \\ i_1',i_2',i_3',i_4'(dist)}}   
\sum_{\substack{j_3, j_3'  \\  j_3 \neq i_3, j_3' \neq i_3'}} 
\eta_{i_2}\eta_{i_3}\eta_{i_4}\eta_{i_2'}\eta_{i_3'}\eta_{i_4'}\mathbb{E}[ W^2_{i_1i_2} W_{i_3j_3} W^2_{i_1'i_2'} W_{i_3'j_3'}]\widetilde{\Omega}_{i_1i_4}\widetilde{\Omega}_{i_1'i_4'}. 
\end{align*} 
Consider the term
\[
 W^2_{i_1i_2} W_{i_3j_3} W^2_{i_1'i_2'} W_{i_3'j_3'}.
\] 
In order for the mean is nonzero, we have three cases 
\begin{itemize} 
\item Case A. $W_{i_1 i_2} = W_{i_3'j_3'}$ and $W_{i_3j_3} = W_{i_1'i_2'}$. 
\item Case B. $ W_{i_3j_3} = W_{i_3'j_3'}$ and $W_{i_1i_2} = W_{i_1' i_2'}$. 
\item Case C. $W_{i_3 j_3} = W_{i_3'j_3'}$ and $W_{i_1i_2} \neq W_{i_1' i_2'}$. 
\end{itemize} 
Consider Case A. In this case, $\{i_1', i_2', i_3'\}$ are three distinct indices in $\{ i_1, i_2, i_3, j_3\}$. In this case, 
\[
 W^2_{i_1i_2} W_{i_3j_3} W^2_{i_1'i_2'} W_{i_3'j_3'} 
= W_{i_1i_2}^3  W_{i_3j_3}^3, 
\] 
where by similar arguments as before 
\[
0 < \mathbb{E}[W_{i_1i_2}^3  W_{i_3j_3}^3] \leq C \Omega_{i_1 i_2} \Omega_{i_3 j_3} \leq C \theta_{i_1} \theta_{i_2} \theta_{i_3} \theta_{j_3}.  
\] 
At the same time, recall that that $0 < \eta_i  \leq C \theta_i$ for any $1 \leq i \leq n$, and that $|\widetilde{\Omega}_{ij}| \leq C \alpha \theta_i \theta_j$ for any $1 \leq i, j \leq n$, $i \neq j$,   where $\alpha = |\lambda_2/\lambda_1|$ with $\lambda_k$ being the $k$-th largest (in magnitude) eigenvalue of $\Omega$, $1 \leq k \leq K$.  By basic algebra, 
\[
|\eta_{i_2}\eta_{i_3}\eta_{i_4}\eta_{i_2'}\eta_{i_3'}\eta_{i_4'} \widetilde{\Omega}_{i_1i_4} \widetilde{\Omega}_{i_1'i_4'}|  \leq C\alpha^2  \theta_{i_1} \theta_{i_2} \theta_{i_3} \theta_{i_4}^2 \theta_{i_1'} 
\theta_{i_2'} \theta_{i_3'} \theta_{i_4'}^2.   
\] 
Note that in the current case, $\{i_1, i_2\} = \{i_3', j_3'\}$ and $\{i_3,j_3\} = \{i_1', i_2'\}$, so for some integers $0 \leq b_1, b_2 \leq 1$ and $b_1 + b_2 = 1$,  
\[
 \theta_{i_1} \theta_{i_2} \theta_{i_3} \theta_{i_4}^2 \theta_{i_1'} 
\theta_{i_2'} \theta_{i_3'} \theta_{i_4'}^2 =  \theta_{i_1}^{1 + b_1}  \theta_{i_2}^{1 + b_2} \theta_{i_3}^2 \theta_{j_3} \theta_{i_4}^2  \theta_{i_4'}^2.
\] 
Recall that $v \asymp \|\theta\|_1^2$. Combining these, the contribution of Case (A) to 
$\mathbb{E}[X_{a1}^2]$ is no greater than 
\[
C \alpha^2 (\|\theta\|_1)^{-6}  \sum_{i_1, i_2, i_3, i_4 (dist)} \sum_{i_4'} \sum_{j_3 (j_3 \neq i_3)}    
\sum_{b_1, b_2 (b_1 + b_2 = 1)}   \theta_{i_1}^{2 + b_1}  \theta_{i_2}^{2 + b_2} \theta_{i_3}^3 \theta_{j_3}^2 \theta_{i_4}^2  \theta_{i_4'}^2, 
\]  
where the right hand side $
\leq C  \alpha^2\cdot\|\theta\|^8\|\theta\|_3^6/\|\theta\|_1^6$.  
This shows that the contribution of Case (A) to $\mathbb{E}[X_{a1}^2]$ is no greater than 
\begin{equation} \label{T21-varXa1-A} 
C\alpha^2\cdot \|\theta\|^8 \|\theta\|_3^6 / \|\theta\|_1^6.  
\end{equation}

Consider Case B.  By similar arguments,    
\[
W^2_{i_1i_2} W_{i_3j_3} W^2_{i_1'i_2'} W_{i_3'j_3'} 
= W_{i_1i_2}^6  W_{i_3j_3}^2,  
\]
where 
\[ \mathbb{E}[W_{i_1i_2}^6  W_{i_3j_3}^2] \leq C \Omega_{i_1 i_2} \Omega_{i_3 j_3} \leq C \theta_{i_1} \theta_{i_2} \theta_{i_3} \theta_{j_3}, 
\] 
Also, by similar arguments,   
\[
|\eta_{i_2}\eta_{i_3}\eta_{i_4}\eta_{i_2'}\eta_{i_3'}\eta_{i_4'} \widetilde{\Omega}_{i_1i_4} \widetilde{\Omega}_{i_1'i_4'}|  \leq C\alpha^2  \theta_{i_1} \theta_{i_2} \theta_{i_3} \theta_{i_4}^2 \theta_{i_1'} 
\theta_{i_2'} \theta_{i_3'} \theta_{i_4'}^2, 
\] 
where as $W_{i_1 i_2} = W_{i_1'i_2'}$ and $W_{i_3 j_3} = W_{i_3' j_3'}$, the right hand side 
\[
\leq C \alpha^2  \theta_{i_1}^2 \theta_{i_2}^2  \theta_{i_3}^{1 + c_1} \theta_{j_3}^{c_2} \theta_{i_4}^2 \theta_{i_4'}^2, 
\] 
where $0 < c_1, c_2 \leq $ are integers satisfying $c_1 + c_2 = 1$. 
Recall $v \sim \|\theta\|_1^2$.  Combining these, the contribution of Case (B) to  $\mathbb{E}[X_{a1}^2]$   
\[
\leq C \alpha^2 (\|\theta\|_1)^{-6}  \sum_{i_1, i_2, i_3, i_4 (dist)} \sum_{i_4'} \sum_{j_3 (j_3 \neq i_3)}    
\sum_{b_1, b_2 (b_1 + b_2 = 1)}   \theta_{i_1}^{3}  \theta_{i_2}^{3} \theta_{i_3}^{2 + c_1} \theta_{j_3}^{1 + c_2} \theta_{i_4}^2  \theta_{i_4'}^2,   
\] 
where by $\|\theta\|^4 \leq \|\theta\|_1 \|\theta\|_3^3$,   the above term 
\[
\leq C \alpha^2 [\|\theta\|^4 \|\theta\|_3^9 / \|\theta\|_1^5, \|\theta\|^8 \|\theta\|_3^6/\|\theta\|_1^6] \leq C  \alpha^2  \|\theta\|^4 \|\theta\|_3^9 / \|\theta\|_1^5.  
\] 
This shows that the contribution of Case (B) to  
$\mathbb{E}[X_{a1}^2]$ is no greater than 
\begin{equation} \label{T21-varXa1-B} 
C \|\theta\|^4 \|\theta\|_3^9 / \|\theta\|_1^5.  
\end{equation}

Consider Case (C). In this case, 
\[
 W^2_{i_1i_2} W_{i_3j_3} W^2_{i_1'i_2'} W_{i_3'j_3'}  =
 W^2_{i_1i_2} W^2_{i_3j_3} W^2_{i_1'i_2'},  
\]
where by similar arguments, 
\[
\mathbb{E}[ W^2_{i_1i_2} W^2_{i_3j_3} W^2_{i_1'i_2'}] \leq C \Omega_{i_1i_2} \Omega_{i_3j_3} \Omega_{i_1'i_2'}   \leq C \theta_{i_1} \theta_{i_2} \theta_{i_3}  \theta_{j_3} \theta_{i_1'} \theta_{i_2'}.   
\] 
Also, by similar arguments, 
\[
|\eta_{i_2}\eta_{i_3}\eta_{i_4}\eta_{i_2'}\eta_{i_3'}\eta_{i_4'} \widetilde{\Omega}_{i_1i_4} \widetilde{\Omega}_{i_1'i_4'}|  \leq C\alpha^2  \theta_{i_1} \theta_{i_2} \theta_{i_3} \theta_{i_4}^2 \theta_{i_1'} 
\theta_{i_2'} \theta_{i_3'} \theta_{i_4'}^2, 
\] 
where as $W_{i_3 j_3} = W_{i_3' j_3'}$, the right hand side 
\[
\leq C \alpha^2  \theta_{i_1} \theta_{i_2} \theta_{i_3}^{1 + c_1} \theta_{j_3}^{c_2} \theta_{i_4}^2 \theta_{i_4'}^2, 
\] 
with the same $c_1, c_2$ as in the proof of Case B. 
Combining these and using $v \asymp \|\theta\|_1^2$,  
we have that under both the null and the alternative, the contribution of Case (C) to $\mathbb{E}[X_{a1}^2]$ 
\[
\leq C\alpha^2  (\|\theta\|_1)^{-6} \sum_{\substack{i_1, i_2, i_3, i_4 (dist) \\ i_1',i_2', i_4' (dist)}}  
\sum_{j_3 (j_3 \neq i_3)}  \theta^2_{i_1} \theta^{2}_{i_2} \theta^{2+c_1}_{i_3}  \theta_{j_3}^{1+c_2} \theta^2_{i_4} \theta^2_{i_1'} \theta^2_{i_2'}\theta^2_{i_4'}, 
\] 
where the right hand size 
\begin{equation} \label{T21-varXa1-C} 
\leq C \alpha^2\cdot[\|\theta\|^{12} \|\theta\|_3^3 / \|\theta\|_1^5 + \|\theta\|^{12} \|\theta\|_3^6/\|\theta\|_1^6] \leq C 
\alpha^2 \|\theta\|^{12}\|\theta\|_3^3 / \|\theta\|_1^5. 
\end{equation} 
Here we have again used $\|\theta\|^4 \leq  \|\theta\|_1 \|\theta\|_3^3$. 

Combining (\ref{T21-varXa1-A}), (\ref{T21-varXa1-B}), and (\ref{T21-varXa1-C})   gives 
\[
\mathbb{E}[X_{a1}^2] \leq C \alpha^2(\|\theta\|^8 \|\theta\|_3^6 / \|\theta\|_1^6 + \|\theta\|^4 \|\theta\|_3^9 / \|\theta\|_1^5 + \|\theta\|^8 \|\theta\|_3^9 /\|\theta\|^5] \leq C  \alpha^2 \|\theta\|^8 \|\theta\|_3^9 /\|\theta\|_1^5, 
\]
where we have used $\|\theta\|^4 \leq \|\theta\|_1 \|\theta\|_3^3$ and $\|\theta\| \goto \infty$.  
This proves (\ref{T21-Varb-1}).

We now consider $\mathbb{E}[X_b^2]$. Write 
\begin{align*} 
\mathbb{E}[X_{b}^2] & = v^{-3} \sum_{\substack{i_1,i_2,i_3,i_4 (dist) \\ i_1',i_2',i_3',i_4'(dist)}}   
\sum_{\substack{j_3, j_3'  \\  j_3 \neq i_3, j_3' \neq i_3'}} 
\eta_{i_2}\eta_{i_3}\eta_{i_4}\eta_{i_2'}\eta_{i_3'}\eta_{i_4'} \\
&\qquad   \mathbb{E}[ W_{i_1j_1} W_{i_2j_2}W_{i_3j_3}  W_{i_1'j_1'} W_{i_2' j_2'}  W_{i_3'j_3'}] \widetilde{\Omega}_{i_1i_4}\widetilde{\Omega}_{i_1'i_4'}. 
\end{align*} 
Consider
\[
W_{i_1j_1}W_{i_2j_2}W_{i_3j_3}, \qquad \mbox{and} \qquad  W_{i_1'j_1'}W_{i_2'j_2'}W_{i_3'j_3'}. 
\] 
Each term has a mean $0$, and two terms are uncorrelated with each other if and only if the two sets of random variables $\{ W_{i_1j_1}, W_{i_2j_2}, W_{i_3j_3}\}$ and $\{ W_{i_1'j_1'}, W_{i_2'j_2'}, W_{i_3'j_3'}\}$ are identical (however, it is possible that $W_{i_1j_1}$ does not equal to $W_{i_1'j_1'}$ but equals to $W_{i_2'j_2'}$, say).   
When this happens, first, $\{i_1, i_2, i_3, j_1, j_2, j_3\} = \{i_1', i_2', i_3', j_1', j_2', j_3'\}$. Recall that $|\widetilde{\Omega}_{ij}| \leq C \alpha \theta_i \theta_j$ for all $1 \leq i, j \leq n$, $i \neq j$, and that $0 < \eta_i \leq C \theta_i$ for all $1 \leq i \leq n$. For 
 integers $a_i \in \{0, 1\}$, $1 \leq i \leq 4$, that satisfy $\sum_{i=1}^6 a_i = 3$,  we have  
\begin{align*} 
|\eta_{i_2}\eta_{i_3}\eta_{i_4}\eta_{i_2'}\eta_{i_3'}\eta_{i_4'} \widetilde{\Omega}_{i_1i_4} \widetilde{\Omega}_{i_1'i_4'}|  & \leq C \eta_{i_1}^{a_1}  \eta_{j_1}^{a_2}   \eta_{i_2}^{1 + a_3} \eta_{j_2}^{a_4}  \eta_{i_3}^{1 + a_5} \eta_{j_3}^{a_6}  \eta_{i_4} \eta_{i_4'} |\widetilde{\Omega}_{i_1i_4}|  |\widetilde{\Omega}_{i_1'i_4'}| \\
&\leq  C \alpha^2   \theta_{i_1}^{1 + a_1}   \eta_{j_1}^{a_2}   \eta_{i_2}^{1 + a_3} \eta_{j_2}^{a_4}  \eta_{i_3}^{1 + a_5} \eta_{j_3}^{a_6}  \eta_{i_4}^2 \eta_{i_4'}^2. 
\end{align*} 
Second, 
\[
\mathbb{E}[ W_{i_1j_1} W_{i_2j_2}W_{i_3j_3}  W_{i_1'j_1'} W_{i_2' j_2'}  W_{i_3'j_3'}] = 
\mathbb{E}[ W_{i_1j_1}^2 W_{i_2j_2}^2 W_{i_3j_3}^2], 
\]  
where by similar arguments, the right hand side 
\[
\leq C \Omega_{i_1j_1} \Omega_{i_2j_2} \Omega_{i_3j_3} \leq C \theta_{i_1} \theta_{j_1} \theta_{i_2} \theta_{j_2} \theta_{i_3}  \theta_{j_3}. 
\] 
Recall that $v \sim \|\theta\|_1^2$. Combining these gives 
\[
\mathbb{E}[X_b^2] \leq C \alpha^2 \|\theta\|_1^{-6}    \sum_{i_1, i_2, i_3, i_4 (dist)} 
\sum_{i_4'}  \sum_{\substack{j_1, j_2, j_3 \\ j_1\neq i_1, j_2\neq i_2,  j_3\neq i_3}} \sum_a
\theta_{i_1}^{2 + a_1}   \eta_{j_1}^{1+a_2}   \eta_{i_2}^{2 + a_3} \eta_{j_2}^{1+a_4}  \eta_{i_3}^{2 + a_5} \eta_{j_3}^{1+a_6}  \eta_{i_4}^2 \eta_{i_4'}^2,  
\] 
where $a = (a_1, a_2, \ldots, a_6)$ as above. By the way $a_i$ are defined, the right hand side  
\[
\leq C \alpha^2    \|\theta\|^4  (\sum_{a} \|\theta\|_{a_1 + 2}^{a_1 + 2}  \cdot \|\theta\|_{a_2 + 1}^{a_2 + 1}  \cdot \|\theta\|_{a_3 + 2}^{a_3 + 2} \cdot  \|\theta\|_{a_4 + 1}^{a_4 + 1}\|\theta\|_{a_5 +2}^{a_5 + 2}\|\theta\|_{a_6 + 1}^{a_6 + 1})/\|\theta\|_1^6,    
\] 
which by $\|\theta\|^4 \leq \|\theta\|_1 \|\theta\|_3^3$, the term in the bracket does not exceed 
\[
C \max\{\|\theta\|^{12},  \|\theta\|_1 \|\theta\|^8 \|\theta\|_3^3,  \|\theta\|_1^2 \|\theta\|^4 \|\theta\|_3^6, \|\theta\|_1^3 \|\theta\|_3^9\} \leq C \|\theta\|_1^3 \|\theta\|_3^9.    
\] 
Combining these gives 
\begin{equation} \label{T21-Varb} 
\mathbb{E}[X_b^2] \leq C \alpha^2 \|\theta\|^4 \|\theta\|_3^9 / \|\theta\|_1^3. 
\end{equation} 

Finally, inserting (\ref{T21-Vara})-(\ref{T21-Varb}) into (\ref{T21-var-split}) gives 
\[
\mathrm{Var}(T_{2a}) \leq C \alpha^2[\|\theta\|^8 \|\theta\|_3^3/\|\theta\|_1^5 +  \|\theta\|^4 \|\theta\|_3^9/\|\theta\|_1^3] \leq C \alpha^2\|\theta\|^4 \|\theta\|_3^9/\|\theta\|_1^3, 
\] 
and (\ref{T21varpf}) follows.  

Consider Item (g) and Item (h). The proof are similar, so we only show Item (g).  The goal is to show (\ref{T23varpf}). 
Recall that 
\beq\label{def:t23-add}
T_{2c} = \sum_{i_1, i_2, i_3, i_4(dist)}\eta_{i_1}\eta_{i_3}\eta_{i_4}\big[ (\eta_{i_2} - \tilde{\eta}_{i_2})^2  (\eta_{i_3} - \tilde{\eta}_{i_3})    \big]\cdot  \widetilde{\Omega}_{i_4i_1},
\eeq
and 
\[
\widetilde{\eta}  - \eta = v^{-1/2} W 1_n. 
\]
Plugging this into $T_{2c}$ gives (note symmetry in $\widetilde{\Omega}$) 
\begin{align*}
T_{2c} &=  - \frac{1}{v^{2/3}} \sum_{i_1, i_2,i_3, i_4 (dist)}\eta_{i_1}\eta_{i_3}\eta_{i_4} \Bigl( \sum_{j_2\neq i_2}W_{i_2j_2}\Bigr)\Bigl( \sum_{\ell_2\neq i_2}W_{i_2\ell_2}\Bigr)\Bigl(\sum_{j_3\neq i_3}W_{i_3j_3}\Bigr)\widetilde{\Omega}_{i_4 i_1}\cr
&= -\frac{1}{v^{3/2}}\sum_{\substack{i_1, i_2,i_3, i_4  (dist)\\ j_1\neq i_1, j_2\neq i_2, j_3\neq i_3}}\eta_{i_1}\eta_{i_3}\eta_{i_4} W_{i_2j_2}W_{i_2\ell_2}W_{i_3j_3}\widetilde{\Omega}_{i_1i_4}.
\end{align*}
By basic combinatorics and careful observations, we have 
\beq \label{2}
W_{i_2j_2}W_{i_2\ell_2}W_{i_3j_3} = 
\begin{cases}
 W_{i_2i_3}^3      ,&\mbox{if }  j_1 = \ell_2 = i_3, j_3 = i_2,  \\
W^2_{i_2j_2}W_{i_3j_3},&\mbox{if } j_1 = \ell_2,  (j_2, j_3) \neq (i_3, i_2), \\
W^2_{i_2i_3}W_{i_2\ell_2}, & \mbox{if }  j_2 = i_3,  j_3 = i_2, \ell_2\neq i_3, \\
W^2_{i_2i_3}W_{i_2j_2}, & \mbox{if }  \ell_2 = i_3,  j_3 = i_2, j_2\neq i_3, \\
W_{i_2j_2}W_{i_2\ell_2}W_{i_3j_3} , & \mbox{otherwise}. 
\end{cases}
\eeq
This allows us to further split $T_{2c}$ into $4$ different terms: 
\begin{equation} \label{T23-split} 
T_{2c}   =  Y_a +  Y_{b1} + Y_{b2} + Y_{b3} + Y_{c}, 
\end{equation} 
\begin{align*}
Y_{a} & =  -\frac{1}{v^{3/2}}\sum_{i_1, i_2, i_3, i_4 (dist)}\sum_{j_3\neq i_3}\eta_{i_1}\eta_{i_3}\eta_{i_4} W^3_{i_2i_3}\widetilde{\Omega}_{i_1i_4},  \cr
Y_{b1} & = -\frac{1}{v^{3/2}}\sum_{i_1, i_2, i_3, i_4 (dist)}\sum_{j_3\neq i_3}\eta_{i_1}\eta_{i_3}\eta_{i_4} W^2_{i_2j_2}W_{i_3j_3}\widetilde{\Omega}_{i_1i_4}, \cr
Y_{b2} & = -\frac{1}{v^{3/2}}\sum_{i_1, i_2, i_3, i_4(dist)}\sum_{j_2\neq i_2}\eta_{i_1}\eta_{i_3}\eta_{i_4} W^2_{i_2i_3}W_{i_2\ell_2}\widetilde{\Omega}_{i_1i_4}, \cr
Y_{b3} & = -\frac{1}{v^{3/2}}\sum_{i_1, i_2, i_3, i_4 (dist)}\sum_{j_1\neq i_1}\eta_{i_1}\eta_{i_3}\eta_{i_4} W^2_{i_2i_3}W_{i_2j_2}\widetilde{\Omega}_{i_1i_4}, \cr
Y_c &=-\frac{1}{v^{3/2}}\sum_{i_1, i_2, i_3, i_4 (dist)}
\sum_{\substack{j_2, \ell_2, j_3 \\ j_2\neq i_2, \ell_2\neq i_2, j_3\neq i_3\\  j_2\neq i_3, \ell_2\neq i_3, j_3\neq i_2  }} 
\eta_{i_1}\eta_{i_3}\eta_{i_4} W_{i_2j_2}W_{i_2\ell_2}W_{i_3j_3} \widetilde{\Omega}_{i_1i_4}. 
\end{align*}
We now show the two claims in (\ref{T23varpf}) separately. Consider the first claim. It is seen that out of the $5$ terms on the right hand side of 
(\ref{T23-split}), the mean of all terms are $0$, except for the first one.  Note that for any $1\leq i, j\leq n, i\neq j$, $\mathbb{E}[W_{ij}^3]\leq C\Omega_{ij}$. Together with $\Omega_{ij}\leq C\theta_i\theta_j$, $\widetilde{\Omega}_{ij}\leq C\alpha\theta_i\theta_j$, $0<\eta_i<C\theta_i$ and $v\sim \|\theta\|_1^2$, it follows \begin{align*}
\mathbb{E}[|Y_a|] & \leq \frac{1}{v^{3/2}}\sum_{i_1, i_2, i_3, i_4 (dist)} \eta_{i_1}\eta_{i_3}\eta_{i_4} \Omega_{i_2i_3}\widetilde{\Omega}_{i_1i_4} \\
& \leq C\alpha\cdot\frac{1}{\|\theta\|_1^3}\sum_{i_1, i_2, i_3, i_4 (dist)} \theta_{i_1}^2\theta_{i_2}\theta_{i_3}^2\eta_{i_4}^2, 
\end{align*}
where the last term is no greater than $
C\alpha\cdot {\|\theta\|^6}/{\|\theta\|_1^3},$ and the first claim of (\ref{T23varpf})  follows.

Consider the second claim of (\ref{T23varpf}).  
By Cauchy-Schwartz inequality,  
\begin{align}  
\mathrm{Var}(T_{2c}) & \leq C (\mathrm{Var}(Y_{a})   + \mathrm{Var}(Y_{b1}) + \mathrm{Var}(Y_{b2}) + \mathrm{Var}(Y_{b3})  + \mathrm{Var}(Y_{c}))  \nonumber \\ 
& \leq  C( \mathrm{Var}(Y_{a}) + \mathbb{E}[Y_{b1}^2] + \mathbb{E}[Y_{b2}^2] + \mathbb{E}[Y_{b3}^2] + 
\mathbb{E}[Y_{c}^2]). \label{T23-var-split} 
\end{align} 
We now study $\mathrm{Var}(Y_{a})$.  Write 
\begin{align*} 
\mathrm{Var}(Y_{a}) & = v^{-3} \sum_{\substack{i_1, i_2, i_3, i_4 (dist) \\ i_1', i_2', i_3', i_4'   (dist) } } 
 \eta_{i_1}  \eta_{i_3}    \eta_{i_4}  \eta_{i_1'}  \eta_{i_3'}  \eta_{i_4'}  \mathbb{E}[(W_{i_2i_3}^3 - \mathbb{E}[W_{i_2i_3}^3])  (W_{i_2'i_3'}^3 - \mathbb{E}[W_{i_2'i_3'}^3]) ]\cdot\widetilde{\Omega}_{i_1i_4}\widetilde{\Omega}_{i_1'i_4'}. 
\end{align*} 
Fix a term $(W_{i_2i_3}^3 - \mathbb{E}[W_{i_2i_3}^3])  (W_{i_2'i_3'}^3 - \mathbb{E}[W_{i_2'i_3'}^3])$. When the mean is nonzero, we must have $\{i_2, i_3\} = \{i_2', i_3'\}$, and when this happens, 
\[
\mathbb{E}[(W_{i_2i_3}^3 - \mathbb{E}[W_{i_2i_3}^3])  (W_{i_2'i_3'}^3 - \mathbb{E}[W_{i_2'i_3'}^3]) ] = \mathrm{Var}(W_{i_2i_3}^3). 
\]
For a random variable $X$, we have $\mathrm{Var}(X)\leq \mathbb{E}[X^2]$, and it follows that \[
\mathrm{Var}(W_{i_2i_3}^3) \leq \mathbb{E}[W_{i_2i_3}^6] \leq \mathbb{E}[W_{i_2i_3}^2], 
\]
where we have used the property that $0\leq W_{i_2i_3}^2\leq 1$. Notice that $\mathbb{E}[W_{i_2i_3}^2]\leq C\theta_{i_2}\theta_{i_3}$, and recall that $v \asymp \|\theta\|_1^2$, $\widetilde{\Omega}_{ij}\leq C\alpha\theta_i\theta_j$ and $0 < \eta_i \leq C \theta_i$ for all $1 \leq i \leq n$. Combining these gives 
\begin{equation} \label{Y3a0-Var}
\mathrm{Var}(Y_{a}) \leq C \alpha^2 (\|\theta\|_1^{-6})  \cdot \sum_{\substack{i_1, i_2, i_3, i_4 (dist) \\ i_1', i_4'   (dist) } }
 \theta_{i_1}^2   \theta_{i_2} \theta_{i_3}^3   \theta_{i_4}^2  \theta_{i_1'}^2 \theta_{i_4'}^2 \leq C \alpha^2 \|\theta\|^8\|\theta\|_3^3/\|\theta\|_1^5. 
\end{equation} 

Additionally, note that 
\begin{itemize}
\item The proof of $Y_{b1}$, $Y_{b2}$, and $Y_{b3}$ is similar to that of $X_{a1}$ in Item (e). 
\item The proof of $Y_c$ is similar to that of $X_b$ in Item (e). 
\end{itemize} 
For these reasons, we skip the proof details, but only to state that,  
both under the null and the alternative,  
\begin{align} 
\mathbb{E}[Y_{b1}^2] & \leq  C \alpha^2\|\theta\|^{8} \|\theta\|_3^3/ \|\theta\|_1  \label{T22-Varb-1},  \\
\mathbb{E}[Y_{b2}^2] & \leq  C\alpha^2 \|\theta\|^{12} \|\theta\|_3^3/\|\theta\|_1^5  \label{T22-Varb-2},  \\ 
\mathbb{E}[Y_{b3}^2] & \leq   C\alpha^2 \|\theta\|^{12} \|\theta\|_3^3/\|\theta\|_1^5 \label{T22-Varb-3}, 
\end{align} 
and therefore,  
\begin{equation} \label{T23-Vara} 
\mathbb{E}[Y_{b1}^2] + \mathbb{E}[Y_{b2}^2] + \mathbb{E}[Y_{b3}^2]  \leq  C \alpha^2\|\theta\|^{8} \|\theta\|_3^3/ \|\theta\|_1. 
\end{equation} 
At the same time, both under the null and the alternative, 
\begin{equation} \label{T23-Varb} 
\mathbb{E}[Y_{c}^2] \leq C\alpha^2\cdot \|\theta\|^{10} \|\theta\|_3^3 / \|\theta\|_1^3. 
\end{equation} 
Inserting (\ref{T23-Vara}) and (\ref{T23-Varb}) into (\ref{T23-var-split}) gives  
\[
\mathbb{E}[T_{2c}^2] \leq C\alpha^2 [ \|\theta\|^8\|\theta\|_3^3/\|\theta\|_1^5 +  \|\theta\|^{8} \|\theta\|_3^3/ \|\theta\|_1+ \|\theta\|^{10} \|\theta\|_3^3 / \|\theta\|_1^3] \leq C \alpha^2\|\theta\|^{8} \|\theta\|_3^3/ \|\theta\|_1. 
\] 
This proves (\ref{T23varpf}). 

Consider Item (i). 
The goal is to show   (\ref{F1varpf}). 
Recall that 
\beq\label{def:F1}
F_{a} =  \sum_{i_1,i_2,i_3,i_4(dist)}\eta_{i_1}\eta_{i_2}\eta_{i_3}\eta_{i_4}\big[(\eta_{i_1} - \tilde{\eta}_{i_1}) (\eta_{i_2} - \tilde{\eta}_{i_2})  (\eta_{i_3} - \tilde{\eta}_{i_3})  (\eta_{i_4} - \tilde{\eta}_{i_4})    \big],
\eeq
and that for any $1 \leq i \leq n$,  
\[
\tilde{\eta}_{i} - \eta_{i}  = v^{-1/2} \sum_{j  \neq i}^n W_{ij}. 
\] 
Inserting it into (\ref{def:F1}) gives 
\[
F_{a} =  \sum_{i_1,i_2,i_3,i_4(dist)}\eta_{i_1}\eta_{i_2}\eta_{i_3}\eta_{i_4}\big[(\eta_{i_1} - \tilde{\eta}_{i_1}) (\eta_{i_2} - \tilde{\eta}_{i_2})  (\eta_{i_3} - \tilde{\eta}_{i_3})  (\eta_{i_4} - \tilde{\eta}_{i_4})    \big],
\] 
By basic combinatorics and basic algebra, we have 
\[
W_{i_1 j_1} W_{i_2 j_2} W_{i_3 j_3} W_{i_4 j_4} 
= 
\left\{ 
\begin{array}{ll} 
W_{i_1i_2}^2 W_{i_3i_4}^2, &\qquad \mbox{if $(i_1, j_1) = (j_2, i_2), (i_3, j_3) = (j_4, i_4)$}, \\
W_{i_1i_3}^2 W_{i_2 i_4}^2, &\qquad   \mbox{if $(i_1, j_1) = (j_3, i_3), (i_2, j_2) = (j_4, i_4)$}, \\
W_{i_1i_4}^2 W_{i_2 i_3}^2,  &\qquad   \mbox{if $(i_1, i_4) = (j_4, i_1), (i_2, j_2) = (j_3, i_3)$},  \\ 
W_{i_1i_2}^2 W_{i_3j_3}W_{i_4 j_4}, &\qquad  \mbox{if $(i_1,j_1) = (j_2, i_2), (j_4, j_3) \neq (i_3, i_4)$},   \\ 
W_{i_1i_3}^2 W_{i_2j_2} W_{i_4j_4}, &\qquad  \mbox{if $(i_1,j_1) = (j_3, i_3), (j_4, j_2) \neq (i_2, i_4)$},    \\
W_{i_1i_4}^2 W_{i_2j_2} W_{i_3j_4}, &\qquad  \mbox{if $(i_1,j_1) = (j_4, i_4), (j_3, j_2) \neq (i_2, i_3)$},  \\ 
W_{i_2i_3}^2 W_{i_1j_1} W_{i_4j_4}, &\qquad  \mbox{if $(i_2,j_2) = (j_3, i_3), (j_4, j_1) \neq (i_1, i_4)$},    \\ 
W_{i_2i_4}^2 W_{i_1j_1} W_{i_3j_3}, &\qquad  \mbox{if $(i_2, j_2) = (j_4, i_4), (j_3, j_1) \neq (i_1, i_3)$},  \\
W_{i_3i_4}^2 W_{i_1j_1} W_{i_2j_2}, &\qquad  \mbox{if $(i_3,j_3) = (j_4, i_4), (j_2, j_1) \neq (i_1, i_2)$}.\\ 
W_{i_1j_1} W_{i_2j_2} W_{i_3 j_3} W_{i_4j_4}, &\qquad \mbox{otherwise}. 
\end{array} 
\right. 
\] 
By symmetry, it allows us to further split  $F_1$ into $3$ different terms: 
\beq\label{F1split}
F_1 = 3 X_a  + 6 X_b + X_c, 
\eeq
where 
\[
X_a = v^{-2} \sum_{i_1, i_2, i_3, i_4 (dist)} \eta_{i_1} \eta_{i_2} \eta_{i_3} \eta_{i_4}   W_{i_1i_2}^2 W_{i_3i_4}^2, 
\] 
\[
X_b = v^{-2} \sum_{i_1, i_2, i_3, i_4 (dist)} \sum_{\substack{j_3, j_4 \\
 (j_3, j_4) \neq (i_4, i_3)}}    \eta_{i_1} \eta_{i_2} \eta_{i_3} \eta_{i_4}  W_{i_1i_2}^2 W_{i_3j_3} W_{i_4 j_4}, 
\]
and
\[
X_c = v^{-2} \sum_{i_1, i_2, i_3, i_4 (dist)} \sum_{\substack{j_1, j_2, j_3, j_4 \\ 
j_k \neq i_{\ell}, k, \ell = 1, 2, 3, 4 }}  
  \eta_{i_1} \eta_{i_2} \eta_{i_3} \eta_{i_4}  W_{i_1j_1}W_{i_2 j_2}  W_{i_3j_3} W_{i_4 j_4}. 
\]

We now show the two claims in  (\ref{F1varpf}) separately.  Consider the first claim of (\ref{F1varpf}). Note that $\mathbb{E}[X_b]  = \mathbb{E}[X_c] = 0$.   Recall that both under the null and the alternative,  for any $i \neq j$, $\mathbb{E}[W_{ij}^2] = \Omega_{ij}(1 - \Omega_{ij}) \leq C \theta_i \theta_j$, and 
that $0 < \eta_i \leq C \theta_i$, and that $v \asymp \|\theta\|_1^2$,  
Therefore,  
\[
0 < \mathbb{E}[X_a] \leq  v^{-2} \sum_{i_1, i_2, i_3, i_4 (dist)} \theta_{i_1}\theta_{i_2} \theta_{i_3}\theta_{i_4} \theta_{i_1} \theta_{i_2} \theta_{i_3} \theta_{i_4} \leq C \|\theta\|^8 / \|\theta\|_1^4. 
\]
Inserting into \eqref{F1split} gives \[
\mathbb{E}[|F_{1}|]\leq  C \|\theta\|^8 / \|\theta\|_1^4,
\]
and the first claim (\ref{F1varpf}) follows.

Consider the second claim (\ref{F1varpf}) next.  By \eqref{F1split} and Cauchy-Schwarz inequality,  
\begin{align}  
\mathrm{Var}(F_1) & \leq C (\mathrm{Var}(X_a) + \mathrm{Var}(X_b) + \mathrm{Var}(X_c))   \leq  C(\mathrm{Var}(X_a)+ \mathbb{E}[X_b^2] +  \mathbb{E}[X_c^2]). \label{F1-var-split} 
\end{align} 
We now consider $\mathrm{Var}(X_a)$, $\mathbb{E}[X_b^2]$, and $\mathbb{E}[X_c^2]$, separately.     
Note that 
\begin{itemize}
\item The proof of $\mathrm{Var}(X_{a})$ is similar to that of $\mathrm{Var}(X_a)$ in Item (a).
\item The proof of $\mathbb{E}[X^2_{b}] $ is similar to that of $\sum_{s=1}^4\mathbb{E}[X^2_{bs}] $ in Item (a).
\item The proof of $\mathbb{E}[X_c^2]$ is similar to that of $\mathbb{E}[X_c^2]$ in Item (a). 
\end{itemize}
For these reasons, we omit the proof details and only state the claims. We have that under both the null and the alternative,
\beq\label{F1a-Var}
\mathrm{Var}(X_{a})\leq C\|\theta\|^8\|\theta\|_3^6/\|\theta\|_1^8.
\eeq
\beq\label{F1b-Var}
\mathrm{Var}(X_b^2) + \mathrm{Var}(Y_{a3})\leq  C\|\theta\|^4\|\theta\|_3^6/\|\theta\|_1^4,
\eeq
\beq\label{F1c-Var}
\mathbb{E}[X_c^2]  \leq C\|\theta\|_3^{12}/\|\theta\|_1^4,
\eeq
Finally, inserting \eqref{F1a-Var},  \eqref{F1b-Var}, and  \eqref{F1c-Var} into  \eqref{F1split} gives that, both under the null and the alternative,
\[
\mathrm{Var}(F_1)\leq C[  \|\theta\|^8\|\theta\|_3^6 / \|\theta\|_1^8 +  \|\theta\|^8\|\theta\|_3^6 / \|\theta\|_1^6  +  \|\theta\|_3^{12} / \|\theta\|_1^4] \leq C\|\theta\|^8\|\theta\|_3^6 / \|\theta\|_1^6, 
\]
where we have used $\|\theta\|\goto\infty$ and $\|\theta\|_3^3\ll\|\theta\|^2\ll\|\theta\|_1$. This gives (\ref{F1varpf}) and completes the proof for Item (i).

 
 Consider Item (j). 
The goal is to show  (\ref{F2varpf}). Recall that 
\[
F_{b} = \sum_{i_1,i_2,i_3,i_4 (dist)} \eta_{i_2}\eta_{i_3}^2\eta_{i_4}\big[(\eta_{i_1} - \tilde{\eta}_{i_1})^2 (\eta_{i_2} - \tilde{\eta}_{i_2}) (\eta_{i_4} - \tilde{\eta}_{i_4})     \big],
\]
and that \[
\widetilde{\eta} - \eta = v^{-1/2} W 1_n.
\] 
Plugging this into $F_b$, we have 
\[
F_b = v^{-2}\sum_{i_1,i_2,i_3,i_4 (dist)} \sum_{\substack{j_1, \ell_1,  j_2, j_4 \\
j_1 \neq i_1, \ell_1 \neq i_1, j_2 \neq i_2, j_4 \neq i_4}}  \eta_{i_2} \eta_{i_3}^2 \eta_{i_4}   
W_{i_1 j_1} W_{i_1 \ell_1} W_{i_2 j_2} W_{i_4 j_4}. 
\]
By basic combinatorics and basic algebra, we have 
\[
W_{i_1 j_1} W_{i_1 \ell_1} W_{i_2 j_2} W_{i_4 j_4}
= 
\left\{ 
\begin{array}{ll} 
W_{i_1i_2}^3 W_{i_4 j_4},  &\qquad    \mbox{if $j_1, \ell_1 = i_2, j_2 = i_1$},  \\ 
W_{i_1i_4}^3 W_{i_2 j_2},  &\qquad    \mbox{if $j_1, \ell_1 = i_4, j_4 = i_1$},  \\
W_{i_1i_2}^2 W_{i_1 i_4}^2, &\qquad   \mbox{if $(j_1, j_2) = (i_2, i_1), (\ell_1, j_4) = (i_4, i_1)$}, \\
W_{i_1i_2}^2 W_{i_1 i_4}^2,  &\qquad    \mbox{if $(\ell_1, j_2) = (i_2, i_1), (j_1, j_4) = (i_4, i_1)$},  \\ 
W_{i_1i_4}^2 W_{i_1 i_2}^2, &\qquad   \mbox{if $(j_1, j_4) = (i_4, i_1), (\ell_1, j_2) = (i_2, i_1)$}, \\
W_{i_1i_4}^2 W_{i_1 i_2}^2,  &\qquad    \mbox{if $(\ell_1, j_4) = (i_4, i_1), (j_1, j_2) = (i_2, i_1)$},  \\ 
W_{i_1j_1}^2 W_{i_2i_4}^2, &\qquad \mbox{if $j_1 =\ell_1, (j_2, j_4) = (i_4, i_2)$}, \\
W_{i_1i_2}^2 W_{i_1j_1} W_{i_4j_4}, &\qquad  \mbox{if $\ell_1 = i_2, j_2 = i_1, j_1\neq i_2, i_4$},    \\
W_{i_1i_2}^2 W_{i_1\ell_1} W_{i_4j_4}, &\qquad  \mbox{if $j_1 = i_2, j_2 = i_1, \ell_1\neq i_2, i_4$},  \\ 
W_{i_1i_4}^2 W_{i_1j_1} W_{i_2j_2}, &\qquad  \mbox{if $\ell_1 = i_4, j_4 = i_1, \ell_1\neq i_2, i_4$},    \\
W_{i_1i_4}^2 W_{i_1\ell_1} W_{i_2j_2}, &\qquad  \mbox{if $j_1 = i_4, j_4 = i_1, j_1\neq i_2, i_4$},    \\  
W_{i_2i_4}^2W_{i_1j_1} W_{i_1\ell_1} , &\qquad  \mbox{if $j_1\neq\ell_1, (j_2, j_4) = (i_4, i_2)$}.\\ 
W_{i_1j_1}^2 W_{i_2j_2}W_{i_4 j_4}, &\qquad  \mbox{if $j_1 = \ell_1, (j_1, j_2) \neq (i_2, i_1), (j_1, j_4) \neq (i_4, i_1)$},   \\ 
W_{i_1 j_1} W_{i_1 \ell_1} W_{i_2 j_2} W_{i_4 j_4}, &\qquad \mbox{otherwise}. 
\end{array} 
\right. 
\] 
By these and symmetry, we can further split $F_b$ into $7$ different terms,

We decompose \beq\label{F2split}
F_b = 2Y_{a1} + 4Y_{a2} + Y_{a3}  + 4Y_{b1} + Y_{b2} + Y_{b3} + Y_{c},
\eeq
where 
\[
Y_{a1} =  v^{-2}\sum_{i_1, i_2, i_3, i_4(dist)}\sum_{j_4, j_4\neq i_4}\eta_{i_2}\eta_{i_3}^2\eta_{i_4}W_{i_1i_2}^3 W_{i_4 j_4}, 
\]
\[
Y_{a2} =  v^{-2}\sum_{i_1, i_2, i_3, i_4(dist)}\eta_{i_2}\eta_{i_3}^2\eta_{i_4}W_{i_1i_2}^2 W_{i_1 i_4}^2, 
\]
\[
Y_{a3} =  v^{-2}\sum_{i_1, i_2, i_3, i_4(dist)}\sum_{j_1, j_1\neq i_1}\eta_{i_2}\eta_{i_3}^2\eta_{i_4}W_{i_1j_1}^2 W_{i_2i_4}^2, 
\]
\[
Y_{b1} =  v^{-2}\sum_{i_1, i_2, i_3, i_4(dist)}\sum_{\substack{j_1, j_4\\ j_1\neq i_1, j_4\neq i_4}}\eta_{i_2}\eta_{i_3}^2\eta_{i_4}W_{i_1i_2}^2 W_{i_1j_1} W_{i_4j_4}, 
\]
\[
Y_{b2} =  v^{-2}\sum_{i_1, i_2, i_3, i_4(dist)}\sum_{\substack{j_1, \ell_1\\ j_1, \ell_1 \neq i_1}}\eta_{i_2}\eta_{i_3}^2\eta_{i_4}W_{i_2i_4}^2 W_{i_1j_1} W_{i_1\ell_1} , 
\]
\[
Y_{b3} =  v^{-2}\sum_{i_1, i_2, i_3, i_4(dist)}\sum_{\substack{j_1, j_2,  j_4\\ j_1\neq i_1, j_2\neq i_2, j_4\neq i_4}}\eta_{i_2}\eta_{i_3}^2\eta_{i_4}W_{i_1j_1}^2 W_{i_2j_2}W_{i_4 j_4}, 
\]
\[
Y_{c} =  v^{-2}\sum_{i_1, i_2, i_3, i_4(dist)}\sum_{\substack{j_1, \ell_1, j_2, j_4\\ j_1, \ell_1\notin\{i_1, i_2, i_4\} \\ j_2\notin\{i_1, i_4\}, j_4\notin\{i_1, i_2\} }}\eta_{i_2}\eta_{i_3}^2\eta_{i_4}W_{i_1 j_1} W_{i_1 \ell_1} W_{i_2 j_2} W_{i_4 j_4}, 
\]
We now consider the two claims in (\ref{F2varpf}) separately.  Consider the first claim. It is seen that only the second and the third terms above have non-zero mean.   Recall that both under the null and the alternative,  for any $i \neq j$, $\mathbb{E}[W_{ij}^2] = \Omega_{ij}(1 - \Omega_{ij}) \leq C \theta_i \theta_j$, $0 < \eta_i \leq C \theta_i$, and that $v \asymp \|\theta\|_1^2$. 
It follows \beq\label{mean-Ya2}
0 < \mathbb{E}[Y_{a2}] \leq  v^{-2} \sum_{i_1, i_2, i_3, i_4 (dist)}\theta_{i_2}\theta^2_{i_3} \theta_{i_4}\cdot \theta_{i_1}^2 \theta_{i_2}  \theta_{i_4} \leq C \|\theta\|^8 / \|\theta\|_1^4. 
\eeq
and
\beq\label{mean-Ya3}
0 < \mathbb{E}[Y_{a3}] \leq  v^{-2} \sum_{i_1, i_2, i_3, i_4 (dist)}\sum_{j_1} \theta_{i_2}\theta^2_{i_3} \theta_{i_4}\cdot \theta_{i_1} \theta_{i_2} \theta_{j_1} \theta_{i_4} \leq C \|\theta\|^6 / \|\theta\|_1^2. 
\eeq
Combining \eqref{mean-Ya2}, \eqref{mean-Ya3} with \eqref{F2split} gives \[
\mathbb{E}[|F_2|] \leq C[\|\theta\|^8 / \|\theta\|_1^4 +\|\theta\|^6 / \|\theta\|_1^2] \leq C \|\theta\|^6 / \|\theta\|_1^2,
\]
where we've used the universal inequality that $\|\theta\|^2\leq\|\theta\|_1$. It follows the first claim of (\ref{F2varpf}). 

We now show the second claim of (\ref{F2varpf}). By Cauchy-Schwarz inequality, 
\begin{align}  
\mathrm{Var}(F_b) & \leq C \big( \mathrm{Var}(Y_{a1}) + \mathrm{Var}(Y_{a2}) + \mathrm{Var}(Y_{a3}) + \mathrm{Var}(Y_{b1}) + \mathrm{Var}(Y_{b2}) + \mathrm{Var}(Y_{b3}) + \mathrm{Var}(Y_{c})  \big) \nonumber \\ 
& \leq C \big( \mathrm{Var}(Y_{a1}) + \mathrm{Var}(Y_{a2}) + \mathrm{Var}(Y_{a3}) + \mathbb{E}[Y^2_{b1}] + \mathbb{E}[Y^2_{b2}]  + \mathbb{E}[Y^2_{b3}]  + \mathbb{E}[Y^2_{c}]  \big). \label{varF2a-add} 
\end{align} 
We now consider $\mathrm{Var}(Y_{a1})$, $ \mathrm{Var}(Y_{a2}) + \mathrm{Var}(Y_{a3}) $, $\mathbb{E}[Y^2_{b1}] + \mathbb{E}[Y^2_{b2}]  + \mathbb{E}[Y^2_{b3}] $, and $\mathbb{E}[Y^2_{c}]$, separately. 
Note that 
\begin{itemize}
\item The proof of $\mathrm{Var}(Y_{a1})$ is similar to that of $\mathrm{Var}(Y_a)$ in Item (b).
\item The proof of $\mathrm{Var}(Y_{a2})$ and $\mathrm{Var}(Y_{a3})$ are similar to that of $\mathrm{Var}(X_a)$ in Item (a).
\item The proof of $\sum_{s=1}^3\mathbb{E}[Y^2_{b_s}] $ is similar to that of $\sum_{s=1}^4\mathbb{E}[X^2_{bs}] $ in Item (a).
\item The proof of $\mathbb{E}[Y_c^2]$ is similar to that of $\mathbb{E}[X_c^2]$ in Item (a). 
\end{itemize}
For these reasons, we omit the proof details and only state the claims. We have that under both the null and the alternative,
\beq\label{varF2a-add1}
\mathrm{Var}(Y_{a1})\leq C\|\theta\|^8\|\theta\|_3^3/\|\theta\|_1^5.
\eeq
\beq\label{varF2a-add3}
\mathrm{Var}(Y_{a2}) + \mathrm{Var}(Y_{a3})\leq  C\|\theta\|^4\|\theta\|_3^6/\|\theta\|_1^4,
\eeq
\beq\label{varF2a-add4}
\sum_{s=1}^3\mathbb{E}[Y^2_{b_s}]  \leq C\|\theta\|^4\|\theta\|_3^6/\|\theta\|_1^2,
\eeq
\beq\label{varF2a-add5}
\mathbb{E}[Y^2_{c}]  \leq C\|\theta\|^6\|\theta\|_3^6 / \|\theta\|_1^4. 
\eeq
Finally, inserting (\ref{varF2a-add1}), (\ref{varF2a-add3}), (\ref{varF2a-add4}), and (\ref{varF2a-add5}) into (\ref{varF2a-add}) gives   
\begin{align*} 
\mathrm{Var}(F_2) & \leq  C[\|\theta\|^8\|\theta\|_3^3/\|\theta\|_1^5 +  \|\theta\|^4\|\theta\|_3^6/\|\theta\|_1^4 + \|\theta\|^4\|\theta\|_3^6/\|\theta\|_1^2 + \|\theta\|^6\|\theta\|_3^6 / \|\theta\|_1^4] \\
& \leq  C\|\theta\|^4 \|\theta\|_3^6 /\|\theta\|_1^4,  \numberthis\label{F2-Vara} 
\end{align*} 
where we have used $\|\theta\|_3^3 \ll \|\theta\|^2 \ll \|\theta\|_1$, $\|\theta\|\goto\infty$ and $\|\theta\|^4 \leq \|\theta\|_1 \|\theta\|_3^3$.   This completes the proof of (\ref{F2varpf}).



 Consider Item (k). 
The goal is to show (\ref{F3meanpf}) and (\ref{F3varpf}). Recall that 
\[
F_{c} = \sum_{i_1,i_2,i_3,i_4 (dist)} \eta_{i_2}^2\eta_{i_4}^2\big[(\eta_{i_1} - \tilde{\eta}_{i_1})^2 (\eta_{i_3} - \tilde{\eta}_{i_3})^2     \big],
\]
and that $\widetilde{\eta} - \eta = v^{-1/2} W 1_n$. Plugging this into $F_3$ gives
\[
F_c = v^{-2}\sum_{i_1,i_2,i_3,i_4 (dist)} \sum_{\substack{j_1, \ell_1,  j_2, j_4 \\
j_1 \neq i_1, \ell_1 \neq i_1, j_3 \neq i_3, \ell_3 \neq i_3 }}  \eta_{i_2}^2\eta_{i_4}^2   
W_{i_1 j_1} W_{i_1 \ell_1} W_{i_3 j_3} W_{i_3 \ell_3}. 
\]
By basic combinatorics and basic algebra, we have 
\[
W_{i_1 j_1} W_{i_1 \ell_1} W_{i_3 j_3} W_{i_3 \ell_3}
= 
\left\{ 
\begin{array}{ll} 
W_{i_1i_3}^4, &\qquad \mbox{if $j_1 =\ell_1 = i_1, j_3 = \ell_3 = i_1$}, \\
W_{i_1i_3}^3 W_{i_1j_1}, &\qquad \mbox{if $j_3 =\ell_3 = i_1, \ell_1 = i_3$}, \\
W_{i_1i_3}^3 W_{i_1\ell_1}, &\qquad \mbox{if $j_3 =\ell_3 = i_1, j_1 = i_3$}, \\
W_{i_1i_3}^3 W_{i_3j_3}, &\qquad \mbox{if $j_1 =\ell_1 = i_3, \ell_3 = i_1$}, \\
W_{i_1i_3}^3 W_{i_3\ell_3}, &\qquad \mbox{if $j_1 =\ell_1 = i_3, j_3 = i_1$}, \\
W_{i_1j_1}^2 W_{i_3j_3}^2, &\qquad \mbox{if $j_1 =\ell_1, j_3 = \ell_3$}, \\
W_{i_1j_1}^2 W_{i_3j_3}W_{i_3 \ell_3}, &\qquad  \mbox{if $j_1 = \ell_1 \neq i_3, j_3\neq \ell_3$},   \\ 
W_{i_3j_3}^2 W_{i_1j_1}W_{i_1 \ell_1}, &\qquad  \mbox{if $j_3 = \ell_3 \neq i_1, j_1\neq \ell_1$},   \\ 
W_{i_1i_3}^2 W_{i_1\ell_1} W_{i_3\ell_3}, &\qquad  \mbox{if $j_1 = i_3, j_3 = i_1$},    \\
W_{i_1i_3}^2 W_{i_1j_1} W_{i_3j_3}, &\qquad  \mbox{if $\ell_1 = i_3, \ell_3 = i_1$},    \\
W_{i_1i_3}^2 W_{i_1j_1} W_{i_3\ell_3}, &\qquad  \mbox{if $\ell_1 = i_3, j_3 = i_1$},    \\
W_{i_1i_3}^2 W_{i_1\ell_1} W_{i_3j_3}, &\qquad  \mbox{if $j_1 = i_3, \ell_3 = i_1$},    \\
W_{i_1 j_1} W_{i_1 \ell_1} W_{i_3 j_3} W_{i_3 \ell_3}, &\qquad \mbox{otherwise}. \\
\end{array} 
\right. 
\] 
By these and symmetry, we can further split $F_3$ into $6$ different terms:  
\begin{equation} \label{F3-split} 
F_c = Z_{a} +  4 Z_{b1} + Z_{b2} +  2 Z_{c1} + 4Z_{c2} + Z_{d},
\end{equation} 
where
\[
Z_{a} =  v^{-2}\sum_{i_1, i_2, i_3, i_4(dist)}\eta_{i_2}^2  \eta_{i_4}^2 W_{i_1i_3}^4, 
\]
\[
Z_{b1} = v^{-2}\sum_{i_1, i_2, i_3, i_4(dist)}\sum_{j_4, j_4\neq i_4}\eta_{i_2}^2  \eta_{i_4}^2W_{i_1i_3}^3 W_{i_3 j_3}, 
\]
\[
Z_{b2} =  v^{-2}\sum_{i_1, i_2, i_3, i_4(dist)}\sum_{j_1, j_1\neq i_1, j_3, j_3\neq i_3}\eta_{i_2}^2  \eta_{i_4}^2 W_{i_1j_1}^2 W_{i_3j_3}^2, 
\]
\[
Z_{c1} =  v^{-2}\sum_{i_1, i_2, i_3, i_4(dist)}\sum_{\substack{j_1, j_3, \ell_3 \\ j_1\notin\{ i_1, i_3\}, j_3, \ell_3    }}\eta_{i_2}^2  \eta_{i_4}^2     W_{i_1j_1}^2 W_{i_3j_3}W_{i_3 \ell_3}, 
\]
\[
Z_{c2} =  v^{-2}\sum_{i_1, i_2, i_3, i_4(dist)}\sum_{\substack{\ell_1, \ell_3 \\ \ell_1\neq i_1, \ell_3\neq i_3}}\eta_{i_2}^2  \eta_{i_4}^2    W_{i_1i_3}^2 W_{i_1\ell_1} W_{i_3\ell_3} , 
\]
\[
Z_{d} =  v^{-2}\sum_{i_1, i_2, i_3, i_4(dist)}\sum_{\substack{j_1, \ell_1, j_3, \ell_3   \\ j_1\neq\ell_1, j_3\neq\ell_3 \\ j_1, \ell_1\neq i_3, j_3, \ell_3\neq i_1 }}\eta_{i_2}^2  \eta_{i_4}^2 W_{i_1 j_1} W_{i_1 \ell_1} W_{i_3 j_3} W_{i_3 \ell_3}. 
\]

We now show  (\ref{F3meanpf}) and (\ref{F3varpf}) separately.  Consider  (\ref{F3meanpf}) first. It is among all the $6$ $Z$-terms, only $Z_a$ and $Z_{b2}$  have non-zero means. 
We now consider $\mathbb{E}[Z_a]$ and $\mathbb{E}[Z_{b2}]$ separately. 

First, consider $\mathbb{E}[Z_a]$. By similar arguments, both under the null and the alternative, 
\[
\mathbb{E}[W_{i_1i_3}^4] \leq C \Omega_{i_1 i_3} \leq C \theta_{i_1} \theta_{i_3}. 
\] 
Recalling that $0 < \eta_i \leq C \theta_i$ and $v \asymp \|\theta\|^2$, it is seen that  
\begin{equation} \label{F3-mean-a} 
\mathbb{E}[Z_a] \leq C (\|\theta\|_1)^{-4} \sum_{i_1, i_2, i_3, i_4(dist)} \theta_{i_2}^2  \theta_{i_4}^2 \theta_{i_1} \theta_{i_3}  \leq C \|\theta\|^4/\|\theta\|_1^2. 
\end{equation} 
Next, consider $\mathbb{E}[Z_{b2}]$. First, recall that under the null,  $\Omega = \theta \theta'$, $v = 1_n' (\Omega - \diag(\Omega)) 1_n$, and $\eta =  v^{-1/2} (\Omega - \diag(\Omega) 1_n$.  It is seen $v \sim \|\theta\|_1^2$, $\eta_i  = (1 + o(1) \theta_i$, $1 \leq i \leq n$, where $o(1) \goto 0$ uniformly for all $1 \leq i \leq n$, and for any $i \neq j$, $\mathbb{E}[W_{ij}^2] = (1 + o(1)) \theta_i \theta_j$, 
where $o(1) \goto 0$ uniformly for all $1 \leq i, j \leq n$. 
It follows 
\begin{equation} \label{F3-mean-b2-null} 
\mathbb{E}[Z_{b2}] =   v^{-2}\sum_{i_1, i_2, i_3, i_4(dist)}\sum_{j_1, j_1\neq i_1, j_3, j_3\neq i_3}\eta_{i_2}^2  \eta_{i_4}^2 \mathbb{E}[W_{i_1j_1}^2 W_{i_3j_3}^2],
\end{equation} 
which  
\[
\sim (\|\theta\|_1)^{-4}\sum_{i_1, i_2, i_3, i_4(dist)}\sum_{j_1, j_1\neq i_1, j_3, j_3\neq i_3}\theta_{i_1}\theta_{i_2}^2  \theta_{i_3}\theta_{i_4}^2\theta_{j_1}\theta_{j_3}  \sim \|\theta\|^4. 
\]
Second, under the alternative, by similar argument, we have that $v \asymp \|\theta\|_1^2$, 
$0 < \eta_i < C \theta_i$ for all $1 \leq i \leq n$, and $\mathbb{E}[W_{ij}^2] \leq C\theta_i \theta_j$ for all $1 \leq i, j \leq n$, $i \neq j$. Similar to that under the null, we have  
\begin{equation} \label{F3-mean-b2-alt}   
0 < |\mathbb{E}[Z_{b2}] |  \leq C \|\theta\|^4. 
\end{equation} 
Inserting (\ref{F3-mean-a}), (\ref{F3-mean-b2-null}), and (\ref{F3-mean-b2-alt}) into (\ref{F3-split}) and recalling that the mean of all other $Z$ terms are $0$,  
\[
\mathbb{E}[F_3]  \sim \|\theta\|^4,  \qquad \mbox{under the null},  
\] 
and 
\[
\mathbb{E}[F_3] \leq C \|\theta\|^4, \qquad \mbox{under the alternative}, 
\]
where we have used $\|\theta\|_1 \goto \infty$. This proves (\ref{F3meanpf}).

We now consider (\ref{F3varpf}). By Cauchy-Schwarz inequality, 
\begin{align}  
\mathrm{Var}(F_c) & \leq C \big(\mathrm{Var}(Z_{a}) +  \mathrm{Var}(Z_{b1}) + \mathrm{Var}(Z_{b2})  + \mathrm{Var}(Z_{c1}) + \mathrm{Var}(Z_{c2})  + \mathrm{Var}(Z_{d})   \big) \nonumber \\ 
& \leq C \big(  \mathrm{Var}(Z_{a}) + \mathbb{E}[Z_{b1}^2] + \mathrm{Var}(Z_{b2}) + \mathbb{E}[Z^2_{c1}] + \mathbb{E}[Z^2_{c2}]   + \mathbb{E}[Z^2_{d}]  \big). \label{F3-var-split} 
\end{align} 
Consider $\mathrm{Var}(Z_a)$. Write 
\begin{align*} 
\mathrm{Var}(Z_{a}) & = v^{-4} \sum_{\substack{i_1, i_2, i_3, i_4 (dist) \\ i_1', i_2', i_3', i_4'   (dist) } } 
 \eta_{i_2}^2    \eta_{i_4}^2  \eta_{i_2'}^2 \eta_{i_4'}^2 \mathbb{E}[(W_{i_1i_3}^4 - \mathbb{E}[W_{i_1i_3}^4])  (W_{i_1'i_3'}^4 - \mathbb{E}[W_{i_1'i_3'}^4]) ]. 
\end{align*} 
Fix a term $(W_{i_1i_3}^4 - \mathbb{E}[W_{i_1i_3}^4])  (W_{i_1'i_3'}^4 - \mathbb{E}[W_{i_1'i_3'}^4])$. When the mean is nonzero, we must have $\{i_1, i_3\} = \{i_1', i_3'\}$, and when this happens, 
\[
\mathbb{E}[(W_{i_1i_3}^4 - \mathbb{E}[W_{i_1i_3}^4])  (W_{i_1'i_3'}^4 - \mathbb{E}[W_{i_1'i_3'}^4]) ] = \mathrm{Var}(W_{i_1i_3}^4). 
\]
For a random variable $X$, we have $\mathrm{Var}(X)\leq \mathbb{E}[X^2]$, and it follows that \[
\mathrm{Var}(W_{i_1i_3}^4) \leq \mathbb{E}[W_{i_1i_3}^8] \leq \mathbb{E}[W_{i_1i_3}^2], 
\]
where we have used the property that $0\leq W_{i_1i_3}^2\leq 1$; note that $\mathbb{E}[W_{i_1i_3}^2]\leq C\theta_{i_1}\theta_{i_3}$. Recall that $v \asymp \|\theta\|_1^2$ and $0 < \eta_i \leq C \theta_i$ for all $1 \leq i \leq n$. Combining these gives 
\begin{equation} \label{F3a-Var}
\mathrm{Var}(Z_{a}) \leq C  (\|\theta\|_1^{-8})  \cdot \sum_{\substack{i_1, i_2, i_3, i_4 (dist) \\ i_2', i_4'   (dist) } }
 \theta_{i_2}^2    \theta_{i_4}^2  \theta_{i_2'}^2 \theta_{i_4'}^2  \theta_{i_1} \theta_{i_3} \leq C \|\theta\|^8/\|\theta\|_1^6. 
\end{equation} 

We now consider all other terms on the right hand side of (\ref{F3-var-split}). Note that 
\begin{itemize}
\item The proof of $\mathbb{E}[Z_{b1}^2]$ is similar to that of $Y_{a1}$ in Item (b). 
\item The proof of $\mathrm{Var}(Z_{b2})$ is similar to that of $X_{a}$ in Item (a). 
\item The proof of $\mathbb{E}[Z_{c1}^2]$ and $\mathbb{E}[Z_{c2}^2]$ are similar to that of $X_{b}$ in Item (a). 
\item The proof of $\mathbb{E}[Z_d^2]$ is similar to that of $X_c$ in Item (a). 
\end{itemize} 
For these reasons, we skip the proof details. 
We have that, under both the null and the alternative, 
\beq\label{F3b1-Var}
\mathbb{E}[Z_{b1}^2]\leq C\|\theta\|^8\|\theta\|_3^3/\|\theta\|_1^5,
\eeq
\beq\label{F3b2-Var}
\mathrm{Var}(Z_{b2})\leq C\|\theta\|^8/\|\theta\|_1^2,
\eeq
\beq\label{F3c-Var}
\mathbb{E}[Z_{c1}^2] + \mathbb{E}[Z_{c2}^2] \leq C\|\theta\|^{10}/\|\theta\|_1^2,
\eeq
and
\beq\label{F3d-Var}
\mathbb{E}[Z_{d}^2]  \leq C\|\theta\|^{12}/\|\theta\|_1^4.
\eeq
Inserting \eqref{F3a-Var}, \eqref{F3b1-Var}, \eqref{F3b2-Var}, \eqref{F3c-Var} and \eqref{F3d-Var} into \eqref{F3-var-split} gives 
\begin{align*}
\mathrm{Var}(F_c) & \leq C[ \|\theta\|^8/\|\theta\|_1^6 + \|\theta\|^8/\|\theta\|_1^2 +\|\theta\|^{10}/\|\theta\|_1^2 + \|\theta\|^{12}/\|\theta\|_1^4 ] \\
& \leq C\|\theta\|^{10}/\|\theta\|_1^2,
\end{align*}
which completes the proof of (\ref{F3varpf}).

\subsubsection{Proof of Lemma~\ref{lem:event}}
Define an event $D$ as 
\[
D =\bigl\{  |V-v|\leq \|\theta\|_1\cdot x_n  \bigr\},  \qquad \mbox{for}\quad \sqrt{\log(\|\theta\|_1)}\ll x_n\ll \|\theta\|_1. 
\] 
We aim to show that
\beq \label{lem-event-0}
\mathbb{E}[(Q_n-Q_n^*)^2\cdot I_{D^c}] = o(\|\theta\|^8). 
\eeq

First, we bound the tail probability of $|V-v|$. Write 
\[
V-v = 2\sum_{i<j}(A_{ij}-\Omega_{ij}).
\]
The variables $\{A_{ij}-\Omega_{ij}\}_{1\leq i<j\leq n}$ are mutually independent with mean zero. They satisfy $|A_{ij}-\Omega_{ij}|\leq 1$ and $\sum_{i<j}\mathrm{Var}(A_{ij}-\Omega_{ij})\leq \sum_{i<j}\Omega_{ij}\leq 1_n'\Omega 1_n/2\leq \|\theta\|_1^2/2$. Applying the Bernstein's inequality, for any $t>0$, 
\[
\mathbb{P}\Bigl(\Bigl|2\sum_{i<j}(A_{ij}-\Omega_{ij}) \Bigr|>t \Bigr)\leq 2\exp\Bigl(  - \frac{t^2/2}{2\|\theta\|_1^2+t/3}\Bigr).
\]
We immediately have that, for some positive constants $C_1, C_2>0$, 
\beq \label{lem-event-tail}
\mathbb{P}(|V-v|>t) \leq \begin{cases}
2\exp\bigl( -\frac{C_1}{\|\theta\|_1^2} t^2\bigr), & \mbox{when } x_n\|\theta\|_1 \leq t \leq \|\theta\|_1^2,\\
2\exp\bigl( -C_2t \bigr), &\mbox{when }t > \|\theta\|_1^2. 
\end{cases} 
\eeq
Especially, letting $t=x_n\|\theta\|_1$, we have
\beq \label{lem-event-PD}
\mathbb{P}(D^c) \leq 2\exp(-C_1x_n^2). 
\eeq

Next, we derive an upper bound of $(Q_n-Q_n^*)^2$ in terms of $V$. Recall that $V$ is the total number of edges and that $Q_n=\sum_{i,j,k,\ell (dist)}M_{ij}M_{jk}M_{k\ell}M_{\ell i}$, where $M_{ij}=A_{ij}-\hat{\eta}_i\hat{\eta}_j$. If one node of $i,j,k,\ell$ has a zero degree (say, node $i$), then $A_{ij}=0$ and $\hat{\eta}_i=0$, and it follows that $M_{ij}=0$ and $M_{ij}M_{jk}M_{k\ell}M_{\ell i}=0$. Hence, only when $(i,j,k,\ell)$ all have nonzero degrees, this quadruple has a contribution to $Q_n$. Since $V$ is the total number of edges, there are at most $V$ nodes that have a nonzero degree. It follows that   
\[
|Q_n|\leq CV^4. 
\]  
Moreover, $Q_n^*=\sum_{i,j,k,\ell (dist)}M^*_{ij}M^*_{jk}M^*_{k\ell}M^*_{\ell i}$, where $M^*_{ij}=\widetilde{\Omega}_{ij}+W_{ij}+\delta_{ij}$. Re-write $M^*_{ij}=A_{ij}-\eta^*_i\eta_j^*+\eta_i(\eta_j-\teta_j)+\eta_j(\eta_j-\teta_j)$. First, since $\eta_i^*\leq C\theta_i$ and $\eta_i\leq C\theta_i$ (see \eqref{eta-bound}), $|M^*_{ij}|\leq A_{ij}+C\theta_i\theta_j+C\theta_i|\eta_j-\teta_j|+C\theta_j|\eta_i-\teta_i|$. Second, note that $\teta_i$ equals to $v^{-1/2}$ times degree of node $i$, where $v\asymp\|\theta\|^2_1$ according to \eqref{v-order}. It follows that $|\eta_i-\teta_i|\leq C(\theta_i +\|\theta\|^{-1}_1V)$. Therefore,
\[
|M^*_{ij}|\leq A_{ij} +C \theta_i\theta_j+ C\|\theta\|_1^{-1}V(\theta_i+\theta_j). 
\]
We plug it into the definition of $Q_n^*$ and note that there are at most $V$ pairs of $(i,j)$ such that $A_{ij}\neq 0$. By elementary calculation, 
\[
|Q_n^*|\leq C(V^4 + \|\theta\|_1^4).
\]
Combining the above gives
\beq \label{lem-event-naivebound}
(Q_n-Q_n^*)^2 \leq 2Q^2_n + 2(Q_n^*)^2 \leq C(V^8 + \|\theta\|_1^{8}). 
\eeq

Last, we show \eqref{lem-event-0}. By \eqref{lem-event-naivebound} and that $V^{8}\leq Cv^{8}+C|V-v|^{8}$, we have
\begin{align}  \label{lem-event-1}
\mathbb{E}[(Q_n-Q_n^*)^2\cdot I_{D^c}] & \leq C\mathbb{E}[|V-v|^8\cdot I_{D^c}] + C(v^8 + \|\theta\|_1^{8})\cdot \mathbb{P}(D^c)\cr
&\leq C\mathbb{E}[|V-v|^8\cdot I_{D^c}] + C \|\theta\|_1^{16} \cdot \mathbb{P}(D^c),
\end{align}
where the second line is from $v\asymp\|\theta\|_1^2$. Note that $x_n\gg\sqrt{\log(\|\theta\|_1)}$. For $n$ sufficiently large, 
$x^2_n\geq 17 C_1^{-1}\log(\|\theta\|_1)$. Combining it with \eqref{lem-event-PD}, we have
\beq \label{lem-event-2}
\|\theta\|_1^{16}\cdot\mathbb{P}(D^c)\leq \|\theta\|_1^{16}\cdot 2e^{-C_1x_n^2} \leq \|\theta\|_1^{16}\cdot 2e^{-17\|\theta\|_1}=o(1). 
\eeq
We then bound $\mathbb{E}[|V-v|^8\cdot I_{D^c}]$. Let $f(t)$ and $F(t)$ be the probability density and CDF of $|V-v|$, and write $\bar{F}(t)=1-F(t)$.  Using integration by part, for any continuously differentiable function $g(t)$ and $x>0$, $\int_{x}^{\infty}g(t)f(t)dt=g(x)\bar{F}(x) + \int_{x}^\infty g'(t)\bar{F}(t)dt$. We apply the formula to $g(t)=t^8$ and $x=x_n\|\theta\|_1$. It yields
\begin{align*}
\mathbb{E}[|V-v|^8\cdot I_{D^c}] &= (x_n\|\theta\|_1)^8\cdot \mathbb{P}(D^c) + \int_{x_n\|\theta\|_1}^\infty 8t^7\cdot \mathbb{P}(|V-v|>t)dt\cr
&\equiv I+II. 
\end{align*}
Consider $I$. By \eqref{lem-event-2} and $x_n\ll\|\theta\|_1$, 
\[
I \ll\|\theta\|_1^{16}\cdot \mathbb{P}(D^c) = o(1). 
\]
Consider $II$. By \eqref{lem-event-tail}, \eqref{lem-event-2}, and elementary probability, 
\begin{align*}
II &\leq  8(\|\theta\|_1^2)^7\cdot \mathbb{P}\bigl(x_n\|\theta\|_1<|V-v|\leq \|\theta\|_1^2\bigr) + \int_{\|\theta\|^2_1} 8t^7 \cdot \mathbb{P}(|V-v|> t)dt\cr
&\leq  C\|\theta\|_1^{14}\cdot\mathbb{P}(D^c) + \int_{\|\theta\|^2_1} 8t^7 \cdot 2e^{-C_2t}dt\cr
&=o(1),
\end{align*}
where in the last line we have used \eqref{lem-event-2} and the fact that $\int_{x}^\infty t^7e^{-C_2t}dt \to 0$ as $x\to\infty$. Combining the bounds for $I$ and $II$ gives
\beq \label{lem-event-3}
\mathbb{E}[|V-v|^8\cdot I_{D^c}]  = o(1). 
\eeq
Then, \eqref{lem-event-0} follows by plugging \eqref{lem-event-2}-\eqref{lem-event-3} into \eqref{lem-event-1}.

\subsubsection{Proof of Lemma~\ref{lem:remainder1}} \label{subsec:remainder1}
There are $175$ post-expansion sums in $(\widetilde{Q}^*_n-Q_n^*)$. They divide into $34$ different types, denoted by $R_1$-$R_{34}$ as shown in Table~\ref{tb:remainder}. 
It suffices to prove that, for each $1\leq k\leq 34$, under the null hypothesis, 
\beq \label{remainder1-Nullgoal}
\bigl| \mathbb{E}[R_k]\bigr| = o(\|\theta\|^4), \qquad \mathrm{Var}(R_k)=o(\|\theta\|^8), 
\eeq
and under the alternative hypothesis, 
\beq \label{remainder1-Altgoal}
\bigl|\mathbb{E}[R_k]\bigr|=o(\alpha^4\|\theta\|^8), \qquad \mathrm{Var}(R_k)=O(\|\theta\|^8+\alpha^6\|\theta\|^8\|\theta\|_3^6). 
\eeq 
\begin{table}[tb!]  
\centering
\caption{The $34$ types of the $175$ post-expansion sums for $(\widetilde{Q}^*_n-Q_n^*)$.}    \label{tb:remainder}
\scalebox{.95}{
\begin{tabular}{lccclcc}
Notation  & $\#$ & $N_{\tilde{r}}$ &  ($N_{\delta}, N_{\widetilde{\Omega}}, N_W)$ &     Examples & $N^*_W$\\
\hline 
$R_1$  & 4  &1&(0, 0, 3)   & $\sum_{i, j, k,\ell (dist)} \tilde{r}_{ij}W_{jk} W_{k\ell} W_{\ell i}$ & 5\\  
$R_2$  & 8  &1&(0, 1, 2)   & $\sum_{i, j, k,\ell (dist)} \tilde{r}_{ij} \widetilde{\Omega}_{jk} W_{k\ell} W_{\ell i}$ & 4\\ 
$R_3$  & 4  & &     & $\sum_{i, j, k,\ell (dist)} \tilde{r}_{ij} W_{jk}\widetilde{\Omega}_{k\ell} W_{\ell i}$ &  4\\ 
$R_4$ & 8 &1&(0, 2, 1)   & $\sum_{i, j, k,\ell (dist)} \tilde{r}_{ij} \widetilde{\Omega}_{jk} \widetilde{\Omega}_{k\ell} W_{\ell i}$ & 3\\  
$R_5$ & 4 &&   & $\sum_{i, j, k,\ell (dist)} \tilde{r}_{ij} \widetilde{\Omega}_{jk}W_{k\ell} \widetilde{\Omega}_{\ell i}$ & 3\\   
$R_6$ & 4 & 1& (0, 3, 0)   & $\sum_{i, j, k,\ell (dist)} \tilde{r}_{ij}\widetilde{\Omega}_{jk} \widetilde{\Omega}_{k\ell} \widetilde{\Omega}_{\ell i}$  & 2\\   
$R_7$   & 8  & 1 & (1, 0, 2)   & $\sum_{i, j, k,\ell (dist)} \tilde{r}_{ij} \delta_{jk} W_{k\ell} W_{\ell i}$ & 5\\  
$R_8$   & 4  &  &   & $\sum_{i, j, k,\ell (dist)} \tilde{r}_{ij} W_{jk}\delta_{k\ell} W_{\ell i}$ & 5\\ 
$R_9$  & 8  &1 &(1, 1, 1)   & $\sum_{i, j, k,\ell (dist)} \tilde{r}_{ij}\delta_{jk}\widetilde{\Omega}_{k\ell}  W_{\ell i}$ & 4   \\ 
$R_{10}$  & 8  & &   & $\sum_{i, j, k,\ell (dist)} \tilde{r}_{ij}\widetilde{\Omega}_{jk}  W_{k \ell}\delta_{\ell i}$ &   
4  \\ 
$R_{11}$ & 8  &&   & $\sum_{i, j, k,\ell (dist)} \tilde{r}_{ij}W_{jk}\delta_{k\ell }\widetilde{\Omega}_{\ell i}  $ &   
4 \\ 
$R_{12}$ & 8 &1& (1, 2, 0)   & $\sum_{i, j, k,\ell (dist)}\tilde{r}_{ij} \delta_{jk}  \widetilde{\Omega}_{k\ell} \widetilde{\Omega}_{\ell i}$ & 3\\ 
$R_{13}$ & 4 & & & $\sum_{i, j, k,\ell (dist)}\tilde{r}_{ij} \widetilde{\Omega}_{jk} \delta_{k\ell}  \widetilde{\Omega}_{\ell i}$ & 3\\
$R_{14}$ & 8 & 1 &(2, 0, 1)   & $\sum_{i, j, k,\ell (dist)} \tilde{r}_{ij}\delta_{jk} \delta_{k\ell} W_{\ell i}$    &  5\\  
$R_{15}$ & 4 &  &   & $\sum_{i, j, k,\ell (dist)} \tilde{r}_{ij}\delta_{jk} W_{k\ell }\delta_{\ell i} $    &  5\\ 
$R_{16}$ & 8 & 1 & (2, 1, 0) & $\sum_{i, j, k,\ell (dist)} \tilde{r}_{ij}\delta_{jk} \delta_{k\ell} \widetilde{\Omega}_{\ell i}$    &  4\\  
$R_{17}$  & 4 &  &  & $\sum_{i, j, k,\ell (dist)} \tilde{r}_{ij}\delta_{jk}\widetilde{\Omega}_{k\ell} \delta_{\ell i}$    & 4 \\  
$R_{18}$ & 4 &  1 & (3,  0,  0)  & $\sum_{i, j, k,\ell (dist)} \widetilde{r}_{ij}\delta_{jk} \delta_{k\ell} \delta_{\ell i}$ & 5\\ 
\hline 
$R_{19}$   & 4  & 2 & (0, 0, 2)   & $\sum_{i, j, k,\ell (dist)} \tilde{r}_{ij} \tilde{r}_{jk} W_{k\ell} W_{\ell i}$ &6 \\ 
$R_{20}$      & 2  &  &    & $\sum_{i, j, k,\ell (dist)} \tilde{r}_{ij} W_{jk}\tilde{r}_{k\ell} W_{\ell i}$ & 6\\  
$R_{21}$   & 4 & 2 & (0, 2, 0)   & $\sum_{i, j, k,\ell (dist)} \tilde{r}_{ij} \tilde{r}_{jk} \widetilde{\Omega}_{k\ell} \widetilde{\Omega}_{\ell i}$ & 4\\ 
 $R_{22}$     & 2  &  &    & $\sum_{i, j, k,\ell (dist)} \tilde{r}_{ij} \widetilde{\Omega}_{jk}\tilde{r}_{k\ell} \widetilde{\Omega}_{\ell i}$ & 4\\
 $R_{23}$     & 4  & 2 & (2, 0, 0)   & $\sum_{i, j, k,\ell (dist)} \tilde{r}_{ij} \tilde{r}_{jk} \delta_{k\ell} \delta_{\ell i}$ & 6\\ 
 $R_{24}$      & 2  &  &    & $\sum_{i, j, k,\ell (dist)} \tilde{r}_{ij} \delta_{jk}\tilde{r}_{k\ell} \delta_{\ell i}$ & 6\\    
 $R_{25}$  & 8  & 2 & (0, 1, 1)   & $\sum_{i, j, k,\ell (dist)} \tilde{r}_{ij} \tilde{r}_{jk} \widetilde{\Omega}_{k\ell} W_{\ell i}$ & 5\\ 
 $R_{26}$      & 4  &  &    & $\sum_{i, j, k,\ell (dist)} \tilde{r}_{ij} \widetilde{\Omega}_{jk}\tilde{r}_{k\ell} W_{\ell i}$ & 5\\  
 $R_{27}$   & 8  & 2 & (1, 1, 0)   & $\sum_{i, j, k,\ell (dist)} \tilde{r}_{ij} \tilde{r}_{jk} \delta_{k\ell} \widetilde{\Omega}_{\ell i}$ & 5\\ 
  $R_{28}$     & 4  &  &    & $\sum_{i, j, k,\ell (dist)} \tilde{r}_{ij} \delta_{jk}\tilde{r}_{k\ell} \widetilde{\Omega}_{\ell i}$ & 5\\
 $R_{29}$      & 8  & 2 & (1, 0, 1)   & $\sum_{i, j, k,\ell (dist)} \tilde{r}_{ij} \tilde{r}_{jk} \delta_{k\ell} W_{\ell i}$ & 6\\ 
 $R_{30}$      & 4  &  &    & $\sum_{i, j, k,\ell (dist)} \tilde{r}_{ij} \delta_{jk}\tilde{r}_{k\ell} W_{\ell i}$ & 6\\  
\hline 
 $R_{31}$   & 4  & 3 & (0, 0, 1)   & $\sum_{i, j, k,\ell (dist)} \tilde{r}_{ij} \tilde{r}_{jk} \tilde{r}_{k\ell} W_{\ell i}$ & 7\\ 
 $R_{32}$   & 4  & 3 & (0, 1, 0)   & $\sum_{i, j, k,\ell (dist)} \tilde{r}_{ij} \tilde{r}_{jk} \tilde{r}_{k\ell} \widetilde{\Omega}_{\ell i}$ & 6\\  
 $R_{33}$  & 4  & 3 & (1, 0, 0)   & $\sum_{i, j, k,\ell (dist)} \tilde{r}_{ij} \tilde{r}_{jk} \tilde{r}_{k\ell} \delta_{\ell i}$ & 7\\   
\hline
 $R_{34}$     & $1$  & 4 & (0, 0, 0)   & $\sum_{i, j, k,\ell (dist)} \tilde{r}_{ij} \tilde{r}_{jk} \tilde{r}_{k\ell} \tilde{r}_{\ell i}$ & 8\\   
\hline
\end{tabular} 
} 
\end{table} 

We need some preparation. First, recall that $\tilde{r}_{ij}= -\frac{v}{V}(\teta_i-\eta_i)(\teta_j-\eta_j)$. It follows that each post-expansion sum has the form
\beq \label{compareXY-0}
\Bigl(\frac{v}{V}\Bigr)^{N_{\tilde{r}}}\sum_{i,j,k,\ell (dist)}a_{ij}b_{jk}c_{k\ell}d_{\ell i},
\eeq
where $a_{ij}$ takes values in $\{\widetilde{\Omega}_{ij}, W_{ij}, \delta_{ij}, -(\teta_i-\eta_i)(\teta_j-\eta_j)\}$ and $b_{jk}, c_{k\ell}, d_{\ell i}$ are similar. The variable $\frac{v}{V}$ has a complicated correlation with each summand, so we want to get rid of it. 
Denote the variable in \eqref{compareXY-0} by $Y$. Write $m=N_{\tilde{r}}$ and 
\beq \label{compareXY-1}
Y =  \Bigl(\frac{v}{V}\Bigr)^m X, \quad\quad \mbox{where}\quad X =  \sum_{i,j,k,\ell (dist)}a_{ij}b_{jk}c_{k\ell}d_{\ell i}. 
\eeq
We compare the mean and variance of $X$ and $Y$. By assumption, $\sqrt{\log(\|\theta\|_1)}\ll \|\theta\|_1/\|\theta\|^2$. Then, there exists a sequence $x_n$ such that
\[
\sqrt{\log(\|\theta\|_1)}\ll x_n\ll \|\theta\|_1/\|\theta\|^2, \qquad \mbox{as}\;\; n\to\infty.
\] 
We introduce an event 
\[
D = \bigl\{|V-v|\leq \|\theta\|_1 x_n  \bigr\}. 
\]
In Lemma~\ref{lem:event}, we have proved $\mathbb{E}[(Q_n-Q_n^*)^2\cdot I_{D^c}] = o(1)$. By similar proof, we can show: as long as $|Y-X|$ is bounded by a polynomial of $V$ and $\|\theta\|_1$, 
\beq \label{compareXY-2}
\mathbb{E}[(Y-X)^2\cdot I_{D^c}]=o(1). 
\eeq
Additionally, on the event $D$, since $v\asymp\|\theta\|_1^2\gg \|\theta\|_1 x_n$, we have $|V-v| = o(v)$. It follows that $\frac{|V-v|}{V}\lesssim \frac{|V-v|}{v}\leq C\|\theta\|^{-1}x_n=o(1)$. For any fixed $m\geq 1$, $(1+x)^m \leq 1+Cx$ for $x$ being close to $0$. Hence, $|1-\frac{v^m}{V^m}|\leq C|1-\frac{v}{V}|\leq C\|\theta\|_1^{-1}x_n=o(\|\theta\|^{-2})$. It implies
\beq \label{compareXY-3}
|Y-X|=o(\|\theta\|^{-2})\cdot |X|, \qquad\mbox{on the event }D. 
\eeq
By \eqref{compareXY-2}-\eqref{compareXY-3} and elementary probability, 
\begin{align*}
|\mathbb{E}[Y-X]| &\leq  |\mathbb{E}[(Y-X)\cdot I_D]| + |\mathbb{E}[(Y-X)\cdot I_{D^c}]|\cr
&\leq o(\|\theta\|^{-2})\cdot \mathbb{E}[|X|\cdot I_D] + \sqrt{\mathbb{E}[(Y-X)^2\cdot I_{D^c}]}\cr
&\leq o(\|\theta\|^{-2}) \sqrt{\mathbb{E}[X^2]} + o(1),
\end{align*}
and
\begin{align*}
\mathrm{Var}(Y) & \leq 2\mathrm{Var}(X) + 2\mathrm{Var}(Y-X)\cr
&\leq 2\mathrm{Var}(X) + 2\mathbb{E}[(Y-X)^2] \cr
&= 2\mathrm{Var}(X) + 2\mathbb{E}[(Y-X)^2\cdot I_D] +2\mathbb{E}[(Y-X)^2\cdot I_{D^c}]\cr
&\leq 2\mathrm{Var}(X) + o(\|\theta\|^{-4})\cdot \mathbb{E}[X^2] +o(1).
\end{align*}
Under the null hypothesis, suppose we can prove that 
\beq \label{remainder1-Nullgoal2}
\mathbb{E}[X^2] = o(\|\theta\|^8). 
\eeq
Since $\mathbb{E}[X^2]=(\mathbb{E}[X])^2+\mathrm{Var}(X)$, it implies $|\mathbb{E}[X]|=o(\|\theta\|^4)$ and $\mathrm{Var}(X)=o(\|\theta\|^8)$. 
Therefore, 
\begin{align*}
& |\mathbb{E}[Y]|\leq |\mathbb{E}[X]| + |\mathbb{E}[Y-X]|= o(\|\theta\|^4), \cr
&\mathrm{Var}(Y)\leq C\mathrm{Var}(X) + o(\|\theta\|^{-4})\cdot \mathbb{E}[X^2] +o(1) =o(\|\theta\|^8).
\end{align*}
Under the alternative hypothesis, suppose we can prove that 
\beq \label{remainder1-Altgoal2}
|\mathbb{E}[X]|=O(\alpha^2\|\theta\|^6), \qquad \mathrm{Var}(X)=o(\|\theta\|^8+\alpha^6\|\theta\|^8\|\theta\|^6_3).
\eeq
Since $\mathbb{E}[X^2]=(\mathbb{E}[X])^2+\mathrm{Var}(X)$, we have $\mathbb{E}[X^2]=O(\alpha^4\|\theta\|^{12})$. 
Then, 
\begin{align*}
& |\mathbb{E}[Y]|\leq O(\alpha^2\|\theta\|^6) + o(\|\theta\|^{-2})\cdot O(\alpha^2\|\theta\|^6)=o(\alpha^4\|\theta\|^8), \cr
&\mathrm{Var}(Y)\leq o(\|\theta\|^8+\alpha^6\|\theta\|^8\|\theta\|^6_3) + o(\|\theta\|^{-4})\cdot O(\alpha^4\|\theta\|^{12}) =o(\|\theta\|^8+\alpha^6\|\theta\|^8\|\theta\|^6_3).
\end{align*}
In conclusion, to prove that $Y$ satisfies the requirement in \eqref{remainder1-Nullgoal}-\eqref{remainder1-Altgoal}, it is sufficient to prove that $X$ satisfies \eqref{remainder1-Nullgoal2}-\eqref{remainder1-Altgoal2}. We remark that \eqref{remainder1-Altgoal2} puts a more stringent requirement on the mean of the variable, compared to \eqref{remainder1-Altgoal}.  

From now on, in the analysis of each $R_k$ of the form \eqref{compareXY-0}, we shall always neglect the factor $(\frac{v}{V})^{N_{\tilde{r}}}$, and show that, after this factor is removed, the random variable satisfies \eqref{remainder1-Nullgoal2}-\eqref{remainder1-Altgoal2}. This is equivalent to pretending
\[
\tilde{r}_{ij} = - (\teta_i-\eta_i)(\teta_j-\eta_j)
\]
and proving each $R_k$ satisfies \eqref{remainder1-Nullgoal2}-\eqref{remainder1-Altgoal2}. Unless mentioned, we stick to this mis-use of notation $\tilde{r}_{ij}$ in the proof below.

Second, we divide 34 terms into several groups using the {\it intrinsic order of $W$} defined below. Note that $\tilde{r}_{ij}=-(\teta_i-\eta_i)(\teta_j-\eta_j)$, $\delta_{ij}=\eta_i(\eta_j-\teta_j)+\eta_j(\eta_i-\teta_i)$, and $\teta_i-\eta_i=\frac{1}{\sqrt{v}}\sum_{s\neq i}W_{is}$. We thus have
\[
\tilde{r}_{ij} = -\frac{1}{v}\Bigl(\sum_{s\neq i}W_{is}\Bigr)\Bigl(\sum_{t\neq j}W_{jt}\Bigr), \qquad \delta_{ij}= -\frac{1}{\sqrt{v}}\eta_i\Bigl(\sum_{t\neq j}W_{jt}\Bigr)- \frac{1}{\sqrt{v}}\eta_j\Bigl(\sum_{s\neq i}W_{is}\Bigr). 
\] 
Each $\tilde{r}_{ij}$ is a weighted sum of terms like $W_{is}W_{jt}$, and each $\delta_{ij}$ is a weighted sum of terms like $W_{jt}$. Intuitively, we view $\tilde{r}$-term as an ``order-2 $W$-term" and view $\delta$-term as ``order-1 $W$-term." It motivates the definition of {\it intrinsic order of $W$} as
\beq \label{intrinsic-order}
N^*_W = N_W + N_\delta + 2N_{\tilde{r}}. 
\eeq 
We group 34 terms by the value of $N^*_W$; see the last column of Table~\ref{tb:remainder}.

\paragraph{Analysis of post-expansion sums with $N^*_W\leq 4$} There are $14$ such terms, including $R_2$-$R_6$, $R_9$-$R_{13}$, $R_{16}$-$R_{17}$, and $R_{21}$-$R_{22}$. They all equal to zero under the null hypothesis, so it is sufficient to show that they satisfy \eqref{remainder1-Altgoal2} under the alternative hypothesis. We prove by comparing  each $R_k$ to some previously analyzed terms. Take $R_9$ for example. Plugging in the definition of $\tilde{r}_{ij}$ and $\delta_{ij}$ gives 
\begin{align*} 
R_9 &=\sum_{i,j,k,\ell (dist)}[(\teta_i-\eta_i)(\teta_j-\eta_j)][(\teta_j-\eta_j)\eta_k+\eta_j(\teta_k-\eta_k)]\widetilde{\Omega}_{k\ell}W_{\ell i}\cr
&= R_{9a}+R_{9b}, 
\end{align*}
where
\begin{align} \label{lem-remainder1-example1}
R_{9a} &= \sum_{i,j,k,\ell (dist)} \eta_k\widetilde{\Omega}_{k\ell}\cdot [ (\teta_i-\eta_i)(\teta_j-\eta_j)^2W_{\ell i}],\cr
R_{9b} &= \sum_{i,j,k,\ell (dist)} \eta_j\widetilde{\Omega}_{k\ell}\cdot[(\teta_i-\eta_i)(\teta_j-\eta_j)(\teta_k-\eta_k)W_{\ell i}].  
\end{align}
At the same time, we recall that $T_1$ in Lemmas~\ref{lem:ProxySgnQ(c)-null}-\ref{lem:ProxySgnQ(c)-alt} is defined as 
\[
T_1 = \sum_{i, j, k, \ell (dist)}\delta_{ij}\delta_{jk}\delta_{k\ell}W_{\ell i}=\sum_{i, j, k, \ell (dist)}\delta_{\ell j}\delta_{jk}\delta_{ki}W_{i\ell}. 
\]
In the proof of the above two lemmas, we express $T_1$ as the weighted sum of $T_{1a}$-$T_{1d}$; see \eqref{T1-split}. Note that $T_{1a}$ and $T_{1d}$ in \eqref{T1-split} can be re-written as
\begin{align} \label{lem-remainder1-example2}
T_{1d} &= \sum_{ i, j, k, \ell (dist)}[\eta_\ell(\teta_j-\eta_j)][(\teta_j-\eta_j)\eta_k][\eta_k(\teta_i-\eta_i)\bigr]W_{i\ell}\cr
&= \sum_{i,j,k,\ell (dist)} \eta^2_k\eta_{\ell}\cdot [ (\teta_i-\eta_i)(\teta_j-\eta_j)^2W_{\ell i}],\cr 
T_{1a} &= \sum_{i, j, k, \ell (dist)} [\eta_\ell(\teta_j-\eta_j)][\eta_j(\teta_k-\eta_k)][\eta_k(\teta_i-\eta_i)\bigr]W_{i\ell}\cr
&= \sum_{i, j, k, \ell (dist)} \eta_j\eta_k\eta_\ell \cdot [(\teta_i-\eta_i)(\teta_j-\eta_j)(\teta_k-\eta_k) W_{i\ell}]. 
\end{align}
Compare \eqref{lem-remainder1-example1} and \eqref{lem-remainder1-example2}. It is seen that $R_{9a}$ and $T_{1d}$ have the same structure, where the non-stochastic coefficients in the summand satisfy $|\eta_k\widetilde{\Omega}_{k\ell}|\leq C\alpha\theta_k^2\theta_\ell$ and $|\eta^2_k\eta_\ell|\leq C\theta_k^2\theta_\ell$, respectively. This means we can bound $|\mathbb{E}(R_{9a})|$ and $\mathrm{Var}(R_{9a})$ in the same way as we bound $|\mathbb{E}[T_{1d}]|$ and $\mathrm{Var}(T_{1d})$, and the bounds have an extra factor of $\alpha$ and $\alpha^2$, respectively. In detail, in the proof of Lemmas~\ref{lem:ProxySgnQ(c)-null}-\ref{lem:ProxySgnQ(c)-alt}, we have shown 
\[
|\mathbb{E}[T_{1d}]|\leq C\|\theta\|^4, \qquad \mathrm{Var}(T_{1d})\leq \frac{C\|\theta\|^6\|\theta\|^3_3}{\|\theta\|_1}.
\]
It follows immediately that 
\[
|\mathbb{E}[R_{9a}]|\leq C\alpha\|\theta\|^4 = o(\alpha^2\|\theta\|^6), \qquad \mathrm{Var}(T_{1d})\leq \frac{C\alpha^2\|\theta\|^6\|\theta\|^3_3}{\|\theta\|_1}=o(\|\theta\|^8).
\]
Similarly, since we have proved
\[
|\mathbb{E}[T_{1a}]|\leq \frac{C\|\theta\|^6}{\|\theta\|_1^2}, \qquad \mathrm{Var}(T_{1a})\leq \frac{C\|\theta\|^4\|\theta\|^6_3}{\|\theta\|^2_1},
\]
it follows immediately that
\[
|\mathbb{E}[R_{9b}]|\leq \frac{C\alpha\|\theta\|^6}{\|\theta\|_1^2} = o(\alpha^2\|\theta\|^6), \qquad \mathrm{Var}(R_{9b})\leq \frac{C\alpha^2 \|\theta\|^4\|\theta\|^6_3}{\|\theta\|^2_1}=o(\|\theta\|^8).
\]
This proves \eqref{remainder1-Altgoal2} for $X=R_{9a}$. 

We use the same strategy to bound all other terms with $N^*_W\leq 4$. The details are in Table~\ref{tb:Order-4}.
In each row of the table, the left column displays a targeting variable $X$, and the right column displays a previously analyzed variable, which we call $X^*$, that has a similar structure as $X$. It is not hard to see that we can obtain upper bounds for $|\mathbb{E}[X]|$ and $\mathrm{Var}(X)$ from multiplying the upper bounds of $|\mathbb{E}[X^*]|$ and $\mathrm{Var}(X^*)$ by $\alpha^m$ and $\alpha^{2m}$, respectively, where $m$ is a nonnegative integer (e.g., $m=1$ in the analysis of $R_9$). Using our previous results, each $X^*$ in the right column satisfies
\[
|\mathbb{E}[X^*]|=O(\alpha^2\|\theta\|^6), \qquad \mathrm{Var}(X^*)=o(\|\theta\|^8+\alpha^6\|\theta\|^8\|\theta\|_3^6). 
\]
So, each $X$ in the left column satisfies \eqref{remainder1-Altgoal2}.


\begin{table}[tb!]  
\centering
\caption{The $14$ types of post-expansion sums with $N^*_W\leq 4$. The right column displays the post-expansion sums defined before which have similar forms as the post-expansion sums in the left column. Definitions of the terms in the right column can be found in \eqref{proof-Z1-decompose}, \eqref{proof-Z2-decompose}, \eqref{proof-Z3-decompose}, \eqref{proof-Z4-decompose}, \eqref{proof-Z5-decompose}, \eqref{T1-split}, \eqref{T2-split}, and \eqref{F-split}. For some terms in the right column, we permute $(i,j,k,\ell)$ in the original definition for ease of comparison with the left column. (In all expressions, the subscript ``$i,j,k,\ell (dist)$" is omitted.)}    \label{tb:Order-4}
\scalebox{.9}{
\begin{tabular}{lc | l c}
  &     Expression  &  & Expression\\
\hline   
$R_2$     & $\sum (\teta_i-\eta_i)(\teta_j-\eta_j) \widetilde{\Omega}_{jk} W_{k\ell} W_{\ell i}$ & $Z_{1b}$ &  $\sum (\teta_i-\eta_i)\eta_j (\teta_j-\eta_j)\eta_k W_{k\ell} W_{\ell i}$\\ 
$R_3$  & $\sum (\teta_i-\eta_i)(\teta_j-\eta_j) W_{jk}\widetilde{\Omega}_{k\ell} W_{\ell i}$ &  $Z_{2a}$ & $\sum \eta_\ell (\teta_j-\eta_j)W_{jk}\eta_k(\teta_i-\eta_i)W_{i\ell}$\\ 
$R_4$ &  $\sum (\teta_i-\eta_i)(\teta_j-\eta_j) \widetilde{\Omega}_{jk} \widetilde{\Omega}_{k\ell} W_{\ell i}$ & $Z_{3d}$ & $\sum (\teta_i-\eta_i)\eta_j(\teta_j-\eta_j)\eta_k\widetilde{\Omega}_{k\ell}W_{\ell  i}$\\  
$R_5$ &  $\sum(\teta_i-\eta_i)(\teta_j-\eta_j) \widetilde{\Omega}_{jk}W_{k\ell} \widetilde{\Omega}_{\ell i}$ &$Z_{4b}$ & $\sum \widetilde{\Omega}_{ij}(\teta_j-\eta_j)\eta_kW_{k\ell}\eta_\ell(\teta_i-\eta_i)$ \\   
$R_6$ &  $\sum (\teta_i-\eta_i)(\teta_j-\eta_j)\widetilde{\Omega}_{jk} \widetilde{\Omega}_{k\ell} \widetilde{\Omega}_{\ell i}$  & $Z_{5a}$ & $\sum \eta_i(\teta_j-\eta_j)\widetilde{\Omega}_{jk}\widetilde{\Omega}_{k\ell}\eta_\ell (\teta_i-\eta_i)$ \\   
$R_9$  & $\sum(\teta_i-\eta_i)(\teta_j-\eta_j)^2\eta_k \widetilde{\Omega}_{k\ell}  W_{\ell i}$ & $T_{1d}$ &   
 $\sum \eta_{\ell}(\teta_j-\eta_j)^2\eta^2_k(\teta_i-\eta_i)W_{i\ell}$\\ 
  & $\sum(\teta_i-\eta_i)(\teta_j-\eta_j)\eta_j (\teta_k-\eta_k)\widetilde{\Omega}_{k\ell}  W_{\ell i}$ & $T_{1a}$ &   
 $\sum \eta_{\ell}(\teta_j-\eta_j)\eta_j(\teta_k-\eta_k)\eta_k(\teta_i-\eta_i)W_{i\ell}$\\ 
$R_{10}$  & $\sum (\teta_i-\eta_i)^2(\teta_j-\eta_j) \widetilde{\Omega}_{jk}  W_{k \ell}\eta_\ell$ &   
$T_{1c}$ & $\sum (\teta_j-\eta_j)\eta_kW_{k\ell}\eta_\ell (\teta_i-\eta_i)^2\eta_j$ \\ 
  & $\sum (\teta_i-\eta_i)(\teta_j-\eta_j) \widetilde{\Omega}_{jk}  W_{k \ell}(\teta_\ell -\eta_\ell)\eta_i$ & $T_{1a}$ & $\sum (\teta_j-\eta_j)\eta_kW_{k\ell}(\teta_\ell-\eta_\ell)\eta_i (\teta_i-\eta_i)\eta_j$  
 \\ 
$R_{11}$  & $\sum (\teta_i-\eta_i)(\teta_j-\eta_j) W_{jk}\eta_k (\teta_\ell-\eta_\ell)\widetilde{\Omega}_{\ell i}$ & $T_{1a}$  & $\sum (\teta_i-\eta_i)\eta_k W_{kj} (\teta_j-\eta_j)\eta_\ell (\teta_\ell-\eta_\ell)\eta_i$   
 \\ 
  & $\sum (\teta_i-\eta_i)(\teta_j-\eta_j)W_{jk}(\teta_k-\eta_k)\eta_\ell \widetilde{\Omega}_{\ell i}  $ & $T_{1b}$ & $\sum \eta_i (\teta_k-\eta_k) W_{kj} (\teta_j-\eta_j)\eta^2_\ell (\teta_i-\eta_i)$   
 \\ 
$R_{12}$ &  $\sum (\teta_i-\eta_i)(\teta_j-\eta_j)^2\eta_k  \widetilde{\Omega}_{k\ell} \widetilde{\Omega}_{\ell i}$ & $T_{2c}$ & $\sum \eta_i(\teta_j-\eta_j)^2\eta_k\widetilde{\Omega}_{k\ell}\eta_\ell(\teta_i-\eta_i)$\\ 
 &  $\sum (\teta_i-\eta_i)(\teta_j-\eta_j)\eta_j(\teta_k-\eta_k) \widetilde{\Omega}_{k\ell} \widetilde{\Omega}_{\ell i}$ & $T_{2a}$ & $\sum \eta_i(\teta_j-\eta_j)\eta_j(\teta_k-\eta_k)\widetilde{\Omega}_{k\ell}\eta_\ell(\teta_i-\eta_i)$\\ 
$R_{13}$ & $\sum (\teta_i-\eta_i)(\teta_j-\eta_j) \widetilde{\Omega}_{jk} (\teta_k-\eta_k)\eta_\ell  \widetilde{\Omega}_{\ell i}$ & $T_{2b}$ & $\sum \eta_i(\teta_j-\eta_j)\widetilde{\Omega}_{jk}(\teta_k-\eta_k)\eta_{\ell}^2(\teta_i-\eta_i)$\\
$R_{16}$ &  $\sum (\teta_i-\eta_i)(\teta_j-\eta_j)^2\eta_k(\teta_k-\eta_k)\eta_{\ell} \widetilde{\Omega}_{\ell i}$    &  $F_b$ &  $\sum \eta_i (\teta_j-\eta_j)^2\eta_k (\teta_k-\eta_k)\eta_\ell^2(\teta_i-\eta_i)$\\    
 &  $\sum (\teta_i-\eta_i)(\teta_j-\eta_j)^2\eta^2_k(\teta_\ell-\eta_\ell)\widetilde{\Omega}_{\ell i}$ & $F_b$ &  $\sum \eta_i (\teta_j-\eta_j)^2 \eta_k^2 (\teta_\ell-\eta_\ell)\eta_\ell(\teta_i-\eta_i)$\\ 
     &  $\sum (\teta_i-\eta_i)(\teta_j-\eta_j)\eta_j(\teta_k-\eta_k)^2\eta_\ell \widetilde{\Omega}_{\ell i}$  &  $F_b$ &  $\sum \eta_i (\teta_j-\eta_j)\eta_j(\teta_k-\eta_k)^2\eta_\ell^2(\teta_i-\eta_i)$\\
    &  $\sum (\teta_i-\eta_i)(\teta_j-\eta_j)\eta_j(\teta_k-\eta_k)\eta_k(\teta_\ell -\eta_\ell) \widetilde{\Omega}_{\ell i}$  &  $F_a$ &  $\sum \eta_i (\teta_j-\eta_j)\eta_j(\teta_k-\eta_k)\eta_k(\teta_\ell-\eta_\ell)\eta_\ell (\teta_i-\eta_i)$\\ 
$R_{17}$  & $\sum (\teta_i-\eta_i)(\teta_j-\eta_j)\eta_j(\teta_k-\eta_k)\widetilde{\Omega}_{k\ell}(\teta_\ell-\eta_\ell)\eta_i$  & $F_a$ &  $\sum \eta_i (\teta_j-\eta_j)\eta_j(\teta_k-\eta_k)\eta_k(\teta_\ell-\eta_\ell)\eta_\ell (\teta_i-\eta_i)$\\ 
  & $\sum (\teta_i-\eta_i)(\teta_j-\eta_j)^2\eta_k \widetilde{\Omega}_{k\ell}(\teta_\ell-\eta_\ell)\eta_i$    & $F_b$ &  $\sum \eta_i (\teta_j-\eta_j)^2 \eta_k^2 (\teta_\ell-\eta_\ell)\eta_\ell(\teta_i-\eta_i)$\\ 
  & $\sum (\teta_i-\eta_i)^2(\teta_j-\eta_j)^2 \eta_k\widetilde{\Omega}_{k\ell} \eta_\ell$    & $F_c$ &  $\sum\eta_\ell(\teta_i-\eta_i)^2\eta_k^2(\teta_j-\eta_j)^2\eta_\ell$\\
$R_{21}$   & $\sum (\teta_i-\eta_i)(\teta_j-\eta_j)^2(\teta_k-\eta_k) \widetilde{\Omega}_{k\ell} \widetilde{\Omega}_{\ell i}$ & $F_b$ &  $\sum \eta_i (\teta_j-\eta_j)^2\eta_k (\teta_k-\eta_k)\eta_\ell^2(\teta_i-\eta_i)$\\
 $R_{22}$     & $\sum (\teta_i-\eta_i)(\teta_j-\eta_j) \widetilde{\Omega}_{jk}(\teta_k-\eta_k)(\teta_\ell-\eta_\ell) \widetilde{\Omega}_{\ell i}$ &$F_a$ &  $\sum \eta_i (\teta_j-\eta_j)\eta_j(\teta_k-\eta_k)\eta_k(\teta_\ell-\eta_\ell)\eta_\ell (\teta_i-\eta_i)$\\ 
\hline
\end{tabular} 
} 
\end{table}

\paragraph{Analysis of post-expansion sums with $N^*_W=5$} There are $10$ such terms, including $R_1$, $R_7$-$R_8$, $R_{14}$-$R_{15}$, $R_{18}$, and $R_{25}$-$R_{28}$. Using the the notation
\[
G_i = \teta_i-\eta_i,
\]
we get the following expressions (note: factors of $(\frac{v}{V})^m$ have been removed; see explanations in \eqref{remainder1-Nullgoal2}-\eqref{remainder1-Altgoal2}):
\begin{align*}
R_1 &= \sum_{i, j, k,\ell (dist)} G_iG_jW_{jk} W_{k\ell} W_{\ell i},\cr
R_7 &=  \sum_{i, j, k,\ell (dist)} G_iG_j\eta_j G_kW_{k\ell} W_{\ell i}+ \sum_{i, j, k,\ell (dist)} G_iG^2_j\eta_k W_{k\ell} W_{\ell i}\cr
&=  \sum_{i, j, k,\ell (dist)} \eta_j (G_iG_jG_kW_{k\ell} W_{\ell i}) + \sum_{i, j, k,\ell (dist)} \eta_k (G_iG_j^2 W_{k\ell} W_{\ell i}),\cr
R_8 &= 2 \sum_{i, j, k,\ell (dist)} G_iG_jW_{jk}\eta_kG_\ell W_{\ell i} = 2 \sum_{i, j, k,\ell (dist)} \eta_k(G_iG_jG_\ell W_{jk}W_{\ell i}),\cr
R_{14} &=  \sum_{\substack{i, j, k,\ell\\ (dist)}} G_iG_j^2\eta_k^2G_\ell W_{\ell i} +2  \sum_{\substack{i, j, k,\ell\\ (dist)}} G_iG_j^2\eta_kG_k \eta_\ell W_{\ell i} +    \sum_{\substack{i, j, k,\ell\\ (dist)}} G_iG_j\eta_jG_k\eta_kG_\ell W_{\ell i} \cr
&=   \sum_{\substack{i, j, k,\ell\\ (dist)}} \eta_k^2(G_iG_j^2G_\ell W_{\ell i}) +2  \sum_{\substack{i, j, k,\ell\\ (dist)}} \eta_k\eta_\ell(G_iG_j^2G_kW_{\ell i}) +   \sum_{\substack{i, j, k,\ell\\ (dist)}} \eta_j\eta_k (G_iG_jG_kG_\ell W_{\ell i}), \cr
R_{15} &=  \sum_{\substack{i, j, k,\ell\\ (dist)}} G_iG_j \eta_jG_k W_{k\ell }G_\ell \eta_i  + 2 \sum_{\substack{i, j, k,\ell\\ (dist)}} G_iG^2_j \eta_k W_{k\ell }G_\ell \eta_i +  \sum_{\substack{i, j, k,\ell\\ (dist)}} G_iG_j^2 \eta_k W_{k\ell }\eta_\ell G_i \cr
&=  \sum_{\substack{i, j, k,\ell\\ (dist)}} \eta_i\eta_j (G_iG_j G_k G_\ell W_{k\ell })   + 2 \sum_{\substack{i, j, k,\ell\\ (dist)}} \eta_i\eta_k (G_iG^2_j G_\ell W_{k\ell }) +  \sum_{\substack{i, j, k,\ell\\ (dist)}} \eta_k\eta_\ell (G^2_iG_j^2 W_{k\ell}), \cr
R_{18} &= 4\sum_{i, j, k, \ell (dist)} \eta_j\eta_k\eta_\ell (G_i^2G_jG_kG_\ell) + 4\sum_{i, j, k, \ell (dist)} \eta_k\eta^2_\ell (G_i^2G^2_jG_k), \\
R_{25} &=\sum_{i, j, k, \ell (dist)} G_iG_j^2G_k\widetilde{\Omega}_{k\ell}W_{\ell i}=\sum_{i, j, k, \ell (dist)} \widetilde{\Omega}_{k\ell} (G_iG_j^2G_kW_{\ell i}),\cr
R_{26} &=\sum_{i, j, k, \ell (dist)} G_iG_j\widetilde{\Omega}_{jk}G_kG_\ell W_{\ell i}=\sum_{i, j, k, \ell (dist)} \widetilde{\Omega}_{jk} (G_iG_jG_kG_\ell W_{\ell i}),\cr
 R_{27} &=\sum_{i, j, k,\ell (dist)} G_iG_j^2 G_k\eta_kG_\ell \widetilde{\Omega}_{\ell i} + \sum_{i, j, k,\ell (dist)} G_iG_j^2 G^2_k\eta_\ell \widetilde{\Omega}_{\ell i}\cr
 &=\sum_{i, j, k,\ell (dist)} \eta_k \widetilde{\Omega}_{\ell i} (G_iG_j^2 G_kG_\ell)  + \sum_{i, j, k,\ell (dist)} \eta_\ell \widetilde{\Omega}_{\ell i}(G_iG_j^2 G^2_k),\cr
R_{28} &= 2\sum_{i, j, k,\ell (dist)} G_iG_j\eta_j G^2_kG_{\ell} \widetilde{\Omega}_{\ell i} =2\sum_{i, j, k,\ell (dist)}   \eta_j \widetilde{\Omega}_{\ell i} (G_iG_jG^2_kG_{\ell}).
\end{align*}
Each expression above belongs to one of the following types: 
\begin{align*}
& J_1 = \sum_{i, j, k,\ell (dist)} G_iG_jW_{jk} W_{k\ell} W_{\ell i}, &&J_2 = \sum_{i, j, k,\ell (dist)} \eta_j (G_iG_jG_kW_{k\ell} W_{\ell i}),\\
&J_3 = \sum_{i, j, k,\ell (dist)} \eta_k(G_iG_jG_\ell W_{jk}W_{\ell i}), && J_4 = \sum_{i, j, k,\ell (dist)} \eta_k (G_iG_j^2 W_{k\ell} W_{\ell i}),\\
& J_5 = \sum_{i, j, k,\ell (dist)} \eta_j\eta_k (G_iG_j G_kG_\ell W_{\ell i}), && J'_5 = \sum_{i, j, k,\ell (dist)} \widetilde{\Omega}_{jk} (G_iG_j G_kG_\ell W_{\ell i}),\\
& J_6 = \sum_{i, j, k,\ell (dist)} \eta_k\eta_\ell (G_iG^2_j G_k W_{\ell i}), && J'_6 = \sum_{i, j, k,\ell (dist)} \widetilde{\Omega}_{k\ell}(G_iG^2_j G_k W_{\ell i}),\\
&J_7 = \sum_{i, j, k,\ell (dist)} \eta^2_k (G_iG^2_j G_\ell W_{\ell i}), && J_8 = \sum_{i, j, k,\ell (dist)} \eta_k\eta_\ell (G_i^2G_j^2W_{k\ell}),\cr
& J_9 = \sum_{i, j, k,\ell (dist)} \eta_k \widetilde{\Omega}_{\ell i} (G_iG_j^2 G_kG_\ell), && J_{10} = \sum_{i, j, k,\ell (dist)} \eta_\ell \widetilde{\Omega}_{\ell i}(G_iG_j^2 G^2_k). 
\end{align*}
Since $|\eta_j\eta_k|\leq C\theta_j\theta_k$ and $|\widetilde{\Omega}_{jk}|\leq C\alpha\theta_j\theta_k$, 
the study of $J_5$ and $J_5'$ are similar. Also, the study of $J_6$ and $J_{6}'$ are similar. We now study $J_1$-$J_{10}$. Consider $J_1$. It is seen that
\[
J_1 = \frac{1}{v}\sum_{i, j, k,\ell (dist)} \Bigl(\sum_{s\neq i}W_{is}\Bigr)\Bigl(\sum_{t\neq j}W_{jt}\Bigr)W_{jk} W_{k\ell} W_{\ell i}= \frac{1}{v}\sum_{\substack{i, j, k,\ell (dist)\\s\neq i, t\neq j}} W_{is}W_{i\ell}W_{jt}W_{jk} W_{k\ell}. 
\]
Since $s$ can be equal to $\ell$ and $t$ can be equal to $k$, there are three different types:
\begin{align*}
& J_{1a} = \frac{1}{v}\sum_{i,j,k,\ell (dist)}W^2_{i\ell}W^2_{jk}W_{k\ell}, \qquad J_{1b}= \frac{1}{v}\sum_{\substack{i, j, k,\ell (dist)\\t\notin\{j, k\}}} W^2_{i\ell}W_{jt}W_{jk} W_{k\ell},\cr
& J_{1c} = \frac{1}{v}\sum_{\substack{i, j, k,\ell (dist)\\s\notin\{ i,\ell\}, t\notin\{j,k\}}} W_{is}W_{i\ell}W_{jt}W_{jk} W_{k\ell}.
\end{align*}
We now calculate $\mathbb{E}[J_{1a}^2]$-$\mathbb{E}[J_{1c}^2]$. Take $J_{1a}$ for example. In order to get nonzero 
$\mathbb{E}[W^2_{i\ell}W^2_{jk}W_{k\ell}W^2_{i'\ell'}W^2_{j'k'}W_{k'\ell'}]$, we need either $W_{k\ell}=W_{k'\ell'}$ or each of the two variables $(W_{k\ell}, W_{k',\ell'})$ equals to another squared-$W$ term. The leading term of $\mathbb{E}[J_{1a}^2]$ comes from the first case. In this case, we have $W_{k\ell}=W_{k'\ell'}$ but allow for $W_{i\ell}\neq W_{i'\ell'}$ and $W_{jk}\neq W_{j'k'}$. It has to be the case of either $(k',\ell')=(k,\ell)$ or $(k',\ell')=(\ell,k)$. Therefore, we have $\mathbb{E}[W^2_{i\ell}W^2_{jk}W_{k\ell}W^2_{i'\ell'}W^2_{j'k'}W_{k'\ell'}]=\mathbb{E}[W^2_{i\ell}W^2_{jk}W^2_{i'\ell'}W^2_{j'k'}W^2_{k\ell}]$. Using similar arguments, we have the following results, where details are omitted, as they are similar to the calculations in the proof of Lemmas~\ref{lem:ProxySgnQ(a)-null}-\ref{lem:ProxySgnQ(c)-alt}. 
\begin{align*}
\mathbb{E}[J^2_{1a}] &\leq \frac{C}{v^2}\sum_{\substack{i,j,k,\ell \\i', j'}}  \mathbb{E}[W^2_{i\ell}W^2_{jk}W^2_{i'\ell}W^2_{j'k}W^2_{k\ell}]\leq \frac{C}{\|\theta\|_1^4}\sum_{\substack{i,j,k,\ell \\i', j'}}\theta_i\theta_j\theta_k^3\theta^3_\ell\theta_{i'}\theta_{j'}\leq C\|\theta\|_3^6,\cr
\mathbb{E}[J^2_{1b}] &\leq \frac{C}{v^2} \sum_{\substack{i, j, k,\ell, t\\i'}} \mathbb{E}[W^2_{i\ell}W^2_{i'\ell}W^2_{jt}W^2_{jk} W^2_{k\ell}]\leq \frac{C}{\|\theta\|_1^4} \sum_{\substack{i, j, k,\ell, t\\i'}} \theta_i\theta_j^2\theta_k^2\theta_\ell^3\theta_t\theta_{i'}\leq \frac{C\|\theta\|^4\|\theta\|_3^3}{\|\theta\|_1},\cr
\mathbb{E}[J^2_{1c}]  &\leq \frac{C}{v^2}\sum_{i, j, k,\ell,s,t} \mathbb{E}[W^2_{is}W^2_{i\ell}W^2_{jt}W^2_{jk} W^2_{k\ell}] \leq \frac{C}{\|\theta\|_1^4} \sum_{i,j,k,\ell,s,t}\theta_i^2\theta_j^2\theta_k^2\theta_\ell^2\theta_s\theta_t\leq \frac{C\|\theta\|^8}{\|\theta\|^2_1}.
\end{align*}
The right hand sides are all $o(\|\theta\|^8)$. It follows that
\[
\mathbb{E}[J_1^2] = o(\|\theta\|^8), \qquad \mbox{under both hypotheses}. 
\]
Consider $J_2$-$J_4$. By definition,
\begin{align*}
J_2&=\frac{1}{v\sqrt{v}}\sum_{\substack{i, j, k,\ell (dist)\\s\neq i,t\neq j,q\neq k}} \eta_j W_{is}W_{jt}W_{k q}W_{k\ell}W_{\ell i},\qquad & J_3=\frac{1}{v\sqrt{v}}\sum_{\substack{i, j, k,\ell (dist)\\s\neq i,t\neq j,q\neq \ell}} \eta_k W_{is}W_{jt}W_{\ell q}W_{jk}W_{\ell i},\cr
J_4 &= \frac{1}{v\sqrt{v}}\sum_{\substack{i, j, k,\ell (dist)\\s\neq i,t\neq j,q\neq j}} \eta_k W_{is}W_{jt}W_{jq}W_{k\ell}W_{\ell i}. 
\end{align*}
The analysis is summarized in Table~\ref{tb:J2to4}. In the first column of this table, we study different types of summands. For example, in the expression of $J_2$, $W_{is}W_{kq}W_{k\ell}W_{\ell i}$ have four different cases: (a) $W_{k\ell}^2W^2_{\ell i}$, (b) $W^2_{k\ell}W_{\ell i}W_{is}$ or $W_{k\ell}W^2_{\ell i}W_{kq}$, (c) $W_{k\ell}W_{\ell i}W^2_{ik}$, and (d) $W_{k\ell}W_{\ell i}W_{is}W_{kq}$. In cases (b) and (d), $W_{is}$ or $W_{kq}$ may further equal to $W_{jt}$. Having explored all variants and considered index symmetry, we end up with 6 different cases, as listed in the first column of Table~\ref{tb:J2to4}. In the second column, we study the mean of the squares of the sum of each type of summands. Take the first row for example. We aim to study
\[
\mathrm{E}\Bigl[\Bigl(\sum_{\substack{i,j,k,\ell (dist)\\t\neq j}} \eta_j (W_{k\ell}^2W_{\ell i}^2)W_{jt} \Bigr)\Bigr].
\] 
The naive expansion gives the sum of $\eta_j\eta_{j'}\, \mathbb{E}[W_{k\ell}^2W_{\ell i}^2W_{jt}W_{k'\ell'}^2W_{\ell' i'}^2W_{j't'}]$ over $(i,j,k,\ell,t,i',j',k',\ell',t')$. However, for this term to be nonzero, all single-$W$ terms have to be paired (either with another single-$W$ term or with a squared-$W$ term). The main contribution is from the case of $W_{jt}=W_{j't'}$. This is satisfied only when $(j', s')=(j,s)$ or $(j', s')=(s, j)$. By calculations which are omitted here, we can show that $(j',s')=(j,s)$ yields a larger bound. Therefore, it reduces to the sum of $\eta_j^2\, \mathbb{E}[(W_{jt}^2)W_{k\ell}^2W_{\ell i}^2W_{k'\ell'}^2W^2_{\ell' i'}]$ over $(i,j,k,\ell,t,i',k',\ell')$, which is displayed in the second column of the table. In the last column, we sum the quantity in the second column over indices; it gives rise to a bound for the mean of the square of sum. See the table for details. Recall that the definition of $J_2$-$J_4$ contains a factor of $\frac{1}{v\sqrt{v}}$ in front of the sum, where $v\asymp\|\theta\|_1^2$ by \eqref{v-order}. Hence, to get a desired bound, we only need that each row in the third column of Table~\ref{tb:J2to4} is
\[
o(\|\theta\|^8\|\theta\|_1^6). 
\]  
This is true. We thus conclude that 
\[
\max\big\{ \mathbb{E}[J^2_2],\;  \mathbb{E}[J^2_3],\;  \mathbb{E}[J^2_4]\bigr\} = o(\|\theta\|^8), \qquad \mbox{under both hypotheses}.
\]

\begin{table}[htb!]  
\caption{Analysis of $J_2$-$J_4$. In the second column, the variables in brackets are paired $W$ terms.}    \label{tb:J2to4}
\centering
\scalebox{.95}{
\begin{tabular}{l  l l c }
  & Types of summand  &  Terms in mean-squared & Bound \\
\hline
 \multirow{6}{*}{$J_2$} & $\eta_j (W_{k\ell}^2W_{\ell i}^2)W_{jt}$ & $\eta_j^2\,\mathbb{E}[(W_{jt}^2)W_{k\ell}^2W_{\ell i}^2W_{k'\ell'}^2W^2_{\ell' i'}]\leq \theta_i \theta^3_j\theta_k\theta^2_\ell\theta_t\theta_{i'} \theta_{k'}\theta_{\ell'}^2$ & $\|\theta\|^4\|\theta\|_3^3\|\theta\|_1^5$ \\
  &  $\eta_j (W_{k\ell}W_{\ell i}W_{ik}^2)W_{jt}$ & $\eta_j^2\, \mathbb{E}[ (W^2_{k\ell}W^2_{\ell i}W^2_{jt})W_{ik}^4]\leq C\theta^2_i\theta^3_j\theta^2_k\theta^2_\ell\theta_t$  & $\|\theta\|^6\|\theta\|_3^3\|\theta\|_1$ \\
 & $\eta_j (W_{k\ell}^2W_{\ell i}W_{is})W_{jt}$ & $\eta_j^2\,\mathbb{E}[(W^2_{\ell i}W^2_{is}W^2_{jt})W^2_{k\ell}W^2_{k'\ell}]\leq C\theta^2_i\theta^3_j\theta_k\theta^3_\ell\theta_s\theta_t\theta_{k'}$ & $\|\theta\|^2\|\theta\|_3^6\|\theta\|_1^4$\\
& $\eta_j (W_{k\ell}^2W_{\ell i})W_{ij}^2$ & $\eta_j\eta_{j'}\, \mathbb{E}[(W^2_{\ell i})W^2_{k\ell}W^2_{ij}W^2_{k'\ell}W^2_{i j'}]\leq C\theta^3_i\theta^2_j\theta_k\theta^3_\ell\theta^2_{j'}\theta_{k'}$ & $\|\theta\|^4\|\theta\|_3^6\|\theta\|_1^2$\\
& $\eta_j (W_{k\ell}W_{\ell i}W_{kq}W_{is})W_{jt}$ & $\eta_j^2\,\mathbb{E}[(W^2_{k\ell}W^2_{\ell i}W^2_{kq}W^2_{is}W^2_{jt})]\leq C\theta^2_i\theta^3_j\theta^2_k\theta^2_\ell\theta_s\theta_t\theta_q$ & $\|\theta\|^6\|\theta\|_3^3\|\theta\|_1^3$\\
& $\eta_j (W_{k\ell}W_{\ell i})W_{kq}W_{ij}^2$ & $\eta_j\eta_{j'}\,\mathbb{E}[(W^2_{k\ell}W^2_{\ell i}W^2_{kq})W^2_{ij}W^2_{ij'}]\leq C\theta^3_i\theta^2_j\theta^2_k\theta^2_\ell\theta_q\theta^2_{j'}$ & $\|\theta\|^8\|\theta\|^3_3\|\theta\|_1$\\
\hline
\multirow{8}{*}{$J_3$}& $\eta_k W^3_{\ell i}W_{jk}^2$ & $\eta_k\eta_{k'}\,\mathbb{E}[W^3_{\ell i}W^2_{jk}W^3_{\ell'i'}W^2_{j'k'}]\leq C\theta_i\theta_j\theta^2_k\theta_\ell\theta_{i'}\theta_{j'}\theta^2_{k'}\theta_{\ell'}$ & $\|\theta\|^4\|\theta\|_1^6$\\
& $\eta_k W^3_{\ell i}(W_{jk}W_{jt})$ & $\eta_k^2\,\mathbb{E}[(W^2_{jk}W^2_{jt})W^3_{\ell i}W^3_{\ell' i'}]\leq C\theta_i\theta^2_j\theta^3_k\theta_\ell\theta_t\theta_{i'}\theta_{\ell'}$ & $\|\theta\|^2\|\theta\|_3^3\|\theta\|_1^5$\\
& $\eta_k (W^2_{\ell i}W_{is})W_{jk}^2$ & $\eta_k\eta_{k'}\, \mathbb{E}[(W^2_{is})W^2_{\ell i}W^2_{jk}W^2_{\ell' i}W^2_{j'k'}]\leq C \theta^3_i\theta_j\theta^2_k\theta_\ell\theta_s\theta_{j'}\theta^2_{k'}\theta_{\ell'}$ & $\|\theta\|^4\|\theta\|_3^3\|\theta\|_1^5$ \\
& $\eta_k (W^2_{\ell i}W_{is})W_{jk}W_{jt}$ & $\eta^2_k\, \mathbb{E}[(W^2_{is}W^2_{jk}W^2_{jt})W^2_{\ell i}W^2_{\ell'i}]\leq C\theta^3_i\theta^2_j\theta^3_k\theta_\ell\theta_s\theta_t\theta_{\ell'}$ & $\|\theta\|^2\|\theta\|_3^6\|\theta\|_1^4$\\
& $\eta_k W^2_{\ell i}W^2_{ij}W_{jk}$ & $\eta^2_k\, \mathbb{E}[(W^2_{jk})W^2_{\ell i}W^2_{ij}W^2_{\ell' i'}W^2_{i'j}]\leq C \theta^2_i\theta^3_j\theta^3_k\theta_\ell\theta^2_{i'}\theta_{\ell'}$ & $\|\theta\|^4\|\theta\|_3^6\|\theta\|_1^2$\\
& $\eta_k (W_{\ell i}W_{is}W_{\ell q})W_{jk}^2$ & $\eta_k\eta_{k'}\, \mathbb{E}[ (W^2_{\ell i}W^2_{is}W^2_{\ell q})W_{jk}^2W_{j'k'}^2]\leq C \theta^2_i\theta_j\theta^2_k\theta^2_\ell\theta_s\theta_q\theta_{j'}\theta^2_{k'}$ & $\|\theta\|^8\|\theta\|_1^4$\\
& $\eta_k (W_{\ell i}W_{is}W_{\ell q})W_{jk}W_{jt}$ & $\eta_k^2\, \mathbb{E}[(W^2_{\ell i}W^2_{is}W^2_{\ell q}W^2_{jk}W^2_{jt})]\leq C \theta^2_i\theta^2_j\theta^3_k\theta^2_\ell\theta_s\theta_t\theta_q$ & $\|\theta\|^6\|\theta\|_3^3\|\theta\|_1^3$\\
& $\eta_k W_{\ell i}W^2_{ij}W_{\ell q}W_{jk}$ & $\eta^2_k\, \mathbb{E}[(W^2_{\ell i}W^2_{\ell q}W^2_{jk})W^4_{ij}]\leq C\theta^2_i\theta^2_j\theta^3_k\theta^2_\ell\theta_q$ & $\|\theta\|^6\|\theta\|_3^3\|\theta\|_1$\\
\hline
\multirow{6}{*}{$J_4$} & $\eta_k(W_{k\ell}W^2_{\ell i})W_{jt}^2$  & $\eta^2_k\, \mathbb{E}[ (W^2_{k\ell}) W^2_{\ell i}W_{jt}^2 W^2_{\ell i'}W_{j' t'}^2 ]\leq C \theta_i\theta_j\theta^3_k\theta^3_\ell\theta_t\theta_{i'}\theta_{j'}\theta_{t'}$ & $\|\theta\|_3^6\|\theta\|_1^6$\\
& $\eta_k(W_{k\ell}W^2_{\ell i})W_{jt}W_{jq}$ & $\eta^2_k\, \mathbb{E}[(W^2_{k\ell}W^2_{jt}W^2_{jq})W^2_{\ell i}W^2_{\ell i'}] \leq C\theta_i\theta^2_j\theta^3_k\theta^3_\ell\theta_t\theta_q\theta_{i'}$ & $\|\theta\|^2\|\theta\|_3^6\|\theta\|_1^4$\\
& $\eta_k(W_{k\ell}W_{\ell i}W_{is})W^2_{jt}$ & $\eta^2_k\, \mathbb{E}[( W^2_{k\ell}W^2_{\ell i}W^2_{is} )W^2_{jt}W^2_{j't'}]\leq C\theta^2_i\theta_j\theta^3_k\theta^2_\ell\theta_s\theta_t\theta_{j'}\theta_{t'}$ & $\|\theta\|^4\|\theta\|_3^3\|\theta\|_1^5$\\
& $\eta_k W_{k\ell}W_{\ell i}W^3_{ij}$ & $\eta^2_k\, \mathbb{E}[(W^2_{k\ell}W^2_{\ell i})W^3_{ij}W^3_{ij'}]\leq C\theta^3_i\theta_j\theta^3_k\theta^2_\ell\theta_{j'}$ &  $\|\theta\|^2\|\theta\|_3^6\|\theta\|_1^2$\\
& $\eta_k(W_{k\ell}W_{\ell i}W_{is})W_{jt}W_{jq}$ & $\eta_k^2\, \mathbb{E}[(W^2_{k\ell}W^2_{\ell i}W^2_{is}W^2_{jt}W^2_{jq}] \leq C\theta^2_i\theta^2_j\theta^3_k\theta^2_\ell\theta_s\theta_t\theta_q $ & $\|\theta\|^6\|\theta\|_3^3\|\theta\|_1^3$\\
& $\eta_k W_{k\ell}W_{\ell i}W^2_{ij}W_{jq}$ & $\eta^2_k\, \mathbb{E}[(W^2_{k\ell}W^2_{\ell i}W^2_{jq})W^4_{ij}]\leq C\theta^2_i\theta^2_j\theta^3_k\theta^2_\ell\theta_q $ & $\|\theta\|^6\|\theta\|_3^3\|\theta\|_1$\\
\hline
\end{tabular} 
} 
\end{table}

Consider $J_5$-$J_8$. It is seen that
\begin{align*}
 J_5 &= \frac{1}{v^2}\sum_{i, j, k,\ell (dist)} \eta_j\eta_kW_{is}W_{jt}W_{kq}W_{\ell m} W_{\ell i}, && J_6 = \frac{1}{v^2} \sum_{i, j, k,\ell (dist)} \eta_k\eta_\ell W_{is}W_{jt}W_{jq} W_{km} W_{\ell i}, \cr
J_7 &= \frac{1}{v^2} \sum_{i, j, k,\ell (dist)} \eta^2_k W_{is}W_{jt}W_{jq} W_{\ell m} W_{\ell i}, && J_8 = \frac{1}{v^2}\sum_{i, j, k,\ell (dist)} \eta_k\eta_\ell W_{is}W_{it}W_{jq}W_{jm}W_{k\ell},\cr
\end{align*}
The analysis is summarized in Table~\ref{tb:J5to8}. We note that $J_7$ can be written as
\[
J_7 = \frac{1}{v^2} \sum_{i, j, \ell (dist)} \beta_{ij\ell} W_{is}W_{jt}W_{jq} W_{\ell m} W_{\ell i}, \qquad\mbox{where}\quad\beta_{ij\ell}\equiv \sum_{k\notin\{i,j,\ell\}}\eta^2_k. 
\]
Although the values of $\beta_{ij\ell}$ change with indices, they have a common upper bound of $C\|\theta\|^2$. We treat $\beta_{ij\ell}$ as $\|\theta\|^2$ in Table~\ref{tb:J5to8}, as this doesn't change the bounds but simplifies notations.  Recall that the definition of $J_5$-$J_8$ contains a factor of $\frac{1}{v^2}$ in front of the sum, where $v\asymp\|\theta\|_1^2$ by \eqref{v-order}. Hence, to get a desired bound, we only need that each row in the third column of Table~\ref{tb:J2to4} is
\[
o(\|\theta\|^8\|\theta\|_1^8). 
\]  
This is true. We thus conclude that 
\[
\max\big\{ \mathbb{E}[J^2_5],\;  \mathbb{E}[J^2_6],\;  \mathbb{E}[J^2_7],\;  \mathbb{E}[J^2_8]\bigr\} = o(\|\theta\|^8), \qquad \mbox{under both hypotheses}.
\]

\begin{table}[htb!]  
\centering
\caption{Analysis of $J_5$-$J_8$. In the second column, the variables in brackets are paired $W$ terms.}    \label{tb:J5to8}
\scalebox{.95}{
\begin{tabular}{l  l l c }
  & Types of summand  &  Terms in mean-squared & Bound \\
\hline
 \multirow{8}{*}{$J_5$} & $\eta_j\eta_k W_{\ell i}^3W_{jk}^2$ & $\eta_j\eta_k\eta_{j'}\eta_{k'}\, \mathbb{E}[W_{\ell i}^3W_{jk}^2W_{\ell' i'}^3W_{j'k'}^2]\leq C\theta_i\theta^2_j\theta^2_k\theta_{\ell}\theta_{i'}\theta^2_{j'}\theta^2_{k'}\theta_{\ell'}$ & $\|\theta\|^8\|\theta\|_1^4$\\
&$\eta_j\eta_k W_{\ell i}^3(W_{jt}W_{kq})$ & $\eta^2_j\eta^2_k\, \mathbb{E}[ (W^2_{jt}W^2_{kq})W_{\ell i}^3W_{\ell' i'}^3 ]\leq C\theta_i\theta^3_j\theta^3_k\theta_\ell\theta_t\theta_q\theta_{i'}\theta_{\ell'}$ & $\|\theta\|_3^6\|\theta\|_1^6$\\
& $\eta_j\eta_k (W_{\ell i}^2W_{is})W_{jk}^2$ & $\eta_j\eta_k\eta_{j'}\eta_{k'}\, \mathbb{E}[ (W^2_{is})W_{\ell i}^2W_{jk}^2W_{\ell' i}^2W_{j'k'}^2]\leq C\theta^3_i\theta^2_j\theta^2_k\theta_\ell\theta_s\theta^2_{j'}\theta^2_{k'}\theta_{\ell'}$ & $\|\theta\|^8\|\theta\|_3^3\|\theta\|_1^3$\\
& $\eta_j\eta_k (W_{\ell i}^2W_{is})(W_{jt}W_{kq})$ & $\eta^2_j\eta^2_k\, \mathbb{E}[ (W^2_{is})W^2_{jt}W^2_{kq}) W^2_{\ell i}W^2_{\ell' i} ]\leq C\theta^3_i\theta^3_j\theta^3_k\theta_\ell\theta_s\theta_t\theta_q\theta_{\ell'}$ & $\|\theta\|_3^9\|\theta\|_1^5$\\
& $\eta_j\eta_k W_{\ell i}^2W_{ij}^2 W_{kq}$ & $\eta_j\eta^2_k\eta_{j'}\, \mathbb{E}[ (W^2_{kq})W_{\ell i}^2W_{ij}^2W_{\ell' i'}^2W_{i'j'}^2 ]\leq C\theta^2_i\theta^2_j\theta^3_k\theta_\ell\theta_q\theta^2_{i'}\theta^2_{j'}\theta_{\ell'}$ & $\|\theta\|^8\|\theta\|_3^3\|\theta\|_1^3$\\
& $\eta_j\eta_k (W_{\ell i}W_{is}W_{\ell m})W_{jk}^2$ & $\eta_j\eta_k\eta_{j'}\eta_{k'}\, \mathbb{E}[ (W^2_{\ell i}W^2_{is}W^2_{\ell m}) W_{jk}^2W_{j'k'}^2]\leq C\theta^2_i\theta^2_j\theta^2_k\theta^2_\ell\theta_s\theta_m\theta^2_{j'}\theta^2_{k'}$ & $\|\theta\|^{12}\|\theta\|_1^2$
\\
& $\eta_j\eta_k (W_{\ell i}W_{is}W_{\ell m})(W_{jt}W_{kq})$ & $\eta^2_j\eta^2_k\, \mathbb{E}[ (W^2_{\ell i}W^2_{is}W^2_{\ell m}W^2_{jt}W^2_{kq})]\leq C\theta^2_i\theta^3_j\theta^3_k\theta^2_\ell\theta_s\theta_t\theta_q\theta_m$ & $\|\theta\|^4\|\theta\|_3^6\|\theta\|_1^4$\\
& $\eta_j\eta_k W_{\ell i}W_{ij}^2W_{\ell m}W_{kq}$ & $\eta_j\eta^2_k\eta_{j'}\, \mathbb{E}[ (W^2_{\ell i}W^2_{\ell m}W^2_{kq})W_{ij}^2W_{ij'}^2]\leq C\theta^3_i\theta^2_j\theta^3_k\theta^2_\ell\theta_q\theta_m\theta^2_{j'} $  & $\|\theta\|^6\|\theta\|_3^6\|\theta\|_1^2$\\
\hline
 \multirow{12}{*}{$J_6$} & $\eta_k\eta_\ell W_{\ell i}^2W_{jt}^2W_{km}$ & $\eta^2_k\eta_\ell\eta_{\ell'}\,\mathbb{E}[ (W^2_{km})W_{\ell i}^2W_{jt}^2W_{\ell' i'}^2W_{j't'}^2 ]\leq C\theta_i\theta_j\theta^3_k\theta^2_\ell\theta_t\theta_m\theta_{i'}\theta_{j'}\theta^2_{\ell'}\theta_{t'}$ & $\|\theta\|^4\|\theta\|_3^3\|\theta\|_1^7$\\
 & $\eta_k\eta_\ell W_{\ell i}^2W_{jk}^3$ & $\eta_k\eta_\ell\eta_{k'}\eta_{\ell'}\,\mathbb{E}[W_{\ell i}^2W_{jk}^3W_{\ell' i'}^2W_{j'k'}^3 ]\leq C\theta_i\theta_j\theta^2_k\theta^2_\ell\theta_{i'}\theta_{j'}\theta^2_{k'}\theta^2_{\ell'}$ & $\|\theta\|^8\|\theta\|_1^4$\\
 & $\eta_k\eta_\ell W_{\ell i}^2(W_{jt}W_{jq})W_{km}$ & $\eta^2_k\eta_\ell\eta_{\ell'}\,\mathbb{E}[ (W^2_{jt}W^2_{jq}W^2_{km}) W_{\ell i}^2W_{\ell' i'}^2 ]\leq C\theta_i\theta^2_j\theta^3_k\theta^2_\ell\theta_t\theta_q\theta_m\theta_{i'}\theta^2_{\ell'}$ & $\|\theta\|^6\|\theta\|_3^3\|\theta\|_1^5$
\\ 
 & $\eta_k\eta_\ell W_{\ell i}^2W^2_{jk}W_{jq}$ & $\eta_k\eta_\ell\eta_{k'}\eta_{\ell'}\,\mathbb{E}[ (W^2_{jq})W_{\ell i}^2W^2_{jk}W_{\ell' i'}^2W^2_{jk'} ] \leq C\theta_i\theta^3_j\theta^2_k\theta^2_\ell\theta_q\theta_{i'}\theta^2_{k'}\theta^2_{\ell'}$ & $\|\theta\|^8\|\theta\|_3^3\|\theta\|_1^3$\\ 
 & $\eta_k\eta_\ell (W_{\ell i}W_{is})W_{jt}^2W_{km}$ & $\eta^2_k\eta^2_\ell\,\mathbb{E}[ (W^2_{\ell i}W^2_{is}W^2_{km})W_{jt}^2W_{j't'}^2 ] \leq C\theta^2_i\theta_j\theta^3_k\theta^3_\ell\theta_s\theta_t\theta_m\theta_{j'}\theta_{t'}$ & $\|\theta\|^2\|\theta\|_3^6\|\theta\|_1^6$
\\
  & $\eta_k\eta_\ell W_{\ell i}W_{ij}^3W_{km}$ & $\eta^2_k\eta^2_\ell\,\mathbb{E}[ (W^2_{\ell i}W^2_{km})W_{ij}^3W_{ij'}^3 ] \leq C\theta^3_i\theta_j\theta^3_k\theta^3_\ell\theta_m\theta_{j'}$ & $\|\theta\|_3^9\|\theta\|_1^3$\\
    & $\eta_k\eta_\ell W_{\ell i}W_{is}W_{jk}^3$ & $\eta_k\eta^2_\ell\eta_{k'}\,\mathbb{E}[ (W^2_{\ell i}W^2_{is})W_{jk}^3W_{j'k'}^3 ] \leq C\theta^2_i\theta_j\theta^2_k\theta^3_\ell\theta_s\theta_{j'}\theta^2_{k'}$ & $\|\theta\|^6\|\theta\|_3^3\|\theta\|_1^3$\\
  & $\eta_k\eta_\ell W_{\ell i}W_{ik}^2 W_{jt}^2$ & $\eta_k\eta^2_\ell\eta_{k'}\,\mathbb{E}[ (W^2_{\ell i})W_{ik}^2 W_{jt}^2W_{ik'}^2 W_{j't'}^2 ] \leq C\theta^3_i\theta_j\theta^2_k\theta^3_\ell\theta_t\theta_{j'}\theta^2_{k'}\theta_{t'}$ & $\|\theta\|^4\|\theta\|_3^6\|\theta\|_1^4$\\
 & $\eta_k\eta_\ell (W_{\ell i}W_{is})(W_{jt}W_{jq})W_{km}$ & $\eta^2_k\eta^2_\ell\,\mathbb{E}[ (W^2_{\ell i}W^2_{is}W^2_{jt}W^2_{jq}W^2_{km})] \leq C\theta^2_i\theta^2_j\theta^3_k\theta^3_\ell\theta_s\theta_t\theta_q\theta_m$ & $\|\theta\|^4\|\theta\|_3^6\|\theta\|_1^4$\\ 
 & $\eta_k\eta_\ell W_{\ell i}W^2_{ij}W_{jq}W_{km}$ & $\eta^2_k\eta^2_\ell\,\mathbb{E}[ (W^2_{\ell i}W^2_{jq}W^2_{km})W^4_{ij}] \leq C\theta^2_i\theta^2_j\theta^3_k\theta^3_\ell\theta_q\theta_m$ & $\|\theta\|^4\|\theta\|_3^6\|\theta\|_1^2$\\ 
 & $\eta_k\eta_\ell W_{\ell i}W_{is}W^2_{jk}W_{jq}$ & $\eta_k\eta^2_\ell\eta_{k'}\,\mathbb{E}[ (W^2_{\ell i}W^2_{is}W^2_{jq}) W^2_{jk}W^2_{jk'}] \leq C\theta^2_i\theta^3_j\theta^2_k\theta^3_\ell\theta_s\theta_q\theta^2_{k'}$ & $\|\theta\|^6\|\theta\|_3^6\|\theta\|_1^2$\\ 
& $\eta_k\eta_\ell W_{\ell i}W_{ik}^2 W_{jt}W_{jq}$ & $\eta_k\eta^2_\ell\eta_{k'}\,\mathbb{E}[ (W^2_{\ell i}W^2_{jt}W^2_{jq})W_{ik}^2W_{ik'}^2 ] \leq C\theta^3_i\theta^2_j\theta^2_k\theta^3_\ell\theta_t\theta_q\theta^2_{k'}$ & $\|\theta\|^6\|\theta\|_3^6\|\theta\|_1^2$\\
\hline
\multirow{11}{*}{$J_7$} & $\|\theta\|^2 W_{\ell i}^3W^2_{jt}$ & $\|\theta\|^4\,\,\mathbb{E}[W_{\ell i}^3W^2_{jt}W_{\ell' i'}^3W^2_{j't'}] \leq C\|\theta\|^4\theta_i\theta_j\theta_\ell\theta_t\theta_{i'}\theta_{j'}\theta_{\ell'}\theta_{t'}$ & $\|\theta\|^4 \|\theta\|_1^8$\\
& $\|\theta\|^2 W_{\ell i}^3(W_{jt}W_{jq})$ & $\|\theta\|^4\,\mathbb{E}[ (W^2_{jt}W^2_{jq})W_{\ell i}^3W_{\ell' i'}^3 ] \leq C\|\theta\|^4 \theta_i\theta^2_j\theta_\ell \theta_t\theta_q\theta_{i'}\theta_{\ell'}$ & $\|\theta\|^6\|\theta\|_1^6$ \\ 
& $\|\theta\|^2 (W_{\ell i}^2W_{is}) W^2_{jt}$ & $\|\theta\|^4 \,\mathbb{E}[ (W^2_{is})W_{\ell i}^2W^2_{jt}W_{\ell' i}^2W^2_{j't'} ] \leq C\|\theta\|^4 \theta^3_i\theta_j\theta_\ell\theta_s\theta_t\theta_{j'}\theta_{\ell'}\theta_{t'}$ & $\|\theta\|^4\|\theta\|_3^3\|\theta\|_1^7$\\
& $\|\theta\|^2 W_{\ell i}^2W_{ij}^3$ & $\|\theta\|^4 \,\mathbb{E}[W_{\ell i}^2W_{ij}^3W_{\ell' i'}^2W_{i'j'}^3 ] \leq C\|\theta\|^4 \theta^2_i\theta_j\theta_\ell\theta^2_{i'}\theta_{j'}\theta_{\ell'}$ & $\|\theta\|^8\|\theta\|_1^4$\\
& $\|\theta\|^2 (W_{\ell i}^2W_{is})(W_{jt}W_{jq})$ & $\|\theta\|^4 \,\mathbb{E}[ (W^2_{is}W^2_{jt}W^2_{jq})W_{\ell i}^2W_{\ell' i}^2 ] \leq C\|\theta\|^4 \theta^3_i\theta^2_j\theta_\ell\theta_s\theta_t\theta_q\theta_{\ell'}$ & $\|\theta\|^6\|\theta\|_3^3\|\theta\|_1^5$\\ 
& $\|\theta\|^2 W_{\ell i}^2W_{ij}^2W_{jq}$ & $\|\theta\|^4 \,\mathbb{E}[ (W^2_{jq})W_{\ell i}^2W_{ij}^2W_{\ell' i'}^2W_{i'j}^2 ] \leq C\|\theta\|^4 \theta^2_i\theta^3_j\theta_\ell \theta_q\theta^2_{i'}\theta_{\ell'}$ & $\|\theta\|^8\|\theta\|_3^3\|\theta\|_1^3$
\\ 
& $\|\theta\|^2 (W_{\ell i}W_{is}W_{\ell m}) W^2_{jt}$ & $\|\theta\|^4 \,\mathbb{E}[ (W^2_{\ell i}W^2_{is}W^2_{\ell m})W^2_{jt}W^2_{j't'} ] \leq C\|\theta\|^4 \theta^2_i\theta_j\theta^2_\ell\theta_s\theta_t\theta_m\theta_{j'}\theta_{t'}$ & $\|\theta\|^8\|\theta\|_1^6$\\
& $\|\theta\|^2 W_{\ell i}W_{ij}^3 W_{\ell m}$ & $\|\theta\|^4 \,\mathbb{E}[ (W^2_{\ell i}W^2_{\ell m})W_{ij}^3W_{ij'}^3 ] \leq C\|\theta\|^4 \theta^3_i\theta_j\theta^2_\ell\theta_m\theta_{j'}$ & $\|\theta\|^6\|\theta\|_3^3\|\theta\|_1^3$\\
& $\|\theta\|^2 (W_{\ell i}W_{is}W_{\ell m}) (W_{jt}W_{jq})$ & $\|\theta\|^4 \,\mathbb{E}[ (W^2_{\ell i}W^2_{is}W^2_{\ell m}W^2_{jt}W^2_{jq})] \leq C\|\theta\|^4 \theta^2_i\theta^2_j\theta^2_\ell\theta_s\theta_t\theta_q\theta_m$ & $\|\theta\|^{10}\|\theta\|_1^4$\\ 
& $\|\theta\|^2 W_{\ell i}W^2_{ij}W_{\ell m}W_{jq}$ & $\|\theta\|^4 \,\mathbb{E}[ (W^2_{\ell i}W^2_{\ell m}W^2_{jq})W^4_{ij}] \leq C\|\theta\|^4 \theta^2_i\theta^2_j\theta^2_\ell\theta_q\theta_m$ & $\|\theta\|^{10}\|\theta\|_1^2$\\ 
& $\|\theta\|^2 W_{\ell i}W^2_{ij}W^2_{\ell j}$ & $\|\theta\|^4 \,\mathbb{E}[ (W^2_{\ell i})W^2_{ij}W^2_{\ell j}W^2_{ij'}W^2_{\ell j'} ] \leq C\|\theta\|^4 \theta^3_i\theta^2_j\theta^3_\ell\theta^2_{j'}$ & $\|\theta\|^8\|\theta\|_3^6$\\ 
\hline
\multirow{4}{*}{$J_8$} & $\eta_k\eta_\ell W^4_{ij}W_{k\ell}$ & $\eta^2_k\eta^2_{\ell}\,\mathbb{E}[ (W^2_{k\ell})W^4_{ij}W^4_{i'j'} ] \leq C\theta_i\theta_j\theta^3_k\theta^3_\ell\theta_{i'}\theta_{j'}$ & $\|\theta\|_3^6\|\theta\|_1^4$\\
& $\eta_k\eta_\ell (W^3_{ij}W_{is})W_{k\ell}$ & $\eta^2_k\eta^2_{\ell} \,\mathbb{E}[ (W^2_{is}W^2_{k\ell})W^3_{ij}W^3_{ij'}] \leq C\theta^3_i\theta_j\theta^3_k\theta^3_\ell\theta_s\theta_{j'}$ & $\|\theta\|_3^9\|\theta\|_1^3$\\
& $\eta_k\eta_\ell (W^2_{ij}W_{is}W_{jq})W_{k\ell}$ & $\eta^2_k\eta^2_{\ell} \,\mathbb{E}[ (W^2_{is}W^2_{jq}W^2_{k\ell})W^4_{ij} ] \leq C\theta^2_i\theta^2_j\theta^3_k\theta^3_\ell\theta_s\theta_q$ & $\|\theta\|^4\|\theta\|_3^6\|\theta\|_1^2$\\
& $\eta_k\eta_\ell (W_{is}W_{it}W_{jq}W_{jm})W_{k\ell}$ & $\eta^2_k\eta^2_{\ell} \,\mathbb{E}[ (W^2_{is}W^2_{it}W^2_{jq}W^2_{jm}W^2_{k\ell}) ] \leq C\theta^2_i\theta^2_j\theta^3_k\theta^3_\ell\theta_s\theta_t\theta_q\theta_m$ & $\|\theta\|^4\|\theta\|_3^6\|\theta\|_1^4$\\
& $\eta_k \eta_\ell W^2_{is}W_{jq}W_{jm}W_{k\ell}$ & $\eta_k^2\eta_\ell^2 \mathbb{E}[(W_{jq}^2W_{jm}^2W^2_{k\ell})W^2_{is}W^2_{i's'}]\leq C\theta_i\theta^2_j\theta^3_k\theta^3_\ell\theta_s\theta_q\theta_m\theta_{i'}\theta_{s'}$ & $\|\theta\|^2\|\theta\|_3^6\|\theta\|_1^6$\\
& $\eta_k \eta_\ell W^2_{is}W^2_{jq}W_{k\ell}$ & $\eta_k^2\eta_\ell^2 \mathbb{E}[(W^2_{k\ell})W^2_{is}W^2_{jq}W^2_{i's'}W^2_{j'q'}]\leq C\theta_i\theta_j\theta^3_k\theta^3_\ell\theta_s\theta_q\theta_{i'}\theta_{j'}\theta_{s'}\theta_{q'}$ & $\|\theta\|_3^6\|\theta\|_1^8$\\
\hline
\end{tabular} 
} 
\end{table}

Consider $J_9$-$J_{10}$. They can be analyzed in the same way as we did for $J_1$-$J_8$. To save space, we only give a simplified proof for the case of $\|\theta\|\gg \alpha [\log(n)]^{5/2}$. For $1\ll \|\theta\|\leq C\alpha [\log(n)]^{5/2}$, the proof is similar to those in Tables~\ref{tb:J2to4}-\ref{tb:J5to8}, which is omitted. For a constant $C_0>0$ to be decided, we introduce an event
\beq  \label{eventE}
E = \cap_{i=1}^n E_i, \qquad \mbox{where}\quad E_i= \big\{ \sqrt{v} |G_i|\leq C_0\sqrt{\theta_i\|\theta\|_1\log(n)} \bigr\}. 
\eeq
Recall that $\sqrt{v}G_i = \sqrt{v}(\teta_i-\eta_i) = \sum_{j\neq i}(A_{ij}-\mathbb{E}A_{ij})$. The variables $\{A_{ij}\}_{j\neq i}$ are mutually independent, satisfying that $|A_{ij}-\mathbb{E}A_{ij}|\leq 1$ and $\sum_j \mathrm{Var}(A_{ij})\leq \sum_j \theta_i\theta_j\leq \theta_i\|\theta\|_1$. By Bernstein's inequality, for large $n$, the probability of $E_i^c$ is $O(n^{-C_0/4.1})$. Applying the probability union bound, we find that the probability of $E^c$ is $O(n^{-C_0/2.01})$. Recall that $V=\sum_{i,j: i\neq j}A_{ij}$. On the event $E^c$, if $V=0$ (i.e., the network has no edges), then $\widetilde{Q}^*_n=Q^*_n=0$; otherwise, $V\geq 1$ and $|\widetilde{Q}^*_n-Q^*_n|\leq n^4$. Combining these results gives
\[
\mathbb{E}\bigl[|\widetilde{Q}^*_n-Q^*_n|^2\cdot I_{E^c}\bigr] \leq  n^4\cdot O(n^{-C_0/2.01}). 
\]  
With an properly large $C_0$, the right hand side is  $o(\|\theta\|^8)$. Hence, it suffices to focus on the event $E$. On the event $E$,
\begin{align*}
|J_9|&\leq \sum_{i,j,k,\ell}|\eta_k\widetilde{\Omega}_{\ell i}||G_iG_j^2G_kG_\ell|\cr
&\leq C\sum_{i,j,k,\ell}(\alpha\theta_i \theta_k\theta_\ell)\frac{\sqrt{\theta_i\theta_j^2\theta_k\theta_\ell\|\theta\|_1^5 [\log(n)]^5}}{\sqrt{v^5}} \cr
&\leq \frac{C\alpha[\log(n)]^{5/2}}{\sqrt{\|\theta\|_1^5}} \Bigl(\sum_i \theta^{3/2}_i\Bigr)\Bigl(\sum_j \theta_j\Bigr)\Bigl(\sum_k \theta^{3/2}_k\Bigr)\Bigl(\sum_\ell \theta^{3/2}_\ell\Bigr) \cr
&\leq \frac{C\alpha [\log(n)]^{5/2}}{\sqrt{\|\theta\|_1^3}} \Bigl(\sum_i \theta^{3/2}_i\Bigr)^3\cr
&\leq \frac{C\alpha [\log(n)]^{5/2}}{\sqrt{\|\theta\|_1^3}} \Bigl(\sum_i\theta_i^2\Bigr)^{3/2}\Bigl(\sum_i\theta_i\Bigr)^{3/2}\cr
&\leq C\alpha [\log(n)]^{5/2}\|\theta\|^3,
\end{align*}
where the second last line is from the Cauchy-Schwarz inequality. Since $\|\theta\|\gg \alpha[\log(n)]^{5/2}$, the right hand side is $o(\|\theta\|^4)$, which implies that $|J_9|^2=o(\|\theta\|^8)$. Similarly, on the event $E$, 
\begin{align*}
|J_{10}| &\leq  \sum_{i, j, k,\ell} |\eta_\ell \widetilde{\Omega}_{\ell i}| |G_iG_j^2 G^2_k|\cr
&\leq C\sum_{i,j,k,\ell}(\alpha\theta_i\theta_\ell^2)\frac{\sqrt{\theta_i\theta_j^2\theta^2_k\|\theta\|_1^5[\log(n)]^5}}{\sqrt{v^5}}\cr
&\leq \frac{C\alpha[\log(n)]^{5/2}}{\sqrt{\|\theta\|_1^5}}\Bigl(\sum_i \theta_i^{3/2}\Bigr)\Bigl(\sum_j \theta_j\Bigr)\Bigl(\sum_k\theta_k\Bigr)\Bigl(\sum_\ell \theta_\ell^2\Bigr)\cr
&\leq \frac{C\alpha[\log(n)]^{5/2}}{\sqrt{\|\theta\|_1^5}}\bigl(\|\theta\|\sqrt{\|\theta\|_1}\bigr)\|\theta\|^2_1\|\theta\|^2\cr
&\leq C\alpha[\log(n)]^{5/2}\|\theta\|^3;
\end{align*} 
again, the right hand side is $o(\|\theta\|^4)$. Combining the above gives
\[
\max\big\{ \mathbb{E}[J^2_9],\;  \mathbb{E}[J^2_{10}]\bigr\} = o(\|\theta\|^8), \qquad \mbox{under both hypotheses}.
\] 
So far, we have proved: for each $R_k$ with $N^*_W=5$, it satisfies $\mathbb{E}[R_k^2]=o(\|\theta\|^8)$. This is sufficient to guarantee \eqref{remainder1-Nullgoal2}-\eqref{remainder1-Altgoal2} for $X=R_k$.

\paragraph{Analysis of post-expansion sums with $N^*_W=6$} There are $7$ such terms, including $R_{19}$-$R_{20}$, $R_{23}$-$R_{24}$, $R_{29}$-$R_{30}$, and $R_{32}$. We plug in the definition of $\tilde{r}_{ij}$ and $\delta_{ij}$ and neglect all factors of $\frac{v}{V}$ (see the explanation in \eqref{remainder1-Nullgoal2}-\eqref{remainder1-Altgoal2}). It gives ($G_i=\teta_i-\eta_i$):
\begin{align*}
R_{19} & = \sum_{i, j, k,\ell (dist)} G_iG_j^2G_k W_{k\ell} W_{\ell i},\cr
R_{20} & = \sum_{i, j, k,\ell (dist)} G_iG_j W_{jk}G_kG_\ell W_{\ell i}, \cr
R_{23} & = \sum_{i, j, k,\ell (dist)} G_iG_j^2G_k( \eta_kG_\ell^2\eta_i+2G_k\eta_\ell G_\ell \eta_i + G_k\eta_\ell^2 G_i) \cr
& = \sum_{i, j, k,\ell (dist)} \eta_i \eta_kG_iG_j^2G_kG_\ell^2 + 2 \sum_{i, j, k,\ell (dist)} \eta_i\eta_\ell G_iG_j^2G^2_kG_\ell + \sum_{i, j, k,\ell (dist)} \eta_\ell^2 G^2_iG_j^2G^2_k \cr
&= 3\sum_{i, j, k,\ell (dist)} \eta_i \eta_kG_iG_j^2G_kG_\ell^2 + \sum_{i, j, k,\ell (dist)} \eta_\ell^2 G^2_iG_j^2G^2_k, \cr
R_{24} &= \sum_{i, j, k, \ell (dist)} G_iG_j(\eta_jG_k+G_j\eta_k)G_kG_\ell(\eta_\ell G_i+G_\ell \eta_i)\cr
&= 4\sum_{i, j, k, \ell (dist)} \eta_j\eta_\ell G^2_iG_jG^2_kG_\ell,\cr
R_{29} &= \sum_{i, j, k,\ell (dist)} G_iG_j^2 G_k (\eta_kG_\ell + G_k\eta_\ell)W_{\ell i}\cr
&= \sum_{i, j, k,\ell (dist)} \eta_k G_iG_j^2 G_kG_\ell W_{\ell i} +  \sum_{i, j, k,\ell (dist)} \eta_\ell G_iG_j^2 G^2_k W_{\ell i},\cr
R_{30} &= 2 \sum_{i, j, k,\ell (dist)} G_iG_j (\eta_jG_k) G_kG_\ell W_{\ell i} = 2 \sum_{i, j, k,\ell (dist)} \eta_j G_iG_j G^2_k G_\ell W_{\ell i}, \cr
R_{32} &= \sum_{i, j, k,\ell (dist)} \widetilde{\Omega}_{\ell i}   G_iG_j^2G_k^2G_\ell. 
\end{align*}
Each expression above belongs to one of the following types: 
\begin{align*}
& K_1 = \sum_{i, j, k,\ell (dist)} G_iG_j^2G_k W_{k\ell} W_{\ell i}, &&K_2 = \sum_{i, j, k,\ell (dist)} G_iG_j G_kG_\ell W_{jk}W_{\ell i},\cr
& K_3 = \sum_{i, j, k,\ell (dist)} \eta_k G_iG_j^2 G_kG_\ell W_{\ell i}, && K_4=  \sum_{i, j, k,\ell (dist)} \eta_\ell G_iG_j^2 G^2_k W_{\ell i},\cr
&K_5 = \sum_{i, j, k,\ell (dist)} \eta_i \eta_kG_iG_j^2G_kG_\ell^2, && K_5' = \sum_{i, j, k,\ell (dist)} \widetilde{\Omega}_{ik} G_iG_j^2G_kG_\ell^2, \cr
& K_6 = \sum_{i, j, k,\ell (dist)} \eta_\ell^2 G^2_iG_j^2G^2_k. 
\end{align*}
Since $|\eta_i\eta_k|\leq C\theta_i\theta_k$ and $|\widetilde{\Omega}_{ik}|\leq C\alpha\theta_i\theta_k$, 
the study of $K_5$ and $K_5'$ are similar; we thus omit the analysis of $K'_5$. We now study $K_1$-$K_6$. 

Consider $K_1$. Re-write
\[
K_1 = \frac{1}{v^2}\sum_{\substack{i,j,k,\ell (dist)\\s\neq i,t\neq j,q\neq j,m\neq k}} W_{is}W_{jt}W_{jq}W_{km}W_{k\ell}W_{\ell i}.  
\]
Note that $W_{km}W_{k\ell}W_{\ell i}W_{is}$ has four different cases: (a) $W_{k\ell}^2W^2_{\ell i}$, (b) $W^2_{k\ell}W_{\ell i}W_{is}$, (c) $W_{k\ell}W_{\ell i}W_{ik}^2$, and (d) $W_{k\ell}W_{\ell i}W_{km}W_{is}$. At the same time, $W_{jt}W_{jq}$ has two cases: (i) $W_{jk}^2$ and (ii) $W_{jt}W_{jq}$. This gives at least $4\times 2=8$ cases. Each case may have sub-cases, e.g., for $(W^2_{k\ell}W_{\ell i}W_{is})W^2_{jt}$, if $(s, t)=(j, i)$, it becomes $W^2_{k\ell}W_{\ell i}W^3_{ij}$. By direct calculations, all possible cases of the summand are as follows:
\begin{align}  \label{K1summand}
& (W^2_{k\ell}W^2_{\ell i})W^2_{jt}, \quad (W^2_{k\ell}W^2_{\ell i})(W_{jt}W_{jq}), \quad  (W^2_{k\ell}W_{\ell i}W_{is})W^2_{jt},\cr
& W^2_{k\ell}W_{\ell i}W^3_{ij}, \quad (W^2_{k\ell}W_{\ell i}W_{is})(W_{jt}W_{jq}), \quad  W^2_{k\ell}W_{\ell i}W^2_{ij}W_{jq}, \cr
 &(W_{k\ell}W_{\ell i}W_{ik}^2)W_{jt}^2, \quad (W_{k\ell}W_{\ell i}W_{ik}^2)(W_{jt}W_{jq}), \cr
 &(W_{k\ell}W_{\ell i}W_{km}W_{is})W^2_{jt}, \quad W_{k\ell}W_{\ell i}W_{km}W^3_{ij}, \cr
 &(W_{k\ell}W_{\ell i}W_{km}W_{is})(W_{jt}W_{jq}), \quad W_{k\ell}W_{\ell i}W_{km}W^2_{ij}W_{jq},\cr
& W_{k\ell}W_{\ell i}W_{kj}^2W^2_{ij}. 
\end{align}
Take the second type for example. We aim to bound $\mathbb{E}[(\sum_{i,j,k,\ell,t,q}W^2_{k\ell}W^2_{\ell i}W_{jt}W_{jq})^2]$, which is equal to
\[
\sum_{\substack{i,j,k,\ell,t,q\\i',j',k',\ell',t',q'}}\mathbb{E}[ W^2_{k\ell}W^2_{\ell i}W_{jt}W_{jq}W^2_{k'\ell'}W^2_{\ell' i'}W_{j't'}W_{j'q'} ]. 
\]
For the expectation to be nonzero, each single $W$ term has to be paired with another term. The main contribution comes from the case that $W_{j't'}W_{j'q'}=W_{jt}W_{jq}$. It implies $(j', t', q')=(j, t, q)$ or $(j', t', q')=(j, q,t)$. Then, the expression above becomes
\begin{align*}
\sum_{\substack{i,j,k,\ell,t,q\\i',k',\ell'}}\mathbb{E}[(W^2_{jt}W^2_{jq})W^2_{k\ell}W^2_{\ell i}W^2_{k'\ell'}W^2_{\ell' i'}]
& \leq C\sum_{\substack{i,j,k,\ell,t,q\\i',k',\ell'}} \theta_i\theta^2_j\theta_k\theta^2_\ell\theta_t\theta_q\theta_{i'}\theta_{k'}\theta^2_{\ell'} \cr
&\leq C\|\theta\|^6\|\theta\|_1^6. 
\end{align*}
There are a total of $9$ indices in this sum, which are $(i,j,k,\ell,t,q,i',k',\ell')$. Similarly, for each type of summand, when we bound the expectation of the square of its sum, we count how many indices appear in the ultimate sum. This number equals to twice of the total number of indices appearing in the summand, minus the total number of indices appearing in single $W$ terms. For the above example, all indices appearing in the summand are $(i,j,k,\ell,t,q)$, while indices appearing in single $W$ terms are $(j,t,q)$; so, the aforementioned number is $2\times 6 - 3=9$. If this number if $m_0$, then the expectation of the square of sum of this type is bounded by $C\|\theta\|_1^{m_0}$. We note that $K_1$ has a factor $\frac{1}{v^2}$ in front of the sum, which brings in a factor of $\frac{C}{\|\theta\|^8_1}$ in the bound. Therefore, for any type of summand with $m_0\leq 8$, the expectation of the square of its sum is $O(1)$, which is $o(\|\theta\|^8)$. As a result, among the types in \eqref{K1summand}, we only need to consider those with $m_0\geq 9$. We are left with
\[
(W^2_{k\ell}W^2_{\ell i})W^2_{jt}, \qquad (W^2_{k\ell}W^2_{\ell i})(W_{jt}W_{jq}), \qquad (W^2_{k\ell}W_{\ell i}W_{is})W^2_{jt}. 
\]
We have proved that the expectation of the square of sum of the second type of summands is bounded by $C\|\theta\|^2\|\theta\|_1^6=o(\|\theta\|^8\|\theta\|_1^8)$. For the other two types, by direct calculations,
\begin{align*}
\mathbb{E}\biggl[ \Bigl(\sum_{\substack{i,j,k,\ell (dist)\\t\neq j}}W^2_{k\ell}W^2_{\ell i}W^2_{jt}\Bigr)^2  \biggr]& \leq \sum_{\substack{i,j,k,\ell,t\\i',j',k',\ell',t'}}\mathbb{E}[ W^2_{k\ell}W^2_{\ell i}W^2_{jt}W^2_{k'\ell'}W^2_{\ell' i'}W^2_{j't'} ]\cr
&\leq \sum_{\substack{i,j,k,\ell,t\\i',j',k',\ell',t'}} \theta_i\theta_j\theta_k\theta^2_{\ell}\theta_t\theta_{i'}\theta_{j'}\theta_{k'}\theta^2_{\ell'}\theta_{t'}\cr
&\leq C\|\theta\|^4\|\theta\|_1^8 = o(\|\theta\|^8\|\theta\|_1^8),\cr
\mathbb{E}\biggl[ \Bigl(\sum_{\substack{i,j,k,\ell (dist)\\s\notin\{i,\ell\}, t\neq j,\\ (s,t)\neq (j,i)}}W^2_{k\ell}W_{\ell i}W_{is}W^2_{jt} \Bigr)^2  \biggr]& \leq \sum_{\substack{i,j,k,\ell,s,t\\j',k',t'}}\mathbb{E}[ (W^2_{\ell i}W^2_{is})W^2_{k\ell}W^2_{jt}W^2_{k'\ell}W^2_{j't'} ]\cr
&\leq C\sum_{\substack{i,j,k,\ell,s,t\\j',k',t'}}\theta^2_i\theta_j\theta_k\theta^3_{\ell}\theta_s\theta_t\theta_{j'}\theta_{k'}\theta_{t'}\cr
&\leq C\|\theta\|^2\|\theta\|_3^3\|\theta\|_1^7 = o(\|\theta\|^8\|\theta\|_1^8). 
\end{align*}
Combining the above gives
\[
\mathbb{E}[K_1^2] = o(\|\theta\|^8), \qquad \mbox{under both hypotheses}. 
\]

Consider $K_2$. Re-write
\[
K_2 = \frac{1}{v^2}\sum_{\substack{i, j, k,\ell (dist)\\s\neq i,t\neq j,q\neq k,m\neq \ell}} W_{is}W_{jt}W_{kq}W_{\ell m}W_{jk}W_{\ell i}.
\]
Note that $W_{qk}W_{kj}W_{jt}$ has three cases: (a) $W_{kj}^3$, (b) $W_{kj}^2W_{jt}$ (or $W_{qk}W_{kj}^2$), and (c) $W_{qk}W_{kj}W_{jt}$. Simiarly, $W_{m\ell}W_{\ell i}W_{is}$ has three cases: (a) $W_{\ell i}^3$, (b) $W_{\ell i}^2W_{is}$ (or $W_{m\ell}W^2_{\ell i}$), and (c) $W_{m\ell}W_{\ell i}W_{is}$. By index symmetry, this gives $3+2+1=6$ different cases. Some case may have sub-cases, due to that $(s,t)$ may equal to $(j,i)$, say. By direct calculations, all possible cases of the summand are as follows:
\begin{align*} 
&W_{kj}^3W^3_{\ell i}, \quad W^3_{kj}(W^2_{\ell i}W_{is}), \quad W^3_{kj}(W_{m\ell}W_{\ell i}W_{is}), \quad  (W^2_{kj}W_{jt})(W^2_{\ell i}W_{is}),\cr
& W^2_{kj}W_{ji}^2W^2_{\ell i}, \quad (W^2_{kj}W_{jt})(W_{m\ell}W_{\ell i}W_{is}), \quad W^2_{kj}W^2_{ji}W_{m\ell}W_{\ell i}, \cr
& (W_{qk}W_{kj}W_{jt})(W_{m\ell}W_{\ell i}W_{is}), \quad W_{qk}W_{kj}W_{ji}^2W_{m\ell}W_{\ell i}, \quad W_{kj}W_{ji}^2W^2_{k\ell}W_{\ell i}. 
\end{align*}
As in the analysis of \eqref{K1summand}, we count the effective number of indices, $m_0$, which equals to twice of the total number of indices appearing in the summand minus the total number of indices appearing in all single-$W$ terms. For the above types of summand, $m_0$ equals to $8, 8, 8, 8, 8, 8, 7, 8, 6, 4$, respectively. None is larger than $8$. We conclude that the expectation of the square of sum of each type of summand is bounded by $C\|\theta\|_1^8$. We immediately have
\[
\mathbb{E}[K_2^2] = \frac{1}{v^4}\cdot C\|\theta\|_1^8=O(1)= o(\|\theta\|^8), \qquad \mbox{under both hypotheses}. 
\]

Consider $K_3$. Re-write
\[
K_3 = \frac{1}{v^2\sqrt{v}}\sum_{\substack{i,j,k,\ell (dist)\\s\neq i,t\neq j,q\neq j,m\neq k,p\neq \ell}}\eta_k W_{is}W_{jt}W_{jq}W_{km}W_{\ell p}W_{\ell i}
\] 
Note that $W_{jt}W_{jq}W_{km}$ has four cases: (a) $W^3_{jk}$, (b) $W^2_{jk}W_{jt}$ (or $W^2_{jk}W_{jq}$), (c) $W^2_{jt}W_{km}$, and (d) $W_{jt}W_{jq}W_{km}$. At the same time, $W_{is}W_{\ell p}W_{\ell i}$ has three cases: (a) $W^3_{\ell i}$, (b) $W^2_{\ell i}W_{is}$ (or $W^2_{\ell i}W_{\ell p}$), and (c) $W_{\ell i}W_{is}W_{\ell p}$. This gives $4\times 3=12$ different cases. Each case may have sub-cases. For example, in the case of $\eta_k(W_{jk}^2W_{jt})(W^2_{\ell i}W_{is})$, if $(s, t)=(j,i)$, it becomes $\eta_kW^2_{jk}W^2_{ji}W^2_{\ell i}$. By direct calculations, we obtain all possible cases of summands as follows:
\begin{align*}
& \eta_k W_{jk}^3W_{\ell i}^3, \quad \eta_k W_{jk}^3(W^2_{\ell i}W_{is}), \quad \eta_k W^3_{jk}(W_{\ell i}W_{is}W_{\ell p}), \quad \eta_k (W_{jk}^2W_{jt})W_{\ell i}^3, \cr
&\eta_k (W_{jk}^2W_{jt})(W^2_{\ell i}W_{is}), \quad \eta_k W^2_{jk}W^2_{ji}W^2_{\ell i}, \quad \eta_k (W_{jk}^2W_{jt})(W_{\ell i}W_{is}W_{\ell p}),\cr
& \eta_k W_{jk}^2W^2_{ji}W_{\ell i}W_{\ell p}, \quad \eta_k (W_{jt}^2W_{km})W_{\ell i}^3, \quad \eta_k (W_{jt}^2W_{km})(W^2_{\ell i}W_{is}),\quad \eta_k W^2_{jt}W^2_{ki}W^2_{\ell i},\cr
&\eta_k (W_{jt}^2W_{km})(W_{\ell i}W_{is}W_{\ell p}),\quad \eta_k W_{jt}^2W^2_{ki}W_{\ell i}W_{\ell p},\quad 
\eta_k (W_{jt}W_{jq}W_{km})W_{\ell i}^3,\cr
&\eta_k (W_{jt}W_{jq}W_{km})(W^2_{\ell i}W_{is}),\quad \eta_k W_{jt}W^2_{ji}W_{km}W^2_{\ell i},\quad \eta_k W_{jt}W_{jq}W^2_{ki}W^2_{\ell i},\cr
& \eta_k (W_{jt}W_{jq}W_{km})(W_{\ell i}W_{is}W_{\ell p}), \quad  \eta_k W_{jt}W^2_{ji}W_{km}W_{\ell i}W_{\ell p}, \quad \eta_k W_{jt}W_{jq}W^2_{ki}W_{\ell i}W_{\ell p}. 
\end{align*}
Same as before, let $m_0$ be the effective number of indices for each type of summand, which equals to twice of number of distinct indices appearing in the summand minus the number of distinct indices appearing in single-$W$ terms (see \eqref{K1summand} and text therein). By direct calculations, $m_0\leq 10$ for all types above. It follows that, for each type of summand, the expectation of the square of their sums is bounded by
\[
\frac{1}{(v\sqrt{v})^2}\cdot C\|\theta\|_1^{m_0}\leq C \|\theta\|_1^{m_0-10}=O(1)=o(\|\theta\|^8). 
\]
We immediately have
\[
\mathbb{E}[K_3^2] =  o(\|\theta\|^8), \qquad \mbox{under both hypotheses}. 
\]

%
%

Consider $K_4$. Re-write
\[
K_4 = \frac{1}{v^2\sqrt{v}}\sum_{\substack{i,j,k,\ell (dist)\\s,t,q,m,p}}\eta_\ell W_{is}W_{jt}W_{jq}W_{km}W_{kp}W_{\ell i}. 
\]
Note that $W_{is}W_{\ell i}$ has two cases: (a) $W_{\ell i}^2$ and (b) $W_{\ell i}W_{is}$. Moreover, there are a total of six cases for $W_{jt}W_{jq}W_{km}W_{kp}$: 
(a) $W^4_{jk}$, (b) $W^3_{jk}W_{jt}$, (c) $W^2_{jk}W_{jt}W_{km}$, (d) $W^2_{jt}W^2_{km}$, (e) $W_{jt}W_{jq}W^2_{km}$, and (f) $W_{jt}W_{jq}W_{km}W_{kp}$. It gives $2\times 6=12$ different cases. Each case may have some sub-cases. It turns out all different types of summand are as follows:
\begin{align*}
& \eta_\ell W^2_{\ell i} W^4_{jk}, \quad \eta_\ell W^2_{\ell i}(W^3_{jk}W_{jt}),  \quad \eta_\ell W^2_{\ell i}(W^2_{jk}W_{jt}W_{km}), \quad  \eta_\ell W^2_{\ell i}(W^2_{jt}W^2_{km}),\cr
& \eta_\ell W^2_{\ell i}(W_{jt}W_{jq}W^2_{km}),\quad \eta_\ell W^2_{\ell i}(W_{jt}W_{jq}W_{km}W_{kp}), \quad \eta_\ell (W_{\ell i}W_{is})W^4_{jk}, \cr
&\eta_\ell (W_{\ell i}W_{is})(W^3_{jk}W_{jt}),\quad \eta_\ell W_{\ell i}W^3_{jk}W^2_{ji},\quad
  \eta_\ell (W_{\ell i}W_{is})(W^2_{jk}W_{jt}W_{km}),\cr
  &\eta_\ell W_{\ell i}W^2_{jk}W^2_{ji}W_{km},\quad \eta_\ell (W_{\ell i}W_{is})(W^2_{jt}W^2_{km}), \quad \eta_\ell W_{\ell i}W^3_{ij}W^2_{km},\cr
  &\eta_\ell  (W_{\ell i}W_{is})(W_{jt}W_{jq}W^2_{km}),\quad \eta_\ell  W_{\ell i}W^2_{ij}W_{jq}W^2_{km}, \quad \eta_\ell W_{\ell i}W_{jt}W_{jq}W^3_{ki}, \cr
  & \eta_\ell (W_{\ell i}W_{is})(W_{jt}W_{jq}W_{km}W_{kp}),\quad \eta_\ell W_{\ell i}W^2_{ij}W_{jq}W_{km}W_{kp}. 
\end{align*}
Same as before, for each type, let $m_0$ be the effective number of indices. It suffices to focus on cases where $m_0\geq 11$. We are left with 
\[
\eta_\ell W^2_{\ell i}(W^2_{jt}W_{km}^2), \qquad \eta_\ell W^2_{\ell i}(W_{jt}W_{jq}W^2_{km}), \qquad \eta_\ell(W_{\ell i}W_{is})(W^2_{jt}W_{km}^2). 
\] 
By direct calculations,
\begin{align*}
\mathbb{E}\biggl[ \Bigl(\sum_{\substack{i,j,k,\ell (dist)\\t\neq j, m\neq k}}\eta_\ell W^2_{\ell i}W^2_{jt}W^2_{km} \Bigr)\biggr]&\leq \sum_{\substack{i,j,k,\ell,t,m\\i',j',k',\ell',t',m'}} \eta_{\ell}\eta_{\ell'}\, \mathbb{E}[W^2_{\ell i}W^2_{jt}W^2_{km}W^2_{\ell' i'}W^2_{j't'}W^2_{k'm'}]\cr
&\leq C\sum_{\substack{i,j,k,\ell,t,m\\i',j',k',\ell',t',m'}}\theta_{i}\theta_j\theta_k\theta^2_\ell \theta_t\theta_m\theta_{i'}\theta_{j'}\theta_{k'}\theta^2_{\ell'} \theta_{t'}\theta_{m'}\cr
&\leq C\|\theta\|^4\|\theta\|_1^{10} = o(\|\theta\|^8\|\theta\|_1^{10}),\cr
\mathbb{E}\biggl[ \Bigl(\sum_{\substack{i,j,k,\ell (dist)\\t\neq j, q\neq j, m\neq k\\t\neq q}}\eta_\ell W^2_{\ell i}W_{jt}W_{jq}W^2_{km} \Bigr)\biggr] & \leq \sum_{\substack{i,j,k,\ell,t,q,m\\i',k',\ell',m'}}\eta_{\ell}\eta_{\ell'}\,\mathbb{E}[(W^2_{jt}W^2_{jq})W^2_{\ell i}W^2_{km}W^2_{\ell'i'}W^2_{k'm'}]\cr
&\leq C\sum_{\substack{i,j,k,\ell,t,q,m\\i',k',\ell',m'}}\theta_i\theta^2_j\theta_k\theta^2_\ell\theta_t\theta_q\theta_m\theta_{i'}\theta_{k'}\theta^2_{\ell'}\theta_{m'}\cr
&\leq C\|\theta\|^6\|\theta\|_1^{8} = o(\|\theta\|^8\|\theta\|_1^{10}),\cr
\mathbb{E}\biggl[ \Bigl(\sum_{\substack{i,j,k,\ell (dist)\\s\neq i, t\neq j, m\neq k\\(s,t)\neq (j,i), (s,m)\neq (k,i)}}\eta_\ell W_{\ell i}W_{is}W^2_{jt}W_{km}^2 \Bigr)\biggr] &\leq C\sum_{\substack{i,j,k,\ell,s,t,m\\j',k',t',m'}}\eta^2_{\ell}\, \mathbb{E}[(W^2_{\ell i}W^2_{is})W^2_{jt}W^2_{km}W^2_{j't'}W^2_{k'm'}]\cr
&\leq C\sum_{\substack{i,j,k,\ell,s,t,m\\j',k',t',m'}}\theta^2_i\theta_j\theta_k\theta^3_\ell\theta_s\theta_t\theta_m\theta_{j'}\theta_{k'}\theta_{t'}\theta_{m'}\cr
&\leq C\|\theta\|^2\|\theta\|_3^3\|\theta\|_1^9 = o(\|\theta\|^8\|\theta\|_1^{10}). 
\end{align*}
It follows that
\[
\mathbb{E}[K_4^2]\leq \frac{1}{(v^2\sqrt{v})^2}\cdot o(\|\theta\|^8\|\theta\|_1^{10})=o(\|\theta\|^8), \qquad \mbox{under both hypotheses}. 
\]

Consider $K_5$-$K_6$. To save space, we only present the proof for the case of $\|\theta\|\gg[\log(n)]^{3/2}$. When $1\ll \|\theta\|\leq C[\log(n)]^{3/2}$, we can bound $\mathbb{E}[K_5^2]$ and $\mathbb{E}[K_6^2]$ in the same way as in the study of $J_1$-$J_8$, so the proof is omitted.  Let $E$ be the event defined in \eqref{eventE}. We have argued that it suffices to focus on the event $E$. On this event, 
$|G_i|\leq C\sqrt{\theta_i\|\theta\|_1\log(n)/v}$. It follows that
\begin{align*}
|K_5| &\leq C\sum_{i,j,k,\ell} (\theta_i\theta_k)\frac{\sqrt{\theta_i\theta_j^2\theta_k\theta_\ell^2}\|\theta\|_1^3[\log(n)]^3}{v^3}\cr
&\leq \frac{C[\log(n)]^3}{\|\theta\|_1^3}\Bigl(\sum_{i}\theta_i^{3/2}\Bigr)\Bigl(\sum_j\theta_j\Bigr)\Bigl(\sum_k \theta_k^{3/2}\Bigr)\Bigl(\sum_\ell\theta_\ell\Bigr) \cr
&\leq \frac{C[\log(n)]^3}{\|\theta\|_1^3}\bigl(\|\theta\|\sqrt{\|\theta\|_1}\bigr)^2\|\theta\|_1^2\cr
&\leq C[\log(n)]^3\|\theta\|^2,
\end{align*}
where we have used the Cauchy-Schwarz inequality $(\sum_{i}\theta_i^{3/2})\leq \|\theta\|\sqrt{\|\theta\|_1}$. Similarly, 
\begin{align*}
|K_6| &\leq C\sum_{i,j,k,\ell} \theta_{\ell}^2\cdot \frac{\theta_i\theta_j\theta_k\|\theta\|_1^3[\log(n)]^3}{v^3}\cr
&\leq \frac{C[\log(n)]^3}{\|\theta\|_1^3}\sum_{i,j,k,\ell}\theta_i\theta_j\theta_k\theta_\ell^2 \cr
&\leq C[\log(n)]^3\|\theta\|^2. 
\end{align*}
When $\|\theta\|\gg[\log(n)]^{3/2}$, both right hand sides are $o(\|\theta\|^4)$. We immediately have
\[
\max\bigl\{\mathbb{E}[K_5^2], \;\mathbb{E}[K_6^2]\bigr\} = o(\|\theta\|^8). 
\]

We have proved: Each $R_k$ with $N^*_W=6$ satisfies $\mathbb{E}[R_k^2]=o(\|\theta\|^8)$. This is sufficient to guarantee \eqref{remainder1-Nullgoal2}-\eqref{remainder1-Altgoal2} for $X=R_k$.

\paragraph{Analysis of terms with $N^*_W\geq 7$} There are $3$ such terms, $R_{31}$, $R_{33}$ and $R_{34}$. Consider $R_{31}$. By definition,
\[
R_{31} = \sum_{i,j,k,\ell (dist)}G_iG_j^2G_k^2G_\ell W_{\ell i} = \frac{1}{v^3}\sum_{\substack{i,j,k,\ell (dist)\\s\neq i,t\neq j,q\neq j,\\m\neq k,p\neq k,y\neq \ell}}W_{is}W_{jt}W_{jq}W_{km}W_{kp}W_{\ell y}W_{\ell i}. 
\]
We note that $W_{\ell i}W_{is}W_{\ell y}$ has three cases: (a) $W^3_{\ell i}$, (b) $W^2_{\ell i}W_{is}$, and (c) $W_{\ell i}W_{is}W_{\ell y}$. Moreover, $W_{jt}W_{jq}W_{km}W_{kp}$ has six cases: (a) $W^4_{jk}$, (b) $W^3_{jk}W_{jt}$, (c) $W^2_{jk}W_{jt}W_{km}$, (d) $W^2_{jt}W^2_{km}$, (e) $W_{jt}W_{jq}W^2_{km}$, and (f) $W_{jt}W_{jq}W_{km}W_{kp}$. This gives $3\times 6=18$ different cases. Since each case may have sub-cases, we end up with the following different types:
\begin{align*}
& W^3_{\ell i}W^4_{jk}, \quad W^3_{\ell i}(W^3_{jk}W_{jt}), \quad W^3_{\ell i}(W^2_{jk}W_{jt}W_{km}), \quad W^3_{\ell i}(W^2_{jt}W^2_{km}),\cr
& W^3_{\ell i}(W_{jt}W_{jq}W^2_{km}),\quad  W^3_{\ell i}(W_{jt}W_{jq}W_{km}W_{kp}),\quad
(W^2_{\ell i}W_{is})W^4_{jk}, \cr
& (W^2_{\ell i}W_{is})(W^3_{jk}W_{jt}), \quad W^2_{\ell i}W^3_{jk}W^2_{ji},\quad
\quad (W^2_{\ell i}W_{is})(W^2_{jk}W_{jt}W_{km}),\cr
& W^2_{\ell i}W^2_{jk}W^2_{ji}W_{km}, \quad (W^2_{\ell i}W_{is})(W^2_{jt}W^2_{km}), \quad W^2_{\ell i}W^3_{ij}W^2_{km}, \cr
& (W^2_{\ell i}W_{is})(W_{jt}W_{jq}W^2_{km}),\quad W^2_{\ell i}W^2_{ij}W_{jq}W^2_{km}, \quad 
 W^2_{\ell i}W_{jt}W_{jq}W^3_{ki},\cr
& (W^2_{\ell i}W_{is})(W_{jt}W_{jq}W_{km}W_{kp}),\quad W^2_{\ell i}W^2_{ij}W_{jq}W_{km}W_{kp}, \cr
& (W_{\ell i}W_{is}W_{\ell y})W^4_{jk},\quad (W_{\ell i}W_{is}W_{\ell y})(W^3_{jk}W_{jt}),\quad  W_{\ell i}W_{\ell y}W^3_{jk}W^2_{ji},\cr
 & (W_{\ell i}W_{is}W_{\ell y})(W^2_{jk}W_{jt}W_{km}),\quad W_{\ell i}W_{\ell y}W^2_{jk}W^2_{ji}W_{km}, \quad W_{\ell i}W^2_{jk}W^2_{ji}W^2_{k\ell}, \cr
 &(W_{\ell i}W_{is}W_{\ell y})(W^2_{jt}W^2_{km}),\quad W_{\ell i}W_{\ell y}W^3_{ji}W^2_{km},\quad W_{\ell i}W_{ji}^3W^3_{k\ell}, \cr
&(W_{\ell i}W_{is}W_{\ell y})(W_{jt}W_{jq}W^2_{km}), \quad W_{\ell i}W_{\ell y}W^2_{ji}W_{jq}W^2_{km},
 \quad W_{\ell i}W_{\ell y}W_{jt}W_{jq}W^3_{ki},\cr
 &W_{\ell i}W^2_{ji}W_{jq}W^3_{ki},\quad  (W_{\ell i}W_{is}W_{\ell y})(W_{jt}W_{jq}W_{km}W_{kp}),\cr
 & W_{\ell i}W_{\ell y}W^2_{ji}W_{jq}W_{km}W_{kp}, \quad W_{\ell i}W^2_{ji}W_{jq}W^2_{k\ell}W_{kp}. 
\end{align*}
For each type, we count $m_0$, the effective number of indices. It equals to twice of the number of distinct indices in the summand, minus the number of distinct indices appearing in all single-$W$ terms. It turns out that $m_0\leq 12$ for all types above. By similar arguments as in \eqref{K1summand}, we conclude that
\[
\mathbb{E}[R_{31}^2]\leq \frac{1}{v^6}\cdot C\|\theta\|_1^{m_0}\leq C\|\theta\|_1^{m_0-12}=O(1)=o(\|\theta\|^8), \qquad \mbox{under both hypotheses}. 
\]

Consider $R_{33}$-$R_{34}$. We only give the proof when $\|\theta\|^6\gg[\log(n)]^7$, as it is much simpler. In the case of $1\ll\|\theta\|^6\leq C[\log(n)]^7$, we can follow similar steps above to obtain desired bounds, where details are omitted. On the event $E$ (see \eqref{eventE} for definition), 
\begin{align*}
|R_{33}|& \leq \sum_{i,j,k,\ell}|\eta_\ell ||G^2_iG_j^2G_k^2G_\ell|\cr
&\leq C\sum_{i, j, k, \ell} \theta_\ell \frac{\sqrt{\theta_i^2\theta_j^2\theta_k^2\theta_\ell \|\theta\|_1^7[\log(n)]^7}}{(\sqrt{v})^7}\cr
&\leq \frac{C[\log(n)]^{7/2}}{\sqrt{\|\theta\|_1^7}} \Bigl(\sum_{i}\theta_i\Bigr)\Bigl(\sum_{j}\theta_j\Bigr)\Bigl(\sum_k\theta_k\Bigr)\Bigl(\sum_{\ell}\theta_\ell^{3/2}\Bigr)\cr
&\leq \frac{C[\log(n)]^{7/2}}{\sqrt{\|\theta\|_1^7}}\cdot \|\theta\|_1^3\bigl(\|\theta\|\sqrt{\|\theta\|_1}\bigr)\cr
&\leq C[\log(n)]^{7/2}\|\theta\|,
\end{align*} 
where we have used the Cauchy-Schwarz inequality $\sum_\ell \theta_\ell^{3/2}\leq \|\theta\|\sqrt{\|\theta\|_1}$ in the second last line. When $\|\theta\|^6\gg [\log(n)]^7$, the right hand side is $o(\|\theta\|^4)$. Similarly, 
\begin{align*}
|R_{34}|& \leq \sum_{i,j,k,\ell}|G^2_iG_j^2G_k^2G^2_\ell|\cr
&\leq C\sum_{i, j, k, \ell} \frac{\theta_i\theta_j\theta_k\theta_\ell \|\theta\|_1^4[\log(n)]^4}{v^4}\cr
& \leq C[\log(n)]^4. 
\end{align*} 
When $\|\theta\|^6\gg [\log(n)]^7$, the right hand side is $o(\|\theta\|^4)$. As we have argued in \eqref{eventE}, the event $E^c$ has a negligible effect. It follows that
\[
\max\bigl\{\mathbb{E}[R_{31}^2],\; \mathbb{E}[R^2_{33}], \; \mathbb{E}[R^2_{34}]\bigr\} = o(\|\theta\|^8), \qquad \mbox{under both hypotheses}. 
\]
This is sufficient to guarantee \eqref{remainder1-Nullgoal2}-\eqref{remainder1-Altgoal2} for $R_k$.

We have analyzed all $34$ terms in Table~\ref{tb:remainder}. The proof is now complete.

\subsubsection{Proof of Lemma~\ref{lem:remainder2}}  \label{subsec:remainder2}
Consider an arbitrary post-expansion sum of the form
\beq \label{lastlemma-1}
\sum_{i,j,k,\ell (dist)} a_{ij}b_{jk}c_{k\ell}d_{\ell i}, \qquad \mbox{where}\quad a,b,c,d\in \{\widetilde{\Omega}, W,\delta, \tilde{r}, \epsilon\}.
\eeq
Let $(N_{\widetilde{\Omega}}, N_W, N_\delta, N_{\tilde{r}}, N_{\epsilon})$ be the number of each type in the product, where these numbers have to satisfy $N_{\widetilde{\Omega}}+N_W+N_\delta+N_{\tilde{r}}+N_{\epsilon}=4$. 
As discussed in Section~\ref{suppC3}, $(Q_n-Q_n^*)$ equals to the sum of all post-expansion sums such that $N_{\epsilon}>0$. Recall that 
\[
\epsilon_{ij} = (\eta_i^*\eta_j^*-\eta_i\eta_j)+(1-\frac{v}{V})\eta_i\eta_j-(1-\frac{v}{V})\delta_{ij}. 
\]
Define 
\[
\epsilon_{ij}^{(1)} = \eta_i^*\eta_j^*-\eta_i\eta_j, \quad \epsilon_{ij}^{(2)} =  (1-\frac{v}{V})\eta_i\eta_j, \quad \epsilon^{(3)}_{ij}= -(1-\frac{v}{V})\delta_{ij}. 
\]
Then, $\epsilon_{ij}=\epsilon_{ij}^{(1)}+\epsilon_{ij}^{(2)}+\epsilon_{ij}^{(3)}$. It follows that each post-expansion sum of the form \eqref{lastlemma-1} can be further expanded as the sum of terms like 
\beq  \label{lastlemma-2}
\sum_{i,j,k,\ell (dist)} a_{ij}b_{jk}c_{k\ell}d_{\ell i}, \quad \mbox{where}\;\; a,b,c,d\in \{\widetilde{\Omega}, W,\delta, \tilde{r}, \epsilon^{(1)}, \epsilon^{(2)}, \epsilon^{(3)}\}.  
\eeq
Let $(N_{\widetilde{\Omega}}, N_W, N_\delta, N_{\tilde{r}})$ have the same meaning as before, and let $N^{(m)}_{\epsilon}$ be the number of $\epsilon^{(m)}$ term in the product, for $m\in\{1,2,3\}$. These numbers have to satisfy 
$N_{\widetilde{\Omega}}+N_W+N_\delta+N_{\tilde{r}}+N^{(1)}_{\epsilon}+N^{(2)}_{\epsilon}+N^{(3)}_{\epsilon}=4$. Now, $(Q_n-Q_n^*)$ equals to the sum of all post-expansion sums of the form \eqref{lastlemma-2} with
\beq \label{lastlemma-3}
N^{(1)}_{\epsilon}+N^{(2)}_\epsilon+ N^{(3)}_{\epsilon} \geq 1. 
\eeq
Fix such a post-expansion sum and denote it by $Y$. We shall bound $|\mathbb{E}[Y]|$ and $\mathrm{Var}(Y)$. 

We need some preparation. First, we derive a bound for $|\epsilon^{(1)}_{ij}|$. By definition, $\eta_i =(1/\sqrt{v})\sum_{j\neq i}\Omega_{ij}$ and $\eta^*_i = (1/\sqrt{v_0})\sum_{j}\Omega_{ij}$. It follows that
\[
\eta^*_i = \frac{\sqrt{v}}{\sqrt{v_0}}\eta_i + \frac{1}{\sqrt{v_0}}\Omega_{ii}. 
\]
We then have
\[
\eta_i^*\eta_j^* = \frac{v}{v_0}\eta_i\eta_j + \frac{\sqrt{v}}{v_0}(\eta_i\Omega_{jj}+\eta_j\Omega_{ii}) + \frac{1}{v_0}\Omega_{ii}\Omega_{jj}.  
\]
Note that $v=\sum_{i\neq j}\Omega_{ij}$ and $v_0=\sum_{ij}\Omega_{ij}\asymp \|\theta\|_1^2$. It follows that $v_0-v=\sum_i\Omega_{ii}\leq \sum_i\theta_i^2\leq \|\theta\|^2$. Therefore,
\begin{align*}
|\eta_i^*\eta_j^* - \eta_i\eta_j| &\leq \Bigl|1-\frac{v}{v_0} \Bigr|\eta_i\eta_j + \frac{\sqrt{v}}{v_0}(\eta_i\Omega_{jj}+\eta_j\Omega_{ii}) + \frac{1}{v_0}\Omega_{ii}\Omega_{jj}\cr
&\leq \frac{C\|\theta\|^2}{\|\theta\|_1^2}\cdot \theta_i\theta_j + \frac{C}{\|\theta\|_1}(\theta_i\theta_j^2+\theta_j\theta_i^2) + \frac{C}{\|\theta\|_1^2}\cdot \theta_i^2\theta_j^2\cr
&\leq C\theta_i\theta_j\cdot\Bigl( \frac{\|\theta\|^2}{\|\theta\|_1^2} + \frac{\theta_i+\theta_j}{\|\theta\|_1} +\frac{\theta_i\theta_j}{\|\theta\|_1^2} \Bigr).
\end{align*}
Since $\|\theta\|^2\leq \theta_{\max}\|\theta\|_1$, the term in the brackets is bounded by $C\theta_{\max}/\|\theta\|_1$. We thus have
\beq  \label{lastlemma-prep}
|\epsilon_{ij}^{(1)}| \leq \frac{C\theta_{\max}}{\|\theta\|_1}\cdot \theta_i\theta_j, \qquad \mbox{for all }1\leq i\neq j\leq n. 
\eeq

Second, in Lemmas~\ref{lem:IdealSgnQ-null}-\ref{lem:remainder1}, we have studied all post-expansion sums of the form 
\[
Z\equiv \sum_{i,j,k,\ell (dist)} a_{ij} b_{jk} c_{k\ell} d_{\ell i}, \qquad \mbox{where}\quad a, b, c, d\in \{\widetilde{\Omega}, W, \delta, \tilde{r}\},
\]
where $(N_{\widetilde{\Omega}}, N_W, N_\delta, N_{\tilde{r}})$ are the numbers of each type in the product. We hope to take advantage of these results. Using the proved bounds for $|\mathbb{E}[Z]|$ and $\mathrm{Var}(Z)$, we can get
\beq  \label{lastlemma-4}
\mathbb{E}[Z^2]\leq C(\alpha^2)^{N_{\widetilde{\Omega}}}\cdot f(\theta; N_{\widetilde{\Omega}}, N_W, N_\delta, N_{\tilde{r}}), 
\eeq
where $\alpha=|\lambda_2|/\lambda_1$ and $f(\theta; m_1, m_2, m_3, m_4)$ is a function of $\theta$ whose form is determined by $(m_1,m_2,m_3,m_4)$. For example,
\[
\begin{cases}
f(\theta; 0, 4, 0, 0)\, \text{=}\, \|\theta\|^8, &\mbox{by claims of $X_1$ in Lemmas~\ref{lem:IdealSgnQ-null}\&\ref{lem:IdealSgnQ-alt}};\\
f(\theta; 4, 0, 0, 0)\, \text{=}\, \|\theta\|^{16}, &\mbox{by claims of $X_6$ in Lemma~\ref{lem:IdealSgnQ-alt}};\\
f(\theta; 3, 1, 0, 0)\, \text{=}\,\|\theta\|^8\|\theta\|_3^6, &\mbox{by claims of $X_5$ in Lemma~\ref{lem:IdealSgnQ-alt}};\\
f(\theta; 1, 2, 1, 0) \, \text{=}\, \|\theta\|^4\|\theta\|_3^6, &\mbox{by claims of $Y_2$, $Y_3$ in Lemma~\ref{lem:ProxySgnQ(a)-alt}};\\
f(\theta; 1,1,1,1) \, \text{=}\, \|\theta\|^8, &\mbox{by claims of $R_{9}$-$R_{11}$ in the proof of Lemma~\ref{lem:remainder1}}. 
\end{cases}
\]
If there are more than one post-expansion sum that corresponds to the same $(N_{\widetilde{\Omega}}, N_W, N_\delta, N_{\tilde{r}})$, we use the largest bound to define $f(\theta; N_{\widetilde{\Omega}}, N_W, N_\delta, N_{\tilde{r}})$. Thanks to previous lemmas, we have known the function $f(\theta; m_1, m_2, m_3,m_4)$ for all possible $(m_1, m_2, m_3, m_4)$.

We now show the claim. 
Recall that $Y$ is the post-expansion sum in \eqref{lastlemma-2}. The key is to prove the following argument: For any sequence $x_n$ such that $\sqrt{\log(\|\theta\|_1)}\ll x_n\ll \|\theta\|_1$,
\begin{align} \label{lastlemma-key}
\mathbb{E}[Y^2] & \leq C(\alpha^2)^{N_{\widetilde{\Omega}}}\times \Bigl( \frac{\theta^2_{\max}}{\|\theta\|^2_1}\Bigr)^{N^{(1)}_\epsilon} \times \Bigl(\frac{x^2_n}{\|\theta\|^2_1}\Bigr)^{N^{(2)}_\epsilon+N_\epsilon^{(3)}}\cr
&\qquad \times  f(\theta; m_1, m_2, m_3, m_4)\bigg|_{\substack{m_1= N_{\widetilde{\Omega}}+N_{\epsilon}^{(1)}+N_{\epsilon}^{(2)},\;\; m_2= N_W,\\ m_3= N_\delta + N_{\epsilon}^{(3)},\;\; m_4 = N_{\tilde{r}},}} 
\end{align}
where $(N_{\widetilde{\Omega}}, N_W, N_\delta, N_{\tilde{r}}, N_{\epsilon}^{(1)}, N_{\epsilon}^{(2)}, N_{\epsilon}^{(3)})$ are the same as in \eqref{lastlemma-2}-\eqref{lastlemma-3}, and $f(\theta; m_1,m_2,m_3,m_4)$ is the known function in \eqref{lastlemma-4}. 

We prove \eqref{lastlemma-key}. Let $D$ be the event 
\[
D = \{ |V-v|\leq \|\theta\|_1x_n\}. 
\]
In Lemma~\ref{lem:event}, we have proved $\mathbb{E}[(Q_n-Q_n^*)^2\cdot I_{D^c}] = o(1)$. By similar proof, we can show: when $|Y|$ is bounded by a polynomial of $V$ and $\|\theta\|_1$ (which is always the case here), 
\[
\mathbb{E}[Y^2\cdot I_{D^c}]=o(1). 
\]
It follows that
\beq \label{lastlemma-5}
\mathbb{E}[Y^2] \leq \mathbb{E}[Y^2\cdot I_D] + o(1). 
\eeq
We then bound $\mathbb{E}[Y^2\cdot I_D]$. 
In the definition of $Y$, each $\epsilon^{(2)}$ term introduces a factor of $(1-\frac{v}{V})$, and each $\epsilon^{(3)}$ term introduces a factor of $-(1-\frac{v}{V})$. We bring all these factors to the front and re-write the post-expansion sum as
\[
Y=(-1)^{N_\epsilon^{(3)}}\Bigl(1 - \frac{v}{V}\Bigr)^{N^{(2)}_{\epsilon}+N^{(3)}_{\epsilon}} X, \qquad X\equiv \sum_{i,j,k,\ell (dist)} a_{ij}b_{jk}c_{k\ell}d_{\ell i}. 
\]
After the factor $(1-\frac{v}{V})$ is removed, $\epsilon^{(2)}$ becomes $\eta_i\eta_j$; similarly, $\epsilon^{(3)}$ becomes $\delta_{ij}$. Therefore, in the expression of $X$, 
\beq \label{lastlemma-6}
\begin{cases}
a_{ij}, b_{ij}, c_{ij}, d_{ij} \in \{\widetilde{\Omega}_{ij}, W_{ij}, \delta_{ij}, \tilde{r}_{ij}, \epsilon_{ij}^{(1)}, \eta_i\eta_j\}, \\
\mbox{number of $\eta_i\eta_j$ in the product is $N^{(2)}_{\epsilon}$},\\
\mbox{number of $\delta_{ij}$ in the product is $N_{\delta}+N_{\epsilon}^{(3)}$},\\
\mbox{number of any other term in the product is same as before}. 
\end{cases}
\eeq
On the event $D$, $|1-\frac{v}{V}|\leq \frac{x_n\|\theta\|_1}{C\|\theta\|_1^2} = O(\frac{x_n}{\|\theta\|_1})$. Hence, 
\[
|Y|\leq C\Bigl(\frac{x_n}{\|\theta\|_1}\Bigr)^{N^{(2)}_\epsilon+N^{(3)}_{\epsilon}}|X|, \qquad \mbox{on the event $D$}.
\]
It follows that  
\beq  \label{lastlemma-7}
\mathbb{E}[Y^2\cdot I_D]\leq C\Bigl(\frac{x^2_n}{\|\theta\|_1^2}\Bigr)^{N^{(2)}_{\epsilon}+N^{(3)}_{\epsilon}}\cdot \mathbb{E}[X^2]. 
\eeq
To bound $\mathbb{E}[X^2]$, we compare $X$ and $Z$. In obtaining \eqref{lastlemma-4}, the only property of $\widetilde{\Omega}$ we have used is
\[
|\widetilde{\Omega}_{ij}|\leq \alpha\cdot C\theta_i\theta_j. 
\]   
In comparison, in the expression of $X$, we have (by \eqref{lastlemma-prep} and \eqref{eta-bound})
\beq \label{lastlemma-nonrandom}
|\widetilde{\Omega}_{ij}|\leq \alpha\cdot C\theta_i\theta_j, \qquad|\epsilon^{(1)}_{ij}|\leq \frac{\theta_{\max}}{\|\theta\|_1}\cdot C\theta_i\theta_j, \qquad |\eta_i\eta_j|\leq C\theta_i\theta_j.
\eeq
If we consider $(\alpha^{N_{\widetilde{\Omega}}}\cdot (\frac{\theta_{\max}}{\|\theta\|_1})^{N^{(1)}_{\epsilon}}\cdot 1^{N^{(2)}_{\epsilon}})^{-1}X$ and $(\alpha^{N_{\widetilde{\Omega}}})^{-1}Z$, we can derive the same upper bound for the second moment of both variables, except that the effective $N_\delta$ in $X$ should be $N_\delta+N^{(3)}_{\epsilon}$ and the effective $N_{\widetilde{\Omega}}$ in $X$ should be $N_{\widetilde{\Omega}} + N^{(1)}_{\epsilon} + N^{(2)}_{\epsilon}$. It follows that
\begin{align} \label{lastlemma-8}
\mathbb{E}[X^2] & \leq C(\alpha^2)^{N_{\widetilde{\Omega}}}\times \Bigl( \frac{\theta^2_{\max}}{\|\theta\|^2_1}\Bigr)^{N^{(1)}_\epsilon} \cr
&\qquad \times  f(\theta; m_1, m_2, m_3, m_4)\bigg|_{\substack{m_1= N_{\widetilde{\Omega}}+N_{\epsilon}^{(1)}+N_{\epsilon}^{(2)},\;\; m_2= N_W,\\ m_3= N_\delta + N_{\epsilon}^{(3)},\;\; m_4 = N_{\tilde{r}}.}} 
\end{align}
We plug \eqref{lastlemma-8} into \eqref{lastlemma-7}, and then plug it into \eqref{lastlemma-5}. It gives \eqref{lastlemma-key}.

Next, we use \eqref{lastlemma-key} to prove the claims of this lemma. Under our assumption, we can choose a sequence $x_n$ such that $\sqrt{\log(\|\theta\|_1)}\ll x_n\ll \|\theta\|_1/\|\theta\|^2$. Also, note that $\|\theta\|_1\geq \theta_{\max}^{-1}\|\theta\|^2\gg \|\theta\|^2$. Then,  
\beq \label{lastlemma-2order}
\frac{\theta_{\max}}{\|\theta\|_1} = o(\|\theta\|^{-2}), \qquad \frac{x_n}{\|\theta\|_1} = o(\|\theta\|^{-2}). 
\eeq
As a result, since $N_{\epsilon}^{(1)}+N_{\epsilon}^{(2)}+N_{\epsilon}^{(3)}\geq 1$,  \eqref{lastlemma-key} implies
\beq \label{lastlemma-9}
\mathbb{E}[Y^2] = o(\|\theta\|^{-4})\cdot f(\theta; m_1, m_2, m_3, m_4),
\eeq
for $m_1= N_{\widetilde{\Omega}}+N_{\epsilon}^{(1)}+N_{\epsilon}^{(2)}$, $m_2= N_W$, $m_3= N_\delta + N_{\epsilon}^{(3)}$ and $m_4 = N_{\tilde{r}}$. We then extract $f(\theta; m_1, m_2, m_3, m_4)$ from previous lemmas. Recall the following facts:
\begin{itemize}
\item Under the null hypothesis, for any previously analyzed post-expansion sum $Z$, $|\mathbb{E}[Z]|\leq C\|\theta\|^4$ and $\mathrm{Var}(Z)\leq C\|\theta\|^8$.
\item Under the alternative hypothesis, except $\sum_{i,j,k,\ell (dist)}\widetilde{\Omega}_{ij}\widetilde{\Omega}_{jk}\widetilde{\Omega}_{k\ell}\widetilde{\Omega}_{\ell i}$, for all previously analyzed post-expansion sum $Z$ , $|\mathbb{E}[Z]|\leq C\alpha^2\|\theta\|^6$ and $\mathrm{Var}(Z)\leq C\|\theta\|^8+C\alpha^6\|\theta\|^8\|\theta\|_3^6$.  
\end{itemize}
Therefore, under both hypotheses, except for $(m_1,m_2,m_3,m_4)=(4,0,0,0)$, 
\beq \label{lastlemma-10}
f(\theta; m_1, m_2, m_3, m_4)\leq C(\|\theta\|^8+\|\theta\|^{12}+\|\theta\|^8\|\theta\|^{6}_3)\leq C\|\theta\|^{12}. 
\eeq
Consider two cases for $Y$. The first case is $N_{\widetilde{\Omega}}+N_{\epsilon}^{(1)}+N_{\epsilon}^{(2)}\neq 4$. Combining 
\eqref{lastlemma-9}-\eqref{lastlemma-10} gives
\[
\mathbb{E}[Y^2] = o(\|\theta\|^{-4})\cdot C\|\theta\|^{12} = o(\|\theta\|^8). 
\]
The claims follow immediately. The second case is $N_{\widetilde{\Omega}}+N_{\epsilon}^{(1)}+N_{\epsilon}^{(2)}= 4$. In this case, 
\[
f(\theta; m_1, m_2, m_3, m_4) = f(\theta; 4,0,0,0)=\|\theta\|^{16}. 
\]
If $N_{\epsilon}^{(1)} + N_{\epsilon}^{(2)}\geq 2$, then by \eqref{lastlemma-key} and \eqref{lastlemma-2order},
\[
\mathbb{E}[Y^2]=o(\|\theta\|^{-8})\cdot C\|\theta\|^{16} = o(\|\theta\|^8). 
\]
The claims follow. It remains to consider $N_{\epsilon}^{(1)} + N_{\epsilon}^{(2)}= 1$ (and so $N_{\widetilde{\Omega}}=3$). Write for short $S=1-\frac{v}{V}$. By \eqref{lastlemma-6}, 
\[
Y=S^{N^{(2)}_{\epsilon}}\cdot X, \qquad \mbox{where}\quad X=\sum_{i,j,k,\ell (dist)}a_{ij}b_{jk}c_{k\ell}d_{\ell i}, 
\]
and $a_{ij}, b_{ij}, c_{ij}, d_{ij}$ can only take values from $\{ \widetilde{\Omega}_{ij}, \epsilon^{(1)}_{ij}, \eta_i\eta_j\}$. So, $X$ is a non-stochastic number. Using \eqref{lastlemma-nonrandom}, we can easily show
\[
|X| \leq C\alpha^{N_{\widetilde{\Omega}}} \Bigl( \frac{\theta_{\max}}{\|\theta\|_1} \Bigr)^{N^{(1)}_\epsilon} \|\theta\|^8.  
\]
When $(N_{\epsilon}^{(1)}, N_{\epsilon}^{(2)})=(1,0)$, we have $Y=X$. By \eqref{lastlemma-2order}, $\frac{\theta_{\max}}{\|\theta\|_1}=o(\|\theta\|^{-2})$. It follows that 
\[
\mathrm{Var}(Y)=0, \qquad  |\mathbb{E}[Y]| =|X|\leq C \alpha^3 \cdot o(\|\theta\|^{-2}) \cdot \|\theta\|^8=o(\alpha^4\|\theta\|^8).  
\]
This gives the desired claims. 
When $(N_{\epsilon}^{(1)}, N_{\epsilon}^{(2)})=(0,1)$, we have $Y=S\cdot X$. So, 
\[
|Y|= |X| \cdot |S|\leq C\alpha^3\|\theta\|^8\cdot |S|. 
\]
Note that $S=1-\frac{v}{V}$, where $v=\mathbb{E}[V]$. Using the tail bound \eqref{lem-event-tail}, we can prove $\mathbb{E}[S^2]\leq C\|\theta\|_1^{-2}$. Therefore,
\[
\mathbb{E}[Y^2]\leq \frac{C\alpha^6\|\theta\|^{16}}{\|\theta\|_1^2}\leq C\alpha^6\|\theta\|^8\|\theta\|_3^6,
\]
where the last inequality is due to $\|\theta\|^4\leq \|\theta\|_1\|\theta\|_3^3$ (Cauchy-Schwarz). The claims follow immediately. \qed

\end{document}